

\magnification = 1200
\input amssym.def
\input amssym.tex
\def \qed {\hfill $\square$}
\def \R {\Bbb R}
\def \Z {\Bbb Z}

\def \S {{\cal S}}
\def \H {{\cal H}}
\def \ms {\medskip}
\def \ssi {\smallskip\noindent}
\def \pari {\par\noindent}
\def \msi {\medskip\noindent}
\def \sm {\setminus}
\def \wt {\widetilde}
\def \wh {\widehat}

\def \dsp {\displaystyle}
\overfullrule=0pt
\def \d {\partial}
\def \diam{\mathop{\rm diam}\nolimits}
\def \dist{\mathop{\rm dist}\nolimits}

\def \Min{\mathop{\rm Min}\nolimits}
\def \Max{\mathop{\rm Max}\nolimits}

\centerline{\bf LOCAL REGULARITY PROPERTIES }

\centerline{\bf OF ALMOST- AND QUASIMINIMAL SETS}

\centerline{\bf WITH A SLIDING BOUNDARY CONDITION}

\vskip 0.5cm
\centerline{ Guy David
\footnote{$^1$}
{Univ Paris-Sud, Laboratoire de Math\'{e}matiques, UMR 8658 Orsay, F-91405}
\footnote{$^2$}
{CNRS, Orsay, F-91405}
\footnote{$^3$}
{Institut Universitaire de France}
}
\vskip 1cm

\noindent
{\bf R\'{e}sum\'{e}.}
On s'int\'{e}resse \`{a} la r\'{e}gularit\'{e} jusqu'\`{a} la
fronti\`{e}re des ensembles presque minimaux et quasiminimaux
sous une condition de glissement. Les comp\'{e}titeurs d'un ensemble
$E$ y sont de la forme $F = \varphi_1(E)$, o\`{u}
$\{ \varphi_t \}$ est une famille \`{a} un param\`{e}tre d'applications
continues d\'{e}finies sur $E$, et qui pr\'{e}servent des ensembles 
fronti\`{e}res donn\'{e}s \`{a} l'avance. On g\'{e}n\'{e}ralise 
des r\'{e}sultats connus \`{a} l'int\'{e}rieur, et on d\'{e}montre notamment
l'Ahlfors r\'{e}gularit\'{e}, la rectifiabilit\'{e} et parfois 
l'uniforme rectifiabilit\'{e} locales des ensembles quasiminimaux, 
la stabilit\'{e} des classes consid\'{e}r\'{e}es par limites, 
et la presque monotonie de la densit\'{e} des ensembles presque minimaux 
sur des boules centr\'{e}es \`{a} la fronti\`{e}re.

\bigskip \noindent
{\bf Abstract.}
We study the boundary regularity of almost minimal and 
quasiminimal sets that satisfy sliding boundary conditions.
The competitors of a set $E$ are defined as $F = \varphi_1(E)$,
where $\{ \varphi_t \}$ is a one parameter family of continuous
mappings defined on $E$, and that preserve a given collection of 
boundary pieces. We generalize known interior regularity
results, and in particular we show that the quasiminimal
sets are locally Ahlfors-regular, rectifiable, and some times 
uniformly rectifiable, that our classes are stable under limits,
and that for almost minimal sets the density of Hausdorff measure in 
balls centered on the boundary is almost nondecreasing.

\medskip \noindent
{\bf AMS classification.}
49K99, 49Q20.
\medskip \noindent
{\bf Key words.}
Minimal sets, Almost minimal sets, Almgren restricted or quasiminimal sets, 
Sliding boundary condition, Hausdorff measure.
 
\vskip 0.5cm

\bigskip
\centerline{CONTENTS}
\msi
\centerline{PART I : INTRODUCTION AND DEFINITIONS}
\ssi
1. Introduction
\pari
2. Generalized sliding quasiminimal sets
\pari
3. Coral GSAQ and Lipschitz retractions on the $L_j$
\msi
\centerline{PART II : AHLFORS REGULARITY AND RECTIFIABILITY}
\ssi
4. Local Ahlfors regularity of quasiminimal sets 
\pari
5. Lipschitz mappings with big images, projections, and rectifiability
\msi
\centerline{PART III : UNIFORM RECTIFIABILITY OF QUASIMINIMAL SETS}
\ssi
6. Local uniform rectifiability in some cases
\par\hskip.3cm
6.a. How we want to proceed, following [DS4] 
\par\hskip.3cm
6.b. Some bad news
\par\hskip.3cm
6.c. What we can say anyway
\pari
7. The local uniform rectifiability of $E^\ast$ and bilateral 
weak geometric lemmas
\pari
8. Big projections and big pieces of Lipschitz graphs
\pari
9. Extension to the Lipschitz assumption
\msi
\centerline{PART IV : LIMITS OF QUASIMINIMAL SETS}
\ssi
10. Limits of quasiminimal sets: the main statement,
rectifiability, and l.s.c.
\pari
11. Construction of a stabler deformation:
the initial preparation
\pari
12. Step 2 of the construction: the places where $f$ is many-to-one
\par\hskip.3cm
Step 2.a. We cover $X_N(\delta_4)$ by balls $B_j = B(x_j,t)$, $j\in J_1$
\par\hskip.3cm
Step 2.b. We cover  $Y_N(\delta_4)$ by balls $D_l$, $l\in {\cal L}$
\par\hskip.3cm
Step 2.c. The collection of disks $Q$, $Q\in {\cal F}_l$
\par\hskip.3cm
Step 2.d. Where do the tangent planes go?
\par\hskip.3cm
Step 2.e. We choose disks $Q_j$, $j\in J_1$
\par\hskip.3cm
Step 2.f. We construct mappings $g_j$, $j\in J_1$
\par\hskip.3cm
Step 2.g. The mappings $g_j$, under the Lipschitz assumption
\par\hskip.3cm
Step 2.h. We glue the mappings $g_j$, $j\in J_1$,
and get a first mapping $g$
\pari
14. Step 3. Places where $f$ has a very contracting 
direction, and the $B_j$, $j\in J_2$
\pari
15. Step 4. The remaining main part of $X_0$
\pari
16. The modified function $g$, and a deformation for $E$
\pari
17. Magnetic projections onto skeletons,
and a deformation for the $E_k$
\par\hskip.3cm
17.a. Magnetic projections onto the faces
\par\hskip.3cm
17.b. A stable deformation for the $E_k$
\par\hskip.3cm
17.c. A last minute modification of our deformation
\pari
18. The final accounting and the proof of Theorem 10.8
in the rigid case
\pari
19. Proof of Theorem~10.8, and variants, under the Lipschitz assumption
\msi
\centerline{PART V : ALMOST MINIMAL SETS AND OTHER THEOREMS ABOUT 
LIMITS}
\ssi
20. Three notions of almost minimal sets
\pari
21. Limits of almost minimal sets and of minimizing sequences
\pari 
22. Upper semicontinuity of $H^d$ along sequences of almost 
minimal sets
\pari
23. Limits of quasiminimal and almost minimal sets in variable domains
\pari
24. Blow-up limits
\msi
\centerline{PART VI : OTHER NOTIONS OF QUASIMINIMALITY}
\ssi
25. Elliptic integrands; the main lower semicontinuity result
\pari
26. Limits of $f$-quasiminimal sets associated to elliptic integrands
\pari
27. Smooth competitors
\msi
\centerline{PART VII : MONOTONE DENSITY}
\ssi
28. Monotone density for minimizers
\pari
29. Minimal sets with constant density are cones
\pari
30. Nearly constant density and approximation by minimal cones
\pari
31. Where should we go now?

\bigskip
\centerline{FREQUENTLY USED NOTATION} 
\pari
$B(x,r) = \big\{ y \, ; \, |y-x| < r \big\}$
is the open ball centered at $x$ and with radius $r>0$.
\pari
$\H^d$ is the $d$-dimensional Hausdorff measure. 
See [Fe] or [Ma]. 
\pari
$GSAQ=GSAQ(U,M,\delta,h)$ is a class of quasiminimal sets;
see Definition 2.3.
\pari
$W_t = \big\{ y \in E \cap B \, ; \varphi_t(y) \neq y \big\}$
and $\widehat W = \bigcup_{0< t \leq 1} W_t \cup \varphi_t(W_t)$;
see (2.1).
\pari
$E^\ast = \big\{ x\in E \, ; \, \H^d(E\cap B(x,r)) > 0
\hbox{ for every } r>0 \big\}$ is the core of $E$; see (3.2).
\pari
$d_{x,r}(E,F)$ is almost a normalized Hausdorff distance in $B(x,r)$; see (10.5).
\pari
$\dagger \  \dagger$ delimits a proof or comment that concerns the Lipschitz 
assumption only.
\pari
$W_f = \big\{ x\in \R^n \, ; \, f(x) \neq x \big\}$; see (11.19).
\pari
$\wt f(x) = \psi(\lambda f(x))$ (used in Part IV, in the Lipschitz case); 
see (11.50), (12.75).  
\pari
$B_j = B(x_j,t)$, $j\in J_1$, is our first collection of balls (Part IV); see (12.8)-(12.9).
\pari
$B_j = B(x_j,r_j)$, $j\in J_2$, is the second collection of balls; see Lemma 14.6.
\pari
$D_j = B(y_j,r_j)$, $j\in J_3$, balls in the image, are used with the $B_{j,x}$; see (15.12)-(15.14).
\pari
$B_{j,x}$, $x\in Z(y_j)$, is our third collection of balls; see (15.19) and (15.1).
\pari
$h(r)$ is a gauge function that measures almost minimality; see (20.1) and Definition 20.2.
\pari
${\cal I}(U,a,b)$, ${\cal I}_l(U,a,b)$, and ${\cal I}^+(U,a,b)$ are classes of elliptic integrands; 
see Definition 25.3, Claim 25.89, and (25.94).

\bigskip
\centerline{INDEX} 
\pari
Ahlfors regularity: Prop.4.1, Prop. 4.74,
\pari
Almost minimal sets: ($A_+$, $A*$, $A'$): Def.20.2, 
Prop.20.9
\pari
BPBI (big pieces of bilipschitz images): Prop.6.41, Thm.7.7,
Prop.7.85, Prop.9.4
\pari
BPLG (big pieces of Lipschitz graphs): Thm.6.1, (8.4), 
Thm.8.5, Thm.9.81
\pari
Big projections: (8.3), Thm.8.5
\pari
BWGL (bilateral weak geometric lemma): Thm.7.7, Lem.7.8, Lem.7.72
\pari
Blow-up limits: (24.6), Thm.24.13, Cor.29.53
\pari
Carleson measures (7.5) and similar estimates in Sections 7, 8, 9
\pari
Concentration property [DMS]: Cor.8.55, Cor.9.103, 
Prop.10.82
\pari
Core, Coral set: Def.3.1, Prop.3.3, Prop.3.27
\pari
Density: (28.3), Thm.28.4, Thm.28.7, Thm.28.15, 
Thm.29.1, Cor.29.53, Prop.30.3, Prop.30.19
\pari
Dyadic cubes: see below (2.6)
\pari
Federer-Fleming projections: Lem.4.3, Lem.4.39, Prop.5.1, Thm.5.16
\pari
Flat configuration, flat faces: Def.24.8, Def.24.29, 
Prop.24.35
\pari
Gauge function: (20.1) and almost minimal sets
\pari
Integrands: Def.25.3, Thm.25.7, Claim 25.89, (25.94), 
Claim 26.4
\pari
Limits: (10.4)-(10.6), Thm.10.8, Thm.21.3, Cor.21.15,
Thm.23.8
\pari
Lipschitz assumption: Def.2.7
\pari
Lipschitz images: Prop.5.1
\pari
Lipschitz retractions: Lem.3.4, Lem.3.17, Lem.17.1, Lem.17.18
\pari
Lower semicontinuity of $\H^d$: Thm.10.97, Remark 23.23, Thm.25.7
\pari
Minimizing sequences: Remark 21.7, Cor.21.15
\pari
Monotonicity: Thm.28.4, Thm.28.7, Thm.28.15, 
Thm.29.1, Cor.29.53, Prop.30.3, Prop.30.19
\pari
Projections (surjective): Lem.7.38, Lem.9.14; Big projections: (8.3), Thm.8.5
\pari
Rectifiability: Thm.5.16, Lem.10.21
\pari
Rigid assumption: (2.6)
\pari
Technical assumptions (unpleasant) (6.2), (9.6), (9.83), (10.7), 
(19.36), Remark 19.52
\pari
Sliding competitor: Def.1.3
\pari
Smoother competitors: Section 27
\pari
Stopping time argument (painful): Section 6
\pari
Uniform rectifiability: Thm.6.1, Section 7, Thm.9.81
\pari
Upper semicontinuity of $\H^d$: Thm.22.1, Lem.22.3

\bigskip
\centerline{PART I : INTRODUCTION AND DEFINITIONS}
\ms
\noindent {\bf 1. Introduction}
\medskip

The main purpose of this paper is to study the boundary regularity properties
of minimal, almost minimal, and quasiminimal sets, subject to sliding boundary 
conditions that we will explain soon.

A long term motivation is to study various types of Plateau problems,
but where the objects under scrutiny are a priori just sets (rather than 
currents or varifolds), and we want to assume as little structure on them as 
possible. In this respect, the sliding conditions below seem natural to 
the author, and should be flexible enough to allow for a variety 
of applications.

Let us give a very simple example of a Plateau problem that we may 
want to study, and for which we do not have an existence result yet.
Let $\Gamma \i \R^n$ be a smooth closed curve, and let $E_0 \i \R^n$
be a compact set that contains $\Gamma$. For instance,
parameterize $\Gamma$ by the unit circle, extend the parameterization
to the closed unit disk, and let $E_0$ be the image of the disk.
Many other examples are possible, but with this one we should not get
a trivial problem for which the infimum is zero.
Our Plateau problem consists in minimizing $\H^2(E)$ among all
sets $E$ that can be written $E = \varphi_1(E_0)$, where
$\{ \varphi_t \}$, $0 \leq t \leq 1$, is a continuous, one parameter 
family of continuous mappings from $E_0$ to $\R^n$, with
$\varphi_0(x) = x$ for $x\in E_0$ and $\varphi_t(x) \in \Gamma$
for $0 \leq t \leq 1$ when $x\in E_0 \cap \Gamma$.
Thus, along our deformation of $E_0$ by the $\varphi_t$,
we allow the points of $\Gamma$ to move, but only along $\Gamma$;
this is why we shall use the term ``sliding boundary condition".

Minimizers of this problem, if they exist, will be among our simplest 
examples of minimal sets with a sliding boundary condition. 
But solutions of other types of Plateau problems (Reifenberg 
minimizers as in [R1,2], 
[De],  
or [Fa], 
or size minimizing currents under the boundary constraint
$\d T = G$, where $G$ denotes the current of integration along 
$\Gamma$, when they exist, also yield  minimal sets with a sliding 
boundary condition. 
Thus regularity results for sliding minimal sets may be useful
for a variety of problems, and we can also hope that they may help
with existence results.

Let us first give some definitions, and then discuss these issues 
a little more. The sets that we want to study are variants of the
Almgren minimal, almost minimal, or quasiminimal sets 
(he said ``restricted sets"), as in [A2], 
but where we add boundary constraints and are interested
in the behavior of these sets near the boundary.

We work in a closed region $\Omega$ of $\R^n$, which may also
be $\R^n$ itself, and we give ourselves a finite collections of 
closed sets $L_j \i \Omega$, $0 \leq j \leq j_{ max}$, that we call 
boundary pieces. It will make our notation easier to consider $\Omega$
as our first boundary piece, i.e., set
$$
L_0 = \Omega.
\leqno (1.1)
$$
For the elementary Plateau problem suggested above, for instance,
we would work with $L_0 = \Omega = \R^n$ and $L_1 = \Gamma$.

We are also given an integer dimension $d$, with
$0 \leq d \leq n-1$, and we consider closed sets $E \i \Omega$,
whose $d$-dimensional Hausdorff measure is locally finite, i.e., 
such that
$$
\H^d(E \cap B(x,r)) < +\infty
\leqno (1.2)
$$
for $x\in \Omega$ and $r>0$. The next definition explains 
what we mean by a deformation of $E$ that preserves the boundary 
pieces.

\ms\proclaim Definition 1.3. Let $B = \overline B(y,r)$ 
be a closed ball in $\R^n$. We say that the closed set $F \i \Omega$ is 
a \underbar{competitor} for $E$ in $B$, with sliding conditions given by the closed sets
$L_j$, $0 \leq j \leq j_{ max}$, when $F = \varphi_1(E)$ for some 
one-parameter family of functions $\varphi_t$, $0 \leq t \leq 1$, 
with the following properties:
$$
(t,x) \to \varphi_t(x) \hbox{ is a continuous mapping from
$[0,1] \times E$ to $\R^n$,}
\leqno (1.4)
$$
$$
\varphi_t(x) = x \hbox{ for $t=0$ and for $x\in E \sm B$,}
\leqno (1.5)
$$
$$
\varphi_t(x) \in B \hbox{ for $x\in E \cap B$ and $t\in [0,1]$,}
\leqno (1.6)
$$
and, for $0 \leq j \leq j_{ max}$, 
$$
\varphi_t(x) \in L_j \hbox{ when $t\in [0,1]$ and $x\in E \cap L_j \cap B$.}
\leqno (1.7)
$$
We also require that
$$
\varphi_1 \hbox{ be Lipschitz,}
\leqno (1.8)
$$
but with no bounds required.

\ms
We shall sometimes say ``sliding competitor in $B$" instead of 
``competitor for $E$ in $B$, with sliding conditions given by the $L_j$, 
$0 \leq j  \leq j_{ max}$", especially when our choice of $\Omega$ and the 
list of $L_j$ are clear from the context.

We shall soon discuss minimality, almost minimality, and quasiminimality
relative to this notion of sliding competitors, but since the class
of competitors is often the most important part of the definitions, 
a number of general comments on Definition 1.3 will be helpful.

It is important here that $\varphi_1$ is allowed not to be injective.
So we are allowed to merge different portions of $E$, or contract them
to a point, or pinch them in some other way. This, together with the 
fact that we shall not count measure with multiplicity, is why the 
union of two parallel disks that lie close to each other will not be 
minimal.

We added the last requirement (1.8) because Almgren 
put it in his definitions, and because this will not disturb.
If we drop it, we get more competitors for $E$, which means that the
almost- and quasiminimality properties are harder to get. Hence the 
regularity results proved here are also valid in the context where we 
drop (1.8). On the other hand, (1.8) will often be easy to prove, so 
it does not bother us much. The author suspects that the reason why Almgren
added (1.8) may be the following. Suppose you want to show that the support 
of a size minimizing current $T$ is a minimal set and, to simplify the 
discussion, that you are proceeding locally, in the complement of the 
boundary sets. You are given a deformation $\{ \varphi_t \}$ as in
Definition~1.3, and of course the simplest way to use it is to show 
that pushing $T$ by the $\varphi_t$, and in particular $\varphi_1$,
defines an acceptable competitor for $T$ (with the same boundary 
constraints). The constraint (1.8) just makes it possible to define 
the pushforward of $T$ by $\varphi_1$, so it is convenient.
See [D8] 
for details on this argument and its extension to the boundary. 

In the other direction, J. Harrison and H. Pugh once asked wether requiring 
$\varphi_1$, or even all the $\varphi_t$, to be smooth, would lead to 
the same classes of almost- and quasiminimal sets. The question 
was raised in the local context with no boundaries, but it also
makes sense in the present context.
The answer is yes under suitable conditions on the $L_j$, and if smooth
means $C^1$. For higher regularity, a proof seems to be manageable, 
but quite ugly, and so we only give a very rough sketch of how we would 
proceed, using the construction of Part~IV. This is discussed in Section 27.  


We are allowed to take $\Omega = \R^n$, and then (1.7) for $j=0$
is just empty and if there is no other boundary piece we get a minor 
variant of Almgren's definition of competitors in $\R^n$. 
Of course we can still restrict the list of competitors like he did, 
by requiring that $B$ lies in a fixed open set $U$, 
or that its diameter be less than some $\delta > 0$; 
we shall do this when we discuss our classes of almost- and quasiminimal 
sets, but let us not worry for the moment.

The main difference with Almgren's definition comes from the sliding boundary 
constraint (1.7), and this is also why we insist on the fact that $\varphi_1$ is 
the endpoint of a continuous deformation. If we did not require (1.7),
and we were given a continuous mapping $\varphi_1$ such that
$\varphi_1(x) = x$ for $x\in E\sm B$ and $\varphi_1(x) \in B$ for 
$x\in B$, we could define the $\varphi_t$ by 
$\varphi_t(x) = t \varphi_1(x) + (1-t) x$, and it is easy to check
that (1.4)-(1.6) would hold (because $B$ is convex).
We could also extend $\varphi_1$ to $\R^n$, which fits with the fact 
that $\varphi_1$ is traditionally defined on $\R^n$, not just on $E$.
But in the present situation we want points of the boundary $L_j$
to stay in $L_j$ (hence, (1.7)), and then it seems natural to say that the 
deformation condition in (1.7) only concerns points of $E$:
we do not want to say that the air besides our soap film $E$
is also concerned by the sliding boundary constraint. Notice that the
$\varphi_t$ can be extended to $\R^n$ (but in a way that may not 
preserve the $L_j$), so we do not have to worry about the case where our 
deformations would yield a tearing apart (cavitation) of the air 
besides the soap film.

Notice that with our convention that $L_0 = \Omega$, the set 
$\varphi_t(E)$ stays in $\Omega$, i.e.,
$$
\varphi_t(x) \in \Omega \hbox{ for $x\in E$ and  $t\in [0,1]$,}
\leqno (1.9)
$$
either because $x\in E \sm B$ and $\varphi_t(x) = x \in E \i \Omega$
by (1.5), or else by (1.1) and (1.7) with $j=0$.

The author thinks that Definition 1.3 is a nice way to encode boundary 
constraints, for instance that would be satisfied when $E$ is a soap 
film in a domain. A Plateau boundary constraint could for instance be
associated to one or a few curves $L_j$, but we could also think about
$L_1 = \d\Omega$ (or some other surface) as being a boundary along which the 
soap film may slide (as if loosely attached to a wall). It is quite 
probable that such boundary conditions were studied in the past, 
but the author does not know where. 

Once we have a notion of competitors, we can define a corresponding 
notion of minimal sets. Let us say, for the moment, 
that the closed set $E \i \Omega$ is minimal, 
with the sliding boundary conditions defined by the $L_j$,
$0 \leq j \leq j_{max}$, if $H^d(E) < + \infty$ and
$$
\H^d(E) \leq \H^d(F)
\ \hbox{ whenever $F$ is a sliding competitor for $E$ in some ball $B$,}
\leqno (1.10)
$$
where we allow $B$ to depend on $F$. Many variants of this definition
will be proposed, where one may localize the definition to an open set
$U$, or add a small error term to the right-hand side in (1.10)
(this is how we will define almost minimal sets), or even allow
stronger distortions (this will give rise to quasiminimal sets). 
We shall give the main definitions in Section 2 (for the generalized 
quasiminimal sets) and later in Section 20 (for almost minimal sets),
but for the moment the sliding minimal sets that satisfy (1.10) will
give a fair idea of what we want to study.

Of course our notion of competitors can be used to define Plateau 
problems, as we did earlier with a single curve. Given a collection
of boundary pieces $L_j$, and a closed set $E_0$, we can try to 
minimize $\H^d(E)$ among all the sets $E$ that are sliding 
competitors of $E_0$ (in some ball $B$ that depends on $E$, 
or in some fixed huge ball that contains $\Omega$). 
If $E_0$ is badly chosen (for instance, if some sliding competitors of
$E_0$ are reduced to a point), the problem may not be interesting,
but it is easy to produce lots of examples where the infimum
will be finite and positive. For most of these examples,
we do not have an existence result. But it is clear that
if minimizers for this Plateau problem exist, they are sliding 
minimal sets.

The main point of this paper is to study the general 
(hence often rather weak) regularity properties of the
minimal sets, and their almost minimal and quasiminimal variants,
in particular when we approach the boundary pieces $L_j$.
In practical terms, this means that we will take many interior
regularity results for Almgren minimal (or quasiminimal) sets,
and try to adapt their proofs so that they work all the way
to the boundary. But before we say more about this,  
let us comment a little more on Definition 1.3 and
our motivations.

The word sliding may be misleading in some cases,
as some sets $L_j$ may be reduced to points, where in effect no 
sliding will be allowed. Our assumptions on the $L_j$ will only allow 
a finite number of points where $E$ is fixed. So, for 
instance, we do not consider the case where $\Gamma$ is a simple
curve and we require that $\varphi_t(x) = x$ for every point 
$x\in E\cap \Gamma$. This will not bother us, and probably
such a condition would make it too hard to produce competitors 
and get information on $E$ near $\Gamma$ when $E$ is a
minimal set with these constraints. 
Of course we could always say that $E$ is locally
minimal (for instance) in the domain $U = \R^n \sm \Gamma$, and get
some information from this, but this is not the point of this paper.
On the contrary, the author believes that because we allow our 
competitors to slide along the $L_j$, we will have an amount of
flexibility in the construction of competitors, which we can use 
to prove some decent regularity results. And at the same time
(1.7) looks like a reasonable constraint, for instance, if we want to
model the behavior of soap films.

We believe that in addition to being interesting by themselves,
regularity results for sliding minimal or almost sets could be useful 
to prove existence results (in very simple cases)
for the Plateau problems discussed above, and also 
for other similar problems, because some other types 
of minimizers also yield sliding minimal sets. 
Let us give two examples.

In [R1],  
Reifenberg proposed a Plateau problem where we are given a compact boundary
set $L \i \R^n$ of dimension $d-1$, and we minimize $\H^d(E)$ among compact sets $E$ 
that bound $L$, in the sense that $L \i E$ and the natural map induced by the inclusion, 
from the $(d-1)$-dimensional \v{C}ech homology group of $L$ to the 
$(d-1)$-dimensional \v{C}ech homology group of $E$, is trivial.
He also proves a fairly general existence result, and good interior 
regularity results for the minimizers 
(see [R1,2]). 
These results were generalized by various authors; see for instance
[A1], [De], and more recently [Fa] 
for a quite general existence result. 
Also see [HP] 
for a simpler variant of [R1] in codimension $1$,  
where one replaces the computation of \v{C}ech homology groups
with a simpler linking condition, and which comes with a simpler
proof and is related to differential chains.

It is easy to see that if the boundary set $L$ is not too ugly,
the minimizing sets that are obtained in these papers are
sliding minimal sets associated to $L_0 = \R^n$ and $L_1 = L$. 
See [D8] 
for the rather easy verification, whose main point is
just that if $E$ bounds $L$ and $F$ is a sliding competitor for $E$, 
then $F$ bounds $L$ too.

Reifenberg' homological Plateau problem and its minimizers
are very nice, and give good descriptions of many soap films,
but some people prefer the related problem of size minimizers.
That is, we are given a $(d-1)$-dimensional integral current $S$, 
with $\d S = 0$, and we look for a $d$-dimensional integral current $T$ 
such that $\d T = S$ and whose size (understand, the $H^d$-measure
of the set where the multiplicity is nonzero, but we shall be slightly
sloppy on the definitions) is minimal.
If $d=2$, $L$ is a nice closed curve in $\R^3$, and $S$ is the
current of integration on $L$, 
T. De Pauw showed in [De] 
that the infimum for this problem is the same as for Reifenberg's
homological problem (where \v{C}ech homology is computed over the 
group $\Bbb Z$); but even though De Pauw showed that Reifenberg 
homological minimizers exist,
size minimizers are not known yet to exist.
Anyway, size minimizers, if they exist, are also supported
(under reasonable conditions) on sliding minimal sets. The
point now is that if $T$ is supported by the closed set $E$
and $F$ is a sliding competitor for $E$, then we can use $\varphi_1$
to push $T$ and get another solution of $\d T = S$, which is
supported on $F$. See [D8] 
for the fairly easy verification of this, and variants where
$\d T$ is only required to be homological to $S$ in a boundary set $L$.

So we have at least two potentially interesting other examples of 
sliding minimal sets. To the author's knowledge, not much is known 
on the boundary behavior of these sets, 
and the results in this paper are probably a good start.
A natural question is whether, if we decide to study them by saying
they are sliding minimal sets and forgetting about the initial problem
they solve, we lose important information that we may have used 
profitably.

Most of the results of this paper concern the weaker notion of
sliding quasiminimal sets, but let us make two short remarks
on sliding minimal sets. An important tool that we can still use
in some cases is Allard's regularity theorem from [All], 
which applies to more general stationary varifolds and
goes all the way to the boundary. But this result uses some
initial flatness assumption that we may not want to assume.

In the special case when $L_0 = \R^3$ and $L_1$ is a
nice curve, G Lawler and F. Morgan propose a conjectural list
of 10 boundary behaviors for minimal sets 
bounded by $L_1$; see [LM] and [Mo3], 
and in particular Figure 13.9.3 (on page 137 of the third edition). 
The present paper tries to go in such directions, so far
in more general contexts but with less precision.

\ms
We may now start a description of the results in this paper.
Generally speaking, we shall take local regularity results 
that we like, and try to extend or modify the proofs so that they
work also near the boundary pieces. 

Many of our results are about what we call generalized 
quasiminimal sets, which are
defined in Section 2 (see Definition~2.3). In the special 
case without boundary pieces, the notion is just a little bit more 
general than the quasiminimal sets that 
Almgren studied in [A2] 
under the name of ``restricted sets".
One advantage of the notion is that it is rather weak and quite
flexible. For instance, it is stable under bilipschitz mappings
(the quasiminimality constant $M$ just gets larger), and 
contains minimizers of functionals like $\int_{E} f(x) d\H^d(x)$,
where we just need to know that $f$ is bounded and bounded from below,
and under the same sort of boundary conditions as above.
Thus the graph of any Lipschitz function 
$F : \R^d \to \R^{n-d}$ is locally quasiminimal (with no
boundary condition). Of course this means that we cannot expect 
better regularity than Lipschitz, but this will already be a good 
start, and in effect we shall not get so far from that.

We shall work locally, in an open set $U$, and 
with two set of assumptions on the boundary pieces.
In the first one, which we shall call the rigid assumption,
$U$ is the unit ball $B_0$, we choose a dyadic grid of $\R^n$, 
and we require all the sets $L_j$ to coincide in $U$ with a finite union 
of faces of cubes of our grid. We do not even require all these faces to 
be of the same dimension.

This already gives some choice, but we do not necessarily want
all the faces to be smooth, and we expect some bilipschitz invariance,
so we also allow a weaker set of assumptions, 
which we call the Lipschitz assumption, where
$U$ and the $L_j$ are obtained from the previous case by composing
with a bilipschitz mapping from $B_0$ to $U$.
We even allow an additional dilation that we shall skip here for 
simplicity. See Definition 2.7.
Some times the regularity results in this second case will require
more complicated proofs, but we decided to include them anyway.

Even this set of assumption is not entirely satisfactory, because
for instance it puts some small bounds on the number of faces that 
may touch a given point, but the dyadic combinatorics are pleasant to use, 
and the author was afraid of the complications that may arise in a more 
general case.

Let us give a rough description of our plan.

\ms
Part I deals with the setup and definitions.
After the definitions of Section~2, we check that sliding quasiminimal 
sets in $U$, under the Lipschitz assumption, are just the images by 
our bilipschitz parameterizaton of sliding quasiminimal sets in $B_0$, 
with the rigid assumption.
See Proposition~2.8. 

In Section 3, we introduce the core $E^\ast$ of a closed set $E$ 
(our name for the closed support of the restriction of $\H^d$ to $E$, 
see Definition 3.1), and show that the core of a 
(generalized sliding) quasiminimal set is quasiminimal with the same 
constants. See Proposition 3.3 (and before, Proposition 3.27 in the
simpler rigid setting). The proof is a little unpleasant (because our 
boundary constraint (1.7) does not obviously cooperate with removing 
some parts of $E$), but afterwards we feel better because we can 
forget about the fuzzy set $E \sm E^\ast$, and restrict our attention
to coral sets, i.e., sets such that $E^\ast = E$.

\ms
Part II contains our first regularity results for 
generalized sliding quasiminimal set. 

In Section 4, we show that the core $E^\ast$ of such a set $E$ 
is locally Ahlfors regular. This means that, if $x\in E^\ast$, 
$B(x,2r) \i U$ (the open set where we work),
$r > 0$ is smaller than the scale constant $\delta$ in the definition 2.3
of quasiminimal sets, and the small parameter $h>0$ in Definition 2.3 
is small enough, then
$$
C^{-1} r^d \leq \H^d(E\cap B(x,r)) 
= \H^d(E^\ast\cap B(x,r)) \leq C r^d.
\leqno (1.11)
$$
See Proposition 4.1 (under the rigid assumption) and Proposition 4.74
(for the Lipschitz case). The proof relies on comparison arguments
based on Federer-Fleming projections. It follows the proof
of [DS4] (for the case without boundaries), 
which itself looks a lot like the proof in [A2] 
of almost the same result.  

Section 5 continues along the same lines. Its main result is 
Theorem 5.16, which says that quasiminimal sets (with a small enough
constant $h$) are rectifiable. We still prove this with a
Federer-Fleming projection, and the proof is probably similar to 
Almgren's original proof (away from the boundaries).
The main point is that near a point of density of the unrectifiable
part of $E$, we could project $E$ on a small subset of $d$-faces, 
so small that an additional projection on faces of dimension $d-1$ 
is possible and allows us to make it essentially disappear. 
It is interesting that the rectifiability
of quasiminimal sets (and their limits, see Part IV) was neglected in
[DS4], just because we could prove stronger properties, 
while here we will have to rely more on it in the cases where
we don't get uniform rectifiability.

On a slightly more technical level, Proposition 5.1 says that
for $B(x,r)$ as above (i.e., $E$ is quasiminimal, $x\in E^\ast$, 
$B(x,2r) \i U$, $r < \delta$, and $h$ is small enough),
there is a Lipschitz mapping $F : E \cap B(x,r) \to \R^d$
such that $\H^d(F(E \cap B(x,r))) \geq C^{-1} r^d$; this is a 
technical lemma that can be used in later proofs (typically, for uniform  
rectifiability). Then Proposition 5.7 is a trick from [DS4]  
that allows us to pretend that $F$ is the orthogonal projection
on some $d$-plane. 

\ms
Part III deals with the local uniform rectifiability of the core $E^\ast$
when $E$ is quasiminimal (and $h$ is small enough, as always).

The main result of this part says that if $E$ is a quasiminimal
set (and $h$ is small enough), and if some technical condition
on the dimension of the faces is satisfied, then $E^\ast$ is locally
uniformly rectifiable, with big pieces of Lipschitz graphs.
See Theorem 6.1 under the rigid assumption, and
Theorem 9.81 under the Lipschitz assumption.

Uniformly rectifiable with big pieces of Lipschitz graphs means
that there are constants $A \geq 0$ and $\theta > 0$ such that,
if $B(x,r)$ is as above, then we can find an $A$-Lipschitz graph
$\Gamma$ of dimension $d$ such that
$$
\H^d(E \cap B(x,r)) \cap \Gamma \geq \theta r^d.
\leqno (1.12)
$$
Thus in $B(x,r)$, a substantial part of $E$ lies in the nice
$A$-Lipschitz graph $\Gamma$. Recall that by definition, $\Gamma$ is
the graph of some $A$-Lipschitz function that is defined on some
$d$-plane $P \i \R^n$, and with values in the orthogonal $(n-d)$-plane
$P^\perp$. In this statement, $A$ and $\theta$ depend only on 
the dimensions $n$ and $d$, the quasiminimality constant $M$
in Definition~2.3, and the bilipschitz constant $\Lambda$ in 
Definition 2.7.

Unfortunately, we we only get this under the technical condition (6.2)
(or its analogue (9.83) when we use the Lipschitz assumption).
It is satisfied if, except for the supporting domain $L_0 = \Omega$,
all the boundary pieces $L_j$ are composed of faces of dimensions at 
most $d$. This takes care of many interesting examples, but it is 
nonetheless frustrating that we have to assume this. Of course
we do not have a counterexample; the main problem could even be that 
even in the case without boundary, we have only one proof of uniform 
rectifiability, and this proof is complicated and fails badly when
we deal with boundaries.

The positive point of uniform rectifiability is that it has the right
invariance under bilipschitz mappings, and that it is, to the author's
knowledge, the best very general (weak) regularity result for our 
quasiminimal sets.

Most of Part III is devoted to a proof of Theorems 6.1 and 9.81
on the local uniform rectifiability of $E^\ast$. We essentially take
the long and complicated proof from [DS4],   
try to adapt it, and see where it fails.

At the start, Propositions 5.1 and 5.7 allow us to assume that
for some orthogonal projection $\pi$ on some $d$-plane, 
$\H^d(\pi(E \cap B(x,r))) \geq C^{-1} r^{d}$; the whole proof then 
consists in showing that we can find a large subset of $E\cap B(x,r)$ 
where $\pi$ is bilipschitz. 
Section~6 describes the general scheme of a stopping time
argument which is designed to select the large subset, why it fails in 
general, and why it still works in some limited cases (but really,
not so many new things happen, compared to the previous situation
with no boundary).
We end up, in Proposition 6.41, with a result that says that in
some cases, $E \cap B(x,r)$ contains a significant part which is
bilipschitz-equivalent to a subset of $\R^d$.

In addition to  the stopping time argument described in Section 6,
Theorems 6.1 and 9.81 use some amount of general uniform 
rectifiability theory which is done, when we work under the rigid
assumption, in Sections 7 and 8.

The uniform rectifiability of an Ahlfors regular set $E$ can be defined in lots
of (eventually equivalent) ways, and in Section 7 we discuss two of them. 
The first one, called BPBI, asks for the existence, in each ball 
$B(x,r)$ centered on $E$, of a substantial part of $E \cap B(x,t)$
that can be send to a subset of $\R^d$ by a bilipschitz mapping.
In the case of quasiminimal sets, we first restrict to the core $E^\ast$ 
and work only locally, i.e., on balls such that $B(x,2r) \i U$, but
let us forget these details.
Now Proposition 6.41 gives something like this, but not in enough 
balls $B(x,r)$, so one has to work more, and in effect go through
the BWGL below.

A second definition of uniform rectifiability is by the bilateral 
weak geometric lemma (BWGL), which asks that for most balls $B(x,r)$
(defined in terms of Carleson measures but please don't mind),
there is a $d$-plane $P$ such that $E \cap B(x,r)$ is $\varepsilon r$-close 
to $P\cap B(x,r)$ (in Hausdorff distance, and where $\varepsilon > 0$ is 
a fixed small constant). It turns out that this one is easier to get.

The only place in Section 7 where the quasiminimality of $E$
is used directly (as opposed to, via a regularity result of a previous section)
is to show that if all the points of $E \cap B(x,r)$ lie within $\varepsilon r$
of some $d$-plane $P$, then the converse is also true: 
all the points of $P \cap B(x,3t/2)$ lie within $\varepsilon r$ of 
$E$, and in addition the orthogonal projection from $E$
to $P$ is locally surjective (see (7.46)).
See Lemma 7.38 for a more precise statement that takes into account
the position of the boundary pieces $L_j$, and Lemma 9.14
for a generalization of this first statement.

This lemma helps because it is relatively easy to find balls where 
$E$ stays close to a plane, but the BWGL requires a bilateral
approximation that Lemma 7.38 then provides. 
The rest of Section 7 consists in playing with bad sets of 
balls and various definitions of uniform
rectifiability, to get the BPBI property (for every small ball, not 
just the good ones in Proposition 6.41). See Proposition 7.85.
So $E$ is locally uniformly rectifiable.

In Section 8, we keep the rigid assumption and go from the BPBI to the BPLG, 
i.e., the existence of big pieces of Lipschitz graphs, as in the statement
of Theorem 6.1. For this, the general theory says that we have to 
find big projections (see Theorem 8.5) and, roughly speaking, this is 
provided by the BWGL or even its unilateral version the WGL,
plus another application of Lemma 7.38 (and (7.46) in particular).

In Section 9 we prove the analogue of Theorem 6.1
under the more general, but some times more painful, Lipschitz 
assumption. The relevant statements are now 
Lemma~9.14 (for the generalization of Lemma 7.38)
and Theorem 9.81 (for the main uniform rectifiability result).

A consequence of the uniform rectifiability of $E^\ast$, that has been 
quite useful for the study of limits far from the boundaries, is
the concentration property introduced by Dal Maso, Morel, and 
Solimini [DMS] in the context of the Mumford-Shah functional. 
The point is that for any sequence $\{ F_k \}$ of sets that satisfies
this property (with uniform constants) and converges to $F$
in Hausdorff distance, and any open set $V$,
we have the lower semicontinuity property
$\H^d(F\cap V) \leq \liminf_{k \to +\infty} \H^d(F_k \cap V)$.

We prove this property in Corollary 8.55
(under the Lipschitz assumption) and Corollary 9.103
(under the Lipschitz assumption), as simple consequences of
the local uniform rectifiability, but then with the additional
technical assumption (9.2) or (9.105). Fortunately, there is
another proof of uniform concentration along sequences that does not use 
these assumptions; see Proposition 10.82.

Most of this Part III is not needed for the next ones; the 
failure of Theorems 6.1 and 9.81 in some cases lead the author to 
finding ways to prove the subsequent theorems (and in particular the
results on limits, see Part V) that would not use uniform rectifiability.
So the reader will get something positive out of the weakness of this part.

\ms
Part IV contains our main results on the limits of quasiminimal sets.
The main statement for this part is Theorem 10.8, which says that if
$E$ is the Hausdorff limit (locally in the open set $U$) of the 
sequence $\{ E_k \}$ of coral (i.e., $E_k^\ast = E_k$)
quasiminimal sets which all lie in a class $GSAQ(U,M,\delta,h)$,
with $h$ small enough, then $E$ lies in the same
quasiminimal class $GSAQ(U,M,\delta,h)$ as the $E_k$. 

Here again, when we work under the Lipschitz assumption, we only
prove this under a minor additional regularity assumption on the
faces that compose the $L_j$. Typically, when such a face is more than 
$d$-dimensional, we require the face to be $C^1$ in a
neighborhood of $\H^d$-almost each of its interior points.
See (10.7), or Remark 19.52 for a weaker condition.

The main ingredient for Theorem 10.8 is the lower semicontinuity
estimate in Theorem~10.97, which says that for $\{ E_k \}$ as
above,
$$
\H^d(E\cap V) \leq \liminf_{k \to +\infty} \H^d(E_k \cap V)
\ \hbox{ for every open set } V \i U.
\leqno (1.13)
$$
This is deduced from Dal Maso, Morel, and Solimini's 
result [DMS] 
and the fact that the sets $E_k$ are uniformly concentrated,
as in Proposition 10.82. In turn Proposition 10.82 is obtained a 
little bit like Corollaries 8.55 and 9.103, but instead of uniform
rectifiability, we use the fact that the limit $E$ is rectifiable
(as in Proposition 10.15), and a compactness argument (Proposition 
10.21). The surprising part, at least to the author, 
is the rectifiability of the limit,
which is just proved like Theorem 5.16 (the rectifiability of 
a single $E_k$), with suitable modifications.

Even though Theorem~10.97 is the main ingredient in 
Theorem 10.8, the full proof takes the rest of Part IV
(Sections 11-19). It follows the argument of [D2], 
but unfortunately with many small modification that
force us to give a full proof.

Perhaps we should mention that it is important to prove
limiting results like Theorem~10.8 and Theorem~10.97 
in the context of sets. 
In the context of integral currents, for instance,
the lower semicontinuity of the mass and strong compactness
theorems exist, that have been used very profitably.
Here we get an acceptable substitute for some of that.
Without this, it would be hard to say much about the
blow-up limits of almost minimal sets, for instance.

\ms
In Part V we study the stronger notion of almost minimality,
and extend the stability results of the previous parts to them. 
A few different definitions are possible, 
but let us give a simple one that works when
we do not localize. In addition to the list of boundary pieces
$L_j$ (which we keep as above), we give ourselves a gauge function 
$h : (0,+\infty) \to [0,+\infty]$, such that
$\lim_{r \to 0} h(r) = 0$. Often we also ask $h$ to be nondecreasing
and continuous from the right, and for some results to have a 
sufficient decay near $0$. A typical choice would be to 
pick $\alpha > 0$ and take 
$h(r) = r^\alpha$ for $0 \leq r < \delta$
and $h(r) = +\infty$ for $r \geq \delta$.
A sliding almost minimal set (of type $A'$) in $\R^n$ is then a 
closed set $E$ such that (1.2) holds, and for which 
$$
\H^d(E \cap B) \leq \H^d(F \cap B) + h(r) r^d
\leqno (1.14)
$$
for each closed ball $B = \overline B(x,r)$ and each 
sliding competitor $F$ for $E$ in $B$.
When $h(r) \equiv 0$, we recover the definition of sliding minimal
sets defined by (1.10).
This notion can be localized to an open set $U$, and 
three slightly different types of sliding almost minimal sets 
(called $A_+$, $A$, $A'$) are introduced in Definition 20.2.
Of course we expect better regularity properties for
the sliding almost minimal sets, especially when $h$ is small;
here we shall not really look for such properties, but rather 
prepare the ground with some preliminary results on limits
of almost minimal sets and monotonicity properties for their density.

In Section 20 we give the three definitions of sliding almost minimal
sets (Definition~20.2), but then prove that the two last ones
($A$ and $A'$) are equivalent. This is Proposition 20.9;
the proof follows [D5], 
were similar notions were defined (to try to unify some 
definitions with Almgren's initial ones).

In section 21 we use our limiting theorem on quasiminimal
sets (Theorem 10.8) to show that limits of coral sliding almost minimal sets
(of a given type) with a given gauge function $h$ are also
coral sliding almost minimal sets, of the same type and with the same
gauge function. This is Theorem 21.3. Also see Remark 21.7 and 
Corollary 21.15 that say that the limit of a locally minimizing
sequence of uniformly quasiminimal sets is locally minimal.

In Section 22 we prove an upper semicontinuity result for $\H^d$:
if the sequence $\{ E_k \}$ of coral sliding almost minimal sets
in $U$ converges to $E$ (as in Theorem 21.3), then for each compact set $H$
in $U$,
$$
\H^d(E\cap H) \geq \limsup_{k \to +\infty} \H^d(E_k \cap H);
\leqno (1.15)
$$
see Theorem 22.1, which is specific to the case when $h(r)$ tends
$0$. For quasiminimal sets, we cannot expect such a neat estimate, 
but we still have the less precise 
$$
(1+Ch) M \H^d(E\cap H) \geq \limsup_{k \to +\infty} \H^d(E_k \cap H)
\leqno (1.16)
$$
which is proved in Lemma 22.3 and is often useful too.
Again similar results were proved in [D5] Lemma 13.12, 
and probably many more places before. Surprisingly, the proof only
uses the rectifiability of the limit $E$, some 
covering lemmas, and an application of the definition of 
quasiminimality in some flat balls.

Theorems 10.8 and 21.3 have an obvious defect: in many situations,
such as for blow-up limits with boundaries $L_j$ that are not cones,
we may want to take limits in situations where the domains, and more
importantly the boundary sets $L_j$, change mildly. 
We do this in Theorem 23.8, but rather than redoing the whole proof, 
we reduce to the previous statements by composing with a variable 
change of variables that sends us back to a fixed domain (the limit).
Our proof forces us to restrict to variable domains that are  
close to the limiting domain in the bilipschitz category, which
is probably not optimal.

We apply this in Section 24 to blow-up limits. Under reasonably mild 
flatness conditions on the sets $L_j$ at the origin
(see Definitions 24.8 and 24.29, and Proposition 24.35
that says that the individual flatness of faces (as in Definition 24.29)
is enough), we show that the blow-up limits at the origin of a
sliding almost minimal set, are sliding minimal sets in $\R^n$,
associated to boundary sets $L_j^0$ obtained from the $L_j$ by the
same blow-up. See Theorem 24.13.

\ms
Part VI deals with two extensions of our notions of quasiminimality
and almost minimality. The main one is related to elliptic integrands.
Instead of using the Hausdorff measure $\H^d(E)$ in our various 
definitions, we may want to use slightly distorted versions like
$\int_E f(x) d\H^d(x)$, where $f: \R^n \to [1, M]$ is a continuous 
function, or even 
$$
J_f(E) = \int_E f(x,T_x E) d\H^d(x),
\leqno (1.17)
$$
where $f$ is now defined on $\R^n \times G(n,d)$, $G(n,d)$ denotes the 
Grasssman manifold of vector $d$-planes in $\R^n$, 
and we restrict to rectifiable sets so that the approximate tangent plane
$T_x E$ is defined almost everywhere on $E$
(see (25.2) for a slightly artificial definition, but that would 
also work on $d$-sets that are not rectifiable).
See Definition 25.3 for an acceptable class of elliptic integrands, 
which is just a little larger than the one introduced by 
Almgren [A1], [A3]. 

The main point of Section 25 is that the technique 
of [DMS] also allows us to prove lower semicontinuity 
results like (1.13), but for integrals like $J_f(E)$.
This was noticed by Yangqin Fang [Fa], 
who wanted such a result to extend Reifenberg's
existence theorem for his homological Plateau problem
to the context of elliptic integrands,
and Fang's proof is so simple that it would have been stupid not 
to give it here.

In Theorem 25.7, we prove that if the sequence 
$\{ E_k \}$ of sliding quasiminimal sets in $U$
satisfies the main assumptions of our limiting
Theorem 10.8, and if the integrand $f : \R^n \times G(n,d) \to [a,b]$
satisfies the condition of Definition 25.3, then 
$$
J_f(E\cap V) \leq \liminf_{k \to +\infty} J_f(E_k \cap V)
\ \hbox{ for every open set } V \i U,
\leqno (1.18)
$$
where as usual $E$ is the limit of the $E_k$.

The proof contains the lower semicontinuity result
that we used for Theorems~10.8 and 21.3, 
so the reader that would not be familiar with [DMS] 
can read Section 25 instead and get a slightly more 
direct proof, even for $f \equiv 1$. We still kept the reference
to [DMS] for the other readers, and also because 
this is really where the ideas are coming from.

The notions of quasiminimality and almost minimality
can also be defined in terms of an elliptic integrand $f$
as above. Since $a \leq f \leq b$ for some $a, b > 0$,
the list of quasiminimal sets is the same, only the 
constants are different. This is why we do not need to be careful 
when we state Theorem 25.7. In Section 26 we explain how to extend
Theorem~10.8 to limits of $f$-quasiminimal sets;
see Claim 26.4. The same thing would happen with other
results of Part V, but we omit the details.

We included Section 27 to answer partially a question of J. Harrison
(initially raised far from the boundary), and to say that the question
is probably not as simple as it seems.
Suppose, in the definition of competitors (Definition 1.3),
that we only included competitors for which $\varphi_1$ is smooth;
would the resulting sets of quasiminimal (or almost minimal) 
sets be different? We discuss some partial positive results,
and a possible strategy for further ones, in Section 27.

\ms
Part VII deals with the monotonicity, or near monotonicity,
of the density
$$
\theta(r) = r^{-d} \H^d(E \cap B(x,r))
\leqno (1.19)
$$
for sliding minimal or almost minimal sets, but only
for balls $B(x,r)$ centered on the boundary pieces.

The simplest result is Theorem 28.4, which says that if $x=0$,
$E$ is coral and locally sliding minimal near $0$, and the $L_j$ are cones,
$\theta$ is nondecreasing near $r=0$. When instead 
$E$ is only almost minimal with a gauge function $h$ that
satisfies a Dini condition, and in addition $0 \in E$
(at least, if we deal with $A$-almost minimal sets),
Theorem 28.7 says that $\theta$ is nearly monotone,
i.e., that we can multiply it by a continuous function 
with a nonzero limit at the origin and get a nondecreasing function.

The case when the $L_j$ are not exactly cones centered at $x$ is 
discussed in Remark 28.11 and Theorem 28.15.

The case of equality in Theorem 28.4, i.e., when 
$E$ is minimal, the $L_j$ are cones, and $\theta$ is constant
on some interval, is treated in Section 29.
Theorem 29.1 says that in this case $E$ coincides, in the
corresponding annulus, with a minimal cone with the same sliding
boundary conditions. We use the proof of [D5], 
by lack of a better idea.

We apply this to blow-up limits of coral sliding almost minimal sets 
and show in Corollary 29.53 that, under reasonable assumptions, 
they are sliding minimal cones associated to the blow-up limits of the $L_j$.

We also use the case of equality above, and a compactness argument,
to find situations where, if the function $\theta$ is nearly constant
on an interval, then $E$ can be well approximated by a minimal
cone, both in terms of Hausdorff distance and measure.
See Proposition 30.3 for a general statement with annuli, and 
Proposition 30.19 for a simpler case in a ball.

In a last Section 31, we rapidly discuss a few directions in which 
this work could be continued or used.

The author wishes to thank Thierry De Pauw, Yangqin Fang, Vincent Feuvrier,
Jenny Harrison, Xiangyu Liang, Frank Morgan, and Harison Pugh 
for helpful discussions and remarks concerning this project.
He gladly acknowledges the generous support of
the Institut Universitaire de France, and of the ANR
(programme blanc GEOMETRYA, ANR-12-BS01-0014).

\bigskip
\noindent {\bf 2. Generalized sliding quasiminimal sets}
\medskip

In this section we give the definition of our most general class
of quasiminimizers (the sets for which we shall prove most of our
regularity results), and also describe the two standard sets 
of assumptions on the boundary pieces $L_j$ that will be allowed.

The following notion comes from [D5], 
where it was introduced to generalize both the notion of Almgren 
quasiminimal set (or ``restricted set", see [A2]) 
and some simpler notions of almost minimal sets.

For the next definition, we shall use a quasiminimality constant  
$M \geq 1$, a diameter $\delta\in(0,+\infty]$, and a small number
$h\in [0,1)$. We want to be able to localize our definitions, 
which forces us to work in an open set $U$; but of 
course we are free to take $U = \R^n$.

Given a closed set $E$, with $\H^d(E \cap H) < +\infty$
for every compact set $H \i U$, and a one-parameter family 
$\{ \varphi_t \}$, $0 \leq t \leq 1$, 
such that (1.4)-(1.8) hold, we set
$$
W_t = \big\{ y \in E \cap B \, ; \varphi_t(y) \neq y \big\}
\leqno (2.1)
$$
for $0 < t \leq 1$, and then
$$
\widehat W = \bigcup_{0< t \leq 1} W_t \cup \varphi_t(W_t).
\leqno (2.2)
$$
Note that $W_t \i \widehat W  \i B$, where $B$ is as in  (1.4)-(1.8), but  
they may be smaller and in particular we shall not force $B$ to be 
contained in $U$.

\medskip
\proclaim Definition 2.3. Let $\Omega \i \R^n$ and the 
$L_j$, $0 \leq j \leq j_{ max}$, be as above. Let $M \geq 1$, 
$\delta\in(0,+\infty]$, $h\in [0,1)$, and the open set $U$ be given.  
Let $E \i \Omega$ be a closed set in $U$ such 
$\H^d(E\cap B) < +\infty$ for every closed ball $B \i U$.
We say that $E \in GSAQ(U,M,\delta ,h)$ when, for every choice
of closed ball $B = \overline B(x,r)$ such that $0 < r < \delta$, 
and every one-parameter family $\{ \varphi_t \}$, $0 \leq t \leq 1$, 
such that (1.4)-(1.8) hold and 
$$
\widehat W \i\i  U
\leqno (2.4)
$$
(i.e., $\widehat W$ is contained in a compact subset of $U$), 
we have 
$$
\H^d(W_1) \leq M \H^d(\varphi_1(W_1)) + h r^d,
\leqno (2.5)
$$
where as before
$W_1 = \big\{ y \in E \cap B \, ; \varphi_1(y) \neq y \big\}$. 

\medskip
Here $GSAQ$ stands for generalized sliding Almgren quasiminimal set;
we should probably have mentioned $\Omega$ and the $L_j$ in the
notation, but this could have been too heavy.

Definition 2.3 is the sliding analogue of Definition 2.10 in [D5]. 
The case when $h=0$ corresponds to quasiminimal sets, 
as in [A2] and [DS4], 
except that here we insist that our final deformation $\varphi_1$
comes as the end of a one-parameter family of continuous maps
that satisfy the constraints (1.7). Without these  constraints
and if $U$ were convex, it would not have been necessary to mention this
(because we could take $\varphi_t(x) = (1-t)x + t \varphi_1(x)$),
but here we need to be more careful.

Notice that we allow competitors of $E$ in balls $B$ that are not
necessarily contained in $U$, but only require (2.4). It would have 
been essentially as reasonable to restrict to $B \i U$; 
this would have given an apparently larger classes $GSAQ$, 
and probably our main results are
still valid in that class. Here we opted for the definition which is
closest to Almgren, also because the invariance under changes of 
variables is a little better (Proposition 2.8 below would not work as 
nicely). If the reader ever encounters a $GSAQ$ set for the weaker
version, but not the one we give, she will probably get the desired
results inside $U$ by noticing that it is also a $GSAQ$ set (official
definition) in a slightly smaller open set.

Notice also that our Lipschitz mappings $\varphi_t$ are only 
defined on $E$. If they were allowed
to take values in $\R^n$, this would not matter because we could 
extend them. Here we also require in (1.7) that our set $\varphi_1(E)$
is a deformation of $E$, with the constraints mentioned above, but we 
see no need to require that $\varphi_t$ extends to a mapping from 
$\Omega$ to $\Omega$, for instance, and requiring boundary 
constraints like (1.7) on the $L_j\sm E$ seems really unnatural.
We want to say that the soap is attached in some way to the
boundaries, not that every deformation comes from some
global deformation in space.

When we take $M=1$ and $h$ small, we get a notion which is  closer to 
the notions of almost minimality used in [D5]. 
We are allowed to take $\delta = +\infty$, but often taking $\delta < 
+\infty$ will help. For instance, sliding almost-minimal sets will be 
sets $E$ that lie in $GSAQ(1,\delta ,h(\delta))$ for $\delta$ small,
and with an $h(\delta)$ that tends to $0$ with $\delta$.

The difference between (2.5) and its analogue for quasiminimal sets
(i.e., when $h=0$) is not enormous; the only situations where we expect 
(2.5) to be harder to use are when $\H^d(W_{1})$ and 
$\H^d(\varphi_1(W_{1}))$ are very small, i.e., when $\varphi_1$ 
only moves very few points of $E$. 
The point of a good part of Section 2 in [D5] 
was to show that these situations can be 
avoided when we prove regularity theorems. Here we shall also
need to check that we can adapt the proofs to the case of sliding
boundary conditions.

\bigskip
We shall work with reasonably strong assumptions on
$\Omega$ and the $L_j$, and already this will give us some
notational trouble.
Let us distinguish between two sets of assumptions.

We introduce first a set of assumptions for $\Omega$ and the
$L_j$, which we shall call the ``rigid assumption". 
Its main advantage is its simplicity, and many results will
be proved first under the rigid assumption, and generalized
(some times painfully) to the Lipschitz assumption below.
Again set $\Omega =L_0$, as in (1.1), to simplify the notation.

We shall say that the \underbar{rigid assumption} is satisfied when 
there is an integer $m \geq 0$ such that, for each $0 \leq j \leq j_{ max}$,
$$\eqalign{
&\hbox{$L_j$
coincides in the unit ball with the union of a finite}
\cr& 
\hbox{ number of faces $F_{j,l}$ of dyadic cubes of side length $2^{-m}$.}}
\leqno (2.6)
$$
We shall sometimes refer to the largest $2^{-m}$ such that (2.6)
holds as the \underbar{rigid scale} of the $L_j$. 

Our cubes and faces will always be closed, by convention.
When we say dyadic cube of side length $2^{-m}$, we mean a
set $[0,2^{-m}]^n + 2^{-m} k$, with $k \in \Bbb Z^n$.
The dimensions of the faces $F_{j,l}$ may be anything from
$0$ to $n$, and they may be different from each other, even for a 
fixed $j$. With this definition, it happens that the origin
plays a special role (it lies on the boundary of all the
faces of dimension $\geq 1$ that touch it), but we shall never
need this coincidence (and it will disappear in the next definition).

In terms of combinatorics, this definition still allows a lot
of different possibilities. We also authorize sets $L_j$
that are unions of faces of large dimensions, connected to
each other by lower dimensional faces, for instance, or that
just meet at one point. 

For us the rigid assumption is a toy model for more general
Plateau problems with boundary conditions of mixed dimensions. 
We decided to work with faces of dyadic cubes because this will make
our life much easier in some case, at least in terms of notation
but maybe not only. There are two main objections with this. The first one
is the rigidity of the faces, and the next definition will take care
of this. The second one is that the dyadic structure puts some 
constraints on the combinatorics of our boundary sets (for instance,
it gives a small bound on the number of $2$-dimensional faces that 
touch a given $1$-dimensional face), and this will
not be addressed. See Remarks 2.12 and 2.13. 

So we want to be able to use less rigid faces, which are
fairly smooth but not completely flat. Also, at least as far as 
quasiminimal sets are concerned, we expect some biLipschitz
invariance of our results, so we introduce the following 
weaker ``Lipschitz assumption", where we keep the same structure
for the $L_j$, but allow Lipschitz faces.

\ms\proclaim Definition 2.7.
We say that the \underbar{Lipschitz assumption} is satisfied in the open 
set $U$ when there is a constant $\lambda > 0$ and a 
bilipschitz mapping $\psi : \lambda U \to B(0,1)$
such that the sets $\psi(\lambda (L_j \cap U))$, $0 \leq j \leq j_{max}$,
satisfy the rigid assumption. 

\ms
Obviously, in this definition $\lambda U = \psi^{-1}(B(0,1))$ needs to be 
bilipschitz equivalent to $B(0,1)$, but this will not be a problem. 
In fact, all our conditions and results will be local, 
so even if our initial domain $U$ is not a nice
open set, we can try to apply our results to a smaller domain $V \i U$ 
(such that $\lambda V$ is bilipschitz equivalent to $B(0,1)$), using the fact that 
$E \in GSAQ(V,M,\delta ,h)$ as soon as $E \in GSAQ(U,M,\delta ,h)$.

Notice that the Lipschitz assumption comes with two 
important constants: the bilipschitz constant for $\psi$, 
and the rigid scale $2^{-m}$ above. The last constant $\lambda > 0$
is just a normalization, and should never play a serious role in the
estimates. In fact, we could have decided to take $\lambda = 1$ in 
the definition, and this would only have forced us to apply our 
results to dilations of the considered sets and domains.

\ms
As far as quasiminimal sets are concerned, there will not be too much
difference between our two assumptions; the following proposition will
allow us to to go from the rigid assumption to the Lipschitz 
assumption, at the price of making some constants larger.

\ms\proclaim Proposition 2.8.
Suppose that the Lipschitz assumption is satisfied in the open 
set $U$, and let $\lambda$ and $\psi$ be as in Definition 2.7. 
Also denote by $\Lambda \geq 1$ the bilipschitz constant for $\psi$.
Then, for each $E \in GSAQ(U, M, \delta , h)$, the set
$\psi(\lambda E)$ lies in $GSAQ(B(0,1), \Lambda^{2d}M, 
\Lambda^{-1}\lambda\delta , \Lambda^{2d}h)$. 

\ms
Indeed, let $\{ \varphi_t \}$, $0 \leq t \leq 1$ be as in 
Definition 1.3 (relative to the definition of a competitor for $\psi(\lambda E)$),
and also assume that $\widehat W \i \i B(0,1)$.
Set $\wt\psi (x) = \psi(\lambda x)$; thus $\wt\psi$
is the natural mapping from $U$ to $B(0,1)$. Then set 
$\widetilde\varphi_t = \widetilde\psi^{-1} \circ \varphi_t \circ 
\widetilde \psi$ for $0 \leq t \leq 1$. It is easy to see that
the $\{ \widetilde\varphi_t \}$, $0 \leq t \leq 1$, satisfy
the conditions of Definition 1.3, except that $B$ should be replaced
with $\widetilde\psi^{-1} (B)$, which itself is contained in a ball
$B'$ of radius $\widetilde r \leq  \Lambda \lambda^{-1} r$, where $r$ is 
the radius of $B$.

In addition, the analogue for the $\{ \wt\varphi_t \}$
of $\widehat W$ is $\widetilde W = \widetilde\psi^{-1} (\widehat W)$,
which is compactly contained in $U$ because $\widehat W \i \i B(0,1)$.

If $r < \Lambda^{-1} \lambda \delta$, then $\widetilde r < \delta$,
and the analogue of (2.5) yields
$$
\H^d(\widetilde W_1) 
\leq M \H^d(\widetilde\varphi_1(\widetilde W_1)) 
+ h  \widetilde r^d,
\leqno (2.9)
$$
with $\widetilde W_1 = 
\big\{ y\in \widetilde\psi^{-1}(E) \, ; \, \widetilde\varphi_1(y) 
\neq y \big\} = \widetilde\psi^{-1}(W_1)$.
We apply $\wt\psi$ and get that
$$\eqalign{
\H^d(W_1) 
&= \H^d(\widetilde\psi(\widetilde W_1))
\leq \lambda^d \Lambda^d  \H^d(\widetilde W_1) 
\cr&
\leq \lambda^d \Lambda^d M \H^d(\widetilde\varphi_1(\widetilde W_1)) 
+ \lambda^d \Lambda^d h  \widetilde r^d
\cr& = \lambda^d \Lambda^d M \H^d(\widetilde\varphi_1(\widetilde W_1)) 
+  \Lambda^{2d} h  r^d
}\leqno (2.10)
$$
by (2.9). In addition, 
$\widetilde\varphi_1(\widetilde W_1) = 
\widetilde\psi^{-1} \circ\varphi_1 \circ \widetilde\psi
(\widetilde W_1) = \widetilde\psi^{-1} \circ\varphi_1 (W_1)$, 
so (2.10) says that
$$\eqalign{
\H^d(W_1) &\leq
\lambda^d \Lambda^d M \H^d(\widetilde\psi^{-1} \circ \varphi_1 (W_1))
+  \Lambda^{2d} h  r^d
\leq \Lambda^{2d} M \H^d(\varphi_1 ( W_1))
+  \Lambda^{2d} h  r^d,
}\leqno (2.11)
$$
as needed for Proposition 2.8.
\qed

\ms
Because of Proposition 2.8, we shall sometimes be able to deduce
local regularity properties for the quasiminimal sets under the
Lipschitz assumption from their counterparts under the
rigid assumption. This will work fine for regularity properties
that are invariant under bilipschitz mappings (local Ahlfors 
regularity, rectifiability, or even uniform rectifiability),
but for more sensitive properties, or when we want a precise
dependence on the quasiminimality constants, we shall often
need to conjugate our rigid proofs and check painfully that 
they extend to the Lipschitz assumption.

\ms\noindent{\bf Remark 2.12.}
Our sets of boundaries are not nearly as general as they should be
(for the weak regularity properties that we shall prove).
There should not be anything so special about dyadic cubes, and 
we should probably have considered more general nets constructed
with convex polyhedra, with a lower bound on the angles in the subfaces.
But then the notations would have been somewhat worse, and the author 
was just afraid. Possibly the difficulty is only a matter of
organization, but the reader should be warned that in a few places,
we shall use the description of the $L_j$ with standard dyadic cubes 
to give short proofs, and the author did not even think about how 
these proofs could be adapted to more general nets. We explain about 
this a few times, but when other things become more complicated
(for instance, in Part IV), we simply forget the issue.

Hopefully, the lack of generality of our rigid and Lipschitz
assumptions will be slightly reduced by the fact that we allowed 
bilipschitz images. But on the other hand, we are missing many
simple combinatorial cases. For instance, if we want to allow an $L_j$
where $20$ faces of dimension $2$ bound a single segment, we will have
to adapt the definitions and proofs below, or play a dirty trick 
such as pretending we live in $\R^{10}$.

When we deal with more precise regularity properties
that are not invariant under bilipschitz mappings, 
we may have to choose new sets of assumptions that are
not as restrictive as the rigid assumption (which forces angles
between faces to be multiples of $90^\circ$, for instance),
and not as lenient as the Lipschitz assumption (which allows
ugly Lipschitz faces). Typically, this will happen in
Section 24, when we study blow-up limits, and where we will
allow $C^1$ faces that make different angles.

\ms\noindent{\bf Remark 2.13.}
On the other hand, at first sight it looks like we are
making our life more complicated than needed, by allowing 
large integers $m \geq 1$. Let us discuss this in the 
simple case of the rigid assumption. We are interested in local 
regularity properties of an almost- or quasi minimal set $E$ 
near a point $x_0$. If we concentrate on balls of size smaller than
$2^{-m-2}$, we are reduced to the situation where each $L_j$ is a 
cone, centered at the origin (or at the point of the dyadic grid 
of size $2^{-m}$ that lies closest to $x_0$). This seems simpler
than the situation we described, but in fact the difference is not
enormous because the combinatorics of the intersections of our
cones with a small sphere are not much simpler than the 
combinatorics of the intersections of small dyadic cubes in 
one less dimension. So we would essentially win a dimension,
but we should not expect drastic simplifications in the 
combinatorics. Also, and this is the main reason, 
allowing $m$ to be large will not complicate our proofs.

\ms\noindent{\bf Remark 2.14.}
The following convention may be useful. We shall say that out list of
boundaries $\{ L_j \}$ is complete when
$$
\hbox{for every choice of
$0 \leq i,k \leq j_{ max}$, }
L_i \cap L_k \hbox{ is one of the $L_j$}
\leqno (2.15)
$$
and also
$$
\hbox{each $L_j$ is connected.}
\leqno (2.16)
$$
Replacing the initial list of $L_j$ with a complete one
costs us nothing. Indeed, adding $L_i \cap L_j$ to our
list does not upset (1.7), because (1.7) for $L_i\cap L_k$ 
is an immediate consequence of (1.7) for $L_i$ and (1.7) 
for $L_k$. And (1.7) for $L_j$ is equivalent to (1.7)
for each of its component, because of (1.4); since we shall
only consider sets $L_j$ which have a finite number of connected 
components, the new collection of sets $L_j$ stays finite. 
We may also assume that $\Omega$ is connected,
because otherwise we could study minimal or almost minimal
sets component by component.

\bigskip
\noindent {\bf 3. Coral GSAQ and Lipschitz retractions on the $L_j$}
\medskip

In this section we deal with two technical problems.
First, we shall later find it more reassuring to restrict our attention 
to ``coral" quasiminimal sets, defined as follows.

\ms\proclaim Definition 3.1.
For $E \i \R^n$ closed, with locally finite $\H^d$ measure,
we denote by $E^\ast$ the closed support of the restriction of 
$\H^d$ to $E$; thus
$$
E^\ast = \big\{ x\in E \, ; \, \H^d(E\cap B(x,r)) > 0
\hbox{ for every } r>0 \big\}.
\leqno (3.2)
$$
We say that $E$ is \underbar{coral} when $E^\ast=E$. 

\ms
The definition comes from [D4], 
where $E^\ast$ was also called the core of $E$,  
and we wanted to distinguish coral from a slightly different notion
of ``reduced".
The main goal of this section is to check that  if $E \in GSAQ(M,\delta ,h)$,
then automatically $E^\ast \in GSAQ(M,\delta ,h)$, but since this
unexpectedly does not seem to follow too obviously from the definitions, 
we shall restrict to the Lipschitz setting for the sets $L_j$ that 
was described in Section 2.

\ms\proclaim Proposition 3.3. 
Suppose that $E \in GSAQ(U, M, \delta , h)$
and the Lipschitz assumption is satisfied on the open set $U$. 
Then $E^\ast \in GSAQ(U, M, \delta , h)$.

\ms
Observe that we do not say that $E^\ast$ is a competitor for $E$, and
indeed it is not always true: it may happen that $E$ is a nice 
$d$-dimensional surface, plus 
a $(d-1)$-dimensional handle that cannot be deformed away (or to
a subset of $E^\ast$) inside $\Omega$. The proof of Proposition 3.3 
will be slightly complicated because when some part of $E \sm E^\ast$ 
lies on the $L_j$, it adds some constraints on the competitors that 
we want to use. Put in another way, we have to show that if $E^\ast$
is not a GSAQ set because of some deformation $\{ \varphi_t \}$, we 
cannot add a set of vanishing measure to $E^\ast$, in particular on 
the $L_j$, in such a clever way that we would not be 
able to extend $\varphi_t$ so that (1.7) holds also on $E \sm E^\ast$.
\ms
Before we really start the proof, we want to construct Lipschitz 
retractions from a neighborhood of each $L_j$ onto $L_j$. 
In fact we shall do this for any finite 
union of faces of dyadic cubes of the same side length.

\ms\proclaim Lemma 3.4. 
Let $L$ be a finite union of faces of dyadic cubes of side length $1$,
possibly of different dimensions, and set
$$
L^\eta = \big\{ y\in \R^n \, ; \, \dist(y,L) \leq \eta \big\},
\leqno (3.5)
$$
where in fact we shall take $\eta =1/3$. There is a Lipschitz mapping
$\pi = \pi_L : L^\eta \to L$ such that
$$
\pi(x) = x \hbox{ for } x\in L
\leqno (3.6)
$$
and $\pi(F) \i F$ for each face $F$ (of any dimension) of each
dyadic cube of side length $1$. The Lipschitz constant for $\pi$ 
is less than $C$, where $C$ depends only on $n$.

\ms
We shall construct $\pi$ as a composition of mappings $\rho_m$.

For $m \geq 0$, denote by $A_m$ the set of faces of dimension $n-m$ 
(of dyadic cubes of side length $1$) which touch 
$L$ but are not contained in $L$. 

We may stop  at $m=n-1$, because $A_n = \emptyset$ by 
definition. Then set 
$$
T_m = L \cup \big[ \bigcup_{F\in A_m} (F\cap L^\eta) \big].
\leqno (3.7)
$$
We shall define $\rho_m$ on $T_m$, also as a composition of
simpler mappings. But let us first check a few facts about distances. 
We shall often use the fact that
$$\eqalign{
&\hbox{if $F$, $F'$ are faces of unit dyadic cubes
and $F$ is neither a point}
\cr& \hskip 0.2cm
\hbox{nor contained in $F'$, then 
$\dist(y,F') \geq \dist(y,\partial F)$
for $y\in F$.}
}\leqno (3.8)
$$
Here and below, $\d F$ is the boundary of the face $F$; it is
thus the union of some sub faces of dimension one less (except if
$F$ is a point and $\d F = \emptyset$).
Now (3.8) can be deduced from simple considerations of Euclidean 
geometry; if we were dealing with faces of polyhedra, we would merely
get that $\dist(y,F') \geq \eta_0 \dist(y,\partial F)$,
where $\eta_0$ depends on the smallest angles that adjacent faces 
of polyhedra can make, and on the smallest distance between non 
adjacent faces, and this would only force us to take $\eta$ smaller
in Lemma 3.4.
But let us just check (3.8) for faces of dyadic cubes.

Let $l$ be the dimension of $F$; thus $l \geq 1$.
Without loss of generality, we may assume that $F$
is given by the equations $0 \leq y_j \leq 1$ for
$1 \leq j \leq l$, and $y_j = 0$ for $j > l$.
Since (3.8) is trivial for points of $\d F$,
we just consider points $y\in F$ such that $0 < y_j < 1$
for $j \leq l$. 
Notice that $\dist(y,\d F)$ is the smallest of the
$\Min(y_j,1-y_j)$, $1 \leq j \leq l$.

Let $z\in F'$ minimize $|z-y|$.
If $z_j \neq y_j$ for some $j \leq l$, then
$z_j \notin (0,1)$, because otherwise we could
replace $z_j$ with $y_j$, and get a new point $z'$ that still
lies in $F'$, but is strictly closer to $y$. 
In this case, $|z-y| \geq |z_j-y_j| \geq \Min(y_j,1-y_j)
\geq \dist(y,\d F)$, as needed. So we may assume that 
$z_j = y_j$ for $1 \leq j \leq l$. If $|z_j| \geq 1$
for some  $j > l$, then $|z-y| \geq |z_j-y_j| = |z_j|
\geq 1 \geq \dist(y,\d F)$, which is fine.
Otherwise, we can replace all $z_j$, $j > l$, with $0$,
and get a new point $z' \in F'$. But $z'=y$, hence $y\in F'$.
This  is impossible, because $0 < y_j < 1$ for $1 \leq i \leq l$
(recall that $y\in F \sm \d F$), and this would force $F \i F'$
(because $F$ is the smallest face that contains $y$).
This proves (3.8).

\ms
Let us deduce from (3.8) that when $F \in A_m$,
$$
\dist(y,L) \geq \dist(y,\partial F)
\ \hbox{  for $y\in F$.}
\leqno (3.9)
$$
Indeed, if $z\in L$ and $F'$ is a face of $L$
that contains $z$, we know that $F'$ does not contain $F$
by definition of $A_m$, and also that $F$ is not reduced to one
point because $m < n$ (recall that $A_n = \emptyset$),
so (3.8) says that
$|y-z| \geq \dist(y,F') \geq \dist(y,\partial F)$.
Similarly,
$$
\dist(y,T_m \sm F) \geq \dist(y,\partial F)
\ \hbox{  for $F \in A_m$ and $y\in F$,}
\leqno (3.10)
$$
because if $z\in  T_m \sm F$, then either $z\in L$ and
we can apply (3.9), or else $z$ lies in some other face
$F' \in A_m$, and we can apply (3.8) because $F \neq F'$ and they
have the same dimension.

\ms
For each face $F\in A_m$, denote by $x_F$ the center of $F$
and by $p_F$ the radial projection from $F \sm \{ x_F \}$ to $\d F$.
That is, $p_F(y)$ is the point $z\in \d F$ such that
$y \in [x_F,z]$. By (3.9), $\dist(x_F,L) \geq 1/2$, hence
$p_F$ is defined and Lipschitz on $F \cap L^\eta$.

Extend $p_F$ to $T_m$ by setting $p_F(y)=y$ for $y\notin F$.
This is coherent, because if $F'$ is a different face
of $A_m$, then $F \cap F' \i \d F$ (recall that $F$ and $F'$
have the same dimension), and similarly $L \cap F \i \d F$
by (3.9); hence both definitions yield $p_F(y)=y$ on these sets.

Observe that $p_F$ respects the faces, i.e., $p_F(G\cap T_m) \i G$ 
for every face $G$ of any dimension of a dyadic cube of side length $1$. 
This is clear when $G$ does not meet the interior of $F$, because
then $p_F(y) = y$ on $G$; otherwise, when $G$ meets the 
interior of $F$, $G$ contains $F$ and we just need to know that
$p_F(F) \i F$. Next let us check that
$$
p_F \hbox{ is $6 \sqrt n$-Lipschitz on } T_m.
\leqno (3.11)
$$
Recall that $\dist(x_F,L) \geq 1/2$ by (3.9), so
$\dist(x_F,L^\eta) \geq 1/2-\eta = 1/6$, and hence $p_F$
is $6 \sqrt n$-Lipschitz on $T_m \cap F$. It is trivially
$1$-Lipschitz on $T_m \sm F$, and for $y\in T_m \cap F$
and $z\in T_m \sm F$,
$$\eqalign{
|p_F(y)-p_F(z)| &= |p_F(y)-z| \leq |p_F(y)-y|+|y-z|
\cr&
\leq \sqrt n \dist(y,\partial F)+|y-z|
\cr&
\leq \sqrt n \dist(y,T_m \sm F) + |y-z|
\leq (1+\sqrt n) |y-z|
}\leqno (3.12)
$$
by (3.10). Thus (3.11) holds.

\ms
Now define $\rho_m$ on $T_m$ to be the composition
of all the $p_F$, $F\in A_m$. Notice that since
$p_F$ only moves the interior points of $F$,
which lie out of $L$ by (3.9) and out of the other
$F'\in A_m$ because distinct faces of the same dimension
have disjoint interiors, we see that the order of composition
does not matter (each point is moved at most once), and in
fact $\rho_m(y) = p_F(y)$ on $F$ for each 
$F \in A_m$, and $\rho_m(y)=y$ on $L$.
Also, $\rho_m$ is $C$-Lipschitz, with $C \leq 36n$ 
(refine the proof of (3.11), or brutally observe that 
on $\{ x, y \}$, $\rho_m$ is the composition of two 
$6 \sqrt n$-Lipschitz mappings).

\ms
Next we want to compose the $\rho_m$.
Let us first check that 
$$
\dist(p_F(y),L) \leq \dist(y,L) \leq \eta
\leqno  (3.13)
$$
for $F\in A_m$ and $y\in F \cap L^\eta$.
If we were working with polyhedra, we would use (3.9)
to show that $\dist(p_F(y),L) \leq C\dist(y,L)$, and this would
be fine too, except that we would need to choose a smaller $\eta$
at the end.

In order to prove (3.13), we may assume that 
$$
F = \big\{ y \in \R^n \, ; \, 0 \leq y_i \leq 1 \hbox{ for }
1 \leq i \leq m-n \hbox{ and } y_i=0 \hbox{ for } i>m-n \big\}
\leqno  (3.14)
$$
and, by symmetry, that all the coordinates of $y$ lie in $[0,1/2]$. 
Set $\wt y = p_F(y)$; then $0 \leq \wt y_i \leq y_i$ for $1 \leq i \leq m-n$,
because $0 \leq y_i \leq 1/2$, the coordinate of $x_F$. 

Let $z\in L$ lie closest to $y$; we just want to find $\wt z \in L$
such that $|\wt z - \wt y| \leq |z-y|$.
Since $|z-y|\leq 1/3$, all the coordinates $z_i$ lie in $[-1/3,5/6]$. 
Let $i$ be such that $z_i \leq 0$; we keep $\wt z_i = z_i$, and
obviously $|\wt z_i - \wt y_i| = |z_i - \wt y_i| \leq |z_i-y_i|$.
For the other $i$, we know that $0 < z_i < 5/6$, and we just
set $\wt z_i = \wt y_i$; notice that the point $\wt z$ that we get
this way lies in the same faces as $z$, because we only replaced
some coordinates that lie in $(0,1)$ with other ones in $[0,1]$,
and this operation preserves any face. Thus $\wt z \in L$, just
like $z$, and since by construction $|\wt z_i - \wt y_i| \leq |z_i-y_i|$ 
for all $i$,  we completed the proof of (3.13).

\ms
Since $p_F$ is the identity out of $F$, (3.13) is also valid
for $y\in T_m \i L^\eta$. We claim that
$$
\rho_m(T_m) \i
L \cup \big[ \bigcup_{F\in A_m} (\d F\cap L^\eta) \big]
\i T_{m+1}.
\leqno  (3.15)
$$
Let $w \in \rho_m(T_m)$ be given, and let $y\in T_m$ be 
such that $w = \rho_m(y)$. If $y \in L$, then
$p_F(y) = y$ for all $F$ (because $F \cap L \i \d F$,
by (3.9)), hence $w=y \in L$. If $y\in F$ for some
$F \in A_m$, then $p_F(y) \in \d F$ by construction,
and then all the other $p_{F'}$ preserve $\d F$, because
they preserve every face of every cube; thus $w = \rho_m(y)$
lies in $\d F$ too (recall that we can compose the $p_F$
in any order that we like). Also, $w \in L^\eta$ by
successive applications of (3.13). 

For the second inclusion, let 
$F \in A_m$ and $w \in \d F\cap L^\eta$ be given.
Let $H$ be a $(m-n-1)$-dimensional face of $\d F$
that contains $w$. If $H \i L$, we are happy because
$L \i T_{m+1}$. Otherwise, as soon as we prove that $H$ meets $L$,
we will know that $H \in A_{m+1}$ (by definition of $A_{m+1}$),
hence $w \in H \cap L^\eta \i T_{m+1}$, as needed.
Now $\dist(w,L) \leq \eta$ because $w\in L^\eta$, so
we can find $z \in L$ such that $|z_i-w_i| \leq 1/3$
for all $i$. When $i$ is such that $z_i \neq w_i$, we can
replace both $z_i$ and $w_i$ with some integer $n_i$
which is close to both of them, without changing the
fact that $w\in H$ and $z\in L$; this way we get a point 
of $H\cap L$. This completes the proof of our claim (3.15).

Now we set $\pi = \rho_{n-1}\circ \ldots \circ\rho_0$.
This is a Lipschitz mapping which is defined on
$T_0 = L^\eta$ and takes values in $T_n = L$. 
Thus $\pi(L^\eta) \i L$. Next, (3.6) holds because 
$p_F(y) = y$ on $L$ for al $F$; finally, $\pi$
preserves the faces because it is a composition of mappings 
that preserve the faces of all dimensions. 
Thus $\pi$ is the desired mapping, and Lemma 3.4 follows.
\qed

\ms\noindent{\bf Remark 3.16.}
For $0 \leq \eta \leq 1/3$, we also get a mapping $\pi_L$
as in Lemma 3.4, which is just the restriction to $L^\eta$
of the mapping that we construct with $\eta = 1/3$.
That is, we always use the same formulas, only the domains of definition 
differ.

\ms
The retraction from Lemma 3.4 is the endpoint of 
a deformation; we shall not need this fact before 
Lemma 8.8, but let us check it now before we forget the notation.

\ms\proclaim Lemma 3.17.
Let $L$, $0 < \eta \leq 1/3$, and $L^\eta$ be as in Lemma 3.4.
Then there is a Lipschitz mapping $\Pi_L : L^\eta \times [0,1] \to \R^n$
such that 
$$
\Pi_L(x,t) = x 
\ \hbox{ for $x\in L$ and for } t= 0,
\leqno (3.18)
$$
$$
\Pi_L(x,1) = \pi_L(x) \ \hbox{ for } x\in L^\eta,
\leqno (3.19)
$$
$$
|\Pi_L(x,t) - \Pi_L(x,s)| \leq C \dist(x,L) |t-s|
\ \hbox{ for $x\in L^{\eta}$ and } 0 \leq s,t \leq 1,
\leqno (3.20)
$$
$$
|\Pi_L(x,t) - \Pi_L(y,t)| \leq C |x-y|
\ \hbox{ for $x, y \in L^\eta$ and } 0 \leq t \leq 1,
\leqno (3.21)
$$
and $\Pi_L$ also preserves the faces of all dimensions, i.e.,
$$\eqalign{
\Pi_L(x,t) \in F &\hbox{ whenever $F$ is any face (of any dimension) }
\cr&
\hbox{of a dyadic cube of side length $1$,  $x\in F$, and $0 \leq t \leq 1$.}
}\leqno (3.22)
$$
The constant $C$ in (3.20) and (3.21) depends only on $n$.

\ms
To see this, observe that $\pi_L$ is obtained by composing
a bounded number of Lipschitz mappings $p_F$, where 
$F\in \cup_{m} A_m$ is some face of dyadic cube.
Recall from the definition below (3.10) that
when $F \in A_m$, $F$ is of dimension $n-m$, $p_F$
is defined on the set $T_m$ of (3.7), and is
equal to the identity everywhere, except on $F$
itself, where it is a radial projection on $\d F$.
We easily go from the identity to $p_F$ by setting
$$
p_F(x,t) = t p_F(x) + (1-t) x
\ \hbox{ for $x \in T_m$ and } 0 \leq t \leq 1;
\leqno (3.23)
$$
then the $p_F(\cdot,t)$ also preserve the faces of all dimensions,
are $6 \sqrt n$-Lipschitz like $p_F$, and
$$
|p_F(x,t)-p_F(x,s)| \leq C \eta |t-s|
\ \hbox{ for $x\in T_m$ and } 0 \leq s,t \leq 1,
\leqno (3.24)
$$
because $|p_F(x,t)-x| \leq C \eta$. 

When we used $p_F$, we composed it with a previous mapping $h$, which
maps $L^\eta$ to $T_m$ by (3.15) and because the other
$p_{F'}$, $F' \i A_m$, map $T_m$ to $T_m$. Then
we can go from $h$ to $p_F \circ h$ by setting
$h_F(x,t) = p_F(h(x),t)$ for $x\in L^\eta$ and 
$0 \leq t \leq 1$. The mapping $h_F$ is $C$-Lipschitz in $x$
and $C \eta$-Lipschitz in $t$, because $h$ is $C$-Lipschitz.

We now concatenate all the deformations $h_F$, reparameterize
by the unit interval, and get a mapping $\Pi_L$ that satisfies
(3.18)-(3.22), except that in (3.20) we only get $\eta$
instead of $\dist(x,L)$. 
But Remark 3.16 extends to our mapping $p_F(x,t)$ and $\Pi_l \,$:
the mapping that we would construct on $L^{\eta'}$, with
$\eta' = \dist(x,L)$, is just the restriction to $L^{\eta'}$ 
of the mapping that we constructed here on $L^\eta$.
Therefore, (3.20) is just the same thing as
(3.24) for $x$, but applied to the mapping $\Pi_L$
associated to $\eta'$.
\qed

\ms\noindent{\bf Remark 3.25.}
Of course we can also define $\pi_L$ and $\Pi_L$ when $L$ is a finite
union of faces of dyadic cubes, not necessarily of size one.
That is, if $L$ is a finite union of faces dyadic cubes of size
$2^{-m}$ (as in the definitions of our $L_j$), we define 
$\pi_L$ by
$$
\pi_L(z) = 2^{-m} \pi_{2^m L}(2^m z),
\leqno  (3.26)
$$
and use a similar definition for $\Pi_L$.
When we define the $\pi_{L_j}$ associated to our boundary pieces $L_j$, 
we shall use this convention; if by luck some $L_j$ are also unions of 
dyadic faces of larger diameters, we shall ignore that fact and stay 
with the same $m$.

\ms
Next we want to prove Proposition 3.3 in a simpler setting.
We shall later see how the proof of Proposition 2.8 allows us
to reduce to this case.

\ms\proclaim Proposition 3.27. 
Suppose $E \in GSAQ(B_0, M, \delta , h)$, where we set $B_0 = B(0,1)$, 
and that the rigid assumption is satisfied. Then 
$E^\ast \in GSAQ(B_0, M, \delta , h)$.

\ms
So let $E \in GSAQ(B_0, M, \delta , h)$ be given.
We shall go from $E$ to $E^\ast$ in a finite number of steps, where
each time we remove a set in some $L_j$. 
We may assume that the set of $L_j$ is complete, as in Remark 2.14,
and that when we enumerate the various $L_j$, we start with the largest
ones for the inclusion relation.

We shall first define some intermediate sets $E_j$. Set
$$
Z' = E \sm E^\ast \ \hbox{ and } \ 
Z_j = Z' \sm [\cup_{i \geq j} L_i]
\ \hbox{ for } \ 0 \leq j \leq j_{ max}+1.
\leqno  (3.28)
$$ 
This is a nondecreasing sequence of open subsets of $E$.
The first one is $Z_0 = \emptyset$, because
$L_0 = \Omega$ contains $E$, and the last one is $Z_{j_{ max}+1} = Z'$.
Note that
$$
\H^d(Z') = \H^d(E \sm E^\ast) = 0,
\leqno  (3.29)
$$
by definition of $E^\ast$ (see (8.26) on page 58 of [D4], 
for instance, for the elementary proof).
Then set 
$$
E_j = E \sm Z_j
\ \hbox{ for } \ 0 \leq j \leq j_{ max}+1.
\leqno (3.30)
$$ 
This is a nonincreasing sequence of closed sets, 
with $E_{0} = E$ and $E_{j_{ max}+1} = E^\ast$. 
We want to prove by induction that $E_j$ lies
in $GSAQ(B_0, M, \delta, h)$, just like $E$. 
The induction assumption holds for $j=0$, so let us assume 
that $0 \leq j \leq j_{max}$ and that 
$$
E_j \in GSAQ(B_0, M, \delta, h),
\leqno (3.31)
$$
and prove that 
$$
E_{j+1} \in GSAQ(B_0, M, \delta, h).
\leqno (3.32)
$$
Set 
$$\eqalign{
Z &= E_{j} \sm E_{j+1} = Z_{j+1} \sm Z_j
= \{ Z' \sm [\cup_{i \geq j+1} L_i] \} \sm 
\{ Z' \sm [\cup_{i \geq j} L_i]\}
\cr& = [Z' \cap L_j] \sm [\cup_{i \geq j+1} L_i]
= [E \cap  L_j \sm E^\ast] \sm [\cup_{i \geq j+1} L_i].
}\leqno (3.33)
$$ 

So we want to prove (3.32). We take a sliding competitor 
$F = \varphi_1(E_{j+1})$ for $E_{j+1}$ in some closed ball $B$, 
and we want to prove the analogue of (2.5) for $E_{j+1}$.
It is tempting to use the same one-parameter family $\{ \varphi_t \}$
and apply (2.5) to it, but since it is only defined for $x\in E_{j+1}$,
we have to extend it to $x\in E_j$. The difficult part will be to make
sure that we still have (1.7), in particular at points of $Z$.
At any rate, we want to define $\wt\varphi_t(x)$ for $x\in E_{j}$
and $0 \leq t \leq 1$, and logically we would like to keep
$$
\wt\varphi_t(x) = \varphi_t(x) \hbox{ for } x\in E_{j+1}
\leqno  (3.34)
$$ 
(so that we would only need to define $\wt\varphi_t(x)$ when $x\in Z$); 
we also would like to keep $\wt\varphi_t(x) = x$ when $t=0$ and when 
$x\in E_j \sm B'$, where $B'$ is a closed ball with the same center 
as the ball $B$ of Definition 1.3, but just a tiny bit
larger, so we should mostly worry about $Z \cap B'$.

We extend $\varphi_1$ first, in a Lipschitz way. That is,
we have a Lipschitz mapping $\varphi_1$, defined on $E_{j+1}$,
and we first extend it to $E_{j+1} \cup [E_j \sm B']$ by setting
$$
\varphi_1(x) = x \ \hbox{ for } x \in E_j \sm B'. 
\leqno  (3.35)
$$
This map is still Lipschitz, 
although with a possibly very large constant (but we don't care).
Indeed, since $\varphi_1$ is Lipschitz on $E_{j+1}$ and on 
$E_j \sm B'$, we only need to estimate $|\varphi_1(x)-\varphi_1(y)|$
when $x \in E_{j+1}$ and $y\in E_j \sm B'$ say; if $x\in B$ and
$y\in 2B'$, we say that $|\varphi_1(x)-\varphi_1(y)| \leq
|\varphi_1(x)|+|\varphi_1(y)| \leq 2C \leq 2C(r'-r)^{-1}|x-y|$,
where $C$ is a bound for $\varphi$ on $2B'$, and $r,r'$ denote the 
radii of $B$ and $B'$. If $x \notin B$, $\varphi_1(x)-\varphi_1(y) = 
x-y$ by definition, and the last case when $x\in B$ and $y \notin 2B'$ 
is even easier. So $\varphi_1$ is Lipschitz on $E_{j+1} \cup [E_j \sm B']$.
Now we use the Whitney extension theorem to get a 
Lipschitz extension $\varphi_1 : E_j \to \R^n$.

Next we set $\varphi_0(x) = x$, as we should do, and get a 
mapping $(x,t) \to \varphi_t(x)$ defined on $[E_j \times \{ 0, 1 \}] \cup 
[E_{j+1} \times [0,1]]$. This mapping is continuous; this is clear 
at points $(x,t)$ where $0 < t < 1$; when $x\in E_{j+1}$ and
$t \in \{ 0, 1 \}$, we use the continuity of the previous $\varphi_t$
and the fact that $\varphi_0$ and $\varphi_1$ are Lipschitz;
finally when $x\in E_j \sm E_{j+1}$, we use the fact that $E_{j+1}$
is closed, hence far from $x$. 

We can also set $\varphi_t(x) = x$ for $0 \leq t \leq 1$ and $x\in 
E_j \sm B'$. The two definitions coincide when $t=1$, by (3.35),
when $t=0$ by definition of our first extension, and when 
$x\in E_{j+1}$ by (1.5); we now get a mapping $(x,t) \to \varphi_t(x)$ 
which is defined on $[E_j \times \{ 0, 1 \}] \cup 
[E_{j+1} \times [0,1]] \cup [(E_j \sm B')\times [0,1]]$.
Let us check that this mapping is continuous.
We just need to check this at points $(x,t)$ that lie in the 
intersection of the closures of our two sets (where we already know that
the function is continuous). When $x\in E_{j+1}$, our first extension
was already defined at $x$, with $\varphi_t(x) = x$ because 
$x \in E_{j+1} \cap (E_j \sm B')^- \i E_{j+1} \sm B$, and by (1.5).
Since the second definition also yields $\varphi_t(x) = x$, we get the
continuity at $(x,t)$. Now suppose $x \in E_j \sm E_{j+1}$. As before,
since $E_{j+1}$ is closed, $x$ is far from $E_{j+1}$ and this forces
$t \in \{ 0,1\}$. In this case too, $\varphi_t(x)$ was already defined 
for the first extension, with $\varphi_t(x) = x$, so $\varphi_t$ is 
continuous at $(x,t)$.

We now use the Titze extension theorem. This gives a continuous mapping 
$(x,t) \to \varphi_t(x)$, from $E_j \times [0,1]$ to $\R^n$.
This mapping satisfies the continuity condition in (1.4), (1.5) 
(with $B$ replaced with $B'$), and (1.8). 
Since we also want (1.6), we compose its values on $B'$
with the $1$-Lipschitz radial projection from $\R^n$ to $B'$; this
does not change the values on $E_{j+1}$, by (1.6), and 
does not destroy (1.4), (1.5), or (1.8).
Of course our last property (1.7) is not automatically satisfied, 
so we'll need to modify the $\varphi_t$ again.

\ms
Define a slightly better $f_t$ for $0 \leq t \leq 1$ by
$$
f_t(x) = \varphi_{t \psi(x)}(x),
\leqno (3.36)
$$
where $\psi(x)$ is a Lipschitz function of $d(x) = \dist(x,E_{j+1})$ 
such that
$$\eqalign{
\psi(x) &= 1 \hbox{ when } d(x) \leq \varepsilon,
\cr
0 \leq \psi(x) &\leq 1 \hbox{ when } 
\varepsilon \leq d(x) \leq 2\varepsilon,
\cr 
\psi(x) &= 0 \hbox{ when } d(x) \geq 2\varepsilon,
}\leqno (3.37)
$$
and where $\varepsilon>0$ is a very small number that will be chosen soon.
Near $E_{j+1}$, we do not change anything (because $\psi(x)=1$).
Let $\varepsilon_0 > 0$ be given. Let us check that if $\varepsilon > 0$ 
is small enough, then  
$$
\dist(x,E_{j+1}\cap L_j)\leq \varepsilon_0
\hbox{ for $x\in Z$ such that } d(x) \leq 2\varepsilon.
\leqno (3.38)
$$
Otherwise, we could find a sequence $\{ x_k \}$
in $Z$, with $d(x_k)$ tending to $0$, and that stays
$\varepsilon_0$-far from $E_{j+1}\cap L_j$. Replace 
$\{ x_k \}$ with a subsequence with some limit $x$; then
$x\in L_j$ because $Z \i L_j$ by (3.33),
and $L_j$ is closed. In addition, $x\in E_{j+1}$ because
$d(x) = 0$ and $E_{j+1}$ is closed, a contradiction with
the fact that $\dist(x_k,E_{j+1}\cap L_j) \geq \varepsilon_0$;
(3.38) follows.

Next observe that
$$
\varphi_t(x) \in L_j \ \hbox{ for }
x\in E_{j+1}\cap L_j \hbox{ and } t\in [0,1],
\leqno (3.39)
$$
by (1.7). Hence, by (3.38) and the uniform continuity
of $\varphi_t$ on $[0,1] \times E_j$, we get that for each $\varepsilon_1 >0$,
we can find $\varepsilon > 0$ such that
$$
\dist(\varphi_t(x),L_j) \leq \varepsilon_1
\hbox{ for $t\in [0,1]$ and $x\in Z$ such that } d(x) \leq 2\varepsilon.
\leqno (3.40)
$$
Then we also get that
$$
\dist(f_t(x),L_j) \leq \varepsilon_1
\hbox{ for $t\in [0,1]$ and $x\in Z$ such that } d(x) \leq 2\varepsilon
\leqno (3.41)
$$
by (3.36) and because $t\psi(x) \in [0,1]$. But then
$$
\dist(f_t(x),L_j) \leq \varepsilon_1
\hbox{ for $t\in [0,1]$ and $x\in Z$,}
\leqno (3.42)
$$
because when $x\in Z$ and $d(x) > 2 \varepsilon$,
$\psi(t) = 0$, hence $f_t(x) = \varphi_0(x) = x$, so
$f_t(x) \in L_j$ too.

This proximity to $L_j$ is the reason why $f_t$ is better than 
$\varphi_t$. On the other hand, observe that
$f_1(x) = \varphi_{\psi(x)}(x)$ is not necessarily Lipschitz,
because we did not require $\varphi_t$ to be Lipschitz in $t$,
or even in $x$ when $t<1$, so we shall fix this now and construct 
new functions $g_t$. The first step is, given a small
$\varepsilon_1 > 0$, to choose a Lipschitz
function $\varphi : E_j \times [0,1] \to \R^n$, such that
$$
\varphi(x,t) = x \ \hbox{ for } x\in E_j \sm B'
\hbox{ and for } t=0
\leqno (3.43)
$$
and 
$$
|\varphi(x,t) - \varphi_t(x)| \leq \varepsilon_1
\ \hbox{ for } x\in E_j \cap B' \hbox{ and } 0 \leq t \leq 1.
\leqno (3.44)
$$
This is easy to do, and we just sketch the proof. On the set
$A_1 = [E_j \sm B'] \times [0,1] \bigcup E_j \times \{ 0 \}$,
we simply take $\varphi(x,t) = \varphi_t(x) = x$. Then we use 
the uniform continuity of $(x,t) \to \varphi_t(x)$ on 
$E_j \times [0,1]$ (which is easy because we just need to 
consider the compact set $E_j \cap 2B' \times [0,1]$)
to get $\tau > 0$ such that
$$
|\varphi_t(x)-\varphi_s(y)| < \varepsilon_1
\leqno (3.45)
$$ 
for $x,y \in E_j$ and $t,s \in [0,1]$ such that 
$|x-y|+|t-s| \leq 100\tau$. 
We select a maximal subset $A_2$ of 
$\big\{ (x,t) \in E_j \times [0,1] \, ; \, 
\dist((x,t),A_1) \geq \tau \big\}$, whose points lie at
mutual distances at least $\tau$, and decide that
$\varphi(x,t) = \varphi_t(x)$ for $(x,t) \in A_2$.
Finally we use a partition of unity to complete the definition
of $\varphi$. We get that $\varphi$ is Lipschitz on $E_j \times [0,1]$
by construction, with a very large constant that depends on $\tau$
(but this is all right), and (3.44) holds because for each
$(x,t) \in E_j \times [0,1]$, $\varphi(x,t)$ is an average
of values of $\varphi(y,s) = \varphi_s(y)$ on $A_1 \cup A_2$ at
nearby points $(y,s)$, and by (3.45). If by bad luck $\varphi(x,t)$
falls out of $B'$ for some pairs $(x,t) \in B' \times [0,1]$, compose 
again with the radial Lipschitz projection onto $B'$, without 
altering (3.44) or the Lipschitz constants.

Next we choose a Lipschitz function $h(x)$ of $d(x)=\dist(x,E_{j+1})$,
such that $0 \leq h \leq 1$ everywhere, $h(x) = 1$
when $d(x) =  0$ (i.e., when $x\in E_{j+1}$),
and $h(x) = 0$ when $d(x) \geq \varepsilon/2$. And we set
$$
\varphi'(x,t) =  h(x) \varphi_t(x) + (1-h(x)) \varphi(x,t)
\leqno (3.46)
$$
on $E_j \times [0,1]$.
Notice that $\varphi'(x,t) = \varphi_t(x)$ on $E_{j+1}$
(because $h(x) = 1$ there), and that 
$$
|\varphi'(x,t)  - \varphi_t(x)| \leq \varepsilon_1
\leqno (3.47)
$$ 
by (3.44). Set
$$
g_t(x) = \varphi'(x,t \psi(x)) 
\leqno (3.48)
$$
(compare with (3.36)). Then 
$$
|g_t(x) - f_t(x)| 
= |\varphi'(x,t \psi(x)) -\varphi_{t \psi(x)}(x)|
\leq \varepsilon_1
\ \hbox{ for } (x,t) \in E_j \times [0,1],
\leqno (3.49)
$$
by (3.47). We still have that 
$g_t(x) \in B'$ when $x\in B'$ (because of similar properties
for $\varphi_t$ and $\varphi(x,\cdot)$, and since $B'$ is convex),
and
$$
g_t(x) = x
\ \hbox{ for $x\in E_{j}\sm B'$ and } t=0,
\leqno (3.50)
$$
because $\varphi_t(x) = x$ by construction (see above (3.36)),
$\varphi(x,t) = x$ by (3.43), $\varphi'(x,t) = x$ by (3.46),
and finally $g_t(x) = x$ by (3.48). Next,
$$
g_t(x) = \varphi_t(x)
\ \hbox{ on } E_{j+1} \times [0,1],
\leqno (3.51)
$$
because $d(x) = 0$, hence $h(x) = \psi(x)= 1$
and $\varphi'(x,t) = \varphi_t(x)$ by (3.46).
Note that $(t,x) \to g_t(x)$ is continuous because all the
ingredients in (3.48) and (3.46) are continuous. Let us
also check that 
$$
g_1 \hbox{ is Lipschitz on } E_j. 
\leqno (3.52)
$$
In the region where $d(x) \leq \varepsilon$, (3.37) says that 
$\psi(x) = 1$, so $g_1(x) = \varphi'(x,1) = h(x) \varphi_1(x) 
+ (1-h(x)) \varphi(x,1)$ by (3.48) and (3.46), which is Lipschitz 
in particular because $\varphi_1$ is Lipschitz. 
In the region where $d(x) \geq \varepsilon/2$,
$h(x) = 0$ so $\varphi'(x,t) = \varphi(x,t)$ and 
$g_1(x) = \varphi(x,\psi(x))$, which is Lipschitz by 
definition of $\varphi$. This proves (3.52) because we have more than 
enough room for the gluing.

We still have 
$$
g_0(x) =  \varphi'(x,0) 
= h(x) \varphi_0(x) + (1-h(x)) \varphi(x,0) = x
\ \hbox{ for } x\in E_j,
\leqno (3.53)
$$
by (3.48), (3.46), (3.43), and the definition of the extension
of $\varphi_t$, but we are still missing (1.7). We shall now set
$$
\wt\varphi_t(x) = \pi_j(g_t(x))
\ \hbox{ for } x \in Z \hbox{ and } t\in [0,1],
\leqno (3.54)
$$
where we denote by $\pi_j$ the Lipschitz retraction 
$\pi_{L_j}$ onto $L_j$ that we constructed with the
help of Lemma 3.4, after scaling as in (3.26). 
Recall that $\pi_j$ is defined on a $\eta$-neighborhood of $L_j$, 
where $\eta = 2^{-m}/3$ is now the third of the side length 
of the dyadic cubes that compose $L_j$. 
Note that the definition in (3.54) makes sense because
for $x\in Z$ and $0 \leq t \leq 1$,
$$
\dist(g_t(x),L_j) \leq \dist(f_t(x),L_j) + \varepsilon_1
\leq 2\varepsilon_1 < \eta
\leqno (3.55)
$$
by (3.49), (3.42), and if we choose $\varepsilon_1 < \eta/2$.
Recall from (3.34) that we would like to set 
$\wt\varphi_t(x) = \varphi_t(x)$ on $E_j \sm Z = E_{j+1}$,
but the desired Lipschitzness of $\wt\varphi_1$ at the
interface between $Z$ and $E_{j+1}$ will force us to modify this 
slightly on a small region near $L_j$. Let $\varepsilon_2 > 0$
be very small, to be chosen later, and let $\xi$
be defined by
$$
\left\{
\eqalign{
&\xi(y)=1 \ \hbox{ for }
0 \leq y \leq \varepsilon_2/2, 
\cr&
\xi(y)=0 \ \hbox{ for } y \geq \varepsilon_2, 
\cr&
\xi \ \hbox{ is affine on } [\varepsilon_2/2,\varepsilon_2].
} \right.
\leqno (3.56)
$$
Also set 
$$
d_t(x) = \dist(\varphi_t(x),L_j \cap B')
\ \hbox{ for } x\in E_{j+1} \hbox{ and } t\in [0,1].
\leqno (3.57)
$$
Notice that if $\xi(d_t(x)) \neq 0$, then
$\dist(\varphi_t(x),L_j) \leq d_t(x) \leq \varepsilon_2$ and 
hence (if $\varepsilon_2 < \eta$) $\pi_j(\varphi_t(x))$ is defined. 
This allows us to set
$$
\wt\varphi_t(x) =  t\,\xi \circ d_t(x) \,\pi_j(\varphi_t(x)) 
+(1-t\, \xi \circ d_t(x)) \,\varphi_t(x)
\ \hbox{ for } x\in E_{j+1} \hbox{ and } t\in [0,1].
\leqno (3.58)
$$

We now start checking that the $\wt\varphi_t$ satisfy 
the required properties (1.4)-(1.8) on $E_j$. We first show
that
$$
(x,t) \to \wt\varphi_t(x)
\hbox{ is continuous on } E_j \times [0,1].
\leqno (3.59)
$$
Since it is clearly continuous 
on $Z \times [0,1]$ and on $E_{j+1}\times [0,1]$, 
the only potential problem is at a point 
$(x,t)$ such that $x\in \overline Z \cap E_{j+1}$ 
(recall that $E_{j+1}$ is closed), and it is even enough
to show that our two continuous definitions (3.54) and (3.58)
give the same result at such a point.
But then $x\in L_j$ (by (3.33) and because $L_j$ is closed), 
so $\varphi_t(x) \in L_j$ (by (1.7) and because $x\in E_{j+1}$),
and using (3.54) yields the result
$$
\pi_j(g_t(x)) = \pi_j(\varphi_t(x)) = \varphi_t(x)
\leqno (3.60)
$$
by (3.51), (3.26), and (3.6). But (3.58) also yields
$\varphi_t(x)$, because $\pi_j(\varphi_t(x)) = \varphi_t(x)$
by (3.60). So (3.59) holds.

Next we check that $\wt\varphi_1$ is Lipschitz.
Again $\wt\varphi_1$ is Lipschitz on $Z$ and on $E_{j+1}$,
but we need to be careful about the interface.
That is, we just need to estimate 
$|\wt\varphi_1(x)-\wt\varphi_1(y)|$ when $x\in E_{j+1}$
and $y\in Z$. We'll distinguish between a few cases.

When $d_1(x) \leq \varepsilon_2/2$, $\xi \circ d_1(x) = 1$
and (3.58) says that
$$
\wt\varphi_1(x) = \pi_j(\varphi_1(x)) = \pi_j(g_1(x))
\leqno (3.61)
$$
by (3.51), so $|\wt\varphi_1(x)-\wt\varphi_1(y)|
= |\pi_j(g_1(x))-\pi_j(g_1(y))| \leq C |x-y|$
by (3.54) and the fact that $g_1$ and $\pi_j$
are Lipschitz. 

We claim that if $\varepsilon > 0$ is chosen small enough
(depending on $\varepsilon_2$), we have that
$d_1(x) \leq \varepsilon_2/2$ when $x\in E_{j+1}$
and $y\in Z$ are such that $|x-y| \leq \varepsilon$
and $\dist(x,B') \leq \varepsilon$.
Indeed, otherwise we can find sequences $\{ x_k \}$ in
$E_{j+1}$ and $\{ y_k \}$ in $Z$, such that 
$|x_k-y_k| \leq 2^{-k}$ and $\dist(x_k,B') \leq 2^{-k}$
but $d_1(x_k) \geq \varepsilon_2/2$.
We can extract a subsequence so that $\{ x_k \}$ converges
to a limit $x\in E_{j+1} \cap  B'$ (recall that we chose $B'$ closed). 
Then $x\in L_j$, because all the $y_k$ lie in $Z \i L_j$, so 
$\varphi_1(x) \in L_j$ by (1.7). 
In addition, $\varphi_1(x) \in B'$ by (1.6) (if $x\in B$) or  by (1.5) 
(if $x\in B'\sm B$), and so $d_1(x) = 0$, which contradicts the fact that 
$d_1(x_k) \geq \varepsilon_2/2$ because $d_1$ is continuous.
This proves our claim.

If $|x-y| \geq \varepsilon$, we simply use the fact that
$|\wt\varphi_1(x)-\wt\varphi_1(y)| \leq C \leq C \varepsilon^{-1} |x-y|$
(because $\varphi_1$ is Lipschitz on each piece, hence bounded) 
to get a (very bad) Lipschitz bound. So we may assume that
$|x-y| \leq \varepsilon$ and, since we already treated the case when
$d_1(x) \leq \varepsilon_2/2$, our claim allows us to suppose that
$\dist(x,B') \geq \varepsilon$. Then $\varphi_1(x) = x$ by (1.5),
and $y \in Z \sm B'$ because $|x-y| \leq \varepsilon$.

Let us compute $\wt\varphi_1(y)$. First observe that 
$\varphi_t(y) = y$ for $0 \leq t \leq 1$, 
because $y \in Z \sm B' \i E_j \sm B'$ and by the definition of
$\varphi_t$ below (3.35) (just before we use the Titze extension 
theorem). Also, $\varphi(y,t) = y$ for $0 \leq t \leq 1$, by (3.43),
and hence $g_1(y) = \varphi'(y,\psi(y)) = y$, by (3.48) and (3.46).
And also $\wt\varphi_1(y) = \pi_j(g_1(y)) = \pi_j(y)$ by (3.54).
In addition, observe that
$$
\pi_j(y) = y  \ \hbox{ for } y\in Z
\leqno (3.62)
$$
because $Z \i L_j$ by (3.33), and then by 
(3.26) and (3.6). So here $\wt\varphi_1(y) = \pi_j(y) =y$ and now
$$\eqalign{
|\wt\varphi_1(x)-\wt\varphi_1(y)|
&= |\xi(d_1(x))\,\pi_j(\varphi_1(x)) + (1-\xi(d_1(x)))\,\varphi_1(x)- y|
\cr&
\leq |\pi_j(\varphi_1(x))-y| + |\varphi_1(x)-y|
=  |\pi_j(x)-y| + |x-y|
\cr&
= |\pi_j(x)-\pi_j(y)| + |x-y| \leq C |x-y|
}\leqno (3.63)
$$
by (3.58), because $\varphi_1(x) =  x$ by (1.5), and by (3.62).
So $\wt\varphi_1$ is Lipschitz, which 
takes care of (1.8). Next we check that
$$
\wt\varphi_0(x) = x \ \hbox{ for } x\in E_j.
\leqno (3.64)
$$
Notice that $\varphi_0(x)=x$ by the definition below (3.35), 
$\varphi(x,0)=x$ by (3.43), and $g_0(x) = x$ by (3.46) and (3.48). 
If $x\in Z$, $\pi_j(x) = x$ by (3.62), and (3.54) yields 
$\wt\varphi_0(x)=\pi_j(g_0(x))=x$. If $x\in E_{j+1}$,
$t\xi(d_t(x)) = 0$ because $t=0$, so
(3.58) says that $\wt\varphi_0(x) = \varphi_0(x) = x$;
hence (3.64) holds. 

For (1.5), we'll need to know that
$$
\wt\varphi_t(x) = x
\ \hbox{ for $x\in Z \sm B'$ and } t\in [0,1],
\leqno (3.65)
$$
and indeed, $g_t(x) = x$ by (3.50) and 
$\wt\varphi_t(x)=\pi_j(g_t(x))=x$ by (3.54) and (3.62).
Similarly, let us check that
$$
\wt\varphi_t(x) = x
\ \hbox{ for $x\in E_{j+1} \sm B''$ and } t\in [0,1],
\leqno (3.66)
$$
where 
$$
B'' = \big\{ x\in \R^n \, ; \, \dist(x,B') \leq \varepsilon_2 \big\}
\leqno (3.67)
$$ 
is just a tiny bit larger than $B'$ and $B$. This time
$\varphi_t(x) = x$ by (1.5), hence 
$d_t(x) \geq \dist(x,B') \geq \varepsilon_2$
by (3.57), $\xi \circ d_t(x) = 0$ by (3.56), and
hence $\wt\varphi_t(x) = \varphi_t(x) = x$ by (3.58).
Since $E_j = Z \cup E_{j+1}$, we get that (1.5) holds
for $B''$ (or any larger ball).

Next we check (1.6), but with an even larger ball $\wt B$. 
Let $w \in \R^n$ denote the center of $B$, and let $r$ and
$r'$ denote the respective radii of $B$ and $B'$. Set
$$
\wt r = r' + 4C \varepsilon_1
\ \hbox{ and }
\wt B = \overline B(w,\wt r)
\leqno (3.68)
$$
where $C$ is a bound for the Lipschitz constant for $\pi_j$
and $\varepsilon_1$ is as in (3.44) and (3.55) 
(in fact, any small number chosen in advance), and 
let us check that
$$
\wt\varphi_t(x) \in \wt B
\ \hbox{ for $x\in E_j \cap \wt B$ and } t\in [0,1].
\leqno (3.69)
$$
First suppose that $x\in Z$. If $x\in Z \sm B'$, (3.65)
says that $\wt\varphi_t(x) = x \in B'$, and we are done.
Otherwise, $g_t(x) \in B'$ (see below (3.49)),
and 
$$
|\wt\varphi_t(x)-g_t(x)| = |\pi_j(g_t(x))-g_t(x)| 
\leq 2C \dist(g_t(x),L_j) \leq 4C \varepsilon_1
\leqno (3.70)
$$
by (3.54), because $\pi_j(z)=z$ on $L_j$ and $\pi_j$ is 
$C$-Lipschitz, and by (3.55). Then $\wt\varphi_t(x) \in \wt B$,
as needed.

So we may assume that $x\in E_{j+1}$. If $x\in E_{j+1} \sm B''$, 
(3.66) says that $\wt\varphi_t(x) = x \in B'' \i \wt B$. 
If instead $x\in E_{j+1} \cap B''$, first notice that 
$\varphi_t(x) \in B''$ (by (1.6) or (1.5)). 
If in addition $d_t(x) \geq \varepsilon_2$, then 
$\xi\circ d_t(x) = 0$ by (3.56) and 
$\wt\varphi_t(x) = \varphi_t(x) \in B''$ by (3.58). 
Otherwise 
$$
|\wt\varphi_t(x)-\varphi_t(x)| 
\leq |\pi_j(\varphi_t(x))-\varphi_t(x)| 
\leq 2C\dist(\varphi_t(x),L_j) 
\leq 2C d_t(x)
\leq 2C\varepsilon_2
\leqno (3.71)
$$  
by (3.58), again because $\pi_j(z)=z$ on $L_j$, and by (3.57).
If $\varepsilon_2$ is small enough, depending on $\varepsilon_1$,
we get that $\wt\varphi_t(x) \in \wt B$ because $\varphi_t(x) \in B''$.
So (3.69) and (1.6) hold.

We finally check (1.7). So we pick $0 \leq i \leq j_{ max}$
and want to show that 
$$
\wt\varphi_t(x) \in L_i 
\ \hbox{ for } x\in E_j \cap L_i.
\leqno (3.72)
$$
We start with the case when $x\in E_{j+1}$.
Notice that $\varphi_t(x) \in L_i$, by (1.7) and 
because $x\in E_{j+1} \cap L_i$. Let $F$ be a face of 
$L_i$ that contains $\varphi_t(x)$; then 
$\pi_j(\varphi_t(x)) \in F$ because Lemma 3.4 says that
$\pi_j$ respects the faces of all dimensions.
Now (3.58) says that $\wt\varphi_t(x)$ is a convex 
combination of $\varphi_t(x)$ and $\pi_j(\varphi_t(x))$,
hence $\wt\varphi_t(x) \in F \i L_i$, as needed.

Now suppose that $x\in Z \cap L_i$. By (3.33), 
$$
x \in L_j \sm \bigcup_{k \geq j+1} L_k
\leqno (3.73)
$$
so we have that $i \leq j$.
In addition, if $L_i \cap L_j$ were a proper subset of
$L_j$, our completeness assumption (2.15) 
would say that it is one of the $L_k$, and
since we enumerated our boundaries in nonincreasing order
(see above (3.28)), we would get that $k \geq j+1$, a contradiction 
with (3.73) since $z\in L_k$. So $L_i \cap L_j = L_j$, i.e., $L_j \i L_i$
and it is enough to check that $\wt\varphi_t(x) \in L_j$.
But (3.54) says that $\wt\varphi_t(x) = \pi_j(g_t(x))$,
which lies in $L_j$ by definition of $\pi_j$. This completes our
verification of (1.7).

We already checked (1.4), (1.5), (1.6), and (1.8) before, 
so this completes the verification that $F=\wt\varphi_1(E_j)$ 
is a sliding competitor for $E_j$ in $\wt B$. 
Recall that we may take $\varepsilon_1$ and $r'-r$, and hence 
also $\wt r -r$ in (3.68), as small as we want. 
Now we apply our induction assumption (3.31) and Definition 2.3 
to get that
$$
\H^d(\wt W) \leq M \H^d(\wt\varphi_1 (\wt W)) + \wt r^d h,
\leqno (3.74)
$$
where 
$$
\wt W= \big\{ x\in E_{j} \, ; \, \wt\varphi_1(x) \neq x \big\}.
\leqno (3.75)
$$
We are interested in $\H^d(W)$, where
$$
W = \big\{ x\in E_{j+1} \, ; \, \varphi_1(x) \neq x \big\},
\leqno (3.76)
$$
because we want an analogue of (3.74) for $\varphi_1$,
and we'll cut $W$ into pieces. Write $E_{j+1} = A_0 \cup A$,
where
$$
A_0 = A_0(\varepsilon_2) = \big\{ x\in E_{j+1} \, ; \, 
0 < d_1(x) < \varepsilon_2 \big\}
\leqno (3.77)
$$
and $A = E_{j+1} \sm A_0$. Observe that 
$$
\lim_{\varepsilon_2 \to 0} \H^d(A_0(\varepsilon_2)) = 0
\leqno (3.78)
$$
because $A_0(\varepsilon_2) \i E_{j+1} \cap \wt B$ for $\varepsilon_2$
small, $\H^d(E_{j+1} \cap \wt B) < + \infty$, and because the monotone 
limit of the sets $A_0(\varepsilon_2)$ (when $\varepsilon_2$ tends to $0$) 
is the empty set. So, given $B'$ and a small $\varepsilon_3 > 0$, 
we know that
$$
\H^d(A_0(\varepsilon_2)) < \varepsilon_3
\leqno (3.79)
$$
for $\varepsilon_2$ small enough. Next let us check that
$$
\wt\varphi_1(x) = \varphi_1(x) \ \hbox{ for } x\in A.
\leqno (3.80)
$$
Write $A = A_1 \cup A_2$, with
$$
A_1 = \big\{ x\in E_{j+1} \, ; \, d_1(x) =0 \big\}
\ \hbox{ and } \ 
A_2 = \big\{ x\in E_{j+1} \, ; \, d_1(x) \geq \varepsilon_2 \big\}.
\leqno (3.81)
$$
When $x\in A_1$, $\xi\circ d_1(x) = 1$ by (3.56), so
$\wt\varphi_1(x) = \pi_j(\varphi_1(x))$ by (3.58). 
In addition, $\varphi_1(x)\in L_j$, because $d_1(x) =0$ and by (3.57),
and $\pi_j(\varphi_1(x)) = \varphi_1(x)$ by (3.26) and (3.6), 
as needed for (3.80). When $x\in A_2$, $\xi(x) = 0$ and (3.58) 
directly yields that $\wt\varphi_1(x) = \varphi_1(x)$.
So (3.80) holds. 

If $x\in W\cap A$, then $\wt\varphi_1(x) = \varphi_1(x) \neq x$
by (3.80) and (3.76); hence $W\cap A\i \wt W$ and
$$\eqalign{
\H^d(W) & 
\leq \H^d(A_0) + \H^d(W\cap A) 
\leq \varepsilon_3 + \H^d(W\cap A) 
\cr&
\leq \varepsilon_3 + \H^d(\wt W)
\leq \varepsilon_3 + M \H^d(\wt\varphi_1 (\wt W)) + \wt r^d h
}\leqno (3.82)
$$
because $W \i E_{j+1}$ and $E_{j+1} = A_0 \cup A$ by definition, 
by (3.79), and by (3.74).

Next we estimate $\H^d(\wt\varphi_1 (\wt W))$. Notice that
$E_j = Z \cup E_{j+1} = Z \cup A_0 \cup A$ by (3.33), so
$$
\wt W \i Z \cup \big(\wt W \cap [A_0 \cup A] \big)
\leqno (3.83)
$$
since $\wt W \i E_{j}$. First,
$$
\H^d(\wt\varphi_1 (Z)) = 0
\leqno (3.84)
$$
because $\wt\varphi_1$ is Lipschitz and $Z \i E \sm E^\ast$ 
is negligible (see (3.33), (3.28), and (3.29)).
Next, $\wt\varphi_1(x) = \varphi_1(x)$ on $A$ (by (3.80)), 
so $\wt W \cap A = W \cap A$ and $\wt \varphi_1(\wt W \cap A)
= \varphi_1(W \cap A)$, hence
$$\eqalign{
\H^d(\wt\varphi_1(\wt W)) &= \H^d(\wt\varphi_1(\wt W \cap [A_0\cup A])
\leq \H^d(\wt\varphi_1(\wt W\cap A_0))
+ \H^d(\varphi_1 (W\cap A)) 
\cr&\leq \H^d(\wt\varphi_1(A_0)) + \H^d(\varphi_1 (W)).
}\leqno (3.85)
$$
by (3.83) and (3.84). We have no nice formula for $\wt\varphi_1$
on $A_0$, but let us check that
$$
\wt\varphi_1 \ \hbox{ is $C$-Lipschitz on } A_0,
\leqno (3.86)
$$
with a constant $C$ that may be enormous and depend on 
various Lipschitz constants (in particular for $\varphi_1$ 
and $\pi_j$), but does not depend on $\varepsilon_2$.
Recall that on $A_0$, $\wt\varphi_1$ is given
by (3.58), i.e., 
$$\eqalign{
\wt\varphi_1(x) 
&= \xi \circ d_1(x) \,\pi_j(\varphi_1(x)) 
+(1- \xi \circ d_1(x)) \,\varphi_1(x)
\cr&
= \varphi_1(x) + \xi \circ d_1(x) [\pi_j(\varphi_1(x))-\varphi_1(x)],
}\leqno (3.87)
$$
where $\xi$ and $d_1$ are still given by (3.56) and (3.57),
and only the variations of $\xi\circ d_1(x)$ will be dangerous here
(because they could involve some $\varepsilon_2^{-1}$).
Write, for $x, y \in A_0$,
$$\eqalign{
|\wt\varphi_1(x) - \wt\varphi_1(y)|
&\leq |\varphi_1(x)-\varphi_1(y)| 
+ |\xi \circ d_1(x)-\xi \circ d_1(y)| |\pi_j(\varphi_1(x))-\varphi_1(x)|
\cr&\hskip 2cm
+ \xi \circ d_1(y) \big|\pi_j(\varphi_1(x))-\varphi_1(x)-
\pi_j(\varphi_1(y))+\varphi_1(y) \big|
\cr&
\leq C |x-y| + |\xi \circ d_1(x)-\xi \circ d_1(y)| 
|\pi_j(\varphi_1(x))-\varphi_1(x)|
}\leqno (3.88)
$$
because all the other functions are Lipschitz with estimates that do 
not depend on $\varepsilon_2$. By (3.56), (3.57), and because 
$\varphi_1$ is $C$-Lipschitz,
$$
|\xi \circ d_1(x)-\xi \circ d_1(y)|
\leq 2 \varepsilon_2^{-1} |d_1(x)-d_1(y)|
\leq C \varepsilon_2^{-1} |x-y|,
\leqno (3.89)
$$
while
$$
|\pi_j(\varphi_1(x))-\varphi_1(x)|
\leq C \dist(\varphi_1(x),L_j) \leq C d_1(x) \leq C \varepsilon_2
\leqno (3.90)
$$
because $\pi_j(z)=z$ on $L_j$, by (3.57), and because $x\in A_0$
(see the definition (3.77)). Altogether,
$|\wt\varphi_1(x) - \wt\varphi_1(y)| \leq C |x-y|$
by (3.88), (3.89), and (3.90); this proves (3.86). Now
$$
\H^d(\wt\varphi_1(A_0)) \leq C \H^d(A_0) \leq C \varepsilon_3
\leqno (3.91)
$$
by (3.79), and 
$$\eqalign{
\H^d(W) & 
\leq \varepsilon_3 + M \H^d(\wt\varphi_1 (\wt W)) + \wt r^d h
\cr&
\leq \varepsilon_3 + M \H^d(\wt\varphi_1(A_0)) + M \H^d(\varphi_1 (W)) 
+ \wt r^d h
\cr&
\leq M \H^d(\varphi_1 (W))  + \wt r^d h + (1+M C) \varepsilon_3
}\leqno (3.92)
$$
by (3.82), (3.85), and (3.91). Recall that $\wt r$ can be chosen as 
close to $r$ as we want, and that $\varepsilon_3$ can be chosen 
arbitrarily small. So we get that 
$\H^d(W) \leq M \H^d(\varphi_1 (W))+r^d h$.
That is, (2.5) holds. This completes our proof of (3.32) given (3.31),
and at the same time our proof of Proposition 3.27 by induction
(recall that $E_{j_{ max}+1} = E^\ast$, see below (3.30)).
\qed

\ms\noindent{\bf Proof of Proposition 3.3.}
We now assume that the Lipschitz assumption is satisfied on the 
open set $U$, and want to check that $E^\ast \in GSAQ(U, M, \delta, h)$
as soon as $E \in GSAQ(U, M, \delta, h)$. We cannot use 
Proposition 2.8 directly, because it would give us bad constants,
but we can change variables and apply the proof above. That is,
let $\lambda > 0$ and $\psi : \lambda U \to B(0,1)$ be as in
Definition 2.7, and then define $\wt \psi$ by
$\wt \psi(x) = \psi(\lambda x)$  and set $F = \wt\psi(E)$. This is a
quasiminimal set in $B(0,1)$ (by Proposition 2.8), but we don't
care so much, and the closed support of the 
restriction of $\H^d$ to $F$ is $F^\ast = \wt\psi(E^\ast)$.
We construct a nonincreasing sequence of sets $F_j$, 
$0 \leq j \leq j_{ max}+1$, as we did near (3.30), and we set
$E_j = \wt\psi^{-1}(F_j)$. Thus $E_0 = E$ and $E_{j_{ max}+1} = E^\ast$, 
and we want to show by induction that $E_j \in GSAQ(U, M, \delta, h)$. 

This is the case for $j=0$, and for the induction step, 
we give ourselves a sliding competitor $\varphi_1(E_{j+1})$ for $E_j$ 
in some ball $B$. We consider the mappings 
$f_t = \wt\psi \circ \varphi_t \circ \wt\psi^{-1}$
on $F_{j+1}$, which define a sliding competitor for $F_{j+1}$
in $H = \wt\psi(B)$. Of course $H$ is not a ball, but
we don't really care, we still can use the proof of 
Proposition 3.27 to construct mappings $\wt f_t$, that define a
sliding competitor for $F_{j}$ in a set $H'$ which is just a tiny bit
little larger than $H$. Then the
$\wt\varphi_t = \wt\psi^{-1} \circ \wt f_t \circ \wt\psi$
define a sliding competitor for $E_j$, in the set $\wt\psi^{-1}(H')$
which is a tiny bit larger than $B$. We apply the definition (2.5)
to this competitor and get that
$$
\H^d(\wt W) \leq M \H^d(\wt\varphi_1 (\wt W)) + \wt r^d h,
\leqno (3.93)
$$
where $\wt r$ is the radius of a ball that contains $\wt\psi^{-1}(H')$
(and can be taken as close to $r$ as we want), and
$$
\wt W= \big\{ x\in E_{j} \, ; \, \wt\varphi_1(x) \neq x \big\}.
\leqno (3.94)
$$
Then we estimate $\H^d(W)$, where 
$W= \big\{ x\in E_{j+1} \, ; \, \varphi_1(x) \neq x \big\}$
as we did after (3.76); the error terms, like the ones in 
(3.79), (3.84), and (3.91) become $C$ times larger because we
compose with $\wt\psi$ and $\wt \psi ^{-1}$, but the argument goes 
through.
\qed

\ms\noindent{\bf Remark 3.95.}
Here we defined the rigid assumption, and then the Lipschitz 
assumption, in terms of dyadic cubes, but we could have obtained
similar results if we used a net of convex polyhedra instead, 
with a rotundity assumption where we ask all the angles in the 
faces of all dimensions to be bounded from below. 
The only place where the argument needs to
be modified is in the proof of Lemma 3.4.
See Remark 2.12 and the comments after (3.8) and (3.13).

\ms\noindent{\bf Remark 3.96.}
In Propositions 3.3 and 3.27, we can get a slightly stronger
conclusion under the same assumption, namely that
all the closed sets $F$ such that $E^\ast \i F \i E$
lie in the same $GSAQ(U, M, \delta, h)$ as $E$.
That is, we never use the fact that $E^\ast$ is the 
closed support of $\H^d_{|E}$, but just the fact that
$E^\ast \i E$ and $\H^d(E\sm E^\ast) = 0$.

This could be useful if we tried to extend Proposition 3.27 to
a situation where we only assume that $U$ is covered by a finite
collection of domains where we have the Lipschitz assumption,
and try to go from $E$ to $E^\ast$ in a finite number of steps.
We would need to check what happens to a sliding competitor in a
ball that is not entirely contained in one domain, though. We shall 
not pursue this here, as Proposition 3.27 shall be enough for our 
purposes.

\ms
The conclusion of this section is that we shall feel free
to restrict our attention to coral quasiminimal sets, with no apparent 
loss of generality. We could probably have managed, in the rest of
this paper, not to assume that $E$ is coral, and then prove regularity 
results on $E^\ast$, and we may even do this in some cases (so as not 
to rely on the proof of Propositions 3.3 and 3.27), but we shall find it 
more comfortable to know that we can work work
directly with coral sets. Recall that we do not say that every
quasiminimal set is a competitor for its core $E^\ast$, or the other 
way around (both things are wrong in general), 
but just that they have the same minimizing properties.

\bigskip
\centerline{PART II : AHLFORS REGULARITY AND RECTIFIABILITY}
\ms
In this part we prove basic regularity properties for the (core of)
sliding quasiminimal sets. The main ones are their local Ahlfors regularity
(Proposition 4.1), which is of constant use, and rectifiability 
(Theorem 5.16), which will be important for the theorems of Part IV
on limits.

Most of the results of this part and the next one (where we prove
the local uniform rectifiability in some cases), and their proofs, 
are generalizations of results of [DS4], 
except for rectifiability which was proved in [A2] 
and ignored in [DS4] (because we thought uniform rectifiability 
was better), but for which the proof of [DS4] 
works as well. 

We cannot repeat all the arguments from [DS4] 
(this would be too long), but fortunately many of the intermediate
results there can be used essentially without modification here, and there 
are only a few places where we need to be careful, because a 
competitor for our quasiminimal set is used. We will try to give an 
idea, but all the details, of the arguments that work with only minor
modifications, and be as precise as possible on the differences, i.e., 
places where a competitor is defined. Hopefully this will make the
reading of this text not too unpleasant, probably at the price of 
often believing the author when he says that some old estimates estimates 
still apply.

\bigskip
\noindent {\bf 4. Local Ahlfors regularity of quasiminimal sets.}
\medskip

We start now our long study of regularity properties of sliding 
quasiminimal sets with the very convenient
local Ahlfors regularity of the core $E^\ast$. 
We start with the rigid case.

\ms\proclaim Proposition 4.1. For each choice of $M \geq 1$, we can 
find $h > 0$ and $C_M \geq 1$, depending on $M$ and the dimensions 
$n$ and $d$, so that the following holds. 
Suppose that $E \in GSAQ(B_0, M, \delta , h)$, where we set 
$B_0 = B(0,1)$, and that the rigid assumption is satisfied.
Let $r_0 = 2^{-m} \in (0,1]$ 
denote the side length of the dyadic cubes used to define the 
rigid assumption.
Then if $x\in E^\ast \cap B_0$ and $0 < r < \Min(r_0,\delta)$  
are such that $B(x,2r) \i B_0$, we have that 
$$
C_M^{-1} r^d \leq \H^d(E\cap B(x,r)) \leq C_M r^d.
\leqno (4.2) 
$$

\ms
Recall that $E^\ast$, the core of $E$, is as in (3.2).
We could also have assumed that $E$ is coral, and then obtained that
(4.2) holds for $x\in E \cap B_0$ (instead of $E^\ast \cap B_0$).
Also observe that we can replace $E$ with $E^\ast$ in (4.2),
since $E^\ast \i E$ and $\H^d(E \sm E^\ast) = 0$.

By Proposition 2.8 and the bilipschitz invariance of local 
Ahlfors regularity, this result implies the corresponding one 
when the bilipschitz assumption holds. See Proposition 4.74.
\ms
We want to say that the standard proof given in [DS4] 
goes through in the present setting. We cannot repeat it
entirely (this will make this paper too huge and boring), so we shall
only recall how the proof goes, and concentrate on the minor
modifications that we need to make, in particular in the choices of cubes. 

There are two main differences, compared
to the initial situation in [DS4]. 
First, the accounting in the definition of general quasiminimal
sets is slightly different (and makes the notion of
quasiminimal sets more general) than the one we used in [DS4]. 
This aspect of things is not really important, and was already 
discussed in [D5]. 
The second difference, which really concerns us here, is the additional 
conditions on the competitors that come from the sliding boundary conditions; 
we may have to go all the way to the boundary, and make sure that all the 
competitors that we use in the proof satisfy the sliding condition.
This will force us to choose more carefully the cubes where we do 
Federer-Fleming constructions, and this is why we shall be more 
prudent in the choice of integers $N_k$ below.

We start our reading of [DS4] 
with the Federer-Fleming construction of Lipschitz projections
that is described in Chapter $\ast$3 (we shall use the convention
that $\ast$ calls a reference in [DS4]). 
We need the following slight variant of Proposition $\ast$3.1.
Here and below, cubes are systematically assumed to be closed,
and the $k$-dimensional skeleton of a cube $R$ is the union of all the
faces of dimension $k$ of $R$. Thus the $1$-dimensional skeleton
of a cube in $\R^3$ is a union of $12$ line segments.

\ms\proclaim Lemma 4.3. Let $N \geq 1$  be an integer,
and let $Q \i \R^n$ be a cube of side length $N 2^k$
(for some $k\in \Z$), which is the almost-disjoint union of
$N^n$ dyadic cubes of side length $2^k$.
Denote by ${\cal R}={\cal R}(Q)$ the set of dyadic cubes of 
side length $2^k$ that are contained in $Q$,
and by ${\cal S}_d$ the union of the $d$-dimensional skeletons of 
the dyadic cubes $R \in \cal R$.
Let $E$ be a compact subset of $Q$ such that $\H^d(E) < +\infty$.
Then there is a Lipschitz mapping $\phi : \R^n \to \R^n$
such that
$$
\phi(x)=x \hbox{ for } x\in \R^n \sm Q
\hbox{ and for } x\in {\cal S}_d, 
\leqno (4.4)
$$
$$
\phi(E) \i  {\cal S}_d \cup \partial Q,
\leqno (4.5)
$$
$$
\phi(R) \i R \ \hbox{ for } R \in {\cal R},
\leqno (4.6)
$$
and
$$
\H^d(\phi(E\cap R)) \leq C \H^d(E\cap R)
\hbox{ for } R \in {\cal  R}.
\leqno (4.7)
$$
Here $C$  depends on $n$ and $d$, but not on $N$ or $E$.

\ms
The only difference with Proposition $\ast$3.1 in [DS4] 
is  that $Q$ is not required itself to be a dyadic cube
(and so $N$ is not required to be a power of $2$);
however, this fact that $Q$ is dyadic (or rather that
$N$ is a power of $2$) is never used in [DS4], and 
the proof can be carried out exactly as before.
\qed

Another consequence of the proof, which works by successive 
``radial" projections onto faces, is that in addition to (4.6),
$$
\phi(F) \i F
\hbox{ when $F$ is a face (of any dimension) of a cube
$R\in {\cal R}$}.
\leqno (4.8)
$$
This will be used to prove the stability
conditions (1.7).

\ms
We now turn to Chapter 4 in [DS4], and prove Proposition 4.1. 
We start with the easier upper bound. We want to find $C_0 \geq 1$
(depending also on $M$) such that
$$
\H^d(E \cap Q_0) \leq C_0 l_0^d
\leqno (4.9)
$$
when $Q_0$ is a cube of side length 
$$
l_0 \leq \Min(2^{-m},{\delta \over n})
\leqno (4.10)
$$
which is dyadic (in the same grid that was used in the description 
of the $L_j$ for the rigid assumption), and such that $2Q_0 \i B_0$.

Indeed, if we can prove (4.9) for such cubes, and if 
$x\in B_0$ (we do not need $x\in E^\ast$ for the lower bound) 
and $0 < r < \Min(r_0,\delta)$ are such that
$B(x,2r) \i B_0$, we can easily cover $B(x,r)$ with less than $C$
cubes $Q_0$ as above, with side lengths less than $r$, and then the upper 
bound in (4.2) follows from (4.9).

So we give ourselves a cube $Q_0$ as above, assume that (4.9) fails,
and we shall derive a contradiction if $C_0$ in (4.9) is large enough
(depending on $n$, $d$ and $M$). Here we shall only need to assume that 
$h \leq 1$. We want to construct by induction an increasing sequence 
of cubes $Q_k$, $k \geq 1$, with the same center as $Q_0$, and 
whose side lengths $l_k$ are such that
$$
l_0 < l_k < 2l_0
\ \hbox{ for } k \geq 1.
\leqno (4.11)
$$
At the same time, we shall define large integers
$N_k$, $k \geq 1$, and for each $k\geq 1$ cut $Q_k$ into $N_k^d$ 
cubes of the same side length $N_k^{-1} l_k$; we shall call
${\cal R}(Q_k)$ the collection of these smaller cubes, 
in accordance with the notation of Lemma 4.3.
When we do this, we want to make sure that for $k \geq 1$,
$$\eqalign{
&\hbox{every cube $R \in {\cal R}(Q_k)$ is a dyadic (sub)cube of the grid}
\cr& 
\hbox{that was used to define the $L_j$ in the rigid assumption.}
}\leqno (4.12)
$$
The fact that these cubes are of a smaller size than the $L_j$
follows from (4.11) and the fact that  $l_0 \leq 2^{-m}$,
because $N_k \geq 2$, but we typically want any face of $R$
that intersects the interior of a face of some $L_j$ to be entirely contained 
in that $L_j$. We require this because we want to apply Lemma 4.3
to $Q_k$ to find a competitor for $E$ in $Q_k$; for similar reasons,
we want $Q_{k-1}$ to be obtained from $Q_k$ by removing from 
${\cal R}(Q_k)$ the two exterior layers of cubes, and then taking the
union. In other words, we want to have
$Q_{k-1} = {N_k - 4 \over N_k} \, Q_{k}$, or equivalently
(since the cubes have the same center)
$$
l_{k-1} = \big(1-{4 \over N_{k}}\big) \, l_{k}
\leqno (4.13)
$$
for $k \geq 1$. In fact denote by $A(Q_{k})$ the union of
the cubes $R \in {\cal R}(Q_k)$ that lie on the exterior layer
(or equivalently that meet $\d Q_k$); we wanted to make sure that
$$
Q_{k-1} \i Q_k \sm A(Q_{k}),
\leqno (4.14)
$$
and (4.13) ensures this (we need to remove a extra layer because
$Q_{k-1}$ and $A(Q_{k})$ are both closed).

Next let us assume for the moment that we can choose $N_k$
so that (4.11) and (4.12) hold, and let us use this to control
$$
m_k = \H^d(E\cap Q_k)
\leqno (4.15)
$$
in terms of $m_{k-1}$. First apply Lemma 4.3 to $Q_k$,
the integer $N_k$, and the decomposition coming from
${\cal R}(Q_k)$. This gives a Lipschitz mapping 
$\phi : \R^n \to \R^n$, which we use to define a family $\{\varphi_t \}$
as in Definition 1.3, simply by setting 
$$
\varphi_t(x) = t \phi(x) + (1-t) x
\ \hbox{ for $x\in \R^n$ and $0 \leq t \leq 1$.}
\leqno (4.16)
$$
The properties (1.4), (1.5), (1.6), and (1.8) are easily checked, 
relative to any closed ball that contains $Q_{k}$, for instance the
ball $B$ with the same center $x_0$ as $Q_k$ and $Q_0$, and with
radius 
$$
r = {\sqrt n \, l_k \over 2} \leq \sqrt n \, l_0 < \delta
\leqno (4.17)
$$
by (4.11) and (4.10). We also need to check (1.7). Let
$x\in L_j$ be given, and let $F$ be a face of $L_j$ that contains
$x$. Let $F' \i F$ be a dyadic subface of the same dimension as $F$,
but of side length $N_k^{-1} l^k$, that contains $x$.
We know that such a face exists, because 
$N_k^{-1} l^k < 2^{-m}$ (see the line below (4.12))
and the ratio is a power of $2$. Since the cubes 
of ${\cal R}(Q_k)$ are dyadic in the same grid as $F$
(see (4.12)), we get that $F'$ is a face of the grid defined
by the cubes of ${\cal R}(Q_k)$. Now 
either $x\notin Q_k$, and then $\varphi_t(x) = \phi(x) = x$
by (4.4), or else (4.8) says that $\phi(x)$, and hence also (by 
(4.16) and the convexity of $F'$) $\varphi_t(x)$, lies on $F' \i F \i L_j$. 
This proves (1.7); hence $\phi(E)$ is a sliding competitor for $E$ in 
$B = \overline B(x_0,r)$, with $r$ as in (4.17).

Let us apply Definition 2.3. We get (2.5), with $\varphi_1 = \phi$. 
That is, 
$$
\H^d(W_1) \leq M \H^d(\phi(W_1)) + h r^d
\leq M \H^d(\phi(W_1)) + n^{d/2} l_0^d
\leqno (4.18)
$$
where $W_1 = \big\{ y \in E \, ; \phi(y) \neq y \big\}$,
if $h \leq 1$, and by (4.17).
Observe that $W_1 \i Q_k$ because $\phi(x) = x$ out of $Q_{k}$ by 
(4.4). In addition, (4.5) and (4.6) imply that
$$
\phi(x) \in \S_d \ \hbox{ for } x \in E \cap Q_{k} \sm A(Q_{k}),
\leqno (4.19)
$$
so
$$
\H^d(\phi(W_1 \sm A(Q_{k}))) \leq \H^d(\S_d) 
\leq C N^{n-d}_{k} l_0^d.
\leqno (4.20)
$$
For $E \cap A(Q_{k})$, we decompose $A(Q_{k})$ into cubes of 
${\cal R} = {\cal R}(Q_k)$ and use (4.7) to say that 
$$\eqalign{
\H^d(\phi(W_1 \cap A(Q_{k})))
&\leq \sum_{R \in {\cal R} ; R \i A(Q_{k})} \H^d(\phi(E\cap R))
\cr&
\leq C \sum_{R \in {\cal R} ; R \i A(Q_{k})} \H^d(E\cap R)
\leq 2^n C \H^d(E\cap A(Q_{k}))
}\leqno (4.21)
$$
because a given point of $E\cap A(Q_{k})$ lies in at most $2^n$
cubes $R$. So 
$$
\H^d(W_1) \leq M \H^d(\phi(W_1)) + n^{d/2} l_0^d
\leq C M \H^d(E\cap A(Q_{k})) + C (M +1) N^{n-d}_{k} l_0^d
\leqno (4.22)
$$
by (4.18), (4.20), and (4.21) (and with a constant $C$ 
that does not depend on $M$). Finally
$E \cap Q_{k} \sm A(Q_{k}) \i \S_d \cup W_1$ because
if $x\in E \cap Q_{k} \sm A(Q_{k})$ lies out of $\S_d$,
then (4.19) says that $x\in W_1$. Then
$$
\H^d(E \cap Q_{k} \sm A(Q_{k})) \leq \H^d(\S_d) + \H^d(W_1)
\leq C M \H^d(E\cap A(Q_{k})) + C (M +1) N^{n-d}_{k} l_0^d.
\leqno (4.23)
$$
We add $\H^d(E\cap A(Q_{k}))$ to both sides, recall that $M \geq 1$,
and get that
$$\eqalign{
m_k &= \H^d(E \cap Q_{k})
\leq CM \H^d(E\cap A(Q_{k})) + C'M N^{n-d}_{k} l_0^d
\cr&
\leq CM \H^d(E \cap Q_k \sm Q_{k-1}) + C'M N^{n-d}_{k} l_0^d
\cr&
= CM [m_k - m_{k-1}] + C'M N^{n-d}_{k} l_0^d
}\leqno (4.24)
$$
by (4.15), (4.14), and (4.15) again.  That is,
$$
[CM-1] m_k \geq C M m_{k-1} - C'M N^{n-d}_{k} l_0^d
\leqno (4.25)
$$
or equivalently (dividing by $CM$)
$$
m_k \Big(1- {1 \over CM} \Big) 
\geq  m_{k-1} - {C' \over C} N^{n-d}_{k} l_0^d.
\leqno (4.26)
$$
We shall choose $N_{k}$ so that
$$
{C' \over C} N^{n-d}_{k} l_0^d \leq  {m_{k-1} \over 10CM},
\leqno (4.27)
$$
so that (4.26) implies that
$$
m_k \geq \Big(1- {1 \over CM} \Big)^{-1} \Big(1 - {1 \over 10CM} 
\Big) \, m_{k-1}
\geq m_{k-1} \Big(1 + {1 \over 10CM} \Big)
\leqno (4.28)
$$
and, by induction,
$$
m_k \geq \Big(1 + {1 \over 10CM} \Big)^k m_0
\geq C_0 l_0^d \Big(1 + {1 \over 10CM} \Big)^k
\leqno (4.29)
$$
by (4.15) and because (4.9) is assumed to fail.

Now we want to check that we can choose $C_0$ (large enough) 
and the $N_k$ (by induction), so that (4.11), (4.12), and (4.27)
hold. The desired contradiction will follow, because the fact that
$m_k = \H^d(E\cap Q_k) \leq \H^d(E\cap 2Q_0) < +\infty$ 
(by (4.15), (4.11), and the finite measure condition in 
Definition 2.3) contradicts (4.29) for $k$ large.

In fact, we shall choose the $N_k$, $k \geq 1$, 
so that the following constraints hold. Set
$$
\lambda_k = {l_0^{-d} m_{k} \over 10C'M} 
\leqno (4.30)
$$
for $k \geq 0$, with $C'$ as in (4.27). Thus (4.27) just 
demands that $N^{n-d}_{k+1} \leq \lambda_k$ for $k \geq 0$,
but we shall pick the $N_{k}$ so that
$$
{\lambda_k \over 3^{n-d}} \leq N_{k+1}^{n-d} \leq \lambda_k
\ \hbox{ for } k \geq 0,
\leqno (4.31)
$$
and also such that for $k \geq 0$,
$$
N_{k+1} \geq N_k + 4 
\ \hbox{ and } \ 
{N_{k+1}-4 \over N_{k}} \hbox{ is a power of $2$,} 
\leqno (4.32)
$$
where we set $N_0 = 1$ for $k=0$.

First we want to check that (4.12) for $k+1$ follows from this
and (if $k \geq 1$) (4.12) for $k$.
We do not need to check (4.12) for $k=0$, but recall that $Q_0$ 
was assumed to be dyadic in the usual grid (the one that was used to 
define the $L_j$ in the rigid assumption). 
When $k \geq 1$, denote by $s_k$ the common side length of the cubes of 
${\cal R}(Q_{k})$; thus $s_k = N_k^{-1} l_k$.
Also set $s_0 = l_0$ (the side length of $Q_0$).
Then, for $k \geq 0$,
$$
s_{k+1} = {l_{k+1} \over N_{k+1}}
= \big(1-{4 \over N_{k+1}}\big)^{-1} {l_{k} \over N_{k+1}}
= {l_k \over N_{k+1}-4}
= {N_k \over N_{k+1}-4} \, s_k
\leqno (4.33)
$$
by (4.13). We know (by definition of $Q_0$ or by induction assumption)
that $s_k$ is a dyadic number and $s_k \leq 2^{-m}$ (the size
of the dyadic cubes that we used to define the $L_j$), and now
(4.32) says that $s_{k+1} \leq s_k$ and is a dyadic number too.
So the cubes of ${\cal R}(Q_{k+1})$ have the right size;
we also need to know that they are dyadic (instead of merely 
translations of dyadic cubes), and for this we use the fact that, as
was observed above (4.13), $Q_k$ is obtained, from the 
decomposition of $Q_{k+1}$ into cubes of side length $s_{k+1}$, 
by removing the two exterior layers of cubes. 
By induction assumption, the cubes of side length $s_k$ that compose $Q_k$ 
are dyadic, hence this is also true for the cubes of side length
$s_{k+1}$ that compose $Q_{k+1}$.
This proves (4.12) (if we have (4.32)).

Now let $k \geq 0$ be given, assume that the $N_l$, 
$1 \leq l \leq k$, were chosen so that (4.11), (4.12),
(4.31), and (4.32), hold for $1 \leq l \leq k$
(no condition if $k=0$), and let us choose $N_{k+1}$.
We first check that
$$
\lambda_k \geq (N_k + 4)^{n-d}.
\leqno (4.34) 
$$
Observe that (4.29) holds (because if $k \geq 1$,
(4.27) follows from (4.31) for $k-1$), so
$$
\lambda_k = {l_0^{-d} m_{k} \over 10C'M} 
\geq {C_0 \over 10C'M}\, \Big(1 + {1 \over 10CM} \Big)^k.
\leqno (4.35)
$$
By taking $C_0$ large enough, we can thus make sure that 
$\lambda_k \geq 300^{n-d}$, for instance, hence (by (4.31) for $k-1$)
$N_k \geq 100$ if $k \geq 1$. Note that (4.34) follows trivially
from (4.35) (or directly from (4.30) and the failure of (4.9)) when $k = 0$.
Otherwise, the definition (4.30) and (4.28) say that
$$
\lambda_k = {m_k \over m_{k-1}} \lambda_{k-1}
\geq \Big(1 + {1 \over 10CM} \Big) \lambda_{k-1}
\geq \Big(1 + {1 \over 10CM} \Big) N_k^{n-d}
\leqno (4.36)
$$
by (4.31), and (4.34) will follow as soon as we check that 
$1 + {1 \over 10CM} \geq \big({N_k +4 \over N_k}\big)^{n-d}
= \big(1 + {4 \over N_k}\big)^{n-d}$, or equivalently
${4 \over N_k} \leq \big(1 + {1 \over 10CM} \big)^{1/(n-d)} - 1$.
But (4.31) for $k-1$ and (4.35) say that
$$\eqalign{
4N_k^{-1} &\leq 12 \lambda_{k-1}^{-1/(n-d)} 
\leq 12 \Big( {C_0 \over 10C'M} \Big)^{-1/(n-d)}
\Big(1 + {1 \over 10CM} \Big)^{-(k-1)/(n-d)}
\cr&
\leq 12 \Big( {C_0 \over 10C'M} \Big)^{-1/(n-d)}
< \big(1 + {1 \over 10CM} \big)^{1/(n-d)} - 1
}\leqno (4.37)
$$
if $C_0$ is large enough (depending on $M$, $n$, and $d$),
so (4.34) holds.

Because of (4.34), picking $N_{k+1} = N_k +4$ would already
yield the second half of (4.31). We take for $N_{k+1}$
the largest integer $N$ such that $N \geq N_k +4$,
${N-4 \over N_k}$ is a power of $2$, and $N^{n-d} \leq \lambda_k$
(the second half of (4.31)).
We know from (4.34) that $N_{k+1} \geq N+4$, and so (4.32) holds.
The second half of (4.31) holds by definition. By maximality
of $N_{k+1}$, $N = 2N_{k+1} -  4$ does not work. Since
${N-4 \over N_k} = 2 {N_{k+1}-4 \over N_k}$ is also a power of $2$,
this means that $(2N_{k+1} -  4)^{n-d} > \lambda_k$, which implies 
that $N_{k+1}^{n-d} >{\lambda_k \over 3^{n-d}}$ because $N_{k+1} > 
N_k$ is large, and as needed for (4.31).

We now check that (4.11) holds for $k+1$ with our choice 
of $N_{k+1}$. By repeated uses of (4.13),
$$
l_{k+1} = l_0 \prod_{1 \leq j \leq k+1} \Big(1-{4\over N_j} \Big)^{-1}
\leqno (4.38)
$$
so it is enough to show that 
$\sum_{1 \leq j \leq k+1} {1\over N_j} \leq 10^{-2}$, say,
which follows from the first line of (4.37) and its analogue for 
$j \leq k$, provided that $C_0$ is chosen large enough. This completes the 
verification and the definition of the $N_k$; the expected 
contradiction follows, and shows that (4.9) holds. 
The upper bound in (4.2) follows, as explained below (4.9).

\ms
Next we want to establish the lower bound in
Proposition 4.1. The main step will be the following.

\ms \proclaim Lemma 4.39. Let $a < 1$ be given.
There are constants $\varepsilon_a$ and $C_a$,
that depend on $n$, $a$, and $M$, but not on $h$ or $\delta$, 
with the following property. 
Let $E$ be as in Proposition 4.1, and let $Q$ be a cube
such that $2Q \i B_0$, whose side length $l(Q)$ is such that
$$
l(Q) \leq \min({\delta \over n}, 2^{-m}),
\leqno (4.40)
$$
and for which
$$
\H^d(E \cap Q) \leq \varepsilon_a l(Q)^{d}.
\leqno (4.41)
$$
Then 
$$
\H^d(E \cap {1 \over 100} Q) \leq a \H^d(E\cap Q) + C_a h l(Q)^{d}.
\leqno (4.42)
$$

\ms
Lemma 4.39 will be proved soon, but let us first check that it yields
the lower bound in Proposition 4.1. Let $x\in E^\ast \cap B_0$ and
$r$ be as in the proposition. We know from general geometric
measure theory that there is a constant $c > 0$ 
(depending at most on $n$) such that
$$
\limsup_{\rho \to 0} \rho^{-d }\H^d(E\cap B(x',\rho)) \geq c
\leqno (4.43)
$$
for $\H^d$-almost every point $x'\in E$. 
See for instance [Ma], Theorem 6.2 on page 89. 
Since $x\in E^\ast$, (3.2) says that $\H^d(E\cap B(x,t)) > 0$
for all $t>0$, so we can choose $x'\in E$, very close to $x$, 
such that (4.43) holds.

We choose $a = 200^{-d}$, and Lemma 4.39 yields
constants $\varepsilon_a$ and $C_a$. We can
safely assume that $\varepsilon_a \leq 400^{-d} c$.
Next we choose $h$ so small that
$C_a h < 200^{-d}\varepsilon_a$ in (4.42). The point is that
if $Q$ satisfies the assumptions of Lemma 4.39, then
$$
\H^d(E \cap {1 \over 100} Q) 
\leq 200^{-d} \H^d(E\cap Q) + 200^{-d}\varepsilon_a l(Q)^{d}
\leq 100^{-d}\varepsilon_a l(Q)^{d}
\leqno (4.44)
$$
by (4.42), so ${1 \over 100}Q$ also satisfies the assumptions 
of Lemma 4.39, and recursively
$$
\H^d(E \cap {1 \over 100^k} Q) \leq 100^{-kd}\varepsilon_a l(Q)^{d}
\leqno (4.45)
$$
for $k \geq 0$, by (repeated uses of) (4.44). If $Q$ is centered at $x'$,
this implies that for $100^{-k-1}l(Q) \leq 2\rho \leq 100^{-k}l(Q)$,
$$
\H^d(E\cap B(x',\rho)) \leq \H^d(E \cap {1 \over 100^k} Q)
\leq 100^{-kd}\varepsilon_a l(Q)^{d}
\leq 200^{d} \varepsilon_a \rho^d \leq c \rho^d/2,
\leqno (4.46)
$$
which is incompatible with (4.43). 

Now we may try this with the largest cube $Q$ centered at $x'$,
and such that $2Q \i B(x,r)$ and (4.40) holds. Notice that $l(Q)$
is then comparable to $r$, 
(because $r < \Min(r_0,\delta) = \Min(2^{-m},\delta)$)
and since (4.41) fails by the discussion 
above, we get that
$$
\H^d(E \cap B(x,r)) \geq \H^d(E \cap Q) \geq \varepsilon_a l(Q)^{d}
\geq C^{-1} r^d,
\leqno (4.47)
$$
as needed for (the lower bound in) (4.2). Since we already established
the upper bound, Proposition 4.1 will follow from Lemma 4.39.

\ms
We now prove the lemma. We are given a cube $Q$, and we first 
reduce to dyadic cubes of the usual grid (the one that was used
to define the $L_j$).
Let $l_0$ denote the largest dyadic number such that $l_0 \leq l(Q)/2$.
Then $l_0 \leq 2^{-m}$ by (4.40). Denote by $Q'_0$ any dyadic cube
of side length $l_0$ in the usual grid, and then let $Q_0$
be a translation of $Q_0'$ by an element of $2^{-3}l_0 \Z^{n}$.
We choose $Q_0$ such that, if $x_0$ and $x_Q$ denote the centers
of $Q_0$ and $Q$, the size of every coordinate of $x_0 - x_Q$
is at most $2^{-4}l_0$. Then
$$
{1 \over 100} Q \i {1 \over 2} Q_0 \i Q_0 \i Q.
\leqno (4.48)
$$
Thus it will be enough to show that 
$$
\H^d(E \cap {1 \over 2} Q_0) \leq a \H^d(E\cap Q_0) + C_a h l_0^{d},
\leqno (4.49)
$$
because (4.42) will follow at once.

We shall now proceed a little bit as for the upper bound,
and define by induction a decreasing sequence of concentric cubes $Q_k$,
$k \geq 0$, such that 
$$
m_k = \H^d(E \cap Q_k)
\leqno (4.50)
$$
is rapidly decreasing. That is, up until we stop the process, which 
will in fact happen after a finite number of steps. We take 
$$
Q_{k+1} = \big(1-{6 \over N_{k}} \big) Q_k
\leqno (4.51)
$$
for $k \geq 0$, where $N_{k}$ is a large number that will be chosen 
later. The main point is that for $k \geq 0$, $Q_{k+1}$ is obtained from $Q_k$ 
by the following manipulation.
First we cut $Q_k$ into $N_{k}^n$ equal cubes $R$, $R \in {\cal R}(Q_k)$,
as we did for Lemma~4.3 (with $N=N_{k}$). Then we remove the three 
exterior layers, and $Q_{k+1}$ is the closure of what remains.
We want to make sure that (as in (4.12)), 
$$
\hbox{every cube $R \in {\cal R}(Q_k)$ is a 
dyadic cube of the usual grid,}
\leqno (4.52)
$$
and because of this we require (a little as in (4.32)) that
$N_0$ be a (large) power of $2$ and, for $k \geq 0$,
$$
N_{k+1} \geq N_k - 6 \geq 1
\ \hbox{ and } \ 
{N_{k+1}\over N_{k}-6} \hbox{ is a power of $2$.} 
\leqno (4.53)
$$
As it turns out, we will only need to define a finite sequence
of numbers $N_k$ (after which we shall stop), 
and we shall even manage to take $N_{k+1} = N_k - 6$
for $k \geq 0$, but let us pretend we need more generality and deal 
with (4.53) for the moment.
Let us first check that (4.52) follows if we apply the rule (4.53).
Denote by $s_k$ the common side length of the cubes of 
${\cal R}(Q_{k})$. Thus $s_{k} = N_{k}^{-1} l_k$.
First, $s_0 = N_{0}^{-1} l_0$ is a dyadic number, because $N_0$ is 
large dyadic and $l_0$ is dyadic. For $k \geq 0$,
$$
s_{k+1} = N_{k+1}^{-1} l_{k+1} 
=  N_{k+1}^{-1} \big(1-{6 \over N_{k}}\big) l_{k}
={ N_k - 6 \over N_{k+1}} \, s_{k}
\leqno (4.54)
$$
by (4.51) and because $s_{k} = N_{k}^{-1} l_{k}$. Thus 
$s_{k+1}$ is dyadic if $s_k$ is dyadic, and $s_{k+1} \leq s_k$
by the first part of (4.53). The verification of 
(4.52), i.e., the fact that the cubes also match the dyadic grid
(instead of just having the right size) is now easy, and goes as for 
(4.12) near (4.33).

We apply Lemma~4.3 to $Q_k$ (decomposed as the union of the 
cubes $R \in {\cal R}(Q_k)$), and get a Lipschitz mapping $\phi : \R^n \to \R^n$ 
which preserves the faces as in (4.8). 
This time, we do not use the function $\phi$ directly to produce a 
competitor, but instead try to project once more on $(d-1)$-faces
when this is possible. Suppose that 
$$
m_k \leq c N_k^{-d} l_k^d,
\leqno (4.55)
$$
where $c > 0$ will be chosen soon. Notice that
$\H^d(\phi(E\cap R)) \leq C \H^d(E\cap R) \leq C m_k$
for $R\in {\cal R}(Q_k)$, by (4.7). Hence
$$
\H^d(R \cap \phi(E)) \leq 2^n C m_k
\ \hbox{ for } R \in {\cal R}(Q_k)
\leqno (4.56)
$$
because if $y=\phi(x)$ lies on $R \cap \phi(E)$,
then (4.4) says that $x\in Q_k$, and then, by (4.6),
$x$ lies in $R$ or one of its neighbors of
${\cal R}(Q_k)$. We take $c$ smaller than $10^{-n} C^{-1}$, and this way
(4.56) says that $\phi(E)$ never gets close to filling the central
part of a $d$-dimensional face of a cube $R \in {\cal R}(Q_k)$. In this 
case, the proof of Lemma~4.3 (where we just do an additional 
Federer-Fleming projection on the interior cubes) says that
we can obtain a new mapping $\phi$ such that, in addition to the
properties above, $\phi(E) \cap [Q_k \sm A(Q_k)]$ is contained in a
$(d-1)$-skeleton, where $A(Q_{k})$ still denotes the union of the cubes 
$R \in {\cal R}(Q_k)$ that lie in the exterior layer. Compare with (4.5), 
and see the discussion in [DS4], below ($\ast$4.22). Thus 
$$
\H^d(\phi(E) \cap [Q_k \sm A(Q_k)]) = 0.
\leqno (4.57)
$$
We may now use this $\phi$ to define a family $\{\varphi_t \}$ of
mappings, by the same formula (4.16) as before. 
The properties (1.4)-(1.8) are verified as before (below (4.16)), 
using in particular (4.8), (4.40), and the fact that the cubes 
$R \in {\cal R}(Q_k)$ lie in the usual dyadic grid. 
So $\phi(E)$ is a sliding competitor for $E$ in some 
ball of radius $r = l_0 \sqrt n/2$. 
We can apply (2.5), and we get that
$$
\H^d(W_1) \leq M \H^d(\phi(W_1)) + h r^d
\leqno (4.58)
$$
as in the first part of (4.18), and with
$W_1 = \big\{ x \in E \, ; \phi(x) \neq x \big\}$.
Since $\phi(x) = x$ out of $Q_k$ by (4.4), we get that
$W_1 \i Q_k$, and since $\phi(Q_k) \i Q_k$ by (4.6),
$\phi(W_1) \i Q_k$ too. Hence
$$
\H^d(\phi(W_1)) = \H^d(Q_k \cap \phi(W_1))
= \H^d(A(Q_k) \cap \phi(W_1))
\leqno (4.59)
$$
by (4.57).
Denote by $A_1(Q_k)$ the union of the two exterior layers of
cubes $R \in {\cal R}(Q_k)$. Then
$$
A(Q_k) \cap \phi(W_1) 
\i \bigcup_{R \in {\cal R}(Q_k) \, ; \, R \i A_1(Q_k)} 
\phi(R \cap W_1)
\leqno (4.60)
$$
by (4.4) and (4.6), hence
$$\eqalign{
\H^d(\phi(W_1)) &= \H^d(A(Q_k) \cap \phi(W_1))
\leq \sum_{R \in {\cal R}(Q_k) \, ; \, R \i A_1(Q_k)} 
\H^d(\phi(R \cap W_1)) 
\cr&\leq C \sum_{R \in {\cal R}(Q_k) \, ; \, R \i A_1(Q_k)}
\H^d(R \cap E)
\leq C \H^d(A_1(Q_k) \cap E)
}\leqno (4.61)
$$
by (4.59) and (4.7). 
Note that $Q_{k+1} \i Q_k \sm A_1(Q_k)$ by
construction; also, $\H^d$-almost every point
$x\in E \cap Q_k \sm A_1(Q_k)$ lies in $W_1$, by (4.57)
(just notice that if $x \notin W_1$, then $x = \phi(x) \in 
\phi(E) \cap [Q_k \sm A(Q_k)]$); hence 
$$
\H^d(E \cap Q_{k+1}) \leq \H^d(W_1) \leq M \H^d(\phi(W_1)) + h r^d
\leq CM \H^d(A_1(Q_k) \cap E) + h r^d
\leqno (4.62)
$$
by (4.58) and (4.61). Next,
$$
\H^d(A_1(Q_k) \cap E) \leq \H^d(E \cap Q_k \sm Q_{k+1}) 
= m_k - m_{k+1}
\leqno (4.63)
$$
by (4.50), so (4.62) says that
$$
m_{k+1} \leq CM (m_k - m_{k+1}) + h r^d,
\leqno (4.64)
$$
hence
$$
m_{k+1} \leq {CM \over 1+ CM} \, m_k + {h r^d \over 1+ CM}.
\leqno (4.65)
$$

Now we want choose the $N_k$ and check the various constraints.
We shall first choose $N_0$ dyadic and very large, depending on 
$a$ and $M$. We also set $N = N_0/16$ and
$$
N_k = N_0 - 6k
\ \hbox{ for } 1 \leq k \leq N.
\leqno (4.66)
$$
This is probably far from optimal, but it will work.
Observe that (4.53) is then satisfied, and in fact all 
the sets ${\cal R}(Q_k)$ that we construct will be composed
of dyadic cubes of the same side. Our last cube is
$$
Q_{N+1} = \Big[\prod_{j=0}^N \big(1-{6 \over N_{j}} \big) \Big] Q_0
\leqno (4.67)
$$
by (4.51). Since ${6 \over N_j} \leq {6 \over N_{N}}
\leq {10 \over N_{0}}$ by definition of $N$,
$$
\sum_{j=0}^N \log\big(1-{6 \over N_{j}} \big)
\geq (N+1) \log\big(1-{10 \over N_0} \big)
\geq -{1 \over 10}
\leqno (4.68)
$$
if $N_0$ is large enough, and hence (because $e^{-1/10} > 1/2$)
$Q_{N+1}$ contains ${1 \over 2} Q_0$. Thus
$$
\H^d(E \cap {1 \over 2} Q_0) \leq m_{N+1}.
\leqno (4.69)
$$
Now we check (4.55). Observe that for $k \geq 0$,
$$
\eqalign{
N_k^{d} m_k l_k^{-d} &\leq N_0^{d} m_0 l_k^{-d}
\leq 2^d N_0^{d} m_0 l_0^{-d}
\leq 2^d N_0^{d} \H^d(E \cap Q) l_0^{-d}
\cr&
\leq 2^d N_0^{d} \varepsilon_a l(Q)^d l_0^{-d}
\leq 8^{d} \varepsilon_a N_0^{d}
}\leqno (4.70)
$$
because $N_k \leq N_0$, $m_k \leq m_0$, $l_k \geq l_0/2$
(since $Q_{N+1}$ contains ${1 \over 2} Q_0$), by (4.48) and 
(4.50), by (4.41), and because $l_0 \geq l(Q)/4$ by definition 
of $l_0$ (below (4.47)). If $\varepsilon_a$ is chosen small enough, 
depending on $M$ and $a$ through $N_0$
(see near (4.72) below for the choice of $N_0$), (4.70) implies (4.55),
we can proceed as above, and (4.65) holds for $0 \leq k \leq N$.
That is, if we set $\rho = {CM \over 1+ CM} < 1$ and 
$\tau = {h r^d \over 1+ CM}$, then
$m_{k+1} \leq \rho m_k + \tau$ for $k \geq 0$, hence 
(by induction)
$$
m_{k} \leq \rho^k m_0 + \tau (1+\rho+\rho^2 + \ldots)
\leq \rho^k m_0 + {\tau \over 1-\rho}
= \rho^k m_0 + h r^d 
\leqno (4.71)
$$
for $0 \leq k \leq N+1$. If $N_0$ (and hence also
$N = N_0/16$) is chosen large enough, depending on 
$a$ and $M$, we get that
$$\eqalign{
\H^d(E \cap {1 \over 2} Q_0) &\leq m_{N+1} 
\leq \rho^{N+1} m_0 + h r^d
= \rho^{N+1} \H^d(E \cap Q_0) + h r^d 
\cr&
\leq a \H^d(E \cap Q_0) + C h l_0^d
}\leqno (4.72)
$$
by (4.69) and (4.71), and because $l_0 \geq l(Q)/4 \geq C^{-1} r$.
So (4.49) holds and (4.42) follows (see below (4.49)). This completes
the proof of Lemma 4.39 and also, as was explained just after the
statement of the lemma, of Proposition 4.1.
\hfill$\square$\ $\square$

\ms\noindent{\bf Remark 4.73.}
The author sees no obvious major obstruction to extending the proof above 
to the case where the rigid assumption is defined in terms of a net of
polyhedra with some uniform size rotundity assumption (instead of dyadic 
cubes of size $2^{-m}$). Still, one would need to construct appropriate 
subnets, or at least adapt the construction of Federer-Fleming projections 
to objects that look like thin neighborhoods of a given polyhedron, but
making sure that we preserve the faces of our initial net. The author 
does not claim that this would be pleasant.

Let us now state a local Ahlfors regularity result under the Lipschitz
assumption, which will follow easily from Propositions 4.1 and 2.8.

\ms\proclaim Proposition 4.74. 
For each choice of $\Lambda \geq 1$ and $M \geq 1$, we can 
find $h > 0$ and $C_M \geq 1$, depending on $\Lambda$, $M$,
and the dimensions $n$ and $d$, so that the following holds. 
Suppose that $E \in GSAQ(U, M, \delta , h)$, 
and that the Lipschitz assumption is satisfied on $U$.
Then if $x\in E^\ast \cap U$ and $0 < r < \Min(\lambda^{-1} r_0,\delta)$  
are such that $B(x,2r) \i U$, we have that 
$$
C_M^{-1} r^d \leq \H^d(E\cap B(x,r)) \leq C_M r^d.
\leqno (4.75) 
$$

\ms
As before, $r_0 = 2^{-m} \in (0,1]$ denotes the side length of the dyadic 
cubes used to define the rigid assumption, and $\Lambda$ and $\lambda > 0$ 
are the constants in the Lipschitz assumption (see Definition 2.7).

It is enough to prove (4.75) for slightly smaller balls, i.e. when
$$
0 < r < \Lambda^{-2} \Min(\lambda^{-1} r_0,\delta)
\ \hbox{ and } \ 
B(x,2\Lambda^2 r) \i U, 
\leqno (4.76) 
$$
because if $B(x,r)$ is as in the original statement, 
a lower bound for $\H^d(E\cap B(x, \Lambda^{-2} r))$
implies a lower bound for $\H^d(E\cap B(x,r))$,
and for the upper bound we may cover $B(x,r)$ by less than
$C$ balls (centered on $E^\ast$ if we want) that satisfy the 
stronger condition (4.76).

So let $B(x,r)$ satisfy the stronger condition.
Set $F = \psi(\lambda E)$; by Proposition 2.8,
$F \in GSAQ(B(0,1), \Lambda^{2d}M, \Lambda^{-1} \lambda \delta, 
\Lambda^{2d} h)$, and Proposition 4.1 applies to that set;
we shall just get a larger constant $C_{\Lambda^{2d} M}$ in (4.2).
Set $y= \psi(\lambda x)$. Of course, $y\in F^\ast$ because 
$x\in E^\ast$, and
$$
\dist(y,\d B(0,1)) \geq \lambda \Lambda^{-1} \dist(x,\d U)
\geq 2\lambda \Lambda r,
\leqno (4.77) 
$$
by (4.76) and the bilipschitz property of $\psi$.
Now $\psi(\lambda B(x,r)) \i B(y,\lambda \Lambda r)$;
let us check that we may apply Proposition 4.1 to 
$B = B(y,\lambda \Lambda r)$.
The fact that $B(y,2\lambda \Lambda r) \i B(0,1)$
follows from (4.77), and 
$$
\lambda \Lambda r 
\leq \lambda \Lambda \Lambda^{-2} \Min(\lambda^{-1} r_0,\delta)
\leq \Min(r_0, \Lambda^{-1}\lambda \delta)
\leqno (4.78) 
$$
by (4.76), so we may apply Proposition 4.1 to $B$ (or a smaller ball
centered at $y\in F^\ast$). 
We get that
$$\eqalign{
\H^d(E\cap B(x,r)) 
&\leq \lambda^{-d} \Lambda^{d} \H^d(\psi(\lambda(E\cap B(x,r))))
\leq \lambda^{-d} \Lambda^{d} \H^d(F\cap B(y,\lambda \Lambda r))
\cr&
\leq \lambda^{-d} \Lambda^{d} C_{\Lambda^{2d} M} (\lambda \Lambda r)^d
= \Lambda^{2d} C_{\Lambda^{2d} M} \, r^d
}\leqno (4.79) 
$$
by (4.2). This is the desired upper bound. For the lower bound,
we observe that $\psi(\lambda B(x,r))$ contains 
$B(y,\lambda \Lambda^{-1} r)$, so
$$\eqalign{
\H^d(E\cap B(x,r)) 
&\geq \lambda^{-d} \Lambda^{-d} \H^d(\psi(\lambda(E\cap B(x,r))))
\cr&
\geq \lambda^{-d} \Lambda^{-d} \H^d(F\cap B(y,\lambda \Lambda^{-1}r))
\cr&
\leq \lambda^{-d} \Lambda^{-d} C_{\Lambda^{2d} M}^{-1} (\lambda \Lambda^{-1} r)^d
= \Lambda^{-2d} C_{\Lambda^{2d} M}^{-1} r^d
}\leqno (4.80) 
$$
by the lower bound in (4.2), applied to the smaller ball $B(y,\lambda \Lambda^{-1} r)$.
This completes the proof of Proposition~4.74.
\qed

\bigskip
\noindent {\bf 5. Lipschitz mappings with big images, projections, and rectifiability.}
\medskip

In this section we extend two propositions from [DS4] 
and prove that quasiminimal sets are rectifiable.

The first proposition, Proposition 5.1 in [DS4] and here, 
concerns the existence of Lipschitz functions defined on a quasiminimal set, 
with values in $\R^d$ and with big images. The second one (Proposition 5.7
below) will concern the quasiminimality of a Lipschitz graph 
over a quasiminimal set. Both will be used to prove uniform
rectifiability estimates in the next section.
But we shall only be able to do this last under additional assumptions, 
so it makes sense to prove the plain rectifiability of quasiminimal 
sets here, because we can prove it in full generality
(and the proof is much easier too). See Theorem 5.16 
below. In the standard case without 
boundaries, the rectifiability was known from Almgren [A2], 
and the proof below is probably quite similar.

We start the section with the existence of Lipschitz functions with big 
images.

\ms\proclaim Proposition 5.1. For each choice of $M \geq 1$, we can 
find $h > 0$ and $C_M \geq 1$, depending on $M$ and the dimensions 
$n$ and $d$, so that the following holds. 
Suppose that $E \in GSAQ(B_0, M, \delta , h)$, where we set 
$B_0 = B(0,1)$, and that the rigid assumption is satisfied.
Let $r_0 = 2^{-m} \in (0,1)$ denote the side length of the dyadic cubes 
used to define the rigid assumption.
Then for $x\in E^\ast \cap B_0$ and $0 < r < \Min(r_0,\delta)$  such that
$B(x,2r) \i B_0$, we can find a $C_M$-Lipschitz mapping $F: E \cap 
B(x,r) \to \R^d$ such that 
$$
\H^d(F(E\cap B(x,r))) \geq C_M^{-1} r^d.
\leqno (5.2)
$$

\ms
By the bilipschitz invariance provided by Proposition 2.8, this result
implies the corresponding one when the Lipschitz assumption holds;
the argument, which goes as in Proposition 4.74, is left to the reader.
Also, we could immediately reduce to the case when $E$ is coral, 
by Proposition 3.27 and because any Lipschitz mapping 
$F : E^\ast \cap B(x,r) \to \R^d$ with a big image could easily be 
extended to $E$. But this would not help much anyway.

The proof is a minor variation on what we did for the lower bound in 
Proposition 4.1. Let $N_0$ be a large power of $2$, to be chosen soon,
and let $l_0$ denote the largest power of $2$ such that $l_0 \leq r/2n$.
Let $Q'$ be a dyadic cube of length $l_0$, and choose a translation 
$Q_0$ of $Q'$ by an integer multiple of $N_0^{-1} l_0$, so that
if $x_0$ denotes the center of $Q_0$, then the size of each coordinate of
$|x-x_0|$ is at most $N_0^{-1} l_0$. Thus $Q_0 \i B(x,r)$, and 
(if $h$ is small enough, depending on $M$ and $n$) 
$$
\H^d(E \cap {1 \over 3} \, Q_0) \geq C^{-1} r^d,
\leqno (5.3)
$$
by Proposition 4.1 and where $C$ depends only on $M$ and $n$.
Next choose an integer $N \in [N_0/2,N_0]$ such that
$$
\H^d(E \cap {N \over N_0} \, Q_0 \sm {N-6 \over N_0}\, Q_0) 
\leq {12 \over N_0}\, \H^d(E \cap Q_0) \leq {C r^d \over N_0} \,;
\leqno (5.4)
$$
the last inequality follows from (4.2), and such an $N$ exists 
because we have $N_0/2$ choices of $N$ and no point of $E$ lies in more than six 
annular regions ${N \over N_0} \, Q_0 \sm {N-6 \over N_0}\, Q_0$.
We can apply a variant of Lemma 4.3 to the cube $Q = {N \over N_0} \, Q_0$, 
with its natural decomposition into $N^n$ subcubes $R \in {\cal R}(Q)$,
which are dyadic of side length $N_0^{-1} l_0$ and lie in the usual 
grid (recall that $N_0$ and $l_0$ are powers of $2$, and that 
$l_0 \leq r/2n \leq r_0/2n$). This gives a first mapping 
$\phi : \R^n \to \R^n$, which preserves 
the faces of all dimensions as in (4.8).

We need a variant because this time we want to say that, in addition 
to the properties already mentioned in Lemma 4.3, $\phi$ is
$C$-Lipschitz, with a constant $C$ that depends on $M$, $n$, and $d$,
but not on $E$ or the $B(x,r)$. This can be arranged, because Proposition 4.1
says that $E \cap Q_0$ is semi-regular, which allows us to apply 
Lemma 3.31 in [DS4]. 
The proof is the same as in [DS4], and as before  
we observe that $N$ does not need to be a power of $2$. 
The point of the argument comes when
we need to choose points in various faces of the skeleton
to perform Federer-Fleming projections centered at these points. 
The semi-regularity of $E \cap Q_0$ says that it, and its images by the 
previous Lipschitz mappings that were already constructed, is sufficiently 
far from dense in any face of dimension $\geq d+1$, so that we can find a 
new center that is far from it; then the Federer-Fleming projection is Lipschitz, 
with good bounds, and we can iterate as long as the faces are at least
$(d+1)$-dimensional. Of course our bounds get worse and worse with 
each iteration (when the codimension is large), but this is all right.

If $\phi(E \cap Q)$ contains a full $d$-dimensional face $T$ of side length
$N_{0}^{-1} l_0$, then we can take for $F$ an extension of $\pi \circ \phi$,
where $\pi : \R^n \to \R^d$ is the composition of the orthogonal 
projection onto the vector space $V$ parallel to $T$, and a linear 
isometry from $V$ to $\R^d$. In this case, (5.2) holds just because
$$
\H^d(F(E\cap B(x,r))) \geq \H^d(\pi\circ\phi(E\cap Q)) 
\geq \H^d(T) = N_0^{-d} l_0^d
\geq C^{-1} r^d.
\leqno (5.5)
$$
So we may assume that $\phi(E \cap Q)$ contains no full $d$-dimensional face 
of a cube $R \in {\cal R}(Q)$, and Proposition 5.1 will follow as 
soon as we derive a contradiction. We then proceed as we did near 
(4.56): we compose $\phi$ with an additional Federer-Fleming projection, 
which is obtained by selecting a center $c_T \in T \sm \phi(E \cap Q)$
in each $d$-face $T$, and projecting on $\partial T$  from there. 
This gives a new mapping, which we shall still call $\phi$, and which 
satisfies the conclusions of Lemma 4.3, plus the fact that 
$\phi(E) \cap [Q \sm A(Q)]$ is contained in a $(d-1)$-dimensional
skeleton, where $A(Q)$ is again the exterior layer of $Q$.
The same computations as in (4.57)-(4.62) yield the analogue of 
(4.62). Here ${N-6 \over N_0}\, Q_0 = {N-6 \over N}\, Q$
plays the role of $Q_{k+1} = (1-{6 \over N_k}) \, Q_k$
(see (4.51)). Thus we get that, if $A_1(Q)$ denotes the union of 
the two exterior layers of $Q$,
$$\eqalign{
\H^d(E \cap {N-6 \over N_0}\, Q_0) 
&\leq \H^d(W_1)
\leq CM \H^d(E \cap A_1(Q)) + h r^d
\cr&\leq CM \H^d(E \cap {N \over N_0} \, Q_0 \sm {N-6 \over N_0}\, Q_0) 
+h r^d
\leq {CM r^d \over N_0} + h r^d
}\leqno (5.6)
$$
as before (i.e., because points of $E \cap {N-6 \over N_0}\, Q_0 \sm W_1$ lie
in some $(d-1)$-dimensional skeleton, by quasiminimality (the analogues 
of (4.58) and (4.61)), and by simple geometry), and then by (5.4).

But if $N_0$ is large enough and $h$ is small enough (both depending 
on $M$, $n$, and $d$), this contradicts (5.3). 
Proposition 5.1 follows.
\qed

\ms
It would be nicer if the mapping $F$ provided by Proposition 5.1
were the orthonormal projection onto a $d$-space. 
The second proposition from [DS4] that we generalize here 
gives a trick that sometimes allows us to pretend that this
is the case. 

\ms\proclaim Proposition 5.7.
Let $U \i \R^n$  be an open set,
$E \in GSAQ(U, M, \delta , h)$ a quasiminimal set, and 
$F : U \to \R^m$ a Lipschitz function. Then 
$\wh E = \big\{ (x,F(x)) \, ; \, x\in E \big\} 
\in GSAQ(U nb\times \R^m, C M, \delta , C h)$
for some constant $C$ that depends on $d$ and the Lipschitz constant for 
$F$, and where on $U \times \R^m$, GSAQ is defined with respect 
to the boundaries $\wh L_j = L_j \times \R^m$.

\ms
This is a minor generalization of Proposition 6.1 in [DS4]. 
Explicit values for $C$ could easily be derived from the proof below, 
but we shall not bother to do so.

The main point of the proof is the following. We are given a competitor
for $\wh E$ in a ball $\wh B$ and want to construct a competitor for $E$, 
to which we apply the definition of $GSAQ(U, M, \delta , h)$ to get estimates.
More specifically, the competitor is $\wh \varphi_1(\wh E)$ for some 
one-parameter family of mappings $\wh \varphi_t : \wh E \to \R^{n+m}$, 
with the properties (1.4)-(1.8). We want to define mappings $\varphi_t$,
and we take
$$
\varphi_t(x) = \pi \circ \wh \varphi_t(x,F(x))
\ \hbox{ for } x\in E,
\leqno (5.8)
$$
where $\pi : \R^{n+m} \to \R^n$ denotes the natural projection onto
$\R^n$. 

Let us check that the $\varphi_t$ satisfy (1.4)-(1.8).
First, (1.4) and (1.8) are trivial. For (1.5) and (1.6), we 
take $B= \pi(\wh B)$. If $x \in E \sm B$, then 
$(x,F(x)) \in \wh E \sm \wh B$, so $\wh \varphi_t(x,F(x)) = (x,F(x))$
and $\varphi_t(x) = \pi(x,F(x)) = x$ by (5.8); similarly,
when $t= 0$, $\wh \varphi_t(x,F(x)) = (x,F(x))$ and 
$\varphi_t(x) = \pi(x,F(x)) = x$; thus (1.5) holds.

Next, if $x\in B$, then either $(x,F(x)) \in \wh E \sm \wh B$, 
and then $\varphi_t(x) = \pi \circ \wh \varphi_t(x,F(x))
= \pi(x,F(x)) = x \in B$ by (5.8) and (1.5), or else
$(x,F(x)) \in \wh E \cap \wh B$ and 
$\wh \varphi_t(x,F(x)) \in \wh B$ by (1.6), so $\varphi_t(x) \in B$, 
as needed for (1.6). Finally, if $x \in E \cap L_j \cap B$,
then $(x,F(x)) \in \wh E \cap \wh L_j$, so 
$\wh \varphi_t(x,F(x)) \in \wh L_j$ by (1.7) or (1.5), and 
$\varphi_t(x) = \pi \circ \wh \varphi_t(x,F(x)) \in L_j$, as needed 
for (1.7). So $\varphi_t(E)$ is a sliding competitor for $E$ in $B$.

Now we check that
$$
W_t \i \pi(\wh W_t)
\ \hbox{ and } 
\varphi_t(W_t) \i \pi(\wh\varphi_t(\wh W_t))
\leqno (5.9)
$$
for $0 \leq t \leq 1$, where 
$W_t = \big\{ y \in E \cap B \, ; \varphi_t(y) \neq y \big\}$ 
is as in (2.1), and 
$\wh W_t = \big\{ z \in \wh E \cap \wh B \, ; \wh\varphi_t(z) \neq z \big\}$
is its analogue for $\wh \varphi_t$. If $x\in W_t$, then
$\pi \circ \wh\varphi_t(x,F(x)) \neq x$, so 
in particular $\wh\varphi_t(x,F(x)) \neq (x,F(x))$; then
$(x,F(x)) \in \wh W_t$
and $x \in \pi(\wh W_t)$. Moreover, 
$\varphi_t(x) = \pi \circ \wh\varphi_t(x,F(x)) \in 
\pi(\wh\varphi_t(\wh W_t))$. So (5.9) holds, and the union of the sets
$W_t \cup \varphi_t(W_t)$ is contained in the projection of the
union of the $\wh W_t \cup \wh\varphi_t(\wh W_t)$. If this last union
is relatively compact in $U \times \R^m$ (as in the assumption (2.4)),
then the first union is relatively compact in $U$, which allows 
us to apply Definition 2.3 and get (2.5). That is,
$$
\H^d(W_1) \leq M \H^d(\varphi_1(W_1)) + h r^d.
\leqno (5.10)
$$
Set $H_1 = \big\{ (x,F(x)) \, ; \, x\in W_1 \big\}$. This is a subset
of $\wh W_1$, because $\pi \circ \wh \varphi_1(x,F(x)) = \varphi_1(x) 
\neq x$ for $x\in W_1$. And
$$
\H^d(H_1) \leq C \H^d(W_1) \leq CM \H^d(\varphi_1(W_1)) + Ch r^d
\leq CM \H^d(\wh\varphi_1(\wh W_1)) + Ch r^d
\leqno (5.11)
$$
because $F$ is Lipschitz, by (5.9), and because $\pi$ is $1$-Lipschitz. 
Now consider $H_2 = \wh W_1 \sm H_1$. First observe that
$$
\H^d(H_2) \leq C \H^d(\pi(H_2))
\leqno (5.12)
$$
because $H_2$ is contained in $\wh E$, which
lies on the graph of the Lipschitz function $F$.
In addition, 
$$
\pi \circ \wh \varphi_1 = \pi
\ \hbox{ on } H_2
\leqno (5.13)
$$
because, if $(x,F(x)) \in H_2$, then $x\in E \sm W_1$, so
$x = \varphi_1(x) = \pi \circ \wh \varphi_1(x,F(x))$
by (5.8). Thus
$$
\H^d(H_2) \leq C \H^d(\pi(H_2)) = C \H^d(\pi\circ \wh \varphi_1(H_2))
\leq C \H^d(\wh \varphi_1(H_2)) \leq C \H^d(\wh \varphi_1(\wh W_1))
\leqno (5.14)
$$
by (5.12) and because $H_2 \i \wh W_1$. Finally, 
$$
\H^d(\wh W_1) \leq \H^d(H_1) + \H^d(H_2)
\leq C(M+1) \H^d(\wh\varphi_1(\wh W_1)) + Ch r^d
\leqno (5.15)
$$
by (5.11) and (5.14), which is (2.5) for $\wh \varphi_1(\wh E)$.
Thus $\wh E$ is a quasiminimal set, and this proves Proposition 5.7.
\qed

\ms
We end this section with the fact that quasiminimal sets are rectifiable.
Recall that this means that such a set $E$ is contained in a countable union of
$d$-dimensional Lipschitz graphs (or $C^1$ surfaces, or Lipschitz images of
$\R^d$, if you prefer), plus a set of vanishing $\H^d$-measure.
Thus $E$ is rectifiable if and only if its core $E^\ast$ is rectifiable.

\ms\proclaim Theorem 5.16. For each choice of $M \geq 1$, we can 
find $h > 0$, depending on $M$, the dimensions 
$n$ and $d$, and the bilipschitz constant $\Lambda$ of $\psi$
in the definition 2.7 of the Lipschitz assumption, such that
if the Lipschitz assumption is satisfied in $U$ and 
$E \in GSAQ(U, M, \delta , h)$, then $E$ is rectifiable.

\ms
Let $E \in GSAQ(U, M, \delta , h)$ be as in the statement, and
let $\psi$ be the bilipschitz function in Definition 2.7. 
Since we know that any Lipschitz image of a rectifiable set 
is rectifiable (see 15.3 and 15.21 in [Ma],   
to which we shall refer for anything that concerns 
rectifiability), it will be enough to show that $\psi(\lambda E)$
is rectifiable (where $\lambda > 0$ is also as in Definition 2.7).
But by Proposition 2.8, $\psi(\lambda E) \in GSAQ(B(0,1), \Lambda^{2d}M, 
\Lambda^{-1}\lambda\delta , \Lambda^{2d}h)$, with the rigid 
assumption. So it is enough to prove Theorem 5.16 when $U=B(0,1)$ 
and the rigid assumption holds.

Recall that $E$, just like any other set of locally finite
$\H^d$-measure, can be written as the disjoint union
$E = E_r \cup E_s$ of a rectifiable part $E_r$ and a totally
non rectifiable (or singular) part $E_s$. We just need to show that
$\H^d(E_s) = 0$, because the union of a rectifiable set and an 
$\H^d$-null set is rectifiable too. So we shall assume
that $\H^d(E_s) >0$ and derive a contradiction. 

Since $\H^d(E)$ is locally finite and $E_r$ does not meet $E_s$,
a standard density result 
(see for instance [Ma], Theorem 6.2 (2) on page 89)  
says that
$$
\lim_{\rho \to 0} \rho^{-d} \H^d(E_r \cap B(x,\rho)) = 0
\leqno (5.17)
$$
for $\H^d$-almost every $x\in E_s$. In addition,
$\H^d(E_s \sm E^\ast) \leq \H^d(E \sm E^\ast) = 0$,
so we can pick $x\in E_s \cap E^\ast$ such that (5.17) holds.

We shall proceed as in the proof of Proposition 5.1.
Let $N_0$ be a large power of $2$ (a constant to be chosen soon),
and let $l_0$ be a very small power of $2$ (it will just need to
be small enough, depending on $N_0$ too). 
Let $Q'$ be a dyadic cube of length $l_0$, 
and choose a translation $Q_0$ of $Q'$ by an integer multiple of 
$N_0^{-1} l_0$, so that if $x_0$ denotes the center of $Q_0$, then
the size of each coordinate of $|x-x_0|$ is smaller than 
$N_0^{-1} l_0$ and so $E \cap {1 \over 3} \, Q_0$ contains
$B(x,{l_0 \over 10})$.
As in (5.3), if $h$ is small enough, depending on $M$ and $n$,
$$
\H^d(E \cap {1 \over 3} \, Q_0) 
\geq\H^d(E \cap B(x,{l_0 \over 10})) 
\geq C^{-1} l_0^d,
\leqno (5.18)
$$
by Proposition 4.1 and where $C$ depends only on $M$ and $n$.
In fact, we do not even need Proposition 4.1 here; we could have
chosen a point $x\in E$ where the upper density 
$\limsup_{r \to 0} r^{-d} \H^d(E \cap B(x,r))$ is
larger than a geometric constant, and then taken $l_0$
as small as we want and such that (5.18) holds.
We choose the integer $N \in [N_0/2,N_0]$ such that
$$
\H^d\big(E \cap {N \over N_0} \, Q_0 \sm {N-6 \over N_0}\, Q_0\big) 
\leq {12 \over N_0}\, \H^d(E \cap Q_0) \leq {C l_0^d \over N_0} \, ,
\leqno (5.19)
$$
as in (5.4) and with the same simple proof by Chebyshev.
(And again we could also have obtained the last inequality because 
$\limsup_{r \to 0} r^{-d} \H^d(E \cap B(x,r)) \leq C$ almost 
everywhere on $E$.)

Then we want to apply the proof of Lemma 4.3 to the cube 
$Q = {N \over N_0} \, Q_0$, with its natural decomposition 
into $N^n$ subcubes $R \in {\cal R}(Q)$,
which are all dyadic of side length $N_0^{-1} l_0$ 
and lie in the usual grid if $l_0$ is small enough. 
[This time, we shall not need $\phi$ to be $C$-Lipschitz (as for 
Proposition 5.1), so we do not need the variant that uses the
local Ahlfors-regularity of $E^\ast$.]

So we want to mimic the construction of $\phi$ in Lemma 4.3
(or rather Proposition 3.1 in [DS4]),  
but with a few changes because we want to project away the 
unrectifiable part.
Our mapping $\phi$ will be obtained as the last element
of a sequence $\phi_{n+1}, \phi_n, \phi_{n-1}, \cdots, \phi_d$,
obtained recursively by composing with mappings $\psi_l$.
That is, we shall start from $\phi_{n+1}(z)=z$ and set 
$$
\phi_l = \psi_l \circ \phi_{l+1} \ \hbox{ for }
d \leq l \leq n. 
\leqno (5.20)
$$

For $R \in {\cal R}(Q)$ (the set of dyadic cubes $R \i Q$ of side length
$N_0^{-1} l_0$) and $0 \leq l \leq n$, denote by
${\cal S}_l(R)$ the union of all the $l$-dimensional
faces of $R$. Also set 
${\cal S}_l = \bigcup_{R \in {\cal R}(Q)} {\cal S}_l(R)$.
We intend to choose the $\psi_l$ in such a way that
$$
\phi_{l}(E \cap Q) \i {\cal S}_{l-1} \cup \d Q
\leqno (5.21)
$$
for $d \leq l \leq n$.

Let us say how we do this. We start with
$n \geq l > d$; the case of $l=d$ is a little special, 
and will be treated at the end. Assume that the $\psi_k$,
$k > l$, were already defined, and satisfy (5.21).
This is of course true when $l=n$.
First we decide that
$$
\psi_l(z)=z \ \hbox{ for } z \in \d Q \cup [\R^n \sm Q],
\leqno (5.22)
$$
because we want (4.4) to hold. 
Since $\phi_{l+1}(E \cap Q) \i {\cal S}_{l} \cup \d Q$
(by (5.21) for $l+1$), the main thing to do now, if we want 
a definition of $\phi_l = \psi_l \circ \phi_{l+1}$ on $E$, 
is to define $\psi_l$  on ${\cal S}_{l} \cap Q$.
We shall define $\psi_l$ simultaneously on all the
$l$-dimensional faces $F$, in such a way that
$\psi_l(z) = z$ on $\d F$. Then there will be no
problem about coherence.

So let $F$ be a $l$-dimensional face of a cube $R \in {\cal R}(Q)$. 
If $F \i \d Q$, we keep $\psi_{l}(z)=z$ on $F$ (because of (5.22)). 
Otherwise, we select an origin $x_F \in F \sm \phi_{l+1}(E \cap Q)$, 
near its center. Other constrains will show up soon, but
for the moment let us record that $\H^l$-almost every point
$x_F$ is like this, because $\phi_{l+1}(E \cap Q)$ is at most 
$d$-dimensional (all our mappings are Lipschitz), and $l > d$.
Notice that $x_F \in F \sm \phi_{l+1}(E)$
too, because $\phi_{l+1}(E \sm Q)$ lies far from the center
of $F$ (by iterations of (5.22)).

Pick a small ball $B_F$ centered at
$x_F$ and such that
$$
\dist(B_F,\phi_{l+1}(E)) > 0.
\leqno (5.23)
$$
The small size of $B_F$ will not matter, it will just make
the Lipschitz constant for $\phi$ enormous, but we don't care. 
We decide that
$$\eqalign{
&\hbox{for $z\in F \sm B_F$, $\psi_l(z)$ is the radial projection}
\hbox{ of $z$ on $\d F$, centered at $x_F$.}
}\leqno (5.24)
$$
This last just means that $\psi_l(z) \in \d F$ and
$z$ lies on the segment $[x_F,\psi_l(z)]$. With this choice,
observe that if $z\in E \cap Q$, then by (5.21) for $l+1$,
$\phi_{l+1}(z)$ either lies on $\d Q$ (and then 
$\phi_l(z) = \phi_{l+1}(z) \in \d Q$ by (5.22) and (5.20)),
or else lies on some $F \sm B_F$ (by (5.23)), so
$\phi_l(z) = \psi_l(\phi_{l+1}(z)) \in \d F$ by (5.24).
Thus (5.21) holds for $l$.

Now we extend $\psi_l$ in a Lipschitz way, first to $F$
(so that $\psi_l(F) \i F$), then (after we are finished with all 
the faces $F$) to faces $G$ of higher dimensions (so that 
$\psi_l(G) \i G$ for every face $G$) and eventually 
the whole $Q$. Thus $\psi_l(Q) \i Q$. 
One checks (see the proof in [DS4]) 
that these definitions give rise to Lipschitz mappings $\psi_l$, 
which satisfy (4.4)-(4.6), and also (4.8). For the remaining
estimate (4.7) on the $\H^d(\phi(E\cap  R))$, we need to be 
more careful about the choice of centers $x_F$, and this is also 
where we shall not proceed exactly as in [DS4]. 

We want to treat the rectifiable and singular parts of $E$ 
separately. We still intend to use Lemma 3.22 in [DS4], 
which goes as follows.

\ms\proclaim Lemma 5.25.
Let $F$ is an $l$-dimensional face of cube, with $l > d$, 
and $A \i F$ a closed set such that $\H^d(A) < +\infty$.
For $\xi \in {1 \over 2} F$, denote by
$\theta_{\xi,F} : F \sm \{\xi \} \to \d F$
the radial projection on $\d F$ centered at $\xi$. Then 
$$
\H^l(F)^{-1} \int_{\xi \in {1 \over 2} F \sm A} \H^d(\theta_{\xi,F}(A))
d\H^l(\xi) \leq C \H^d(A).
\leqno (5.26)
$$

\ms
As the proof will show, the lemma stays true if $A$
is merely Borel-measurable, but its closure has a finite
$\H^d$-measure. The main point of the proof is that
for a given $\xi\in {1 \over 2} F \sm \overline{A}$,
$$
\H^d(\theta_{\xi,F}(A)) 
\leq C \int_A |x-\xi|^{-d} \diam(F)^d d\H^d(x),
\leqno (5.27)
$$
which follows from computing the local Lipschitz 
constant of $\theta_{\xi,F}$ near $x$. See (3.24)
and (3.20) in [DS4]. 
We integrate this over $\xi \in {1 \over 2} F \sm A$,
use Fubini, and get that
$$\leqalignno{
\int_{\xi \in {1 \over 2} F \sm A} \H^d(\theta_{\xi,F}(A))
&d\H^l(\xi) 
\leq C \diam(F)^d \int_{x\in A}  
\int_{\xi \in {1 \over 2} F \sm \overline {A}} 
|x-\xi|^{-d} d\H^l(\xi) d\H^d(x)
\cr&
\leq C \diam(F)^d \int_{x\in A}  
\Big\{ \int_{\xi \in F \cap B(x,2\diam(F))} 
|x-\xi|^{-d} d\H^l(\xi) \Big\} d\H^d(x)
\cr&
\leq C \diam(F)^d \int_{x\in A} \diam(F)^{l-d} d\H^d(x)
\leq C \diam(F)^l \H^d(A);
& (5.28)
}
$$
(5.26) and the lemma follow.
\qed

\ms
In [DS4] and for Proposition 4.1, 
Lemma 5.25 is used with $A = F \cap \phi_{l+1}(E)$ to choose 
$x_F$ so that, with the definitions (5.20) and (5.24),
$$
\H^d(\psi_l(F \cap \phi_{l+1}(E)))
\leq C \H^d(F \cap \phi_{l+1}(E)).
\leqno (5.29)
$$
Such a choice is possible, by Fubini. Let us record here the 
fact that, by the proof of Lemma 5.25, we can even get that
$$
\int_{F \cap \phi_{l+1}(E)} |x-x_F|^{-d} \diam(F)^d d\H^d(x)
\leq C \H^d(F \cap \phi_{l+1}(E)),
\leqno (5.30)
$$
which is stronger than (5.29) because of (5.27).

With this choice of $x_F$ for each $F$, we can 
sum over $F$, compose our mappings, and get that
$$
\H^d(\phi_{d+1}(E\cap R)) \leq C \H^d(E\cap R)
\ \hbox{ for } R \in {\cal R}(Q)
\leqno (5.31)
$$
as in (4.7). The proof is the same as in [DS4] 
and for Proposition 4.1.

As we said earlier, here we want to take advantage of
the fact that $x_F$ is chosen by a Fubini argument to
apply Lemma 5.25 with $A = F \cap \phi_{l+1}(E_r)$
and get, in addition to (5.29), that
$$
\H^d(\psi_l(F \cap \phi_{l+1}(E_r)))
\leq C \H^d(F \cap \phi_{l+1}(E_r)),
\leqno (5.32)
$$
and finally obtain, after composing, that
$$
\H^d(\phi_{d+1}(E_r\cap R)) \leq C \H^d(E_r\cap R)
\ \hbox{ for } R \in {\cal R}(Q).
\leqno (5.33)
$$
We sum this and get that
$$
\H^d(\phi_{d+1}(E_r \cap Q))
\leq \sum_{R \in {\cal R}(Q)} \H^d(\phi_{d+1}(E_r \cap R))
\leq C \H^d(E_r \cap Q)
\leqno (5.34)
$$
because the cubes $R \in {\cal R}(Q)$ have bounded overlap.

Let $\varepsilon > 0$ be very small, to be chosen
soon. Because $x$ was chosen so that (5.17) holds,
we deduce from (5.34) that
$$
\H^d(\phi_{d+1}(E_r \cap Q)) \leq \varepsilon l_0^d
\leqno (5.35)
$$
if $l_0$ was chosen small enough.

For the unrectifiable part $E_s$ of $E$, we use 
the following fact, which is proved in Lemma~4.3.3 
on page 111 of [Fv1]  
or Lemma 6 on page 26 of [Fv3]. 
If $F$ is a face of dimension $l > d$, and
if $A \i F$ is such that $\H^d(A) < +\infty$
and $A$ is totally non rectifiable (of dimension $d$),
then for almost every choice of $x_F$, 
$\psi_l(A) \i \d F$ is also totally non rectifiable.
Of course we choose the various $x_F$ so that this 
happens (we had some latitude left to do this); 
then when we compose the $\psi_l$ we get that
$$
\phi_{d+1}(E_s \cap Q) \hbox{ is totally non rectifiable.}
\leqno (5.36)
$$
All this information is valid also on the cubes $R$ that meet
$\d Q$; we concentrated on what happens on faces $F$ that are
not contained in $\d Q$, but on $\d Q$ we simply need to know that
all our mappings are the identity.

Next, if $F$ is any $d$-dimensional face of a cube $R\in {\cal R}(Q)$,
then 
$$
\H^d(F \cap \phi_{d+1}(E \cap Q))
= \H^d(F \cap \phi_{d+1}(E_r \cap Q)) \leq \varepsilon l_0^d,
\leqno (5.37)
$$
because the totally non rectifiable set $\phi_{d+1}(E_s \cap Q)$ 
can only meet the rectifiable set $F$ on a $\H^d$-null set, and by (5.35). 

Thus $\phi_{d+1}(E \cap Q)$ never fills a $d$-face $F$
(if $\varepsilon$ is small enough), and this allows us 
to choose, for each $d$-dimensional face $F$ in $Q$
which is not contained in $\d Q$, a point $x_F$ near
the center of $F$ that does not lie on $\phi_{d+1}(E \cap Q)$.
We then choose $B_F$ and define $\psi_d$ as we did above, near 
(5.23). This gives a last mapping $\phi_d = \psi_d \circ \phi_{d+1}$,
which still satisfies (5.21).

We shall need in a later section to know that if $\tau > 0$ is small enough
(depending also on our choice of mappings $\phi_l$ and their 
bad Lipschitz constants), and if $H \i Q$ is a compact set such that
$$
\dist(z,E) \leq \tau \ \hbox{ for } z\in H,
\leqno (5.38)
$$
then
$$
\phi_{l}(H \cap Q) \i {\cal S}_{l-1} \cup \d Q
\ \hbox{ for } n+1 \geq l \geq d.
\leqno (5.39)
$$
We naturally prove this by descending induction.
Obviously this is true for $l=n+1$, because
$\phi_{n+1}(z) = z$ and ${\cal S}_{n}=Q$.
Let $l \geq d$ be given, and suppose
that (5.39) holds for $l+1$.
Let $z\in H \cap Q$ be given; by induction assumption,
$\phi_{l+1}(z) \in {\cal S}_{l} \cup \d Q$.
If $\phi_{l+1}(z)\in \d Q$, (5.22) and (5.20) say that 
$\phi_{l}(z) = \phi_{l+1}(z) \in \d Q$, so we may assume
that $\phi_{l+1}(z)$ lies in some $l$-face $F$ that is
not contained in $\d Q$. We know, since $\phi_{l+1}$
is Lipschitz (and $\tau$ is as small as we want) that
$\phi_{l+1}(z)$ is arbitrarily close to $\phi_{l+1}(E)$,
so (5.23) says that $\phi_{l+1}(z)$ lies out of $B_F$,
and hence $\psi_{l}(z)$ is given by (5.24). Then 
$\phi_{l}(z) = \phi_{l+1}(z) \in \d F$, as before, and
(5.39) holds for $l$ too. This proves (5.39).
\ms

Return to $E$, and set $\phi^\ast = \phi_d$
(we write $\phi^\ast$ instead of $\phi$ to make sure
that the $\phi^\ast_t$ below will not be confused with 
the $\phi_l$ above). Thus
$$
\phi^\ast(E \cap Q) = \phi_d(E \cap Q) \i {\cal S}_{d-1} \cup \d Q
\leqno (5.40)
$$
by (5.21) with $l=d$.
Now we use the quasiminimality of $E$ to get a contradiction.
It is easy to construct a one-parameter family
$\{ \phi_t^\ast \}$, that satisfies (1.4)-(1.8),
and for which $\phi_1^\ast = \phi^\ast$; the verification
is the same as for Proposition 4.1, for instance near (4.16).
Set
$$
W_t = \big\{ y \in \R^n \, ; \varphi_t^\ast(y) \neq y \big\}
\leqno (5.41)
$$
for $0 < t \leq 1$ and
$$
\widehat W = \bigcup_{0< t \leq 1} W_t \cup \varphi_t^\ast(W_t);
\leqno (5.42)
$$
these are well defined here because the mappings $\varphi_t^\ast$
are defined everywhere. 
We can easily arrange the interpolation between the identity
and $\phi^\ast$ so that $\phi_t^\ast(z) = z$
for $z\in \R^n \sm Q$ and $\phi_t^\ast(Q) \i Q$, and so we 
get that 
$$
\widehat W \i Q \i Q_0 \i \overline B(x,2\sqrt n l_0)
\leqno (5.43)
$$
(see the definition of $Q_0$ and $Q$ near (5.18) and (5.19)).
If $l_0$ is chosen small enough, 
$\overline B(x,2\sqrt n l_0)$ is arbitrarily small
and contained in $U$, so we can apply Definition 2.3.
We get that
$$
\H^d(E \cap W_1) \leq M \H^d(\phi^\ast(E \cap W_1)) + h r^d.
\leqno (5.44)
$$
Denote by ${\cal R}_{ext}$ the collection of
small cubes $R \in {\cal R}(Q)$ that touch $\d Q$
(that is, the ${\cal R}_{ext}$ is the outer layer
of cubes in ${\cal R}(Q)$). Then set
$$
Q' = \bigcup_{ R \in {\cal R}(Q) \sm {\cal R}_{ext}} R.
\leqno (5.45)
$$
Recall from the definition of $Q$ below (5.19) that 
$$
Q = {N \over N_0} \, Q_0 \ \hbox{ for some integer }
N \in [N_0/2,N_0]. 
\leqno (5.46)
$$
Also, the side length of our
cubes $R \in {\cal R}(Q)$ is $N_0^{-1} l_0$, so
$$
Q' = {N-2 \over N} \, Q = {N-2 \over N_0} \, Q_0 
\supset {1 \over 3} \, Q_0\, 
\leqno (5.47)
$$
because $N_0$ is very large. Let us check that
$$
E\cap Q' \sm W_1 \i {\cal S}_{d-1}.
\leqno (5.48)
$$
Let $z\in E\cap Q' \sm W_1$ be given.
Then $z = \phi^\ast(z) \in {\cal S}_{d-1}$ 
by (5.41), (5.40), and because $z\notin \d Q$; (5.48) follows.
Then
$$
\H^d(E \cap W_1)  \geq \H^d(E \cap Q' \cap W_1)
= \H^d(E \cap Q') 
\geq H^d(E \cap {1 \over 3} Q_0) \geq C^{-1} l_0^{d}
\leqno (5.49)
$$
by (5.48), (5.47), and (5.18).
On the other hand, by (5.40) (and the first half of (4.4)),
$$
\H^d(\phi^\ast(E \cap W_1))
= \H^d(\d Q \cap \phi^\ast(E \cap W_1)).
\leqno (5.50)
$$
Let us check that
$$
\d Q \cap \phi^\ast(E \cap W_1)
\i \bigcup_{R \in {\cal R}_{ext}}
\phi^\ast(E \cap R),
\leqno (5.51)
$$
where ${\cal R}_{ext}$ still denotes the outer rim of small 
cubes $R \in {\cal R}(Q)$ that touch $\d Q$.
Let $w\in \d Q \cap \phi^\ast(E \cap W_1)$ be
given, and let $z\in E \cap W_1$ be such that 
$\phi^\ast(z)=w$. Observe that $z$ lies out of $Q'$, 
because (4.6) says that $\phi^\ast(Q') \i Q'$.
So $z \in E \cap R$ for some $R \in {\cal R}_{ext}$,
and (5.51) follows. Next we verify that for
$R \in {\cal R}_{ext}$,
$$
\H^d(\phi^\ast(E \cap R) \sm \phi_{d+1}(E \cap R)) = 0.
\leqno (5.52)
$$
Let $w\in \phi^\ast(E \cap R)$ be given, 
and choose $z \in E \cap R$ such that
$w = \phi^\ast(z)$. Recall that
$w = \phi^\ast(z) = \phi_d(z) = \psi_d(\phi_{d+1}(z))$
by definition of $\phi^\ast$ and (5.20). By (5.21),
$\phi_{d+1}(z) \in {\cal S}_d \cup \d Q$.
If $\phi_{d+1}(z) \in \d Q$, then $\psi_d$
does not move it (by (5.22)), and so 
$w = \psi_d(\phi_{d+1}(z)) = \phi_{d+1}(z)$, 
which is fine for (5.52). Otherwise, $\phi_{d+1}(z)$
lies on some $d$-dimensional face $F$ that is not contained
in $\d Q$, and by construction
its image by $\psi_d$ (that is, $w$) lies on $\d F$, 
which is $(d-1)$-dimensional. So (5.52) holds.
Altogether,
$$\eqalign{
\H^d(\phi^\ast(E \cap W_1))
&\leq \sum_{R \in {\cal R}_{ext}} \H^d(\phi^\ast(E \cap R))
\leq \sum_{R \in {\cal R}_{ext}} \H^d(\phi_{d+1}(E \cap R))
\cr&
\leq C \sum_{R \in {\cal R}_{ext}} \H^d(E\cap R))
\leq C \H^d(E \cap Q \sm {\rm int}(Q'))
}\leqno (5.53)
$$
by (5.50), (5.51), (5.52), (5.31), and the fact that the
cubes $R$ have bounded overlap. Since
$$
Q \sm {\rm int}(Q') = {N \over N_0} \, Q_0
\sm {\rm int}\Big({N-2 \over N_0} \, Q_0 \Big)
\leqno (5.54)
$$
by (5.46) and (5.47), it follows from (5.53), (5.54), and
(5.19) that
$$\eqalign{
\H^d(\phi^\ast(E \cap W_1)) 
&\leq C \H^d(E \cap Q \sm {\rm int}(Q'))
\cr&
\leq C \H^d\Big( E \cap {N \over N_0} \, Q_0
\sm {\rm int}\Big({N-2 \over N_0} \, Q_0 \Big)\Big)
\leq C {l_0^d \over N_0}.
}\leqno (5.55)
$$
If $N_0$ is large enough and $h$ is small enough (depending on $M$ in 
particular), we get a contradiction with (5.44) or (5.49);
thus we could not find our initial point of density $x \in E_s$,
and the rectifiability of $E$ follows.
This completes our proof of Theorem 5.16.
\qed

\bigskip
\centerline{PART III : UNIFORM RECTIFIABILITY OF QUASIMINIMAL SETS}
\ms
This part is largely independent from the next ones,
which is probably a good thing because we shall only be able to 
complete the desired program in some specific cases, depending on the
dimensions of the faces of the $L_j$.

The main goal is to prove that sliding quasiminimal sets are locally
uniformly rectifiable, with big pieces of Lipschitz graphs.

When we wrote the long paper [DS4], 
and even for later results, the author thought that the local uniform
rectifiability of $E$ was an unavoidable main step for many things, 
including the stability of our classes of minimizers under limits
(as in Part IV below). As we shall see later, this is not the case,
and the proof of rectifiability is enough for many purposes.

This is fortunate, because we shall not be able to prove the 
local uniform rectifiability of $E$ in all the interesting cases, and
also because even when it works, the proof is more difficult than
usual. 

We nonetheless include a part on uniform rectifiability here because 
the author cannot deny his past, and it is a nice regularity property. 
It is probably almost the best general result that we can hope to prove 
for quasiminimal sets. That is, because quasiminimality is bilipschitz 
invariant (or directly), Lipschitz graphs are quasiminimal, and
uniformly rectifiable sets are not so different (in terms of regularity)
from Lipschitz graphs. Even for almost minimal or minimal sets, it is 
not so clear how to get better general regularity results (i.e., that 
would hold without assuming some a priori flatness, for instance), 
even though in this case we expect better regularity.

We continue with the same general writing style as in Part II, i.e.,
giving a rapid general description of [DS4], 
except at places where 
modifications are needed (and then we need to be more precise).

\bigskip
\noindent {\bf 6. Local uniform rectifiability in some cases.}
\medskip

So we want to prove that sliding quasiminimal sets are locally
uniformly rectifiable, with big pieces of Lipschitz graphs,
and we shall only be able to do this under an
additional assumptions on the dimensions of the faces of the $L_j$.
The main result of this section and the next two is the following
theorem, and its generalization (Theorem 9.81) under the Lipschitz 
assumption.

\ms\proclaim Theorem 6.1. 
For each choice of $M \geq 1$, we can 
find $h > 0$, $A \geq 0$, and $\theta > 0$, depending on $M$ and 
the dimensions $n$ and $d$, so that the following holds. 
Suppose that $E \in GSAQ(B_0, M, \delta , h)$, where we set 
$B_0 = B(0,1)$, and that the rigid assumption is satisfied.
Let $r_0 = 2^{-m} \leq 1$ 
denote the side length of the dyadic cubes used to define the 
rigid assumption.
Let $x\in E^\ast \cap B_0$ and $0 < r < \Min(r_0,\delta)$
be such that $B(x,2r) \i B_0$. 
Assume in addition that 
$$\eqalign{
&\hbox{if $j \in [0,j_{ max}]$ is such that some face of dimension
(strictly) more than $d$}
\cr&\hskip 1cm
\hbox{of $L_j$ meets $B(x,r)$, then $E^\ast \cap B(x,r) \i L_j$.}
}\leqno (6.2) 
$$
Then we can find a $d$-dimensional $A$-Lipschitz
graph $\Gamma \i \R^n$ such that
$$
\H^d(E\cap \Gamma \cap B(x,r)) \geq \theta r^d. 
\leqno (6.3)
$$

\ms
By $d$-dimensional $A$-Lipschitz graph, we mean a set $\Gamma$ which 
is the image, under an isometry of $\R^n$, of the graph of some
Lipschitz function from $\R^d$ to $\R^{n-d}$ whose Lipschitz norm is
at most $A$. 
Notice that we do not have so much of a restriction on dimensions 
when $d=2$ and $n=3$, which will probably be our main interest 
in the future (but even so we do not allow $L_1$ to be a half space
in which $E$ is not contained).
Also, Theorem 6.1 does not necessarily apply when $d=2$, $n=4$, and some $L_j$
are $3$-dimensional. 

The author does not know whether this additional restriction on the 
dimensions is really needed. 

The restrictions in Theorem 6.1 do not seem too bad, for instance 
because they allow boundary constraints given by sets $L_j$ of 
dimensions at most $d$, and the typical setting for a Plateau problem 
is like this. But in terms of proof, Theorem 6.1 is rather 
disappointing because it does not contain much more information than 
what is readily available from the interior uniform rectifiability
(away from the $L_j$). For instance, if all the $L_j$ are
at most $(d-1)$-dimensional, the local uniform rectifiability of 
$E^\ast$ near the $L_j$ follows from the inside uniform rectifiability
and the local Ahlfors-regularity given by Proposition 4.1
(there is just not enough room near the $L_j$ for a bad behavior).
We will be able to obtain more cases (for instance, increase the 
dimension of the $L_j$ by one) by various general tricks, but the 
center of the proof is still the result from [DS4]. 
That is, a
simpler special case will be obtained in Proposition 6.41,
with a minor modification of the argument of [DS4], 
and then
the extension of this result that we do in Sections 7 and 8 
will mostly use general manipulations of uniform rectifiability 
and Carleson measures, and for instance we shall only construct 
competitors once, in Lemma 7.38 or its generalization Lemma 9.14. 
It would be nice to have a different,
simpler proof of the uniform rectifiability of $E^\ast$ away from the 
$L_j$, but for moment we only know one (very complicated) proof.

Here is our plan for the rest of this part.
We shall start this section with a rapid description of the
proof of local uniform rectifiability given in [DS4]. 
We shall then say (largely for the record) why it does not seem to 
go through with our sliding conditions. 
In the last subsection, we prove a weaker variant of Theorem 6.1,
Proposition 6.41, which is what we can almost directly 
obtain from the proof of [DS4]. 
In Section 7, we shall prove the conclusion of Proposition 6.41
(the existence of a big piece of bilipschitz image of a subset
of $\R^d$) under the weaker assumptions of Theorem~6.1; see
Proposition 7.85. Theorem 6.1 itself will only be proved in 
full in Section~8, with a small additional argument on the existence
of big projections. Finally, Theorem 6.1 will be generalized to
the case of Lipschitz assumption in Section 9.
See Theorem 9.81.

\ms
\noindent{\bf 6.a. How we want to proceed, following [DS4]} 
\ms

So we are given a quasiminimal set $E$ and a ball $B(x,r)$, 
as in the statement of Theorem 6.1.
Because of Proposition 3.3 (which says that $E^\ast$ is quasiminimal 
too), we can assume that $E$ is coral (i.e., $E = E^\ast$); otherwise
just prove and apply the theorem for $E^\ast$.

We first use Proposition 5.1 to find a $C_M$-Lipschitz mapping 
$F_r: E \cap B(x,r) \to \R^d$ such that
$\H^d(F_r(E\cap B(x,r))) \geq C_M^{-1} r^d$, as in (5.2). By Whitney's 
extension theorem, we can extend $F_r$ into a $CC_M$-Lipschitz mapping
defined on $\R^n$.

Next we apply Proposition 5.7, which says that $\wh E$, the
graph of $F_r$ over $E$, is a quasiminimal set in $\R^{n+d}$.
We shall denote by $\pi_x : \R^{n+d} \to \R^n$ and $\pi : \R^{n+d} \to \R^d$ 
the two natural projections, and consider the smaller set
$\wh E_0 = \wh E \cap \pi_x^{-1}(B(x,r))$. Then
$\pi(\wh E_0) = F_r(E\cap B(x,r))$ and 
$$
\H^d(\pi(\wh E_0)) = \H^d(F_r(E\cap B(x,r))) \geq C_M^{-1} r^d.
\leqno (6.4)
$$
Next we want to use a stopping time argument from [D1] 
to find a large piece of $\wh E_0$ where $\pi$ is bilipschitz.
More precisely, we want to find a closed set $\wh\Gamma_0 \i \wh E_0$
such that
$$
\H^d(\wh\Gamma_0) \geq \theta' r^d
\ \hbox{ and } \ 
|y-z| \leq A'|\pi(y)-\pi(z)|
\hbox{ for } y,z \in \wh \Gamma_0,
\leqno (6.5)
$$
where $\theta' > 0$ and $A'$ are constants that depend only on
$n$, $d$, and $M$.

If we do so, this will not directly give a big piece of Lipschitz
graph in $E\cap B(x,r)$, as required in the statement of Theorem 6.1,
but the following weaker conclusion: there is a closed set
$G_0 \i E \cap B(x,r)$ and a mapping $\phi : G_0 \to \R^d$
such that 
$$
\H^d(G_0) \geq \theta r^d
\ \hbox{ and } \ 
C'_M |y-z| \leq |\phi(y)-\phi(z)| \leq C'_M |y-z|
\hbox{ for } y,z \in G_0,
\leqno (6.6)
$$
where $\theta$ and $C'_M$ depend only on $n$, $d$, and $M$.
In other words, instead of a big piece of Lipschitz graph in
$E\cap B(x,r)$, we only find a big piece of bilipschitz image
of a subset of $\R^d$. 

The verification (from (6.5)) is easy: we just try
$G_0 = \pi_x(\wh\Gamma_0)$; then (6.6) follows from
(6.5) because $\pi_x : \wh E \to E$ is bilipschitz.

In the terminology of [DS1] or [DS3], 
(6.6) (for all $x$ and $r$) says that locally, $E$ has big pieces of 
bilipschitz images of $\R^d$ (BPBI), which amounts to saying that 
$E$ is locally uniformly rectifiable, while in the statement
of Theorem 6.1 we claim that if also contains big pieces of 
Lipschitz graphs (BPLG) locally.

Now we can go from BPBI to BPLG by a general argument on uniformly 
rectifiable sets, for which we just need to check that $E$ also has 
``big projections". This will be discussed soon, but anyway the
most important part of Theorem 6.1 is the local uniform rectifiability
provided by (6.6).

Return to $\wh E_0 = \wh E \cap \pi_x^{-1}(B(x,r))$, 
our quest of $\wh \Gamma_0 \i \wh E_0$
such that (6.5) holds, and the stopping time argument
from [D1]. 
We would like to use the proof described in Sections 8 and 9 
of [DS4], which we try to explain now. 

A first ingredient of the proof is the construction of what we call
cubical patchworks on $\wh E$, which are the analogue on $\wh E$ of the
standard dyadic cubes on $\R^n$, and which will be very
useful because we want to run stopping time arguments on $\wh E$.
This construction is done in Section 7 of [DS4], 
and goes through in the present setting because it only uses the 
local Ahlfors-regularity of $E$ near $x$. This last holds as soon 
as $h$ is small enough (depending on $n$ an $M$), by Proposition 4.1
and because we assumed that $E = E^\ast$. 
Naturally, we shall always assume that this ($h$ small enough) is the case.
Let us say what the cubical patchwork is in the situation
of Theorem 6.1. We get a set $F$ and collections $\Sigma_j$,
$j \geq 0$, of so-called dyadic cubes, with the following properties.
First,
$$
\wh E \cap B(\wh x,r/10) \i F \i \wh E \cap B(\wh x,r) \i \wh E_0,
\leqno (6.7)
$$
where we call $\wh x = (x,F_r(x))$ the natural center for $\wh E_0$,
and $F$ also is locally Ahlfors-regular, in the sense that
$$
C^{-1} t^d \leq \H^d(F\cap B(y,t)) \leq C t^d
\ \hbox{ for $y\in F$ and } 0 < t < r.
\leqno (6.8)
$$
For each $j \geq 0$, $\Sigma_j$ is a collection
of measurable subsets $Q$ of $F$, which we shall abusively call cubes, 
such that $F$ is the disjoint union of the cubes $Q$, $Q\in \Sigma_j$.
The cubes have some low regularity properties, and particular
they have a center $c_Q$ such that 
$$
F \cap B(c_Q,C^{-1} 2^{-j}r) \i Q \i F \cap B(c_Q,C 2^{-j}r)
\ \hbox{ for } Q \in \Sigma_j.
\leqno (6.9)
$$
They also have small boundaries (see (7.4) and (7.10) in [DS4]), 
but we shall not use this here. Finally, the 
$\Sigma_j$ have the same structure as for the usual dyadic cubes:
if $i \leq j$, $Q \in \Sigma_i$, and $R \in \Sigma_j$,
then $R \i Q$ or else $R \cap Q = \emptyset$.

The main property that we need to prove if we want to get (6.5)
is a little complicated, and involves a (given) large constant $C_1$,
a (given) small constant $\gamma$, and constants $C_2$
(very large) and $\eta$ (very small), to be chosen (depending
on $C_1$, $\gamma$, $M$, and $n$). For $y\in F$ and $j \geq 0$,
set
$$
T_j(y) = \bigcup_{Q \in \Sigma_j \, ; \, Q \cap B(y,C_2 2^{-j} r) 
\neq \emptyset} \, Q.
\leqno (6.10)
$$
Thus $T_j(y)$ is a little bit like $F \cap B(y,C_2 2^{-j} r)$,
but we prefer to cut neatly along dyadic cubes. The stopping time
argument from [D1]  
that we want to use likes the situations (depending on
$y\in F$ and $j \geq 0$) when our projection $\pi : \R^{n+d} \to \R^d$
has a local surjectivity property, namely when
$$
\pi(T_j(y)) \supset \R^{d} \cap B(\pi(y),C_1 2^{-j} r).
\leqno (6.11)
$$
It also likes it when there is a cube $R\i T_j(y)$ that does significantly
better than average in terms of projections, i.e., when
$$
\hbox{ there exists $R \in \Sigma_j$ such that $R \i T_j(y)$ and }\ \ 
{\H^d(\pi(R)) \over \H^d(R)} \geq (1+2\eta) \, {\H^d(\pi(T_j(y)))
\over \H^d(T_j(y))}.
\leqno (6.12)
$$
The property that makes things work in [DS4]  
is the following.

\ms\proclaim Definition 6.13.
We say that the main lemma holds if for each choice of
$C_1$ and $\gamma > 0$, we can find $C_2$ and $\eta > 0$,
depending on $C_1$, $\gamma$, $M$, $n$, and $d$ (which includes
a dependence on the local Ahlfors-regularity and cubical patchwork 
constants) such that, whenever $y \in F$ and $j \geq 0$
are such that
$$
B(y,2C_22^{-j} r) \i B(\wh x,r/10)
\leqno (6.14)
$$
and
$$
{\H^d(\pi(T_j(y)))
\over \H^d(T_j(y))} \geq \gamma,
\leqno (6.15)
$$
then we have (6.11) or (6.12).

\ms
This property is proved in [DS4], as Main Lemma 8.7. 
The fact that it allows us to apply a theorem from [D1] 
and get a graph $\wh\Gamma_0$ as in (6.5) is proved in 
Section 8 of [DS4],  
and the proof goes through without major modification in 
the present context. [Again, it only uses the local Ahlfors regularity 
properties of $E$, and no construction of competitors.]

\ms
Thus we want to know whether the main lemma holds in the context
of sliding minimizers, and we study the proof given in 
Section 9 of [DS4]. 

We assume that we can find $y \in F$ such that (6.14) and (6.15)
hold, but not (6.11) or (6.12), and we want to reach a contradiction
(for a correct choice of $C_2$ and $\eta$). That is, we want to
construct an appropriate deformation of $\wh E$ (which is a 
quasiminimal set by Proposition 5.7), for which most of the measure
near $x$ disappears. A first step in the verification, which is done
in Section 9-2 of [DS4], 
consists in obtaining the following description of $F$ 
(or equivalently $\wh E$) near $y$.

As in (9.62) of [DS4], 
we apply a dilation to all our sets so that 
$$
2C_1 2^{-j} r = 1;
\leqno (6.16)
$$
this normalization will allow us to work with (standard!) dyadic cubes 
of unit side length in the $d$-plane $P = \R^d$. We can also assume 
that $y=0$. Still denote by $\pi$ the orthogonal projection on
$P = \R^d$, and by $\pi_x$ the orthogonal projection on $V = \R^{n}$
(in [DS4] it is called $h$, 
but we want to avoid a conflict of notation here).

We shall restrict our attention to the box $P_0 \times V_0$,
where $P_0 = [-N,N]^d \i P$ for some large integer $N$,
and $V_0 = V \cap \overline B(0,\rho_0)$ for some $\rho_0 \in [N,CN]$.
Here $N$ will be chosen very large, depending on 
$C_1$, $\gamma$, $M$, $n$, and $d$, and $C$ is so large
(depending on the same constants), that a Chebyshev argument
allows us to choose $\rho_0\in [N,CN]$ so that 
$$
\dist(P_0 \times (V \cap \d B(0,\rho_0)), F) \geq N.
\leqno (6.17)
$$
[See (9.77) in [DS4], and we won't   
need to modify this part of the argument.] 
Later on, we shall choose $C_2$ (depending also on $N$), so large that 
$$
P_0 \times V_0 \i B(y,C_2 2^{-j-1} r) \i B(\wh x,r/11);
\leqno (6.18)
$$ 
the first part is easy to arrange (because $y=0$ and by the normalization 
(6.16)), and the second inclusion comes from (6.14). Now set
$$
F_0 = F \cap (P_0 \times V_0) = \wh E \cap (P_0 \times V_0) \i \wh 
E_0,
\leqno (6.19)
$$
where the last part comes from (6.7) and (6.18). Notice that
$$
F_0 \i T_j(y)
\leqno (6.20)
$$
by (6.18) and (6.10), so it will be easy to use our assumption 
that (6.11) and (6.12) fail.

Denote by $A_i$, $i\in I_0$, the collection of cubes in $P$,
contained in $P_0$, that are obtained from the unit cube $[0,1]^d$ 
by an integer translation in $\Z^d$. That is, we cut $P_0$ into 
$(2N)^{d}$ dyadic unit cubes (and the point of the normalization 
above is that we can use unit cubes).

We finally come to our description of $F_0$. First set
$$
I_1 = \big\{ i\in I_0 \, ; \, 
\hbox{ there is an $x_i \in int(A_i)$ such that }
F_0 \cap \pi^{-1}(x_i) = \emptyset \big\},
\leqno (6.21)
$$
and let us check that
$$
I_1 \hbox{ is not empty.}
\leqno (6.22)
$$
Recall that (6.11) fails, so we can find 
$w\in P \cap B(y,C_1 2^{-j}r) \sm \pi(T_j(y))$.
By (6.16) and because $y=0$, $w\in B(0,1/2)$; 
by (6.20), $w \in P \sm \pi(F_0)$.
By (6.19), $F_0$ is compact, so a whole neighborhood of
$w$ in $P$ lies in $P \sm \pi(F_0)$. This neighborhood 
contains an interior point of some $A_i$, $i\in I_0$,
and by definition this $i$ lies in $I_1$. This proves (6.22).

Next, for each $i\in I_0$ there is a finite set 
$\Xi(i) \in F_0 \cap \pi^{-1}(A_i)$, with at most $C$ elements, 
and such that
$$
\dist(z,\Xi(i)) \leq 1
\ \hbox{ for every } z\in F_0 \cap \pi^{-1}(A_i).
\leqno (6.23)
$$
This is checked in [DS4], 
and the same proof applies here; see the verification 
of ($\ast$9.3) (understand, (9.3) in [DS4]) 
below ($\ast$9.89), which relies a lot on
Lemma $\ast$9.83. Let us just say here why this is not 
surprising. 

Let ${\cal R} = \big\{ R \in \Sigma_j \, ; \, R \i T_j(y) \}$
denote the set of cubes that compose $T_j(y)$; these cubes are
disjoint by definition of $\Sigma_j$. Set
$a_0 = {\H^d(\pi(T_j(y))) \over \H^d(T_j(y))}$; 
thus $a_0 \geq \gamma$ by (6.15).
Also set $a(R) = {\H^d(\pi(R)) \over \H^d(R)}$
for $R \in {\cal R}$, and 
$$
a_1 = \H^d(T_j(y))^{-1} \sum_{R \in {\cal R}} \H^d(\pi(R))
= \H^d(T_j(y))^{-1} \sum_{R \in {\cal R}} a(R) \H^d(R);
\leqno (6.24)
$$
thus $a_1$ is a weighted average of the $a(R)$
(with the weights $\H^d(R)$, $R \in {\cal R}$), and at the same time
$$
a_1 = \H^d(T_j(y))^{-1} \sum_{R \in {\cal R}} \H^d(\pi(R))
\geq \H^d(T_j(y))^{-1} \H^d(\pi(T_j(y))) = a_0
\leqno (6.25)
$$
because $\pi(T_j(y))$ is the union of the $\pi(R)$.
Now (6.12) fails, so 
$a(R) \leq (1+2\eta) a_0 \leq (1+2\eta) a_1$ for $R \in {\cal R}$.
If $\eta$ is small, this forces all the $a(R)$ to be very close
to their average $a_1$, and also $a_1-a_0$ to be very small.
A first consequence is that $a(R) \geq \gamma/2$, and hence
$\H^d(\pi(R)) \geq \gamma \H^d(R)/2 \geq C^{-1}$, 
for $R \in {\cal R}$. But also, the various $\pi(R)$, 
$R \in {\cal R}$, are almost disjoint,
because when $R_1$, $R_2 \in {\cal R}$ are different,
$$
\H^d(T_j(y)) (a_1-a_0)
= \H^d(\pi(T_j(y))) - \sum_{R \in {\cal R}} H^d(\pi(R))
\geq \H^d(\pi(R_1)\cap \pi(R_2))
\leqno (6.26)
$$
by the proof of (6.25). Then there are at most $C$
such cubes $R \in {\cal R}$ such that $\pi(R)$ falls near 
a given $A_i$, and the existence of $\Xi(i)$ as in 
(6.23) follows reasonably easily. 
The bound $C$ on the cardinal of the $\Xi(i)$ depends on 
$\gamma$ (and the local Ahlfors-regularity constants, as usual).

Finally set $I_2 = I_0 \sm I_1$. We also prove that for each $i\in I_2$,
there is a point $z_i \in F_0$ such that
$\pi(z_i) \in int(A_i)$ and
$$
|z-z_i| \leq 1 \ \hbox{ for all $z\in F_0$ 
such that } \pi(z) = \pi(z_i).
\leqno (6.27)
$$
This would be obvious if the $\pi(R)$ were disjoint, because 
$\diam(R) \leq 2C2^{-j}r < 1$ by (6.9), by the normalization (6.16),
and because we can assume that $C_1 > 2C$, where $C$,
the constant from (6.9), depends only on the Ahlfors regularity 
constant. Here the $\pi(R)$ are merely nearly disjoint, 
so we have to work a little more, i.e., use Chebyshev.
The verification is done in [DS4],  
below ($\ast$9.89), proof of ($\ast$9.5).
This completes our description of $F_0$.

\ms
The next stage in our proof is a deformation lemma
(Proposition 9.6 in [DS4]) 
that sends most of $F_0$ to a $(d-1)$-dimensional set.
The proposition concerns a more arbitrary closed set in $\R^{n+d}$, 
but we apply it to $F_0$, and the main assumptions are (6.22), (6.23)
and (6.27), that we just obtained. It yields the following.

\ms\proclaim Lemma 6.28.
There exists a family $\{\phi_t\}$, $0 \leq t \leq 1$,
of Lipschitz mappings of $\R^{n+d}$, with the following properties:
$$
(t,z) \to \phi_t(z) \hbox{ is Lipschitz (from
$[0,1] \times \R^{n+d} \to \R^{n+d}$);}
\leqno (6.29)
$$
$$
\phi_t(z) = z
\hbox{ for $t=0$ and when $\dist(z,P_0 \times V_0) \geq d+3$;}
\leqno (6.30)
$$
$$\eqalign{
&\hbox{if $\phi_t(z) \neq z$ for some $z\in \R^{n+d}$ and
$t\in [0,1]$, then}
\cr& \hskip-0.3cm
\hbox{$\dist(z,F_0) \leq C$ and $\dist(\phi_s(z),F_0) \leq C$ 
for $0 \leq s \leq 1$}
}\leqno (6.31)
$$
(so our deformation may move some points a lot, but only close to
$F_0$ and somehow along $F_0$);
$$
\phi_t(F_0) \i P_0 \times V_0
\hbox{ for all $t\in [0,1]$;}
\leqno (6.32)
$$
$$
H^{d-1}(\phi_1(F_0)) < +\infty
\leqno (6.33)
$$
(the main point: we essentially make $F_0$ disappear);
$$
\phi_t \hbox{ is $C$-Lipschitz on } \R^{n+d} \sm (P_0 \times V_0).
\leqno (6.34)
$$

\ms
As usual, the constant $C$ in this statement depends on 
$n$, $M$, $C_1$, and $\gamma$ through the constants that arise in the 
description of $F_0$ above. The last property (6.34)
is useful because we do not want to lose what
we win by (6.32) by increasing too much the Hausdorff measure
of $F = \wh E$ near the boundary of $P_0 \times V_0$.

Lemma 6.28 still holds in our case, with the same proof, but
differences will occur in the way it is applied.

As the reader may have guessed, the mappings $\phi_t$
are used in [DS4] to produce a deformation of $\wh E$  
which contradicts its quasiminimality. The reader should not worry 
about the way the Hausdorff measure estimates go,
because it will be the same as in [DS4], 
but let us just say a few words to explain some of our choices.
For instance, (6.34) goes with some control on the size of the set
$$
H = \big\{ z \in \wh E \, ; \, 
0 < \dist(z,P_0 \times V_0) \leq d+3 \big\},
\leqno (6.35)
$$
where the $\phi_t$ may differ from the identity, but we cannot use 
(6.33). And we required in (6.17) that
$\dist(P_0 \times (V \cap \d B(0,\rho_0)), F) \geq N$
to get an easier control on $H$. Indeed,
let $z\in H$ be given. Then $z\in F$ by (6.7), the first 
part of (6.18), and (6.14). 
Write $z = (\pi(z),\pi_x(z)) \in P \times V$, then
$\dist(\pi(z),P_0) \leq d+3$ by (6.35), so 
$\dist(\pi_x(z), V \cap \d B(0,\rho_0))) \geq N-d-3 > N/2$
by (6.17), and since $V_0 = V \cap \overline B(0,\rho_0)$ and
$\dist(\pi_x(z),V_0) \leq d+3$ by (6.35),
we get that $\pi_x(z) \in V \cap B(0,\rho_0-N/2)$.
Altogether, $H$ is contained in the simpler set
$$
H' = \big\{ z \in \wh E \, ; \, \pi_x(z) \in B(0,\rho_0-N/2)
\, \hbox{ and } \, 0 < \dist(\pi(z), P_0) \leq d+3 \big\},
\leqno (6.36)
$$
which is easier to control because we can use the fact that $P$ is
$d$-dimensional to cover $H'$ by a $C N^{d-1}$ balls of radius $1$,
using something like (6.23). Near the end of the argument, $N$ is 
chosen so large that the contribution of $H$ to the $\H^d$-measure
of the image $\phi_1(\wh E)$, which is less than $C N^{d-1}$, 
is negligible compared to the mass of $\wh E$ that we save by (6.33), 
which larger than $C^{-1} N^d$.
So it is important that $C$, in particular in (6.34), does not depend on $N$.
But again the reader should not worry too much, the computations 
are done in [DS4].  

Up to now, we did not need to worry, because all our constructions 
relied on the general properties of $E$ (local Ahlfors-regularity,
existence of Lipschitz mappings with a big image), and not on the 
definition of a quasiminimal set. This changes now.

In [DS4], the $\phi_t$ define a competitor for $\wh E$, 
which is significantly better than $\wh E$ (by (6.33) and (6.34),
and observations as above) and this leads to the desired contradiction. 
The computations are done at the end of Section $\ast$9.2, near ($\ast$9.90) 
and below.
In the case of generalized quasiminimal sets, we need to add 
a small term $C h r_0^d$ (coming from the $hr^d$ in (2.5))
to the right-hand side of ($\ast$9.93). But if $h$ is chosen small enough,
this does not upset the end of the proof: the additional term is 
small compared to $M \H^d(\wh E \cap W)$ in the right-hand side,
because this last is larger than $C^{-1} r^d$ by ($\ast$9.105) and the 
line before.

In the present situation, the difficulty will come from the fact that
the $\phi_t$ may fail to define a competitor for $\wh E$, because
we don't know whether they respect the boundaries $L_j$ as in (1.7).
Note however that the other constraints (1.4)-(1.6) and (1.8) are satisfied,
with the same verification as in [DS4]. 

\ms
\noindent{\bf 6.b. Some bad news}
\ms
Let us try to continue, and see whether the $\phi_t$ defined
in [DS4] satisfy the last constraint (1.7), or  
can be modified so that (1.7) holds. The main point of this short subsection 
is to explain why the author thinks there is a serious difficulty for the
brutal extension of the proof of [DS4]. 

There is a first obvious reason. Suppose $n=3$, $\Omega = L_0$
is the half space $\{x_1 \geq 0 \big\}$, $L_1 = \d \Omega$
is the vertical plane $\{x_1 = 0 \big\}$, and $E = P \cap \Omega$
for some $2$-plane $P$ perpendicular to $\{x_1 = 0 \big\}$. For 
instance, $E = \big\{ x_1 \geq 0 \hbox{ and } x_3 = 0 \big\}$.
In this case, we don't need the trick of replacing $E$ with
$\wh E$, because the projection $\pi$ over $P$ already has a big
image, but if we did it, we would just replace $\R^3$ with
$\R^5 = \R^3 \times P$, and $E = P \cap \Omega$ with the slightly 
tilted half plane $\big\{ (x_1,x_2,0,x_1,x_2) \, ; \, x_1 \geq 0 
\big\}$, and the discussion would stay the same as below.

Pick a small ball $B(x,r)$ centered at $x=0$, and try to apply the
proof above. Also pick $y=0$; we would like to say that the main lemma
from Definition 6.13 holds, but we can't. Indeed (6.12) never holds, 
because $\H^d(\pi(R))$ is just proportional to $\H^d(R)$, (6.15) holds 
for the same reason (the proportionality constant is not small),
and (6.11) fails because $\pi(E) = P \cap \Omega$ only covers half
the desired ball.

In [DS4], this would never happen, because we would be allowed to 
deform points of the boundary $\big\{ x_1 = x_3 = 0 \big\}$ along $E$
into the domain, thus making a good piece of $E$ disappear and 
contradicting the quasiminimality of $E$. And indeed the 
proof of [DS4] does something like that, which is not allowed here  
because of (1.7) for $L_1$.

Let us say a little more about how these things happen
in the proof of [DS4]; 
the reader may also skip the following discussion and turn to the proof of 
Theorem 6.1 which starts in in the next subsection. 

There are three phases in the construction of the $\phi_t$
in [DS4]. 
In the first one we we move points horizontally 
(i.e., with trajectories parallel to $P$), 
independently in each $\pi^{-1}(A_i)$, so as to project on
$\pi^{-1}(\d A_i)$ whenever this is possible. That is, let
$\phi_{1/3}$ denote the endpoint of this first phase, and let
$z \in \pi^{-1}(A_i)$ be given. When $i\in I_1$, we manage
to obtain that $\pi(\phi_{1/3}(z))$ is the radial projection 
of $\pi(z)$ on $\partial A_i$, centered at the point $x_i$
of (6.21). When $i\in I_2$, we manage to obtain that 
$\pi(\phi_{1/3}(z))$ is the the radial projection 
of $\pi(z)$ on $\partial A_i$, centered at the point $\pi(z_i)$
of (6.27), but only when $z \in \pi^{-1}(\d A_i)$
and when $z$ lies far from $z_i$ (more precisely, 
when $|\pi_x(z)-\pi_x(z_i)| \geq 2$).

We can keep the first phase as it is, because since we only move
points horizontally and the boundaries for $\wh E$ are the sets
$\wh L_j = P \times L_j$, the condition (1.7) is automatically
satisfied. At the end of this first stage, $\phi_{1/3}(F_0)$,
seen from far, looks a little like a piece of graph (over the
union of the $A_i$, $i\in I_2$), plus some uncontrolled junk
above the $\partial A_i$.

For the second phase (corresponding to $\phi_t$, $1/3 \leq t \leq 2/3$), 
we move points vertically, so as to merge the various points 
of $\pi^{-1}(x) \cap \phi_{1/3}(F_0)$ into a single point when this
is possible. The difficulty is to make this in a Lipschitz way 
with respect to $x$. Between $\phi_{1/3}$ and $\phi_{2/3}$, 
we move the points linearly, i.e., we take
$$
\phi_t(z) = (2-3t) \phi_{1/3}(z)
+ (3t-1) \phi_{2/3}(z)
\ \hbox{ for } 1/3 \leq t \leq 2/3,
\leqno (6.37)
$$
with, in coordinates, 
$$
\phi_{2/3}(z) = (\pi(\phi_{1/3}(z)), \varphi(\phi_{1/3}(z))
\in P \times V
\leqno (6.38)
$$
for some $\varphi : \R^{n+d} \to V$ that describes the vertical
motion. Observe that our notation here is slightly different; 
what we denote by $(\pi(z),\varphi(z))$ now was 
called $\phi_2$ in [DS4]); then (6.38) here corresponds 
to ($\ast$9.48) there. 

It is not so important to describe the precise definition of
$\varphi$ and $\phi_{2/3}$ here. Let us just say that this is done with 
partitions of unity, and that the main point is that the resulting set 
$F_2 = \phi_{2/3}(F_0)$ has the following nicer property.

Recall that for $i\in I_1$,
$\phi_{1/3}(F_0) \cap \pi^{-1}(int(A_i)) = \emptyset$,
so 
$$
F_2 \cap \pi^{-1}(int(A_i)) = \emptyset
\hbox{ for } i\in I_1,
\leqno (6.39)
$$
just because our second phase moves points vertically.
[Also see ($\ast$9.43) in [DS4]]. 
When $i\in I_2$, we only know that
$|\pi_x(z)-\pi_x(z_i)| \leq 2$ for all 
$z \in \pi^{-1}(int(A_i)) \cap \phi_{1/3}(F_0)$.
But by our our vertical motion, we make sure that
$$
F_2 \cap \pi^{-1}(int(A_i)) \i \Gamma_i
\hbox{ for } i\in I_2,
\leqno (6.40)
$$
where $\Gamma_i$ is the graph over $int(A_i)$ of some Lipschitz
function. See ($\ast$9.44) in [DS4]. 
So the point of the vertical motion is to merge
the various points of $\pi^{-1}(w)$, $w\in int(A_i)$; the partitions 
of unity help us do this in a nice Lipschitz way.

Our control on the sets $F_2 \cap \pi^{-1}(\d A_i)$ is a little less 
precise, but still ($\ast$9.45) in [DS4] says that 
each of them is contained in a finite union of Lipschitz
graphs over $\d A_i$, so their total $\H^d$-measure is null.

In the present situation, there would be a way to modify the 
construction of $\phi_{2/3}$ and $F_2$, so that we also have
the preservation (1.7) of the boundary pieces $L_j$. In other
words, the serious problem is not here yet. 
The idea is to try to favor choices of points with integer 
coordinates in $V$ in the description $\Xi_i$, but
let us not be more precise, because more serious 
problems will arise in the third phase.

\ms
In the third and last phase of the construction of [DS4], 
points move a lot more. The mappings $\phi_t$, $2/3 \leq t \leq 1$, 
are obtained by composing successive deformations, 
each time occurring on $\pi^{-1}(A_i \cup A_j)$ for some pair of adjacent 
cubes in $P_0$. That is, we set $t_k = 2/3 + 2^{-k}/6$ and construct
recursively $\phi_{t}$, $t\in I_k = [t_k,t_{k+1}]$.
At the start our set $F(k) = \phi_{t_k}(F_0)$ is
composed of a certain number of Lipschitz graphs $\Gamma_i$,
$i\in I(k) \i I_2$ over the corresponding open squares $int(A_i)$,
plus a set $Z(k)$ of finite $H^{d-1}$-measure. Notice that we have such a 
description for $F(0) = F_2$, where $I(k) = I_2$ and the small
set $Z(0)$ lies above the $\d A_i$.

If $I(k) = \emptyset$, we stop, and keep $\phi_t = \phi_{t_k}$ for
$t_k \leq t \leq 1$. Otherwise, we select an $i\in I(k)$
and a $j \in I_0 \sm I(k)$ that are contiguous, i.e., shares
a face $S$ of dimension $d-1$. Such a pair exists, because 
$I_1 \i I_0 \sm I(k)$ is not empty. We first construct our 
deformation on $\Gamma_i$, so that it moves points inside 
$\Gamma_i$ so that the final image lies in 
$\Gamma_i \cap \pi^{-1}(\d A_i)$ (and even in
$\Gamma_i \cap \pi^{-1}(\d A_i \sm int(S))$)
and fixes every point of $\Gamma_i \cap \pi^{-1}(\d A_i \sm int(S))$.
We just use $\pi^{-1}(int(S)) \cap \Gamma_i$ as a base to push
the points in the direction of $\pi^{-1}(\d A_i)$.

Then we extend our Lipschitz deformation into a Lipschitz deformation
of $\R^{n+d}$, which leaves $\pi^{-1}(P \sm A_i \cup A_j)$ alone, and
it is easy to see that $F_{k+1} = \phi_{t_{k+1}}(F_0)$ satisfies the 
induction assumption. At the end of the construction, 
$I(k) = \emptyset$, $H^{d-1}(F_k) < +\infty$ as in (6.33), and we are happy.

Unfortunately, we cannot arrange (1.7) for the mappings that we just 
constructed. The main difficulty is when $\pi^{-1}(A_i)$ contains some
points of some $L_j$, say, for the only initial index $i\in I_1$. 
In our construction, these points get pushed to $\pi^{-1}(\d A_i)$,
and then to other boxes. Along the way, they have to stay close to 
$\wh E$, and this may well be incompatible with (1.7), for instance
if $\wh E$ gets away from $L_j$.

Now we could hope to be lucky, and have a sequence of indices
$i\in I_0$, that can be removed in the corresponding order, and
such that the list of sets $L_j$ that touch $\pi^{-1}(A_i)$
is a nondecreasing function of time, or some similar condition that
seems hard to get in practice. But in view of the counterexample 
(a half plane) given at the beginning of the subsection, this
hope looks very optimistic.

\ms\noindent
{\bf 6.c. What we can say anyway} 
\ms
There is one special case when the proof of [DS4] 
can easily be adapted, and which we record now.

Recall that $B_0$ is the unit open ball and that
$r_0 = 2^{-m}$ is the scale of the dyadic cubes in the description
of the $L_j$. 

\ms \proclaim Proposition 6.41. 
For each choice of $M \geq 1$, we can 
find $h > 0$, $\theta > 0$, and $C_M \geq 1$, 
depending on $M$ and $n$, so that the following holds. 
Suppose that $E \in GSAQ(B_0, M, \delta , h)$ and that the 
rigid assumption is satisfied, and let 
$x\in E^\ast \cap B_0$ and $0 < r < \Min(r_0,\delta)$ be such 
that $B(x,2r) \i B_0$. Also assume that
$$
E \cap B(x,r) \i L_j \ \hbox{ for every $j$ such that
$L_j$ meets $B(x,r)$}.
\leqno (6.42)
$$
Then there is a closed set
$G_0 \i E^\ast \cap B(x,r)$ and a mapping $\phi : G_0 \to \R^d$
such that 
$$
\H^d(G_0) \geq \theta r^d
\ \hbox{ and } \ 
C_M^{-1} |y-z| \leq |\phi(y)-\phi(z)| \leq C_M |y-z|
\hbox{ for } y,z \in G_0.
\leqno (6.43)
$$

\ms
Thus $E^\ast \cap B(x,2r)$ contains a big piece of bilipschitz image
of a subset of $\R^d$. The case when no $L_j$ meets $B(x,r)$
corresponds to the result of [DS4]. 

We shall find it more convenient to assume (6.42) as it is, but
it would be enough to assume that $E^\ast \cap B(x,r) \i L_j$ when
$L_j$ meets $B(x,r)$; the stronger corresponding statement is
simply deduced by applying Proposition 6.41 to the set $E^\ast$,
which also lies in $GSAQ(B_0, M, \delta , h)$ by Proposition 3.3.

We shall see later how to deduce Theorem 6.1
from the proposition, but for the moment let us prove it.
We follow the scheme of [DS4], as explained 
above, and in particular get mappings $\phi_t$ as in 
Lemma 6.28; the only thing that we need to do is modify
them so that the satisfy (1.7) in addition.

We are only interested in the $L_j$ that meet $B(x,r/5)$. Indeed,
for the other ones, the constraint (1.7) is trivially satisfied 
because we shall only consider competitors for $\wh E$ in $B(\wh x,r/10)$. 
With the $\phi_t$ that we have so far, this follows from (6.30) and 
(6.18), and this will stay true after we modify the $\phi_t$ below.

Let $J_0$ denote the set of indices $j$ such that $L_j$ meets 
$B(x,r/5)$. We may assume that $J_0$ is not empty, because
otherwise there is nothing to check since the $\phi_t$ from 
Lemma~6.28 do the job as in [DS4]. 
Set
$$
L = \bigcap_{j\in J_0} L_j \, ;
\leqno (6.44)
$$
observe that 
$$
E \cap B(x,r) \i L,
\leqno (6.45)
$$
by (6.42), so $L$ is not empty. Apply Lemma 3.4 to $L$ 
(after a dilation of factor $r_0^{-1}$, as in (3.26), used to 
return to faces of unit size). This gives a Lipschitz 
retraction 
$$
\pi_{L} : L^\eta \to L, \hbox{ with } 
L^\eta = \big\{ y\in \R^n \, ; \, \dist(y,L) \leq \eta \big\}. 
\leqno (6.46)
$$
In Lemma 3.4 we could take $\eta = 1/3$, but here,
since $L$ is composed of faces of size $r_0$, we take
$\eta = r_0/3$ because we conjugate with a dilation.
Retraction means that $\pi_{L}(y) = y$ on $L$, and Lemma 3.4 also
says that $\pi_{L}$ preserves the faces of size $r_0$ of any
dimension. Finally define the analogue of $\pi_L$ on $\R^{n+d}$ by
$$
\Pi(z) = (p,\pi_L(v)) = (\pi(z),\pi_L(\pi_x(z)))
\ \hbox{ for } z = (p,v) \in P \times V.
\leqno (6.47)
$$

We are ready to define the mappings $\phi_t^\ast$
that will replace the $\phi_t$ from Lemma 6.28.
We want to set
$$
\phi_t^\ast(z) = \Pi(\phi_t(z))
\ \hbox{ for } z\in F,
\leqno (6.48)
$$
so let us check that this makes sense. 
If $\phi_t(z) = z$, then its $V$-coordinate $\pi_x(z)$ lies in 
$E\cap B(x,r)$ (by (6.7)), hence also in $L$ by (6.45), so 
$\pi_L(\pi_x(z)) = \pi_x(z)$ and $\Pi(\phi_t(z))$ is not only 
defined, but equal to $z$. For the record,
$$
\phi_t^\ast(z) = z
\ \hbox{ when $z\in F$ and $\phi_t(z) = z$.}
\leqno (6.49)
$$
If $\phi_t(z) \neq z$, (6.31) says that $\dist(\phi_t(z),F_0) \leq C$.
Now $F_0 \i L$ by (6.7) and (6.45), and we claim that $C$ in (6.31) is 
much smaller than $\eta = r_0/3$. 
Indeed, the normalization (6.16) says that $2C_1 2^{-j}r = 1$, and (6.18) 
implies that $C_2 2^{-j-1} r \leq r/11 \leq r_0/11$ 
(by an assumption of Proposition 6.41). Now $C_2$
is huge, much larger than $C_1$, which proves our claim.
Thus $\phi_t(z) \in L^\eta$, $\Pi(\phi_t(z))$ is defined, and (6.48) 
makes sense.

We select a very small number $\rho > 0$ (which will even depend
on the Lipschitz constant for the $\phi_t$), and keep
$$
\phi_t^\ast(z) = \phi_t(z)
\ \hbox{ when } \dist(z,F) \geq \rho
\leqno (6.50)
$$
and 
$$
\phi_t^\ast(z) = \phi_t(z) = z 
\ \hbox{ when $t=0$, or $\dist(z,P_0 \times V_0) \geq d+3$,
or $\dist(z,F_0) \geq C$},
\leqno (6.51)
$$
where $C$ is as in (6.31). The fact that $\phi_t(z) = z$
comes from (6.30) in the first two cases, and (6.31) in the last one.
Notice that (6.51) and (6.48) are compatible, by (6.49).
In addition, $(t,z) \to \phi_t^\ast(z)$ is Lipschitz (with a very 
large constant that depends on $\rho$) on the set where we defined 
it so far, by (6.29) and because $\Pi$ is Lipschitz. 
Indeed, we need to estimate $|\phi_t^\ast(z)-\phi_s^\ast(z')|$,
and the only case where we do not already know the Lipschitz estimate
is when we use two different definitions, i.e., when $z \in F$
and $\dist(z',F) \geq \rho$, or the other way around.

Next we extend $\phi_t^\ast$ to $\R^{n+d}$ in a Lipschitz way,
using the standard proof with Whitney cubes (here their size is 
at most $\rho$ because $\phi_t^\ast(z)$ was defined when $z\in F$ and 
when $\dist(z,F) \geq \rho$) and partitions of unity. 
Our extension $\phi_t^\ast$ thus 
satisfies (6.29) (by construction) and (6.30) (by (6.51)).

Let us now check (6.31), and so let $z\in \R^{n+d}$ be such that 
$\phi_t^\ast(z) \neq z$ for some $t\in [0,1]$. 
Because of (6.51), we know that $\dist(z,F_0) \leq C$.
If $\dist(z,F) \geq \rho$, (6.50) says that 
$\phi_t^\ast(z) = \phi_t(z)$, and then (6.31) says that
$\dist(z,F_0) \leq C$ and $\dist(\phi_s^\ast(z),F_0) = \dist(\phi_s(z),F_0)
\leq C$ for $0 \leq s \leq 1$. 

So we may assume that $\dist(z,F) < \rho$, and we let $z_0 \in F$ 
be such that $|z-z_0| \leq \rho$. 
By construction, every $\phi_s^\ast(z)$ is a convex
combination of various values of $\phi_s(w)$ or 
$\Pi(\phi_s(w))$, where $w\in B(z_0,10\rho)$.
Since $\phi_s$ and $\Pi \circ \phi_s$ are Lipschitz,
$|\Pi(\phi_s(w)) - \Pi(\phi_s(z_0))|
\leq C |\phi_s(w) - \phi_s(z_0)| \leq C \rho$, with a very large constant
$C$ but that does not depend on $\rho$. 

If in addition $\dist(z_0, F_0) \geq d+3$, then $\phi_s(z_0) = z_0$ 
by (6.30), hence $\Pi(\phi_s(z_0)) = \phi_s^\ast (z_0) =z_0$ by (6.49)
and (6.48). Then $|\phi_s^\ast(z)-z_0| \leq C \rho$, and
$$
\dist(\phi_s^\ast(z),F_0) \leq \dist(z_0,F_0) + C \rho 
\leq \dist(z,F_0)+C\rho \leq C+1
\leqno (6.52)
$$
for $0 \leq s \leq 1$, and if $\rho$ is small enough
(recall that $\dist(z,F_0) \leq C$ because $\phi_t^\ast (z) \neq z$
for some $t$).

Otherwise, if $\dist(z_0, F_0) < d+3$, then 
$\dist(\phi_s(z_0), F_0) \leq C$ for $0 \leq s \leq 1$, by (6.31)
or (if $\phi_s(z_0) = z_0$) simply because $\dist(z_0, F_0) \leq d+3$. 
Since $\Pi$ coincides with the identity on $F_0$, this implies that
$|\Pi(\phi_s(z_0))-\phi_s(z_0)| \leq C' \dist(\phi_s(z_0), F_0)
\leq C''$, where now $C'$ and $C''$ also depend on the Lipschitz constant
for $\Pi$, which is all right because this Lipschitz constant
depends only on the geometry of $L$. In this case all the 
$\phi_s(w)$ and $\Pi(\phi_s(w))$ lie within $C'' + C \rho$
of $\phi_s(z_0)$, and
$$
\dist(\phi_s^\ast(z),F_0) \leq 
\dist(\phi_s(z_0),F_0) + C'' + C \rho
\leq C + C'' + 1
\leqno (6.53)
$$
if $\rho$ is small enough.
That is, the $\phi_t^\ast$ satisfy (6.31), even though with a larger
geometric constant.

For the analogue of (6.32) we need to check that 
$\phi_t^\ast(z) \in P_0 \times V_0$ when $z\in F_0$. 
We already know from (6.32) that 
$\phi_t(z) \in P_0 \times V_0$, and (6.48) says that 
$\phi_t^\ast(z) = \Pi(\phi_t(z))$. 
Write $\phi_t(z) = (p,v)$, with $p = \pi(\phi_t(z))$ 
and $v = \pi_x(\phi_t(z))$.
Then $\phi_t^\ast(z) = (p,\pi_L(v))$ by (6.47). We know
that $p \in P_0$, so we just need to check that $\pi_L(v) \in V_0$.

Recall that $\dist(\phi_t(z),F_0) \leq C$, either by (6.31) or else
because $\phi_t(z) = z \in F_0$. Choose $w\in F_0$ such that 
$|w-\phi_t(z)| \leq C$. 
Observe that $\pi_x(w) \in E \cap B(x,r) \i L$, by (6.7) and (6.45), and so 
$\pi_L(\pi_x(w)) = \pi_x(w)$ by definition of $\pi_L$. Now
$$
|\pi_L(v)-v| \leq |\pi_L(v) - \pi_L(\pi_x(w))| + |\pi_x(w)-v| 
\leq C |\pi_x(w)-v| \leq C |w-\phi_t(z)| \leq C
\leqno (6.54)
$$
because $\pi_L$ is Lipschitz and $v = \pi_x(\phi_t(z))$,
and where $C$ is a geometric constant that depends on 
the constant in (6.31) and the Lipschitz constant for $\pi_L$.
So we still can choose $N$ in the definition of $P_0$
and $V_0$ (see the discussion below (6.16)) much larger
than this. 

Now $w\in F_0 = F \cap (P_0 \times V_0)$ (by (6.19)),
and (6.17) says that 
$$
\dist(w, P_0 \times (V \cap \d B(0,\rho_0)) \geq N,
\leqno (6.55)
$$
so $\dist(\pi_x(w),V \cap \d B(0,\rho_0)) \geq N$.
Recall that $V_0 = V \cap B(0,\rho_0)$, so 
$\pi_x(w) \in V \cap B(0,\rho_0)$, hence in fact
$\pi_x(w) \in V \cap B(0,\rho_0 - N)$.
Since
$$\eqalign{
|\pi_L(v)-\pi_x(w)| &\leq |\pi_L(v)-v| + |v-\pi_x(w)| \leq C + |v-\pi_x(w)|
\cr&
= C + |\pi_x(\phi_t(z))-\pi_x(w)| \leq C + |\phi_t(z)-w| \leq 2C
}\leqno (6.56)
$$
by (6.54), because $v = \pi_x(\phi_t(z))$ and by definition of $w$,
we get that $\pi_L(v) \in V \cap B(0,\rho_0) = V_0$, as needed. 
This proves (6.32) for the $\phi_t^\ast$.

Observe that (6.33) for $\phi_1^\ast$ follows from its analogue 
for $\phi_1$, by (6.48) and because $\Pi$ is Lipschitz.

As for (6.34), we just need to know that the $\phi_t^\ast$ 
(and in fact $\phi_1^\ast$ is enough) are Lipschitz on $F \sm (P_0 \times V_0)$. 
Indeed, (6.34) (i.e., ($\ast$ 9.13) in [DS4]) 
is only used once, near the end of Section 9.2 of [DS4], 
to prove ($\ast$9.98), and for this we only need the restriction to $F$. 
But then we can use (6.48), and (6.34) for the $\phi_t^\ast$ follow
from (6.34) for the $\phi_t$, because $\Pi$ is $C$-Lipschitz.

This completes the verification of (6.29)-(6.34), but recall that in addition we 
need to check that (1.7) holds, with respect to the quasiminimal set $\wh E$
and the boundaries $\wh L_j = P \times L_j$. That is, we are given 
$z\in \wh E \cap \wh L_j$ and $0 \leq t \leq 1$, and we want to check that 
$\phi_t^\ast(z) \in \wh L_j$. 

The conclusion is trivial if  $\phi_t^\ast(z) = z$, so we may assume that 
$\phi_t^\ast(z) \neq z$. By (6.51), $\dist(z,P_0 \times V_0) \leq d+3$, 
and hence $z \in B(\wh x, r/10)$, by (6.18). In particular,
$\pi_x(z) \in B(x,r/10)$. This excludes the case when $L_j$ does not meet 
$E \cap B(x,r/10)$, because $\pi_x(z) \in E \cap L_j$. In the remaining 
case, $j \in J_0$ and $L \i L_j$ by (6.44).

Return to our $z\in \wh E \cap \wh L_j$. We know that 
$z \in B(\wh x, r/10)$, so (6.7) says that $z\in F$,
and $\phi_t^\ast(z) = \Pi(\phi_t(z))$ by (6.48). Therefore 
$\phi_t^\ast(z) \in P \times L \i P \times L_j = \wh L_j$
by (6.47) and the fact that $\pi_L$ maps $L^\eta$ to $L$
by definition, and as needed for (1.7).

We finally completed our list of verifications; now we can apply 
the fact that $\wh E$ is quasiminimal, compute as in [DS4], 
and get the same conclusion as Proposition $\ast$8.15 there, 
which happen to be the same as in Proposition 6.41, which follows.
\qed

\ms\noindent
{\bf 7. The local uniform rectifiability of $E^\ast$
and bilateral weak geometric lemmas} 
\ms

Our next goal for this section is to extend Proposition 6.41 to the case
when only the assumptions of Theorem 6.1 are satisfied; see
Proposition 7.85 below. Then, in Section 8, we shall take care of the 
difference between big pieces of bilipschitz images 
and big pieces of Lipschitz graphs, and prove Theorem 6.1.

For our first verification, we shall mostly use Proposition 6.41
itself, the smallness or regularity of the faces that compose the $L_j$,
and general knowledge on uniformly rectifiable sets; Lemma 7.38 
will be the only place where we use the quasiminimality of $E$ in this 
argument, to show that a quasiminimal set that stays very close to 
the interior of a $d$-face does not have big holes there.

We shall use a characterization of uniform rectifiability by
the so-called bilateral weak geometric lemma. We are given a locally
Ahlfors-regular set $F$ of dimension $d$, which we want to study; 
our main example will be $E^\ast$, or a piece of $E^\ast$.
First define the standard P. Jones numbers $\beta(x,r)$ by
$$
\beta(x,r) = \beta_F(x,r) = \inf_{P}\Big\{ {1 \over r}\sup_{y\in F\cap B(x,r)}
\dist(y,P) \Big\},
\leqno (7.1)
$$
where $x\in F$, $r > 0$, and
the infimum is taken over all the affine $d$-planes $P$ through $x$.
It is just as convenient here to restrict to planes that contain $x$,
even though the other option could be used too, and would give 
equivalent results. We shall also use the bilateral variant
$$
b\beta(x,r) = b\beta_F(x,r) = \inf_{P}\Big\{ 
{1 \over r}\sup_{y\in F\cap B(x,r)}\dist(y,P) 
+{1 \over r}\sup_{y\in P\cap B(x,r)}\dist(y,F)
\Big\},
\leqno (7.2)
$$
where we also account for big holes in the middle of $F$.
We shall be interested in the size of the bad sets
$$
{\cal B}(\varepsilon) = {\cal B}_F(\varepsilon) = 
\big\{ (x,r) \in F \times (0,+\infty) \, ; \,
b\beta(x,r) > \varepsilon \big\}
\leqno (7.3)
$$
when $\varepsilon > 0$ is small enough, because the local 
uniform rectifiability of (Ahlfors regular) sets turns out to be 
equivalent to Carleson measure estimates on the
${\cal B}(\varepsilon)$. 

The following result concerns unbounded Ahlfors-regular sets; 
it will need to be adapted to the present situation, 
but gives an idea of what we want to do. We consider a closed 
(unbounded) Ahlfors-regular set $F$ of dimension $d$. 
This last means that there is a constant $C_0 \geq 1$ such that
$$
C_0^{-1} r^d \leq \H^d(F \cap B(x,r)) \leq C_0 r^d
\ \hbox{ for $x\in F$ and } 0 < r < +\infty.
\leqno (7.4)
$$
We shall say that $F \in BWGL(\varepsilon, C(\varepsilon))$
(or that $F$ satisfies a bilateral weak geometric lemma, with the
constants $\varepsilon$ and $C(\varepsilon)$)
when 
$$
\int_{y\in F \cap B(x,r)}\int_{0 < t < r} 
{\bf 1}_{{\cal B}(\varepsilon)}(y,t) \, {d\H^d(y) dt \over t}
\leq C(\varepsilon) r^d
\leqno (7.5)
$$
for $x\in F$ and $0 < r < +\infty$.

We say that $F\in BPBI(\theta,C_1)$ (for big pieces of bilipschitz 
images) when for all $x\in F$ and $r > 0$, 
we can find a closed set $G_0 \i F \cap B(x,r)$ and a Lipschitz mapping 
$\phi : G_0 \to \R^d$ such that 
$$
\H^d(G_0) \geq \theta r^d
\ \hbox{ and } \ 
C_1^{-1} |y-z| \leq |\phi(y)-\phi(z)| \leq C_1 |y-z|
\hbox{ for } y,z \in G_0.
\leqno (7.6)
$$
Thus this is the same property as in Proposition 6.41, except
that there we restricted to the pairs $(x,r)$ such that
$0 < r < \Min(r_0,\delta)$ and $B(x,2r) \i B_0$, and there was
an extra assumption to get it.

\ms\proclaim Theorem 7.7.
Let $F \i \R^n$ be a closed Ahlfors-regular set of dimension
$d$. If $F\in BPBI(\theta,C_1)$ for some choice of 
$\theta > 0$ and $C_1 \geq 1$, then for every $\varepsilon > 0$,
$F \in BWGL(\varepsilon, C(\varepsilon))$ for some
$C(\varepsilon)$ that depends only on $n$, $C_0$ (the regularity
constant in (7.4)), $\theta$, $C_1$, and $\varepsilon$.
Conversely, there exists $\varepsilon > 0$, that depends only
on $n$ and $C_0$, such that if 
$F \in BWGL(\varepsilon, C(\varepsilon))$ for some 
$C(\varepsilon) \geq 0$, then there exist 
$\theta > 0$ and $C_1 \geq 1$,
that depend only on $n$, $\varepsilon$, and $C(\varepsilon)$, such 
that $F\in BPBI(\theta,C_1)$.

\ms
Notice that $BWGL(\varepsilon, C(\varepsilon))$ is smaller
when $\varepsilon$ is smaller, so that in the converse statement,
assuming that $F \in BWGL(\varepsilon', C(\varepsilon'))$ 
for some $\varepsilon' \leq \varepsilon$ would be enough too.

Theorem 7.7 follows from Theorems 2.4 and 1.57 in [DS3]; 
see the remark above Theorem 2.4 (for the fact that only one
small $\varepsilon$ is needed), Definitions 2.2 and 1.69 (for the
BWGL), and (1.60) and (1.61) (for two of the equivalent conditions in 
Theorem 1.57, one clearly stronger and one clearly weaker than 
our BPBI here).

We shall use both parts of the equivalence here. We start with
the direct part.

\ms\proclaim Lemma 7.8.
Let $E$, $x\in E^\ast$, and $r$ be as in Proposition 6.41.
In particular, assume that (6.42) holds.
Then for each $\varepsilon > 0$
$$
\int_{y\in E^\ast \cap B(x,r/8)}\int_{0 < t < r/8} 
{\bf 1}_{{\cal B}_{E^\ast}(\varepsilon)}(y,t) \, {d\H^d(y) dt \over t}
\leq C(\varepsilon) r^d,
\leqno (7.9)  
$$
where $C(\varepsilon)$ depends only on $n$, $M$, and $\varepsilon$.

\ms
We want to deduce Lemma 7.8 from Proposition 6.41 and Theorem 7.7,
but a small localization argument is needed. 
To this effect, we apply Proposition 7.6 in [DS4]. 
This is the same proposition that we used to get the set
$F$ in (6.7), but here we shall need to use it in a more precise 
way. We apply it to the set $E^\ast$ and the ball
$B(x,r/2)$; the assumptions follow from the local 
Ahlfors-regularity of $E^\ast$ that we proved in 
Proposition 4.1.

What we get from Proposition 7.6 in [DS4] 
is a (bounded) Ahlfors-regular set $F$ such that
$$
E^\ast \cap B(x,r/4) \i F \i E^\ast \cap B(x,r/2),
\leqno (7.10)
$$
and a cubical patchwork for $F$, i.e., collections
of decompositions of $F$ into pseudocubes $Q$, $Q \in \Sigma_j$,
like the one we used near (6.9). But this time we shall find it
convenient to use the fact that our patchwork is adapted
to $E$, in the sense that for every cube $Q \in \cup_{j \geq 0} 
\Sigma_j$ and every (small) $\tau > 0$,
$$\eqalign{
&\H^d\big(\big\{ w\in Q \, ; \, \dist(w,E \sm Q) \leq \tau \diam(Q) \big\}\big)
\leq C \tau^{1/C} \diam(Q)^d.
}\leqno (7.11)
$$
This is part of ($\ast$7.10) in [DS4], 
and this is a little more precise than the usual ``small boundary
condition" on cubes of $F$, because it also controls the difference 
between $E$ and $F$. Here $C$ in (7.11) depends on $n$ and $M$
(through the local Ahlfors-regularity constant), but not on $\tau$.

We want to apply Theorem 7.7, and since $F$ is not unbounded,
we replace it with $F' = F \cup P$, where $P$ is any affine $d$-plane
such that $\dist(x,P) = 2r$. It is easy to see that $F'$ is 
Ahlfors-regular (as in (7.4)), and now we want to check that
$F'\in BPBI(\theta,C_1)$, for some choice of $\theta > 0$ and
$C_1 \geq 1$ which will be just a little worse than the constants of 
Proposition 6.41.

So we pick $y\in F'$ and $t>0$ and try to find a big bilipschitz piece 
$G_0$ in $F' \cap B(y,t)$, as in (7.6). If $t > 3r$ or $y\in P$, 
$F'$ contains a $d$-dimensional disk of radius $t/3$ (contained in 
$P$), which is a nice choice of $G_0$ for (7.6). 
So assume that $y\in F$ and $t<3r$, and try to find
$G_0 \i F$. Let $Q$ be a cube of our patchwork such that 
$$
y \in Q \i F \cap B(y,t/10),
\leqno (7.12)
$$
and which is maximal with these properties.
Thus ($\ast$7.2) in [DS4] 
says that $\diam(Q) \geq t/C$, where $C$ depends on
$n$ and $M$ (through the local Ahlfors-regularity constant in 
Proposition~4.1). Also recall from ($\ast$7.2) that
$\H^d(Q) \geq C^{-1} \diam(Q)^d$; then choose the constant 
$\tau \in (0,1/10)$ so small (again depending on $n$ and $M$ only) 
that the right-hand side of (7.11) is smaller than $\H^d(Q)$.
This allows us to find $w\in Q$ such that 
$$
\dist(w,E^\ast \sm F) \geq \dist(w,E^\ast \sm Q) \geq \tau \diam(Q).
\leqno (7.13)
$$

We want to apply Proposition 6.41 to the pair $(w,\tau \diam(Q))$,
so we check the hypotheses. First, $w\in E^\ast$ because $F \i E^\ast$
(see (7.10)). Next, 
$$
\tau \diam(Q) \leq \tau t/5  \leq 3 \tau r/5 \leq 3r/50
\leqno (7.14)
$$
by (7.12), because $t<3r$, and because we chose $\tau < 1/10$. 
In particular $\tau \diam(Q) < r < \Min(r_0,\delta)$, and also 
$$\eqalign{
B(w,2\tau \diam(Q)) &\i B(y,2 \tau \diam(Q) + t/10) 
\i B(x,2 \tau \diam(Q) + t/10 + r/2)
\cr&
\i B(x,6r/50+3r/10+r/2) \i B(x,r) \i B_0
}\leqno (7.15)
$$
by (7.12), because $y\in F$, and by (7.10) and (7.14). Finally, (6.42)
for $B(w,\tau \diam(Q))$ follows from (6.42) for $B(x,r)$,
simply because $B(w,\tau \diam(Q)) \i B(x,r)$ by (7.15).
So Proposition 6.41 applies, and gives a set $G_0 \i E^\ast \cap 
B(w,\tau \diam(Q))$, such that
$$
\H^d(G_0) \geq \theta \tau^d \diam(Q)^d
\ \hbox{ and } \ 
C_M^{-1} |y-z| \leq |\phi(y)-\phi(z)| \leq C_M |y-z|
\hbox{ for } y,z \in G_0.
\leqno (7.16)
$$
This set works in the definition (7.6), because $G_0 \i F$ by
(7.13). Note that since $\diam(Q) \geq C^{-1}t$ (see below (7.12)), 
$\theta \tau^d \diam(Q)^d \geq \theta' t^d$
for some $\theta'$ that depends only on $n$ and $M$.
Thus $F'\in BPBI(\theta',C_M)$, as needed, and now Theorem 7.7
says that for every choice of $\varepsilon > 0$, we can find
$C(\varepsilon)$ so that (7.5) holds for $F'$ (and any ball
centered on $F'$). We just apply this to the ball
$B(x,r/8)$ (which is centered on $F'$ by (7.10), and get that
$$
\int_{y\in F' \cap B(x,r/8)}\int_{0 < t < r/8} 
{\bf 1}_{{\cal B}_{F'}(\varepsilon)}(y,t) \, {d\H^d(y) dt \over t}
\leq 8^d C(\varepsilon) r^d.
\leqno (7.17)
$$
But $F'$ coincides with $E^\ast$ on $B(x,r/4)$ by (7.10),
so (7.17) is just the same as (7.9). Lemma~7.8 follows.
\qed.

\ms
We slowly return to the extension of Proposition 6.41 to the
situation of Theorem 6.1. We fix a small $\varepsilon > 0$,
and want to control the size of ${\cal B}_{E^\ast}(\varepsilon)$,
by cutting it into smaller pieces that we can control.
Denote by
$$
{\cal A} = \big\{ (y,t) \in E^\ast \times \Min(r_0,\delta)
\, ; \, B(y,2t) \i B_0 \big\}
\leqno (7.18)
$$
the set of balls that we like to consider.
We first get rid of the balls that lie close to a  
face of dimension at most $d-1$ in our initial net.

\ms\proclaim Lemma 7.19. Denote by $F_1$ the union of all the faces of 
dimension at most $d-1$ of cubes from the dyadic net that was used 
to define the $L_j$, and set
$$
{\cal B}_1 =  \big\{ (y,t) \in {\cal A} \, ; \, 
B(y,10t) \hbox{ meets } F_1 \big\}.
\leqno (7.20)  
$$
Then 
$$
\int_{y \in E^\ast \cap B(x,r)}\int_{0 < t < r} 
{\bf 1}_{{\cal B}_1}(y,t) \, {d\H^d(y) dt \over t}
\leq C_1 r^d
\leqno (7.21)
$$
for $(x,r) \in {\cal A}$, with a constant $C_1$
that depends only on $n$ and $M$.

\ms
For $0 < t \leq r$, cover $F_1 \cap B(x,20r)$ with balls $B_i$, 
$i\in I(t)$, of the same radius $10t$.
We can do this with less than $C (r/t)^{d-1}$ balls, i.e., so that
$\sharp I(t) \leq C (r/t)^{d-1}$. Then the local Ahlfors-regularity
given by Proposition 4.1 yields
$$\leqalignno{
\int_{y \in E^\ast \cap B(x,r)} \int_{0 < t < r}& 
{\bf 1}_{{\cal B}_1}(y,t) \, {d\H^d(y) dt \over t}
= \int_{0 < t < r} \H^d(\{ y\in E^\ast \cap B(x,r); \dist(y,F_1) < 10t \})
\, {dt \over t}
\cr&\hskip0.5cm
\leq \int_{0 < t < r} \sum_{i \in I(t)} \H^d(\{ E^\ast \cap B(x,r) \cap 2B_i)) 
\, {dt \over t}
\leq C \int_{0 < t < r} \sharp I(t) \, t^d \, {dt \over t}
\cr& \hskip0.5cm
\leq C r^{d-1} \int_{0 < t < r} dt  \leq C r^d,
& (7.22) 
}
$$
where the fact that $\H^d(\{ E^\ast \cap B(x,r) \cap 2B_i)) \leq C t^d$
comes from Proposition 4.1 and the fact that $B(x,2r) \i B_0$; this 
proves (7.21).
\qed

\ms
Next we consider the pairs $(y,t) \in {\cal A} \sm {\cal B}_1$
such that $B(y,2t)$ meets a $d$-dimensional face, without staying too 
close to it.

\ms\proclaim Lemma 7.23.
Denote by $F_2$ the union of all the $d$-dimensional 
faces of cubes from the dyadic net that was used to define 
the $L_j$, and set
$$\eqalign{
{\cal B}_2 &=  \big\{ (y,t) \in {\cal A} \sm {\cal B}_1 \, ; \, 
B(y,2t) \hbox{ meets $F_2$ but there exists }
\cr&\hskip 1cm
w\in E^\ast \cap B(y,2t) \hbox{ such that } \dist(w,F_2) \geq \varepsilon t
\big\}.
}\leqno (7.24)  
$$
Then
$$
\int_{y \in E^\ast \cap B(x,r)}\int_{0 < t < r/10} 
{\bf 1}_{{\cal B}_2}(y,t) \, {d\H^d(y) dt \over t}
\leq C_2(\varepsilon) r^d
\leqno (7.25)
$$
for $(x,r) \in {\cal A}$, with a constant $C_2(\varepsilon)$
that depends only on $n$, $M$, and $\varepsilon$.

\ms
Let $(x,r) \in {\cal A}$ be given; we want to estimate
the left-hand side of (7.25), which we write as
$$
\Lambda = \int\int_{(y,t) \in {\cal B}_2(x,r)}  \, {d\H^d(y) dt \over t}
\leqno (7.26)
$$
where ${\cal B}_2(x,r) = \big\{ (y,t) \in {\cal B}_2  
\, ; \, y\in B(x,r) \hbox{ and } 0 < t < r/10 \big\}$. For each
$(y,t) \in {\cal B}_2(x,r)$ we use the definition (7.24) to pick 
$w\in E^\ast \cap B(y,2t)$ such that $\dist(w,F_2) \geq \varepsilon t$, 
and we set
$$
Z(y,t) = E^\ast \cap B(w,\varepsilon t/2).
\leqno (7.27)
$$
Obviously
$$
\dist(z,F_2) \geq \varepsilon t/2
\ \hbox{ for } z \in Z(y,t).
\leqno (7.28)
$$
Also, $|w-y| \leq 2t$ since $w\in B(y,2t)$, then $|w-x| \leq 2t+r 
\leq 12r/10$,
because $y\in B(x,r)$, so 
$$
B(w,\varepsilon t) \i B(x,13r/10) \i B_0;
\leqno (7.29)
$$
then we can apply Proposition 4.1 to $B(w,\varepsilon t)$
and get that
$$
\H^d(Z(y,t)) \geq C^{-1} \varepsilon^d t^d.
\leqno (7.30)
$$
By choosing $w$ out of a fixed countable dense subset of $E^\ast$,
we can make sure that the relation ``$z \in Z(y,t)$" is measurable
in all variables. Then 
$$
\Lambda \leq C \varepsilon^{-d} \int\int_{(y,t) \in {\cal B}_2(x,r)}
\int_{z\in Z(y,t)} t^{-d} \, {d\H^d(y) d\H^d(z) dt \over t}.
\leqno (7.31)  
$$
Let us use Fubini. In the domain of integration,
$z$ lies in $E^\ast \cap B(x,13r/10)$ by (7.29), and
then $|y-z| \leq |y-w| + |w-z| \leq 2t+\varepsilon t$. So
$y\in E^\ast \cap B(z,3t)$, whose
$\H^d$-measure is less than $C t^d$, by Proposition 4.1
and because $B(z,6t) \i B(z,6r/10) \i B(x,2r) \i B_0$. 
In addition, $\varepsilon t/2 \leq \dist(z,F_2)$ by (7.28),
and $\dist(z,F_2) \leq |z-w| + |w-y| + \dist(y,F_2)
\leq \varepsilon t/2 + 2t + 2t \leq 5 t$ by (7.24) in particular.
Therefore, setting $d(z) = \dist(z,F_2) > 0$ to save space, 
$$\eqalign{
\Lambda &\leq C \varepsilon^{-d}
\int_{z \in E^\ast \cap B(x,13r/10)}
\int_{5^{-1}d(z) \leq t \leq 2 \varepsilon^{-1}d(z)} 
t^{-d} \H^d(E^\ast \cap B(z,3t)) \, {d\H^d(z) dt \over t}
\cr&
\leq C \varepsilon^{-d} \int_{z \in E^\ast \cap B(x,13r/10)}
\int_{5^{-1}d(z) \leq t \leq 2 \varepsilon^{-1}d(z)} {d\H^d(z) dt \over t}
\cr&
\leq C \varepsilon^{-d} \log(10/\varepsilon) \, \H^d(E^\ast \cap B(x,13r/10)) 
\leq C \varepsilon^{-d} \log(10/\varepsilon) \, r^d,
}\leqno (7.32)  
$$
by Proposition 4.1, this time applied to a few balls $B_i$
of radius $r/10$ that cover $B(x,13r/10)$, to make sure 
that the $2B_i$ are contained in $B_0$. This completes the 
proof of Lemma 7.23.
\qed

\ms\proclaim Lemma 7.33. Suppose that $\varepsilon$ is small
enough, depending on $n$ and $M$.  
Let $(y,t) \in {\cal A} \sm ({\cal B}_1 \cup{\cal B}_2)$
be such that $t < r_0/10$ and 
$B(y,2t)$ meets $F_2$. Then
$b\beta_{E^\ast}(y,t) \leq 4\varepsilon$.

\ms
Let $(y,t)$ be as in the statement. Since 
$(y,t) \notin {\cal B}_2$, we know that 
$$
\dist(w,F_2) < \varepsilon t \hbox{ for every } w\in E^\ast \cap B(y,2t).
\leqno (7.34)  
$$
In particular, $\dist(y,F_2) < \varepsilon t$ and there is a
$d$-dimensional face $F$ from our usual net such that
$\dist(y,F) < \varepsilon t$. Let us check that
$$
\dist(y,F') \geq 8t \hbox{ for every face $F'$ of the usual net
such that $F \not\i F'$.}
\leqno (7.35)  
$$
Let $\d F$ denote the boundary of $F$; this is a union
of $(d-1)$-dimensional faces, and 
$\dist(y,\d F) \geq \dist(y,F_1) \geq 10t$
by definition of $F_1$ and because $(y,t) \notin {\cal B}_1$.
Also let $f\in F$ be such that $|y-f| \leq \varepsilon t$; then
(3.8) says that
$$
\dist(f,F') \geq \dist(f,\d F) \geq \dist(y,\d F) - \varepsilon t
\geq 10t - \varepsilon t \geq 9t
\leqno (7.36) 
$$
(in this lemma, we may assume that $\varepsilon$ is as small
as we want, but anyway $\varepsilon >1$ does not make sense because 
all the $\beta$-numbers are $\leq 1$); (7.35) follows at once.
Of course this would be easy to adapt to polyhedral nets.

Notice that (7.35) applies to every $d$-dimensional face $F' \neq F$ 
of our net, so $\dist(y,F_2 \sm F) \geq 8t$ and (7.34) implies that
$$
\dist(w,F) < \varepsilon t \hbox{ for every } w\in E^\ast \cap B(y,2t).
\leqno (7.37)  
$$
Now we we need to use the quasiminimality of $E$ to prove that 
$E$ has no apparent hole in $B(y,t)$; the next lemma is just a 
little more general than what we need; also, we shall need to return 
to its proof and generalize it in Lemma 9.14.

\ms\proclaim Lemma 7.38.
Let $C_0 \geq 1$ be given. 
Let $E \in GSAQ(B_0, M, \delta , h)$, and suppose that the 
rigid assumption is satisfied and that $h$ is small enough,
depending on $n$, $M$, and $C_0$. Let $y\in E^\ast$ and $t>0$
be such that $0 < t < \Min(r_0,\delta)$ and $B(y,2t) \i B_0$.
Let $P \i \R^n$ be a $d$-plane, and assume that
$$
\dist(w,P) < \varepsilon t \hbox{ for } w\in E^\ast \cap B(y,2t)
\leqno (7.39)  
$$
for some $\varepsilon > 0$ that we assume to be small enough,
depending on $n$, $M$, and $C_0$. 
Also suppose that 
$$
P \cap B(y,2t) \i L_j \ \hbox{ for every $j$ such that $L_j$ meets $B(y,2t)$}
\leqno (7.40) 
$$
and that we have a Lipschitz
function $h : E^\ast \cap B(y,2t)\times [0,1] \to \R^n$
such that 
$$
h(w,0) = w \hbox{ and } h(w,1) = \pi(w)
\ \hbox{ for } \ w\in E^\ast \cap B(y,2t),
\leqno (7.41) 
$$
where $\pi$ denotes the orthogonal projection on $P$,
$$
|h(w,s) - h(w,s')| \leq C_0 \varepsilon t |s-s'|
\ \hbox{ for } \ w\in E^\ast \cap B(y,2t) \hbox{ and }
0 \leq s,s' \leq 1,
\leqno (7.42) 
$$
$$
|h(w,s) - h(w',s)| \leq C_0 |w-w'|
\ \hbox{ for } \ w,w' \in E^\ast \cap B(y,2t) \hbox{ and }
0 \leq s \leq 1,
\leqno (7.43) 
$$
and finally
$$
h(w,s) \in L_j  \hbox{ for } 0 \leq s \leq 1 \ \hbox{ whenever }
w \in E^\ast \cap L_j \cap B(y,2t).
\leqno (7.44) 
$$
Then
$$
\dist(p,E^\ast) \leq \varepsilon t \hbox{ for } 
p\in P \cap B(y,3t/2)
\leqno (7.45)  
$$
and 
$$
\pi(E^\ast \cap B(y,5t/3)) \hbox{ contains } 
P \cap B(\pi(y),3t/2).
\leqno (7.46)  
$$

\ms
Of course the simplest choice of path $h$ is to take
$h(w,s) = s \pi(w) + (1-s) w$, but it does not always
work, because it may be more efficient to follow the faces
of dyadic cubes to stay in the $L_j$ and get (7.44).
We keep this type of issues for the next sections. 

\ms
Let us first check that Lemma 7.38 implies Lemma 7.33. Let
$(y,t)$ and $F$ be as in Lemma 7.33, let $P$ be the $d$-plane 
that contains $F$; then (7.39) follows from (7.37).

Next let $j \leq j_{ max}$ be such that $L_j$ meets $B(y,2t)$.
Let $F'$ be a face of $L_j$ that meets $B(y,2t)$; by (7.35),
$F'$ contains $F$. Also, $\dist(y,F) \leq \varepsilon t$
by definition of $F$, and $\dist(y,\d F) \geq \dist(y,F_1)
\geq 10t$ because $(y,t) \notin {\cal B}_1$ or by the end of (7.36),
so $P \cap B(y,2t) \i F \i F' \i L_j$, and (7.40) holds.

We take $h(w,s) = s \pi(w) + (1-s) w$, then $h$ is Lipschitz,
(7.41) holds trivially, (7.42) and (7.43) are true with
$C_0 = 1$ because $|\pi(w)-w| \leq \varepsilon t$
for $w\in E^\ast \cap B(y,2t)$, by (7.39), and because $\pi$ is 
$1$-Lipschitz.

We finally check (7.44), i.e., that $[w,\pi(w)] \i L_j$ when
$w \in E^\ast \cap L_j \cap B(y,2t)$.
Observe that $P$ is defined by some equations
$w_i = n_i r_0$, with $n_i \in {\Bbb Z}$, and that $\pi(w)$ is obtained
from $w$ by replacing the $w_i$ such that $w_i \neq n_i r_0$
with $n_i r_0$. When this happens, $|w_i - n_i r_0| \leq \varepsilon t$,
because $|\pi(w)-w| \leq \varepsilon t$. Then $[w,\pi(w)]$ is 
contained in any face of any dimension that contains $w$
(we just replace some noninteger coordinates of $r_0^{-1} w$ with other 
ones that lie in the same dyadic intervals). We apply this to any face
of $L_j$ that contains $w$ and get that $[w,\pi(w)] \i L_j$,
as needed for (7.44). 

So Lemma 7.38 applies. If we could use $P$ in the definition of
$b\beta_{E^\ast}(y,t)$, (7.39) and (7.45) would imply that 
$b\beta_{E^\ast}(y,t) \leq 2\varepsilon$. We cannot exactly,
because maybe $P$ does not contain $y$, but 
$\dist(y,P) \leq \varepsilon t$ by (7.39), so we can use a small
translation of $P$ that goes through $y$, and we get that
$b\beta_{E^\ast}(y,t) \leq 4\varepsilon$, as needed.
Hence Lemma 7.33 will follow from Lemma 7.38 as soon as we prove it.

\ms
Lemma 7.38 is a variant of Lemma 10.10 in [DS4], 
but we need to modify some things because of the boundary 
constraints (1.7). 
We define a first family of deformations $\varphi_s$.
First let $\psi : [0,+\infty) \to [0,1]$ be
such that 
$$\eqalign{
&\psi(\rho) = 1 \hbox{ for } 0 \leq \rho \leq {5t \over 3} + (C_0+1)\varepsilon t,
\cr& \psi(\rho) = 0 \hbox{ for } \rho \geq {5t \over 3} + (C_0+2)\varepsilon t, 
\hbox{ and }
\cr&
\psi \hbox{ is affine on } [{5t \over 3} 
+ (C_0+1)\varepsilon t,{5t \over 3} + (C_0+2)\varepsilon t].
}\leqno (7.47)  
$$
Then set
$$
\varphi_s(w) = h(w, s \psi(|w-y|))
\ \hbox{ for $w\in E^\ast$ and } 0 \leq s \leq 1;
\leqno (7.48)  
$$
the fact that $h(w,s)$ is only defined for $w \in E^\ast \cap B(y,2t)$
is not a problem, because we can set 
$$
\varphi_s(w) = w
\ \hbox{ for $w\in E^\ast \sm B(y,11t/6)$ and } 0 \leq s \leq 1,
\leqno (7.49)  
$$
where the two definitions make coincide on $B(y,2t) \sm B(y,11t/6)$, 
if $\varepsilon$ is so small that ${5t \over 3} + (C_0+2)\varepsilon t < 11t/6$, 
because $s \psi(|w-y|) = 0$ there. 
To see that $(s,w) \to \varphi_s(w)$ is Lipschitz, we observe that 
the two definitions yield Lipschitz functions and coincide in 
$B(y,2t) \sm B(y,11t/6)$.

We shall not use the $\varphi_t$ as they are,
but let us check that they satisfy the properties (1.4)-(1.8),
with the closed ball
$$
B = \overline B(y,11t/6)
\leqno (7.50) 
$$ 
and with respect to the set $E^\ast$. 
First observe that we just checked (1.4), and that (1.5) and (1.8) 
are very easy consequences of the definition.

For (1.6), let $w\in B$ and $0 \leq s \leq 1$ be given; we 
want to check that $\varphi_s(w) \in B$.
This is trivial if $\varphi_s(w) = w$, so we may assume that
$w \in B(y,2t)$ and $\varphi_s(w)$ is given by (7.48).
Then $s\psi(|w-y|) \neq 0$, and hence 
$|w-y| < {5t \over 3} + (C_0+2)\varepsilon t$.
Notice that
$$
|\varphi_s(w)-w| = |h(w, s \psi(|w-y|))-w|\leq C_0 \varepsilon t
\ \hbox{ for $w\in E^\ast$ and } 0 \leq s \leq 1
\leqno (7.51) 
$$ 
by (7.41) and (7.42) if $w \in B(y,2t)$, and because
$\varphi_s(w) = w$ otherwise. If $\varepsilon$ is small enough, 
$\varphi_s(w)\in B$ when 
$|w-y| < {5t \over 3} + (C_0+2)\varepsilon t$, as needed for (1.6).

Finally (1.7) holds because if $x\in E^\ast \cap L_j \cap B$ 
and $s\in [0,1]$, then 
$x\in B(y,2t)$ and $\varphi_s(x) \in L_j$ by (7.48) and (7.44).

\ms
Now we shall assume that (7.46) fails, use this to construct a 
deformation that completes the $\varphi_t$ and makes
$E^\ast \cap B(y,t)$ essentially vanish, and get a contradiction. 
So let us assume that we can find
$$
p \in P \cap B(\pi(y),3t/2) \sm \pi(E^\ast \cap B(y,5t/3)).
\leqno (7.52)  
$$  
Observe that
$$
\pi(\varphi_1(w)) \hbox{ lies out of $B(\pi(y),3t/2)$ for } 
w \in E^\ast \cap B(y,2t) \sm B(y,5t/3)
\leqno (7.53)  
$$
just because 
$|\pi(\varphi_1(w))-\pi(y)| 
\geq |w-y| - |\pi(\varphi_1(w)) - w| - |\pi(y)-y|
\geq {5t \over 3} - |\pi(\varphi_1(w)) - w| - \varepsilon t$ by (7.39)
and $|\pi(\varphi_1(w)) - w| \leq |\pi(\varphi_1(w)) - \pi(w)|+
|\pi(w)-w| \leq |\varphi_1(w) - w|+|\pi(w)-w| \leq (C_0+1) \varepsilon t$
by (7.51) and (7.39).
Thus
$$
p \in P \cap B(\pi(y),3t/2) \sm \pi(E^\ast \cap B),
\leqno (7.54)  
$$
where $B = \overline B(y,11t/6)$ as before, by (7.52) and (7.53).
Since $E^\ast \cap B$ is compact, we can find $\tau > 0$
(possibly extremely small) such that
$$
P \cap B(p,\tau) \hbox{ does not meet } \pi(E^\ast \cap B).
\leqno (7.55)  
$$
Define $g : P \cap B(\pi(y),{5t \over 3}) \sm B(p,\tau) \to P \cap \d 
B(\pi(y),{5t \over 3})$ as the radial projection centered at $p$, i.e., 
by the fact that 
$$
g(w) \in \d B(\pi(y),{5t \over 3}) \hbox{ and }
w \in [p,g(w)]. 
\leqno (7.56)  
$$
Also set $g(z) = z$ for $z \in \d B(\pi(y),{5t \over 3})$;
this gives a Lipschitz mapping defined on the union, and with values
in $\overline B(\pi(y),{5t \over 3})$. We extend $g$ to 
$\overline B(\pi(y),{5t \over 3})$ in a Lipschitz way, with values 
in $\overline B(\pi(y),{5t \over 3})$ (use the Whitney extension theorem, and
compose with the radial projection on $\overline B(\pi(y),{5t \over 3})$
if needed). Finally extend $g$ to $\R^n$ by setting 
$$
g(z) = z \ \hbox{ for }
z \in \R^n \sm B(\pi(y),{5t \over 3}).
\leqno (7.57)  
$$
This yields a Lipschitz mapping defined on $\R^n$, which we still call $g$.

Now we define the $\varphi_s$, $1 \leq s \leq 2$, by
$$
\varphi_s(w) = (2-s) \varphi_1(w) + (s-1) g(\varphi_1(w))
\ \hbox{ for $w\in \R^n$ and } 1 \leq s \leq 2.
\leqno (7.58)  
$$
Let us check the analogue of (1.4)-(1.8) for the 
$\varphi_s$, $0 \leq s \leq 2$, with the same choice of 
$B = \overline B(y,11t/6)$ and again with respect to $E^\ast$.

The mapping $(s,w) \to \varphi_s(w)$ is Lipschitz on
$[1,2] \times E^\ast$, so (1.4) and (1.8) hold. We already know
that $\varphi_0(w) = w$ for $w\in \R^n$. 

Next, if $w\in E^\ast \sm B$, we know from our earlier 
verification of (1.5) that $\varphi_s(w) = w$ for $0 \leq s \leq 1$, 
and in particular $\varphi_1(w) = w \notin B$, 
hence $\varphi_1(w) \notin B(\pi(y),{5t \over 3})$
(recall that $|\pi(y)-y| \leq \varepsilon t$), and
$g(\varphi_1(w)) = \varphi_1(w)$ by (7.57). Then
$\varphi_s(w) = w$ for $1 \leq s \leq 2$ by (7.58),
and the analogue of (1.5) holds.

If $w\in B$, we know that $\varphi_s(w) \in B$ for $0 \leq s \leq 1$;
then $g(\varphi_1(w)) \in B$ because $g(B(\pi(y),{5t \over 3})) \i 
B(\pi(y),{5t \over 3})$ and $g(z) = z$ out  of $B(\pi(y),{5t \over 3})$.
So (1.6) holds because the $\varphi_s(w)$, $s \geq 1$, lie on
the segment $[\varphi_1(w),g(\varphi_1(w))] \i B$.

We are left with (1.7) to check, and again it is nice to do this 
relatively to $E^\ast$ (and not the full $E$).
Let $j$ and $w\in E^\ast \cap L_j \cap B$ be given; we want to show 
that $\varphi_s(w) \in L_j$ for $1 \leq s \leq 2$ (we already know 
this for $s \leq 1$). First assume that
$$
w \in B(y,{5t \over 3} + (C_0+1)\varepsilon t).
\leqno (7.59)
$$
Then $\psi(|w-y|) = 1$ by (7.47) and 
$\varphi_1(w) = h(w,1) = \pi(w)$ by (7.48) and (7.41).
In particular, $\varphi_1(w) \in P \sm B(p,\tau)$,
by (7.55)). If $\varphi_1(w) \in B(\pi(y),{5t \over 3})$,
then $g(\varphi_1(w))$ is the radial projection of $\varphi_1(w)$ on 
$\d B(\pi(y),{5t \over 3})$ (as in (7.56)); otherwise,
$g(\varphi_1(w)) = \varphi_1(w)$ by (7.57); in both cases,
$\varphi_s(y) \in [\varphi_1(w),g(\varphi_1(w))] \i P \cap B$
(recall that $\varphi_s(y) \in B$ by (1.6)). Now 
$w\in L_j \cap B \i B(y,2t)$, so $P \cap B \i P \cap B(y,2t) \i L_j$
by (7.40) and  $\varphi_s(y) \in L_j$ when (7.59) holds.

If (7.59) fails, $w \in B \sm B(y,{5t \over 3} + (C_0+1)\varepsilon t)$.
Then $\varphi_1(w)$ lies out of $B(y,{5t \over 3} + \varepsilon t)$
by (7.51), hence also out of $B(\pi(y),{5t \over 3})$. In this case,
$g(\varphi_1(w)) = \varphi_1(w)$, hence (7.58) and 
(1.7) for $0 \leq s \leq 1$ yield
$\varphi_s(w) = \varphi_1(w) \in L_j$, as needed.

This completes our proof of (1.7) for the $\varphi_s$, 
$0 \leq s \leq 2$. Note also that (2.4) holds, 
because $\wh W \i B \i\i B_0$ since we assumed that $B(y,2t) \i B_0$.
We can now apply (2.5), because Proposition 3.3
says that $E^\ast$ is quasiminimal just like $E$. 
This is one instance where we use Proposition 3.3 for real;
of course we could also have assumed that (7.41)-(7.45) hold
with the whole $E_k$, or worked more here to extend our $\varphi_t$ 
correctly. Anyway, we get that
$$
\H^d(W_2) \leq M \H^d(\varphi_2(W_2)) + h r^d,
\leqno (7.60)
$$
where we set 
$W_2 = \big\{ w\in E^\ast \cap B \, ; \, \varphi_2(w) \neq w \big\}$ 
as in Definition 2.3.

Let us first control $\varphi_2$ on
$A_1 = \big\{ w\in E^\ast\cap 2B \, ; \, 
\varphi_1(w) \in B(\pi(y),{5t \over 3}) \big\}$.
We claim that
$$
\varphi_2(A_1) \i P \cap \d B(\pi(y),5t/3).
\leqno (7.61) 
$$
Indeed, let $w\in A_1$ be given.
Recall that $|\pi(y)-y| \leq \varepsilon t$ by (7.39)
and $|\varphi_1(w)-w| \leq C_0\varepsilon t$ by (7.51),
so $w \i B(y,{5t \over 3}+ (C_0+1)\varepsilon t)$
because $\varphi_1(w) \in B(\pi(y),{5t \over 3})$. Then
$\psi(|w-y|) = 1$ by (7.47), so $\varphi_1(w) = h(w,1) = \pi(w)$ 
by (7.48) and (7.41). Also, $\pi(w) \in P \sm B(p,\tau)$ 
by (7.55), so altogether $\varphi_1(w) 
= \pi(w)\in P \cap B(\pi(y),{5t \over 3})\sm B(p,\tau)$.
This is the case when $g(\varphi_1(w)) = g(\pi(w))$ 
is the radial projection of $\pi(w)$ on $\d B(\pi(y),{5t \over 3})$,
as in (7.56).
But  $\varphi_2(w) = g(\varphi_1(w))$ by (7.58), so
$\varphi_2(w) \in P \cap \d B(\pi(y),{5t \over 3})$, as needed
for (7.61). 

We like (7.61) because it immediately implies that
$$
\H^d(\varphi_2(A_1)) = 0.
\leqno (7.62)
$$
Also,
$$
E^\ast \cap B(y,{5t \over 3}- (C_0+1)\varepsilon t) \i A_1 \cap  W_2
\leqno (7.63)
$$
because if $w\in E^\ast \cap B(y,{5t \over 3}- (C_0+1)\varepsilon t)$,
then $\varphi_1(w) \in B(\pi(y),{5t \over 3})$
(again by (7.51) and because $|\pi(y)-y| \leq \varepsilon t$),
hence $w\in A_1$; in addition $w\in B(\pi(y),{5t \over 3})$ and
$\varphi_2(w) \in \d B(\pi(y),{5t \over 3})$ by (7.61), so
$\varphi_2(w) \neq w$ and $w\i W_2$.

Next consider 
$$
A_2 = E^\ast \cap B(y,{5t \over 3}+ (C_0+2)\varepsilon t) \sm A_1.
\leqno (7.64)
$$
Notice that $A_2$ is fairly small, because it is contained
in $E^\ast \cap B(y,{5t \over 3}+ (C_0+2)\varepsilon t) 
\sm B(y,{5t \over 3} -(C_0+1)\varepsilon t)$
(by (7.63)), and in an $\varepsilon t$-neighborhood of $P$ by (7.39). 
So we can cover $A_2$ by less than $C \varepsilon^{-d+1}$ balls $B_l$ 
of radius $(C_0+10)\varepsilon t$, centered on the $(d-1)$-dimensional sphere
$P \cap \d B(y,{5t \over 3})$. Proposition~4.1 says that
$H^d(E \cap B_l) \leq C (C_0\varepsilon t)^d$ for each $l$
(recall that $C_0 \geq 1$), so
$$
\H^d(A_2) \leq C C_0^d \varepsilon t^d.
\leqno (7.65)
$$
But we mostly need to control $\varphi_2(A_2)$, so let us prove that
$$
\hbox{$\varphi_2$ is $2 C_0$-Lipschitz on $A_2$.}
\leqno (7.66)
$$
First we check that
$$
\varphi_2(w) = \varphi_1(w) = h(w,\psi(|w-y|))
\ \hbox{ for } w \in A_2.
\leqno (7.67)
$$
Indeed $\varphi_1(w) \notin B(\pi(y),{5t \over 3})$ 
since $w \notin A_1$, then $g(\varphi_1(w)) = \varphi_1(w)$ by (7.57),
and so $\varphi_2(w) = g(\varphi_1(w)) = \varphi_1(w)$ by (7.58).
The last identity comes from (7.48).

Now let $w, w' \i A_2$, be given,
and set $a = \psi(|w-y|)$ and $a' = \psi(|w'-y|)$, where 
$\psi$ is still as in (7.47) and (7.48).
Thus $|a'-a| \leq (\varepsilon t)^{-1}  |w'-w|$. Now
$$\eqalign{
|\varphi_2(w) - \varphi_2(w')| &= |h(w,a)-h(w',a')|
\cr&
\leq |h(w,a)-h(w,a')| + |h(w,a')-h(w',a')|
\cr&
\leq C_0 \varepsilon t |a'-a| + C_0 |w'-w| \leq 2C_0 |w'-w|
}\leqno (7.68)
$$
by (7.67), (7.42), and (7.43).
This proves (7.66).

\ms
Next we check that 
$$
\varphi_2(w) = w \ \hbox{ for }
w = E^\ast \sm B(y,{5t \over 3}+ (C_0+2)\varepsilon t)
\leqno (7.69)
$$
Indeed $\psi(|w-y|) = 0$ by (7.47), hence $\varphi_1(w) = w$ 
by (7.48) or (7.49), 
and so $\varphi_2(w) = g(\varphi_1(w)) = g(w)$ by (7.58).
But $w \i \R^n \sm B(y,{5t \over 3})$, so $g(w)=w$ by (7.57)
and as needed for (7.69).

By (7.69) and (7.64), 
$W_2 = \big\{ w\in E^\ast \cap B \, ; \, \varphi_2(w) \neq w \big\}
\i A_1 \cup A_2$, and 
$$
\H^d(\varphi_2(W_2)) \leq \H^d(\varphi_2(A_2)) \leq 2^d C_0^d \H^d(A_2)) 
\leq C C_0^{2d} \varepsilon t^d
\leqno (7.70)
$$
by (7.62), (7.66), and (7.65). On the other hand,
$$
\H^d(W_2) \geq \H^d(E^\ast \cap B(y,{5t \over 3}- (C_0+1)\varepsilon t))
\geq C^{-1} t^d
\leqno (7.71)
$$
by (7.63) and Proposition 4.1, and so (7.70) and (7.71) 
contradict (7.60) if $h$ and $\varepsilon$ are chosen small enough, 
depending on $M$, $n$, and $C_0$. So we were wrong to assume that 
there exists $p$ so that (7.52) holds, and this proves (7.46).

Now (7.45) follows from (7.46), because for $p\in P \cap B(y,3t/2)$,
we can find $w\in E^\ast \cap B(y,5t/3)$ such that $\pi(w) = p$,
and $|p-w| = |\pi(w)-w| \leq \varepsilon t$ because
$\dist(w,P) \leq \varepsilon t$ by (7.39).
Lemma 7.38 follows, and also Lemma 7.33 (see the comments
below the statement of Lemma 7.38).
\hfill $\square$ $\square$

\ms 
We are finally in position to gather the estimates on the various 
bad sets and resume our proof of Theorem 6.1. 
We start with a control on the
bad sets ${\cal B}_{E^\ast}(\varepsilon)$ of (7.3).

\ms\proclaim Lemma 7.72.
Let $E$, $x\in E^\ast$, and $r$ be as in Theorem 6.1.
Then for each $\varepsilon > 0$
$$
\int_{y\in E^\ast \cap B(x,r/4)}\int_{0 < t < r/10} 
{\bf 1}_{{\cal B}_{E^\ast}(\varepsilon)}(y,t) \, {d\H^d(y) dt \over t}
\leq C(\varepsilon) \, r^d,
\leqno (7.73)  
$$ 
where $C(\varepsilon)$ depends only on $n$, $M$, and $\varepsilon$.

\ms
Obviously it will be enough to prove (7.73) for 
${\cal B}_{E^\ast}(4\varepsilon)$ instead of 
${\cal B}_{E^\ast}(\varepsilon)$. Also, we may as well suppose 
that $\varepsilon$ is small (depending on $n$ and $M$), because 
${\cal B}_{E^\ast}(4\varepsilon)$ is larger when $\varepsilon$ is smaller
(see the definition (7.3)). 

Let  $x$, $r$, and $\varepsilon$ be as in
the statement, and set
$$
{\cal B} = {\cal B}(4\varepsilon, x, r) = 
\big\{ (y,t) \in {\cal B}_{E^\ast}(4\varepsilon)
\, ; \, y\in E^\ast \cap B(x,r/4) \hbox{ and } 0 < t < r/10 \big\};
\leqno (7.74) 
$$
then (7.73) is the same as
$$
\int\int_{(y,t) \in {\cal B}}  \, {d\H^d(y) dt \over t}
\leq C(\varepsilon) r^d.
\leqno (7.75)  
$$
Clearly ${\cal B} \i {\cal A}$, with
$\cal A$ as in (7.18), because $t \leq r < \Min(r_0,\delta)$
and $B(y,2t) \i B(x,2r) \i B_0$. 
The set ${\cal B} \cap {\cal B}_1$
is taken care of by Lemma 7.19, and similarly
${\cal B} \cap {\cal B}_2$ is controlled by Lemma 7.23. So we just need
to show that
$$
\int\int_{(y,t) \in {\cal B}'}  \, {d\H^d(y) dt \over t}
\leq C(\varepsilon) r^d,
\leqno (7.76)  
$$
with ${\cal B}' = {\cal B} \sm ({\cal B}_1 \cup {\cal B}_2)$. 
Notice that if $(y,t) \in {\cal B}'$, then 
$b\beta_{E^\ast}(y,t) > 4\varepsilon$ because 
$(y,t)  \in {\cal B}_{E^\ast}(4\varepsilon)$ 
(see the definition (7.3)), and then Lemma 7.33 says that 
$B(y,2t)$ does not meet $F_2$. That is,
$$
{\cal B}' \i \big\{(y,t)\in {\cal A} \, ; \, y\in E^\ast \cap B(x,r/4),
\, t < r/10, \hbox{ and $B(y,2t)$ does not meet }F_2 \big\}.
\leqno (7.77) 
$$

At this point, we want to cut ${\cal B}'$ 
into smaller sets for which we can apply Lemma 7.8, and
for this a covering of $E^\ast \cap B(x,r/4) \sm F_2$
will be useful. For $z \in E^\ast \cap B(x,r/4) \sm F_2$, set
$$
d(z) = \Min(r,\dist(z,F_2)) > 0
\leqno (7.78) 
$$
and $B_z = B(y,{d(z) \over 100})$. Then select a maximal set
$Z \i E^\ast \cap B(x,r/4) \sm F_2$ such that
$$
\hbox{the $B_z$, $z\in Z$, are disjoint.}
\leqno (7.79) 
$$
For each $y\in E^\ast \cap B(x,r/4) \sm F_2$,
we select $z = z(y) \in Z$ so that $B_z$ meets $B_y$;
such a $z$ exists by maximality of $Z$, and it is
easy to select $z(y)$ in a measurable way, because
$Z$ is at most countable. Then cut ${\cal B}'$ as 
$$
{\cal B}' = \bigcup_{z\in Z} {\cal B}'(z),
\leqno (7.80) 
$$
where
$$
{\cal B}'(z) = \big\{ (y,t)\in {\cal B}' \, ; \, z(y) = z \big\}.
\leqno (7.81) 
$$

Fix $z\in Z$ for the moment. We want to apply Lemma 7.8
to the quasiminimal set $E^\ast$ and the
pair $(z, d(z)/2)$, so let us check the hypotheses.
We know from Proposition 3.3 that $E^\ast \in GSAQ(B_0,M,\delta,h)$,
just like $E$, but with $E^\ast$ (6.42) will be easier to check.
First recall that $z \in E^\ast \cap B(x,r/4)$ and 
$d(z) \leq r$; hence the first assumptions that 
$z\in E^\ast \cap B_0$, $0 < d(z)/2 < \Min(r_0,\delta)$,
and $B(z,d(z)) \i B_0$ follow from the similar assumptions
for $(x,r)$. Now we check the main assumption
(6.42). Let $j$ be such that $L_j$ meets $B(z,d(z)/2)$;
we want to show that $E^\ast \cap B(z,d(z)/2) \i L_j$.

Recall from (7.78) that $\dist(z,F_2) \geq d(z)$, where 
$F_2$ denotes the union of all the $d$-dimensional 
faces of cubes from our dyadic grid (see Lemma 7.23).
This means that the faces of $L_j$ that meet $B(z,d(z)/2)$
are at least $(d+1)$-dimensional. We know that there is at least
one face like this, because we assume that $L_j$ meets $B(z,d(z)/2)$.
But $B(z,d(z)/2) \i B(x,r)$ (because $d(z) \leq r$); 
then the main assumption (6.2) says that $E^\ast \cap B(x,r) \i L_j$, 
which is enough for (6.42). 
So we may apply Lemma 7.8, and we get that 
$$
\int_{y\in E^\ast \cap B(z,d(z)/16)}\int_{0 < t < d(z)/16} 
{\bf 1}_{{\cal B}_{E^\ast}(4\varepsilon)}(y,t) \, {d\H^d(y) dt \over t}
\leq C(4\varepsilon) d(z)^d.
\leqno (7.82)  
$$

Return to ${\cal B}'(z)$. If $(y,t) \in {\cal B}'(z)$, then 
$100|z-y| < d(z) + d(y)$ because $B_z$ meets $B_y$ when $z=z(y)$; since
$d(y) \leq d(z) + |z-y|$ by (7.78), we get that
$100|z-y| < 2 d(z) + |z-y|$, and hence $|z-y| < d(z)/49$
and also $d(y) \leq d(z) + |z-y| \leq {50 \over 49}\, d(z)$.

If in addition $t < d(z)/16$, $(y,t)$ lies in the domain
of integration of (7.82) (see (7.74) and the definition
of ${\cal B}'$ below (7.76)) and it will be taken care of by (7.82). 
Otherwise, observe that $B(y,2t)$ does not meet $F_2$ (by (7.77)).
If $d(y) = \dist(y,F_2)$, this shows that 
$t \leq d(y)/2 \leq d(z)/3$. Otherwise $d(y) = r$,
so $d(z) \geq {49d(y) \over 50} \geq {49r \over 50}$.
In this case too $t \leq d(z)/3$, because $t \leq r/10$ when 
$(y,t) \in {\cal B}'$. Altogether
$$\leqalignno{
\int\int_{(y,t) \in {\cal B}'(z)} \, {d\H^d(y) dt \over t}
&\leq C(4\varepsilon) d(z)^d + 
\int_{y\in E^\ast \cap B(z,d(z)/49)} \int_{d(z)/16 < t \leq d(z)/3}
\, {d\H^d(y) dt \over t}
\cr&
\leq C(4\varepsilon) d(z)^d + 
\ln(16/3) \, \H^d(E^\ast \cap B(z,d(z)/49)) \leq C'(\varepsilon) d(z)^d
& (7.83) 
}
$$
by Proposition 4.1 (recall that $B(z,d(z)) \i B(x,2r)$, so 
$B(z,d(z)/49)$ is not too large).

We now use (7.80) and sum over $z \,$:
$$\eqalign{
\int\int_{(y,t) \in {\cal B}'} \, {d\H^d(y) dt \over t}
&= \sum_{z\in Z} \int\int_{(y,t) \in {\cal B}'(z)} \, {d\H^d(y) dt \over t}
\leq  C'(\varepsilon) \sum_{z\in Z} d(z)^d
\cr&
\leq C C'(\varepsilon) \sum_{z\in Z} \H^d(E^\ast \cap B_z)
\leq C C'(\varepsilon) r^d
}\leqno (7.84)  
$$
by Proposition 4.1, and because the $B_z$ are disjoint (by (7.79))
and contained in $B(x,r)$. This is (7.76), and Lemma 7.72 follows.
\qed

\ms
Next we use Theorem 7.7 to prove the analogue of Proposition 6.41
under the (weaker) assumptions of Theorem 6.1. A more direct approach
to Theorem 6.1 is also possible, using the weak geometric lemma and
big projections now, but the next proposition is really a logical
consequence of Lemma 7.72.

\ms\proclaim Proposition 7.85.
Let $E$, $x\in E^\ast$, and $r$ be as in Theorem 6.1
(again with $h$ small enough, depending on $M$ and $n$).
Then there is a closed set $G_0 \i E^\ast \cap B(x,r)$ 
and a mapping $\phi : G_0 \to \R^d$ such that 
$$
\H^d(G_0) \geq \theta r^d
\ \hbox{ and } \ 
C_M^{-1} |y-z| \leq |\phi(y)-\phi(z)| \leq C_M |y-z|
\hbox{ for } y,z \in G_0,
\leqno (7.86)
$$
where $\theta > 0$ and $C_M$ depend only on $M$ and $n$.

\ms
In other words, we can find, in $E^\ast \cap B(x,r)$,
a big piece of bilipschitz image of a subset of $\R^d$. 
Proposition 7.85 really follows from Lemma 7.72 and the 
proof of Theorem~7.7, but we need a localization argument 
so that we can use the statement of Theorem 7.7 as it is.

Let $E$, $x\in E^\ast$, and $r$ be as in Proposition 7.85
or Theorem 6.1; let use again Proposition 7.6 in [DS4], 
as we did for the proof of Lemma 7.8 but applied to a slightly
smaller radius, to find a bounded Ahlfors-regular set $F$ such that
$$
E^\ast \cap B(x,r/16) \i F \i E^\ast \cap B(x,r/8).
\leqno (7.87)
$$
Since we want an unbounded Ahlfors-regular set,
we consider the set $F' = F \cup H$, where $H$
is a $d$-plane such that $\dist(x,H) = r$. We want 
to use Theorem 7.7, so let us prove that for 
every small $\varepsilon > 0$, there is a constant 
$C(\varepsilon)$, that depends only on $n$, $M$, and $\varepsilon$,
such that $F' \i BWGL(\varepsilon,C(\varepsilon))$.

So we let $(x_1,r_1) \i F' \times (0,+\infty)$ be given;
we want to prove that
$$
\int\int_{(y,t) \in {\cal B}} \, {d\H^d(y) dt \over t} \leq 
C(\varepsilon) r_1^d ,
\leqno (7.88)
$$
where ${\cal B} = \big\{ (y,t) \in (F' \cap B(x_1,r_1))\times (0,r_1) \, ; \,
b\beta_{F'}(y,t) > \varepsilon \big\}$. We shall need to cut
${\cal B}$ into many pieces, to control various interfaces,
but let us start with the most interesting case when 
$x_1 \in B(x,r/8)$ and $r_1 \leq r/40$. Then the main piece is
$$
{\cal B}_1 = \big\{ (y,t) \in {\cal B} \, ; \, 
\dist(w,F) \leq {\varepsilon t \over 2} \hbox{ for } w\in E^\ast \cap B(y,2t) 
\big\}.
\leqno (7.89)
$$
We claim that
$b\beta_{E^\ast}(y,2t) > \varepsilon/4$ for $(y,t) \in {\cal B}_1$.
Notice that $y\in E^\ast$ because $F \i E^\ast$, so at least 
$b\beta_{E^\ast}(y,2t)$ was officially defined in (7.2). If 
$b\beta_{E^\ast}(y,2t) \leq \varepsilon/4$, there is a $d$-plane $P$ 
through $y$ such that
$$
\sup_{w\in E^\ast\cap B(y,2t)}\dist(w,P) 
+\sup_{w\in P\cap B(y,2t)}\dist(y,E^\ast)
\leq \varepsilon t/2;
\leqno (7.90)
$$
but
$$
\sup_{w\in F'\cap B(y,t)}\dist(w,P) \leq \sup_{w\in E^\ast\cap 
B(y,2t)}\dist(w,P)
\leqno (7.91)
$$
because $F'\cap B(y,t) = F\cap B(y,t) \i E^\ast \cap B(y,t)$, and 
$$
\sup_{w\in P\cap B(y,t)}\dist(y,F') \leq \sup_{w\in P\cap B(y,t)}\dist(y,F)
\leq \sup_{w\in P\cap B(y,t)}\dist(y,E^\ast) + {\varepsilon t \over 2}
\leqno (7.92)
$$
because $F \i F'$ and by definition of ${\cal B}_1$, which contradicts 
the fact that $(y,t) \in {\cal B}$. Now
$$
\int\int_{(y,t) \in {\cal B}_1} \, {d\H^d(y) dt \over t} 
\leq \int_{y\in F\cap B(x_1,r_1)} \int_{0 < r < r_1} 
{\bf 1}_{{\cal B}_{E\ast}(\varepsilon/4)}(y,2t) \, {d\H^d(y) dt \over t}
\leq C(\varepsilon) r_1^d,
\leqno (7.93)
$$
where the last inequality comes from Lemma 7.72, applied to the
pair $(x_1, 20r_1)$. 

Then we need to control ${\cal B} \sm {\cal B}_1$.
To $(y,t)\in {\cal B} \sm {\cal B}_1$ we associate 
$z = z(y,t) \in E^\ast \cap B(y,2t)$ such that 
$\dist(z,F) \geq {\varepsilon t \over 3}$, and the set
$A(y,t) = E^\ast \cap B(z,{\varepsilon t \over 6})$. 
We can choose $z(y,t)$ and $A(y,t)$ in such a way that
$(y,t,w) \to {\bf 1}_{A(y,t)}(w)$ is measurable,
for instance by choosing the first available
$z$ from a sufficiently dense countable set in $E^\ast$.
Notice that if $w\in A(y,t)$, then $|w-y| \leq 3t$ and
$t \leq 6\varepsilon^{-1} \dist(w,F)$. In addition,
$t \geq |w-y|/3 \geq \dist(w,F)/3$ because $y\in F$,
so that 
$$
t \in T(w) = (0,r_1] \cap [\dist(w,F)/3,6\varepsilon^{-1} 
\dist(w,F)]
\leqno (7.94)
$$ 
(recall that $0 < t \leq r_1$ for $(y,t)\in {\cal B}$).
Finally $w \in B(x_1,4r_1)$ because $y \in B(x_1,r_1)$.
Now multiple uses of Proposition 4.1 yield
$$\eqalign{
&\hskip-0.4cm
\int\int_{(y,t) \in {\cal B} \sm {\cal B}_1} \, {d\H^d(y) dt \over t}
\leq C
\int\int_{(y,t) \in {\cal B}\sm {\cal B}_1}
\int_{w\in A(y,t)} \, (\varepsilon t)^{-d} \, {d\H^d(y) d\H^d(w)dt \over t}
\cr&
\leq C(\varepsilon) \int_{w\in E^\ast \cap B(x_1,4r_1)}
\int_{t \in T(w)} 
\int_{y\in F\cap B(w,3t)}
\, {\H^d(y) d\H^d(w)dt \over t^{d+1}}
\cr&
\leq C(\varepsilon) \int_{w\in E^\ast \cap B(x_1,4r_1)}
\int_{t \in T(w)} \, {d\H^d(w)dt \over t}
\cr&
\leq C(\varepsilon) \ln(18/\varepsilon) \, 
\int_{w\in E^\ast \cap B(x_1,4r_1)}  d\H^d(w)
\leq C(\varepsilon) r_1^d,
}\leqno (7.95)
$$
as needed for (7.88).

This proves (7.88) when $x_1 \in B(x,r/8)$ and $r_1 \leq r/40$.
When $x_1 \in H$ and $r_1 \leq r/40$, $F'$ coincides with
$H$ on $B(x_1,2r_1)$, and ${\cal B} = \emptyset$.
This settles the case when $r_1 \leq r/40$ because
$F \i B(x,r/8)$ by (7.87). So assume now that $r_1 > r/40$.

Note that $b\beta_{F'}(y,t) \leq \varepsilon$ for $y\in F'$
and $t \geq 10 \varepsilon^{-1} r$, simply because we can use
$H$ in the definition of $b\beta_{F'}(y,t)$, and by (7.87).
Thus ${\cal B} = {\cal B}_2 \cup {\cal B}_3$, where 
${\cal B}_2 = \big\{ (y,t) \in {\cal B} \, ; \, 0 < t < r/40 \big\}$
and ${\cal B}_3 = \big\{ (y,t) \in {\cal B} \, ; \,
r/40 \leq t < 10 \varepsilon^{-1} r \big\}$.
Since $b\beta_{F'}(y,t) = 0$ when $y\in H$ and $0 < t < r/40$,
$$\eqalign{
\int\int_{(y,t) \in {\cal B}_2} \, {d\H^d(y) dt \over t} 
&\leq \int_{y\in F} \int_{0 < t < r/40} 
{\bf 1}_{{\cal B}_{F}(\varepsilon)}(y,t) \, {d\H^d(y) dt \over t} 
\cr&
\leq C(\varepsilon) r^d  \leq C(\varepsilon) 40^d r_1^d,
}\leqno (7.96)
$$
where the main inequality comes from the case of $x_1 = x$ and
$r_1 = r/40$, which was treated before. And
$$\eqalign{
\int\int_{(y,t) \in {\cal B}_3} \, {d\H^d(y) dt \over t} 
&= \int_{y\in F' \cap B(x_1,r_1)} \int_{r/40 \leq t < 10 \varepsilon^{-1} r} 
\, {d\H^d(y) dt \over t} 
\cr&\leq C(\varepsilon) \H^d(F' \cap B(x_1,r_1))
\leq C(\varepsilon) r_1^d.
}\leqno (7.97)
$$

This completes our proof of (7.88), from which we deduce that
$F' \i BWGL(\varepsilon,C(\varepsilon))$ and, by Theorem 7.7,
that $F' \in BPBI(\theta,C_M)$ for some choice of $\theta > 0$
and $C_M \geq 1$ that depend only on $n$ and $M$.
The reader should not be surprised not to see $\varepsilon$ any more;
recall that only one value of $\varepsilon$ is needed, that depends 
on $M$ and $n$ through the Ahlfors-regularity constants for $F'$.
We apply the definition (7.6) of $BPBI(\theta,C_M)$ to the pair
$(x,r)$ and get a set $G_0 \i F \cap B(x,r) \i E^\ast \cap B(x,r)$
(by (7.87)); then $G_0$ satisfies (7.86), and this completes our 
proof of Proposition 7.85.
\qed

\ms\noindent
{\bf 8. Big projections and big pieces of Lipschitz graphs} 
\ms

We shall complete the proof of Theorem 6.1 in this section,
with the help of yet another theorem of uniform rectifiability
that we state now, in its initial context of unbounded Ahlfors-regular sets.
First we need to define the weak geometric lemma, big projections, and 
big pieces of Lipschitz graphs for such sets.

Let $F \i \R^n$ be an (unbounded) Ahlfors-regular set of dimension $d$;
this means that $F$ is closed and (7.4) holds. Let $\beta_F(x,r)$
be the P. Jones number defined in (7.1), and set
$$
{\cal M}_F(\varepsilon) = 
\big\{ (x,r) \in F \times (0,+\infty) \, ; \,
\beta(x,r) > \varepsilon \big\}
\leqno (8.1)
$$
for $0 < \varepsilon < 1$. We say that 
$F \in WGL(\varepsilon, C(\varepsilon))$ 
(or that $F$ satisfies a weak geometric lemma, with the
constants $\varepsilon$ and $C(\varepsilon)$) when 
$$
\int_{y\in F \cap B(x,r)}\int_{0 < t < r} 
{\bf 1}_{{\cal M}_F(\varepsilon)}(y,t) \, {d\H^d(y) dt \over t}
\leq C(\varepsilon) r^d
\leqno (8.2)
$$
for $x\in F$ and $0 < r < +\infty$.
This is the same thing as the bilateral weak geometric lemma,
but with the smaller functions $\beta_F(x,r)$ that only check
whether every point of $F \cap B(x,r)$ lies near a plane.

We say that $F \in BP(\alpha)$ (for big projections) if for
every choice of $x\in F$ and $0 < r < +\infty$, we can find a 
$d$-plane $P$ such that
$$
\H^d(\pi(F\cap B(x,r))) \geq \alpha r^d,
\leqno (8.3)
$$
where $\pi$ denotes the orthogonal projection on $P$.

Finally, we say that $F\in BPLG(\theta,A)$ (for big pieces Lipschitz graphs) 
when for all $x\in F$ and $r > 0$, we can find a $d$-dimensional 
$A$-Lipschitz graph $\Gamma \i \R^n$ (which means, the graph of some 
Lipschitz function which is defined on a $d$-plane $P$, with values in $P^\perp$, 
and with a Lipschitz constant at most $A$) such that 
$$
\H^d(F \cap \Gamma \cap B(x,r))) \geq \theta r^d.
\leqno (8.4)
$$
This property is slightly, but in general strictly stronger than uniform 
rectifiability (or equivalently $BPBI$). The following is 
Theorem 1.14 on page 857 of [DS2]. 

\ms\proclaim Theorem 8.5.
Let $F \i \R^n$ be a closed unbounded Ahlfors-regular set 
of dimension $d$. If $F \in BP(\alpha)$ for some $\alpha > 0$
and $F \in WGL(\varepsilon, C(\varepsilon))$ 
for some small enough $\varepsilon$ (depending only on 
$n$, $\alpha$, and the constant $C_0$ in (7.4)), then
$F\in BPLG(\theta,A)$, where $\theta > 0$ and
$A \geq 0$ depend only on $n$, $\alpha$, $\varepsilon$, 
and $C(\varepsilon)$.

\ms\noindent{\bf Proof of Theorem 6.1.}
Let $E$, $x$, and $r$ be as in Theorem 6.1, and let $F$ and 
$F' = F \cup H$ be the Ahlfors-regular sets that we 
already used for Proposition 7.85. Thus 
$E^\ast \cap B(x,r/16) \i F \i E^\ast \cap B(x,r/8)$ 
as in (7.87), and $H$ is a $d$-plane such that $\dist(x,H) = r$.

We want to say that $F' \in BPLG(\theta,A)$ for some choice
of $\theta$ and $A$, so let us check that
$F' \in WGL(\varepsilon, C(\varepsilon)) \cap BP(\alpha)$.
In fact, we already know that $F'\in WGL(\varepsilon, C(\varepsilon))$
with $\varepsilon$ as small as we want, because we even checked
in (7.88) that $F'\in BWGL(\varepsilon, C(\varepsilon))$, which
is obviously stronger. So we just need to check that 
$F' \in BP(\alpha)$ for some $\alpha >0$, that depends only
on $n$ and $M$.

So let $x_1 \i F'$ and $r_1 > 0$ be given; we want to find
$P = P(x_1,r_1)$ such that, as in (8.3),
$$
\H^d(\pi(F'\cap B(x_1,r_1))) \geq \alpha r_1^d,
\leqno (8.6)
$$
where $\pi$ denotes the orthogonal projection on $P$.
The case when $x_1 \in H$ or $x_1\in F$ but $r_1 \geq 4r$
is trivial, because we can take $P = H$, so let us assume
that
$$
x_1 \in F \ \hbox{ and } 0 < r_1 < 4 r.
\leqno (8.7)
$$
The idea is to find a pair $(x_2,r_2)$ to which
we can apply Lemma 7.38, and deduce (8.6) from (7.46).
Let us state what we need.

\ms\proclaim Lemma 8.8.
Theorem 6.1 will follow as soon as we prove the following.
For each $\varepsilon > 0$, there is a small constant $c_\varepsilon$,
which depend only on $M$, $n$, and $\varepsilon$, such that
if $x\in E^\ast$, $r$, $F$, $x_1 \in F$, and $r_1 < 4r$ are as above, 
then we can find $x_2\in F \cap B(x_1,r_1/10)$, 
$r_2 \in [c_\varepsilon r_1, r_1/10]$, 
and a $d$-plane $P$ such that the assumptions of Lemma 7.38 
(relative to the pair $(x_2,r_2)$) are
satisfied, and
$$
H^d(E^\ast \cap B(x_2,2r_2) \sm F) \leq \varepsilon r_2^d.
\leqno (8.9)
$$

\ms
To prove the lemma, we assume the existence of $(x_2,r_2)$ and check (8.6). 
Since Lemma~7.38 applies (if $h$ is small enough, depending on $M$ and $n$), 
we get a plane $P$ such that (7.46) holds, and then
$$\eqalign{
\H^d(\pi(F'\cap B(x_1,r_1))) 
&\geq \H^d(\pi(F\cap B(x_2,2r_2)))
\geq \H^d(\pi(E^\ast\cap B(x_2,2r_2))) - \varepsilon r_2^d
\cr& \geq \H^d(P\cap B(x_2,3r_2/2))) - \varepsilon r_2^d
\geq C^{-1} r_2^d \geq C(\varepsilon) r_1^d
}\leqno (8.10)
$$
by (8.9), if $\varepsilon$ is small enough, and 
because $r_2 \geq c_\varepsilon r_1$. Then (8.6) holds, $F'$ has big 
projections, and Theorem 8.5 says that it contains big pieces of
Lipschitz graphs too. We apply the definition (8.4) to the
ball $B(x,r/16)$, and the Lipschitz graph $\Gamma$ such that
$\H^d(F' \cap \Gamma \cap B(x,r/16))) \geq 16^{-d} \theta r^d$
also works for Theorem 6.1, except that we need to divide $\theta$
by $16^d$.
\qed

\ms
So we are left with the assumption of Lemma 8.8 to check. 
That is, we are given $x_1 \in F$ and $r_1 < 4r$ as above, and
we want to find a pair $(x_2,r_2)$. We shall produce first an 
intermediate pair $(y_0,t_0)$, with better regularity 
properties than $(x_1,r_1)$. These properties will be stated
in terms of a very small $\varepsilon_0 < \varepsilon$, which
will be chosen later, depending on $n$, $M$, and $\varepsilon$.

Recall that $F_1$ denotes the union of all the faces of 
dimension at most $d-1$ of cubes from our usual dyadic net
(see Lemma 7.19), and $F_2$ is the union of all the $d$-dimensional 
faces of cubes from that net, as in Lemma 7.23.

\ms\proclaim Lemma 8.11. There exists a constant
$c(\varepsilon_0) > 0$, that depends on $n$, $M$, 
and $\varepsilon_0$, such that for $(x_1,r_1)$ as above,
we can find $(y_0,t_0)$ with the following properties:
$$
y_0 \in F \cap B(x_1,r_1/200)
\ \hbox{ and } c(\varepsilon_0) r_1 \leq t_0 \leq r_1/100,
\leqno (8.12)
$$
$$
\dist(y_0,F_1) \geq 10 t_0,
\leqno (8.13)
$$
$$
\hbox{if $B(y_0,2t_0)$ meets $F_2$, then }
\dist(w,F_2) \leq \varepsilon_0 t_0
\ \hbox{ for } w\in E^\ast \cap B(y_0,2t_0),
\leqno (8.14)
$$
$$
b\beta_{E^\ast}(y_0,t_0) \leq \varepsilon_0,
\leqno (8.15)
$$
and 
$$
\H^d(E^\ast \cap B(y_0,2t_0)\sm F) \leq \varepsilon_0 t_0^d.
\leqno (8.16)
$$

\ms
The lemma will follow from more Carleson packing estimates. Let
$$
{\cal A}_0 = \big\{ (y,t) \, ; \, y \in F \cap B(x_1,r_1/200)
\ \hbox{ and } 0 < t \leq r_1/200 \big\}
\leqno (8.17)
$$
be a little smaller than the set of pairs that we want to pick from; 
notice that it is contained in the set $\cal A$ of (7.18).
Let ${\cal B}_1$ and ${\cal B}_2$ be as in (7.20) and (7.24),
except that we replace $\varepsilon$ with $\varepsilon_0$ in (7.24);
we know from Lemmas 7.19 and 7.23 that
$$
\int\int_{(y,t) \in {\cal A}_0 \sm {\cal B}_1 \cup {\cal B}_2} 
\, {d\H^d(y) dt \over t} \leq C(\varepsilon_0) r_1^d.
\leqno (8.18)
$$
Next set
$$
{\cal B}_3 = \big\{ (y,t) \in {\cal A}_0 \, ; \, 
b\beta_{E^\ast}(y,t) > \varepsilon_0 \big\}.
\leqno (8.19)
$$
Since $x_1 \in F \i E^\ast \cap B(x,r/8)$ and $0 < r_1 < 4 r$
by (8.7), the pair $(x_1,r_1/8)$ satisfies the assumptions of
Lemma 7.72 (in particular because (6.2) is hereditary), so we get 
from that lemma that
$$
\int\int_{(y,t) \in {\cal B}_3} 
\, {d\H^d(y) dt \over t} 
\leq C(\varepsilon_0) \, r_1^d.
\leqno (8.20)
$$
Next we want to check that
$$
\int\int_{(y,t) \in {\cal B}_4} 
\, {d\H^d(y) dt \over t} 
\leq C(\varepsilon_0) \, r_1^d
\leqno (8.21)
$$
for the set
$$
{\cal B}_4 = \big\{ (y,t) \in {\cal A}_0 \, ; \, 
\H^d(E^\ast \cap B(y,2t)\sm F) > \varepsilon_0 t^d \big\}.
\leqno (8.22)
$$
Let us first find, for each pair $(y,t) \in {\cal B}_4$, a point
$z(y,t)$ such that
$$
z(y,t) \in E^\ast \cap B(y,2t) \ \hbox{ and } \ 
\dist(z(y,t),F) \geq 2\gamma t,
\leqno (8.23)
$$
where the very small constant $\gamma > 0$ will chosen soon, depending
on $\varepsilon_0$ as well.

For this we need to know that $E^\ast\sm F$ is not too dispersed,
and we are going to use again the dyadic patchwork that is provided
for us by Proposition 7.6 of [DS4], 
and that we rapidly described near (6.7) and (7.10).
This time, we shall need the other part of the small relative boundary 
condition ($\ast$7.10) (that is, (7.10) in [DS4]), 
namely the fact that for every cube $Q$ of our cubical patchwork
for $F$, 
$$\H^d \big(\big\{ w\in E^\ast \sm Q \, ; \, \dist(x,Q) \leq \tau 
\diam(Q) \big\}\big) \leq C \tau^{1/C} \diam(Q)^d
\leqno (8.24)
$$
for every (small) $\tau > 0$, with a constant $C$ that depend only
on $M$ and $n$.

Return to the pair $(y,t) \in {\cal B}_4$, and denote by $Q_i$, $i\in I$,
the collection of maximal cubes of the patchwork such that
$Q_i \cap B(y,3t) \neq \emptyset$ and $\diam(Q_i) \leq t$.
These sets are disjoint (by ($\ast$7.3) and maximality),
and since $\diam(Q_i) \geq t/C$ (by maximality and the size
property ($\ast$7.2)), we get that $\H^d(Q_i) \geq C^{-1} t^d$
(by the second part of ($\ast$7.2)) for $i\in I$.
But $H^d(\cup_{i\in I} Q_i) \leq \H^d(F \cap B(y,4t)) \leq Ct^d$
by Proposition 4.1, and hence $I$ has at most $C$ elements.

Suppose we cannot find $z(y,t)$ such that (8.23) holds.
Then each $z\in E^\ast \cap B(y,2t) \sm F$ lies within
$2\gamma t$ of $F \cap B(y,3t)$, so $\dist(z,Q_i) \leq 2\gamma t$
for some $i\in I$. Recall that $\diam(Q_i) \geq C^{-1} t$,
so $z$ lies in the set of (8.24), with $\tau = 2 C \gamma t$.
We get that
$$
\H^d(E^\ast \cap B(y,2t) \sm F)
\leq \sum_{i\in I} C \tau^{1/C} \diam(Q)^d
\leq C' \tau^{1/C} t^d < \varepsilon_0 t^d
\leqno (8.25)
$$
by (8.24), and if $\tau$ is chosen small enough, 
depending on $\varepsilon_0$. [And this is how we choose $\gamma$.]
This contradiction with the definition (8.22) shows that we can find
$z(y,t)$ as in (8.23). We may now prove the Carleson measure
estimate (8.21) in the usual way. Denote by $A$ the right-hand side 
of (8.21); then
$$\eqalign{
A  &\leq C \int\int_{(y,t) \in {\cal B}_4} 
(\gamma t)^{-d} \int_{w\in E^\ast \cap B(z(y,t),\gamma t)} 
\, {d\H^d(y)d\H^d(w) dt \over t}
}\leqno (8.26)
$$
because $\H^d(E^\ast \cap B(z(y,t),\gamma t) 
\geq C^{-1} (\gamma t)^{d}$ by Proposition 4.1.
Then apply Fubini. Notice that $w\in E^\ast \cap B(x_1,r_1/10)$
because $z(y,t) \in B(y,2t)$, $y \in F \cap B(x_1,r_1/200)$,
and $t \leq r_1/200$ by (8.17). Next
$2\gamma t \leq \dist(z(y,t),F) \leq \dist(z(y,t),y) \leq 2t$ 
by (8.23) and because $y\in F$, so $\gamma t \leq \dist(w,F) \leq 3t$
and $t \in T(w) = (0,r_1/200] \cap [\dist(w,F)/3,\dist(w,F)/\gamma]$.
Finally,
$|y-w| \leq |y-z(y,t)|+ \gamma t \leq 3t$, and (8.26) yields
$$
\eqalign{
A &\leq C \gamma^{-d} \int_{w \in E^\ast \cap B(x_1,r_1/10)} 
\int_{t \in T(w)} t^{-d} \int_{y\in F \cap B(w,3t)} \, {d\H^d(w) dt d\H^d(y) \over t}
\cr&
\leq C \gamma^{-d} \int_{w \in E^\ast \cap B(x_1,r_1/10)}  
\int_{t \in T(w)} \, {d\H^d(w) dt \over t}
\cr&
\leq C \gamma^{-d} \, \ln(3/\gamma) \, \H^d(E^\ast \cap B(x_1,r_1/10))
\leq C(\varepsilon_0) \, r_1^d
}\leqno (8.27)
$$
by Proposition 4.1 again. This proves (8.21). 

Let us finally use our Carleson estimates to find a pair $(y_0,t_0)$
such that
$$
(y_0,t_0) \i {\cal A}_0 \sm \big[ {\cal B}_1 \cup {\cal B}_2 \cup
{\cal B}_3 \cup {\cal B}_4 \big]
\ \hbox{ and } \ 
t_0 \geq c(\varepsilon_0) r_1
\leqno (8.28)
$$
where as usual $c(\varepsilon_0) > 0$ depends on $M$, $n$, and 
$\varepsilon_0$. Observe that
$$\eqalign{
\int\int_{(y,t) \in {\cal A}_0 \, ; \, 
t_0 \geq c(\varepsilon_0) r_1} \, {d\H^d(y) dt \over t} 
&= \H^d(F \cap B(x_1,r_1/200)) 
\int_{c(\varepsilon_0) r_1 \leq t \leq r_1/200} \, {dt \over t}
\cr&
\geq C^{-1} r_1^d \ln(200/c(\varepsilon_0))
}\leqno (8.29)
$$
by (8.17) and Proposition 4.1. If we choose $c(\varepsilon)$ small 
enough, the right-hand side of (8.23) is larger than the sum of the
right-hand sides of (8.18), (8.20), and (8.21), and then the domain
of integration in (8.29) is not contained in 
${\cal B}_1 \cup {\cal B}_2 \cup {\cal B}_3 \cup {\cal B}_4$, which 
means that we can choose $(y_0,t_0)$ as in (8.28).
The pair $(y_0,t_0)$ satisfies (8.12) by (8.17) and (8.28),
and (8.13)-(8.16) by definition of the ${\cal B}_j$; Lemma 8.11
follows.
\qed
 
\ms
We now use the new pair $(y_0,t_0)$ to find the pair $(x_2,r_2)$
demanded by Lemma 8.8.

A first possibility is that $B(y_0,2t_0)$ meets $F_2$.
Then $\dist(w,F_2) \leq 4\varepsilon_0 t_0$
for $w\in E^\ast \cap B(y_0,2t_0)$ (almost as in (7.34)), and 
the proof of Lemma 7.33 (both before and just after the statement
of Lemma 7.38) applies, and says that we can apply Lemma 7.38 to the pair 
$(y,t) = (y_0,t_0)$ itself, with $C_0 = 1$ and
if we took $4\varepsilon_0 \leq \varepsilon$.
In this case we may take $(x_2,r_2) = (y_0,t_0)$ in the statement
of Lemma 8.8, and (8.9) holds by (8.16).

So we may assume that
$$
B(y_0,2t_0) \cap F_2 = \emptyset.
\leqno (8.30)
$$
Set 
$$
J = \big\{ j \in [0,j_{ max}] \, ; \, L_j \cap B(y_0,2t_0) \neq 
\emptyset \big\} \ \hbox{ and } L = \cap_{j\in J} L_j.
\leqno (8.31)
$$ 
Notice that $J \neq \emptyset$ because $E \i L_0 = \Omega_0$.
Let $j \in J$, and let $A$ be any face of our grid that is contained 
in $L_j$ and meets $B(y_0,2t_0)$. Such a face exists by definition of
$J$, and by (8.30) $A$ is more than $d$-dimensional. Since
$L_j \cap B(x,r) \supset A \cap B(y_0,2t_0) \neq \emptyset$,
(6.2) says that $E^\ast \cap B(y_0,2t_0) \i E^\ast \cap B(x,r) \i L_j$.
Hence $E^\ast \cap B(y_0,2t_0) \i L_j$ for $j\in J$. That is,
$$
E^\ast \cap B(y_0,2t_0) \i L.
\leqno (8.32)
$$
Since there are at most $C$ faces contained in $L$ and that meet $B(y_0,2t_0)$,
we can find such a face $A \i L$, so that
$$
\H^d(E^\ast \cap B(y_0,t_0/2) \cap A) 
\geq C^{-1} \H^d(E^\ast \cap B(y_0,t_0/2))
\geq C^{-1} t_0^d
\leqno (8.33)
$$
by Proposition 4.1. We claim that we can find $d+1$ points
$$
w_0, \ldots, w_d \in E^\ast \cap B(y_0,t_0/2) \cap A, 
\leqno (8.34)
$$
so that
for $1 \leq l \leq d$,
$$
\dist(w_l, P(w_0, \ldots, w_{l-1})) \geq c t_0,
\leqno (8.35)
$$
where $P(w_0, \ldots, w_{l-1})$ denotes the affine subspace 
of dimension $l-1$ spanned by $w_0, \ldots, w_{l-1}$ and $c>0$
is a constant that depends only on $M$ and $n$.
Indeed, if we cannot find some $w_l$, the whole 
$E^\ast \cap B(y_0,t_0/2)\cap A$ lies within $c t_0$ of 
$P(w_0, \ldots, w_{l-1}) \cap B(y_0,t_0)$ and can be covered by 
less than $C c^{-l+1}$ balls of radius $2ct_0$. Then
$\H^d(E^\ast \cap B(y_0,t_0/2) \cap A) \leq C c^{-l+1} (ct_0)^d
\leq C c t_0^d$ by Proposition 4.1, and this contradicts 
(8.33) if $c$ is small enough; this proves the claim.

Set $P = P(w_0, \ldots, w_{d})$, and denote by 
$W$ the convex hull of the $w_l$; thus 
$$
W \i P \cap A,
\leqno (8.36)
$$
because $A$ is convex and contains the $w_l$. 

Recall that $b\beta_{E^\ast}(y_0,t_0) \leq \varepsilon_0$ by
(8.15); so there is a $d$-plane $P'$ through $y_0$ such that 
$$\eqalign{
&\sup\big\{  \dist(w,P') \, ; \, w\in  E^\ast \cap B(y_0,t_0) \big\}
\cr&\hskip 3cm
+ \sup\big\{  \dist(p,E^\ast) \, ; \, w\in  P'\cap B(y_0,t_0) \big\}
\leq \varepsilon_0 t_0.
}\leqno (8.37)
$$
In particular, we can find $d+1$ points $p_l \in P'$,
$0 \leq l \leq d$, such that $|p_l - w_l| \leq \varepsilon_0 t_0$.

By a simple but slightly unpleasant verification using the affine 
independence estimate (8.35), we get that $P$ is $C\varepsilon_0t_0$-close to $P'$
in $B(y_0,10r_0)$, i.e., that
$$\eqalign{
\dist(p,P') &\leq C \varepsilon_0 t_0 \hbox{ for } p\in P\cap B(y_0,10r_0)
\cr&\hskip 0.5cm \hbox{ and }
\dist(p',P) \leq C \varepsilon_0 t_0 \hbox{ for } p'\in P'\cap B(y_0,10r_0).
}\leqno (8.38)
$$
The idea is simply that we can compute $P$ and $P'$ from the 
position of the $w_l$ and the $p_l$, in a stable way. 
If we were to do the computation, we would take 
coordinates where $P = \big\{ x_{d+1} = \cdots = x_n = 0 \big\}$,
observe that the $p_l-p_0$, $1 \leq l \leq d$, form a basis of
${\rm Vect}(P')$, write any unit vector $v\in {\rm Vect}(P')$ as 
$v = \sum_l a_l (p_l-p_0)$, with bounded coefficients $a_l$,
get that the coordinates $v_l$, $l \geq d+1$
are all less than $C \varepsilon_0$, and then conclude.
From (8.38) and (8.36) we deduce that 
$$
\dist(w,P) \leq C \varepsilon_0 t_0 
\ \hbox{ for } w\in  E^\ast \cap B(y_0,t_0).
\leqno (8.39)
$$
Also set $\displaystyle w = {1 \over d+1} \sum_{l=0}^d w_l \,$; 
obviously 
$$
w\in P \cap B(y_0,t_0/2)
\leqno (8.40)
$$
by (8.34), (8.36), and convexity, but another easy consequence of (8.35) 
that we leave to the reader is that
$$
B(w,c_0t_0) \cap P \i W,
\leqno (8.41)
$$
where $c_0 \in (0,1)$ is another small constant that depends 
on $M$ and $n$, but not on $\varepsilon$ or $\varepsilon_0$.
We want to choose
$$
x_2 \in F \cap B(w,c_0t_0/100)
\leqno (8.42)
$$
so we pick $p\in P'$ such that $|p-w| \leq C \varepsilon_0 t_0$,
(which exists because $w\in P \cap B(y_0,t_0/2)$ and by (8.38)),
and then $\xi \in E^\ast$ such that $|\xi-p| \leq \varepsilon_0 t_0$
(by (8.37)). By Proposition 4.1 and (8.16),
$$
\H^d(E^\ast \cap B(\xi,c_0t_0/200)) 
\geq C^{-1} c_0^dt_0^d
> \varepsilon_0 t_0^d
\geq \H^d(E^\ast \cap B(y_0,2t_0)\sm F)
\leqno (8.43)
$$
if $\varepsilon_0$ is small enough,
so we can find $x_2 \i F \cap B(\xi,c_0t_0/200)$; then (8.42)
holds because $|\xi-w| \leq C \varepsilon_0 t_0 < c_0t_0/200$.

We choose $x_2$ as in (8.42), and $r_2 = c_0 t_0/10$, and now we need 
to check that it satisfies the conditions mentioned in Lemma 8.8.
First notice that
$$
c_0 c(\varepsilon_0) r_1/10 \leq c_0 t_0/10 = r_2 
\leq t_0/10 \leq r_1/1000
\leqno (8.44)
$$
by (8.12), so the size of $r_2$ is correct. Next, 
$$
|x_2-y_0| \leq |x_2-w|+|w-y_0| \leq c_0t_0/100 + t_0/2 
\leq {2t_0 \over 3} \leq {2 r_1 \over 300}
\leqno (8.45)
$$
by (8.42), (8.40), and (8.12), and since in addition
$|y_0-x_1| \leq r_1/200$ by (8.12), we easily get that
$x_2 \i B(x_1,r_1/10)$ (as needed). 

Next we prove (8.9). Notice that
$$
B(x_2,2r_2) \i B\big(y_0,2r_2+{2t_0 \over 3}\big)
\i B(y_0,{9t_0 \over 10})
\leqno (8.46)
$$
by the beginning of (8.45) and because $r_2 = c_0 t_0/10$, so
$$
H^d(E^\ast \cap B(x_2,2r_2) \sm F) \leq 
H^d(E^\ast \cap B(y_0,t_0) \sm F) \leq \varepsilon_0 t_0^d 
= \varepsilon_0 10^d c_0^{-d} r_2^d 
\leq \varepsilon r_2^d
\leqno (8.47)
$$
by (8.16) and if $\varepsilon_0$ is chosen small enough; 
so (8.9) holds.

We still need to check that $(x_2,r_2)$ satisfies the hypotheses of
Lemma 7.38, with the same $P$ as above. First, (7.39) follows from 
(8.39) if $\varepsilon_0 \leq c_0\varepsilon/10$, because 
$r_2 = c_0 t_0/10$ and $B(x_2,2r_2) \i B(y_0,t_0)$ by (8.46).

Next, suppose that $L_j$ meets $B(x_2,2r_2)$. Then $j\in J$ 
(see (8.46) and (8.31)), and
$$
P \cap B(x_2,3r_2) \i P \cap B(w,c_0t_0) \i W \i A \i L \i L_j
\leqno (8.48)
$$ 
by (8.42) and because $r_2 = c_0 t_0/10$, then by (8.41), (8.36),
the definition of $A$, and (8.31). This proves (7.40).

Now we construct $h : E^\ast \cap B(x_2,2r_2) \times [0,1] \to \R^n$,
with the properties (7.41)-(7.44). 

Let $w \in E^\ast \cap B(x_2,2r_2)$ be given. 
By (8.39) and (8.46), $|\pi(w)-w| \leq C \varepsilon_0 t_0$. 
Then 
$$
\pi(w) \in P \cap B(x_2,2r_2+C \varepsilon_0 t_0) \i A
\leqno (8.49)
$$
by (8.48), so $\dist(w,A) \leq |\pi(w)-w| \leq C \varepsilon_0 t_0$.
With the notation of (3.5), $w \in A^\eta$ with
$\eta = C \varepsilon_0 t_0$.

Let us apply Lemma 3.17 to $L= r_0^{-1} A$
(we need to normalize as in Remark 3.25). We get a mapping
$\Pi_A$, obtained from the one we get from Lemma 3.17 by a formula
like (3.26). For $w\in E^\ast \cap B(x_2,2r_2)$, set
$$
h(w,s) = \Pi_A(w,2s)
\ \hbox{ for } 0 \leq s \leq 1/2
\leqno (8.50)
$$
and
$$\eqalign{
h(w,s) &= (2s-1) \pi(w) + (2-2s) h(w,1/2)
\cr&
= (2s-1) \pi(w) + (2-2s) \Pi_A(w,1)
\hskip0.2cm \hbox{ for } 1/2 \leq s \leq 1.
}\leqno (8.51)
$$

The fact that $h(w,0) = w$ follows from (3.18),
and $h(w,1) = \pi(w)$ from (8.51), so (7.41) holds.
Next (3.20) says that for $0 \leq s,s' \leq 1/2$,
$$
|h(w,s)-h(w,s')| = |\Pi_A(w,2s)-\Pi_A(w,2s')| 
\leq C \dist(w,A) |s-s'|
\leq C \varepsilon_0 t_0 |s-s'|
\leqno (8.52)
$$
(the two renormalizations cancel). For $1/2 \leq s,s' \leq 1$,
$$\eqalign{
|h(w,s)-h(w,s')| &= 2|s-s'| |\pi(w)-h(w,1/2)|
\cr& \leq 2|s-s'| \big(|\pi(w)-h(w,0)| + |h(w,0)- h(w,1/2)| \big)
\cr&
= 2|s-s'| \big(|\pi(w)-w)| + |h(w,0)- h(w,1/2)| \big)
\cr&
\leq C |s-s'| \varepsilon_0 t_0 
}\leqno (8.53)
$$
because $h(w,0) = w$ and $|\pi(w)-w| \leq C \varepsilon_0 t_0$,
and by (8.52). Altogether 
$$
|h(w,s)-h(w,s')| \leq C \varepsilon_0 t_0 |s-s'|
= 10 C c_0^{-1} \varepsilon_0  r_2 |s-s'|
\leq \varepsilon r_2 |s-s'|
\leqno (8.54)
$$
for $0 \leq s,s' \leq 1$, and if $\varepsilon_0$
is small enough. So we get (7.42) with $C_0 = 1$.

Next each $h(\cdot,s)$ is $C$-Lipschitz because the
$\Pi_A(\cdot,2s)$ are $C$-Lipschitz, so (7.43) holds
for some $C_0$ that depends only on $n$.

Finally we need to check that $h(w,s) \in L_j$
when $w \in E^\ast \cap L_j \cap B(x_2,2r_2)$.
For $0 \leq s \leq 1/2$, we use the fact that
$\Pi_A$ preserves the faces (as in (3.22)), so that
$h(w,s) = \Pi_A(w,2s)$ lies in any of the faces of $L_j$
that contains $w$. For $s \geq 1/2$, we use the convexity
of $A$, the fact that 
$\Pi_A(w,1) = \pi_A(w) \in A$ by (3.19) and Lemma 3.4,
and the fact that $\pi(w) \in A$ by (8.49), to get that
$h(w,s) \in A \i L \i L_j$. 

Thus (7.44) holds too, the pair $(x_2,r_2)$ satisfies the hypotheses of
Lemma 7.38, and this completes the verification of the hypothesis of 
Lemma 8.8. Finally Theorem 6.1 follows, by Lemma 8.8.
\qed

\ms
It is now very easy to prove that, still under the rigid
assumption, our quasiminimal sets have the property of concentration,
as in the following corollary. Incidentally, this corollary will be improved 
in Corollary 9.103 and Proposition 10.82. We include it here because it
is easy to get now, but the proof of Section 10 is both globally simpler and
more general. 
The property of concentration, introduced in [DMS], 
is a very nice tool to get the lower semicontinuity of $\H^d$, 
restricted to sequences of quasiminimal sets; 
this will be discussed again in Section 10, 
and even generalized slightly in Section 25. 

\ms\proclaim Corollary 8.55. For each choice of $n$, $M \geq 1$, and 
$\varepsilon > 0$, we can find $h > 0$ and $c_\varepsilon > 0$
such that the following holds. 
Suppose that $E \in GSAQ(B_0, M, \delta , h)$, with $B_0 = B(0,1)$, 
and that the rigid assumption is satisfied.
Let $r_0 = 2^{-m} \leq 1$ 
denote the side length of the dyadic cubes of the usual grid.
Let $x\in E^\ast \cap B_0$ and $0 < r < \Min(r_0,\delta)$
be such that $B(x,2r) \i B_0$. 
Assume in addition that (6.2) holds.
Then we can find a pair $(y,t)$, such that
$y \i E^\ast \cap B(x,r/100)$, $c_\varepsilon r \leq t \leq r/100$,
and
$$
\H^d(E^\ast\cap B(y,t)) \geq (1-\varepsilon) \omega_d t^d,
\leqno (8.56)
$$
where $\omega_d$ denotes the $d$-dimensional Hausdorff measure
of the unit ball in $\R^d$.

\ms
Indeed, let $(x,r)$ be as in the statement; we want to follow the 
end of the proof of Theorem 6.1. We do not need to define $F$ and
$F'$, as we did after the statement of Theorem~8.5; instead, we start 
with $x_1 = x$ and $r_1 = r$, and check that we can find $(x_2,r_2)$
as in Lemma 8.8. That is, the proof of Lemma 8.11 gives a pair
$(y_0,t_0)$, with the properties (8.12)-(8.15), where we do not need 
(8.16) and we can replace $F$ with $E^\ast$ in (8.12),
and then we use the pair $(y_0,t_0)$, as we did after Lemma 8.11, to find 
a pair $(x_2,r_2)$ that satisfies the assumptions of Lemma 7.38. 
We take $y = x_2$ and $t=r_2$. Thus there is a $d$-plane $P$ such that
$$
\dist(w,P) < \varepsilon t \hbox{ for } w\in E^\ast \cap B(y,2t)
\leqno (8.57)  
$$
because (7.39) holds, and 
$$
\pi(E^\ast \cap B(y,5t/3)) \hbox{ contains } 
P \cap B(\pi(y),3t/2),
\leqno (8.58)  
$$
by (7.46) and where $\pi$ still denotes the orthogonal
projection on $P$. For each $p \in P \cap B(\pi(y),(1-2\varepsilon)t)$,
(8.58) gives a point $w\in E^\ast \cap B(y,5t/3)$ such that
$\pi(w) = p$. But $|p-w| = |\pi(w)-w| = \dist(w,P) \leq \varepsilon t$
by (8.57) and similarly $|\pi(y)-y|\leq\varepsilon t$, 
so $|w-y| \leq |p-\pi(y)|+ |p-w| + |\pi(y)-y| \leq |p-\pi(y)| + 2 
\varepsilon t \leq t$ and $w \in B(y,t)$. 
So $P \cap B(\pi(y),(1-2\varepsilon)t) \i \pi(E^\ast \cap B(y,t))$, hence
$$\eqalign{
\H^d(E^\ast \cap B(y,t)) &\geq \H^d(\pi(E^\ast \cap B(y,t)))
\cr&
\geq \H^d(P \cap B(\pi(y),(1-2\varepsilon)t)) 
\geq \omega_d(1-2\varepsilon)^d t^d,
}\leqno (8.59)  
$$
which implies (8.56), even though only for the slightly larger
$\varepsilon'$ such that $1-\varepsilon' = (1-2\varepsilon)^d$;
but $\varepsilon'$ is as small as we want and Corollary 8.55 follows.
\qed

\ms\noindent
{\bf 9. Extension to the Lipschitz assumption.} 
\ms

In this section we intend to generalize Theorem 6.1
and Corollary 8.55 to the case when we only have the 
Lipschitz assumption. So we shall assume, throughout this
section, that 
$$
E \in GSAQ(U,M,\delta,h)
\leqno (9.1)
$$
for an open set $U \i \R^n$, and that
$$
\hbox{ the Lipschitz assumption is satisfied in $U$.}
\leqno (9.2)
$$
Recall from Definition 2.7 that this means that
there is a constant $\lambda > 0$ and a bilipschitz 
mapping $\psi : \lambda U \to B(0,1)$ such that the sets 
$\psi(\lambda L_j \cap U)$, $0 \leq j \leq j_{ max}$ 
satisfy the rigid assumption described near (2.6).

Recall also that this last comes with the rigid scale
$r_0 = 2^{-m}$, which is the side length of the cubes
in the usual dyadic grid.

We shall denote by $\Lambda \geq 1$ the bilipschitz constant
for $\psi$; thus
$$
\Lambda^{-1} |x-y| \leq |\psi(x)-\psi(y)| \leq \Lambda |x-y|
\ \hbox{ for } x,y \in U
\leqno (9.3)
$$
and we expect our main constants to depend on $M$,
$n$, and now $\Lambda$.

We shall also systematically assume that $h$ in (9.1) is sufficiently
small, depending on the various constants at hand (such as
$M$, $n$, $\Lambda$) for the proofs to work.

We start with an extension of Proposition 7.85, which will be easy because
the conclusion of Proposition 7.85 (in essence, uniform 
rectifiability) is essentially invariant under bilipschitz mappings.

\ms\proclaim Proposition 9.4.
We can find constants 
$\theta > 0$ and $C(M,\Lambda) \geq 1$, 
that depend only on $n$, $M$ and $\Lambda$, 
such that if $h$ is small enough, depending on 
$n$, $M$ and $\Lambda$,
$$
x\in E^\ast \cap U, \hskip0.2cm
0 < r < \Min(\lambda^{-1} r_0,\delta), \hskip0.2cm 
B(x,2r) \i U
\leqno (9.5)
$$
and
$$
\eqalign{
&\hbox{if $j \in [0,j_{ max}]$ is such that some face of dimension
(strictly) more than $d$}
\cr&\hskip 1.8cm
\hbox{of $L_j$ meets $B(x,r)$, then $E^\ast \cap B(x,r) \i L_j$,}
}\leqno (9.6) 
$$
then there is a closed set $G_0 \i E^\ast \cap B(x,r)$ 
and a mapping $\phi : G_0 \to \R^d$ such that 
$$
\H^d(G_0) \geq \theta r^d
\ \hbox{ and } \ 
C(M,\Lambda)^{-1} |y-z| \leq |\phi(y)-\phi(z)| \leq C(M,\Lambda) |y-z|
\hbox{ for } y,z \in G_0.
\leqno (9.7)
$$

\ms
By definition, the faces of the $L_j$ are the images
by $\lambda^{-1} \psi^{-1}$ of the faces of the boundaries
$$
\wt L_j = \psi(\lambda L_j).
\leqno (9.8)
$$
Of course we want to use Proposition 2.8, which says that
$$
\psi(\lambda E) \in GSAQ(B(0,1), \Lambda^{2d}M, 
\Lambda^{-1}\lambda\delta , \Lambda^{2d}h). 
\leqno (9.9)
$$
Let $(x,r)$ be as in the statement; we want to apply Proposition 7.85
to the set $\wt E = \psi(\lambda E)$ and the pair $(\wt x,\wt r)$, where
$\wt x = \psi(\lambda x)$ and $\wt r = \Lambda^{-1}\lambda r$.
First observe that $(\wt E)^\ast = \psi(\lambda E^\ast)$ (see the 
definition (3.2)), so $\wt x \in (\wt E)^\ast \cap B(0,1)$; it is 
clear that $0 < \wt r < \Min(r_0,\Lambda^{-1}\lambda\delta)$ by (9.5)
and also 
$$
B(\wt x,2\wt r) = B(\psi(\lambda x),2\Lambda^{-1}\lambda r)
\i \psi(B(\lambda x,2\lambda r)) \i \psi(\lambda U) = B(0,1).
\leqno (9.10)
$$
Next we check (6.2). 
Let $j$ be such that $B(\wt x,\wt r)$ meets some face of
dimension larger than $d$ of $\wt L_j = \psi(\lambda L_j)$. Since
$B(\wt x,\wt r) \i \psi(B(\lambda x,\lambda r))$ by the proof of (9.10),
$B(\lambda x,\lambda r)$ meets $\lambda L_j$, so $B(x,r)$ meets $L_j$
and (9.6) says that $E^\ast \cap B(x,r) \i L_j$ and then
$$\eqalign{
(\wt E)^\ast \cap B(\wt x,\wt r) &= \psi(\lambda E^\ast) \cap B(\wt x,\wt r)
\i \psi(\lambda E^\ast) \cap \psi(B(\lambda x,\lambda r))
\cr&= \psi(\lambda E^\ast \cap B(\lambda x,\lambda r))
\i \psi(\lambda L_j) = \wt L_j,
}\leqno (9.11)
$$
as needed for (6.2). 

Now Proposition 7.85 gives a closed set $\wt G_0 \i (\wt E)^\ast \cap 
B(\wt x, \wt r)$ and a mapping $\wt\phi : \wt G_0 \to \R^d$ such that 
$$
\H^d(\wt G_0) \geq \wt\theta \, \wt r^d
\ \hbox{ and } \ 
\wt C^{-1} |y-z| \leq |\wt\phi(y)-\wt\phi(z)| \leq \wt C |y-z|
\hbox{ for } y,z \in \wt G_0,
\leqno (9.12)
$$
where the constants $\wt\theta$ and $\wt C$ depend
only on $M$ and $n$. We set $G_0 = \lambda^{-1} \psi^{-1}(\wt G_0)$
and $\phi(y) = \lambda^{-1}\wt\phi(\psi(\lambda x))$ for $y\in G_0$.
Then $G_0 \i E^\ast \cap B(x,r)$ (because $\lambda^{-1} 
\psi^{-1}(B(\wt x, \wt r)) \i \lambda^{-1} B(\psi^{-1}(\wt x),\Lambda 
\wt r) \i B(x,r)$), and $\phi$ is $\Lambda \wt C$-bilipschitz.
Finally,
$$
\H^d(G_0) = \lambda^{-d} \H^d(\psi^{-1}(\wt G_0))
\geq \lambda^{-d} \Lambda^{-d} \H^d(\wt G_0)
\geq\lambda^{-d} \Lambda^{-d} \wt\theta \, \wt r^d = \Lambda^{-2d} \wt\theta r^d.
\leqno (9.13)
$$
So (9.6) holds with $C(M,\Lambda) = \Lambda \wt C$ and $\theta = 
\Lambda^{-2d} \wt\theta$; Proposition 9.4 follows. 
\qed

\ms
For the extension of Theorem 6.1 itself, we need to work a bit more,
because the existence of big projections, or of big pieces of Lipschitz 
graphs, is not bilipschitz invariant. We need the following extension of 
Lemma 7.38.

\ms\proclaim Lemma 9.14. There exist $C_0 \geq 1$, 
that depends only on $n$ and $\Lambda$, 
and small constants $\eta \in (0,1)$ and $\overline\varepsilon >0$, 
that depend only on $n$, $M$, and $\Lambda$, such that the following holds
if $h$ is small enough, depending only on $n$, $M$, and $\Lambda$.
Let $y\in E^\ast$ and $t>0$ be such that 
$$
0 < t < C_0^{-1} \Min(\lambda^{-1} r_0,\delta)
\hbox{ and } B(y,(C_0+1)t) \i U.
\leqno (9.15)
$$
Set 
$$
J = \big\{ j\in[0,j_{max}] \, ; \,  L_j \hbox{ meets } B(y,2t) \big\}
\ \ \hbox{ and } \ \ 
L = \bigcap_{j \in J} L_j
\leqno (9.16)  
$$
and suppose that
$$
\dist(w,L) \leq \eta t
\ \hbox{ for } w\in E^\ast \cap B(y,2t).
\leqno (9.17)  
$$
Finally let $P$ be a $d$-plane, and suppose that
$$
\dist(w,P) < \varepsilon t \ \hbox{ for } w\in E^\ast \cap B(y,2t)
\leqno (9.18)  
$$
for some $\varepsilon \leq \overline\varepsilon$.
Then
$$
\dist(p,E^\ast) \leq \varepsilon t \hbox{ for } 
p\in P \cap B(y,3t/2)
\leqno (9.19) 
$$
and, if we denote by $\pi$ the orthogonal projection onto $P$,
$$
\pi(E^\ast \cap B(y,5t/3)) \hbox{ contains } 
P \cap B(\pi(y),3t/2).
\leqno (9.20)  
$$

\ms
We could have taken $L = \R^n$ when $J = \emptyset$, but the simplest 
is to observe that  this does not happen, as $L_0 = \Omega$ meets $B(y,2t)$
because it contains $E$.
We still want to proceed as in Lemma 7.38, but we shall need to
be more careful about the way we move points around. Set 
$$
\wt L = \bigcap_{j \in J} \wt L_j = \psi(\lambda L)
\ \ \hbox{ and } \ \ 
\wh L = r_0^{-1} \wt L.
\leqno (9.21)  
$$
Notice that $\wh L$ is composed of faces of dyadic cubes of unit size.
Denote by $\wh\Pi$ the deformation that we get when we apply
Lemma 3.17, with $\eta = 1/3$, to $\wh L$. 
Thus $\wh\Pi(\wh w,s)$ is defined when $\wh w \in \wh L^{1/3}$
and $0 \leq s \leq 1$.
Next use Remark 2.25 to define a similar deformation onto 
$\wt L$, by
$$
\wt\Pi(\wt w,s) = r_0 \wh\Pi(r_0^{-1} \wt w, s)
\ \hbox{ for $\wt w \i \wt L^{r_0/3}$ and } 0 \leq s \leq 1.
\leqno (9.22)  
$$
Observe that 
$\wh\pi = \wh\Pi(\cdot,1)$ is a Lipschitz retraction
from $\wh L^{1/3}$ to $\wh L$, and 
$\wt\pi = \wt\Pi(\cdot,1)$ is a Lipschitz retraction 
from $\wt L^{r_0/3}$ to $\wt L$.

We conjugate once more to get a deformation onto
$L$. Set $\eta_0 = \Lambda^{-1}\lambda^{-1}r_0/3$, to make sure that
$$
\psi(\lambda w) \in \wt L^{r_0/3}
\ \hbox{ when } w \in L^{\eta_0} \cap U,
\leqno (9.23)  
$$
and then define $\pi_L$ and $\Pi$ by
$$
\pi_L(w) = \lambda^{-1} 
\psi^{-1}\big(\wt\pi(\psi(\lambda w))\big)
\ \hbox{ for } w \in L^{\eta_0} \cap U
\leqno (9.24)  
$$
and 
$$
\Pi(w,s) = \lambda^{-1} 
\psi^{-1}\big(\wt\Pi(\psi(\lambda w),s)\big)
\ \hbox{ for } w \in L^{\eta_0} \cap U
\hbox{ and } 0 \leq s \leq 1.
\leqno (9.25)  
$$
Notice also that 
$$
B(y,3t) \i L^{\eta_0} \cap U;
\leqno (9.26)
$$
the first inclusion holds because $L$ meets $B(y,2t)$ (by (9.16))
and $t \leq C_{0}^{-1} \lambda^{-1}r_0 \leq \eta_0/10$
by (9.15) and if $C_0$ is large enough; the second one holds because 
$B(y,(C_0+1)t) \i U$ by (9.15).

\ms
We shall now assume that (9.20) fails, combine $\Pi$ and a variant of the 
deformation used in Lemma 7.38 to build mappings $\varphi_s$, $0\leq s \leq 3$,
apply the definition of a quasiminimizer, and get a contradiction.

We start with a first stage, where we try to go from $w \in E^\ast$
to $\pi_L(w)$, but a first cut-off function $\psi_1$ will be required.
Set
$$
a = C_1 \Lambda^2 (\varepsilon + \eta) t, 
\leqno (9.27) 
$$
where the geometric constant $C_1 \geq 1$, which depends only on $n$
(through the constants of Lemma 3.17), will be chosen soon.
Then define $\xi_1 : [0,+\infty) \to [0,1]$ by
$$\eqalign{
&\xi_1(\rho) = 1 \hbox{ for } 0 \leq \rho \leq {5t \over 3} + 5a,
\cr& \xi_1(\rho) = 0 \hbox{ for } \rho \geq {5t \over 3} + 6a, 
\hbox{ and }
\cr&
\xi_1 \hbox{ is affine on } [{5t \over 3} + 5a,{5t \over 3} + 6a],
}\leqno (9.28)  
$$
and set
$$
\varphi_s(w) = \Pi(w,s\xi_1(|w-y|))
\ \hbox{ for $w \in E^\ast \cap B(y,2t)$ and } 0 \leq s \leq 1.
\leqno (9.29)  
$$
We also set
$$
\varphi_s(w) = w
\ \hbox{ for $w \in E^\ast \sm B(y,{5t \over 3} + 6a)$ and }
0 \leq s \leq 1.
\leqno (9.30)  
$$
We shall take $\overline \varepsilon$ and $\eta$ so small that
when $\varepsilon \leq \overline\varepsilon$, 
$B(y,{5t \over 3} + 10a) \i B(y,2t)$. Notice then that
the two definitions above coincide when 
$w \in E^\ast \cap B(y,2t)\sm B(y,{5t \over 3} + 6a)$, and
hence $\varphi_s(w)$ is a Lipschitz function
of $w$ and $s$. In addition,
$$
\varphi_0(w) = w
\ \hbox{ for } w \in E^\ast 
\leqno (9.31)  
$$
and 
$$
\varphi_1(w) = \Pi(w,1) = \pi_L(w) \in L
\ \hbox{ for } w \in E^\ast \cap B(y,{5t \over 3} + 5a),
\leqno (9.32)  
$$
by (3.19), Lemma 3.4, and the conjugations.

We want to estimate $|\varphi_s(w)-\varphi_{s'}(w)|$ 
when $w\in E^\ast \cap B(y,2t)$, and it will be convenient to set
$$
\wt w = \psi(\lambda w)
\hbox{ and } \wh w = r_0^{-1} \wt w
\ \hbox{ for } w\in E^\ast \cap B(y,2t).
\leqno (9.33) 
$$
Notice that $w \in L^{\eta_0} \cap  U$ by (9.26), so 
$\wt w \i \wt L^{r_0/3}$ and $\wh w \i \wh L^{1/3}$
by (9.23) and (9.21). Also set $\alpha = s\xi_1(|w-y|)$
and $\alpha' = s'\xi_1(|w-y|)$.
With these notations,
$$\leqalignno{
|\varphi_s(w)-\varphi_{s'}(w)| &= |\Pi(w,s\xi_1(|w-y|))-\Pi(w,s'\xi_1(|w-y|))|
= |\Pi(w,\alpha)-\Pi(w,\alpha')|
\cr&
= \lambda^{-1} \big|\psi^{-1}\big(\wt\Pi(\psi(\lambda w),\alpha)\big)
- \psi^{-1}\big(\wt\Pi(\psi(\lambda w),\alpha')\big)\big|
\cr& 
\leq \lambda^{-1} \Lambda 
|\wt\Pi(\psi(\lambda w),\alpha)- \wt\Pi(\psi(\lambda w),\alpha')|
& (9.34)
\cr&
= \lambda^{-1} \Lambda
|\wt\Pi(\wt w,\alpha)- \wt\Pi(\wt w,\alpha')|
= \lambda^{-1} \Lambda r_0 |\wh\Pi(\wh w,\alpha)
- \wh\Pi(\wh w,\alpha')|,
}
$$
by (9.32), (9.25), and (9.22).
We now apply (3.20), with
$$
\eta' =: \dist(\wh w,\wh L) = r_0^{-1} \dist(\wt w,\wt L)
\leq \lambda \Lambda r_0^{-1} \dist(w,L) 
\leqno (9.35)
$$
and also $\eta' \leq 1/3$ because $\wh w \i \wh L^{1/3}$.
We get that for
$w\in E^\ast \cap B(y,2t)$ and $0 \leq s \leq 1$,
$$\eqalign{
|\varphi_s(w)-\varphi_{s'}(w)| 
&\leq C \lambda^{-1} \Lambda r_0 \eta' |\alpha-\alpha'|
\leq C \Lambda^2 \dist(w,L) |\alpha-\alpha'|
\cr&\leq C \Lambda^2 \eta t |\alpha-\alpha'|
\leq C \Lambda^2 \eta t |s-s'|
}\leqno (9.36)
$$
by (9.35) and (9.17). 

Next, let $0 \leq s \leq 1$ and
$w, w' \in E^\ast \cap B(y,2t)$ be given; we want
to estimate
$$\eqalign{
|\varphi_s(w)-\varphi_{s}(w')|
&= |\Pi(w,s\xi_1(|w-y|))-\Pi(w',s\xi_1(|w'-y|))|
\cr&
\leq |\Pi(w,s\xi_1(|w-y|))-\Pi(w,s\xi_1(|w'-y|))|
\cr& \hskip 2cm
+ |\Pi(w,s\xi_1(|w'-y|))-\Pi(w',s\xi_1(|w'-y|))|.
}\leqno (9.37)
$$
If we set $\alpha = s\xi_1(|w-y|)$ and 
$\alpha' = s\xi_1(|w'-y|)$, the proof of (9.36)
yields
$$
\eqalign{
|\Pi(w,s\xi_1(|w-y|)) &-\Pi(w,s\xi_1(|w'-y|))|
= |\Pi(w,\alpha)-\Pi(w,\alpha')|
\cr&
\leq C \Lambda^2 \eta t |\alpha-\alpha'|
\leq C \Lambda^2 \eta t |\xi_1(|w-y|)-\xi_1(|w'-y|)|
\cr&
\leq C \Lambda^2 {\eta t \over a} |w-w'|
\leq C_1^{-1} C |w-w'| \leq C |w-w'|
}\leqno (9.38)
$$
by (9.28) and (9.27). The last term in (9.37) is less
than $C \Lambda^2 |w-w'|$ because the $\wh\Pi(\cdot,s)$ are 
$C$-Lipschitz and we conjugate with $\psi$ and two dilations. 
So
$$
|\varphi_s(w)-\varphi_{s}(w')| \leq C \Lambda^2 |w-w'|
\leqno (9.39)
$$
for $0 \leq s \leq 1$ and $w, w' \in E^\ast \cap B(y,2t)$.
Let us finally record that
$$
|\varphi_s(w)-w| \leq C \Lambda^2 \eta t < a
\ \hbox{ for $w\in E^\ast$ and } 0 \leq s \leq 1,
\leqno (9.40) 
$$ 
either by (9.36) (if $w \in E^\ast \cap B(y,2t)$)
or trivially by (9.30), and (for the second part)
by (9.27) and because we choose $C_1$ is large 
enough now.

\ms
We are ready to start the second stage of the deformation.
We do not want to change anything outside of $B(y,{5t \over 3} + 4a)$,
so let immediately set
$$
\varphi_s(w) = \varphi_1(w)
\ \hbox{ for $w \in E^\ast \sm B(y,{5t \over 3} + 4a)$
and } 1 \leq s \leq 2.
\leqno (9.41) 
$$
Define a second cut-off function $\xi_2 : [0,+\infty) \to [0,1]$ by
$$\eqalign{
&\xi_2(\rho) = 1 \hbox{ for } 0 \leq \rho \leq {5t \over 3} + 3a,
\cr& \xi_2(\rho) = 0 \hbox{ for } \rho \geq {5t \over 3} + 4a, 
\hbox{ and }
\cr&
\xi_2 \hbox{ is affine on } [{5t \over 3} + 3a,{5t \over 3} + 4a].
}\leqno (9.42)  
$$
This time, we try to go from $\pi_L(w)$ to $\pi_L(\pi(w))$.
Set
$$
\overline\varphi_s(w) = 
(s-1) \xi_2(|w-y|) \pi(w) + [1-(s-1)\xi_2(|w-y|)] w,
\leqno (9.43)
$$
and then
$$
\varphi_s(w) = \pi_L(\overline\varphi_s(w))
\leqno (9.44)  
$$
for $w \in E^\ast \cap B(y,{5t \over 3} + 5a)$
and $1 \leq s \leq 2$.

First observe that when 
$w \in E^\ast \cap B(y,{5t \over 3} + 5a) \sm B(y,{5t \over 3} + 4a)$,
$\xi_2(|w-y|) = 0$ by (9.42), so $\overline\varphi_s(w) = w$
and $\varphi_s(w) = \pi_L(w) = \varphi_1(w)$
by (9.32). So the two definitions of $\varphi_s(w)$ coincide on 
$E^\ast \cap B(y,{5t \over 3} + 5a) \sm B(y,{5t \over 3} + 4a)$.

Similarly, when $w \in E^\ast \cap B(y,{5t \over 3} + 5a)$,
(9.32) says that $\varphi_1(w) = \pi_L(w)$, and (9.44) gives the
same result because $\overline\varphi_1(w) = w$.
Altogether $\varphi_s(w)$ is a Lipschitz function of $s$ and $w$.
But more precisely, if $1 \leq s \leq 2$ and
$w, w' \in E^\ast \cap B(y,{5t \over 3} + 5a)$,
$$\eqalign{
|\varphi_s(w)-\varphi_s(w')|& = |\pi_L(\overline\varphi_s(w))
-\pi_L(\overline\varphi_s(w'))|
\leq C \Lambda^2 |\overline\varphi_s(w)-\overline\varphi_s(w')|
}\leqno (9.45)
$$
and, if we set $\alpha = (s-1) \xi_2(|w-y|)$
and $\alpha' = (s-1) \xi_2(|w'-y|)$,
$$\eqalign{
|\overline\varphi_s(w)-\overline\varphi_s(w')|
&= |\alpha \pi(w)+ (1-\alpha) w - \alpha' \pi(w') - (1-\alpha') w'|
\cr&
= \big|(\alpha-\alpha')(\pi(w)-w) + \alpha'(\pi(w)-\pi(w')) 
+ (1-\alpha')(w-w') \big|
\cr&
\leq |\alpha-\alpha'| |\pi(w)-w| + \alpha' |\pi(w)-\pi(w')| + 
(1-\alpha') |w-w'|
\cr&
\leq |\alpha-\alpha'| \varepsilon t  + |w-w'|
\cr&
= (s-1) \big|\xi_2(|w-y|)-\xi_2(|w'-y|)\big| \varepsilon t + |w-w'|
\cr&
\leq {\varepsilon t \over a} |w-w'| + |w-w'|
\leq 2 |w-w'|
}\leqno (9.46)
$$
by (9.18), (9.42), and (9.27). Therefore
$$
|\varphi_s(w)-\varphi_s(w')| \leq C \Lambda^2 |w-w'|
\ \hbox{ for $w, w' \in E^\ast \cap B(y,{5t \over 3} + 5a)$
and } 1 \leq s \leq 2.
\leqno (9.47)
$$
The variations in $s$ are easier, since for $1 \leq s, s' \leq 2$ and
$w \in E^\ast \cap B(y,{5t \over 3} + 5a)$,
$$\eqalign{
|\varphi_s(w)-\varphi_{s'}(w)|
&\leq C \Lambda^2 |\overline\varphi_s(w)-\overline\varphi_{s'}(w)|
= C \Lambda^2 \xi_2(|w-y|) |s-s'||\pi(w)-w|
\cr&
\leq C \Lambda^2 |s-s'||\pi(w)-w|
\leq C \Lambda^2 |s-s'|\varepsilon t 
}\leqno (9.48)
$$
by (9.44), because $\pi_L$ is $C\Lambda^2$-Lipschitz, and by (9.43)
and (9.18). The case when $w \in E^\ast \sm B(y,{5t \over 3} + 5a)$
is even more trivial, because $\varphi_s(w)=\varphi_{s'}(w)=\varphi_1(w)$
by (9.41), so
$$
|\varphi_s(w)-\varphi_1(w)| \leq C \Lambda^2 \varepsilon t |s-s'|  < a |s-s'|  
\ \hbox{ for $1 \leq s, s' \leq 2$ and } w \in E^\ast,
\leqno (9.49) 
$$
where the last inequality comes from (9.27) (if $C_1$ is large 
enough).

Let us also record the fact that
$$
\overline\varphi_2(w) = \pi(w)
\hbox{ and }
\varphi_2(w) = \pi_L(\pi(w))
\ \hbox{ for $w \in E^\ast \cap B(y,{5t \over 3} + 3a)$},
\leqno (9.50)
$$
because $\xi_2(|w-y|) = 1$ by (9.42), and
by the definitions (9.43) and (9.44).

\ms
We are now ready for the third stage where we try to move points along
$P$ to a lower-dimensional sphere. Let us first decide that
$$
\varphi_s(w) = \varphi_2(w)
\ \hbox{ for $w \in E^\ast \sm B(y,{5t \over 3} + 3a)$
and } 2 \leq s \leq 3.  
\leqno (9.51) 
$$

We now assume that (9.20) fails. This means that we can find
$$
p\in P \cap B(\pi(y),3t/2) \sm \pi(E^\ast \cap B(y,5t/3)).
\leqno (9.52)
$$
Observe that for $w \in E^\ast \cap B(y,2t) \sm B(y,5t/3)$,
$\pi(w)$ lies outside of $B(\pi(y),3t/2)$ anyway, because
$|\pi(w)-w| \leq \varepsilon t$ by (9.18). So in fact
$p$ lies out of $\pi(E^\ast \cap B(y,2t))$.

The slightly smaller compact set $\pi(E^\ast \cap \overline B(y,{11t \over 6}))$ 
does not contain $p$ either, so we can find a very small $\tau > 0$ such that
$$
P \cap B(p,\tau) \hbox{ does not meet $\pi(E^\ast \cap \overline B(y,{11t \over 6}))$.}
\leqno (9.53)
$$
We intend to move points inside
$$
B_1 = \overline B(\pi(y),5t/3).
\leqno (9.54)
$$
First define $g : P \cap B_1 \sm B(p,\tau) \to P \cap \d B_1$,
to be the radial projection on $\d B_1$, centered at $p$.
Thus $g(z)$ is characterized by the fact that 
$$
g(z) \in P \cap \d B_1 \hbox{ and } z \in [p,g(z)]
\ \hbox{ for } z \in P \cap B_1 \sm B(p,\tau).
\leqno (9.55)
$$
We also set
$$
g(z) = z 
\ \hbox{ for } z \in \R^n \sm B_1,
\leqno (9.56)
$$
observe that this defines a Lipschitz mapping
on $[P \cap B_1 \sm B(p,\tau)] \cup [\R^n \sm B_1]$,
and extend it to $\R^n$ in a Lipschitz way, so that
$g(B_1) \i B_1$.
The Lipschitz constant is very large, because we do
not control $\tau$, but we don't care.
Now set
$$
\overline\varphi_s(w) = (s-2) g(\pi(w)) + (3-s) \pi(w)
\leqno (9.57)
$$
and 
$$
\varphi_s(w) = \pi_L(\overline\varphi_s(w))
\leqno (9.58)
$$
for $w \in E^\ast \cap B(y,{5t \over 3} + 3a)$
and $2 \leq t \leq 3$. Notice that
$\pi_L(\overline\varphi_s(w))$ is well defined, because
$\overline\varphi_s(w) \in B(y,2t) \i L^{\eta_0}\cap U$
by (9.26) (also see the line below (9.30).
For such $w$, (9.50) yields $\overline\varphi_2(w) = \pi(w)$
and $\varphi_2(w) = \pi_L(\pi(w))$, so the two
definitions of $\varphi_2(w)$ coincide. 

If $w \in E^\ast \cap B(y,{5t \over 3} + 3a) \sm B(y,{5t \over 3} + a)$, 
then in addition
$$
|\pi(w)-\pi(y)| \geq |w-y|-2\varepsilon t 
\geq {5t \over 3} + a -2\varepsilon t \geq {5t \over 3} + a/2,
\leqno (9.59)
$$
by (9.18), (9.27), and if $C_1 \geq 4$; 
then $\pi(w) \in P \sm B_1$ and (9.56) says that $g(\pi(w)) = \pi(w)$.
In this case $\overline\varphi_s(w) = \pi(w)$ and $\varphi_s(w) = 
\pi_L(\pi(w)) = \varphi_2(w)$.
We claim that because of this,
$$
\varphi_s(w) = \varphi_2(w)
\ \hbox{ for $w \in E^\ast \sm B(y,{5t \over 3} + a)$
and } 2 \leq s \leq 3.
\leqno (9.60) 
$$
We just checked this when $w \in E^\ast \cap B(y,{5t \over 3} + 3a)$,
but otherwise this is just (9.51). 

We have two Lipschitz definitions of the $\varphi_s$ 
(by the formulas (9.51) and (9.57)) that overlap on the annulus 
$B(y,{5t \over 3} + 3a) \sm B(y,{5t \over 3} + a)$, 
hence $(w,s) \to \varphi_s(w)$ is Lipschitz on $E^\ast \times [2,3]$.

We just constructed Lipschitz mappings 
$\varphi_s : E^\ast \times [0,3] \to \R^n$, and we want to check the
properties (1.4)-(1.8), for the longer interval $[0,3]$, and
with respect to the ball 
$$
B = \overline B(y,C_0 t),
\leqno (9.61)
$$
where the value of $C_0 \geq 2$ will be decided soon.
We already know that $(w,s) \to \varphi_s(w)$ is Lipschitz,
so (1.8) holds. Also, $\varphi_0(w) = w$ by (9.31). Next,
$$
\varphi_s(w) = w
\ \hbox{ for $w \in E^\ast \sm B(y,{5t \over 3} + 6a)$ and }
0 \leq s \leq 3;
\leqno (9.62)  
$$
by (9.30), (9.41), and (9.51). This takes care of (1.5)
because $\overline B(y,{5t \over 3} + 6a) \i B(y,2t) \i B$ if
$\varepsilon$ and $\eta$ are small enough (see (9.27)).

For (1.6), we just need to check that 
$$
\varphi_s(w) \in B \ \hbox{ for 
$w\in E^\ast \cap B(y,{5t \over 3} + 6a)$ and }
0 \leq s \leq 3
\leqno (9.63)  
$$
because $\varphi_s(w) = w \in B$ when 
$w \in B \sm B(y,{5t \over 3} + 6a)$. For
$0 \leq s \leq 2$, (9.40) and (9.49) say that 
$|\varphi_s(w)-w| \leq 2a$, and then
$\varphi_s(w) \in B(y,{5t \over 3} + 8a) \i B$.

So we may assume that $s \geq 2$, and that 
$\varphi_s(w) \neq \varphi_2(w)$. This implies that
$w \in B(y,{5t \over 3} + a)$, by (9.60),
and that $g(\pi(w)) \neq \pi(w)$ (because (9.57)
and (9.58) apply). Then $\pi(w) \in B_1$ by (9.56), 
and also $\pi(w) \in P \sm B(p,\tau)$ by (9.53). 
Hence (9.55) applies and $g(\pi(w)) \in P \cap \d B_1$. By (9.57), 
$\overline\varphi_s(w) \in [\pi(w),g(\pi(w))] \i P \cap B_1$
and 
$$\eqalign{
|\varphi_s(w)-\varphi_2(w)| 
&= |\pi_L(\overline\varphi_s(w))-\pi_L(\overline\varphi_2(w))|
\leq C \Lambda^2 |\overline\varphi_s(w))-\overline\varphi_2(w)|
\leq 4 C \Lambda^2 t
}\leqno (9.64) 
$$
because $\diam(B_1) \leq 4t$. 
Since $\varphi_2(w) \in B(y,{5t \over 3} + 8a) \i B(y,2t)$, 
we simply choose $C_0 \geq 2 + 4C\Lambda^2$, and
(9.63) follows from (9.64). 

Finally we check (1.7). Let $j \leq j_{max}$ and 
$w \in E^\ast \cap L_j$ be given; we want to check that
$\varphi_s(w) \in L_j$ for $0 \leq s \leq 3$.
We may assume that $w \in E^\ast \sm B(y,{5t \over 3} + 6a)$,
because otherwise (9.62) says that $\varphi_s(w) = w$. 

We first consider $s \leq 1$. Then 
$\varphi_s(w) = \Pi(w,s\xi_1(|w-y|))$ by (9.29).
Set $s' = s\xi_1(|w-y|)$, and recall that the 
mapping $\wh \Pi(\cdot,s')$ of Lemma 3.17 preserve all the 
faces of unit dyadic cubes. This is also true for $\Pi(\cdot,s')$ and 
the faces of the $L_j$.
Then $\varphi_s(w)$ lies is in any of the faces of $L_j$ that
contains $x$.

Next consider $s \in (1,2]$. If 
$w \in E^\ast \sm B(y,{5t \over 3} + 4a)$,
then $\varphi_s(w) = \varphi_1(w) \in L_j$
by (9.41) and the previous case. Otherwise,
$\varphi_s(w) = \pi_L(\overline\varphi_s(w)) \in L$
by (9.44) and the definition of $\pi_L$. Since 
$w \in E^\ast \cap L_j \cap B(y,{5t \over 3} + 4a)$,
(9.16) says that $j \in J$ and $L \i L_j$,
so $\varphi_s(w) \in L_j$.

We are left with the case when $s > 2$.
If $w \in E^\ast \sm B(y,{5t \over 3} + 3a)$, (9.51) says that
$\varphi_s(w) = \varphi_2(w) \in L_j$. Otherwise, 
$\varphi_s(w) = \pi_L(\overline\varphi_s(w)) \in L \i L_j$
by (9.58) and the same argument as above.

This completes the verification of (1.4)-(1.8),
relative to the set $E^\ast$ and the ball
$B = \overline B(y,C_0 t)$.
The condition (2.4) is also satisfied, since the set $\wh W$ 
of (2.2) is contained in $B \i\i U$ by (9.61) and (9.15). 
In addition $C_0 t < \delta$ by (9.15). 
Finally recall that $E^\ast \in GSAQ(U,M,\delta,h)$, by (9.1) and 
Proposition 3.3. By Definition 2.3, (2.5) holds, i.e.,
$$
\H^d(W) \leq M \H^d(\varphi_3(W)) + h r^d,
\leqno (9.65)
$$
where 
$$
W = \big\{ w \in E^\ast \cap B \, ; \varphi_3(y) \neq y \big\}.
\leqno (9.66)
$$
First we consider 
$$
A_1 = E^\ast \cap B(y,{5t \over 3}- a)
\leqno (9.67)
$$
and
$$
A_2 = \big\{ w\in E^\ast \cap B(y,{5t \over 3} + 3a) \, ; \,
\varphi_3(w) \neq \varphi_2(w) \big\}.
\leqno (9.68)
$$
Let us check that
$$
\varphi_3(w) \in \pi_L(P \cap \d B_1)
\ \hbox{ for } w \in A_1 \cup A_2.
\leqno (9.69)
$$
First let $w \in A_1$ be given.
Notice that $|\pi(y)-\pi(w)| \leq |y-w| + |\pi(y)-y| +|\pi(w)-w|
\leq |y-w| + 2 \varepsilon t < {5t \over 3}$
by (9.18), so $\pi(w) \in B(\pi(y),{5t \over 3}) \i B_1$.
But $\pi(w) \in P$ by definition of $\pi$,
and $\pi(w) \notin P \cap B(p,\tau) \i B_1$ by (9.53),
so (9.55) says that $g(\pi(w)) \in P \cap \d B_1$.
Then
$$
\varphi_3(w) = \pi_L(\overline\varphi_3(w)) = \pi_L(g(\pi(w)))
\in \pi_L(P \cap \d B_1)
\leqno (9.70)
$$
by (9.58) and (9.57). Similarly, let $w \in A_2$ be given.
Since $\varphi_2(w)$ and $\varphi_3(w)$ are still
given by (9.58) and (9.57) in this case, the fact that
they are different implies that $g(\pi(w)) \neq \pi(w)$.
Then $\pi(w) \in B_1$ by (9.56). Again, 
$\pi(w) \notin P \cap B(p,\tau)$ by (9.53),
so $g(\pi(w)) \in P \cap \d B_1$ by (9.55), and (9.70)
holds as above. This proves (9.69).

Recall that $\H^d(P \cap \d B_1) = 0$
(this is a $(d-1)$-dimensional set), so
$$
\H^d(\varphi_3(A_1 \cup A_2)) 
\leq \H^d(\pi_L(P \cap \d B_1)) = 0,
\leqno (9.71)
$$
by (9.69) and because $\pi_L$ is Lipschitz. Since by (9.66)
$A_1 \sm W \i \varphi_3(A_1)$,
we get that
$$
\H^d(A_1 \sm W) = 0
\leqno (9.72)
$$
by (9.71). Next set
$$
A_3 = E^\ast \cap B(y,{5t \over 3} + 6a) \sm [A_1 \cup A_2],
\leqno (9.73)
$$
and let $w\in A_3$ be given. If
$w\in E^\ast \cap B(y,{5t \over 3} + 3a)$, 
then $\varphi_3(w) = \varphi_2(w)$ because $w\notin A_2$.
Otherwise, $\varphi_3(w) = \varphi_2(w)$ by (9.60).
Thus
$$
\varphi_3(w) = \varphi_2(w) \ \hbox{ for } w\in A_3.
\leqno (9.74)
$$
We cut $A_3$ into two pieces. On $A_{31} = A_3 \sm B(y,{5t \over 3} + 4a)$,
(9.41) says that $\varphi_2(w) = \varphi_1(w)$,
and (9.39) says that $\varphi_1$ is $C \Lambda^2$-Lipschitz.
On $A_{32} = A_3 \cap B(y,{5t \over 3} + 4a)$, (9.47)
says that $\varphi_2$ is $C \Lambda^2$-Lipschitz. Altogether,
$$
\H^d(\varphi_3(A_3)) = \H^d(\varphi_2(A_3))
= \H^d(\varphi_1(A_{31})) + \H^d(\varphi_2(A_{32}))
\leq C \Lambda^{2d} \H^d(A_3).
\leqno (9.75)
$$
Notice that $A_3 \i E^\ast \cap B(y,{5t \over 3} + 6a) \sm B(y,{5t \over 3}- a)$
(by (9.67)), which by (9.18) is contained in a $(6a+\varepsilon t)$-neighborhood of
$P \cap \d B(y,{5t \over 3})$. Also recall that $\varepsilon t < a$, 
by (9.27). Thus we can cover $A_3$ by less than $C ({t \over a})^{d-1}$ balls $B_i$ 
of radius $10a$. By Proposition 4.1, 
$\H^d(A_3 \cap B_i) \leq \H^d(E^\ast \cap B_i) \leq C a^d$ because $A_3 \i E^\ast$.
Altogether,
$$
\H^d(A_3) \leq \sum_{i} \H^d(A_3 \cap B_i) 
\leq C \big({t \over a}\big)^{d-1} a^d
= C \, {a \over t} \, t^d
= C C_1 \Lambda^2 (\varepsilon + \eta) \, t^d
\leqno (9.76)
$$
by (9.27), and (9.75) yields
$$
\H^d(\varphi_3(A_3)) \leq C \Lambda^{2d} \H^d(A_3)
\leq C \Lambda^{2d} C_1 \Lambda^2 (\varepsilon + \eta) \, t^d.
\leqno (9.77)
$$
Notice that if $w \in W$, (9.62) says that 
$w \in E^\ast \cap B(y,{5t \over 3} + 6a)$. Thus
$W \i A_1 \cup A_2 \cup A_3$, by (9.73).

We may now return to (9.65). First observe that
$$\eqalign{
\H^d(\varphi_3(W))
&\leq \H^d(\varphi_3(A_1 \cup A_2)) + \H^d(\varphi_3(W \cap A_3))
\cr&
=  \H^d(\varphi_3(W \cap A_3))
\leq C C_1 \Lambda^{2d+2} (\varepsilon + \eta) \, t^d
}\leqno (9.78)
$$
by (9.71) and (9.77). On the other hand, Proposition 4.1 yields
$$
C^{-1} t^d \leq \H^d(A_1) = \H^d(A_1 \cap W) \leq \H^d(W)
\leqno (9.79)
$$
by (9.67) and (9.72). We now apply (9.65) and get that
$$
C^{-1} t^d \leq \H^d(W)
\leq M \H^d(\varphi_3(W)) + h r^d
\leq C C_1 M \Lambda^{2d+2} (\varepsilon + \eta) \, t^d
+ h r^d.
\leqno (9.80)
$$
If $\eta$, $\varepsilon$, and $h$ are small enough,
depending on $M$, $n$, and $\Lambda$ (recall that $C_1$
depends only on $n$), we get the desired contradiction 
that proves (9.20).

We still need to prove (9.19), but it follows from
(9.20) and (9.18), with the same short proof as in Lemma 7.38,
a little below (7.71).

This completes our proof of Lemma 9.14.
\qed

\ms
Let us now state and prove the generalization of Theorem 6.1.

\ms\proclaim Theorem 9.81.
For each choice of constants $n$, $M \geq 1$ and $\Lambda \geq 1$, 
we can find $h > 0$, $A \geq 0$, and $\theta > 0$, 
depending only on $n$, $M$, and $\Lambda$, such that if
$E\in GSAQ(U,M,\delta,h)$ is a quasiminimal set in $U \i \R^n$
(as in (9.1)), and if the pair $(x,r)$ is such that
$$
x\in E^\ast \cap U, \hskip0.2cm
0 < r < \Min(\lambda^{-1} r_0,\delta), \hskip0.2cm 
B(x,2r) \i U,
\leqno (9.82)
$$
and
$$
\eqalign{
&E^\ast \cap B(x,r) \i L_j
\ \hbox{Êfor every $j \in [0,j_{ max}]$ such that some face}
\cr&\hskip0.4cm
\hbox{of $L_j$, of dimension (strictly) more than $d$, meets $B(x,r)$,}
}\leqno (9.83) 
$$
then  we can find a $d$-dimensional $A$-Lipschitz
graph $\Gamma \i \R^n$ such that
$$
\H^d(E\cap \Gamma \cap B(x,r)) \geq \theta r^d. 
\leqno (9.84)
$$

\ms
By $d$-dimensional $A$-Lipschitz graph, we still mean the image, 
under an isometry of $\R^n$, of the graph of some
$A$-Lipschitz function from $\R^d$ to $\R^{n-d}$. 

We want to copy the proof that we did for Theorem 6.1
(see below Theorem 8.5).
We localize as usual: we define an unbounded Ahlfors-regular set
$F' = F \cup H$ such that $E^\ast \cap B(x,r/16) \i F 
\i E^\ast \cap B(x,r/8)$, as we did near (7.87). This,
and most of our Carleson estimates, depend only on the local 
Ahlfors-regularity of $E^\ast$ (Proposition 4.1 in the rigid case,
Proposition 4.74 in the Lipschitz case).

Next we check that $F'$ is uniformly rectifiable, and
more precisely that $F' \in BPBI(\theta,C)$ (see (7.6)
for the definition), for some constants $\theta$ and $C$ 
that are a little worse than those of Proposition 9.4, 
but still depend only on $n$, $M$, and $\Lambda$. 
For this we shall reduce to Proposition 9.4,
as in the proof of Lemma 7.8; we sketch the argument,
but the reader may return to that proof for additional detail.

We are given a ball $B(y,t)$ centered on $F'$,
and we look for a bilipschitz image of a piece of $\R^d$
inside $F \cap B(y,t)$. 
When $y \in H$ or $y \in F$ and $t \geq 3 r$, 
we can easily find this bilipschitz image inside $H$,
so we may assume that $y \in F$ and $t \leq 3r$.
Then we find a reasonably large cube $Q$, of the cubical
patchwork for $F$, such that $y \in Q \i F \cap B(y,t/10)$
(as in (7.12)). And inside $Q$, we find a point $w \in F$ 
such that $\dist(w,E^\ast \sm F) \geq \tau \diam(Q)$
(for some constant $\tau > 0$ that depends on the local
Ahlfors regularity constant), as in (7.13). 
This part uses the fact that our cubical patchwork 
for $F$ is adapted to the set $E^\ast$, as in (7.11).
Then we apply Proposition 9.4 to the ball $B(w,\tau \diam(Q))$,
and find a set $G \i E^\ast \cap B(w,\tau \diam(Q))$ which is 
the bilipschitz image of a piece of $\R^d$; 
this set is contained in $F \cap B(y,t)$, in particular because 
$\dist(w,E^\ast \sm F) \geq \tau \diam(Q)$,
and it is large enough because $\diam(Q) \geq C^{-1} t$;
see (7.14)-(7.16) for the verification. Thus 
$F' \in BPBI(\theta,C)$, as announced.

By Theorem 7.7, $F' \in BWGL(\varepsilon,C(\varepsilon))$
for every $\varepsilon$, where as usual $C(\varepsilon)$
depends also on $n$, $M$, and $\Lambda$. As in the rigid
case, it is enough to show that $F'$ has big projections, 
because then Theorem 8.5 will say that $F'$ contains big pieces
of Lipschitz graphs, and we can use one of these pieces
(contained in $F \i E^\ast \cap B(x,r/8)$) in the statement
of Theorem~9.81.

Again we are given a ball  $B(x_1,r_1)$ centered on $F'$,
and now we want to find a $d$-plane $P$ such that
$$
\H^d(\pi(F' \cap B(x_1,r_1))) \geq \alpha r_1^d,
\leqno (9.85)
$$
where $\pi$ denotes the orthogonal projection on $P$.
When $x_1 \in H$ or $x_1 \in F$ and $r_1 \geq 3 r$, 
we easily get this with $P=H$, because $\H^d(\pi(F' \cap B(x_1,r_1)))
\geq \H^d(H \cap B(x_1,r_1)) \geq C^{-1} r_1^d$. So we may
assume that $x_1 \in F$ and $r_1 \leq 3 r$.

Let $\varepsilon_0 > 0$ be very small, to be chosen later.
We proceed as in Lemma 8.11 to find a pair
$(y_0,t_0)$, with the properties (8.12)-(8.16), except
that in (8.12) we replace $r_1/100$ with the smaller 
$(100C_0)^{-1} r_1$, with $C_0$ as in (9.15),
that in (8.13) we require that 
$\dist(y_0,F_1) \geq 10 \Lambda^2 t_0$,
and that in (8.15) and (8.19) we replace
$b\beta_{E^\ast}$ with $b\beta_{F}$.
The reader recalls that in Lemma 8.11, the
pair $(x_1,r_1)$ was also such that $x_1 \in F$ and 
$r_1 \leq 3 r$. Most of the proof of Lemma 8.11
can be repeated here, because it relies on Carleson measure 
computations based on the local Ahlfors regularity of 
$E^\ast$ and $F$, and simple distance estimates with faces
of our dyadic grid on $U$. There is just one exception,
which is the Carleson measure estimate (8.20) on the bad set
${\cal B}_3$ where $E^\ast$ is not flat. Here we replace 
$b\beta_{E^\ast}$ with $b\beta_{F}$ in (8.15) and (8.19),
and we get the analogue of (8.20) because $F$ is uniformly 
rectifiable, hence satisfies a bilateral weak geometric lemma 
(i.e, $F' \in BWGL(\varepsilon,C(\varepsilon))$ as above).

So we get the pair $(y_0,t_0)$, and we now check that
Lemma 9.14 applies to the pair $(y_0,t_0/3)$.
The first part of (9.15) holds because
$$
{t_0 \over 3} \leq  {r_1 \over 100 C_0} \leq {4 r \over 100 C_0}
< {4 \over 100 C_0} \Min(\lambda^{-1} r_0, \delta)
\leqno (9.86)
$$
by the modified (8.12), (8.7), and (9.82).
For the second part, observe that
$$
|y_0-x| \leq |y_0-x_1|+|x_1-x| \leq {r_1 \over 200}+{r\over 8}
\leq {r \over 4}
\leqno (9.87)
$$ 
by (8.12) and because $x_1 \in F \i B(x,r/8)$ and $r_1 \leq 3r$; 
then (9.86) and (9.82) yield
$$
B(y_0,{(C_0+1)t_0 \over 3}) \i B(x,{r\over 4}+{(C_0+1)t_0 \over 3}) 
\i B(x,r/2) \i U.
\leqno (9.88)
$$
Recall that $J$ in (9.16) is not empty (because $y \in \Omega = L_0$), 
and let us check (9.17). A first possibility is that for each $j\in J$, 
$L_j$ has a face of dimension larger than $d$ that meets $B(y_0,2t_0/3)$. 
Since $B(y_0,2t_0/3) \i B(x,r)$, (9.83) says that 
$E^\ast \cap B(y_0,2t_0/3) \i E^\ast \cap B(x,r) \i L_j$ for each $j$, 
hence $E^\ast \cap B(y_0,2t_0/3) \i L$
and (9.17) holds with $\eta = 0$.

So let us assume that for some $j\in J$, no face
of $L_j$ of dimension larger than $d$ meets $B(y_0,2t_0/3)$.
Since $j \in J$, $L_j$ meets $B(y_0,2t_0/3)$, which means
that some face $H$ of $L_j$ meets $B(y_0,2t_0/3)$. Then
$H$ is of dimension at most $d$.

It will be easier to make our metric computations with standard
dyadic cubes, so we set $h(y) = r_0^{-1} \psi(\lambda y)$
for $y\in \R^n$, and observe that $h(H)$ is a standard
unit dyadic face by construction. 
Notice that
$$
\dist(h(y_0),h(H)) 
\leq r_0^{-1} \lambda \Lambda \dist(y_0,H)
\leq r_0^{-1} \lambda \Lambda t_0
\leq {1 \over 10},
\leqno (9.89)
$$
where the last inequality comes from (9.86)
if we took $C_0 \geq \Lambda$.
But by (8.13) with the new constant $10 \Lambda^2$,
$$ 
\dist(h(y_0), h(F_1))
\geq r_0^{-1} \lambda \Lambda^{-1} \dist(y_0, F_1)
\geq 10 r_0^{-1} \lambda \Lambda t_0,
\leqno (9.90)
$$
so (by (9.89)) $H$ is not contained in $F_1$, and it is $d$-dimensional.
By (8.14),
$$
\dist(h(w),h(F_2)) \leq r_0^{-1} \lambda \Lambda \dist(w,F_2)
\leq r_0^{-1} \lambda \Lambda \varepsilon_0 t_0
\leqno (9.91)
$$
for $ w\in E^\ast \cap B(y_0,2t_0)$. We use (9.89) to find 
$\xi\in h(H)$ such that $|h(y_0)-\xi| \leq r_0^{-1} \lambda \Lambda t_0$;
then if $\wh H$ is any other $d$-dimensional face of the dyadic net
(i.e., $\wh H \neq h(H)$),
$$\eqalign{
\dist(h(y_0), \wh H) &\geq \dist(\xi, \wh H) - r_0^{-1} \lambda \Lambda t_0
\geq \dist(\xi, \d(h(H)) - r_0^{-1} \lambda \Lambda t_0
\cr&
\geq \dist(h(y_0), \d(h(H)) - 2r_0^{-1} \lambda \Lambda t_0
\geq 8 r_0^{-1} \lambda \Lambda t_0,
}\leqno (9.92)
$$
where the main inequalities come from (3.8) and (9.90).
If $w\in E^\ast \cap B(y_0,2t_0)$, this yields
$$
\dist(h(w), \wh H) \geq \dist(h(y_0), \wh H) - 2r_0^{-1} \lambda \Lambda t_0
\geq 6 r_0^{-1} \lambda \Lambda t_0.
\leqno (9.93)
$$
In other words, all the other faces that compose $h(F_2)$
are too far, and (9.91) implies that
$$
\dist(h(w),h(H)) \leq r_0^{-1} \lambda \Lambda \dist(w,F_2)
\leq r_0^{-1} \lambda \Lambda \varepsilon_0 t_0
\ \hbox{ for } w\in E^\ast \cap B(y_0,2t_0).
\leqno (9.94)
$$
In fact, (9.93) also implies that $\dist(w,H') \geq 6 t_0$
for every $d$-dimensional face $H'= \lambda^{-1}\psi^{-1}(r_0 \wh H)$ 
of dimension $d$ of our (distorted) dyadic grid on $U$, other than $H$,
and (8.14) (or the second half of (9.91)) now says that
$$
\dist(w,H) = \dist(w,F_2) \leq \varepsilon_0 t_0
\ \hbox{ for } w\in E^\ast \cap B(y_0,2t_0).
\leqno (9.95)
$$

Let us now check that
$$
L_i \hbox{ contains $H \, $ for } i \in J.
\leqno (9.96)
$$
Let $G$ be a be a face of $L_i$ that meets $B(y_0,2t_0/3)$;
such a face exists by definition of $J$. Also let 
$\xi \in h(H)$ be, as above, such that 
$|h(y_0)-\xi| \leq r_0^{-1} \lambda \Lambda t_0$.
If (9.96) fails, $G$ does not contain $H$, 
$h(G)$ does not contain $h(H)$, $h(H)$ is not reduced to one point
because it is $d$-dimensional, and so (3.8) says that 
$$\eqalign{
\dist(h(y_0), h(G)) 
&\geq \dist(\xi,h(G)) - r_0^{-1} \lambda \Lambda t_0
\geq \dist(\xi, \d(h(H)) - r_0^{-1} \lambda \Lambda t_0
\cr&
\geq \dist(h(y_0), \d(h(H)) - 2 r_0^{-1} \lambda \Lambda t_0
\geq 8 r_0^{-1} \lambda \Lambda t_0,
}\leqno (9.97)
$$
by (9.92) or directly (9.90). This is impossible, because
$$
\dist(h(y_0), h(G)) \leq r_0^{-1} \lambda \Lambda\dist(y_0, G)
\leq r_0^{-1} \lambda \Lambda t_0
\leqno (9.98)
$$ 
since $G$ meets $B(y_0,t_0)$. So (9.95) holds.

By (9.95) and the definition (9.16), $L$ contains $H$. Then 
$\dist(w,L) \leq \dist(w, H) \leq \varepsilon_0 t_0$
for $w\in E^\ast \cap B(y_0,2t_0)$, by (9.95).
This proves (9.17), with $\eta = 3\varepsilon_0$.

Finally we need to check (9.18) for some $P$; we know from the 
modified (8.15) that $b\beta_{F}(y_0,t_0) \leq \varepsilon_0$, 
so there is a plane $P$ through $y_0$ such that in particular
$$
\dist(w,P) \leq \varepsilon_0 t_0
\hbox{ for } w\in F \cap B(y_0,t_0).
\leqno (9.99)
$$
If (9.18) fails for this $P$ (and the pair $(y_0,t_0/3)$),
we can find $w \in E^\ast \cap B(y_0, 2t_0/3)$ such that
$\dist(w,P) \geq \varepsilon t_0/3$. If 
$\varepsilon_0 \leq \varepsilon/6$, (9.99) implies
that $\dist(w,F) \geq \varepsilon t_0/3$.
But then
$$
\H^d(E^\ast \cap B(y_0,2t_0) \sm F)
\geq \H^d(E^\ast \cap B(w,\varepsilon t_0/3))
\geq C^{-1} \varepsilon^d t_0^d
\leqno (9.100)
$$
by Proposition 4.74. This contradicts (8.16) if $\varepsilon_0$
is small enough (depending on $\varepsilon$, $M$, and $\Lambda$);
thus (9.18) holds.

This completes the verification of the hypotheses of Lemma 9.14
for the pair $(y_0,t_0/3)$. We apply Lemma 9.14 and get
that 
$$
\pi(E^\ast \cap B(y_0, 5t_0/9) \hbox{ contains }
P \cap B(\pi(y),t_0/2), 
\leqno (9.101)
$$
as in (9.20). Then we use the modified (8.16) to get that
$$\leqalignno{
\H^d(\pi(F \cap B(y_0, 5t_0/9)) 
&\geq \H^d(\pi(E^\ast \cap B(y_0, 5t_0/9)))
-\H^d(E^\ast \cap B(y_0, 5t_0/9) \sm F)
\cr&
\geq \H^d(P \cap B(\pi(y),t_0/2)) - \varepsilon_0 t_0^d
\geq C^{-1} t_0^d
&(9.102)
}
$$
if $\varepsilon_0$ is small enough.
But $B(y_0, 5t_0/9) \i B(x,r_1)$ by (8.12),
so (9.102) implies that 
$\H^d(\pi(F \cap B(x_1,r_1)) \geq C^{-1} t_0^d$,
and (9.85) follows because $t_0 \geq c(\varepsilon_0) r_1$
by (8.12).

Thus $F'$ has big projections, hence contains big pieces
of Lipschitz graphs, and we know that Theorem~9.81 follows.
\qed

\ms
Let us also state the property of concentration under
the Lipschitz assumption. The following is a generalization 
of Corollary 8.55; it will be further generalized in 
Section 10, where we shall (simplify the proof and)
remove the unnatural assumption (9.105).  See Proposition 10.82.

\ms\proclaim Corollary 9.103. For each choice of $n$, $M \geq 1$, 
$\Lambda \geq 1$ and $\varepsilon > 0$, we can find $h > 0$ and 
$c_\varepsilon > 0$ such that the following holds. 
Suppose that $E \in GSAQ(U,M,\delta,h)$ for some open set
$U \i \R^n$, and that the Lipschitz assumption is satisfied,
with the constants $\lambda$ and $\Lambda$ (as in (9.3)).
Also denote by $r_0 = 2^{-m} \leq 1$ the side length of the 
dyadic cubes of the usual grid. Then let $(x,r)$ be such that
$$
x\in E^\ast \cap U, \hskip0.2cm
0 < r < \Min(\lambda^{-1} r_0,\delta), \hskip0.2cm 
B(x,2r) \i U,
\leqno (9.104)
$$
and 
$$
\eqalign{
&E^\ast \cap B(x,r) \i L_j
\ \hbox{Êfor every $j \in [0,j_{ max}]$ such that some face}
\cr&\hskip0.4cm
\hbox{of $L_j$, of dimension (strictly) more than $d$, meets $B(x,r)$.}
}\leqno (9.105) 
$$
Then we can find a pair $(y,t)$, such that
$y \i E^\ast \cap B(x,r/100)$, $c_\varepsilon r \leq t \leq r/100$,
and
$$
\H^d(E^\ast\cap B(y,t)) \geq (1-\varepsilon) \omega_d t^d,
\leqno (9.106)
$$
where $\omega_d$ denotes the $d$-dimensional Hausdorff measure
of the unit ball in $\R^d$.

\ms
The proof is the same as for Corollary 8.55, except that we use
Lemma 9.14 instead of Lemma 7.38; the point is that given 
a small constant $\varepsilon > 0$ and a pair 
$(x,r)$ as in the statement, we can find a new pair $(y,t)$,
with $y\in E^\ast \cap B(x,r/2)$ and $c(\varepsilon) \leq t \leq r/10$,
which satisfies the assumptions of Proposition 9.14.
For instance, we can proceed as in the final part of the proof 
of Theorem~9.81 above and pick the pair $(y_0,t_0/3)$ 
associated to $x_1=x$ and $r_1=r$.

Then we just need to use the properties (9.18)-(9.20) to prove (9.106) 
for this pair (unfortunately with a slightly larger constant 
$\varepsilon'$, but this does not matter), exactly as we did in 
(8.57)-(8.59).
\qed

\bigskip
\centerline{PART IV : LIMITS OF QUASIMINIMAL SETS}
\ms
In this part, we want to generalize results that come mainly
from [D2], that concern  
limiting properties of quasiminimal (or almost minimal, or
minimal) sets. The main result of this part (Theorem 10.8) is that 
our various classes of quasiminimal sets are stable under limits, and 
the main reason why it holds is the lower semicontinuity of
$\H^d$, restricted to a sequence of quasiminimal sets with uniform 
quasiminimality constants; cf. Theorem~10.97. See the general 
introduction for a discuss of the interest of these results.

In turn the main ingredient in the proof of Theorems 10.97 and 10.8
is the fact that our quasiminimal sets satisfy a concentration 
property that was introduced by Dal Maso, Morel, and Solimini [DMS] 
in the different, but related context of minimizers for
the Mumford-Shah functional in image processing (see Proposition 
10.82 below). They used this property to prove lower semicontinuity 
results for $\H^d$ (on some minimizing sequences) and get an existence 
result for minimal segmentations of the functional. At the same time,
E. De Giorgi, M. Carriero, and A. Leaci [DeCL] 
obtained the same existence result, using the weak form of the 
functional and a compactness result of Ambrosio in the class $SBV$
of special bounded variation functions.

For minimal sets and surfaces, it seems that the idea of
using the concentration property to obtain existence results was not 
considered before [D2], 
probably because people were very happy with the quite strong
compactness properties of integral currents and varifolds.
Most often, when a limiting property for minimal sets was needed,
people would revert to currents or varifolds, take a limit there, 
and return to sets. Nonetheless, it is good to have
limiting theorems like Theorems 10.97 and 10.8, both because this 
looks more direct and, for instance, because
some minimal sets are hard to describe as supports of currents, typically 
for orientation reasons or because multiplicities could become too 
large.

Recall that we already proved the concentration property in some cases, 
as a consequence of the uniform rectifiability of the quasiminimal 
sets; see Corollaries 8.55 and 9.103. But the very surprising thing,
at least for the author, is that there is a more direct route to this, 
through the fact that limits of quasiminimal sets 
(with uniform quasiminimality constants) are rectifiable 
(Proposition 10.15 below), which gives a simpler and general proof.
Also, we shall give a slightly more direct proof of Theorem 10.97
(still based on the same ideas but improved by Y. Fang) in Section 25,
that also works when $\H^d$ is multiplied by some elliptic integrands.

We prove Theorem 10.97 in the next section, but the proof of 
Theorem 10.8 (the quasiminimality of limits) will be quite long, and 
we split in into smaller pieces (Sections 11-19). The difficulty is
that given a sequence $\{ E_k \}$ of quasiminimal sets that converge 
to the set $E$, and a competitor $F = \varphi_1(E)$ for $E$, the
obvious competitor $\varphi_1(E_k)$ may be very bad, and we have to 
spend some energy to make it better, typically by pinching parallel
leaves of $\varphi_1(E_k)$ to diminish their total mass.
Unfortunately, unlike what happens with uniform rectifiability,
we essentially have to redo most of the proof of [D2]. 
Note that it is far from impossible (since the author worked by himself for
all of this) that a better proof exists. But in the mean time we seem
to be stuck with a long, technical, but not so inventive proof.

\msi
{\bf 10. Limits of quasiminimal sets: the main statement,
rectifiability, and l.s.c.}

\ms 
In this section and the next ones, we take a sequence of sets $E_k$, 
which are quasiminimal in a same domain, with sliding conditions 
with respect to the same boundary pieces $L_j$, and with uniform
quasiminimality constants, and we try to prove that if the cores 
$E_k^\ast$ converge (in local Hausdorff distance) to $E$, then $E$ 
is quasiminimal with the same constants.
It is also natural to make the $L_j$ vary as well, 
but it will be simpler for us not to do this until Section 23. 

In this section we shall take care of the (plain) rectifiability of 
$E$, the uniform concentration property for the $E_k$, and the
lower semicontinuity of $\H^d$ along our sequence.

Let us describe our assumptions for the next few sections.
We fix an open set $U \i \R^n$, and boundary pieces $L_j$, 
$0 \leq j \leq j_{max}$, and we assume that
$$
\hbox{ the Lipschitz assumption is satisfied in $U$.}
\leqno (10.1)
$$
Recall from Definition 2.7 that this assumption comes 
with a positive dilation constant $\lambda > 0$,
and a bilipschitz mapping $\psi : \lambda U \to B(0,1)$;
we shall denote by $\Lambda$ (a bound for) the bilipschitz 
constant of $\psi$, as in (9.3).

Next we are given a sequence $\{ E_k \}$ of closed sets 
in $U$. By this we mean that $E_k$ is contained in $U$ and relatively 
closed in $U$, but it would make no difference if we just assumed that 
$E_k \i \R^n$ and that its intersection with $U$ is 
closed in $U$, because anyway we shall never look at points that 
lie outside of $U$. 
We assume that there are constants $M \geq 1$,
$\delta \in (0,+\infty]$, and $h >0$ (systematically assumed
to be small enough, depending on $n$, $M$, and $\Lambda$) 
such that for all $k$,
$$
E_k \in GSAQ(U,M,\delta,h),
\leqno (10.2)
$$
and
$$
\hbox{$E_k$ is coral, i.e., } E_k^\ast = E_k 
\leqno (10.3)
$$
(see Definitions 2.3 and 3.1).
This time there is no point in trying to avoid the assumption (10.3): 
if the $E_k$ are not coral, the sets $E_k^\ast \sm E_k$ may
converge to anything, even though they have a vanishing measure.
This means that in concrete problems where the sets 
$E_k^\ast \sm E_k$ have some meaning, one may need to do
something special about them, probably after taking care of the
$E_k^\ast$. 

We also assume that there is a closed set $E$ in $U$, such that
$$
\lim_{k \to +\infty} E_k = E
\ \hbox{ locally in } U,
\leqno (10.4)
$$
where the limit is defined as follows. For each choice of $x\in U$ and
$r> 0$ such that $\overline B(x,r) \i U$, and of two sets
$E$, $F$, which we can assume to be closed in $U$, we set
$$\eqalign{
d_{x,r}(E,F) &= {1 \over r} \, \sup \big\{ \dist(y,F) \, ; \, y\in E \cap B(x,r) \big\}
\cr& \hskip2cm
+ {1 \over r} \, \sup \big\{ \dist(y,E) \, ; \, y\in F \cap B(x,r) \big\},
}\leqno (10.5)
$$
where by convention $\sup \big\{ \dist(y,F) \, ; \, y\in E \cap B(x,r) \big\}
= 0$ when $E \cap B(x,r)$ is empty, and similarly 
$\sup \big\{ \dist(y,E) \, ; \, y\in F \cap B(x,r) \big\}=0$
when $F \cap B(x,r)$ is empty. But when $E \cap B(x,r)$ is nonempty
but $F$ is empty, for instance, we set $d_{x,r}(E,F) = +\infty$;
we include these cases to be able to say that sets that go away to the 
boundary tend to the empty set, but in fact this situation does not 
interest us (because we know that $\emptyset$ is minimal, for 
instance). Now (10.4) means that
$$
\lim_{k \to +\infty} d_{x,r}(E_k,E) = 0
\leqno (10.6)
$$
for all choices of $x\in U$ and $r> 0$ such that $\overline B(x,r) \i U$.
This is easily seen to be equivalent to other ways of defining 
(10.4), for instance where we would replace our family of balls $B(x,r)$ with 
an exhaustion of $U$ by compact subsets. The main point of using this 
definition is that given any sequence $\{ E_k \}$ of closed sets in $U$,
we can always find a subsequence that converges to some closed set.

There is a small technical assumption that we want to make
when the Lipschitz assumption holds:
$$\eqalign{
&\hbox{for each $0 \leq j \leq j_{max}$ and each face $F \i U$
of our net such that}
\cr&
\hbox{${\rm dimension}(F) > d$ and $F \i L_j$, but the interior of $F$ 
(as a face)}
\cr&
\hbox{is not contained in the interior of $L_j$ 
(as a subset of $\R^n$), we have}
\cr& 
\hbox{that for $\H^d$-almost every interior point
$y$ of $F$, we can find $t>0$}
\cr&
\hbox{such that the restriction of $\psi$ 
to $\lambda F \cap B(\lambda y,t)$ is of class $C^1$.}
}\leqno (10.7)
$$
This condition is a little strange so let us explain a little.
Notice that we require the exceptional set to be small for $\H^d$,
regardless of the dimension of the face $F$. The interior of $L_j$
is really taken with the topology of $\R^n$, not with respect to
the dimension of some faces: we add this constraint on $F$ because
we don't want to put regularity conditions on the faces $F$ that lie in 
the middle of out initial domain $\Omega = L_0$, for instance. But we
require some control on the boundary of $\Omega$.

Let us state (10.7) in a slightly different way. For each $y\in U$,
denote by $F(y)$ the smallest face of our grid that contains $y$; thus
$y$ is an interior point of $F(y)$ (when $F(y) = \{ y \}$,
we may say that ${\rm int}(F(y)) = \{ y \}$, but this case will be 
rapidly dismissed anyway). We require that for $0 \leq j \leq j_{max}$ and
for $\H^d$-almost every point of $L_j$, if ${\rm dimension}(F(y)) > d$
and $y$ does not lie in the $n$-dimensional interior of $L_j$,
we can find $t>0$ such that the restriction of $\psi$ 
to $\lambda F(y) \cap B(\lambda y,t)$ is of class $C^1$.

In this second condition also, we only really exclude the case when 
${\rm dimension}(F(y)) = d$, because the lower dimensional skeletons
have vanishing $\H^d$-measure anyway.

Let us check that the two conditions are equivalent.
If (10.7) holds and $y$ lies in none of the exceptional sets 
associated to faces $F$, and if $y \in U$ does not lie in the
$n$-dimensional interior of $L_j$ and is such that 
${\rm dimension}(F(y)) > d$, we can apply (10.7) to $F = F(y)$
(which is contained in $L_j$ because $y\in L_j$ and because 
$L_j$ is composed of faces), and we get the desired $t>0$
because $y$ does not lie in the exceptional set of $F$.

If our second condition holds and the face $F$ is such that
$F \i L_j$, ${\rm dimension}(F) > d$, and $F$ is not contained
in the $n$-dimensional interior of $L_j$, observe that for
every interior point $y$ of $F$, we have that $F(y) = F$.
Thus, for all the points $y \in {\rm int}(F)$ that do not lie in the
exceptional set of the second condition, we can find $t>0$ as
needed.

We shall state later a weaker (but a little more complicated to state) 
condition that works as well (see Remark 19.52), 
but we did not find any obvious way to get rid of it entirely.
Notice that (10.7) does not require anything in faces $F$
of dimension $d$, and that is trivially satisfied 
under the rigid assumption, or if the bilipschitz mapping 
$\psi : \lambda U \to B(0,1)$ is of class $C^1$.

Here is our main result about limits.

\ms\proclaim Theorem 10.8.
Let $U$, $\{ E_k \}$, and $E$ satisfy the hypotheses above
(including (10.7) if the Lipschitz assumption holds).
Also suppose that $h$ is small enough, depending only on $n$, 
$M$, and $\Lambda$. Then $E$ is coral, and
$$
E \in GSAQ(U,M,\delta,h),
\leqno (10.9)
$$
with the same constants $M$, $\delta$, and $h$.

\ms
See Remark 19.52 for a statement where (10.7) 
is slightly weakened.

The smallness of $h$ and the additional assumption
(10.7) will only be used in the various limiting arguments,
but will have no impact on the constants in the conclusion. 

In addition, the smallness of $h$ is only used to show that
the $E_k$ have, uniformly in $k$, some good regularity properties
(that imply, in particular, the lower semicontinuity estimate (10.98)),
and then we don't need it any more. So, if we have (10.2)
for some combination of $M$, $\delta$, and $h$ for which $h$ is
small enough (as required), and (10.2) also holds for some other
combination of $M$, $\delta$, and $h$ (this time, with no constraint),
then our conclusion (10.9) holds for both combinations.

Theorem 10.8 generalizes Theorem 4.1 on page 126 of [D2]; 
we shall try to follow the proof, but since many modifications 
will be needed in the middle of the construction, we shall need to explain
things with more detail than in the previous part.

\ms
The main goal of the rest of this section
is to prove the lower semicontinuity of $\H^d$,
when we restrict our attention to sequences of
quasiminimal sets that satisfy the assumptions (10.1)-(10.4) of 
Theorem 10.8 (we shall not need (10.7) for quite some time).
This is Theorem 10.97 below, which in a way is the 
main tool for our proof of Theorem~10.8.

We intend to deduce Theorem 10.97 from a result of 
Dal Maso, Morel, and Solimini [DMS] 
which says that $\H^d$ is lower semicontinuous along
uniformly concentrated sequences, but we shall 
follow a different route to the concentration property.

In [D2] we deduced it from the local 
uniform rectifiability of the $E_k$, but here we were only 
able to get this under additional (and not too natural) 
assumptions on the dimensions of the faces of the $L_j$. 
That is, Corollaries 8.55 and 9.103 have some unnatural 
assumptions that we want to avoid.

So we shall prove a weaker regularity condition, 
the existence of reasonably large balls where 
a given quasiminimal set is approximated by a $d$-plane,
and then show that we can use it to prove the desired
concentration property. See Lemma 10.21 for the approximation 
property, and Proposition 10.82 for the concentration property.
Also see the later Section~25 for a more direct proof of 
Theorem 10.97 that works with some elliptic integrands.

\ms
Before we get to Lemma 10.21, we shall consider a sequence
of quasiminimal sets that satisfy the assumptions (10.1)-(10.4)
(again, we do not need (10.7) in this section),
and prove various simple properties.

The first ones are the local Ahlfors regularity of
the limit (see (10.11)), very rough lower and upper
semicontinuity properties of $\H^d$ along the sequence
((10.12) and (10.14)), and the rectifiability of the limit 
(Proposition 10.15).

This last, which is not a direct consequence of the rectifiability
of the $E_k$, but essentially follows from its proof,
is useful because we use it to prove the approximation
Lemma~10.21 (through a compactness argument), and because we
shall use the rectifiability of $E$ to construct the
competitors in the next sections.

The construction of a competitor for the main part of 
the argument will only start in the next section,
and estimates will continue up until Section 19. 

\ms
So let $\{ E_k \}$ and its limit $E$ satisfy our assumptions 
(10.1)-(10.4); we want to derive a few simple properties.

Let us first observe that the $E_k$ are locally Ahlfors-regular,
with uniform estimates. This means that there exists a constant
$C_M$, that depends only on $n$, $M$ and $\Lambda$, such that 
$$
C_M^{-1} r^d \leq \H^d(E_k \cap B(x,r)) \leq C_M r^d
\leqno (10.10)
$$
for every pair $(x,r)$ such that $x\in E_k$,
$0 < r < \Min(\lambda^{-1} r_0,\delta)$,
and $B(x,2r) \i U$. This is an easy consequence
of Propositions 4.1 and 4.74, that we already used a lot in the
last sections. We deduce from this that
$$
C_M^{-1} r^d \leq \H^d(E \cap B(x,r)) \leq C_M r^d
\leqno (10.11)
$$
when $x\in E$ and $0 < r < \Min(\lambda^{-1} r_0,\delta)$
are such that $B(x,2r) \i U$, with a possibly larger constant $C_M$, but
that still depends only on $n$, $M$ and $\Lambda$. This is easy to check,
because the (local) Ahlfors regularity of $E$ follows from the 
existence of a locally finite Borel measure $\mu$ on $E$ that satisfies
(10.11) (where we would replace $\H^d(E \cap B(x,r))$ with 
$\mu(E \cap B(x,r))$), and such a measure is easy to obtain as
a weak limit of the restriction of $\H^d$ to $E_k$, or for a subsequence of 
$\{ E_k \}$. More detail can be found in the proof of
Lemma 4.2 in [D2]. 
Also notice that $E$ is clearly coral because of (10.11).

Let us deduce from (10.10) and (10.11) that
$$
\H^d(E\cap V) \leq C_M\liminf_{k \to +\infty} \, \H^d(E_k \cap V)
\ \hbox{ for every open set $V \i U$,}
\leqno (10.12)
$$
where again $C_M$ depends only on $n$, $M$ and $\Lambda$. To do this,
cover $E\cap V$ by balls $B(x_j,r_j)$
such that $x_j \in E$, $0 < r_j < 10^{-1}\Min(\lambda^{-1} r_0,\delta)$,
and $B(x_j,10r_j) \i V \i U$. Then use the usual $5$-covering
lemma to cover $\H^d(E\cap V)$ with a family 
$B(x_j,5r_j)$, $j\in J$, such that the $B(x_j,r_j)$ are disjoint.

For each finite subset $J_0$ of $J$, we deduce from (10.4) that
for $k$ large enough, every $B(x_j,r_j/2)$ contains a point $y_j \in E_k$.
Then
$$\eqalign{
\sum_{j\in J_0} \H^d(E\cap B(x_j,5r_j))
&\leq C \sum_{j\in J_0} r_j^d
\leq C \sum_{j\in J_0} \H^d(E_k\cap B(y_j,r_j/2))
\cr&
\leq C \sum_{j\in J_0} \H^d(E_k\cap B(x_j,r_j))
\leq C \H^d(E_k\cap V)
}\leqno (10.13)
$$
by (10.11) and (10.10). Thus 
$\sum_{j\in J_0} \H^d(E\cap B(x_j,5r_j))
\leq C\liminf_{k \to +\infty} \, \H^d(E_k \cap V)$
for each finite set $J_0 \i J$. We take the supremum,
observe that $\H^d(E\cap V) \leq \sum_{j\in J} \H^d(E\cap 
B(x_j,5r_j))$ because the $B(x_j,5r_j)$ cover $E\cap V$,
and get (10.12). A similar argument shows that
$$
\limsup_{k \to +\infty} \H^d(E_k \cap H)
\leq C_M \H^d(E \cap H)
\leqno (10.14)
$$
whenever $H$ is a compact subset of $U$, and where $C_M$ depends only
on $n$, $M$, and $\Lambda$. We skip the details, because
this is the same as (3.11) in [D2], 
and the proof applies here.

In the present context, we cannot really hope for
(10.14) to hold with $C_M = 1$:
$\H^d(E)$ could be smaller than the limit of the 
$\H^d(E_k)$ (even if it exists). For instance, with our assumptions,
$E_k$ could be the graph of $f_k(x) = 2^{-k} \cos(2^k x)$,
which is somewhat longer than its limit (a straight line).
Nonetheless, we shall see in Lemma 22.3 that (10.14) holds
with the more precise constant $C_M = (1+Ch)M$.

And in the more restricted context of almost-minimal sets, 
we will have a much better control on the upper semicontinuity defect,
and show that (10.14) holds holds with the optimal constant $C_M=1$.
See Theorem 22.1.

Surprisingly, both result only use very little information:
the rectifiability of $E^\ast$, a covering argument, 
and an application of the definition of quasiminimality.

Fortunately, the situation is different for (10.12),
for which $C_M$ can be removed, even for quasiminimal sets.
See Theorem 10.97 below, which is the main goal of this section.

\ms
We shall use the fact that $E$ is rectifiable. Notice that
in general, limits of rectifiable sets are not always rectifiable,
but in the present situation the proof of rectifiability given in 
Section 5 will kindly pass to  the limit.

\ms\proclaim Proposition 10.15.
Let $U$, $\{ E_k \}$, and $E$ satisfy the hypotheses (10.1), (10.2), 
(10.3), and (10.4). Also suppose that $h$ is small enough, 
depending only on $n$, $M$, and $\Lambda$.
Then $E$ is rectifiable. 

\ms
Of course the rectifiability of $E$ is compatible with
Theorem 10.8 and the fact that quasiminimal sets are 
rectifiable, but we need to prove it first.
If we are in a position to apply Theorem 6.1
or Theorem 9.81 to the $E_k$, we can deduce the rectifiability of $E$
from its uniform rectifiability (because unlike plain rectifiability,
uniform rectifiability (with uniform estimates) goes to the limit).
But this gives a much longer proof, and that does not even always 
works.

We start the proof like we did for Theorem 5.16.
Let $\lambda > 0$ and $\psi : \lambda U \to B(0,1)$
be as in Definition 2.7. Since the bilipschitz
mapping $\psi$ preserves rectifiable sets, it
is enough to show that $E' = \psi(\lambda E)$ is
rectifiable. But $E'$ is the limit (locally in 
$B(0,1) = \psi(\lambda U)$) of the sets
$E'_k = \psi(\lambda E_k)$, and by Proposition 2.8
each $E'_k$ lies in $GSAQ(B(0,1), \Lambda^{2d}M, 
\Lambda^{-1}\lambda\delta , \Lambda^{2d}h)$. So it is
enough to prove the proposition under the rigid assumption.

As before, we assume that $E$ is not rectifiable to get
a contradiction. Let $N_0 \geq 0$ be a large integer.
We first find an origin $x \in E_s$  
(the singular part of $E$) such that (5.17) holds, 
pick a cube $Q_0$ and an integer $N \in [N_0/2,N_0]$
such that (5.18) and (5.19) hold, 
and construct a mapping $\phi^\ast = \phi_d$ as in
(5.20)-(5.40).

Now we cannot use the fact that $E$ is quasiminimal,
but instead we shall apply the definition of 
quasiminimality to $E_k$ for $k$ large, with 
the same family $\{ \phi_t^\ast \}$, $0 \leq t \leq 1$,
as before. The $\{ \phi_t^\ast \}$ satisfy (1.4)-(1.8)
with respect to $E_k$ as well as $E$; the main point is
that (1.7) only uses the fact that the various mappings
$\psi_l$ that we composed to get $\phi^\ast$ all preserve 
every face of every cube $R \in {\cal R}(Q)$, as in (4.8), 
and the verification is still the same as below (4.16).  

We still can apply Definition 2.3, because of (5.43)
(which is the same for $E_k$ as for $E$), and we get that
$$
\H^d(E_k \cap W_1) \leq M \H^d(\phi^\ast(E_k \cap W_1)) + h r^d
\leqno (10.16)
$$
(compare with (5.44)). Next observe that
for $k$ large enough, (5.38) holds for $H = E_k \cap Q$, so 
we get (5.39). That is, 
$$
\phi_l(E_k \cap Q) \i {\cal S}_{l-1} \cup \d Q
\leqno (10.17)
$$
for $n+1 \geq l \geq d$, which will be a good
replacement for (5.21) and (5.40).

We may now continue as in Section 5, and get that
$$
\H^d(E_k \cap W_1)  
\geq H^d(E_k \cap {1 \over 3} Q_0) 
\geq C^{-1} l_0^{d}
\leqno (10.18)
$$
as in (5.49), and where the last inequality now follows 
from the first half of (5.18), and the facts that $E_k$ is 
locally Ahlfors-regular with uniform bounds and 
$\lim_{k \to +\infty} \dist(x,E_k) = 0$.

The proof of (5.53) also goes through (it just uses 
(5.21) and (5.40)), and yields
$$
\H^d(\phi^\ast(E_k \cap W_1)) 
\leq C \H^d(E_k \cap Q \sm {\rm int}(Q')),
\leqno (10.19)
$$
where the slightly smaller $Q'$ is defined in (5.47)
(we removed from $Q$ the exterior layer of small cubes).
Let us now apply (10.14) to $H = Q \sm {\rm int}(Q')$;
we get that for each $\varepsilon > 0$, 
$$
\H^d(E_k \cap Q \sm {\rm int}(Q'))
\leq C \H^d(E \cap Q \sm {\rm int}(Q')) + \varepsilon
\leq C {l_0^d \over N_0} + \varepsilon
\leqno (10.20)
$$
for $k$ large, and by the last part of (5.55).
We choose $\varepsilon$ small and $N_0$ large, as 
we did in Section~5, and get a contradiction with (10.16) 
or (10.18). Proposition 10.15 follows.
\qed

\ms
The following lemma is a slightly more uniform version of our
rectifiability result for quasiminimal sets; we shall deduce it 
from Proposition 10.15.

\ms\proclaim Lemma 10.21.
For each choice of $n$, $M \geq 1$, $\Lambda$, 
and $\varepsilon > 0$, we can find $h> 0$ and
$c_\varepsilon > 0$ such that if 
$U \i \R^n$ is open, $E \in GSAQ(U,M,\delta,h)$
is a  sliding quasiminimal set, and if
the pair $(x,r)$ is such that
$$
x\in E^\ast, \hskip0.2cm
0 < r < \Min(\lambda^{-1} r_0,\delta), \hskip0.1cm \hbox{and } 
B(x,2r) \i U,
\leqno (10.22)
$$
then we can find $y\in E^\ast \cap B(x,r/2)$, 
$t \in [c_\varepsilon r, r/2]$, and a $d$-plane $P$ through $y$, 
such that 
$$
\dist(z,P) \leq \varepsilon t
\ \hbox{ for } z\in E^\ast \cap B(y,t).
\leqno (10.23)
$$

\ms
Here $\Lambda$, $\lambda$, and $r_0$ are as in the definition of the
rigid and Lipschitz assumptions. Observe that we say that the 
constant $c_\varepsilon$ does not depend on $r_0$, $\delta$, $\lambda$,
or the precise list of boundary pieces $L_j$,
which is natural but will cost us a slightly more complicated
compactness argument.

The conclusion of Lemma 10.21 is a little more quantitative than
rectifiability, and could easily be deduced from the local uniform 
rectifiability of $E^\ast$ if we could prove it 
(we would use Lemma 7.8 and a small Chebyshev argument 
to find the pair $(y,t)$).
But it is not as strong as local uniform rectifiability,
for which we would need to know that for most pairs
$(y,t) \in E^\ast \cap B(x,r/2) \times (0,r/2]$
(in the sense that the complement satisfies a Carleson packing 
condition), we can find $P$ such that (10.23) holds.

The gap could possibly be filled, for instance
if we could prove a regularity result that says that
if $\varepsilon$ is small enough and (10.23) holds for the pair 
$(y,t)$, then it also holds for all pairs 
$(w,s) \in E^\ast \cap B(y,t/2) \times (0,t/2]$
(with a possibly different, but arbitrarily small $\varepsilon'$, 
and where the plane $P$ may depend on $(w,s)$).
Such a regularity result exists in the standard case
without boundaries, and could for instance be deduced from
Allard's theorem [All], 
but the proof is not easy.

Note also that because our proof will use a compactness 
argument, we shall not get any computable lower bound for 
$c_\varepsilon$.

\ms
Since Proposition 3.3 says that $E^\ast \in GSAQ(U,M,\delta,h)$, 
it will be enough to prove Lemma 10.21 when $E$ is coral, 
i.e., when $E^\ast = E$.

We shall assume that Lemma 10.21 is false and derive a contradiction.
Let $n$, $d$, $M$, $\Lambda$, and $\varepsilon > 0$ be such that 
the statement fails, and for each $k \geq 0$, choose a domain $U_k$, 
a coral set $E_k \in GSAQ(U_k,M,\delta_k,h_k)$, and a pair
$(x_k,r_k)$ that satisfy the hypotheses of the lemma, 
with $h_k = 2^{-k}$, but for which but we cannot find $y_k$, $t_k$, 
and $P_k$ such that
$$
y_k\in E_k \cap B(x_k,r_k/2), \, t_k \in [2^{-k}r_k, r_k/2],
\leqno (10.24)
$$
and
$$
\dist(z,P_k) \leq \varepsilon t_k
\hbox{ for } z\in E_k \cap B(y_k,t_k).
\leqno (10.25)
$$

The general scheme is quite simple: we want to take a limit,
obtain a limiting set which is rectifiable (by Proposition 10.15),
find a tangent plane to the limit at some point, and use it to get
a contradiction. But our quasiminimality assumption involves a
$\Lambda$-bilipschitz mapping $\psi_k : \lambda_k U_k \to B(0,1)$, 
a scale $\delta_k$, a choice of basic size $r_{0,k}$ for our dyadic grid, 
and a collection of boundary pieces $L_{j,k}$, that may all depend
on $k$; our argument will be a little complicated by the fact that
we shall need to make two or three changes of variables 
so as to be able to apply Proposition 10.15 to a fixed ball
and with fixed boundaries.

Here is the first change of variables. Set
$$
\wt E_k = \psi_k(\lambda_k E_k)
\ \hbox{ and } \
\wt x_k = \psi_k(\lambda_k x_k) \in \wt E_k
\leqno (10.26)
$$
(because $x_k \in E_k$), and observe that
$$
\dist(\wt x_k,\d B(0,1))
\geq \Lambda^{-1} \dist(\lambda_k x_k,\d (\lambda_k U))
\geq 2\Lambda^{-1} \lambda_k r_k
\leqno (10.27)
$$
because $\psi_k$ is $\Lambda$-bilipschitz
and by (10.22). Let $a_n$ denote the integer power of $2$ 
such that
$$
\sqrt n \leq a_n < 2 \sqrt n
\leqno (10.28)
$$
(we prefer to use dyadic numbers here), and
let $m_k$ be the integer such that
$$
2^{-m_k} \leq {\lambda_k r_k \over 8 a_n \Lambda} < 2^{-m_k+1}.
\leqno (10.29)
$$
Then choose a new origin $o_k \in (2^{-m_k} \Bbb Z)^n$ 
such that
$$
|o_k-\wt x_k| \leq 2^{-m_k}\sqrt n \leq 2^{-m_k} a_n
\leq {\lambda_k r_k \over 8 \Lambda}.
\leqno (10.30)
$$
Also set
$$
\gamma_k = 2^{-m_k+2} a_n \ \hbox{ and }
\wt B = B(o_k,\gamma_k).
\leqno (10.31)
$$
Observe that
$$
B(\wt x_k, 2^{-m_k} a_n) 
\i B(o_k, 2^{-m_k+1} a_n) 
= B(o_k, \gamma_k/2) 
\i \wt B
\leqno (10.32)
$$
and
$$
\wt B \i B(\wt x_k, 2^{-m_k+3} a_n)
\i B(\wt x_k, \Lambda^{-1}\lambda_k r_k )
\i B(0,1)
\leqno (10.33)
$$
(by (10.29) and (10.27)). Our second change of variable is 
the dilation $A_k$ given by
$$
A_k(z) = \gamma_k^{-1} \, (z-o_k)
= (2^{-m_k+2} a_n)^{-1} \, (z-o_k)
\ \hbox{ for } z\in \R^n
\leqno (10.34)
$$
which is of course chosen so that
$$
A_k(\wt B) = B(0,1).
\leqno (10.35)
$$
We set
$$
\wt E_k^\sharp = A_k(\wt E_k) \cap B(0,1)
= A_k(\wt E_k \cap \wt B)
= A_k(\psi_k(\lambda_k E_k) \cap \wt B)
\leqno (10.36)
$$
and
$$
\wt x_k^\sharp = A_k(\wt x_k) = A_k(\psi_k(\lambda_k x_k))
\in \wt E_k^\sharp
\leqno (10.37)
$$
(because $\wt x_k \in \wt E_k \cap \wt B$, and by (10.26)).

Let us check that $\wt E^\sharp$ is quasiminimal.
Recall that $E_k \in GSAQ(U_k,M,\delta_k,h_k)$; then
Proposition 2.8 says that
$$
\wt E_k = \psi_k(\lambda_k E_k) \in
GSAQ(B(0,1),\wt M, \wt \delta_k, \wt h_k),
\leqno (10.38)
$$
with $\wt M = \Lambda^{2d} M$, 
$\wt \delta_k = \Lambda^{-1} \lambda_k \delta_k$, 
and $\wt h_k = \Lambda^{2d} h_k$. For $\wt E_k$, we have 
the rigid assumption, and the boundary sets are the
$\wt L_{j,k} = \psi_k(\lambda_k L_{j,k})$, 
which by assumption are composed of faces of dyadic cubes 
of side length $r_{0,k}$.

We claim that $\wt E_k^\sharp$ is also quasiminimal,
and more precisely that
$$
\wt E_k^\sharp \in GSAQ(B(0,1),\wt M, \wt\delta_k^\sharp, \wt h_k),
\ \hbox{ with } 
\wt\delta_k^\sharp = \gamma_k^{-1} \wt \delta_k.
\leqno (10.39)
$$
This comes directly from the definitions, using the fact that by 
(10.38) $\wt E_k \cap \wt B$ is quasiminimal in $\wt B \i B(0,1)$. 
We just need to multiply $\wt \delta_k$ with the dilation factor
$\gamma_k^{-1}$; the other constants stay the same when we dilate 
everything. Also, the boundary constraints are now given by the sets
$$
\wt L_{j,k}^{\sharp} = A_k(\wt L_{j,k}) \cap B(0,1)
= A_k(\psi_k(\lambda_k L_{j,k})) \cap B(0,1).
\leqno (10.40)
$$
Let us verify that the $A_k(\wt L_{j,k})$ lie in
an acceptable grid. As was said above, the 
$\wt L_{j,k}$ are composed of faces of dyadic cubes 
of side length $r_{0,k}$. Our center $o_k$ lies in 
$(2^{-m_k} \Bbb Z)^n$, and
$$
2^{-m_k} \leq {\lambda_k r_k \over 8 a_n \Lambda}
\leq {r_{0,k} \over 8 a_n \Lambda} \leq r_{0,k}
\leqno (10.41)
$$
by (10.29) and (10.22), so the sets $\wt L_{j,k} - o_k$ 
are composed of faces of dyadic cubes of size $2^{-m_k}$. 
By (10.34), $A_k$ maps our cubes into dyadic cubes of 
side length $(2^{-m_k+2} a_n)^{-1} 2^{-m_k} = (4a_n)^{-1}$
That is, for $\wt E_k^\sharp$ we have the rigid assumption, 
with a scale $\wt r_0^\sharp = (4a_n)^{-1}$ that does not
depend on $k$.

This is good, because there is only a finite number of 
possibilities for the $\wt L_{j,k}^{\sharp}$, and modulo 
replacing $\{ E_k \}$ with a subsequence we may assume that the 
$\wt L_{j,k}^{\sharp}$ are always the same. 
The scale constant $\wt\delta_k^\sharp$ will not create trouble
either, because
$$\eqalign{
\wt\delta_k^\sharp &= \gamma_k^{-1} \wt\delta_k
= \gamma_k^{-1} \Lambda^{-1} \lambda_k \delta_k
= (2^{-m_k+2} a_n)^{-1} \Lambda^{-1} \lambda_k \delta_k
\cr&
\geq (4a_n)^{-1} \, {8 a_n \Lambda \over \lambda_k r_k}
\, \Lambda^{-1} \lambda_k \delta_k
= {2\delta_k \over r_k} \geq 2
}\leqno (10.42)
$$
by (10.31), (10.29), and our assumption (10.22).
Altogether (10.39) simplifies into
$$
\wt E_k^\sharp \in GSAQ(B(0,1),\wt M, 2, \wt h_k),
\hbox{ with $\wt M = \Lambda^{2d} M$, $\wt h_k = \Lambda^{2d} h_k$,}
\leqno (10.43)
$$
with the rigid assumption at the fixed scale $(4a_n)^{-1}$ and 
with a fixed set of boundaries.

Let us replace $\{ E_k \}$ with a new subsequence,
so that 
$$
\wt E_k^\sharp \hbox{ converges to a limit } \wt E_\infty^\sharp
\leqno (10.44)
$$
(locally in $B(0,1)$, as in (10.4)-(10.6)).
Recall that in (10.43), 
$\wt h_k = \Lambda^{2d} h_k$ tends to $0$,
because we assumed that $h_k = 2^{-k}$.
Thus we can apply Proposition 10.15 and get that
$$
\wt E_\infty^\sharp \hbox{ is rectifiable.}
\leqno (10.45)
$$

Then extract again a subsequence, so that
$\wt x_k^\sharp$ tends to some limit $\wt x^\sharp_{\infty}$.
Recall that $\wt x_k^\sharp = A_k(\wt x_k)\in \wt E_k^\sharp$ 
is defined in (10.37), and let us check that
$$
\wt x^\sharp_{\infty} = \lim_{k \to +\infty} \wt x_k^\sharp
\in \wt E_\infty^\sharp \cap \overline B(0,1/4).
\leqno (10.46)
$$
First of all,
$$\eqalign{
|\wt x_k^\sharp| &= |A_k(\wt x_k)| = |A_k(\wt x_k)-A_k(o_k)|
= (2^{-m_k+2}a_n)^{-1} |\wt x_k-o_k|
\cr&
\leq (2^{-m_k+2}a_n)^{-1} \, 2^{-m_k}a_n = {1 \over 4}
}\leqno (10.47)
$$
by (10.34) and (10.30). So $\wt x_k^\sharp \in \overline B(0,1/4)$,
hence $\wt x^\sharp_{\infty} = \lim_{k \to +\infty} \wt x_k^\sharp
\in \overline B(0,1/4)$ too.
Finally, $\wt x^\sharp_{\infty} \in \wt E_\infty^\sharp$
because $\wt x_k^\sharp \in \wt E_k^\sharp$,
$\wt E_k^\sharp$ converges to $\wt E_\infty^\sharp$ locally
in $B(0,1)$, and our points stay in $\overline B(0,1/4)$.
This proves (10.46).

\ms
We need a last change of coordinates.
Recall that $A_k^{-1}(B(0,1)) = \wt B \i B(0,1)$
by (10.35) and (10.33), so we can define $\theta_k$ on $B(0,1)$ by
$$
\theta_k(z) = \gamma_k^{-1} \psi_k^{-1}(A_k^{-1}(z)).
\leqno (10.48)
$$

Observe that $\theta_k$ is a $\Lambda$-bilipschitz mapping
(we conjugate $\psi_k^{-1}$ by a dilation, and translate). 
We can extend $\theta_k$
to the closed ball $\overline B(0,1)$, either because it is Lipschitz
on $B(0,1)$, or simply using (10.48) and the fact that 
$A_k^{-1}(\overline B(0,1)) \i B(0,1)$. By Arzel\`{a}-Ascoli,
the collection of mappings 
$\theta_k(\cdot) - \theta_k(0) : \overline B(0,1) \to \R^n$
is totally bounded, so we can extract another subsequence so that
$$
\hbox{the functions $\theta_k(\cdot) - \theta_k(0)$ converge,
uniformly on $\overline B(0,1)$, to a limit $\theta_\infty$.}
\leqno (10.49)
$$
We need to remove the constant $\theta_k(0)$ 
because the domains $U_k$ could go away to infinity very fast;
we could also avoid this minor problem by translating $U_k$, $E_k$,
and $x_k$ by a same vector $v_k$, with the effect of
precomposing $\psi_k$ with a translation, adding 
$\lambda_k v_k$ to $\psi_k^{-1}$, keeping $\wt E_k$ and 
$A_k$ as they are, adding $\gamma_k^{-1} \lambda_k v_k$
to $\theta_k$, and making $\theta_k(0) = 0$
(if $v_k$ is chosen correctly).

Note that $\theta_\infty$ is $\Lambda$-bilipschitz on $B(0,1)$, 
because the $\theta_k$ are. Next set 
$$
E_k^\sharp = \theta_k(\wt E_k^\sharp), \
E_\infty^\sharp = \theta_\infty(\wt E_\infty^\sharp),
\, \hbox{ and } \, 
x_\infty^\sharp = \theta_\infty(\wt x_\infty^\sharp)
\in E_\infty^\sharp
\leqno (10.50)
$$
(because $\wt x_\infty^\sharp \in \wt E_\infty^\sharp$
by (10.46)).
We expect correct translations of $E_k^\sharp$
to converge to $E_\infty^\sharp$, but in fact we shall
find it more convenient to work directly on $\wt E_k^\sharp$ 
and $\wt E_\infty^\sharp$. Obviously, 
$$
E_\infty^\sharp \hbox{ is rectifiable,}
\leqno (10.51)
$$
by (10.45) and because $\theta_\infty$ is bilipschitz.
We need the last change of variable because we want $E_k$ to be close 
to a plane, not $\wt E_k$, and for this (10.51) will be more useful.

Let $\varepsilon^\sharp$ be very small (to be chosen near the end);
we shall use $\varepsilon^\sharp$ to measure various small
quantities that are not necessarily connected to each other.
We claim that we can choose 
$$
y_\infty^\sharp \in E_\infty^\sharp \cap B(x_\infty^\sharp,\varepsilon^\sharp) 
\leqno (10.52)
$$
so that
$$
E_\infty^\sharp \hbox{ has a tangent plane $P^\sharp$
at } y_\infty^\sharp.
\leqno (10.53)
$$
Indeed, $E_\infty^\sharp$ is rectifiable by (10.51)
so it has an approximate tangent plane at $\H^d$-almost
every point. In addition, recall from (10.44) that 
$\wt E_\infty^\sharp$ is the local limit in $B(0,1)$
of a sequence of reduced quasiminimal sets $\wt E_k^\sharp$,
with uniform constants (see (10.43)). When $k$ is large,
$\wt h_k = \Lambda^{2d} h_k = \Lambda^{2d} 2^{-k}$
is as small as we want; then the assumptions (10.1)-(10.4) 
are satisfied, and so $\wt E_\infty^\sharp$ is locally Ahlfors-regular
in $B(0,1)$, as in (10.11). Its bilipschitz 
image $E_\infty^\sharp = \theta_\infty(\wt E_\infty^\sharp)$
is also locally Ahlfors-regular. Then the
approximate tangent planes are automatically 
tangent planes (see for instance Exercise 41.21 on
page 277 of [D4]),  
and so $E_\infty^\sharp$ admits a (true) tangent plane at
$\H^d$-almost every point $y_\infty^\sharp$. Finally, 
$\H^d(E_\infty^\sharp \cap B(x_\infty^\sharp,\varepsilon^\sharp)) > 0$
(again by local Ahlfors-regularity, and because 
$x_\infty^\sharp \in E_\infty^\sharp$),
so we can choose $y_\infty^\sharp$ in
$B(x_\infty^\sharp,\varepsilon^\sharp)$, as needed.

By (10.53), we can find $\rho \in (0,1)$ such that 
$$
\dist(z,P^\sharp) \leq \varepsilon^\sharp \rho
\ \hbox{ for } 
z\in E_\infty^\sharp \cap B(y_\infty^\sharp,\rho).
\leqno (10.54)
$$
We want to use $y_\infty^\sharp$ and $P^\sharp$ to find,
for $k$ large enough, a pair $(y_k,t_k)$ that satisfies
(10.24) and (10.25); the desired contradiction will
ensue.

Since $y_\infty^\sharp \in E_\infty^\sharp = \theta_\infty(\wt E_\infty^\sharp)$
(by (10.52) and (10.50)), there exists $\wt y_\infty^\sharp \in \wt E_\infty^\sharp$
such that $\theta_\infty(\wt y_\infty^\sharp) = y_\infty^\sharp$.
Notice that 
$$
|\wt y_\infty^\sharp-\wt x_\infty^\sharp|
\leq \Lambda |y_\infty^\sharp - x_\infty^\sharp|
\leq \Lambda \varepsilon^\sharp
\leqno (10.55)
$$
because $\theta_\infty$ is $\Lambda$-bilipschitz,
$x_\infty^\sharp = \theta_\infty(\wt x_\infty^\sharp)$
by (10.50), and by (10.52). Then 
$\wt y_\infty^\sharp \in \wt E_\infty^\sharp \cap B(0,1/2)$ by (10.46), 
and by (10.44) we can find points $\wt y_k ^\sharp \in \wt E_k^\sharp$ 
so that
$$
\wt y_\infty^\sharp = \lim_{k \to +\infty} \wt y_k ^\sharp .
\leqno (10.56)
$$
Set $y_k^\sharp = \theta_k(\wt y_k ^\sharp) \in E_k^\sharp$
(by (10.50)), 
$\wt y_k = A_k^{-1}(\wt y_k ^\sharp) \in \wt E_k \cap \wt B$
(by (10.36) and (10.35)), and
$y_k = \lambda_k^{-1}\psi_k^{-1}(\wt y_k) \in E_k$ 
(by (10.33) and (10.26)). By (10.48),
$$
y_k^\sharp = \theta_k(\wt y_k ^\sharp) =
\gamma_k^{-1} \psi_k^{-1}(A_k^{-1}(\wt y_k ^\sharp))
= \gamma_k^{-1} \lambda_k y_k.
\leqno (10.57)
$$
For (10.24) we first need to check that for $k$ large,
$$
y_k \in E_k \cap B(x_k,r_k/2).
\leqno (10.58)
$$
But $y_k = \lambda_k^{-1} \psi_k^{-1}(A_k^{-1}(\wt y_k ^\sharp))$,
and similarly $x_k = \lambda_k^{-1}\psi_k^{-1}(A_k^{-1}(\wt x_k ^\sharp))$
by (10.37), so 
$$\eqalign{
|y_k-x_k| &= |\lambda_k^{-1} \psi_k^{-1}(A_k^{-1}(\wt y_k ^\sharp))
- \lambda_k^{-1}\psi_k^{-1}(A_k^{-1}(\wt x_k ^\sharp))|
\cr&
\leq \lambda_k^{-1} \Lambda |A_k^{-1}(\wt y_k ^\sharp)-A_k^{-1}(\wt x_k ^\sharp)|
= \lambda_k^{-1} \Lambda 2^{-m_k+2} a_n |\wt y_k ^\sharp - \wt x_k ^\sharp|
\cr&
\leq 4 \lambda_k^{-1} \Lambda \, {\lambda_k r_k \over 8 a_n \Lambda}
\, a_n |\wt y_k ^\sharp - \wt x_k ^\sharp|
= {r_k \over 2} \, |\wt y_k ^\sharp - \wt x_k ^\sharp|
}\leqno (10.59)
$$
because $\psi_k$ is $\Lambda$-bilipschitz, and 
by (10.34) and (10.29). Now $\wt y_k ^\sharp$ tends 
to $\wt y_\infty ^\sharp$ by (10.56), $\wt x_k ^\sharp$ 
tends to $\wt x_\infty ^\sharp$ by (10.46), and so
$$
|\wt y_k ^\sharp - \wt x_k ^\sharp|
\leq |\wt y_\infty ^\sharp - \wt x_\infty ^\sharp|
+ |\wt y_k ^\sharp -\wt y_\infty ^\sharp|
+ |\wt x_k ^\sharp- \wt x_\infty ^\sharp|
\leq |\wt y_\infty ^\sharp - \wt x_\infty ^\sharp| + \varepsilon^\sharp
\leq (\Lambda+1) \varepsilon^\sharp
\leqno (10.60)
$$
for $k$ large, and by (10.55). By (10.59),
$$
|y_k-x_k| \leq {(\Lambda+1) \varepsilon^\sharp \, r_k \over 2} 
\leqno (10.61)
$$
for $k$ large, and of course (10.58) follows.
We choose 
$$
\dsp t_k = {\rho r_k \over 20 \Lambda^2}
\leq {r_k \over 20\Lambda^2}, 
\leqno (10.62)
$$
with $\rho \in (0,1)$ as in (10.54). 
Obviously the pair $(y_k,t_k)$ satisfies (10.24)
for $k$ large, so we just need to find a $d$-plane
$P_k$ through $y_k$ such that (10.25) holds.

\ms
So let $z\in E_k \cap B(y_k,t_k)$ be given. 
We first want to define
$$
\wt z = \psi_k(\lambda_k z) \in \wt E_k \, ,
\ 
\wt z^\sharp = A_k(\wt z) \in \wt E_k^\sharp
\, \hbox{ and } \,
z^\sharp = \theta_k(\wt z^\sharp) \in E_k^\sharp \, ,
\leqno (10.63)
$$
because (10.54) gives us some control on
$E_\infty^\sharp$ and hence probably on $E_k^\sharp$.

We start with $\wt z = \psi_k(\lambda_k z)$,
which is defined because $z \in B(y_k,t_k) \i U_k$
(recall that $B(y_k,t_k) \i B(x_k,r_k)$ by
(10.61) and (10.62), and use (10.22)).
Next, $\wt z \in \wt E_k$ by (10.26).
Before we switch to $\wt z^\sharp$, observe that
$$\eqalign{
|\wt z - \wt x_k| 
&= |\psi_k(\lambda_k z)-\psi_k(\lambda_k x_k)|
\leq \Lambda |\lambda_k z - \lambda_k x_k|
\leq \Lambda \lambda_k (|z-y_k|+|y_k - x_k|)
\cr&
\leq \Lambda \lambda_k 
(t_k + {(\Lambda+1) \varepsilon^\sharp \, r_k \over 2})
\leq {\lambda_k r_k \over 19 \Lambda} 
< 2^{-m_k} a_n
}\leqno (10.64)
$$
by (10.26), (10.61), if $\varepsilon^\sharp$ is small enough,
and by (10.62) and (10.29). Thus 
$$
\wt z \in B(\wt x_k,2^{-m_k} a_n) 
\i B(o_k, \gamma_k/2) = {1 \over 2} \wt B
\leqno (10.65)
$$
by (10.32) and (10.31). By (10.34) or (10.35), 
$A_k(B(o_k, \gamma_k/2)) = B(0,1/2)$, so
$$
\wt z^\sharp = A_k(\wt z) \in \wt E_k^\sharp \cap B(0,1/2),
\leqno (10.66)
$$
because (10.36) says that $\wt E_k^\sharp = A_k(\wt E_k \cap \wt B)$.
Now $z^\sharp = \theta_k(\wt z^\sharp)$ is defined and lies in $E_k^\sharp$
by the definition (10.50) and because $\wt z^\sharp \in \wt E_k^\sharp$.
This completes our verification of (10.63).

Notice that
$$
z^\sharp = \theta_k(\wt z^\sharp) = 
\gamma_k^{-1} \psi_k^{-1}(A_k^{-1}(\wt z^\sharp))
= \gamma_k^{-1} \psi_k^{-1}(\wt z)
= \gamma_k^{-1} \lambda_k z
\leqno (10.67)
$$
by (10.63) and (10.48). Set
$$
\varepsilon_k = 
\sup_{w \in \wt E_k^\sharp \cap \overline B(0,1/2)}
\dist(w,\wt E_\infty^\sharp);
\leqno (10.68)
$$
then $\varepsilon_k$ tends to $0$ because 
$\overline B(0,1/2)$ is a compact subset of $B(0,1)$,
and by the local convergence of $\wt E_k^\sharp$ to $\wt E_\infty^\sharp$
(see (10.44)). Similarly,
$$
\varepsilon'_k = 
\sup_{w\in B(0,1)} |\theta_k(w)-\theta_k(0)-\theta_\infty(w)|
\leqno (10.69)
$$
tends to $0$, by the uniform convergence in (10.49).

Return to $\wt z^\sharp$, and choose 
$\wt \xi^\sharp \in \wt E_\infty^\sharp$ such that 
$$
|\wt \xi^\sharp-\wt z^\sharp|
\leq \dist(\wt z^\sharp,\wt E_\infty^\sharp) \leq \varepsilon_k
\leqno (10.70)
$$
(compare (10.68) with (10.66)). Set
$\xi^\sharp = \theta_\infty(\wt \xi^\sharp)$; then 
$$
\xi^\sharp \in E_\infty^\sharp = \theta_\infty(\wt E^\sharp_\infty)
\leqno (10.71)
$$
(see (10.50)). Also,
$$\eqalign{
|\xi^\sharp - y_\infty^\sharp|
&= |\theta_\infty(\wt \xi^\sharp)- \theta_\infty(\wt y_\infty^\sharp)|
\leq |\theta_k(\wt \xi^\sharp)- \theta_k(\wt y_\infty^\sharp)|
+ 2\varepsilon'_k
\cr&
\leq |\theta_k(\wt z^\sharp)- \theta_k(\wt y_k^\sharp)|
+ \Lambda \big(|\wt \xi^\sharp-\wt z^\sharp| 
+ |\wt y_\infty^\sharp - \wt y_k^\sharp|\big)
+ 2\varepsilon'_k
\cr&
\leq |\theta_k(\wt z^\sharp)- \theta_k(\wt y_k^\sharp)|
+ \Lambda \varepsilon_k 
+ \Lambda |\wt y_\infty^\sharp - \wt y_k^\sharp|
+ 2\varepsilon'_k
\cr&
\leq |\theta_k(\wt z^\sharp)- \theta_k(\wt y_k^\sharp)|
+ \varepsilon^\sharp \rho
}\leqno (10.72)
$$
because $y_\infty^\sharp = \theta_\infty(\wt y_\infty^\sharp)$
(see above (10.55)) and $\theta_k$ is $\Lambda$-Lipschitz,
by (10.70), then by (10.56), because $\varepsilon_k$ and  
$\varepsilon'_k$ tend to $0$, and if $k$ is large enough.
Next
$$
\theta_k(\wt z^\sharp)- \theta_k(\wt y_k^\sharp)
= \gamma_k^{-1} \lambda_k (z - y_k)
\leqno (10.73)
$$
by (10.67) and (10.57); since
$$
|z - y_k| \leq t_k = {\rho r_k \over 20 \Lambda^2}
\leqno (10.74)
$$
because $z\in E_k \cap B(y_k,t_k)$ and by (10.62),
we deduce from (10.73) that
$$\eqalign{
|\theta_k(\wt z^\sharp)- \theta_k(\wt y_k^\sharp)|
&\leq \gamma_k^{-1} \lambda_k \, {\rho r_k \over 20 \Lambda^2}
= (2^{-m_k+2}a_n)^{-1} \lambda_k \, {\rho r_k \over 20 \Lambda^2}
\cr&
\leq {1 \over 2} \, {8a_n \Lambda \over \lambda_k r_k} \, a_n^{-1} \lambda_k
\, {\rho r_k \over 20 \Lambda^2}
= {\rho \over 5\Lambda}
}\leqno (10.75)
$$
by (10.31) and (10.29). Therefore
$$
|\xi^\sharp - y_\infty^\sharp|
\leq |\theta_k(\wt z^\sharp)- \theta_k(\wt y_k^\sharp)|
+ \varepsilon^\sharp \rho
\leq {\rho \over 5\Lambda} + \varepsilon^\sharp \rho
 < {\rho \over 4 \Lambda} 
\leqno (10.76)
$$
by (10.72) and (10.75), and if $\varepsilon^\sharp$ is small enough. 
Since $\xi^\sharp \in E_\infty^\sharp$ by (10.71), we get that
$\xi^\sharp \in E_\infty^\sharp \cap B(y_\infty^\sharp,\rho)$,
and (10.54) says that 
$\dist(\xi^\sharp,P^\sharp) \leq \varepsilon^\sharp \rho$.
Set 
$$
P'_k = \gamma_k \lambda_k^{-1}[ P^\sharp + \theta_k(0)]; 
\leqno (10.77)
$$
then
$$\eqalign{
\dist(z,P'_k) 
&= \gamma_k \lambda_k^{-1} 
\dist(\gamma_k^{-1} \lambda_k z,P^\sharp + \theta_k(0))
= \gamma_k \lambda_k^{-1} \dist(z^\sharp,P^\sharp + \theta_k(0))
\cr&
\leq \gamma_k \lambda_k^{-1}
\big[\dist(\xi^\sharp + \theta_k(0),P^\sharp+ \theta_k(0)) 
+ |\xi^\sharp + \theta_k(0)-z^\sharp| \big]
\cr&
= \gamma_k \lambda_k^{-1}
\big[\dist(\xi^\sharp,P^\sharp) + |\xi^\sharp + \theta_k(0)-z^\sharp| \big]
\cr&
\leq \gamma_k \lambda_k^{-1} \big[ \varepsilon^\sharp \rho
+ |\xi^\sharp +\theta_k(0)-z^\sharp| \big]
}\leqno (10.78)
$$
by (10.67). Recall that $\xi^\sharp = \theta_\infty(\wt \xi^\sharp)$
(see above (10.71)) and $z^\sharp = \theta_k(\wt z^\sharp)$
(see (10.63)), so
$$\eqalign{
|\xi^\sharp+\theta_k(0)-z^\sharp|
&= |\theta_\infty(\wt \xi^\sharp)+\theta_k(0)-\theta_k(\wt z^\sharp)|
\cr&
\leq |\theta_\infty(\wt \xi^\sharp)+\theta_k(0)-\theta_k(\wt \xi^\sharp)|
+ |\theta_k(\wt \xi^\sharp)-\theta_k(\wt z^\sharp)|
\cr&
\leq \varepsilon'_k + \Lambda |\wt \xi^\sharp-\wt z^\sharp|
\leq \varepsilon'_k + \Lambda \varepsilon_k
}\leqno (10.79)
$$
by (10.69) and (10.70). We combine this with (10.78) and get that
$$\eqalign{
\dist(z,P'_k) &\leq \gamma_k \lambda_k^{-1}
\big[ \varepsilon^\sharp \rho
+ \varepsilon'_k + \Lambda \varepsilon_k]
= 2^{-m_k+2} a_n \lambda_k^{-1} \big[ \varepsilon^\sharp \rho
+ \varepsilon'_k + \Lambda \varepsilon_k]
\cr&
\leq 4 \, {\lambda_k r_k \over 8 a_n \Lambda}
a_n \lambda_k^{-1} \big[ \varepsilon^\sharp \rho
+ \varepsilon'_k + \Lambda \varepsilon_k]
= {r_k \over 2  \Lambda} \,
\big[ \varepsilon^\sharp \rho
+ \varepsilon'_k + \Lambda \varepsilon_k]
\cr&
= {10 \Lambda t_k \over \rho}
\big[ \varepsilon^\sharp \rho
+ \varepsilon'_k + \Lambda \varepsilon_k]
}\leqno (10.80)
$$
by (10.31), (10.29), and (10.62).
Recall that $\varepsilon_k$ and $\varepsilon'_k$ tend to $0$;
then, for $k$ large enough (depending also on $\rho$,
but this is all right), we deduce from (10.80) that
$$
\dist(z,P'_k) \leq 20 \Lambda \varepsilon^\sharp t_k
\leq {\varepsilon t_k \over 2}
\leqno (10.81)
$$
if $\varepsilon^\sharp$ is chosen small enough.

Now all this is true for all $k$ large (not depending on $z$), 
and all $z\in E_k \cap B(y_k,t_k)$. In particular,
$z=y_k$ yields $\dist(y_k,P'_k) \leq {\varepsilon t_k \over 2}$.
We choose for $P_k$ a translation of $P'_k$ that contains
$y_k$; this is required for (10.23) and (10.25),
but fortunately we just need to translate by at most
${\varepsilon t_k \over 2}$. Then
$\dist(y_k,P_k) \leq \varepsilon t_k$
for $z\in E_k \cap B(y_k,t_k)$, by (10.81)
and as needed for (10.25). 

We finally found a plane $P_k$ through $y_k$ that satisfies 
(10.25); as announced earlier, its existence contradicts the definition
of our sequence $\{ E_k \}$; this completes our proof of Lemma 10.21 
by contradiction.
\qed

\ms
Our next preparatory result is a (simplified) generalization of 
Corollaries 9.103 and 8.55; it says that if $E$ is a quasiminimal
set, its core $E^\ast$ is concentrated, with uniform bounds.
The terminology comes from [DMS]  
(and is justified by that fact that (10.84) below says that
$E^\ast$ is almost as concentrated in $B(y,t)$ as a $d$-plane
through $y$), and the result is interesting because it will soon 
allow us to prove the lower semicontinuity of $\H^d$ along 
convergent sequences of uniformly quasiminimal sets.

\ms
\ms\proclaim Proposition 10.82. 
For each choice of constants $n$, $M \geq 1$, 
$\Lambda \geq 1$ and $\varepsilon > 0$, we can find $h > 0$ and 
$d_\varepsilon > 0$ such that the following holds. 
Suppose that $E \in GSAQ(U,M,\delta,h)$ for some open set
$U \i \R^n$, and that the Lipschitz assumption are satisfied,
with the constants $\lambda$ and $\Lambda$ (as in (9.3)).
Also denote by $r_0 = 2^{-m} \leq 1$ the side length of the 
dyadic cubes of the usual grid. Then let $(x,r)$ be such that
$$
x\in E^\ast, \hskip0.2cm
0 < r < \Min(\lambda^{-1} r_0,\delta), \hskip0.2cm 
B(x,2r) \i U.
\leqno (10.83)
$$
Then we can find a pair $(y,t)$, such that
$y \i E^\ast \cap B(x,r/2)$, $d_\varepsilon r \leq t \leq r/4$,
and
$$
\H^d(E^\ast\cap B(y,t)) \geq (1-\varepsilon) \omega_d t^d,
\leqno (10.84)
$$
where $\omega_d$ denotes the $d$-dimensional Hausdorff measure
of the unit ball in $\R^d$.

\ms
The proof will be similar to the proof of 
Corollaries 9.103 and 8.55, but we shall rely on Lemma 10.21
rather than the uniform rectifiability of our quasiminimal
sets, which we do not know how to prove with enough
generality. This slightly different approach is new,
and even in the case of standard quasiminimal
sets without boundaries, it has the advantage of
not using our complicated proof of uniform rectifiability
(with the tough stopping time argument on the projections).
But one more compactness argument is used, and we loose
an ``explicit" control on $d_\varepsilon$.

Compared with Corollary 9.103, we just get rid of the 
unpleasant additional assumption (9.105) on the dimensions
of some faces. Recall that Corollary 8.55 works under
the rigid assumption, and also has the unpleasant dimensionality
assumption (6.2).

For the proof, first notice that by Proposition 3.3,
$E^\ast \in GSAQ(U,M,\delta,h)$, so it is enough to prove
Proposition 10.82 when $E$ is coral, i.e., when $E^\ast = E$.

Let $E$, $x$, and $r$ be as in the statement.
Our goal is to apply Lemma 9.14 to some pair $(y,t)$,
because surjective projections will help us find
lower bounds on $\H^d(E\cap B(y,t))$.
Let $C_0 \geq 1$, $\eta >0$, and $\overline \varepsilon >0$
be as in the statement of that lemma, and recall that
they depend only on $n$, $M$, and $\Lambda$.

Let $\varepsilon_0$ be very small, to be chosen near
the end, and apply Lemma 10.21, with the constant 
$\varepsilon_0$, to the pair $(x,(2C_0)^{-1}r)$;
the hypotheses for Lemma 10.21 are the same as for the present
proposition, so the pair $(x,(2C_0)^{-1}r)$ satisfies them. 
We get $(y_0,t_0)$ such that
$$
y_0\in E \cap B(x,{r \over 4C_0})
\ \hbox{ and } \
{c_{\varepsilon_0} r \over 2C_0} \leq t_0 \leq {r \over 4C_0},
\leqno (10.85)
$$
and a plane $P$ through $y_0$ such that
$$
\dist(z,P) \leq \varepsilon_0 t_0
\ \hbox{ for } z\in E \cap B(y_0,t_0).
\leqno (10.86)
$$
We would be happy to apply Lemma 9.14 directly to $(y_0,t_0)$,
but the unpleasant assumption (9.17) on the proximity to some
boundaries $L_j$ may not be satisfied. As in (9.16), set
$$
J(y,t) = \big\{ j\in[0,j_{max}] \, ; \,  L_j \hbox{ meets } B(y,2t) \big\}
\ \ \hbox{ and } \ \ 
L(y,t) = \bigcap_{j \in J(y,t)} L_j
\leqno (10.87)  
$$
for $y\in E$ and $t > 0$. 
Recall that $J(y,t) \neq \emptyset$ because $y\in E \i L_0 = \Omega$.
We want to find pairs $(y,t)$ such that
$$
\dist(w,L(y,t)) \leq \eta t
\ \hbox{ for } w\in E \cap B(y,2t),
\leqno (10.88) 
$$
as in (9.17). We shall restrict to pairs $(y,t)$ such that
$B(y,2t) \i B(y_0,t_0)$, near which (10.86) says that $E$
stays very close to $P$.

We shall define a (finite) sequence of pairs 
$(y_k,t_k)$. Naturally, we start with $(y_0,t_0)$.

Suppose we already defined $(y_k,t_k)$.
If the pair $(y_k,t_k/2)$ satisfies (10.88),
we stop the construction. Otherwise, we define 
$(y_{k+1},t_{k+1})$ as follows.

If $J(y_k,t_k/4) \neq J(y_k,t_k)$ (which by (10.87)
means that it is strictly smaller), set
$(y_{k+1},t_{k+1}) = (y_k,t_k/4)$. We are happy, because
$J(y_{k+1},t_{k+1})$ is strictly contained in $J(y_k,t_k)$
and this cannot happen too often.

If $J(y_k,t_k/4) = J(y_k,t_k)$ and
the pair $(y_k,t_k/4)$ satisfies (10.88),
set $(y_{k+1},t_{k+1}) = (y_k,t_k/2)$. This time
we are happy too because we know that we will 
stop next time.

In the remaining case, the failure of (10.88)
for $(y_k,t_k/4)$ gives a point $w\in E \cap B(y_k,t_k/2)$ 
such that 
$$
\dist(w,L(y_k,t_k)) =
\dist(w,L(y_k,t_k/4)) \geq \eta t_k/4, 
\leqno (10.89) 
$$
where the first identity comes from the fact
that $J(y_k,t_k/4) = J(y_k,t_k)$. In this last case, we set 
$(y_{k+1},t_{k+1}) = (w, a t_k)$, for some small constant
$a \in (0,1/2)$ that will be chosen soon. Notice that
$B(y_{k+1},t_{k+1}) \i B(y_k,t_k)$ in all cases, 
so that we know that 
$$
B(y_k,t_k) \i B(y_0,t_0)
\ \hbox{for all $k \geq 1$.}
\leqno (10.90) 
$$

This completes our definition of the pairs
$(y_k,t_k)$. Now we want to show that the construction
stops after at most $j_{max} + 2$ steps
(where $j_{max}+1$ still denotes the number of boundary 
sets $L_j$), and for this it will be enough to show that
in our last case, 
$$
J(w,a t_k) \hbox{ is strictly contained in } J(y_k,t_k).
\leqno (10.91) 
$$
Notice that $J(w,at_k) \i J(y_k,t_k)$ by (10.87) and because 
$B(w, 2a t_k) \i B(w, t_k) \i B(y_k,2t_k)$, so (10.91) just
means that $J(w,a t_k) \neq J(y_k,t_k)$.
Let us suppose that (10.91) fails, and derive a contradiction. 

Then $J(w,a t_k) = J(y_k,t_k)$, and
for each $j \in J(y_k,t_k)$, the definition (10.87) says
that $L_j$ meets $B = B(w,2at_k)$. Choose a face
$F_j \i L_j$ that meets $B$, and set
$$
F = \bigcap_{j \in J(y_k,t_k)} F_j \i L(y_k,t_k),
\leqno (10.92) 
$$
where the inclusion comes from (10.87).
We want to show that $F$ is nonempty and meets a larger ball. 
We return to the standard grid because this will make the computations easier.

Set $w' = \psi(\lambda w)$, 
$r' = 2\lambda \Lambda at_k$, and $B' = B(w',r')$.
Then $\psi(\lambda B) = \psi(B(\lambda w, 2\lambda at_k))
\i B(w',2\lambda \Lambda at_k) = B'$,
just because $\psi$ is $\Lambda$-lipschitz.
Now $B'$ is not very large, because
$r' = 2\lambda \Lambda at_k \leq 2\lambda \Lambda at_0
\leq \lambda \Lambda a r \leq \Lambda a r_0$
by construction, (10.85), and (10.83).

Also set $F'_j = \psi(\lambda F_j)$ for $j \in J(y_k,t_k)$.
The $F'_j$ are now real dyadic faces of side length $r_0$,
and they all meet $B'$ because the $F_j$ meet $B$. 

We need to know the following geometrical fact about our net.
We have a collection of faces, that all meet a small ball
$B'$, and we want to know that their intersection
meets $CB'$. This is probably true with general
polyhedral networks, but here again let us cheat and
use the fact that we have a cubical network. 

Write things in coordinates. Each $F'_j$ is
given by the equations $z_i \in I_{i,j}$, 
$1 \leq i \leq n$, where each $I_{i,j}$ is either
a point or a dyadic interval of size $r_0$.
Let $w'_i$ denote the $i$-th coordinate of $w'$.
Since $B'$ meets $F'_j$, we get that
$\dist(w'_i,I_{i,j}) < r'$ for all $i$.
If $a < (3\Lambda)^{-1}$, then $r' < r_0/3$;
then for each $i$, either all the $I_{i,j}$
are equal, or else they all have a common endpoint $\xi_i$ 
which in addition is such that
$|w'_i-\xi_i| \leq r'$ (easy proof by induction on $j$). 

In all cases, we get $\xi_i \in \bigcap_{j\in J(y_k,t_k)} I_{i,j}$
such that $|w'_i-\xi_i| \leq r'$.
Now the point $\xi$ with coordinates $\xi_i$ lies
in $\bigcap_{j\in J(y_k,t_k)} F'_j$, and $|\xi-w'| \leq \sqrt n r'$.
Set $\zeta = \lambda^{-1} \psi^{-1}(\xi)$; then
$\zeta \in F$, and 
$$
|\zeta-w| = |\lambda^{-1} \psi^{-1}(\xi)-\lambda^{-1} \psi^{-1}(w')|
\leq \lambda^{-1} \Lambda |\xi-w'| \leq \lambda^{-1} \Lambda \sqrt n r'
= 2 \Lambda^2 \sqrt n at_k
\leqno (10.93) 
$$
Choose $a = {\eta \over 10 \Lambda^2 \sqrt n}$;
then $|\zeta - w| < \eta t_k /4$. We also know that
$\zeta \in F$, so $\zeta \in L(y_k,t_k)$, by (10.92).
This contradicts (10.89).

So (10.91) holds, and our construction stops after at most 
$j_{max} + 2$ steps. Let $(y_k,t_k)$ be the last pair,
where we stop. Set $y = y_k$ and $t = t_k/2$.
By definition of stopping, $(y,t)$ satisfies (10.88).

Let us try to apply Lemma 9.14 to the pair $(y,t)$.
First we need to check (for (9.15)) that 
$0 < t < C_0^{-1} \Min(\lambda^{-1} r_0, \delta)$,
but this is true because 
$t \leq t_k \leq t_0 \leq {r \over 4C_0}$
(by (10.85)), and by (10.83). Also,
$B(y,(C_0+1)t) \i 2C_0 B(y,t) \i 2C_0 B(y_0,t_0)
= B(y_0,2C_0t_0) \i B(y_0,r/2) \i B(x,r) \i U$
because $B(y_k,t_k) \i B(y_0,t_0)$ by (10.90),
and by (10.85) and (10.83); this proves (9.15).

The ugly condition (9.17) is now satisfied, precisely
because $(y,t)$ satisfies (10.88). For the last condition
(9.18), let us take the same plane $P$ as in (10.86),
and show that
$$
\dist(z,P) \leq A\varepsilon_0 t
\ \hbox{ for } z \in E \cap B(y,2t),
\leqno (10.94) 
$$
where $A = 2a^{-j_{max} - 2}$ is just another geometric 
constant, that only depends on $n$ and $\Lambda$.

Let $z\in E \cap B(y,2t) = E \cap B(y_k,t_k)$ be given.
Then $z \in E \cap B(y_0,t_0)$ (again by (10.90)), and
$\dist(z,P) \leq \varepsilon_0 t_0$ by (10.86).
So we just need to check that $t_0 \leq A t$.

But during our construction, we always took 
$t_{l+1} \geq a t_l$. Therefore,
$t = {1 \over 2} \, t_k \geq  {1 \over 2} \, a^{j_{max} + 2} t_0
= A^{-1} t_0$, and (10.94) follows.

We just proved that (9.17) is satisfied, with the constant
$A\varepsilon_0$. We shall of course choose $\varepsilon_0$
so small that $A\varepsilon_0 \leq \overline \varepsilon$,
where $\overline \varepsilon$ is the threshold in Lemma 9.14;
then the lemma applies, and we get that (9.20) holds,
i.e., that
$$
\pi(E \cap B(y,5t/3)) \hbox{ contains } 
P \cap B(\pi(y),3t/2),
\leqno (10.95)  
$$
where we denote by $\pi$ the orthogonal projection onto $P$.

We shall now conclude as in the other corollaries.
We want to check that our pair $(y,t)$ satisfies the 
conclusions of Proposition 10.82. We know that
$y \i E^\ast \cap B(x,r/2)$ because $y=y_k \in B(y_0,t_0)$
(by (10.90)) and by (10.85). Similarly, $t \leq t_0 \leq r/4$
by (10.85), and $t \geq A^{-1} t_0 
\geq {c_{\varepsilon_0} r \over 2AC_0}$ 
(by (10.85) again). So we will be able to take
$d_\varepsilon = {c_{\varepsilon_0} \over 2AC_0}$, and
we just need to check (10.84).

For each $p \in P \cap B(y,(1-A\varepsilon_0)t)$,
(10.95) gives a point $z\in E \cap B(y,3t/2)$ such that
$\pi(w) = p$. Since 
$|p-w| = |\pi(w)-w| = \dist(w,P) \leq A\varepsilon_0 t$
by (10.94), $w \in B(y,t)$. So $P \cap B(y,(1-A\varepsilon_0)t)
\i \pi(E \cap B(y,t))$, and
$$\eqalign{
\H^d(E \cap B(y,t)) &\geq \H^d(\pi(E\cap B(y,t)))
\cr&
\geq \H^d(P \cap B(y,(1-A\varepsilon_0)t)) 
\geq \omega_d(1-A\varepsilon_0)^d t^d,
}\leqno (10.96)  
$$
where $\omega_d$ is the same as in (10.84).
We choose $\varepsilon_0$ so small, depending on $\varepsilon$,
that $(1-A\varepsilon_0)^d \leq 1 - \varepsilon$, and then deduce
(10.84) from (10.96). This completes our proof of 
Proposition~10.82.
\qed

\ms
We finally come to the lower semicontinuity of $\H^d$ 
along convergent sequences of uniformly quasiminimal sets,
which we will deduce from Proposition~10.82 and the 
lower semicontinuity result of
Dal Maso, Morel, and Solimini [DMS]  
for the uniformly concentrated set.

\ms\proclaim Theorem 10.97.
Let $U$, $\{ E_k \}$, and $E$ satisfy the hypotheses (10.1), 
(10.2), (10.3), and (10.4). Also suppose that $h$ is small enough, 
depending only on $n$, $M$, and $\Lambda$.
Then 
$$
\H^d(E\cap V) \leq \liminf_{k \to +\infty} \, \H^d(E_k \cap V)
\ \hbox{ for every open set $V \i U$.}
\leqno (10.98)
$$

\ms
See Theorem 25.7 for an extension of Theorem 10.97 
where we also prove the lower semicontinuity of
$\int_{E_k\cap V} f(x) d\H^d(x)$ for some continuous
functions $f$ or even elliptic integrands
(where $f$ may also depend on the tangent plane to $E_k$ at $x$).
The proof of Theorem 25.7, which is based on a recent result
of Y. Fang [Fa], can thus be used as 
an alternative to [DMS]  
by readers that would not already be familiar with it.

\ms
Here again, we do not need (10.7). 
Of course the major difference with (10.12) is that we removed
the ugly constant $C_M$.

The (fairly short) proof of Theorem 10.97 is the same as for 
Theorem 3.4 in [D2]: 
the conclusion of Proposition~10.82  is stronger than 
what we prove in Lemma 3.6 of [D2],  
which was already more than enough to apply the results
of [DMS].  
\qed

\ms
Theorem 10.97 will lie at the center of our proof of Theorem 10.8, 
even though many complications will occur, 
both in the definition of the competitors
(in particular because we have to follow the sliding boundary rules)
and in the accounting (because Almgren's definition 
of quasiminimal sets does not cooperate too well with 
deformation mappings that are not injective).

\ms\noindent
{\bf 11. Construction of a stabler deformation:
the initial preparation} 
\ms
In this section and the next ones, 
we continue with the notation and assumptions of Theorem 10.8, 
except that we don't yet need to assume (10.7), which will only 
be needed for the final Hausdorff measure computations. 

Also, the construction of our main deformation will
be a little more unpleasant when we work under 
the Lipschitz assumption, so in most section
we shall first describe the construction under 
the simpler rigid assumption, and explain the necessary
modifications for the general case (in the best cases, this is just
a conjugation of some mappings with our bilipschitz mapping $\psi$,
but in some cases more work is needed) to the end of sections or
subsections, so that the reader may easily skip them.
We even put the corresponding text between daggers ($\dagger$)
to make the skipping easier (but it would be a shame).

Most of the next sections consists in describing the construction 
of a deformation that was done in [D2],  
and adapting it to the sliding boundary conditions.
After the construction itself, we shall complete
the argument with some Hausdorff measure estimates.
The proof will finally be competed in Sections 18 (under the rigid 
assumption) and 19 (in the Lipschitz case). We return to an
almost self-contained mode, because so many modifications
are needed from the original proof in [D2]  
(after all, most of that paper is the construction of a competitor).

So let $\{ E_k \}$ be a sequence of quasiminimal sets, such that
(10.1)-10.4) hold for some relatively closed set $E \i U$.
Since we want to show that $E$ is quasiminimal, we give ourselves a
one-parameter family of functions $\varphi_t$, 
$0 \leq t \leq 1$, such that (1.4)-(1.8) hold
for some closed ball $B$ and relative to $E$; we assume that
$$
B = \overline B(X_0,R_0), \ 
0 < R_0 < \delta,
\ \hbox{ and } \ 
\wh W \i\i U
\leqno (11.1)
$$
(as in (2.4), and where $\wh W$ is as in (2.2)). Very often, we 
shall replace our Lipschitz assumption (10.1) with the rigid 
assumption, so $U$ will be the unit ball, but this does not 
matter yet.

We want to prove (2.5), and naturally we would like to use 
the $\varphi_t$ to construct a competitor
for $E_k$ for $k$ large, apply (2.5) to $E_k$, 
get some information, and take a limit. 
Our first task will be to extend the $\varphi_t$ to $\R^n$, 
because they are not defined on $E_k$ yet.
We know that we shall probably need to modify the extension slightly, 
because we want to have (1.7) for $E_k$ and not just $E$. 
But even in the standard case when we have no $L_j$,
we cannot use the $\varphi_t$, or their extension,
directly, because of complications with the multiplicities
that will be explained soon and will be our main source of trouble.

Anyway, let us first define extensions, that we shall also
call $\varphi_t$. We shall find it better to use a specific extension
algorithm, because this way it will be easier to derive estimates.

We first set $\varphi_t(x) = x$ near $\d U$. That is, set
$$
\delta_0 = \dist(\wh W, \R^n \sm U) > 0
\ \hbox{ and } \ 
U_{ext} = \big\{ x\in \R^n \, ; \, \dist(x, \R^n \sm U) \leq \delta_0/2 \}
\leqno (11.2)
$$
(we know that $\delta_0 > 0$ because of (11.1)). We set
$$
\varphi_t(x) = x \ \hbox{ for $x\in U_{ext}$ and } 0 \leq t \leq 1.
\leqno (11.3)
$$
At this point, we have a definition of $\varphi_t$ on
$E \cup U_{ext}$, and it satisfies the properties (1.4), (1.5), (1.6)
and (1.8), with $E$ replaced by $E \cup U_{ext}$. For instance,
$\varphi_1$ is Lipschitz on $E \cup U_{ext}$ because $\varphi_1(x) = x$
on $E \sm \wh W$, $\dist(U_{ext},\wh W) \geq \delta_0 > 0$, and $\varphi_1(x) - x$
is bounded on $E \cap \wh W$.

Now we extend all these mappings to $\R^n$.
We shall use the Whitney algorithm, as it is described in
Chapter IV.2 of [St], for instance,  
and we refer to this book for details on the construction that follows. 
We cover $U \sm (E \cup U_{ext})$ by Whitney cubes 
$Q_j \i \R^n \sm (E \cup U_{ext})$, $j\in J$, with disjoint 
interiors, and such that
$$
10\diam(Q_j) \leq \dist(Q_j,E \cup U_{ext}) \leq 21\diam(Q_j).
\leqno (11.4)
$$
This easy to do (use maximal dyadic cubes with the first inequality),
and the point is that then the cubes $3Q_j$ have bounded overlap.

Also choose, for each $j$, a point $\xi_j \in E \cup U_{ext}$ such that
$\dist(\xi,Q_j) = \dist(Q_j,E \cup U_{ext})$, and construct a partition of unity 
subordinate to the $Q_j$, which means a collection of smooth functions
$\chi_j \geq 0$ such that
$$
\sum_j \chi_j = {\bf 1}_{\R^n \sm (E \cup U_{ext})},
\leqno (11.5)
$$
$$
0 \leq \chi_j \leq {\bf 1}_{2Q_j}
\ \hbox{ for each $j$,}
\leqno (11.6)
$$
and
$$
||\nabla \chi_j||_\infty \leq C \diam(Q_j)^{-1}.
\leqno (11.7)
$$
We use the Whitney extension formula to extend
the function $\varphi_t(x) - x$ from $E \cup U_{ext}$ to $\R^n$. That is, we set
$$
\varphi_t(x) = x + \sum_j \chi_j(x) [\varphi_t(\xi_j)- \xi_j]
\ \hbox{ for } x\in \R^n \sm (E \cup U_{ext}).
\leqno (11.8)
$$
Naturally we keep $\varphi_t$ as it was on the set $E \cup U_{ext}$.

It will be useful to know that if $\chi_j(x) \neq 0$, for some
$j$ such that $\xi_j \in E$, then $x\in 2Q_j$ and hence
$$
\dist(x,E) \leq |\xi_j - x| = \dist(Q_j,E)
\leq \dist(x,E) + \diam(Q_j)
\leqno (11.9)
$$
and hence, since (11.4) says that $10 \diam(Q_j) \leq \dist(Q_j,E)$,
$$
\dist(x,E) \geq \dist(Q_j,E) - \diam(Q_j) \geq 9 \diam(Q_j),
\leqno (11.10)
$$
which we can plug back in (11.9) to get that
$$
|\xi_j - x| \leq \dist(x,E) + \diam(Q_j)
\leq {10 \over 9} \dist(x,E).
\leqno (11.11)
$$

The advantage of extending $\varphi_t(x) - x$ is that we more
easily spot places where it vanishes. Let us check that
$$
\varphi_t(x) = x
\ \hbox{ when } \dist(x,E) < {9 \over 10}\dist(x,W_t),
\leqno (11.12)
$$
where as before
$$
W_t = \big\{ y\in E \, ; \, \varphi_t(y) \neq y \big\}.
\leqno (11.13)
$$
This is clear when $x\in E$ (because then $\dist(x,W_t) > 0$), 
and when $x\in U_{ext}$ (by (11.3)),
so let us consider $x\in \R^n \sm (E \cup U_{ext})$. 
If $\varphi_t(x) \neq x$, (11.8) says that we can find $j$
such that $\chi_j(x) \neq 0$ and $\varphi_t(\xi_j) \neq \xi_j$.
Then $\xi_j \in E$ by (11.3), and 
$|\xi_j - x| \leq {10 \over 9} \dist(x,E)$ by (11.11); 
but we assumed that $\dist(x,W_t) > {10 \over 9} \dist(x,E)$,
so $\xi_j \notin W_t$, hence $\varphi_t(\xi_j) = \xi_j$,
a contradiction which proves (11.12).

We also have that
$$
\varphi_0(x) = x \ \hbox{ for } x\in \R^n,
\leqno (11.14)
$$
just because all the $\varphi_0(\xi_j)-\xi_j$ vanish, that
$$
(t,x) \to \varphi_t(x) 
\ \hbox{ is continuous on } [0,1] \times \R^n
\leqno (11.15)
$$
(in particular because each fixed $\varphi_t(\xi_j)$ is a 
continuous function of $t$), and similarly
$$
\varphi_1 : \R^n \to \R^n
\ \hbox{ is Lipschitz,}
\leqno (11.16)
$$
because we used the standard formula for the Whitney 
extension theorem.

Our extensions $\varphi_t$ will be better when we stay close to $E$,
for instance because otherwise (11.12) does not give much of a 
control, but this is all right because we only need them on sets
$E_k$ that tend to $E$.
But also, the $\varphi_t$ have some defects that we'll need to fix.
The main one, which will be explained soon, is that $\varphi_1$
may be one-to-one on the $E_k$, while it is many-to-one on $E$,
which possibly makes the $\varphi_1(E_k)$ much worse competitors
than $\varphi_1(E)$.

The truth is that we mostly care about $\varphi_1$.
We need to keep the $\varphi_t$, $0 \leq t < 1$,
because eventually we shall need to check (1.7) on the $E_k$,
but the main point of the argument concerns estimates
like (2.5), where only $\varphi_1$ counts.

Except for the fact that we have to worry about (1.7) and have
a slightly more complicated way to define competitors
(we now have a ball $B$ and the open set $U$), we will
mostly follow the construction of [D2].  
Even though we will change many little things, it will some times
be convenient to refer to [D2]  
for small independent things.

So we are interested in the values of $f = \varphi_1$ near $E$.
Even in the standard case with no $L_j$, if by bad luck there
is a region where $E$ and $E_k$ are composed of many pieces, 
$f$ maps all the pieces of $E$ to a same small disk, say, 
but maps all the $E_k$ to parallel, but disjoint little disks, 
$f(E)$ may be a quite reasonable competitor 
(because the measure of the single disk is small),
but not $f(E_k)$. Then, when we apply (2.5) to $E_k$ and
$f(E_k)$, we won't get much information, not enough to control
$f(E)$. What we intend to do in this case, when $f$ sends
many pieces to parallel and nearby disks, is to modify $f$
(typically, by composing it with a projection) 
so that it sends all these pieces to a single disk.
This will make $f(E_k)$ a much better competitor, and then we
have a good chance to run the usual lower semicontinuity
arguments and get the desired inequality (2.5).

A good part of the construction below consists in doing such grouping,
but obviously this will require some nontrivial amount of
cutting and pasting. As the reader may have guessed, we shall use
the fact that $E$ is rectifiable (to show that it has a tangent
at most points), the fact that $\varphi_1$ is Lipschitz on $E$
(to show that it is often close to its differential), and lots of
covering arguments (to reduce to situations where $\varphi_1$
is almost affine and $E$ is almost flat). The word stability in the
title of the section refers to the fact that after this grouping,
the total measure of $f(E_k)$ will be much less dependent on $k$.

\ms\noindent{\bf 11.17. Remark about the many constants.}
Since there will be lots of constants in this argument, 
let us announce here in which order we intend to choose them,
so that the reader may more easily check that we do not cheat.
We shall systematically denote by $C$ constants that depend
only on $n$ $M$, and $\Lambda$ (when we work under
the Lipschitz assumption). This includes the local 
Ahlfors-regularity constants for $E$.

Next observe that from now on, $U$, $B$, the $\varphi_t$,
and in particular $f = \varphi_1$ are fixed, so we shall not mind
if our constants depend on $r_0$, $\lambda$, of $\psi$ in the Lipschitz 
assumption, or on $f$, typically through its Lipschitz
constant $|f|_{lip}$. Similarly, we can let our constants 
depend on the number $\H^d(\{ x\in E \, ; \, f(x) \neq x \})$. 
In both cases, we shall often indicate this dependence, 
mostly to comfort the reader, but will not be a real issue.

A first string of constants is
$$\eqalign{
&\hbox{$\gamma > 0$ (small), $a < 1$ (close to $1$), $\alpha > 0$ 
(small), $N$ (large)}, 
\cr&\hskip2cm
\hbox{$\eta > 0$ (small), and $\varepsilon > 0$ 
(small)}
}\leqno (11.18)
$$
to be chosen in this order. 
Our small constants $\delta_i$, $1 \leq \delta_i \leq 9$, 
will be allowed to depend on these constants 
(even though the first ones don't), 
and are thus chosen after $\eta$ and $\varepsilon$. 
Typically, they are chosen smaller and smaller.
They act as scale parameters, and force us to work with balls 
that are small enough for some properties to hold, but they 
should not have an incidence on the estimates.

Finally, we will choose two last small constant $\varepsilon_0$
and $\varepsilon_\ast$. 
They also determine the distance that we allow between 
$E$ and the $E_k$, and so they could also have been called
$\delta_{10}$ and $\delta_{11}$, but we decided to revert
to the name $\varepsilon$ to show that they are even smaller.

Let us also mention that our estimates will only be valid for $k$
large, depending on the various constants, and in particular
the $\delta_i$ and $\varepsilon_0$. 

We cut the construction into a few smaller steps; only the first
one will be completed in this section. As explained above, we shall
first carry the construction under the rigid assumption, and then we 
shall explain how to modify things under the Lipschitz assumption.

\ms\noindent{\bf Step 1. We remove a few small bad sets.}
Before we define balls $B_j$ and modify $f$ on them,
we remove some small bad sets from $E$, where
$E$ or $f$ is not regular enough. 
Set
$$
W_f = \big\{ x\in \R^n \, ; \, f(x) \neq x \big\}
= \big\{ x\in \R^n \, ; \, \varphi_1(x) \neq x \big\};
\leqno (11.19) 
$$
our star starting set is
$$
X_0 = E \cap W_f = \big\{ x\in E \, ; \, f(x) \neq x \big\}
\i \wh W \i \i U,
\leqno (11.20)
$$
(by (2.2) and (11.1))
and we immediately replace it with a compact subset $X_1 \i X_0$,
such that
$$
\H^d(X_0 \sm X_1) \leq \eta,
\leqno (11.21)
$$
where $\eta$ is some very small positive number,
that will tend to $0$ at the end of the argument. Set
$$
\delta_1 = \dist(X_1,\R^n \sm W_f) > 0;
\leqno (11.22)
$$
the fact that $\delta_1 > 0$ (because $X_1$ is compact and $W_f$ 
is open) will make it easier to stay inside $W_f$ when we 
cover $X_1$ by small balls.

Some manipulations will be easier if we force $f(x)$ to stay far 
from the boundaries of faces, because it will make the smallest
face that contains $f(x)$ locally constant.
For $0 \leq l \leq n$, and under the rigid assumption,
denote by ${\cal S}_l$ the union of all the 
faces of dimension $l$ of the dyadic cubes of side length $r_0$ of
our usual grid. When we work under the Lipschitz assumption, 
we shall just call ``cubes" the images of standard dyadic cubes by 
$\lambda^{-1}\psi^{-1}$, and similarly for faces, and we shall define the 
${\cal S}_l$ in terms of these distorted faces. Set
$$
X_{1,\delta}(l) = \big\{ x\in X_1 \, ; \, 
f(x) \in {\cal S}_l \sm {\cal S}_{l-1} \ \hbox{ and } \ 
\dist(f(x),{\cal S}_{l-1}) \geq \delta \big\}
\leqno (11.23)
$$
for $\delta > 0$ and $1 \leq l \leq n$; for  $l=0$, simply put
$$
X_{1,\delta}(0) = \big\{ x\in X_1 \, ; \, 
f(x) \in {\cal S}_0 \big\}.
\leqno (11.24)
$$
Then set
$$
X_{1,\delta} = \bigcup_{0 \leq l \leq n} X_{1,\delta}(l).
\leqno (11.25)
$$
Since 
$X_1$ is the monotone union of the $X_{1,\delta}$ 
(when $\delta$ decreases to $0$), we can choose $\delta_2 >0$ 
so small that $\H^d(X_1 \sm  X_{1,\delta_2}) \leq \eta$. Then we set
$$
X_2 = X_{1,\delta_2} = \bigcup_{0 \leq l \leq n} X_{1,\delta_2}(l)
\ \hbox{ and so } \H^d(X_1 \sm  X_2) \leq \eta.
\leqno (11.26)
$$
For each $x\in X_2$, denote by $F(f(x))$ the smallest face of our dyadic 
grid that contains $f(x)$. We claim that (under the rigid assumption)
$$
F(f(x)) = F(f(y)) \ \hbox{ for $x,y \in X_2$ such that } |f(x)-f(y)| 
< \delta_2.
\leqno (11.27)
$$
Indeed, let $l = l(x)$ be such that $x\in X_{1,\delta_2}(l)$, and define
$l(y)$ similarly. By symmetry, we may assume that $l \geq l(y)$.
By definition, $f(x)\in {\cal S}_l$ and (if $l \geq 0$)
$\dist(f(x),{\cal S}_{l-1}) \geq \delta_2$. In particular,
$f(y)$ cannot lie on ${\cal S}_{l-1}$, so $l(y) \geq l$ and,
by our symmetry assumption, $l(y)=l$. Now $F(f(x))$ and $F(f(y))$ are two 
faces of the same dimension $l$. Suppose for a moment that they are different.
If $l=0$, this is impossible because $|f(x)-f(y)| < \delta_2$ and so
(if we chose $\delta_2 \leq r_0$) $f(x)$ and $f(y)$ cannot both lie on
${\cal S}_0$. Otherwise, (3.8) says that
$$
\dist(f(x),\d F(f(x))) \leq \dist(f(x),F(f(y))) \leq |f(x)-f(y)| < 
\delta_2,
\leqno (11.28)
$$
which contradicts the fact that $\dist(f(x),{\cal S}_{l-1}) \geq 
\delta_2$. This proves (11.27) in the rigid case. 

$ \dagger$ We can do the same argument under the Lipschitz assumption, 
where faces are only bilipschitz images of faces of a true rigid grid; 
we just need to require that
$|f(x)-f(y)| < \Lambda^{-2} \delta_2$ in (11.27), because we may lose
a constant $\Lambda^2$ in the first inequality of (11.28). $\dagger$

\ms
Next we want to use the rectifiability of $E$,
which we deduce from Proposition 10.15.
By Theorem 15.21 on page 214 of [Ma], 
we can find a countable collection of $C^1$ embedded 
submanifolds $\Gamma_s$ (or, if you prefer, images by 
rotations of $C^1$ graphs) of dimension $d$, such that
$$
\H^d(E \sm \bigcup_{s} \Gamma_s) = 0.
\leqno (11.29)
$$
To be fair, we don't really need $C^1$ submanifolds,
and Lipschitz graphs would not have required much more work, 
but on the other hand, we should not pretend that we do not use
strong results, when in the proof of Proposition 10.15
we used the much stronger fact that unrectifiable $d$-sets 
have negligible projections in almost all directions.

Select a finite set $S$ of indices, so that if we set
$X_3 = X_2 \cap \Big[\bigcup_{s\in S} \Gamma_s \Big]$, then
$$
\H^d(X_2 \sm X_3) \leq \eta.
\leqno (11.30)
$$
Put a total order on $S$, and set
$$
X_3(s) = X_3 \cap \Gamma_s \sm \bigcup_{s'<s} \Gamma_{s'}
\leqno (11.31)
$$
for $s\in S$. Thus the $X_3(s)$ are disjoint and cover $X_3$.
We know that 
$$
\lim_{r \to 0} r^{-d} H^d(\Gamma_s \cap B(x,r)\sm X_3(s)) = 0
\leqno (11.32)
$$
for $\H^d$-almost every $x\in X_3(s)$, just because 
$\H^d(\Gamma_s)$ is locally finite and $\Gamma_s \sm X_3(s)$
does not meet $X_3(s)$; see for instance
Theorem 6.2 on page 89 of [Ma].  
Similarly,
$$
\lim_{r \to 0} r^{-d} H^d(E \cap B(x,r)\sm X_3(s)) = 0
\leqno (11.33)
$$
for $\H^d$-almost every $x\in X_3(s)$, this time because
$X_3(s)$ has a neighborhood where $\H^d(E)$ is finite
(recall that $X_3(s) \i X_1$ and $X_1$ is a compact subset of 
$\wh W \i\i U$, by (11.20)).

Let us check that if $x\in X_3(s)$ is such that (11.33) holds, 
and if $P_x$ denotes the tangent plane to $\Gamma_s$ at $x$, 
then $P_x$ is also a tangent plane to $E$ at $x$, which means that 
$$
\lim_{y \to x \, ; \, y\in E} \  |y-x|^{-1} \dist(y,P_x) = 0.
\leqno (11.34)
$$
Indeed, if $d(y) = \dist(y,X_3(s))$, then $d(y) \leq |y-x|$ trivially,
and
$$\eqalign{
C^{-1} d(y)^{d} &\leq \H^d(E \cap B(y,d(y))
= \H^d(E \cap B(y,d(y))\sm X_3(s))
\cr&\leq H^d(E \cap B(x,2|y-x|)\sm X_3(s))
= o(|y-x|^d)
}\leqno (11.35)
$$
because $E$ is locally Ahlfors-regular and 
$B(y,d(y)) \i B(x,2|y-x|)$, and by (11.33).
But then $\dist(y,\Gamma_s) \leq d(y) = o(|x-y|)$ because
$X_3(s) \i \Gamma_s$, and now (11.34) holds because 
$\Gamma_s$ is tangent to $P_x$ at $x$.
Notice that without surprise,  $E$ and $\Gamma_s$
share the same tangent plane at $x$ (regardless of
the rest of $E$, which anyway does not matter because
of (11.33)); the uniqueness of the tangent plane to $E$ 
follows from the local Ahlfors-regularity of $E$, but
we could also deduce it from the fact that by (11.32),
most points of $\Gamma_s$ lie in $X_3(s)$, and hence in $E$.

We also want to throw out the points of $X_3(s)$ where
$f$ is not differentiable in the direction of $P_x$.
That is, denote by $X_4(s)$ the set of points $x\in X_3(s)$
such that (11.32) and (11.33) hold, and there exists
an affine function 
$A_x : \R^n \to \R^n$, of rank at most $d$ and with Lipschitz norm 
$$
|DA_x| \leq |f|_{lip},
\leqno (11.36)
$$
and such that 
$$
\lim_{y \to x \, ; \, y \in \Gamma_s} {|f(y)-A_x(y)| \over |y-x|} = 0.
\leqno (11.37) 
$$
Let us check that if $x\in X_3(s)$ satisfies (11.32) and (11.33),
and if the restriction of $f$ to $\Gamma_s$ is differentiable at $x$,
then $x\in X_4(s)$; clearly this will imply that
$$
\H^d(X_3(s) \sm X_4(s)) = 0.
\leqno (11.38)
$$
Let $x\in X_3(s)$ be like this. There is a small neighborhood of $x$
where $\Gamma_s$ is the graph of some $C^1$ function $F$, defined on $P_x$
and with values in $P_x^\perp$, the $(n-d)$-vector space 
perpendicular to $P_x$. That is, $\Gamma_s$ is locally equal to the
set of points $z+F(z)$, $z\in P_x$, and with these conventions
$F(x)=0$ and $Df(x) = 0$. 

Set $g(z) = f(z+F(z))$ for $z\in P_x$ near $x$.
By assumption, $g$ is differentiable at $z=x$. 
Denote by $P'_x$ the vector space parallel to $P_x$ and by 
$A : P'_x \to \R^n$ the differential in question.
Obviously, the rank of $A$ is at most $d$, and the norm of $A$
is at most $|f|_{lip}$, because for each $\lambda > |f|_{lip}$,
there is a neighborhood of $x$ in $P_x$ where $g$ is
$\lambda$-Lipschitz.

Set $A_x(y) = f(x)+A(\pi(y-x))$ for $y\in \R^n$, and where
$\pi$ denotes the orthogonal projection on $P'_x$.
Then $A_x$ is affine, with rank at most $d$, and (11.36) holds.
For $y\in \Gamma_s$ close enough to $x$, write $y = z+F(z)$,
where $z = x + \pi(y-x)$ is the projection of $y$ on $P_x$; then
$$\eqalign{
|f(y)-A_x(y)| &= |g(z)-A_x(y)| = |g(z)-f(x)-A(\pi(y-x))|
\cr&= |g(z)-f(x)-A(z-x)| = o(|z-x|) = o(|y-x|)
}\leqno (11.39)
$$
by definition of $A_x$ and $A$ (and because $f(x)=g(x)$ since
$F(x)=0$). So (11.37) holds and $x\in X_4(s)$, 
as needed for (11.38). 

Let us also check that 
$$
\lim_{y \to x \, ; \, y \in E \cup P_x} {|f(y)-A_x(y)| \over |y-x|} = 0
\ \hbox{ when } x\in X_4(s).
\leqno (11.40) 
$$
For $y\in E \cup P_x$ near $x$, choose $w\in \Gamma_s$ such that
$|w-y| \leq 2\dist(y,\Gamma_s)$; we know from (11.34) or the fact that
$P_x$ is tangent to $\Gamma_s$ that $|w-y| = o(|y-x|)$, and then
$$\eqalign{
|f(y)-A_x(y)| &\leq |f(w)-A_x(w)|+ |f(w)-f(y)|+|A_x(w)-A_x(y)|
\cr&
\leq |f(w)-A_x(w)|+2|w-y| |f|_{lip} 
\cr&
= o(|w-x|) + o(|y-x|) = o(|y-x|),
}\leqno (11.41) 
$$
by (11.37); (11.40) follows.

We want to have all the properties above with some uniformity.
So we let $\varepsilon > 0$ be small, and denote by $X_5(s)$
the set of points $x\in X_4(s)$ such that the following
properties hold. First, there is a $C^1$ mapping
$F_x : P_x \to P_x^\perp$, such that
$$
F_x(x) = 0 \, , \,  DF_x(x) = 0 \, ,  \hbox{ and } 
||DF_x||_\infty \leq \varepsilon
\leqno (11.42)
$$
and
$$
z + F_x(z) \in \Gamma_s
\ \hbox{ for } z\in P_x \cap B(x,\delta_3),
\leqno (11.43)
$$
where $\delta_3$ is another small constant, that will
be chosen soon, depending on $\varepsilon$. This is just a 
quantified version of the description of $\Gamma_s$ near $x$
that we used below (11.38). Next,
$$
\H^d(B(x,r) \cap [\Gamma_s \cup E] \sm X_3(s)) \leq \varepsilon r^d
\ \hbox{ for } 0 < r \leq \delta_3
\leqno (11.44) 
$$
(compare with (11.32) and (11.33)), and also
$$
\dist(y,P_x) \leq \varepsilon |y-x|
\ \hbox{ for } y\in E \cap B(x,\delta_3)
\leqno (11.45)
$$
as in (11.34). Finally,
$$
|f(y)-A_x(y)| \leq \varepsilon |y-x|
\ \hbox{ for } y\in [E \cup P_x \cup \Gamma_s] \cap B(x,\delta_3),
\leqno (11.46)
$$
as in (11.37) and (11.40). 

Each set $X_4(s)$ is the monotone union, when $\delta_3$ goes
to $0$, of the corresponding sets $X_5(s)$ that we just defined.
So we can choose $\delta_3 > 0$ so small that if we set
$$
X_4 = \bigcup_{s\in S} X_4(s)
\ \hbox{ and }
X_5 = \bigcup_{s\in S} X_5(s) \i X_4,
\leqno (11.47)
$$
then $\H^d(X_4 \sm X_5) \leq \eta$,
and hence
$$
\H^d(X_0 \sm X_5) \leq 4 \eta
\leqno (11.48)
$$
by (11.21), (11.26), (11.30), (11.31), and (11.38).

\ms\noindent $\dagger$ {\bf Remark 11.49.}
In this first step that we just finished, the flatness of the
faces does not show up. Under the Lipschitz assumption, 
we can proceed exactly as we did, except that cubes and
faces are more complicated. 
Things will become more unpleasant
when we try to use the existence of the tangent planes $P_x$
to derive information on the closeness of $E$ to the (no longer flat) 
faces of cubes, and in fact we shall need to require the 
equivalent of (11.36), (11.37), and (11.40) also for the mapping 
$\wt f$ defined by
$$
\wt f(x) = \psi(\lambda f(x)) \in B(0,1)
\ \hbox{ for } x\in E
\leqno (11.50)
$$
(which makes sense because we know that $f = \varphi_1 : E \to U$ and 
$\psi : \lambda U \to B(0,1)$); see near (12.36). $\dagger$

\bigskip\noindent
{\bf 12. Step 2 of the construction: the places where $f$ is many-to-one} 
\ms

In this section we continue the construction of Section 11, 
and modify the function $f$ in some balls $B_j$, $j\in J_1$,
where we can make $f$ highly non injective.
Let $N$ be a large number, and set
$$
Y_N = \big\{ y\in f(X_5) \, ; \, X_5 \cap f^{-1}(y)
\hbox{ contains at least $N$ distinct points}\big\}.
\leqno (12.1)
$$
As before, it will be easier to (demand some uniformity
and) control the set 
$$\eqalign{
Y_N(\delta_4) = &\big\{ y\in f(X_5) \, ; \, X_5 \cap f^{-1}(y)
\hbox{ contains at least }
\cr&\hskip1cm
\hbox{ $N$ distinct points at mutual distances $\geq \delta_4$}\big\}
}\leqno (12.2)
$$
for some small $\delta_4 >0$. Since $f^{-1}(Y_N)$ is the
monotone union of the $f^{-1}(Y_N(\delta_4))$, we can choose
$\delta_4 > 0$ so small, depending on $\eta$ and other constants, that
if we set
$$
X_N(\delta_4) =  X_5 \cap f^{-1}(Y_N(\delta_4)),
\leqno (12.3)
$$
then
$$
\H^d\big([X_5 \cap f^{-1}(Y_N)] \sm X_N(\delta_4)\big)
\leq \eta. 
\leqno (12.4)
$$

We shall need a covering of $X_N(\delta_4)$.

\msi{\bf Step 2.a. We cover $X_N(\delta_4)$ by balls $B_j = B(x_j,t)$, $j\in J_1$.}
\ms
We want to cover $X_N(\delta_4)$ with small balls of the same
very small radius $t/2$, but let us first say how small we want
$t$ to be. Set
$$
\delta_5 = \inf\big\{ |f(x)-x| \, ; \, x\in X_1 \big\};
\leqno (12.5)
$$ 
notice that $\delta_5 > 0$ because $X_1$ is compact and $f(x) \neq x$
for $x \in X_1 \i X_0$ (by (11.20)). Recall that we set 
$$
\delta_0 = \dist(\wh W,\R^n \sm U) > 0
\leqno (12.6)
$$
in (11.2). Pick
$$
\delta_6 < {1 \over 10 \Lambda^2 (1+|f|_{lip})} \,
\Min\big(\lambda^{-1} r_0, \delta, \delta_0, \delta_1,
\delta_2,\delta_3,\delta_4,\delta_5 \big)
\leqno (12.7)
$$ 
(where $\lambda = \Lambda = 1$ in the rigid case)
and any $t$ such that
$$
0 < t < \delta_6.
\leqno (12.8)
$$
Then pick a maximal collection
$\{ x_j \}$, $j\in J_1$, of points in $X_N(\delta_4)$, that 
lie at mutual distances at least $t/3$. Thus
$$
X_N(\delta_4) \i \bigcup_{j \in J_1} B(x_j,t/2)
\leqno (12.9)
$$
by maximality, and we claim that 
$$
\hbox{$J_1$ has at most $C t^{-d} \H^d(X_0)$ elements.}
\leqno (12.10)
$$
Indeed, $x_j \in X_5 \i X_1$ for $j\in J_1$, so
$$
t < \delta_6 < {1 \over 10} \, \delta_1 
= {1 \over 10} \, \dist(X_1,\R^n \sm W_f)
\leq {1 \over 10} \, \dist(x_j,\R^n \sm W_f),
\leqno (12.11)
$$
which implies that $E \cap B(x_j,t/6) \i E \cap W_f = X_0$
(see the definition (11.20)). Let us also show that 
$$
\H^d(E\cap B(x_j,t/6)) \geq C^{-1} t^d.
\leqno (12.12)
$$ 
We want to apply Proposition 4.1 or Proposition 4.74, 
so we just need to check that $t/6 < \Min(\lambda^{-1} r_0,\delta)$ 
and $B(x_j,t/3) \i U$.
The first one follows from (12.7) and (12.8), and the second one
holds because $x\in X_0 \i \wh W$ and $t < \delta_6 \leq \delta_0$.
So Proposition 4.1 or 4.74 applies and gives (12.12). 

The $E \cap B(x_j,t/6)$, $j\in J_1$ are disjoint (by definition of $J_1$),
and contained in $X_0$, so (12.10) follows from (12.11).
 
We agree that (12.10) is not a very good bound, 
but a large choice of $N$, depending on $\H^d(X_0)$
and $\eta$, will compensate. At least, notice that 
$\H^d(X_0) < +\infty$ (because $X_0 \i \wh W \i\i U$), and will be
taken as a constant (it only depends on $E$ and $f$).

\msi{\bf Step 2.b. We cover  $Y_N(\delta_4)$ by balls $D_l$, $l\in {\cal L}$.}
\ms
We also want a covering of $Y_N(\delta_4)$. So we take a
maximal set of points $y_l$, $l\in {\cal L}$, in $Y_N(\delta_4)$,
at mutual distances at least $t/2$. Then 
$$
\hbox{the balls $D_l = B(y_l,t)$, $l\in {\cal L}$,
cover } Y_N(\delta_4).
\leqno (12.13)
$$
Let us prove that the cardinality of ${\cal L}$ is
$$
|{\cal L}| \leq C N^{-1} (1+|f|_{lip})^d t^{-d} \H^d(X_0).
\leqno (12.14)
$$
For each $l\in {\cal L}$, select $N$ points
$x_{l,j} \in X_5$, $1 \leq j \leq N$, at mutual
distances at least $\delta_4$, such that $f(x_{l,j})=y_l$. 
Such points exist, by the definition (12.2), and they lie in
$X_N(\delta_4)$, by (12.3).

Set $s = {t \over 4(1+|f|_{lip})}$ and $B_{l,j} = B(x_{l,j},s)$;
we claim that
$$
\hbox{$B_{l,j}$ is disjoint from $B_{l',j'}$ when $(l,j) \neq (l',j')$.}
\leqno (12.15)
$$
When $l=l'$ and $j \neq j'$, this is because 
$s \leq t/4 \leq \delta_4/4 \leq |x_{l,j}-x_{l,j'}|/4$ (by (12.8),
(12.7), and our choice of points $x_{l,j}$); 
when $l \neq l'$, this is because if $x\in B_{l,j}$ and $x' \in B_{l',j'}$, 
then
$$
|f(x)-f(x')| 
\geq |f(x_{l,j})-f(x_{l',j'})| - 2 s |f|_{lip} 
= |y_{l}-y_{l'}| - 2 s |f|_{lip}
\geq t/2 - t/4 > 0.
\leqno (12.16)
$$
With the same verification as for (12.12), $\H^d(E \cap B_{l,j}) 
\geq C^{-1} s^d \geq C^{-1} (1+|f|_{lip})^{-d} t^{d}$.
Also, $E \cap B_{l,j} \i X_0$ because 
$s \leq t/4 \leq \dist(x_{l,j}, \R^n \sm W_f)$ 
(because $x_{l,j}\in X_5$ and by (11.20) and the proof of (12.11)).
Since there are $N |{\cal L}|$ balls $B_{l,j}$, which are disjoint by 
(12.15) and have a total mass at most $\H^d(X_0)$, we get (12.14).

\msi{\bf Step 2.c. The collection of disks $Q$, $Q\in {\cal F}_l$.}
\ms
For each $l\in {\cal L}$, we want to select a reasonably large collection 
of affine subspaces $P$, and then disks $Q \i P$. The general idea is 
that when we modify $f$ on the $B_j$, $j\in J_1$, we want to send 
points preferably to these disks $Q$, which will make $f$ more
stably many-to-one.

This is a stage of the construction where we need 
to be careful about the boundary pieces $L_j$, so we shall 
diverge slightly from [D2])  
and also we shall need to modify some things when we work under the 
Lipschitz assumption, but let us first describe what we do when
the rigid assumption hold. 

Let us fix $l\in {\cal L}$. Denote by $F_l$ the smallest face
of our usual grid that contains $y_l$; this makes sense, because
if $y_l$ lies on two faces $F$ and $F'$, then $F\cap F'$ is also
a face that contains $y_l$. Denote by $d(l)$ the dimension of
$F_l$, and by $W(y_l)$ the $d(l)$-dimensional affine subspace
that contains $F_l$. 

If $d(l) \leq d$, we just choose one affine subspace, namely
$P = W(y_l)$, and one disk $Q$, namely
$$
Q = P \cap B(y_l, 3 (1+|f|_{lip}) t).
\leqno (12.17) 
$$
If $d(l) > d$, we choose a whole collection of affine $d$-planes $P$ 
through $D_l$, all of them contained in $W(y_l)$, with a density 
property that we shall explain soon. 
For each $P$ that we choose, we still define $Q$ by (12.17);
this gives a collection ${\cal F}_l$ of disks $Q$, 
$Q \in {\cal F}_l$, (which is just composed of one disk when
$d(l) \leq d$). 

The density property is the following.
Let $\alpha > 0$ be a very small constant (to be chosen later);
we demand that if $d(l) > d$, then for each
affine $d$-plane $P'$ through $D_l$ 
which is contained in $W(y_l)$, 
we can find $Q \in {\cal F}_l$ such that
$$
\dist(z,Q) \leq \alpha t
\ \hbox{ for } z\in P' \cap B(y_l, 3(1+|f|_{lip}) t). 
\leqno (12.18) 
$$
So we choose the set ${\cal F}_l$ like this, but with not
too many elements, so that 
$$
|{\cal F}_l| \leq C(\alpha,f)
\leqno (12.19)
$$
for some constant $C(\alpha,f)$ that depends on $\alpha$
and $|f|_{lip}$. Now we set
${\cal F} = \cup_{l \in {\cal L}} {\cal F}_l$ 
and observe that
$$
|{\cal F}| \leq C(\alpha,f) |{\cal L}|
\leq C(\alpha,f) N^{-1} t^{-d} \H^d(X_0)
\leqno (12.20)
$$
by (12.14), and with new constants $C(\alpha,f)$ that depend on $\alpha$
and $|f|_{lip}$. Let us record the fact that by 
construction,
$$
Q \i W(y_l) \cap B(y_l, 3 (1+|f|_{lip}) t)
\ \hbox{ for $l \in {\cal L}$ and } Q \in {\cal F}_l.
\leqno (12.21)
$$ 
We shall see why this is important when we check the boundary 
condition (1.7). 

\ms
$\dagger$
When we work under the Lipschitz assumption, we proceed
almost the same way.
Recall that the dyadic cubes and faces are now obtained
from the standard ones in the unit ball, using (the inverse of) the 
bilipschitz mapping $\psi : \lambda U \to B(0,1)$; thus $F_l$, 
the smallest face that contains $y_l$, is the image 
$F_l = \lambda^{-1}\psi^{-1}(\wt F_l)$ 
of some true dyadic face $\wt F_l$ (the smallest one that contains 
$\wt y_l = \psi(\lambda y_l)$).
When $d(l) \leq d$, we just choose one affine space $\wt P$, namely
the $d(l)$-dimensional affine subspace $W(\wt y_l)$ that contains
$\wt F_l$, and set
$$
\wt Q = \wt P \cap B(\wt y_l, 20\lambda \Lambda^2 (1+|f|_{lip}) t)
\ \hbox{ and } \ 
Q = \lambda^{-1}\psi^{-1}(\wt Q).
\leqno (12.22) 
$$
When $d(l) > d$, we still define $\wt Q$ and $Q$ by (12.22), but we let
$\wt P$ run through a fairly large collection of $d$-planes $\wt P$,
such that
$$
\wt P \i W(\wt y_l) \ \hbox{ and } \ 
\wt P \hbox{ meets } 
B(\wt y_l, 2 \lambda \Lambda (1+|f|_{lip}) t).
\leqno (12.23) 
$$
We choose this collection so dense that,
if $\wt P'$ is any other affine $d$-plane such that
$$
\wt P' \i W(\wt y_l) \ \hbox{ and } \ 
\wt P' \hbox{ meets } 
B(\wt y_l, 2 \lambda \Lambda (1+|f|_{lip}) t)),
\leqno (12.24) 
$$
then there is a $d$-plane $\wt P$ in the collection such that
$$
\dist(z,\wt P) \leq \lambda \Lambda^{-1} \alpha t
\ \hbox{ for } z\in \wt P' 
\cap B(\wt y_l, 10\lambda \Lambda^2 (1+|f|_{lip}) t). 
\leqno (12.25) 
$$
We get a collection of plates $Q$, $Q \in {\cal F}_l$, and again we 
can manage so that each plane $\wt P$ we choose satisfies (12.23), and
$|{\cal F}_l| \leq C(\alpha,f)$, where here
$C(\alpha,f)$ also depends on $\Lambda$. 
This way we still get (12.19) and (12.20), and the 
analogue of (12.21), namely the fact that
$$
\wt Q = \psi(\lambda Q) \i W(\wt y_l) \cap 
B(\wt y_l, 20\lambda \Lambda^2 (1+|f|_{lip}) t) 
\ \hbox{ for $l \in {\cal L}$ and } Q \in {\cal F}_l.
\leqno (12.26)
$$
Incidentally, notice that none of our main constants will
depend on $\lambda$, and in fact we could easily get rid of $\lambda$ 
with a simple dilation of $E$ and the $E_k$; we shall not do
this and keep mentioning $\lambda$ in our estimates, but the
reader could also decide not to bother and make $\lambda=1$ everywhere.
$\dagger$

\msi{\bf Step 2.d. Where do the tangent planes go?}
\ms
The following lemma will help us soon in our choice
of approximate tangent planes in the image, but we state it
independently. As usual, we start under the rigid assumption.

\ms\proclaim Lemma 12.27.
For $x \in X_5$, let $P_x$ and $A_x$ be as in (11.34) and (11.40),
denote by $F(f(x))$ the smallest face of our usual grid 
that contains $f(x)$ 
and by $W(f(x))$ the smallest affine subspace 
that contains $F(f(x))$. Then
$$
A_x(P_x) \i W(f(x)).
\leqno (12.28)
$$

\ms
Fix $x$ as in the statement, and denote by $P$ and $W$
the vector spaces parallel to $P_x$ and $W(f(x))$ respectively.
Since we already know that $A_x(x)=f(x) \in W(f(x))$, we just need to check
that $A(P) \i W$, where $A$ denotes the linear part of $A_x$.

First we want to find, for each small $\rho > 0$, a
collection of affinely independent points 
$w_k^\rho \in X_3 \cap B(x,\rho)$, $1 \leq k \leq d$, 
and more precisely such that if we set $w_0 = x$, 
then for $1 \leq k \leq d$,
$$
\dist(w_k^\rho, P(w_0^\rho, \ldots, w_{k-1}^\rho)) \geq c \rho,
\leqno (12.29)
$$
where $P(w_0^\rho, \ldots, w_{k-1}^\rho)$ denotes the affine subspace 
of dimension $k-1$ spanned by $w_0^\rho, \ldots, w_{k-1}^\rho$ and $c>0$
is a constant that depends only on $M$ and $n$.

The proof of existence is the same as for (8.35); we just need to 
know that $E$ is locally Ahlfors-regular and that for $\rho > 0$ 
small enough,
$$
\H^d(X_3 \cap B(x,\rho)) \geq C^{-1} \rho^d.
\leqno (12.30)
$$
This last follows from Proposition 4.1 
(the local Ahlfors-regularity of $E$), (11.33),
and the fact that $X_3(s) \i X_3$ by (11.31). 
Alternatively (if you don't like (8.35)), we could take $d$ 
affinely independent points in $P_x$ (with a property like (12.29),
even with $c=1/2$), recall that $P_x$ was initially defined as the 
tangent plane to the $\Gamma_s$ that contains $x$, and use (11.32) 
to find points of $X_3$ that lie close enough to these points.
Either way we can find the $w_k^\rho$.

Take a sequence of radii $\rho$ that tends to $0$, and for which 
$$
\hbox{each $\rho^{-1}(w_k^\rho-x)$, $0 \leq k \leq d$,
has a limit $w_k$.}
\leqno (12.31)
$$
Then $w_0 = 0$ and (12.29) yields
$$
\dist(w_k, P(0, \ldots, w_{k-1})) \geq c 
\leqno (12.32)
$$
for $1 \leq k \leq d$, and so the $w_k$, $1 \leq k \leq d$,
are linearly independent. In addition, $w_k \in P$, 
because
$$\eqalign{
\dist(\rho^{-1}(w_k^\rho-x), P) &= \rho^{-1} \dist(w_k^\rho-x,P)
= \rho^{-1} \dist(w_k^\rho,P_x) 
\cr&
= \rho^{-1} o(|w_k^\rho-x|) = o(1)
}\leqno (12.33)
$$
by (11.34) (which holds because $x\in X_3$)
and because $w_k^\rho \in X_3 \i E$.
So the $w_k$, $1 \leq k \leq d$, are a basis of $P$, 
and it is enough to check that $A(w_k) \in W$
for $1 \leq k \leq d$. But
$$
\eqalign{
\dist(A(w_k^\rho-x), W)
&= \dist(A_{x}(w_k^\rho)- A_{x}(x), W)
\cr&
= \dist(A_{x}(w_k^\rho), W+A_{x}(x))
= \dist(A_{x}(w_k^\rho), W(f(x)))
\cr&
\leq \dist(f(w_k^\rho), W(f(x))) + |f(w_k^\rho)-A_{x}(w_k^\rho)|
}\leqno (12.34)
$$
because $A_x$ is affine and $A_{x}(x) = f(x) \in W(f(x))$,
and by definition of $W$ (the vector space parallel to $W(f(x))$).
But $|f(w_k^\rho)-f(x)| \leq \rho |f|_{lip} < \delta_2$ 
if $\rho$ is small enough, so $F(f(w_k^\rho)) = F(f(x))$ by (11.27) 
(recall that $x\in X_5$ and $w_k^\rho \in X_3$, 
so they both lie in $X_2$), and then 
$f(w_k^\rho) \in F(f(w_k^\rho)) = F(f(x)) \i W(f(x))$,
by definitions. Hence (12.34) yields
$$
\dist(A(w_k^\rho-x), W) \leq |f(w_k^\rho)-A_{x}(w_k^\rho)|
= o(|w_k^\rho-x|) = o(\rho)
\leqno (12.35)
$$
by (11.40). We divide by $\rho$, take the limit, use (12.31), 
and get that $A(w_k) \in W$, as needed; Lemma~12.27 follows.
\qed

\ms
$\dagger$
Under the Lipschitz assumption, we shall use an analogue of 
Lemma~12.27 where $f$ is replaced with the mapping $\wt f$ 
defined by
$$
\wt f(x) = \psi(\lambda f(x)) \in B(0,1)
\ \hbox{ for } x\in E
\leqno (12.36)
$$
(which is the same as (11.50)). Recall that $f : E \to U$
(because $f = \varphi_1$ on $E$, and because $\wh W \i \i U$),
and then $\psi : \lambda U \to B(0,1)$ by definition of the Lipschitz 
assumption. In fact, $\wt f$ is even defined in a neighborhood of $E$
in $U$.

In the definition of $X_4(s)$, and in addition to the 
differentiability at $x$ of the restriction of $f$ to
$\Gamma_s$ (see below (11.37)), let us also require
that the restriction of $\wt f$ to $\Gamma_s$ be 
differentiable at $x$. This is also true for $\H^d$ 
almost every $x\in X_3(s)$, and the proof of (11.36) and (11.37)
shows that there is an affine function $\wt A_x : \R^n \to \R^n$,
of rank at most $d$, such that
$$
|D \wt A_x| \leq \lambda \Lambda |f|_{lip},
\leqno (12.37)
$$
and 
$$
\lim_{y \to x \, ; \, y \in \Gamma_s} {|\wt f(y)- \wt A_x(y)| \over |y-x|} = 0,
\leqno (12.38) 
$$
where $\wt f$ is defined, near $x$, by the same formula (12.36) as 
above. Then the proof of (11.40) yields 
$$
\lim_{y \to x \, ; \, y \in E \cup P_x} {|\wt f(y)- \wt A_x(y)| \over |y-x|} = 0
\leqno (12.39) 
$$
when $x\in X_4(s)$, if we add this requirement to the definition of $X_4(s)$.
Then we get the following variant of Lemma~12.27.

\ms\proclaim Lemma 12.40.
For $x \in X_5$, let $P_x$ and $\wt A_x$ be as in (11.34) and (12.39),
denote by $F(\wt f(x))$ the smallest face of our usual rigid dyadic grid 
that contains $\wt f(x) = \psi(\lambda f(x))$ 
and by $W(\wt f(x))$ the smallest affine subspace 
that contains $F(\wt f(x))$. Then
$$
\wt A_x(P_x) \i W(\wt f(x)).
\leqno (12.41)
$$

\ms
The proof is the same as for Lemma~12.27; we just replace $f$
with $\wt f$ and (11.40) with (12.39) in (12.35).
\qed $\dagger$

\msi{\bf Step 2.e. We choose disks $Q_j$, $j\in J_1$.}
\ms
For each $j\in J_1$, we set $B_j = B(x_j,t)$ (see Step 2.a).
We want to choose a $Q_j$ in our large collection $\cal F$, 
and not so far from $f(E \cap B_j)$
so that composing $f$ with a projection on $Q_j$
will not move points too much. As before, we start with the
easier rigid assumption.

So fix $j\in J_1$. Recall that 
$x_j \in X_N(\delta_4) = X_5 \cap f^{-1}(Y_N(\delta_4))$
(see the line below (12.8), and (12.3)). Then
$f(x_j) \in Y_N(\delta_4)$ and by (12.13) we can find
$l = l(j) \in {\cal L}$ such that $f(x_j) \in D_l = B(y_l,t)$.
Since $y_l \in Y_N(\delta_4) \i f(X_5)$ by (12.1),
we can find $x(l) \in X_5$ such that $y_l = f(x(l))$.

Still denote by $F_l$ the smallest face that contains
$y_l$ and by $W(y_l)$ the affine space spanned by $F_l$.
With the notation of (11.27), $F_l = F(y_l) =
F(f(x(l)))$, and (11.27) says that $F_l = F(f(x_j))$ as well,
since $x_j$ and $x(l)$ both lie in $X_5 \i X_2$
and $|f(x_j)-f(x(l))| = |f(x_j)-y_l| \leq t < \delta_6
< \delta_2$ (see (12.8) and (12.7)). 
Then $W(y_l) = W(f(x_j))$ (the affine space spanned by $F(f(x_j))$
as well.

Set $P_j = P_{x_j}$; then Lemma 12.27 says that
$$
A_{x_j}(P_j) \i W(f(x_j)) = W(y_l).
\leqno (12.42)
$$
We claim that that we can find a disk
$Q_j \in {\cal F}_l$ such that
$$
\dist(z,Q_j) \leq \alpha t
\ \hbox{ for } z\in A_{x_j}(P_j) \cap B(y_l,3 (1+|f|_{lip}) t). 
\leqno (12.43) 
$$
Indeed, if $W(y_l)$ is of dimension at most $d$, 
${\cal F}_l$ is composed of a single element 
$Q = W(y_l) \cap B(3 (1+|f|_{lip}) t)$, 
which satisfies (12.43) by (12.42) and (12.17). 
Otherwise, (12.42) (and the fact that $A_{x_j}(P_j)$
is at most $d$-dimensional) allows us to pick an affine $d$-plane $P'$
such that $A_{x_j}(P_j) \i P' \i W(y_l)$. Observe that $P'$
goes through $D_l$ because $A_{x_j}(x_j) = f(x_j) \in D_l$, 
so we can choose $Q_j \in {\cal F}_l$ so that (12.18) holds; 
(12.43) follows because $A_{x_j}(P_j) \i P'$.

Let us check that
$$
\dist(f(z),Q_j) \leq \big(2\varepsilon +2\varepsilon |f|_{lip} +\alpha \big) \, t
\ \hbox{ for } z\in E \cap 2B_j = E \cap B(x_j,2t).
\leqno (12.44)
$$
Indeed, $x_j \in X_5$, so (11.46) holds for $x=x_j$.
We get that
$$
|f(z)-A_{x_j}(z)| \leq \varepsilon |z-x_j| \leq 2\varepsilon t
\ \hbox{ for } z\in E \cap 2B_j,
\leqno (12.45)
$$
because $|z-x_j| \leq 2t < 2\delta_6 < \delta_3$
by (12.8) and (12.7). Denote by $\pi(z)$ the projection of $z$ 
on $P_j$; then
$|\pi(z)-z| \leq 2\varepsilon t$ by (11.45) (again applied with 
$x=x_j$ and valid because $|z-x_j| < \delta_3$). By (11.36),
$$
|A_{x_j}(z)-A_{x_j}(\pi(z))| \leq |\pi(z)-z| \, |f|_{lip}
\leq 2\varepsilon t \, |f|_{lip}.
\leqno (12.46)
$$
In addition,
$$\eqalign{
|A_{x_j}(\pi(z))-y_l|
&\leq |A_{x_j}(\pi(z))-f(x_j)|+|f(x_j)-y_l|
\cr&
= |A_{x_j}(\pi(z))-A_{x_j}(x_j)|+|f(x_j)-y_l|
\cr&
\leq |f|_{lip} |\pi(z)-x_j| + t
\leq (2|f|_{lip} + 1) t
}\leqno (12.47)
$$
because $A_{x_j}(x_j) = f(x_j)$, $f(x_j) \in D_l = B(y_l,t)$, 
and $|\pi(z)-x_j| \leq |z-x_j| \leq 2t$. Since
$A_{x_j}(\pi(z)) \in A_{x_j}(P_j)$ trivially,
(12.47) allows us to apply (12.43) to it; this yields 
$\dist(A_{x_j}(\pi(z)), Q_j) \leq \alpha t$, and now
(12.44) follows from (12.45) and (12.46).

\ms
$\dagger$
Under the Lipschitz assumption, we still can define 
$l = l(j) \in {\cal L}$ such that $f(x_j) \in D_l$, 
and $x(l)\in X_5$ such that $y_l = f(x(l))$. 
The smallest face $F_l$ that contains
$y_l$ is still the same as for $f(x_j)$; equivalently,
the smallest face $\wt F_l$ of the rigid grid that
contains $\wt y_l = \psi(\lambda y_l)$
is the same as for $\wt f(x_j) =  \psi(\lambda f(x_j))$.
Then the smallest affine space $W(\wt y_l)$ that contains
$\wt F_l$ is the same as for $\wt f(x_j)$.
This time, we use Lemma 12.40, which says that
$$
\wt A_{x_j}(P_{x_j}) \i W(\wt f(x_j)) = W(\wt y_l),
\leqno (12.48)
$$
and we want to use this to find a close enough 
$Q_j \in {\cal F}_l$. If $F(y_l)$ is of dimension at 
most $d$, we pick the only element $Q$ defined by (12.22)
with $\wt P = W(\wt y_l)$. Otherwise, we want to use 
the rule (12.25) to select a $\wt P$. First use (12.48) 
to choose $\wt P'$, of dimension $d$, such that 
$\wt A_{x_j}(P_{x_j}) \i \wt P' \i W(\wt y_l)$,
and observe that $\wt P'$ meets
$B(\wt y_l, \lambda \Lambda (1+|f|_{lip}) t))$
because $\wt f(x_j) = \wt A_{x_j}(x_j) \in \wt A_{x_j}(P_{x_j})$
and 
$$
|\wt f(x_j)-\wt y_l| = |\psi(\lambda f(x_j))-\psi(\lambda y_l)|
\leq \lambda \Lambda |f(x_j)-y_l| \leq \lambda \Lambda t
\leqno (12.49)
$$
because $f(x_j) \in D_l$. Thus (12.24) holds, and we can find 
a $d$-plane $\wt P$ in our collection, so that (12.25) holds.
We call that plane $\wt P_j$.
Hence
$$
\dist(z,\wt P_j) \leq \lambda \Lambda^{-1} \alpha t
\ \hbox{ for } z\in \wt A_{x_j}(P_{x_j})
\cap B(\wt y_l, 10\lambda \Lambda^2 (1+|f|_{lip}) t),
\leqno (12.50) 
$$
since $\wt A_{x_j}(P_{x_j}) \i \wt P'$.
Then let $Q_j\in {\cal F}_l$ be the set $Q$ defined by (12.22)
with this choice of $\wt P = \wt P_j$; we want to check that 
$$
\dist(f(z),Q_j) \leq 
10 \Lambda^3 (1+|f|_{lip}) \varepsilon  t + \alpha \, t
\ \hbox{ for } z\in E \cap B(x_j,10\Lambda t),
\leqno (12.51)
$$
as in (12.44). For this, we assume that in the definition of $X_5(s)$
(near (11.42)), we added to the requirement (11.46) its analogue
for $\wt f$. That is, we first take $\delta_3$ so small that
$f(y) \in U$ for $x\in X_1$ and $y\in B(x,\delta_3)$
(this is easy, because $f$ is Lipschitz, $X_1$ is compact, and 
$f(X_1) \i U$); this way, we can define 
$\wt f(y) = \psi(\lambda f(y))$ for $x\in X_1$ and $y\in B(x,\delta_3)$, 
as in (12.36). But more importantly, we take $\delta_3$ so small that
with this definition,
$$
|\wt f(y)- \wt A_x(y)| \leq \lambda\varepsilon |y-x|
\ \hbox{ for } y\in [E \cup P_x \cup \Gamma_s] \cap B(x,\delta_3).
\leqno (12.52)
$$
This is possible, for the same reason as for (11.46).

Now let $z\in E \cap B(x_j,10\Lambda t)$ be given; then
$$
|\wt f(z)-\wt A_{x_j}(z)| 
\leq \lambda \varepsilon |z-x_j| 
\leq 10 \lambda \Lambda\varepsilon t
\leqno (12.53)
$$
by (12.52) and because 
$|z-x_j| \leq 10 \Lambda t < 10 \Lambda \delta_6 < \delta_3$
by (12.8) and (12.7). Again the projection $\pi(z)$ of $z$ on $P_j$
is such that $|\pi(z)-z| \leq 10 \Lambda\varepsilon t$, by (11.45), so
$$
|\wt A_{x_j}(z)- \wt A_{x_j}(\pi(z))| \leq \lambda \Lambda |f|_{lip} |\pi(z)-z| 
\leq 10 \lambda \Lambda^2 |f|_{lip} \varepsilon t 
\leqno (12.54)
$$
by (12.37), and also
$$\eqalign{
|\wt A_{x_j}(\pi(z))- \wt y_l|
&\leq |\wt A_{x_j}(\pi(z))- \wt f(x_j)|+|\wt f(x_j)- \wt y_l|
\cr&
= |\wt A_{x_j}(\pi(z))- \wt A_{x_j}(x_j)|
+|\psi(\lambda f(x_j))- \psi(\lambda y_l)|
\cr&
\leq \lambda \Lambda |f|_{lip} |\pi(z)-x_j| 
+ \lambda \Lambda |f(x_j)- y_l|
\cr&
\leq \lambda \Lambda |f|_{lip} |z-x_j| 
+ \lambda \Lambda t
\leq 10\lambda\Lambda^2 (1+|f|_{lip}) t
}\leqno (12.55)
$$
because $\wt A_{x_j}(x_j) = \wt f(x_j)$, by (12.37) again,
because $x_j = \pi(x_j)$ since $P_j$ goes through $x_j$,
and because $f(x_j) \in D_l = B(y_l,t)$ and $z\in B(x_j,10\Lambda t)$.

Now $\wt A_{x_j}(\pi(z)) \in \wt A_{x_j}(P_j)$, so by
(12.55) we may apply (12.50) to it; we get that
$\dist(\wt A_{x_j}(\pi(z)),\wt P_j) \leq \lambda \Lambda^{-1} 
\alpha t$. Since
$$
|\wt f(z)-\wt A_{x_j}(\pi(z))| \leq 10 \lambda \Lambda^2 
(1+|f|_{lip}) \varepsilon t
\leqno (12.56)
$$
by (12.53) and (12.54), we get that
$$
\dist(\wt f(z),\wt P_j) 
\leq 10 \lambda \Lambda^2 (1+|f|_{lip}) \varepsilon t 
+\lambda \Lambda^{-1} \alpha t
\ \hbox{ for } z\in E \cap B(x_j,10\Lambda t).
\leqno (12.57)
$$
In addition, 
$|\wt f(z)-\wt y_l| \leq 11\lambda \Lambda^2 (1+|f|_{lip}) t$
by (12.55) and (12.56), so the point of $\wt P_j$ that
minimizes the distance to $\wt f(z)$ automatically lies in 
$\wt Q = \psi(\lambda Q_j)$ by (12.22). Finally,
$$\eqalign{
\dist(f(z),Q_j) &= \dist(f(z),\lambda^{-1}\psi^{-1}(\wt Q))
\leq \lambda^{-1} \Lambda \dist(\wt f(z),\wt Q)
\cr&
= \lambda^{-1} \Lambda \dist(\wt f(z),\wt P_j)
\leq 10 \Lambda^3 (1+|f|_{lip}) \varepsilon t  + \alpha t
}\leqno (12.58)
$$
by (12.22) again and (12.57); (12.51) follows. 
$\dagger$

\msi{\bf Step 2.f. We construct mappings $g_j$, $j\in J_1$.}
\ms
Return to the rigid assumption.
For each $j\in J_1$, we now define a Lipschitz mapping
$g_j : U \to \R^n$. We use a new constant $a \in (0,1)$
quite close to $1$. Set 
$$
g_j(x) = f(x)
\ \hbox{ for } x \in U \sm B(x_j,{1+a \over 2}\, t)
\leqno (12.59)
$$
and
$$
g_j(x) = \pi_j(f(x))
\ \hbox{ for } x \in aB_j = B(x_j,at),
\leqno (12.60)
$$
where $\pi_j$ denotes the orthogonal projection onto
the affine plane spanned by $Q_j$. In the middle, interpolate
linearly as usual, by setting
$$
g_j(x) = {2|x-x_j|-2at \over (1-a)t}\, f(x) +
{(1+a)t -2|x-x_j|\over (1-a)t}\, \pi_j(f(x))
\leqno (12.61)
$$
for $x \in B(x_j,{1+a \over 2}\, t) \sm B(x_j,at)$.
This gives a Lipschitz function $g_j$, with a quite large Lipschitz 
norm that we don't want to compute, and such that
$$\eqalign{
||g_j - f||_\infty 
&\leq \sup_{x\in B_j} |\pi_j(f(x))-f(x)|
\leq \sup_{x\in B_j} \dist(f(x),Q_j)
\cr&
\leq \dist(f(x_j),Q_j) + \sup_{x\in B_j} |f(x)-f(x_j)|
\leq \dist(f(x_j),Q_j) + |f|_{lip} \, t 
\cr&
\leq \big(2\varepsilon +2\varepsilon |f|_{lip} +\alpha \big) \, t
+ |f|_{lip} t 
\leq (1+ |f|_{lip}) t
}\leqno (12.62)
$$
by (12.44) and if $\varepsilon$ and $\alpha$ are small enough. 
Fortunately, the estimates get better near $E$. Set
$$
E^{\varepsilon t} = \big\{ x\in U \, ; \, \dist(x,E) \leq 
\varepsilon t \big\};
\leqno (12.63)
$$
we claim that
$$
|g_j(x) - f(x)| \leq (2\varepsilon +3|f|_{lip} \, \varepsilon + \alpha)\,  t
\ \hbox{ for } x\in E^{\varepsilon t}.
\leqno (12.64)
$$
Indeed, by (12.59) we can assume that 
$x \in B(x_j,{1+a \over 2}\, t)$; choose $z \in E$
such that $|z-x| \leq \varepsilon t$; then $z\in 2B_j$, and
$$\eqalign{
|g_j(x) - f(x)| &\leq |\pi_j(f(x))-f(x)|
\leq \dist(f(x),Q_j)
\cr&
\leq \dist(f(z),Q_j)+ |f(x)-f(z)|
\cr&
\leq \big(2\varepsilon +2\varepsilon |f|_{lip} +\alpha \big) \, t
+ |f|_{lip} \varepsilon t
= (2\varepsilon +3 \varepsilon |f|_{lip} + \alpha)\, t
}\leqno (12.65)
$$
as above, and by (12.44); the claim follows. 
Similarly, let us check that
$$
g_j \ \hbox{ is $(1+|f|_{lip})$-Lipschitz on } E^{\varepsilon t}.
\leqno (12.66)
$$
Let $x, y \in E^{\varepsilon t}$, and use (12.59)-(12.61)
to write 
$$
g_j(x) = \beta(x) f(x) + (1-\beta(x)) \pi_j(f(x))
\leqno (12.67)
$$
for some $\beta(x) \in [0,1]$, and similarly
$g_j(y) = \beta(y) f(y) + (1-\beta(y)) \pi_j(f(y))$.
Write 
$$
g_j(y) = \beta(x) f(y) + (1-\beta(x)) \pi_j(f(y))
- [\beta(x)-\beta(y)][f(y)-\pi_j(f(y))], 
\leqno (12.68)
$$
and then subtract (12.68) from (12.67) to get that
$$\eqalign{
g_j(x)-g_j(y)
&= \beta(x) [f(x)-f(y)] +(1-\beta(x)) [\pi_j(f(x))-\pi_j(f(y))]
\cr&\hskip 4.5cm
+ [\beta(x)-\beta(y)][f(y)-\pi_j(f(y))].
}\leqno (12.69)
$$
The first part is at most
$$
\beta(x) |f(x)-f(y)| +(1-\beta(x)) |\pi_j(f(x))-\pi_j(f(y))|
\leq |f(x)-f(y)| \leq |f|_{lip} |x-y|
$$
because $\pi_j$ is $1$-Lipschitz. 
If $y \in B(x_j,{1+a \over 2}\, t)$, the proof of (12.65)
shows that 
$$
|\pi_j(f(y))-f(y)| \leq 
(2\varepsilon +3 \varepsilon |f|_{lip} + \alpha)\, t
\leqno (12.70)
$$
and since $\displaystyle |\beta(x)-\beta(y)| \leq {2 |x-y| \over (1-a)t}$ 
by (12.59)-(12.61), we get that
$$
|\beta(x)-\beta(y)||f(y)-\pi_j(f(y))|
\leq (2\varepsilon +3 \varepsilon |f|_{lip}
+ \alpha) \, {2 |x-y| \over (1-a)} \leq |x-y|
\leqno (12.71)
$$
if $\varepsilon$ and $\alpha$ are small enough, depending on
$|f|_{lip}$, but also on $a$. (This is all right, see Remark 11.17.)
Altogether
$|g_j(x)-g_j(y)| \leq (1+|f|_{lip}) |x-y|$, as needed 
for (12.66) in this first case.

If $x \in B(x_j,{1+a \over 2}\, t)$, the same argument, with $x$ and
$y$ exchanged from the start (i.e., also in (12.68) and (12.69)) gives
the desired result. Finally, if both $x$ and $y$ lie
out of $B(x_j,{1+a \over 2}\, t)$, then $\beta(x)=\beta(y) = 1$ by 
(12.59), and $|g_j(x)-g_j(y)| \leq |f|_{lip} |x-y|$ directly by (12.70).
This completes our proof of (12.66).

\msi $\dagger${\bf Step 2.g. The mappings $g_j$, under the Lipschitz assumption.}$\dagger$
\ms

$\dagger$ Now we do the same thing under the Lipschitz assumption.
We shall try to do the linear algebra and convex combinations 
on $B(0,1)$, because we want to preserve the faces when we can, 
but (later on, when we interpolate between the $g_j$) 
we shall still use partitions of unity defined on $U$.

Before we define mappings $\wt g_j$ we need to extend our definition of
the $\wt f$ of (12.36). We are particularly interested in the set
$$
U_{int} = \big\{ x\in U \, ; \, \dist(x,X_0) \leq {\delta_0 \over 
4(1+|f|_{lip})} \big\},
\leqno (12.72)
$$
because we shall see that it is so small that (12.36) makes sense on 
it, and sufficiently large to contain the $2B_j$, $j\in J_1$.
Let us first check that
$$
\dist(f(x),\wh W) \leq {\delta_0 \over 3}
\ \hbox{ for } x\in U_{int}.
\leqno (12.73)
$$
For $x\in U_{int}$, pick $y\in X_0$ such that
$|y-x| \leq {\delta_0 \over 3(1+|f|_{lip})}$; then $f(y) \in \wh W$ 
by (11.20), (2.1), and (2.2), so
$$
\dist(f(x),\wh W) \leq |f(x)-f(y)|
\leq |f|_{lip} |y-x| \leq {\delta_0 \over 3}
\leqno (12.74)
$$
as needed. For such $x$, $f(x) \in U$ because
$\delta_0 = \dist(\wh W, \R^n \sm U)$ by (12.6). Hence we can define 
$\psi(\lambda f(x))$. So we can extend the definition (12.36), and set 
$$
\wt f(x) = \psi(\lambda f(x)) 
\ \hbox{ for } x\in U_{int}\cup E.
\leqno (12.75)
$$
Note that $\wt f(x) \in B(0,1)$ automatically, because 
$\psi(\lambda U) = B(0,1)$. It will also be good to know that
$$
2B_j = B(x_j,2t) \i U_{int} \ \hbox{ for } j\in J_1,
\leqno (12.76)
$$
which is true because $x_j \in X_5 \i X_0$ and 
$t < \delta_6 \leq {\delta_0 \over 10(1+|f|_{lip})}$ 
by (12.8) and (12.7).

Next we define intermediate mappings $\wt g_j$. 
Recall that $Q_j$ is defined by (12.22) for some affine plane 
$\wt P = \wt P_j$; we denote by $\wt \pi_j$ the orthogonal 
projection onto $\wt P_j$, and then set
$$
\wt g_j(x) = \wt f(x)
\ \hbox{ for } x \in U_{int} \sm B(x_j,{1+a \over 2}\, t)
\leqno (12.77)
$$
(a little as in (12.59)),
$$
\wt g_j(x) = \wt\pi_j(\wt f(x))
\ \hbox{ for } x \in aB_j = B(x_j,at),
\leqno (12.78)
$$
(as in (12.60)), and
$$
\wt g_j(x) = {2|x-x_j|-2at \over (1-a)t}\, \wt f(x) +
{(1+a)t -2|x-x_j|\over (1-a)t}\, \wt\pi_j(\wt f(x))
\leqno (12.79)
$$
for $x \in B(x_j,{1+a \over 2}\, t) \sm B(x_j,at)$
(as in (12.61)). The simplest for us will be not to define 
$\wt g_j$ in $U \sm U_{int}$.

Let us concentrate on what happens in $2B_j = B(x_j,2t)$. Recall that
$2B_j \i U_{int}$ by (12.76), so $\wt g_j$ is defined on $2B_j$.
Then, for $x\in 2B_j$,
$$\eqalign{
|\wt f(x)-\wt y_l| &= |\psi(\lambda f(x)) - \psi(\lambda y_l)|
\leq \lambda \Lambda |f(x)-y_l|
\cr&
\leq \lambda \Lambda \big(|f(x)-f(x_j)|+|f(x_j)-y_l|\big)
\cr&
\leq \lambda \Lambda \big(|f|_{lip}|x-x_j|+t\big)
\leq \lambda \Lambda \big(2|f|_{lip}+1\big)\, t
}\leqno (12.80)
$$
by (12.75), the definition of $\wt y_l$ above (12.48),
and the fact that $f(x_j) \in D_l$. Hence
$$\leqalignno{
|\wt g_j(x) - \wt f(x)|
&\leq |\wt\pi_j(\wt f(x)) - \wt f(x)|
= \dist(\wt f(x), \wt P_j)
\leq  \dist( \wt y_l, \wt P_j) + |\wt f(x)- \wt y_l|
\cr&
\leq  2 \lambda \Lambda (1+|f|_{lip}) \, t 
+ |\wt f(x)-\wt y_l|
\leq 4 \lambda \Lambda (1+|f|_{lip}) \, t
& (12.81)
}
$$
by (12.77)-(12.79), because $\wt P_j$ was chosen
(near (12.49)) so that (12.23) holds, and by (12.80). 
Let us also check that
$$
\wt g_j(x) \in B(0,1) \ \hbox{ for } x\in 2B_j.
\leqno (12.82)
$$
By (12.80) and (12.81),
$|\wt g_j(x) - \wt y_l| \leq 6\lambda \Lambda (1+|f|_{lip})\, t$.
But $y_l \in X_0 \i W_1$, so 
$$\eqalign{
\dist(\wt y_l,\d B(0,1))
&= \dist(\psi(\lambda y_l),\psi(\lambda \d U))
\geq \lambda \Lambda^{-1} \dist(y_l,\d U)
\cr&
\geq \lambda \Lambda^{-1} \dist(\wh W,\R^n \sm U) 
= \lambda \Lambda^{-1} \delta_0
\cr&
\geq 10 (1+|f|_{lip})\, \lambda \Lambda \delta_6
\geq 10 (1+|f|_{lip})\, \lambda \Lambda t
}\leqno (12.83)
$$
because $\wt y_l = \psi(\lambda y_l)$ (see above (12.22)),
because $\psi$ has a bilipschitz extension from the closure of 
$\lambda U$ to $\overline B(0,1)$, and by (12.6)-(12.8). Then 
$$
\dist(\wt g_j(x),\R^n \sm B(0,1)) 
\geq \dist(\wt y_l,\d B(0,1)) - |\wt g_j(x) - \wt f(x)|
\geq 4(1+|f|_{lip})\, \lambda \Lambda t > 0,
\leqno (12.84)
$$
and (12.82) holds.

We may now define $g_j$ on $2B_j$ by
$$
g_j(x) = \lambda^{-1}\psi^{-1}(\wt g_j(x))
\ \hbox{ for } x\in 2B_j,
\leqno (12.85)
$$
because $\wt g_j(x) \in B(0,1)$ and hence $\psi^{-1}(\wt g_j(x))$
is defined. We also get that $g_j(x) \in U$, because 
$\psi : \lambda U \to B(0,1)$.

Observe that when $x\in 2B_j \sm B(x_j,{1+a \over 2}\, t)$, (12.85), 
(12.77) and (12.75) yield $g_j(x) = \lambda^{-1}\psi^{-1}(\wt g_j(x))
= \lambda^{-1}\psi^{-1}(\wt f(x)) = f(x)$. So we can safely set
$$
g_j(x) = f(x)
\ \hbox{ for } x\in U \sm B(x_j,{1+a \over 2}\, t),
\leqno (12.86)
$$
the two definitions coincide on $2B_j \sm B(x_j,{1+a \over 2}\, t)$,
and $g_j$ is Lipschitz on $U$.

Return to $x\in 2B_j$. Since $f(x) = \lambda^{-1}\psi^{-1}(\wt f(x))$
by (12.75), we see that
$$
||g_j - f||_{L^\infty(2B_j)} 
\leq \lambda^{-1} \Lambda ||\wt g_j - \wt f||_{L^\infty(2B_j)}
\leq 4 \Lambda^2 (1+|f|_{lip}) \, t,
\leqno (12.87)
$$
by (12.85) and (12.81). This will be a good enough analogue for (12.62).

We also need better estimates when $x\in E^{\varepsilon t}$
(the small neighborhood of $E$ defined in (12.63)). Let 
$x \in E^{\varepsilon t} \cap B(x_j,{1+a \over 2}t)$
be given, and pick $z\in E$ such that $|z-x| \leq \varepsilon t$; then
$$\eqalign{
|\wt g_j(x) - \wt f(x)|
&\leq  |\wt \pi_j(\wt f(x))- \wt f(x)|
\leq  \dist( \wt f(x), \wt P_j)
\cr&
\leq  \dist( \wt f(z), \wt P_j) +  |\wt f(x)- \wt f(z)|
\leq  \dist( \wt f(z), \wt P_j)  +  |x- z| \, |\wt f|_{lip} 
\cr&
\leq  10 \lambda \Lambda^2 (1+|f|_{lip}) \varepsilon t 
+\lambda \Lambda^{-1} \alpha t  + \varepsilon t \lambda \Lambda |f|_{lip}
\cr&
\leq 11 \lambda \Lambda^2 (1+|f|_{lip}) \varepsilon t 
+\lambda \Lambda^{-1} \alpha t 
}\leqno (12.88)
$$
by (12.77)-(12.79), because $\wt \pi_j$ denotes the projection
on the plane $\wt P_j$ that was used to construct $Q_j$
(see near (12.22) and (12.77)), by (12.57), and 
because $|\wt f|_{lip} \leq \lambda \Lambda |f|_{lip}$ by (12.75). 
Then (12.85) yields
$$
|g_j(x) - f(x)| \leq \lambda^{-1} \Lambda |\wt g_j(x) - \wt f(x)|
\leq 11 \Lambda^3 (1+|f|_{lip}) \varepsilon t + \alpha t,
\leqno (12.89)
$$
which is a good replacement for (12.65). But $g_j(x) = f(x)$ for 
$x \in E^{\varepsilon t} \sm \cap B(x_j,{1+a \over 2}\, t)$ 
(by (12.86), so 
$$
|g_j(x) - f(x)| \leq 11 \Lambda^3 (1+|f|_{lip}) \varepsilon t + \alpha t
\ \hbox{ for } x \in E^{\varepsilon t},
\leqno (12.90)
$$
which is an acceptable analogue of (12.64).

Next we copy the proof of (12.66). Let us estimate 
$|\wt g_j(x)- \wt g_j(y)|$ when $x, y \in E^{\varepsilon t} \cap 2B_j$.
Then we can use (12.77)-(12.79); as before, we write $\wt g_j(x)$ 
as a linear combination of $\wt f(x)$ and $\wt \pi_j(\wt f(x))$, 
and similarly for $\wt g_j(y)$, and then we compute as in (12.67)-(12.71).

As in (12.69), we get that
$\wt g_j(x)- \wt g_j(y) = A+B$, with
$$
A = \beta(x) [\wt f(x)- \wt f(y)] +(1-\beta(x)) 
[\wt \pi_j( \wt f(x))- \wt \pi_j(\wt f(y))]
\leqno (12.91)
$$
and 
$$
B = [\beta(x)-\beta(y)][\wt f(y)- \wt\pi_j(\wt f(y))].
\leqno (12.92)
$$
As before, $|A| \leq |\wt f(x)- \wt f(y)| 
\leq \lambda\Lambda |f|_{lip} |x-y|$
because $\wt\pi_j$ is $1$-Lipschitz and $\wt f$
is $\lambda\Lambda |f|_{lip}$ -Lipschitz.

In the first case when 
$y \in E^{\varepsilon t} \cap B(x_j,{1+a \over 2}\, t)$,
the proof of the second part of (12.88) yields
$$
|\wt f(y)-\wt\pi_j(\wt f(y))| = \dist( \wt f(y), \wt P_j)
\leq 11 \lambda \Lambda^2 (1+|f|_{lip}) \varepsilon t 
+\lambda \Lambda^{-1} \alpha t;
\leqno (12.93)
$$
since $\displaystyle |\beta(x)-\beta(y)| \leq {2 |x-y| \over (1-a)t}$ 
as before, we get that
$$
|B| \leq \big\{11 \lambda \Lambda^2 (1+|f|_{lip}) \varepsilon t 
+\lambda \Lambda^{-1} \alpha t \big\}  \, {2 |x-y| \over (1-a) t} 
\leq \lambda \Lambda^{-1} |x-y|
\leqno (12.94)
$$
if $\varepsilon$ and $\alpha$ are small enough. In this case
$$
|\wt g_j(x)- \wt g_j(y)| \leq |A|+|B| 
\leq \lambda \Lambda |f|_{lip} |x-y| + \lambda \Lambda^{-1}|x-y|.
\leqno (12.95)
$$
The other two cases are treated as before, and we get that
$$
|\wt g_j(x)- \wt g_j(y)| \leq \lambda \, {1+\Lambda^2 |f|_{lip} \over \Lambda} |x-y|
\ \hbox{ for } x, y \in E^{\varepsilon t} \cap 2B_j.
\leqno (12.96)
$$
By this and (12.85), we get that
$$
\hbox{$g_j$ is } (1+\Lambda^2 |f|_{lip}) \hbox{-Lipschitz on } 
E^{\varepsilon t} \cap 2B_j.
\leqno (12.97)
$$

This is a good enough analogue of (12.66), which was the last 
estimate of Step 2.f. $\dagger$

\bigskip\noindent
{\bf 13. Step 2.h. We glue the mappings $g_j$, $j\in J_1$,
and get a first mapping $g$.}

\ms
We start a new section, but continue to take care of the places where
$f$ is very many-to-one. Now we want to use the functions $g_j$ 
that we just built to modify $f$ on a subset of 
$$
V = \bigcup_{j \in J_1} B(x_j,t) = \bigcup_{j \in J_1} B_j. 
\leqno (13.1)
$$
Part of the difficulty will be that the balls $B(x_j,t)$ are not disjoint 
(but fortunately they have the same radius). We shall apply
the same trick as in [D2], based on 
adapted partitions of unity. Set
$$
R_j = B_j \sm a B_j
= \big\{ z\in \R^n \, ; \, at \leq |z-x_j| < t \big\}
\leqno (13.2)
$$
for $j\in J_1$, and
$$
R = \bigcup_{j\in J_1} R_j \, .
\leqno (13.3)
$$
Then let $\varphi_j$ be a smooth function such that 
$$
0 \leq \varphi_j(x) \leq 1
\ \hbox{ and }
|\nabla \varphi_j(x)| \leq {2 \over (1-a) t}
\ \hbox{ for } x\in \R^n,
\leqno (13.4)
$$
$$
\varphi_j(x) = 1 \hbox{ for } x\in a B_j,
\hbox{ and } \ \varphi_j(x) = 0 \hbox{ for } x\in \R^n \sm B_j.
\leqno (13.5)
$$
Next put an arbitrary order on the set $J_1$, and set
$$
\psi_j(x) = \sup_{i\in J_1 ; \, i \leq j} \varphi_i(x)
- \sup_{i\in J_1 ; \, i < j} \varphi_i(x)
\leqno (13.6)
$$
for $j\in J_1$ and $x\in \R^n$ (and where the empty sup is zero). Clearly
$$
\sum_{i\in J_1 ; \, i \leq j} \psi_i(x) 
= \sup_{i\in J_1 ; \, i \leq j} \varphi_i(x).
\leqno (13.7)
$$
Finally set 
$$
\psi(x) = \sum_{i\in J_1} \psi_i(x) = \sup_{i\in J_1} \varphi_i(x);
\leqno (13.8)
$$
then 
$$
0 \leq \psi(x) \leq 1
\ \hbox{ for } x\in \R^n,
\leqno (13.9)
$$
$\psi$ is $C [(1-a) t]^{-1}$-Lipschitz (because the $B_j$ have
bounded overlap), and $\psi = 1$ on $\bigcup_{j\in J_1} aB_j$.
Because of our particular choice of functions, we get that
$$
\psi(x) = 0 \hbox{ for } x\in \R^n \sm \big[\bigcup_{j\in J_1} B_j \big]
= \R^n \sm V
\leqno (13.10)
$$
and 
$$\eqalign{
&\hbox{if $x\in a B_j \sm R$ for some $j\in J_1$, then one of the 
$\psi_i(x)$ is equal to $1$}
\cr&\hskip2cm
\hbox{ and the other ones are equal to $0$}
}\leqno (13.11)
$$
(where in fact $i$ is the first index in $J_1$ such that $x\in aB_i$, 
or equivalently $x\in B_i$).

\ms
We now use the $\psi_j$ to construct a mapping $g : U \to \R^n$. 
As usual, we first do the description under the rigid assumption. We set 
$$
g(x) = f(x) + \sum_{j \in J_1} \psi_j(x) [g_j(x)-f(x)],
\leqno (13.12)
$$
which we write like this because $\sum_{j \in J_1} \psi_j(x)$
may be smaller than $1$ at some points. Notice that
$$
||g - f||_\infty \leq (1+ |f|_{lip}) \, t
\leqno (13.13)
$$
by (12.62) and (13.9), and 
$$
|g(x) - f(x)| \leq \sum_{j\in J_1} \psi_j(x) |g_j(x)-f(x)|
\leq (2\varepsilon +3\varepsilon |f|_{lip} + \alpha)\,  t
\ \hbox{ for } x\in E^{\varepsilon t},
\leqno (13.14)
$$
by (12.64). We also claim that because of (12.66),
$$
g \ \hbox{ is $(2+3|f|_{lip})$-Lipschitz on } E^{\varepsilon t}.
\leqno (13.15)
$$
Indeed, since $f$ is Lipschitz, it is enough to show that $g-f$ 
is $(2+2|f|_{lip})$-Lipschitz on $E^{\varepsilon t}$ and estimate
$$
\Delta(x,y) = (g-f)(x)-(g-f)(y)
= \sum_{j \in J_1} \Big\{\psi_j(x) [g_j(x)-f(x)]
-\psi_j(y) [g_j(y)-f(y)] \Big\}
\leqno (13.16)
$$
for $x, y \in E^{\varepsilon t}$.
Since $\Delta(x,y) \leq 2 (2\varepsilon +3\varepsilon |f|_{lip} + 
\alpha)\,  t$ by the $L^\infty$ bound in (13.14), we may assume that 
$|x-y| \leq t/10$.

Let $j \in J_1$ be such that $g_j(x)-f(x) \neq 0$; then
$x\in B(x_j,{1+a \over 2}\, t)$ by (12.59), and $y\in {3 \over 2} B_j$
because $|x-y| \leq t/10$. Similarly, if $g_j(y)-f(y) \neq 0$,
then $y\in B(x_j,{1+a \over 2}\, t)$ and $x\in {3 \over 2} B_j$.
So $x,y \in {3 \over 2} B_j$ for every $j\in J_1$ that has a 
contribution to the right-hand of (13.16). There are at most $C$
indices $j$ like this (recall that $B_j = B(x_j,t)$ and that
$|x_i - x_j| \geq t/3$ when $i \neq j$; see above (12.9)),
and each contribution is estimated as follows. We write
$$\leqalignno{
&\big|\psi_j(x) [g_j(x)-f(x)] -\psi_j(y) [g_j(y)-f(y)] \big|
\cr& \hskip 1cm
\leq \psi_j(x) |g_j(x)-f(x)-g_j(y)+f(y)|
+ |\psi_j(x)-\psi_j(y)| \, |g_j(y)-f(y)|
\cr& \hskip 1cm
\leq \psi_j(x) |g_j(x)-f(x)-g_j(y)+f(y)|
+ 4 (1-a)^{-1} t^{-1} |x-y| \, |g_j(y)-f(y)|
&(13.17)
\cr& \hskip 1cm
\leq \psi_j(x) (1+2|f|_{lip}) |x-y|
+ 4 (1-a)^{-1} |x-y| (2\varepsilon +3\varepsilon |f|_{lip} + \alpha)
}
$$
because $\psi_j$ is $4[(1-a)t]^{-1}$-Lipschitz
(by (13.6) and (13.4)), and by (12.66), (12.64), and (13.14).

When we sum (13.17) over $j$, the first term gives a total contribution
of at most $(1+2|f|_{lip}) |x-y|$, by (13.8) and (13.9),
and the second one of at most 
$|x-y|$, if $\varepsilon$ and $\alpha$ are chosen small enough,
depending on $n$, $|f|_{lip}$ and $a$. So $|\Delta(x,y)| \leq 
(2+2|f|_{lip}) |x-y|$ by (13.16), and our claim (13.15) follows.

Let us record the fact that, by (13.12) and (12.59),
$$
g(x) = f(x) \ \hbox{ on } \R^n \sm \bigcup_{j\in J_1} B(x_j,{1+a \over 2}\, t).
\leqno (13.18)
$$
Also, we claim that
$$
g\Big(\bigcup_{j\in J_1} B_j \sm R\Big) \i \bigcup_{l\in {\cal L}}
\bigcup_{Q \in {\cal F}_l} Q = \bigcup_{Q \in {\cal F}} Q.
\leqno (13.19)
$$
Let $j\in J_1$ and $x\in B_j \sm R$ be given. Then 
$x\in aB_j$ (see (13.2) and (13.3)), and by (13.11) exactly one 
$\psi_i(x)$ is equal to $1$, and the other ones are equal to $0$. 
For this $i$, $x\in a B_i$ (see below (13.11)) and $g_i(x) = \pi_i(f(x))$
by (12.60), so $g(x) = g_i(x) = \pi_i(f(x))$ by (13.12), and hence (13.19)
will follow as soon as we prove that $\pi_i(f(x)) \in Q_i$.

Obviously $\pi_i(f(x))$ lies on the affine subspace $P$ spanned
by $Q_i$ (by definition of $\pi_i$ below (12.60)), so by (12.17)
it is enough to show that 
$$
\pi_i(f(x)) \in B(y_l,3(1+|f|_{lip}) t),
\leqno (13.20)
$$
where $l =l(i)$ is the index that we used in the definition of $Q_i$, 
above (12.42). But 
$$
|f(x) - y_l| \leq |f(x)-f(x_i)| + |f(x_i)-y_l|
\leq |f|_{lip} t + t
\leqno (13.21)
$$ 
because $x\in B_i$ and $f(x_i) \in D_l \,$. 
In addition, if $\pi'_i$ denotes the orthogonal projection onto 
the affine plane through $y_l$ parallel to $P$, 
then $||\pi'_i-\pi_i||_\infty \leq t$ because $P$ goes through 
$D_l = B(y_l,t)$ (by definition of ${\cal F}_l$; see below (12.17)). 
Then
$$
|\pi_i(f(x))-y_l| 
\leq |\pi'_i(f(x))-y_l| + t
\leq |f(x)-y_l| + t \leq (2+|f|_{lip}) t
\leqno (13.22)
$$
by (13.21), and now (13.20) and (13.19) follow.

Since by definition (12.17), $\H^d(Q) \leq C (1+|f|_{lip})^d t^d$ 
for all $Q\in {\cal F}$, (13.19) and (12.20) imply that
$$
\H^d\Big(g\Big(\bigcup_{j\in J_1} B_j \sm R\Big)\Big)
\leq \sum_{Q \in {\cal F}} \H^d(Q)
\leq C(\alpha,f) N^{-1} \H^d(X_0),
\leqno (13.23)
$$
where $C(\alpha,f)$ depends on $\alpha$ and $|f|_{lip}$.
This is still good, because $N$ will be chosen very large, 
depending on $f$, $\H^d(X_0)$, $\alpha$, and $\eta$ 
in particular.

\ms
Because of (13.23), we shall not need to worry too
much about what happens in $\bigcup_{j\in J_1} B_j \sm R$.
The set $R$ will not disturb much either, because $E \cap R$
is small. Indeed, we claim that
$$
\H^d\big(E \cap \overline R \big) 
= \H^d\big(E \cap \bigcup_{j\in J_1} [\overline B_j \sm aB_j]\big) 
\leq C (1-a) \H^d(X_0),
\leqno (13.24)
$$
where $C$ depends only on $M$ and $n$. First fix $j\in J_1$,
and observe that 
$$
\dist(x,P_{x_j}) \leq \varepsilon |x-x_j| 
\leq \varepsilon t 
\ \hbox{ for } x\in E \cap \overline B_j,
\leqno (13.25)
$$
because $|x-x_j| \leq t < \delta_3$ by (12.7) and (12.8), because
$x_j \in X_5$, and by (11.45). By elementary geometry, we can cover 
$P_{x_j} \cap [\overline B_j \sm aB_j]$ by less than 
$C (1-a)^{-d+1}$ balls of radius $(1-a) t$. 
Then the double balls cover $E \cap \overline B_j \sm aB_j$
(if $\varepsilon$ is small enough compared to $1-a$), and the 
local Ahlfors-regularity of $E$ (with the same justification
as for (12.12)) yields
$$
\H^d(E \cap \overline B_j \sm aB_j)
\leq C (1-a)^{-d+1} [(1-a) t]^d
= C (1-a) t^d.
\leqno (13.26)
$$
Next $\overline R \i \bigcup_{j\in J_1} \overline R_j
= \bigcup_{j\in J_1} [\overline B_j \sm aB_j]$ 
by (13.2) and (13.3). Also recall from (12.10) that 
$J_1$ has at most $C t^{-d} \H^d(X_0)$ elements; 
then (13.24) follows from (13.26).

\ms
$\dagger$
We now switch to the Lipschitz assumption. Set 
$$
V' = \bigcup_{j \in J_1} 2B_j \i U_{int},
\leqno (13.27)
$$
where $U_{int}$ is defined in (12.72) and the inclusion follows
from (12.76). We keep the same functions $\psi_j$ as above
(not to be confused with our bilipschitz mapping $\psi$),
and use the definition of $\wt f$ in (12.75) to set
$$
\wt g(x) = \wt f(x) + \sum_{j \in J_1} \psi_j(x) [\wt g_j(x)- \wt f(x)]
\ \hbox{ for } x\in V'
\leqno (13.28)
$$
(compare with (13.12); we still want to do the linear algebra on $B(0,1)$
and the partitions of unity on $U$). We intend to set
$$
g(x) = \lambda^{-1}\psi^{-1}(\wt g(x))
\ \hbox{ for } x\in V',
\leqno (13.29)
$$
so we need to check that $\wt g(x) \in B(0,1)$. Notice that
$$
\wt g(x) = \wt f(x) 
\ \hbox{ when } x \in V' \sm \bigcup_{j\in J_1} B(x_j,{1+a \over 2}\, t),
\leqno (13.30)
$$
because (12.77) says that $\wt g_j(x) = \wt f(x)$ for all $j$,
and by (13.28).
For such an $x$, $\wt g(x) = \wt f(x) = \psi(\lambda f(x))$
by (12.75), $\wt g(x)\in B(0,1)$ because $\psi$ maps to $B(0,1)$, 
and so (13.29) makes sense and we get that
$$
g(x) = f(x)
\ \hbox{ when } x \in V' \sm \bigcup_{j\in J_1} B(x_j,{1+a \over 2}\, t).
\leqno (13.31)
$$
Next suppose that $x$ lies in some $B(x_j,{1+a \over 2}\, t)$.
Obviously $\wt f(x)$ and $\wt g(x)$ are defined because $x\in V'$. 
Also observe that in fact, (12.77) says that $x \in B(x_j,{1+a \over 2}\, t)$ 
for all the indices $j$ such that $\wt g_j(x)-\wt f(x) \neq 0$, so
$$\eqalign{
|\wt g(x) - \wt f(x)|
&\leq \sum_{j \in J_1} \psi_j(x) |\wt g_j(x)- \wt f(x)|
\cr&
\leq 4\lambda\Lambda (1+ |f|_{lip}) \, t \sum_{j \in J_1} \psi_j(x)
\leq 4\lambda\Lambda (1+ |f|_{lip}) \, t
}\leqno (13.32)
$$
by (13.28), (12.81), (13.8) and (13.9). 
In addition,
$$\leqalignno{
\dist(\wt g(x),\R^n \sm B(0,1)) 
&\geq \dist(\wt y_l,\d B(0,1)) - |\wt g(x) - \wt y_l|
\cr&
\geq 10 (1+|f|_{lip})\, \lambda \Lambda t
- |\wt g(x) - \wt f(x)| - |\wt f(x) - \wt y_l|
\cr&
\geq 10 (1+|f|_{lip})\, \lambda \Lambda t
- 4(1+|f|_{lip})\, \lambda \Lambda t - 2(1+|f|_{lip})\, \lambda \Lambda t
&(13.33)
\cr&
\geq 4 (1+|f|_{lip})\, \lambda \Lambda t
}
$$
because $\wt y_l = \psi(\lambda y_l) \in B(0,1)$, by (12.83),
(13.32), and (12.80). In this case too, 
$\wt g(x) \in B(0,1)$ and we can define $g(x)$ as in (13.29).
This completes the legitimation of (13.29).

We decide to set directly
$$
g(x) = f(x)
\ \hbox{ when } x \in U \sm V';
\leqno (13.34)
$$
since $2B_j \i V'$ for all $j$, we see that
$\dist(U \sm V', B(x_j,{1+a \over 2}\, t)) > t$ and (13.31) 
gives us a large enough transition region where the two
definitions of $g$ give the same result. But in fact we shall
never use that definition outside of $V'$.

By (13.31) and (13.32), $||\wt g - \wt f||_{L^\infty(V')}
\leq 4\lambda\Lambda (1+ |f|_{lip}) \, t$, and then 
(by (13.29) and (13.34))
$$
||g - f||_\infty \leq 4 \Lambda^2 (1+ |f|_{lip}) \, t.
\leqno (13.35)
$$

Next we restrict to $E^{\varepsilon t}$ and check that
$$
|g(x) - f(x)| 
\leq 11 \Lambda^3 (1+|f|_{lip}) \varepsilon t + \alpha t
\ \hbox{ for } x\in E^{\varepsilon t}.
\leqno (13.36)
$$
By (13.31) and (13.34), we can assume that 
$x\in \bigcup_{j\in J_1} B(x_j,{1+a \over 2}\, t)) \i V'$. 
As for (13.32),
$$
|\wt g(x) - \wt f(x)| 
\leq \sum_{j\in J_1} \psi_j(x) |\wt g_j(x)-\wt f(x)|,
\leqno (13.37)
$$
and the only indices $j\in J_1$ that contribute are such that
$x \in E^{\varepsilon t} \cap B(x_j,{1+a \over 2}\, t)$ (use (13.28) 
and (12.77)). For these $j$, (12.88) applies and says that 
$|\wt g_j(x) - \wt f(x)| \leq 11 \lambda \Lambda^2 (1+|f|_{lip}) \varepsilon t 
+\lambda \Lambda^{-1} \alpha t$. We sum in $j$, use the
fact that $\sum_{j\in J_1} \psi_j(x) \leq 1$ by (13.8) and (13.9),
and get that $|\wt g(x) - \wt f(x)| \leq 
11 \lambda \Lambda^2 (1+|f|_{lip}) \varepsilon t 
+\lambda \Lambda^{-1} \alpha t$.
Now (13.36) follows from (13.29).

We also want to check that
$$
\wt g \ \hbox{ is } \lambda \, {2+3\Lambda^2 |f|_{lip} 
\over \Lambda}\hbox{-Lipschitz on } E^{\varepsilon t} \cap V'.
\leqno (13.38)
$$
We follow the proof of (13.15); given 
$x, y \in E^{\varepsilon t} \cap V'$, we set
$$\eqalign{
\wt \Delta(x,y) &= (\wt g-\wt f)(x)-(\wt g-\wt f)(y)
\cr&
= \sum_{j \in J_1} \Big\{\psi_j(x) [\wt g_j(x)-\wt f(x)]
-\psi_j(y) [\wt g_j(y)-\wt f(y)] \Big\}
=:  \sum_{j \in J_1} \Delta_j(x,y)
}\leqno (13.39)
$$
as in (13.16); since $\wt f$ is $\lambda\Lambda$-Lipschitz,
we just need to show that
$$
|\wt \Delta(x,y)| \leq \lambda \, {2+2\Lambda^2 |f|_{lip} 
\over \Lambda} \, |x-y|.
\leqno (13.40)
$$
When $|x-y| \geq t/10$, (13.40) holds because 
$|\wt \Delta(x,y)| \leq 22 \lambda\Lambda^4 (1+|f|_{lip}) \varepsilon t 
+ 2\lambda \Lambda \alpha t$ by (13.36) and (13.9) (and because we 
can choose $\varepsilon$ and $\alpha$ very small), 
so we may assume that $|x-y| \leq t/10$. By the same argument as above,
$x,y \in {3 \over 2} B_j$ for every $j\in J_1$ such that 
$\Delta_j(x,y) \neq 0$ in (13.39), and there are
at most $C$ indices $j$ for which this happens. For such $j$,
we proceed as in (13.17) and get that
$$\leqalignno{
|\Delta_j(x,y)|
& \leq \psi_j(x) |\wt g_j(x)-\wt f(x)-\wt g_j(y)+\wt f(y)|
+ |\psi_j(x)-\psi_j(y)| \, |\wt g_j(y)-\wt f(y)|
\cr& 
\leq \psi_j(x) |\wt g_j(x)-\wt f(x)-\wt g_j(y)+\wt f(y)|
+ 4 (1-a)^{-1} t^{-1} |x-y| \, |\wt g_j(y)-\wt f(y)|
\cr& 
\leq \lambda \, {1+2\Lambda^2 |f|_{lip} \over \Lambda} \, \psi_j(x) |x-y|
+ 4 (1-a)^{-1} |x-y| (11 \lambda\Lambda^2 (1+|f|_{lip}) \varepsilon 
+ \lambda \Lambda^{-1} \alpha)
&(13.41)
}
$$
because $\psi_j$ is still $4[(1-a)t]^{-1}$-Lipschitz,
$\wt f$ is $\lambda \Lambda$-Lipschitz,
$\wt g_j$ is $\lambda \, {1+\Lambda^2 |f|_{lip} \over \Lambda}$-Lipschitz
on $E^{\varepsilon t} \cap 2B_j$ (by (12.96)), and by (12.88)
(if $y \in B(x_j,{1+a \over 2}t)$; otherwise $\wt g_j(y)=\wt f(y)$
directly by (12.77)).

When we sum this over $j$, the first term gives a total contribution
which is bounded by $\lambda \, {1+2\Lambda^2 |f|_{lip} \over 
\Lambda} \, |x-y|$, and the second one contributes at most 
$\lambda \Lambda^{-1} |x-y|$, if $\varepsilon$ and $\alpha$ are 
small enough (depending on $a$); (13.40) and (13.38) follow.

We deduce from (13.38) and (13.29) that 
$$
g \ \hbox{ is } (2+3\Lambda^2 |f|_{lip})\hbox{-Lipschitz on } 
E^{\varepsilon t} \cap V'.  
\leqno (13.42)
$$

Return to what we did in the rigid case.
We still have (13.18) in the present Lipschitz case (see (13.31) and 
(13.34)). Let us check now that (13.19) also holds now, 
i.e., that
$$
g\Big(\bigcup_{j\in J_1} B_j \sm R\Big) \i \bigcup_{l\in {\cal L}}
\bigcup_{Q \in {\cal F}_l} Q = \bigcup_{Q \in {\cal F}} Q.
\leqno (13.43)
$$
As before, any $x\in \bigcup_{j\in J_1} B_j \sm R$
lies in some $aB_j$ (the first one), for which
$\wt g(x) = \wt g_j(x) = \wt \pi_j(\wt f(x))$
(by (13.11) and (13.28)), and it is enough to check
that $g(x) \in Q_j$, or equivalently that
$\wt \pi_j(\wt f(x)) \in \wt Q_j$, since 
$Q_j = \lambda^{-1}\psi^{-1}(\wt Q_j)$ by (12.22) 
and $g(x) = \lambda^{-1}\psi^{-1}(\wt g(x))$ by (13.29).

Recall from the definition above (12.77) that $\wt \pi_j$ is the 
orthogonal projection onto the plane $\wt P_j$ that was defined
below (12.49), subject to the constraint (12.23) for some
$y_l$ such that $f(x_j) \in D_l$ (see above (12.48)). That is,
$$
\wt P_j \hbox{ meets } B(\wt y_l, 2 \lambda \Lambda (1+|f|_{lip}) t)).
\leqno (13.44) 
$$
Since $\wt Q_j = \wt P_j \cap B(\wt y_l, 20\lambda \Lambda^2 
(1+|f|_{lip}) t)$ by (12.22), it is enough to check that
$$
|\wt \pi_j(\wt f(x))-\wt y_l| < 20\lambda \Lambda^2 
(1+|f|_{lip}) t.
\leqno (13.45) 
$$
But 
$$\eqalign{
|\wt f(x)-\wt y_l| &= |\psi(\lambda f(x))-\psi(\lambda y_l)|
\leq \lambda \Lambda |f(x)-y_l|
\cr&
\leq \lambda \Lambda (|f(x)-f(x_j)|+|f(x_i)-y_l|)
\cr&
\leq \lambda \Lambda ( a |f|_{lip} t +  t)
\leq \lambda \Lambda ( 1 + |f|_{lip})\, t
}\leqno (13.46) 
$$
by (12.75) and the definition of $\wt y_l = \psi(\lambda y_l)$ above (12.22), 
and because $x\in aB_i$. 

Let $\pi$ denote the orthogonal projection onto the plane
parallel to $\wt P_j$, but through $\wt y_l$; then
$||\pi-\wt\pi_j||_{\infty} \leq 2 \lambda \Lambda (1+|f|_{lip}) t)$
by (13.44), and
$$\eqalign{
|\wt \pi_j(\wt f(x))-\wt y_l| 
&\leq ||\pi-\wt\pi_j||_{\infty} + |\pi(\wt f(x))-\wt y_l| 
\leq ||\pi-\wt\pi_j||_{\infty} + |\wt f(x)-\wt y_l|
\cr&
\leq 2 \lambda \Lambda (1+|f|_{lip}) t) 
+ \lambda \Lambda ( 1 + |f|_{lip})\, t
}\leqno (13.47) 
$$
by (13.46). This is better than (13.45), and (13.43) follows.

Clearly (13.23) still holds, even though with a larger constant 
$C(\alpha,f)$, now by (13.43), (12.22), and (as before) (12.20).

The last estimates (13.24)-(13.26) stay the same; they do not
even involve $\psi$.
$\dagger$

\ms\noindent
{\bf 14. Step 3. Places where $f$ has a very contracting 
direction, and the $B_j$, $j\in J_2$.} 
\ms
At the beginning of Section 12, we were left with a set
$X_5 \i X_0$, such that $\H^d(X_0 \sm X_5) \leq 4 \eta$
by (11.48). Set (as in (13.1))
$$
V = \bigcup_{j\in J_1} B_j = \bigcup_{j\in J_1} B(x_j,t),
\leqno (14.1)
$$
which contains $X_5 \cap X_N(\delta_4)$ by (12.9).
In principle, we already took good care of $V$
in Sections 12 and 13, by (13.23) and (13.26). We also know 
from (12.3) and (12.4) that
$$
\H^d\big([X_5 \cap f^{-1}(Y_N)] \sm X_N(\delta_4))\big)
= \H^d\big([X_5 \cap f^{-1}(Y_N)] \sm f^{-1}(Y_N(\delta_4))\big) \leq \eta,
\leqno (14.2)
$$
where $Y_N$ is as in (12.1). Next consider
$$
X_6 = X_5 \sm \Big[ f^{-1}(Y_N) \cup V \Big].
\leqno (14.3)
$$
If $x\in X_5 \sm [V \cup X_6]$, then it lies in $f^{-1}(Y_N)$
($V$ is not allowed) and, since it does not lie in 
$X_5 \cap X_N(\delta_4)$ (which is contained in $V$ too), it
lies in the set of (14.2). So
$$
\H^d(X_0 \sm [V \cup X_6])
\leq \H^d(X_0 \sm X_5) + \H^d(X_5 \sm [V \cup X_6])
\leq 5 \eta,
\leqno (14.4)
$$
and we may now turn to $X_6$.

\ms
Our next target is the set of points $x\in X_6$ where $A_x$ has
a very contracting direction along $P_x$. That is, we want to control
the set
$$
X_7 = \big\{ x\in X_6 \, ; \, \hbox{ there is a unit vector
$\nu \in P_x'$ such that } |DA_x(\nu)| \leq \gamma \big\},
\leqno (14.5)
$$
where $P_x'$ denotes the vector space parallel to $P_x$,
$DA_x$ is the differential of $A_x$, and $\gamma <1$ is another
very small positive constant, to be chosen later. 

The following is very similar to Lemma 4.60 in [D2], 
whose fairly standard proof applies here too (so we skip it).

\ms\proclaim Lemma 14.6.
We can find a finite collection of balls $B_j = B(x_j,r_j)$, 
$j\in J_2$, with the following properties:
$$
x_j \in X_7 \,\hbox{ and } \, 0 < r_j \leq \delta_6
\hbox{ for } j\in J_2,
\leqno (14.7)
$$
where $\delta_6$ is as in (12.7), 
$$
\hbox{ the $\overline B_j$, $j\in J_2$ are disjoint, and do not meet }
\bigcup_{j\in J_1} B\big(x_j,{(1 +a) t \over 2}\big),
\leqno (14.8)
$$
and 
$$
\H^d\big(X_7 \sm \bigcup_{j\in J_2}\overline B_j \big)
\leq \eta.
\leqno (14.9)
$$

\ms
This time, since the $\overline B_j$ are disjoint, we shall not need a
subtle partition of unity as before, and we can define functions $g_j$
independently. Also, what we intend to do here in the $B_j$, 
$j\in J_2$, will be independent of what we did in the 
$B\big(x_j,{(1 +a) t \over 2}\big)$.
Again we start with the rigid assumption.

We set $P_j = P_{x_j}$, $Q_j = A_{x_j}(P_j)$, 
denote by $\pi_j$ the orthogonal projection on $Q_j$,
and define $g_j$ by the same formulae (12.59)-(12.61) as before
(with $t$ replaced by $r_j$).

Notice that $g_j(x) \in [f(x),\pi_j(f(x))]$; then 
$$
|g_j(x)-f(x)| \leq |\pi_j(f(x))-f(x)| \leq |f|_{lip} |x-x_j|
\leqno (14.10)
$$
because $\pi \circ f - f = (\pi -I) \circ f$ is $|f|_{lip}$-Lipschitz and
vanishes at $x_j$ (recall that $Q_j$ goes through $f(x_j)$ because 
the definition (11.37) says that $A_{x_j}(x_j) = f(x_j)$).
When $x\in B_j$, we get that $|g_j(x)-f(x)| \leq |f|_{lip} \, r_j$.
When $x\in U \sm B_j$, the analogue of (12.59) says that $g_j(x)=f(x)$. 
Altogether,
$$
||g_j - f||_\infty \leq |f|_{lip} \, r_j \leq |f|_{lip} \, \delta_6
\leqno (14.11)
$$
by (14.7) and as in (12.62). Also, the the proof of (12.66)
(which could also be simplified here) says that
$$
g_j \ \hbox{ is $(1+|f|_{lip})$-Lipschitz on } E^{\varepsilon r_j}.
\leqno (14.12)
$$

For $k$ large enough ($k$ is the index in our initial sequence
of quasiminimal sets $E_k$, which converges to $E$), 
$$
\dist(z,P_j \cap a B_j) \leq 2 \varepsilon r_j
\ \hbox{ for } z \in E_k \cap aB_j,
\leqno (14.13)
$$
by (11.45). Set ${\cal E} = A_{x_j}(P_j \cap a B_j) \i Q_j$, 
let $z \in E_k \cap aB_j$ be given, and let $w\in P_j \cap a B_j$ 
be such that $|z-w| \leq 3 \varepsilon r_j$; then
$$\eqalign{
\dist(f(z),{\cal E}) 
&\leq |f(z)-f(w)| + |f(w)-A_{x_j}(w)| + \dist(A_{x_j}(w),{\cal E})
\cr&
= |f(z)-f(w)| + |f(w)-A_{x_j}(w)|
\cr&
\leq 3 \varepsilon r_j |f|_{lip} + \varepsilon |w-x_j|
\leq (1 + 3 |f|_{lip}) \varepsilon r_j 
}\leqno (14.14)
$$
because $A_{x_j}(w) \in A_{x_j}(P_j \cap a B_j) = {\cal E}$,
by (11.46), and again because $w\in P_j \cap a B_j$.
By (12.60), $g_j(z) = \pi_j(f(z))$, so it 
lies in $Q_j$, and (by (14.14)) in a 
$(1 + 3 |f|_{lip}) \varepsilon r_j$-neighborhood ${\cal E}'$
of ${\cal E}$ in $Q_j$.
Now $x_j \in X_7$, so ${\cal E}$ is an ellipsoid, with a shortest axis
of length at most $2 \gamma a r_j$ (by (14.5)), and other axes of 
length at most $2 |f|_{lip} a r_j$. Then 
$$\eqalign{
\H^d(g_j(E_k \cap aB_j)) 
&\leq \H^d({\cal E}')
\leq C (1+|f|_{lip})^{d-1} (\gamma + (1 + 3 |f|_{lip}) \varepsilon) r_j^d
\cr&
\leq C (1+|f|_{lip})^{d-1} \gamma r_j^d
}\leqno (14.15)
$$
again for $k$ large enough and if $\varepsilon$ is small enough,
depending on $\gamma$.

Set $R_j = B_j \sm a B_j$, as before. Then 
$$
\H^d(E \cap \overline R_j)
= \H^d(E \cap \overline B_j \sm aB_j)
\leq C (1-a) r_j^d
\leqno (14.16)
$$
by the same proof as for (13.26). Since
$$
r_j < \delta_6 < {1 \over 10} \delta_1 = {1 \over 10}\dist(X_1,\R^n \sm W_f)) 
\leq \dist(x_j,\R^n \sm W_f))
\leqno (14.17)
$$
by (14.7), (12.7), and (11.22)
(and as in (12.11)), we also get that
$$
r_j^d \leq C \H^d(E \cap B_j)
\ \hbox{ for } j\in J_2,
\leqno (14.18)
$$
by the local Ahlfors-regularity of $E$, and where 
the use of Proposition 4.1 is justified as for (12.12).
Then
$$
\sum_{j\in J_2} r_j^d \leq C \sum_{j\in J_2} \H^d(E \cap B_j)
\leq C\H^d\big(E \cap \bigcup_{j\in J_2} B_j\big)
\leq C \H^d(E \cap W_f)
\leqno (14.19)
$$
by (14.19), (14.8), and (14.17), and now (14.16) implies that
$$\eqalign{
\H^d(E \cap \bigcup_{j\in J_2} \overline R_j)
&\leq \sum_{j\in J_2} \H^d(E \cap \overline B_j \sm aB_j)
\leq C (1-a) \sum_{j\in J_2} r_j^d
\cr&\leq C (1-a) \H^d(E \cap W_f) = C (1-a) \H^d(X_0) 
}\leqno (14.20)
$$
(because $E \cap W_f = X_0$ by (11.20)).

So we should not worry too much about the $R_j$, and
since we have some control on the $aB_j$ by (14.15) and (14.19), 
we shall now concentrate on 
$$
X_8 = X_6 \sm \big[ X_7 \cup \bigcup_{j\in J_2} B_j \big].
\leqno (14.21)
$$

Let $X_9$ be a compact subset of $X_8$ such that $\H^d(X_8 \sm X_9) 
\leq \eta$. Since
$$
X_6 \i X_8 \cup \Big( \bigcup_{j\in J_2} B_j \Big) \cup X_7
\i X_8 \cup \Big( \bigcup_{j\in J_2} \overline B_j \Big) 
\cup \Big( X_7 \sm \bigcup_{j\in J_2} \overline B_j \Big)
\leqno (14.22)
$$
we get that
$$
\H^d\big(X_6 \sm \big[X_9 \cup \big(\bigcup_{j\in J_2} \overline B_j 
\big)\big]\big) \leq \H^d(X_8 \sm X_9) 
+ \H^d\big(X_7 \sm \bigcup_{j\in J_2} \overline B_j\big)
\leq 2 \eta
\leqno (14.23)
$$
by (14.9). Let us deduce from this and (14.4) that
$$
\H^d\big(X_0 \sm \big[X_9 \cup \big(\bigcup_{j\in J_1 \cup J_2} \overline B_j 
\big)\big]\big) \leq 7 \eta.
\leqno (14.24)
$$
Let $Z, Z', Z''$ the sets in the left-hand sides of 
(14.24), (14.23) and (14.4) respectively; we want to check that
$Z \i Z' \cup Z''$. Let $x\in Z\sm Z''$ be given. Then $x\in X_0$,
and so $x\in V \cup X_6$. But $x\in V =  \bigcup_{j\in J_1} B_j$
is impossible because $x\in Z$ (also see the definition (14.1)),
hence $x\in X_6$. Then $x\in Z'$, as needed. So (14.24) holds.

\ms$\dagger$
Under the Lipschitz assumption, we need to modify the definition 
of $g_j$. We still set $P_j = P_{x_j}$, but we consider
$\wt Q_j = \wt A_{x_j}(P_j)$ (where $\wt A_{x_j}$ is the affine 
approximation of $\wt f$, as in (12.37)-(12.39)). 
We denote by $\wt \pi_j$ the orthogonal projection onto
$\wt Q_j$, and will define $\wt g_j$ as we did near (12.77).
Again we first work in the set $U_{int}$ defined by (12.72), 
because this is where we extended $\wt f$ (see (12.75)). 
Notice that $U_{int}$ contains all the $2B_j$, $j\in J_3$, 
by proof of (12.76) (just use (14.7) instead of (12.7)). 

We define $\wt g_j$ on $U_{int}$ with the same formulas
(12.77)-(12.79) as before, with the choice of $\wt \pi_j$ that we just
made, and $t$ replaced with $r_j$. 

We continue, as in Section 12, with estimates for $x\in 2B_j \i U_{int}$.
First observe that 
$$
|\wt g_j(x)-\wt f(x)| \leq |\wt \pi_j(\wt f(x))-\wt f(x)|
\leq \lambda \Lambda |f|_{lip} |x-x_j|
\leqno (14.25)
$$
because $\wt g_j(x) \in [\wt f(x),\wt \pi_j(\wt f(x))]$
by (12.77)-(12.79), and because 
$(\wt \pi_j - I) \circ \wt f$ is $\lambda \Lambda |f|_{lip}$-Lipschitz
and vanishes at $x_j$ by definition of $\wt A_{x_j}$ and $\wt \pi_j$. 
Next 
$$\eqalign{
|\wt g_j(x)-\wt f(x_j)| 
&\leq |\wt g_j(x)-\wt f(x)| + |\wt f(x)-\wt f(x_j)|
\cr& \leq 2\lambda \Lambda |f|_{lip} |x-x_j|
\leq 2\lambda \Lambda |f|_{lip} \, r_j
\leq 2\lambda \Lambda |f|_{lip} \, \delta_6
}\leqno (14.26)
$$
by (14.25), because $\wt f$ is $\lambda \Lambda |f|_{lip}$-Lipschitz,
and by (14.7). But 
$$
\dist(f(x_j),\R^n \sm U) \geq \dist(\wh W, \R^n \sm U) = \delta_0
\geq 10\Lambda^2(1+|f|_{lip}) \delta_6
\leqno (14.27)
$$
because $x_j \in E_0 = E \cap W_f$ (see (11.20) and (11.19)), 
so $f(x_j) \in \wh W$ (see (2.1) and (2.2)), 
and by (12.6) and (12.7). Hence
$$\eqalign{
\dist(\wt g_j(x),\R^n \sm B(0,1))
&\geq \dist(\wt f(x_j),\R^n \sm B(0,1)) - 2\lambda \Lambda |f|_{lip} 
\, \delta_6
\cr&
\geq \lambda\Lambda^{-1} \dist(f(x_j),\R^n \sm U)
- 2\lambda \Lambda |f|_{lip} \, \delta_6
\cr& 
\geq 10\lambda \Lambda(1+|f|_{lip}) \delta_6
- 2\lambda \Lambda |f|_{lip} \, \delta_6
\geq 8\lambda \Lambda(1+|f|_{lip}) \delta_6
}\leqno (14.28)
$$
by (14.26), (12.75), and (14.27). Thus $\wt g_j(x) \in B(0,1)$
when $x\in 2B_j$.

When $x\in U_{int} \sm 2B_j$, and even when 
$x\in U_{int} \sm B(x_j, {1+a \over 2} r_j)$, (12.77) and (12.75)
say that $\wt g_j(x) = \wt f(x) = \psi(\lambda f(x))$; then of course
$\wt g_j(x) \in \psi(\lambda U) = B(0,1)$. Thus in both cases
$\wt g_j(x) \in B(0,1)$, and we can define $g_j$ on $U_{int}$ by
$$
g_j(x) = \lambda^{-1}\psi^{-1}(\wt g_j(x))
\ \hbox{ for } x\in U_{int},
\leqno (14.29)
$$
(compare with (12.85)). This formula yields 
$g_j(x) = f(x)$ for $x\in U_{int} \sm B(x_j, {1+a \over 2} r_j)$, 
and we may even extend it by deciding that
$$
g_j(x) = f(x)
\ \hbox{ for } x\in U \sm B(x_j, {1+a \over 2} r_j)
\leqno (14.30)
$$
(now compare with (12.86)). But in fact the values of $g_j$
on $U \sm U_{int}$, or even $U \sm 2B_j$, will never matter.

Return to the modifications concerning this section.
The analogue of (14.11) is now
$$
||g_j - f||_\infty 
\leq \Lambda^2 |f|_{lip} r_i \leq \Lambda^2 |f|_{lip} \delta_6,
\leqno (14.31)
$$
which follows from (14.30), (14.25), and (14.29).
Then we worry about the Lipschitz estimate (14.12). The fact that
$$
\hbox{$\wt g_j$ is }
\lambda \, {1+\Lambda^2 |f|_{lip} \over \Lambda} 
\hbox{-Lipschitz on } E^{\varepsilon t} \cap 2B_j
\leqno (14.32)
$$
is proved as (12.96) or (12.66) (with some simplifications), 
and implies that
$$
\hbox{$g_j$ is } (1+\Lambda^2 |f|_{lip}) \hbox{-Lipschitz on } 
E^{\varepsilon t} \cap 2 B_j;
\leqno (14.33)
$$
will be good enough to take replace (14.12).

Observe that (14.13) still holds with the same proof.
Next we generalize (14.14) and (14.15). Set 
$\wt {\cal E} = \wt A_{x_j}(P_j \cap a B_j) \i \wt Q_j$
(recall that we set $\wt Q_j = \wt A_{x_j}(P_j)$).
Let $z \in E_k \cap aB_j$ be given, and use (14.13) to find 
$w\in P_j \cap a B_j$ such that $|z-w| \leq 3 \varepsilon r_j$; then
$$\eqalign{
\dist(\wt f(z),\wt{\cal E}) 
&\leq |\wt f(z)-\wt f(w)| + |\wt f(w)-\wt A_{x_j}(w)| 
+ \dist(\wt A_{x_j}(w),\wt{\cal E})
\cr&
= |\wt f(z)-\wt f(w)| + |\wt f(w)-\wt A_{x_j}(w)|
\cr&
\leq 3 \varepsilon r_j |\wt f|_{lip} + \lambda \varepsilon |w-x_j|
\leq (1 + 3 \Lambda |f|_{lip}) \lambda\varepsilon r_j 
}\leqno (14.34)
$$
because $\wt A_{x_j}(w) \in \wt A_{x_j}(P_j \cap a B_j) = \wt {\cal E}$,
by (12.52), by (12.75), and because $w\in P_j \cap a B_j$.

By the analogue of (12.78), $\wt g_j(z) = \wt\pi_j(\wt f(z))$, so it 
lies in $\wt Q_j$ (by definition of $\wt \pi_j$), and (by (14.34)) in a 
$(1 + 3 \Lambda |f|_{lip}) \lambda \varepsilon r_j$-neighborhood 
$\wt {\cal E}'$ of $\wt{\cal E}$ in $\wt Q_j$. So we just checked that
$\wt g_j(E_k \cap aB_j) \i \wt {\cal E}'$ for $k$ large.

Now $\wt {\cal E}$ is an ellipsoid in $\wt Q_j$, its
axes all have lengths smaller than $2\lambda \Lambda |f|_{lip}r_j$
by (12.37), and one of them is much shorter. Indeed,
let $\nu$ be a unit vector in $P_{x_j}'$ such that
$|DA_{x_j}(\nu)| \leq \gamma$; such a vector exists because
$x_j \in X_7$ (see (14.5) and (14.7));
notice that both $f$ and $\wt f$ are differentiable at $x_j$ 
in the direction of $\nu$, and recall that 
$\wt f(x) = \psi(\lambda f(x))$ near $x_j$. Then 
$$\eqalign{
|D\wt A_{x_j}(\nu)|
&= \big|\lim_{t \to 0} t^{-1}
[\wt A_{x_j}(x_j + t\nu)-\wt A_{x_j}(x_j)]\big|
= \big|\lim_{t \to 0} t^{-1} [\wt f(x_j + t\nu)- \wt f(x_j)]\big|
\cr&
\leq \lambda \Lambda \big|\limsup_{t \to 0} 
t^{-1} [f(x_j + t\nu)-f(x_j)]\big|
\cr&
= \lambda \Lambda \big|\limsup_{t \to 0} 
t^{-1} [A_x(x_j + t\nu)-A_x(x_j)]\big|
= \lambda \Lambda |DA_{x_j}(\nu)|
\leq \lambda \Lambda \gamma
}\leqno (14.35)
$$
by (12.39), because $\wt A_{x_j}(x_j) = \wt f(x_j)$, and by (11.40); 
hence the smallest axis of $\wt {\cal E}$ has length at most 
$2\lambda \Lambda\gamma r_j$.
Thus 
$$
\H^d(\wt g_j(E_k \cap aB_j)) 
\leq \H^d(\wt {\cal E}')
\leq C \lambda^d \Lambda^d (1+ \Lambda |f|_{lip})^{d} \gamma r_j^d
\leqno (14.36)
$$
for $k$ large, and hence also
$$
\H^d(g_j(E_k \cap aB_j)) 
\leq \Lambda^d \lambda^{-d} \H^d(\wt g_j(E_k \cap aB_j)) 
\leq C \Lambda^{2d} (1+ \Lambda |f|_{lip})^{d} \gamma r_j^d
\leqno (14.37)
$$
by (14.29), and as in (14.15).

Finally the estimates (14.16)-(14.24) go through with only
minor modifications.
$\dagger$

\ms
\msi{\bf 15. Step 4. The remaining main part of $X_0$}
\ms

Return to the rigid assumption.
We care about $X_9$ now (see near (14.22)). Set $Y_9 = f(X_9)$ and, for
$y \in Y_9$, 
$$
Z(y) = X_9 \cap f^{-1}(y) = \big\{ x\in X_9 \, ; \, f(x)=y \big\}.
\leqno (15.1)
$$
Notice that $Z(y)$ has at most $N$ points, because 
$X_9 \i X_6 \i X_5 \sm f^{-1}(Y_N)$
(see (14.21), the definition of $X_9$ just below (14.21),  (14.3), and (12.1)).
We claim that for each $y\in Y_9$, there is a 
positive radius $r(y)$ such that 
$$
X_9 \cap f^{-1}(B(y,r)) \i \bigcup_{x\in Z(y)} B(x,2\gamma^{-1} r)
\ \hbox{ for } 0 < r \leq r(y).
\leqno (15.2)
$$
Here $\gamma$ is the same as in the definition (14.5) of $X_7$.
The proof uses a small compactness argument (to make sure that
it is enough to control $f$ near the $x\in Z(y)$,
(11.46) (to show that $A_x$ controls $f$ near $x\in Z(y)$), 
and the fact that we excluded $X_7$ in (14.21)
(to exclude points near $x\in Z(y)$ that don't lie
in $B(x,2\gamma^{-1} r)$).
We don't repeat it here because it is the same as in 
Lemma 4.69 in [D2].  

Since $X_9$ is compact and disjoint from
the finite collection of $B_j$, $j\in J_1 \cup J_2$
(see (14.21), (14.3), and (14.1)), the number
$$ 
\delta_7 = \dist\big(X_9, 
\bigcup_{j\in J_1 \cup J_2} {1+a \over 2} \, B_j \big)
\leqno (15.3)
$$
is positive. Just by making $r(y)$ smaller if needed,
we may assume that for each $y\in Y_9$, 
$$
0 < r(y) < {\gamma \over 4} \, \Min\big( \delta_6 , \delta_7, 
\Min\{|x-x'| \, ; \, x, x' \in Z(y), x \neq x'\} \big).
\leqno (15.4)
$$
Then
$$
B(x,2\gamma^{-1}r(y)) \cap {1+a \over 2} \, B_j
= \emptyset
\ \hbox{ for $y\in Y_9$, $x\in Z(y)$, and } j\in J_1 \cup J_2,
\leqno (15.5)
$$
and (for each $y\in Y_9$)
$$
\hbox{the balls $B(x,2\gamma^{-1} r(y))$, $x \in Z(y)$, are disjoint.}
\leqno (15.6)
$$
We'll need some uniformity (i.e., to know that $r(y)$ is not 
too small), so let us choose a new small constant 
$\delta_8 \in (0,\delta_7] > 0$ 
such that if we set 
$$
Y_{10} = \big\{ y\in Y_9 \, ; \, r(y) > \delta_8 \big\}
\ \hbox{ and } \ 
X_{10} = X_9 \cap f^{-1}(Y_{10}), 
\leqno (15.7)
$$
then 
$$
H^d(X_9 \sm X_{10}) \leq \eta.
\leqno (15.8)
$$
As usual, such a $\delta_8$ exists, because
the monotone union of the sets $X_{10}$,
when $\delta_8$ tends to $0$, is $X_9$.
Next set
$$
Y_{11} = \big\{ y\in Y_{10} \, ; \, 
\hbox{ all the affine planes $Q_x = A_x(P_x)$,
$x\in Z(y)$, coincide} \big\}.
\leqno (15.9)
$$
Notice that the $Q_x$ are $d$-planes, because we excluded the case 
when $A_x$ has a very contracting direction in $P_x$. Also set
$$
X_{11} = X_{10} \cap f^{-1}(Y_{11}).
\leqno (15.10)
$$
As we check in (4.77) of [D2],  
$$
\H^d(X_{10} \sm X_{11}) =0.
\leqno (15.11)
$$
The same proof is valid here; we sketch it to prevent the reader
from worrying. 
We use the fact that $f(E)$ is rectifiable, and prove 
that it does not have any approximate tangent plane at points of
$Y_{10} \sm Y_{11}$ (too many tangent directions exist), so
$\H^d(Y_{10} \sm Y_{11}) = 0$. For this last,
we use again the fact that we excluded contracting directions.
For the accounting, we also use the fact that $f$ is at most
$N$-to-$1$ on $X_{10}$, to return from $Y_{10}$ to $X_{10}$
and prove (15.11). Incidentally, (4.77) of in [D2]  
is wrongly referred to as (4.78) at the end of the proof
in [D2] (sorry!). 

\ms
We want to cover $X_{11}$, but it will be more efficient to cover 
$Y_{11}$ first. We choose a finite collection of balls
$D_j = B(y_j,r_j)$, $j\in J_3$, so that 
$$
y_j \in Y_{11} \ \hbox{ and } \ 
0 < r_j < \delta_8 
\ \hbox{ for } j\in J_3, 
\leqno (15.12)
$$
$$
\hbox{ the $\overline D_j$, $j\in J_3$, are disjoint}
\leqno (15.13)
$$
and
$$
\H^d\big(X_{11} \sm f^{-1}\big(\bigcup_{j\in J_3} \overline D_j \big)\big)
\leq \eta.
\leqno (15.14)
$$
See Lemma 4.79 in [D2], where 
(15.14) is deduced form a similar estimate on $Y_9$,
using again the fact that we excluded contracting directions
and $f$ is at most $N$-to-$1$ on $X_{11}$.

Observe that for $j\in J_3$,
$$
r_j < \delta_8 \leq r(y_j) \leq {\gamma \over 4} \, \Min\big( 
\delta_6 , \delta_7, 
\Min\{|x-x'| \, ; \, x, x' \in Z(y_j), x \neq x'\} \big)
\leqno (15.15)
$$
by (15.12), because $y_j \in Y_{10}$, and by (15.7) and (15.4).

We want to modify $f$ on the sets $f^{-1}(D_j)$, as we did
in the balls $B_j$, $j\in J_1 \cup J_2$. We shall be able to 
proceed independently on each $f^{-1}(D_j)$, because the $D_j$ 
are disjoint by (15.13). In fact, for each $j$ we shall only modify $f$
on the $f^{-1}(D_j) \cap B(x,2\gamma^{-1}r_j)$, $x\in Z(y_j)$,
which are disjoint by (15.15) and contain the interesting part 
of $f^{-1}(D_j)$ by (15.2).

Fix $j\in J_3$. By definition, the $d$-planes $A_x(P_x)$, $x\in Z(y_j)$,
are all equal; let us call $Q_j$ this common $d$-plane that we get.
For each $x\in Z(y_j)$, set
$$
E(x) = P_x \cap A_x^{-1}(Q_j \cap D_j).
\leqno (15.16)
$$
This is a $d$-dimensional ellipsoid in $P_j$,
whose axes have lengths between $2 |f|_{lip}^{-1} r_j$ and
$2 \gamma^{-1} r_j$, by (11.36) and the definition (14.5)
of $X_7$ (which we excluded in (14.21)). The analogues of 
$a B_j$, ${1+a \over 2} B_j$, and $B_j$ in the previous sections
will be
$$
B_{j,x}^- = \big\{ z\in \R^n \, ; \, 
\dist(z,a E(x)) < 20^{-1}(1+|f|_{lip})^{-1} (1-a) r_j \big\}
\leqno (15.17)
$$ 
and
$$
B_{j,x}^+ = \big\{ z\in \R^n \, ; \, 
\dist(z,a E(x)) \leq 10^{-1}(1+|f|_{lip})^{-1} (1-a) r_j \big\},
\leqno (15.18)
$$
and 
$$
B_{j,x} = \big\{ z\in \R^n \, ; \,
\dist(z, (2-a)E(x)) \leq (1-a) r_j \big\}.
\leqno (15.19)
$$
Observe that 
$$
B_{j,x}^- \i B_{j,x}^+ \i B_{j,x} \i B(x,{3 \over 2}\gamma^{-1}r_j) 
\leqno (15.20)
$$
if $(1-a)$ is small enough, depending on $|f|_{lip}$ and $\gamma$.
Recall from Remark 11.17 
that $a$ is allowed to depend on $\gamma$.
By (15.20) and (15.15), 
$$
\hbox{for each $j\in J_3$, the $B_{j,x}$, $x\in Z(y_j)$, are disjoint.}
\leqno (15.21)
$$
Let us also check that
$$
f(B_{j,x}^+) \i {1+a \over 2} D_j.
\leqno (15.22)
$$
Let $z\in B_{j,x}^+$ be given, and let $w$ be a point of
$aE(x)$ such that 
$$
|z-w| \leq 10^{-1}(1+|f|_{lip})^{-1} (1-a) r_j.
\leqno (15.23)
$$
Observe that 
$$\eqalign{
|w-x| &\leq |w-z|+|z-x| 
\leq 10^{-1}(1+|f|_{lip})^{-1} (1-a) r_j + {3 \over 2}\gamma^{-1}r_j
\cr&
\leq 2\gamma^{-1}r_j
\leq \delta_6 < \delta_3/10
}\leqno (15.24)
$$
by (15.23), (15.20), if $a$ is small enough, and by
(15.15) and (12.7). In addition, $w\in P_x$, so
(11.46) applies and says that
$$
|f(w)-A_x(w)| \leq \varepsilon |w-x|
\leq 2\varepsilon\gamma^{-1}r_j.
\leqno (15.25)
$$
Now
$$\eqalign{
|f(z) - y_j| &\leq |f(z)-f(w)| + |f(w)-A_x(w)| + |A_x(w)-y_j|
\cr&
\leq |z-w| |f|_{lip}  + 2 \varepsilon\gamma^{-1}r_j + |A_x(w)-y_j|
\cr&
\leq 10^{-1} (1-a) r_j + 2 \varepsilon\gamma^{-1}r_j + a r_j
< {1+a \over 2} \, r_{j}
}\leqno (15.26)
$$
by (15.25), (15.23), because $A_x(w) \in a D_j$
by definition of $E(x)$ (see (15.16)), and if $\varepsilon$
is small enough; (15.22) follows.

Because of (15.22), (15.21), and (15.13),
$$
\hbox{ the $B_{j,x}^+$, $j\in J_3$ and $x\in Z(y_j)$, are all disjoint,}
\leqno (15.27)
$$
even for different values of $j$, because the $D_j$ are disjoint.
Set
$$
R_{j,x} = B_{j,x}^+ \sm B_{j,x}^-
\leqno (15.28)
$$ 
by analogy with the previous constructions. Denote by $\pi_j$
the orthogonal projection onto $Q_j$, and define
$g_{j,x}$ as follows: set
$$
g_{j,x}(z) = \pi_j(f(z))
\ \hbox{ when } z\in B_{j,x}^-,
\leqno (15.29)
$$
$$
g_{j,x}(z) = f(z)
\ \hbox{ when } z\in \R^n \sm B_{j,x}^+,
\leqno (15.30)
$$
and interpolate in the usual linear way in the 
remaining intermediate region $R_{j,x}$. That is, set
$$
g_{j,x}(z) = (1-\beta(z))\pi_j(f(z)) + \beta(z) f(z)
\ \hbox{ for } z\in R_{j,x},
\leqno (15.31)
$$
with 
$$
\beta(z) = {20(1+|f|_{lip}) \dist(z,a E(x)) \over (1-a) r_j} \, - 1.
\leqno (15.32)
$$
We also define a function $g_j$: we set
$g_j = g_{j,x}$ on each $B_{j,x}^+$, $x\in Z(y_j)$,
and $g_j(z) = f(z)$ on the rest of $\R^n$; the definition
is coherent, by (15.30) and (15.21), and we even get a lipschitz 
mapping (possibly with very bad constants). Let us check that
$$
||g_j - f||_\infty \leq r_j \leq \delta_6.
\leqno (15.33)
$$
The second inequality comes from (15.15). 
By (15.31), it is enough to check that 
$|\pi_j(f(z))-f(z)| \leq r_j$
for $z\in B_{j,x}^+$, and this is clear because 
$f(z) \in D_j$ by (15.22), and $D_j$ is centered on $Q_j$.

When $z\in B_{j,x}^+ \cap E^{\varepsilon r_j} = 
\big\{ z\in \R^n \, ; \, \dist(z,E) \leq \varepsilon r_j \big\}$, 
the estimate improves: we can choose $w\in E$
such that $|w-z| \leq \varepsilon r_j$, and, by (11.45),
$p\in P_x$ such that 
$|p-w| \leq \varepsilon |w-x| \leq 2\varepsilon r_j$
(with the same sort of justification as (15.24) for (15.25)).
Thus $|p-z| \leq 3\varepsilon r_j$ and
$|p-x| \leq |p-z|+|z-x| \leq 3 \varepsilon r_j + {3 \over 2}\gamma^{-1} 
r_j \leq 2 \gamma^{-1} r_j$, by (15.20) and if $\varepsilon$ is small enough.
Then
$$\eqalign{
|g_j(z) - f(z)| &\leq |\pi_j(f(z))-f(z)|
\leq \dist(f(z),Q_j)
\cr&
\leq \dist(A_x(p),Q_j) +  |A_x(p)-f(p)| + |f(p)-f(z)|
\cr&
= |A_x(p)-f(p)| + |f(p)-f(z)|
\cr&
\leq \varepsilon |p-x| + |p-z| |f|_{lip} 
\leq (2\gamma^{-1} +3|f|_{lip}) \varepsilon r_j
}\leqno (15.34)
$$
because $A_x(p) \in Q_j$ by definition of $Q_j$ and by (11.46).
Then
$$
g_j \ \hbox{ is $C(1+|f|_{lip})$-Lipschitz on } 
B_{j,x}^+ \cap E^{\varepsilon r_j}.
\leqno (15.35)
$$ 
by the same proof as for (12.66); here again, the small
$\varepsilon$ wins against the large $\gamma^{-1}$, 
$|f|_{lip}$, and $(1-a)^{-1}$.

We shall also need to know that
$$
X_9 \cap f^{-1}\big(\bigcup_{j\in J_3} \overline D_j \big)
\i \bigcup_{j\in J_3}\bigcup_{x\in Z(y_j)} B_{j,x}.
\leqno (15.36)
$$
Indeed let $z\in X_9 \cap f^{-1}\big(\bigcup_{j\in J_3} \overline D_j \big)$ 
be given, and let $j\in J_3$ be such that $f(z) \in \overline D_j$.
Notice that $|f(z)-y_j| \leq r_j < \delta_8 < r(y_j)$
by (15.12), because $y_j \in Y_{11} \i Y_{10}$ (see (15.12) and (15.9)), 
and by (15.7), so (15.2) says that $z \in \overline B(x,2\gamma^{-1}r_j)$ for some
$x\in Z(y_j)$. Now $|z-x| \leq 2\gamma^{-1}r_j \leq \delta_3/10$
by the last part of (15.24), so $\dist(z,P_x) \leq \varepsilon |z-x|$ 
by (11.45).
Let $w\in P_x$ be such that $|z-w| \leq \varepsilon |z-x|$; then
$$\eqalign{
|A_x(w)-y_j| &\leq |A_x(w)-f(w)| + |f(w)-f(z)| + |f(z)-y_j|
\cr&
\leq \varepsilon |w-x| + |f|_{lip} |z-w| + r_j
< r_j + 3 \varepsilon \gamma^{-1} (1+|f|_{lip}) \, r_j
}\leqno (15.37)
$$
by (11.46) and because $|w-x| \leq 3 \gamma^{-1}r_j  < \delta_3/5$. 
Recall that $D_j$ is centered at $y_j = f(x) = A_x(x) \in A_x(P_x) = 
Q_j$, so (15.37) says that $A_x(w) \in Q_j \cap 
(1+3 \varepsilon \gamma^{-1} (1+|f|_{lip})) D_j$.
Then
$$
w \in (1+3 \varepsilon \gamma^{-1} (1+|f|_{lip})) E(x) \i (2-a) E(x)
\leqno (15.38)
$$
by (15.16), because $a<1$, and if $\varepsilon$ is small enough.
Now
$$
\dist(z,(2-a) E(x)) \leq |z-w| \leq \varepsilon |z-x|
\leq 2 \varepsilon \gamma^{-1} r_j < (1-a) r_j
\leqno (15.39)
$$
if $\varepsilon$ is small enough, and hence $z\in B_{x,j}$
(see (15.19)). This proves (15.36).

\ms$\dagger$
When we work under the Lipschitz assumption, 
we need a few modifications to the definitions above.
Surprisingly, we do not modify anything before (15.28).
One could argue that it would be more natural to cover 
$\wt Y_{11} = \psi(\lambda Y_{11})$ instead of $Y_{11}$,
but we prefer to keep the same definitions, 
and we will be able to handle the differences.
In particular, we shall prove that 
$$
\hbox{for each $y\in Y_{11}$, all the affine planes 
$\wt Q_x = \wt A_x(P_x)$, $x\in Z(y)$, coincide.}
\leqno (15.40)
$$
But let us first check that if $y\in Y_{10}$,
$x\in Z(y)$, and $Q_x = A_x(P_x)$, then
the restriction of $\psi$ to $\lambda Q_x$ is differentiable 
at $\lambda y$, with a derivative $D_{\psi}$ such that
$$
\lambda D_{\psi}(DA_x(v)) = D\wt A_x(v)
\ \hbox{ for } v \in P'_x,
\leqno (15.41)
$$
where $P'_x$ denotes the vector space parallel to $P_x$
(and as we would expect from the chain rule). And indeed,
$$\eqalign{
D\wt A_x(v) &= \lim_{t \to 0} t^{-1} [\wt A_x(x+tv)-\wt A_x(x)]
= \lim_{t \to 0} t^{-1} [\wt f(x+tv)-\wt f(x)]
\cr&
= \lim_{t \to 0} t^{-1} [\psi(\lambda f(x+tv))-\psi(\lambda f(x))],
}\leqno (15.42)
$$
by (12.39), and where the last line comes from the convention 
that we used for (12.38) and (12.39), that $\wt f$ is defined by 
the formula (12.36) near $x$. But
$$\eqalign{
f(x+tv) - y & = 
f(x+tv) - f(x) = A_x(x+tv) - f(x) + o(t) 
\cr&
= A_x(x+tv) - A_x(x) + o(t) 
= t DA_x(v) + o(t) 
}\leqno (15.43)
$$
because $x\in Z(y)$, by (11.40), and because $A_x(x) = f(x)$ 
by (11.40) and $A_x$ is affine, so
$$
\psi(\lambda f(x+tv)) 
= \psi[\lambda (y+t DA_x(v) + o(t))]
= \psi[\lambda y + \lambda t DA_x(v)] + o(t)
\leqno (15.44)
$$
because $\psi$ is Lipschitz. So (15.42) says that
$$\eqalign{
D\wt A_x(v) &= \lim_{t \to 0} t^{-1} 
\big(\psi[\lambda y + \lambda t DA_x(v)]-\psi(\lambda f(x))\big)
\cr&
= \lim_{t \to 0} t^{-1} 
\big(\psi[\lambda y + \lambda t DA_x(v)]-\psi(\lambda y) \big).
}\leqno (15.45)
$$
Now let $w$ be any vector in the vector space
$Q_x'$ parallel to $Q_x$, write $w = DA_x(v)$
for some $v\in P'_x$, and observe that (15.45)
says that $\psi$ is differentiable at $\lambda y$ in the
direction $\lambda w$, with a derivative equal to $D\wt A_x(v)$
(and hence that satisfies (15.41)).
We could easily get the differentiability (instead of
the differentiability in each direction), because $\psi$
is Lipschitz, but let us not even bother,
because we just need the formula (15.41) for the 
directional derivatives. Notice however that
since $DA_x : P'_x \to Q'_x$ is a bijection
(because $DA_x$ has no contracting direction because $x\in X_9$; 
see (14.5), (14.21) and the line below it, and (15.1)), 
(15.41) allows us to compute $D_{\psi}$ from $DA_x$ and $D\wt A_x$.

We are now ready to prove (15.40). Let $y\in Y_{11}$
be given; by (15.9), all the affine planes $A_x(P_x)$, $x\in Z(y)$,
are equal to some affine space $Q_y$; in addition, we just checked
that $\psi$ has directional derivatives at $ty$ along $Q_y$,
given by a mapping $D_{\psi}$ that we can compute from the 
values of $DA_x$ and $D\wt A_x$ at some $x\in Z(y)$.
Now (15.41) says that for each $x\in Z(y)$, the vector space
$\wt Q_x'$ parallel to $\wt A_x(P_x)$ is given by
$$
\wt Q_x' = D\wt A_x(P_x') = D_{\psi}(DA_x(P_x')) = D_{\psi}(Q_y'),
\leqno (15.46)
$$
where $Q_y'$ is the vector plane parallel to $Q_y$.
In particular, $\wt Q_x'$ does not depend on $x\in Z(y)$. Since all
the $\wt A_x(P_x)$ go through $\wt f(x) = \psi(\lambda f(x)) = 
\psi(\lambda y)$ by construction, they are all equal, and 
(15.40) follows. 

\ms
Return to the definition of the $g_{j,x}$ near (15.29). 
As before, we first define auxiliary functions $\wt g_{j,x}$
on the set $U_{int}$ defined by (12.72), and on which we
extended $\wt f$ in (12.75). Notice that for $j\in J_3$
and $x\in Z(y_j)$,
$$
B_{j,x}^+ \i B_{j,x} \i  B(x,{3 \over 2}\gamma^{-1}r_j) 
\i B(x,3\gamma^{-1}r_j) \i U_{int}
\leqno (15.47)
$$
by (15.20), because $x \in X_9 \i X_0$ (by (15.1)), since
$r_j \leq {\gamma \over 4} \delta_6 \leq {\gamma \delta_0\over 40 
(1+|f|_{lip})}$ by (15.15) and (12.7), and by the definition 
(12.72).

For $j\in J_3$, we denote by $\wt Q_j$ the common value of the
affine planes $\wt A_x(P_x)$, $x\in Z(y_j)$, and by
$\wt \pi_j$ the orthogonal projection onto $\wt \pi_j$. 
Then we set
$$
\wt g_{j,x}(z) = \wt \pi_j(\wt f(z))
\ \hbox{ when } z\in B_{j,x}^-,
\leqno (15.48)
$$
$$
\wt g_{j,x}(z) = \wt f(z)
\ \hbox{ when } z\in U_{int} \sm B_{j,x}^+,
\leqno (15.49)
$$
and 
$$
\wt g_{j,x}(z) = (1-\beta(z)) \wt\pi_j(\wt f(z)) + \beta(z) \wt f(z)
\ \hbox{ for } z\in R_{j,x},
\leqno (15.50)
$$
with $\beta(z)$ as in (15.32). Also set
$\wt g_j = \wt g_{j,x}$ on each $B_{j,x}^+$, $x\in Z(y_j)$,
and $\wt g_j(z) = \wt f(z)$ on the rest of $E \cup U_{int}$; the definition
is still coherent, for the same reasons as before, and $\wt g_j$
is Lipschitz. The analogue of (15.33) is 
$$
||\wt g_j - \wt f||_{L^\infty(U_{int})} \leq \lambda \Lambda r_j 
\leq \lambda \Lambda \delta_6,
\leqno (15.51)
$$
which we prove as before: the second inequality follows from
(15.15), and for the first one it is enough to observe that
for $z\in B_{j,x}^+$,
$$\eqalign{
|\wt g_{j,x}(z)-\wt f(z)| 
&\leq |\wt \pi_j(\wt f(z))-\wt f(z)|
= \dist(\wt f(z),\wt Q_j)
\leq |\wt f(z)-\wt f(x)|
\cr&
\leq \lambda \Lambda |f(z)-f(x)|
= \lambda \Lambda |f(z)-y| \leq \lambda \Lambda r_j
}\leqno (15.52)
$$
because $\wt Q_j$ goes through $\wt f(x)$ and by (15.22).

Next we want to define $g_j$. We want to set
$$
g_j(z) = \lambda^{-1} \psi^{-1}(\wt g_j(z))
\ \hbox{ for } z \in U_{int}
\leqno (15.53)
$$
so let us check that $\wt g_j(z) \in B(0,1)$.
When $z\in B_{j,x}^+$ for some $j\in J_3$ and $x\in Z(y_j)$,
$$\eqalign{
\dist(\wt g_{j,x}(z), \R^n \sm B(0,1))
&\geq \dist(\wt f(z), \R^n \sm B(0,1)) - \lambda \Lambda r_j
\cr&
\geq \lambda \Lambda^{-1} \dist(f(z), \R^n \sm U) 
- \lambda \Lambda r_j
}\leqno (15.54)
$$
by (15.52) and because $\wt f(z) = \psi(\lambda f(z))$ and 
$\psi : \lambda U \to B(0,1)$ is bilipschitz; then
$$\eqalign{
\dist(f(z), \R^n \sm U) 
&\geq \dist(f(x), \R^n \sm U) - |f(z)-f(x)|
\cr&
= \dist(f(x), \R^n \sm U) - |f(z)-y_j| 
\cr&
\geq \dist(f(x), \R^n \sm U) - r_j
}\leqno (15.55)
$$
because $f(z) \in D_j$ by (15.22), and 
$$
\dist(f(x), \R^n \sm U) \geq \dist(\wh W, \R^n \sm U) = \delta_0
\geq 10\Lambda^2(1+|f|_{lip})\delta_6
\leqno (15.56)
$$
because $x \in X_0$, hence $f(x) \in \wh W$ by (11.20), (2.1), and 
(2.2), and by (12.6) and (12.7). Altogether,
$$\eqalign{
\dist(\wt g_{j,x}(z), \R^n \sm B(0,1))
& \geq \lambda \Lambda^{-1} \dist(f(z), \R^n \sm U) 
- \lambda \Lambda r_j
\cr&
\geq \lambda \Lambda^{-1} [10\Lambda^2(1+|f|_{lip})\delta_6
- r_j ] - \lambda \Lambda r_j
\cr&\geq 8\lambda \Lambda (1+|f|_{lip})\delta_6
}\leqno (15.57)
$$
by (15.54), (15.55), (15.56), and because $r_j \leq \delta_6$ by
(15.15). So $\wt g_{j,x}(z) \in B(0,1)$ and $g_j(z)$ is correctly 
defined in (15.53).

When $z\in U_{int} \sm B_{j,x}^+$ for all $j\in J_3$ and $x\in Z(y_j)$,
we defined $\wt g_j(z) = \wt f(z)$ below (15.50), and $\wt f(z) = 
\psi(\lambda f(z))$ by (12.75), so (15.53) makes sense, and even 
yields $g_j(z) = f(z)$. We can thus extend the definition of $g_i$,
and set $g_j(z) = f(z)$ for $z \in U \sm U_{int}$, but in fact 
we won't even need that. Anyway, we get that
$$
g_j(z) = f(z) \ \hbox{ for }
z \in U_{int} \sm \bigcup_{j\in J_3}\ \bigcup_{x\in Z(y_j)} B_{j,x}^+.
\leqno (15.58)
$$
Notice that
$$
||g_j - f||_{L^\infty(U_{int})} 
\leq \lambda^{-1} \Lambda ||\wt g_j - \wt f||_\infty 
\leq \Lambda^2 r_j 
\leq \Lambda^2 \delta_6,
\leqno (15.59)
$$
by (15.58), (15.53) and (15.51). Let us also check that
$$
|g_j(z) - f(z)| 
\leq \Lambda^2 (2\gamma^{-1}+3 |f|_{lip})\, \varepsilon r_j
\ \hbox{ for } z\in B_{j,x}^+ \cap E^{\varepsilon r_j}.
\leqno (15.60)
$$
We prove this as in (15.34). We can again choose $w\in E$ such that 
$|w-z| \leq \varepsilon r_j$, and (by (11.45)) $p\in P_x$ such that 
$|p-w| \leq \varepsilon |w-x| \leq 2\varepsilon r_j$; 
thus $|p-z| \leq 3\varepsilon r_j$ and $|p-x| \leq |p-z|+|z-x| 
\leq 3\varepsilon r_j+{3 \over 2} \gamma^{-1} r_j 
\leq 2\gamma^{-1} r_j$ by (15.20), and 
$$\eqalign{
|\wt g_j(z) - \wt f(z)| &\leq |\wt \pi_j(\wt f(z))-\wt f(z)|
\leq \dist(\wt f(z), \wt Q_j)
\cr&
\leq \dist(\wt A_x(p), \wt Q_j) + |\wt A_x(p)-\wt f(p)| 
+ |\wt f(p)-\wt f(z)|
\cr&
= |\wt A_x(p)-\wt f(p)| + |\wt f(p)- \wt f(z)|
\cr&
\leq \lambda \varepsilon |p-x| + |p-z| |\wt f|_{lip}
\leq (2\gamma^{-1}+3\Lambda |f|_{lip}) \lambda \varepsilon r_j
}\leqno (15.61)
$$
because $\wt A_x(p) \in \wt Q_j = \wt A_x(P_x)$, and
by (12.52) (with the same justification as for (15.25));
(15.60) follows.

We claim that now 
$$
g_j \ \hbox{ is $C\Lambda^2 (1+|f|_{lip})$-Lipschitz on } 
B_{j,x}^+ \cap E^{\varepsilon r_j},
\leqno (15.62)
$$
with the same proof as for (12.96). Finally, (15.36) still holds 
in the Lipschitz context; its proof only involves $f$ and 
arguments anterior to (15.29) and the definition
of the $g_{j,x}$, so we can keep it.
$\dagger$

\ms\noindent
{\bf 16. The modified function $g$, and a deformation for $E$.}
\ms

We are now ready to define a (new) function $U: \R^n \to \R^n$,
which is a first competitor for the replacement of $f$.
We already defined a function $g$ in Step 2.f, by (13.12)
or (13.29) and (13.31), and we intend to keep it like this on 
$$
V_1 = \bigcup_{j\in J_1} {1 + a \over 2} B_j
= \bigcup_{j\in J_1} B(x_j,{(1 + a)t \over 2}).
\leqno (16.1)
$$
That is, we set
$$
g(z) = f(z) + \sum_{j \in J_1} \psi_j(z) [g_j(z)-f(z)]
\ \hbox{ for } z \in V_1
\leqno (16.2)
$$
in the rigid case, and 
$$
g(z) = \lambda^{-1}\psi^{-1}(\wt g(z))
\ \hbox{ for } z \in V_1,
\leqno (16.3)
$$
with 
$$
\wt g(z) = \wt f(z) + \sum_{j \in J_1} \psi_j(z) [\wt g_j(z)- \wt f(z)]
\leqno (16.4)
$$
under the Lipschitz assumption. We also set
$$
g(z) = g_j(z)
\ \hbox{ for } z \in {1 + a \over 2} B_j = B(x_j,{(1 + a) r_j \over 2})
\leqno (16.5)
$$
when $j \in J_2$, and 
$$
g(z) = g_j(z) = g_{j,x}(z) 
\ \hbox{ for } z \in B_{j,x}^+
\leqno (16.6)
$$
when $j \in J_3$ and $x\in Z(y_j)$. Finally, set
$$
V_1^+ = \Big[\bigcup_{j\in J_1 \cup J_2} {1 + a \over 2} B_j \Big]
\cup \Big[\bigcup_{j\in J_3 \, ; \, x\in Z(y_j)} B_{j,x}^+ \Big];
\leqno (16.7)
$$
we just defined $g$ on $V_1^+$, and we keep
$$
g(z) = f(z) 
\ \hbox{ for } z\in \R^n \sm V_1^+.
\leqno (16.8)
$$

$\dagger$ 
Under the Lipschitz assumption, we
also have a function $\wt g$, defined on $V_1^+$,
and such that $\wt g(z) = \psi(\lambda g(z))$.
On $V_1$, we wrote this explicitly in (16.3)
and (16.4); on the balls ${1 + a \over 2} B_j$,
$j\in J_2$, this comes from the fact that $g_i$
was defined by (14.29) (also recall that $2B_j \i U_{int}$
for $j\in J_2$); on the $B_{j,x}^+$, this comes from
(15.47) and (15.53) (also see the line below (15.50)). 
$\dagger$

Let us check that all these definitions are independent
because the corresponding sets are disjoint. First,
the $B_j$, $j\in J_2$, are disjoint from each other 
and from $\bigcup_{j\in J_1} {1 + a \over 2} B_j$, by (14.8).
The $B_{j,x}^+$ are disjoint from each other by (15.27).
Finally, if $j\in J_3$ and $x\in Z(y_j)$, 
$$
B_{j,x}^+ \i B(x,{3 \over 2}\gamma^{-1}r_j)
\i B(x,{\delta_7\over 2})
\leqno (16.9)
$$
(15.20) and (15.15). This last ball does not meet any
${1 + a \over 2} B_j$, $j\in J_1 \cup J_2$, by the definition
(15.3) of $\delta_7$ and because $x \in X_9$ (by (15.1)).

Next we check that
$$
\hbox{$g$ is Lipschitz on $U$}
\leqno (16.10)
$$
(but possibly with a very bad norm).
Recall that (16.2) or (16.3) would also yield
$g(z) = f(z)$ for $z\in \d V_1$, by (13.18) or (13.31)
(recall that our initial $g$ was Lipschitz).
Similarly, $g_j(z) = f(x)$ for $j\in J_2$ and 
$z \in \d B(x_j,{(1 + a) r_j \over 2})$; in the rigid case,
this is because we still use (12.59) (see above (14.10)), 
and in the Lipschitz case this comes from (12.77) and (14.29), 
or directly from (14.30).
Finally, $g_{j,x}(z) = f(z)$ for $j \in J_3$, $x\in Z(y_j)$,
and $z \in \d B_{j,x}^+$, by (15.30) or (15.49) and (15.53). 
Thus (16.8) does not introduce any discontinuity, 
and (16.10) follows easily, because $g$ is
Lipschitz on the closure of each piece.

Let us finally record that
$$
||g-f||_\infty \leq 4 \Lambda^2 (1+ |f|_{lip}) \delta_6
\leqno (16.11)
$$
in the rigid case by (13.13) and (12.8), (14.11), and (15.33), 
and in the Lipschitz case by (13.35) and (12.8), (14.31), and (15.59).

\ms
We would like to use $g$ to define new competitors,
and a natural first step is to check that $g$ is the
endpoint of a one-parameter family of functions $g_t$, 
that satisfies the conditions (1.4)-(1.8), and in particular 
the boundary conditions (1.7), relative to $E$.

This will not be entirely satisfactory, because 
we would like (1.7) to hold with respect to the $E_k$,
but we shall take care about that in the next section. 

Recall that $f$ itself is defined as $f(x) = \varphi_1(x)$,
for some one-parameter family of functions $\varphi_t$, which we 
extended from $E$ to $\R^n$ at the beginning of Section 11,
and for which (1.4)-(1.8) hold by assumption. 
We start under the rigid assumption and set
$$
g_t(x) = \varphi_{2t}(x)
\ \hbox{ for } 0 \leq t \leq 1/2
\leqno (16.12)
$$
and
$$
g_t(x) = (2-2t) f(x) + (2t-1) g(x)
\ \hbox{ for } 1/2 \leq t \leq 1.
\leqno (16.13)
$$

Recall that (1.4)-(1.8) for the $\varphi_t$ holds
with respect to the ball $B = B(X_0,R_0)$ of (11.1); 
here we shall find it convenient to use a slightly 
larger ball $B'$.

\ms\proclaim Lemma 16.14.
The functions $g_t$, $0 \leq t \leq 1$, 
satisfy (1.4)-(1.8), relative to $E$ and the ball 
$B' = \overline B(X_0,R_0 + 4\Lambda^2 (1+ |f|_{lip}) \delta_6)$.

\ms
We shall need to know that 
$$
\dist(x,X_1) \leq \delta_6 \ \hbox{ and } \ 
\dist(x,\R^n \sm W_f) > \delta_1/2
\ \hbox{ for } x\in V_1^+,
\leqno (16.15)
$$
where we defined $V_1^+$ in (16.7) and 
$\delta_1 = \dist(X_1,\R^n \sm W_f)$ in (11.22).
The second part follows from the first part, because
$\delta_6 < \delta_1/2$ by (12.7). For the first part
there are three similar cases.
When $x\in B_j$ for some $j\in J_1$, this is true because
$x_j \in X_N(\delta_4) \i X_1$ and 
$|x-x_j| \leq t < \delta_6$; 
see the line below (12.8), the various definitions of the $X_j$, 
and (12.8). 
When $x\in B_j$ for some $j\in J_2$, we use 
(14.7) instead. When $x \in B_{j,z}^+$ for some $j\in J_3$
and $z\in Z(y_j)$, we use the fact that $z\in X_9 \i X_1$
by (15.1), and $B_{j,x}^+ \i B(x,{3 \over 2}\gamma^{-1}r_j)
\i B(x,\delta_6)$ by (15.20) and (15.15). So (16.15) holds.

The properties (1.4) and (1.8) hold by construction.
For (1.5), since we know that $g_0(x) = \varphi_0(x) = x$
for $x\in \R^n$, it is enough to check that 
$$
g_t(x)=x \ \hbox{ for $x\in E \sm B$ and } 0 \leq t \leq 1.
\leqno (16.16)
$$
Let $x\in E \sm B$ be given. By (1.5) for the $\varphi_t$,
$\varphi_t(x) = x$ for $0 \leq t \leq 1$, hence by (16.12)
$g_t(x)=x$ for $0 \leq t \leq 1/2$. If $x\in V_1^+$,
(16.15) says that $x\in W_f$ and, since $x\in E$,
this forces $x \in B$ (because $\varphi_1(x) \neq x$ by the
definition (11.19), and by (1.5)); this is impossible.
So $x\in \R^n \sm V_1^+$, and $g_t(x)=f(x)=x$ by (16.13) and (16.8);
this proves (16.16) and (1.5).

For (1.6), we need to check that
$$
g_t(x) \in B' \ \hbox{ when $x\in E \cap B'$ and $0 \leq t \leq 1$.}
\leqno (16.17)
$$
This is trivial when $x\in E \cap B' \sm B$, because $g_t(x) = x \in B$.
If $x\in B$ and $0 \leq t \leq 1/2$, $g_t(x) = \varphi_{2t}(x) \in B$
by (16.12) and (1.6) for the $\varphi_t$. Finally, if
$x\in B$ and $t > 1/2$, $g_t(x)$ lies on the segment
$[f(x),g(x)]$, which is contained in $B'$ because
$f(x) \in B$ and $|g(x)-f(x)| \leq 4\Lambda^2 (1+ |f|_{lip}) \delta_6$
by (16.11).

We still need to check (1.7), i.e., that for $0 \leq k \leq j_{max}$,
$$
g_t(x) \in L_k
\ \hbox{ when $x \in E \cap L_k \cap B'$ and $0 \leq t \leq 1$.}
\leqno (16.18)
$$
[We just used the letter $k$ to avoid a conflict with the notation
for the $B_j$, but of course $k$ is not the index for the sequence
$\{ E_k \}$ here.]
We may assume that $x\in B$, because otherwise $g_t(x) = x \in L_k$,
and that $t \geq 1/2$, because otherwise 
$g_t(x) = \varphi_{2t}(x) \in L_k$ by (1.7) for the $\varphi_t$.
By (16.13), $g_t(x)$ lies on the segment $[f(x),g(x)]$,
so we just need to check that 
$$
[f(x),g(x)] \i L_k
\hbox{ for } x \in E \cap L_k \cap B.
\leqno (16.19)
$$
Since this is trivial when $g(x) = f(x)$, we may 
assume that $x \in V_1^+$. 

First suppose that $x\in V_1$,
and let $j\in J_1$ be such that $x\in B_j$.

Return to the definition of $Q_j$ (Step 2.e, starting above (12.42)). 
Still denote by $x_j$ the the center of $B_j$; we chose 
$l \in {\cal L}$ such that $f(x_j) \in D_l$, and observed that 
we can find $x(l) \in X_5$ such that $y_l = f(x(l))$. 

But $X_5 \i X_2$, so by (11.26) there is an $m \in [0,n]$ such that
$x(l) \in X_{1,\delta_2}(m)$. That is, by (11.23)-(11.24)
$y_l = f(x(l)) \in {\cal S}_m \sm {\cal S}_{m-1}$, 
and (if $m \geq 1$)
$$
\dist(y_l,{\cal S}_{m-1}) \geq \delta_2.
\leqno (16.20)
$$
Still denote by $F_l$ the smallest face
of our grid that contains $y_l$, and by $W(y_l)$ 
the affine plane spanned by $F_l$; obviously 
$F_l$ and $W(y_l)$ are $m$-dimensional.
Also notice that
$$
|f(x)-y_l| \leq |f(x)-f(x_j)| + |f(x_j)-y_l| 
\leq t |f|_{lip} + t
\leq \delta_6 (1+|f|_{lip}) 
\leqno (16.21)
$$
because $f(x_j) \in D_l$ and by (12.8). 
Let us check that
$$
\hbox{any face $F$ of our grid that contains $f(x)$
contains $F_l$ too.}
\leqno (16.22)
$$
We use coordinates and the dyadic structure to prove this, 
but probably polyhedra would work as well.
Also recall that we work under the rigid assumption for the moment.
For $1 \leq i \leq n$, denote by $a_i$ and $b_i$ 
the $i$-th coordinate of $f(x)$ and $y_l$ respectively. 
Thus 
$$
|b_i-a_i| \leq {1 \over 10} \, \min(\delta_2, r_0)
\leqno (16.23)
$$
by (16.21) and (12.7).
Set $I_0 = \big\{ i \in [1,n] \, ; \,  
b_i \notin r_0 {\Bbb Z} \big\}$ (recall that $r_0$
is the scale of our dyadic grid).
For $i\in I_0$, (16.20) says that 
$\dist(b_i,r_0 {\Bbb Z}) > \delta_2$,
so $[a_i,b_i]$ does not meet $r_0 {\Bbb Z}$
(by (16.23)). 

Denote by $w$ the point obtained from $f(x)$
by replacing each $a_i$, $i\in I_0$, with
$b_i$. We get that $w \in F$ too.
And if we want to go from $w$ to $z$, we just 
need to replace each coordinate $a_i$, $i \notin I_0$,
with $b_i$, which by (16.23) and the definition of $I_0$
is the closest point of $r_0 {\Bbb Z}$. 
Then $y_l$ lies in any face that may contain $w$, 
including $F$. Altogether, $y_l \in F$, and since $F_l$ 
is the smallest face that contains $y_l$, we get that $F_l \i F$,
as needed for (16.22).

\ms
Recall that $g_j(x) \in [f(x), \pi_j(f(x))]$
(by (12.59)-(12.61)), where $\pi_j$ is the orthogonal projection 
on the affine plane $\wh Q_j$ spanned by $Q_j$, that $Q_j$
lies in ${\cal F}_l$ (see above (12.43)), and hence goes through $D_l$ 
and is contained in $W(y_l)$ (see above (12.18)). 
Let $\pi$ denote the orthogonal projection
onto the affine plane through $y_l$ parallel to $Q_j$; then 
$$\eqalign{
|\pi_j(f(x))-y_l| &\leq |\pi(f(x))-y_l|+||\pi-\pi_j||_{\infty}
\leq |f(x)-y_l|+ t \leq t |f|_{lip} + 2t
\cr&
\leq 2\delta_6 (1+|f|_{lip}) \leq \delta_2/5
< \dist(y_l,\d F_l)
}\leqno (16.24)
$$
by various parts of (16.21), (12.7), and (16.20).
Also, $\pi_j(f(x)) \in \wh Q_j \i W(y_l)$, 
the affine space spanned by $F_l$; then (16.24)
implies that $\pi_j(f(x)) \in F_l$ because the segment
$[\pi_j(f(x)), y_l] \i W(y_l)$ does not meet $\d F_l$.
Thus
$$
\pi_j(f(x)) \in F_l \i F
\hbox{ for any face $F$ of our grid that contains $f(x)$,}
\leqno (16.25)
$$
by (16.22) for the second part.

By (1.7) for the $\varphi_t$, $f(x) = \varphi_1(x) \in L_k$.
Let $F$ be a face of $L_k$ that contains $f(x)$.
The proof of (16.25) shows that $\pi_i(f(x)) \in F$
for each $i\in J_1$ such that $x\in B_i$ 
(that is, not only for $i=j$), and then
$g_i(x) \in [f(x), \pi_i(f(x))]$
lies in $F$ too (because every face is convex).

By (16.2), (13.8), and (13.9), $g(x)$ lies in the convex hull 
of $f(x)$ and the $g_i(x)$, where $i\in J_1$ 
is such that $\psi_i(x) \neq 0$.
For such $i$, (13.5) and (13.6) imply that $x\in B_i$,
so $g_i(x) \in F$. Altogether, $g(x) \in F$
and $[f(x),g(x)] \i F \i L_k$, as needed for (16.19).

\ms
Our second case for the proof of (16.19) is when
$x\in {1 + a \over 2} B_j$ for some $j\in J_2$.
Set $y_j = f(x_j)$, denote by $F(y_j)$ the smallest
face that contains $y_j$, by $W(y_j)$ the affine 
subspace spanned by $F(y_j)$, and by $m$ their dimension. 
This time $x_j \in X_7 \i X_2$ by (14.7) and various definitions, 
so $y_j = f(x_j) \in {\cal S}_m \sm {\cal S}_{m-1}$
and
$$
\dist(y_l,{\cal S}_{m-1}) \geq \delta_2
\leqno (16.26)
$$
if $m \geq 1$, by the proof of (16.20). Since
$$
|f(x)-y_j| = |f(x)-f(x_j)| \leq r_j |f|_{lip} \leq \delta_6 |f|_{lip}
\leq \delta_2/10
\leqno (16.27)
$$
by (14.7), and (12.7), the same proof as for (16.22)
shows that any face $F$ of our grid that contains $f(x)$
contains $F(y_j)$ too.

Here $g(x) = g_j(x) \in [f(x), \pi_j(f(x))]$
by (16.5) and because $g_j(x)$ is given by (12.59)-(12.61),
and where $\pi_j$ now denotes the orthogonal projection on
$Q_j = A_{x_j}(P_j)$ (see below (14.9)). But Lemma 12.27 says
that $Q_j = A_{x_j}(P_j) \i W(f(y_j))$, and since
$|\pi_j(f(x))-y_j| \leq |f(x)-y_j| \leq \delta_2/10$
because $Q_j$ goes through $y_j$ and by (16.27),
the proof of (16.25) shows that $\pi_j(f(x)) \in F(y_j)$.

As before, the $\varphi_t(x)$ and $f(x) = \varphi_1(x)$ lie in $L_k$.
Let $F$ be a face of $L_k$ that contains $f(x)$; then 
$\pi_j(f(x)) \in F(y_j) \i F$, and 
$g(x) \in [f(x), \pi_j(f(x))]$ lies in $F$
too (by convexity). So 
$[f(x),g(x)] \i F \i L_k$, and (16.19) holds in this case too.

Our last case is when $x$ lies in $B_{j,z}^+$ for some
$j\in J_3$ and $z\in Z(y_j)$ (recall that $x\in V_1^+$
and see the definition (16.7)). We proceed as in the second case,
notice that $x_j \in X_9$ by (15.1), replace (16.27) with 
the fact that $f(x) \in D_j = B(y_j,r_j) \i B(y_j,\delta_6)$
by (15.22) and (15.15) (see the definition of $D_j$ above (15.12)).
Then $g(x) = g_{j,z}(x) \in [f(x), \pi_j(f(x))]$
by (16.6) and (15.29)-(15.32), and where $\pi_j$ 
denotes the orthogonal projection on
$Q_j = A_z(P_z)$ (see above (15.29) and (15.16)), which is again
contained in $W(y_j)$ by Lemma~12.27. The rest of the argument is 
the same. This completes our proof of (16.19) and, by the same
token, of (16.18); this was our last verification;
Lemma 16.14 follows.
\qed

\ms$\dagger$
Under the Lipschitz assumption, we keep $g_t(z) = \varphi_{2t}(z)$
for $0 \leq t \leq 1/2$, as in (16.12), but for $t \geq 1/2$,
we want to preserve the faces when this is possible, and this is
easier to do after the usual change of variable, so we want to set
$$
g_t(z) = \lambda^{-1} \psi^{-1}(\wt g_t(z))
\ \hbox{ for } z\in U_{int},
\leqno (16.28)
$$
where
$$
\wt g_t(z) = (2-2t) \wt f(z) + (2t-1) \wt g(z),
\leqno (16.29)
$$
and $\wt g(z)$ is as in (16.4) when $z\in V_1$, 
$\wt g(z) = \wt g_j(z)$ when $z \in {1 + a \over 2} B_j$
for some $j\in J_2$,
$\wt g(z) = \wt g_j(z) = \wt g_{j,x}(z)$
when $z \in B_{j,x}^+$ for some $j \in J_3$ and $x\in Z(y_j)$,
and $\wt g(z) = \wt f(z)$ when $z\in \R^n \sm V_1^+$.
On $V_1^+$, this definition is the same as in the remark 
below (16.8), which was also based on (16.3)-(16.4)
(also see (13.28) and (13.29)), (14.29), and (15.53).

We need to check that
$$
\wt g_t(z) \in B(0,1)
\ \hbox{for } t \geq 1/2, 
\leqno (16.30)
$$
so that (16.28) makes sense.  This is clear when
$z\in \R^n \sm V_1^+$, because 
$\wt g_t(z)=\wt f(z) = \psi(\lambda f(z)$ by (12.75);
otherwise, we already checked that $\wt g(z) \in B(0,1)$
(typically, when we wanted to define $g$ by 
$g(z) = \lambda^{-1} \psi^{-1}(\wt g(z))$); see (12.82), above
(14.29), and below (15.53). Then $\wt g_t(z)$, which lies on the
segment between $\wt g(z)$ and $\wt f(z)= \psi(\lambda f(z)$,
lies in $B(0,1)$ too. Thus (16.28) makes sense and $g_t(z) \in U$
for $z\in U_{int}$ and $1/2 \leq t \leq 1$.

When $t=1/2$, (16.28) and (16.29) yield $\wt g_t = \wt f$ 
and $g_t=f= \varphi_1$, so $g_t$ is continuous across $t=1/2$. 
When $t=1$, we retrieve $\wt g_1 = \wt g$ and $g_1 = g$.

We only defined $\wt g_t(z)$ and $g_t(z)$ when $z\in U_{int}$;
when $z \in U \sm U_{int}$, we do not define $\wt g_t(z)$
and directly set $g_t(z) = f(z)$, as in (16.8). This does not
create a discontinuity, because $V_1^+$ lies well inside
$U_{int}$ (recall the definition (12.72) and the inclusions in
(12.76), the lines above (14.25), and (15.47)), and because
the definition above also gives $g_t(z) = f(z)$
when $z \in U_{int} \sm V_1^+$.

Now we check that Lemma 16.14 is still valid in the
present case. We do not need to change anything
before the last line of the proof of (16.17), where
we just need to observe that (again for $t \geq 1/2$ and $x\in B$)
$|g_t(x)-f_t(x)| \leq \lambda^{-1} \Lambda |\wt g_t(x)- \wt f(x)|
\leq 4 \Lambda^2 (1+|f|_{lip}) \delta_6$
by (16.28) and the proof of (16.11)
(more precisely, the line above (13.35), (12.8), (14.25),
and (15.51), but if you are ready to loose an extra $\Lambda^2$,
just use (16.11)); so $g_t(x) \in B'$ as before.

Thus we may turn to (1.7), or equivalently (16.18) or, after a change 
of variable, the fact that 
$$
\wt g_t(x) \in \wt L_k = \psi(\lambda L_k)
\ \hbox{ when $x \in E \cap L_k \cap B'$ and $0\leq t \leq 1$.}
\leqno (16.31)
$$
The verification for $0 \leq t \leq 1/2$ is the same as before,
so we may assume that $t \geq 1/2$, and by (16.29) we just need to 
check that
$$
[\wt f(x),\wt g(x)] \i \wt L_k
\hbox{ for } x \in E \cap L_k \cap B
\leqno (16.32)
$$
(compare with (16.19)).

We continue the argument as below (16.19), starting with the 
case when $x\in V_1$ and so $x\in B_j$ for some $j\in J_1$.
Let $l \in {\cal L}$ be as before; thus $f(x_j) \in D_l$ and
$y_l = f(x(l))$ for some $x(l) \in X_5$. We shall also
use $\wt y_l = \psi(\lambda y_l) \in B(0,1)$, and $m \in [0,n]$
such that $y_l \in {\cal S}_m \sm {\cal S}_{m-1}$ 
(just $y_l \in {\cal S}_m$ if $m=0$); then (16.20) holds
as before. Still denote by $F_l$ the smallest face of
the twisted grid that contains $y_l$, set 
$\wt F_l = \psi(\lambda F_l)$ (the smallest rigid face
that contains $\wt y_l$),
and call $\wt W(y_l)$ the affine space spanned by $\wt F_l$.
Next we check that (16.22) holds, or equivalently that
$$
\hbox{any face $\wt F$ of the true grid that contains $\wt f(x)$
contains $\wt F_l$ too.}
\leqno (16.33)
$$
The proof needs to be modified slightly. From (16.20) we deduce that
$$
\dist(\wt y_l,\wt{\cal S}_{m-1}) 
= \dist(\psi(\lambda y_l),\psi(\lambda {\cal S}_{m-1}))
\geq \Lambda^{-1} \lambda \delta_2.
\leqno (16.34)
$$
We still have (16.21), with the same proof, which yields
$$
|\wt f(x)-\wt y_l| = |\psi(\lambda f(x))-\psi(\lambda y_l)|
\leq \lambda \Lambda |f(x)-y_l| \leq \lambda \Lambda \delta_6 (1+|f|_{lip}).
\leqno (16.35)
$$
Denote by $a_i$ and $b_i$ the coordinates of $\wt f(x)$ and $\wt y_l$;
now
$$
|b_i-a_i| \leq \lambda \Lambda \delta_6 (1+|f|_{lip})
\leq {\lambda \over 10 \Lambda} \, \min(\delta_2, \lambda^{-1} r_0)
\leq {1 \over 10}\, \min(\Lambda^{-1} \lambda \delta_2, r_0).
\leqno (16.36)
$$
still by (12.7). From this and (16.34) we deduce the analogue of
(16.22) as before, when we had (16.23) and (16.20).

Now we use the fact that 
$\wt g_j(x) \in [\wt f(x),\wt\pi_j(\wt f(x))]$, by (12.77)-(12.79),
where $\wt\pi_j$ is the orthogonal projection onto the affine
plane $\wt P_j$ that contains $\wt Q_j$; see the description
above (12.77), and recall that $\wt P_j$ satisfies (12.23).

Let $\wt\pi$ be the projection onto the affine plane through
$\wt y_l$ parallel to $\wt P_j$; the analogue of (16.24) is 
$$\eqalign{
|\wt\pi_j(\wt f(x))-\wt y_l| &\leq |\wt\pi(\wt f(x))-\wt y_l|
+||\wt \pi-\wt\pi_j||_{\infty}
\cr&
\leq |\wt f(x)-\wt y_l|+ 
2\lambda \Lambda (1+|f|_{lip})t
\cr&
\leq \lambda \Lambda |f(x)-y_l|
+ 2\lambda \Lambda (1+|f|_{lip})t
\cr&
\leq 3\lambda \Lambda (1+|f|_{lip})t
\leq 3\lambda \Lambda (1+|f|_{lip}) \delta_6 
\cr&
\leq {3 \lambda \delta_2 \over 10 \Lambda}
< \dist(\wt y_l,\d \wt F_l)
}\leqno (16.37)
$$
by (12.23), (16.21), (12.7), and (16.34).

As before, $\wt\pi_j(\wt f(x))$ lies on the affine plane 
$\wt P_j$ that contains $\wt Q_j$, which is contained in
$\wt W(y_l)$ by (12.23); since $\wt W(y_l)$ is the 
affine space spanned by $\wt F_l$, and $\wt y_l \in \wt F_l$, 
we get that $\wt\pi_j(\wt f(x)) \in \wt F_l \i \wt F$ for any
(straight) face $\wt F$ that contains $\wt f(x)$ (by (16.37)).
The rest of the proof of (16.33) (by convexity) goes as before.

The other cases are easier (see near (16.26)); 
we replace Lemma 12.27 with Lemma~12.40 when needed, and
otherwise proceed as above. This completes our 
proof of Lemma 16.14 under the Lipschitz assumption. $\dagger$

\bigskip\noindent 
{\bf 17. Magnetic projections onto skeletons,
and a deformation for the $E_k$.}
\ms

We just checked that $g$ and the $g_t$ define
(a hopefully stabler) acceptable deformation for $E$,
but we still want to modify them so that they work for
the $E_k$, at least for $k$ large.
For this we will need some way to push points back
to the $L_j$ (when they are close to the $L_j$).
The name magnetic for the projections below
was used in [Fv1]  
in a similar context; it is nice because it conveys the idea 
of a strong attraction, but with a very short range.

\msi{\bf 17.a. Magnetic projections onto the faces.}
\ms
We start with a projection on nearby faces of a given dimension, 
and then we shall see how to work in all dimensions at the same time.
In what follows, $m \in [0,n)$ is an integer, and $s$ is a 
small number that plays the role of an attraction range, 
which will later depend on various parameters.
Also recall that ${\cal S}_m$ denotes the $m$-dimensional
skeleton of our usual dyadic grid.

\ms\proclaim Lemma 17.1.
Let a dimension $m \in [0, n[$ and $s\in (0,{r_0 \over 10})$
be given. There is a mapping 
$\Pi = \Pi_{m,s} : \R^n \to \R^n$, with the following
properties:
$$
\Pi_{m,s}(x) = x
\ \hbox{ when $x\in {\cal S}_m$ and when } 
\dist(x,{\cal S}_m) \geq 2s,
\leqno (17.2)
$$
$$
\Pi_{m,s}(x) \in {\cal S}_m 
\ \hbox{ when } \dist(x,{\cal S}_m) \leq s,
\leqno (17.3)
$$
$$
\Pi_{m,s}(x) \hbox{ is a $C$-Lipschitz 
function of $s\in (0,{r_0 \over 10})$ and } x\in \R^n,
\leqno (17.4)
$$
where $C$ depends only on $n$, and
$$
\Pi_{m,s}  \hbox{ preserves all the faces of our usual grid,}
\leqno (17.5)
$$
which means that if $F$ is a face of any dimension, then
$\Pi_{m,s}(x) \in F$ for $x\in F$.

\ms
We start with the (rigid) case when $r_0=1$. Naturally we shall
use Lemma 3.17, with $L = {\cal S}_m$ and $\eta = 1/3$; 
we get a mapping $\Pi_L : L^\eta \times [0,1] \to \R^n$,
with the properties (3.18)-(3.22). Recall that $L^\eta$ is, as in 
(3.5), an $\eta$-neighborhood of $L$. For convenience, we extend 
$\Pi_L$ by setting $\Pi_L(x,0) = x$ for $x\in \R^n$; this is 
compatible with (3.18).

Let $\theta : [0,+\infty) \to [0,1]$ be a smooth cut-off function
such that 
$$
\theta(t) = 1 \hbox{ for } 0 \leq t \leq 1, \ 
0 \leq \theta(t) \leq 1 \hbox{ for } 1 \leq t \leq 2, \ 
\theta(t) = 0 \hbox{ for } t \geq 2,
\leqno (17.6)
$$
and $|\theta'(t)| \leq 2$ everywhere. Set $d(x) = \dist(x,{\cal S}_m)$
for $x\in \R^n$, and then
$$
\Pi_{m,s}(x) = \Pi_L(x,\theta(s^{-1}d(x)))
\ \hbox{ for $x\in \R^n$ and $0 < s \leq 10^{-1}$.}
\leqno (17.7)
$$
First observe that if $x\in \R^n \sm L^\eta$,
then $d(x) \geq \eta = 1/3$, $s^{-1}d(x) \geq 10/3$, 
hence $\theta(s^{-1}d(x))= 0$ and $\Pi_L(x,\theta(s^{-1}d(x)))$
is well defined (and is equal to $x$). 
So $\Pi_{m,s}$ is well defined on $\R^n$.

The second part of (17.2) holds for the same reason: if 
$\dist(x,{\cal S}_m) \geq 2s$, then $\theta(s^{-1}d(x))= 0$
and $\Pi_{m,s}(x) = \Pi_L(x,0) = x$ by (3.18).
Similarly, if $\dist(x,{\cal S}_m) \leq s$, then $\theta(s^{-1}d(x))= 1$
and $\Pi_{m,s}(x) = \Pi_L(x,1) \in  L = {\cal S}_m$, by
(3.19) and the definition of $\pi_L$ in Lemma 3.4, so (17.3)
holds. Finally, if $x\in {\cal S}_m$, 
$\Pi_L(x,t) = x$ for all $t$, by (3.18), so $\Pi_{m,s}(x) =x$
and the first part of (17.2) holds too.

Let us check that $\Pi_{m,s}(x)$ is Lipschitz in $x$.
First consider $x, y \in L^\eta$; then
$$\eqalign{
|\Pi_{m,s}(x)&-\Pi_{m,s}(y)| 
= |\Pi_L(x,\theta(s^{-1}d(x)))-\Pi_L(y,\theta(s^{-1}d(y)))| 
\cr&
\leq |\Pi_L(x,\theta(s^{-1}d(x)))-\Pi_L(x,\theta(s^{-1}d(y)))|
\cr&\hskip4cm
+ |\Pi_L(x,\theta(s^{-1}d(y)))-\Pi_L(y,\theta(s^{-1}d(y)))|
\cr&
\leq C d(x) |\theta(s^{-1}d(x))-\theta(s^{-1}d(y))|
+ C |x-y|
}\leqno (17.8)
$$
by (3.20) and (3.21). Let us check that
$$
d(x) |\theta(s^{-1}d(x))-\theta(s^{-1}d(y))| \leq 6|x-y|.
\leqno (17.9)
$$
If $d(x) \leq 3s$, simply say that
$|\theta(s^{-1}d(x))-\theta(s^{-1}d(y))| \leq 2 s^{-1} |d(x)-d(y)|
\leq 2 s^{-1} |x-y|$, and (17.9) follows. If $d(x) \geq 3s$
and $d(y) \geq 2s$, then $\theta(s^{-1}d(x))=\theta(s^{-1}d(y))=0$
and (17.9) is trivial. In the last case when 
$d(x) \geq 3s$ and $d(y) \leq 2s$, 
$$
d(x) |\theta(s^{-1}d(x))-\theta(s^{-1}d(y))| = d(x) \theta(s^{-1}d(y)) 
\leq d(x) \leq 3|d(x)-d(y)| \leq 3|x-y|,
\leqno (17.10)
$$
and (17.9) holds too. Then 
$|\Pi_{m,s}(x)-\Pi_{m,s}(y)| \leq C |x-y|$,
by (17.8) and (17.9), and this takes care of our first case when 
$x, y \in L^\eta$.

Suppose $x\in L^\eta$ and $y\in \R^n \sm L^\eta$, and let
$z\in [x,y]$ lie on the boundary of $L^\eta$; then
$\Pi_{m,s}(z) = z$ and $\Pi_{m,s}(y)=y$ by (17.2), and
$$\eqalign{
|\Pi_{m,s}(x)-\Pi_{m,s}(y)| 
&\leq |\Pi_{m,s}(x)-\Pi_{m,s}(z)|+|\Pi_{m,s}(z)-\Pi_{m,s}(y)|
\cr&
= |\Pi_{m,s}(x)-\Pi_{m,s}(z)|+|z-y|
\cr&
\leq C |x-z|  + |z-y| \leq C |x-y|
}\leqno (17.11)
$$
by the previous case. The case when $x\in \R^n \sm L^\eta$ and 
$y\in L^\eta$ is similar, and when $x,y\in \R^n \sm L^\eta$ we
simply get that $|\Pi_{m,s}(x)-\Pi_{m,s}(y)|=|x-y|$ by (17.2).
So $\Pi_{m,s}$ is $C$-Lipschitz.

For the Lipschitz dependence on $s$, first let $x\in L^\eta$
and $0 \leq s \leq t \leq 10^{-1}$ be given. Then
$$\eqalign{
|\Pi_{m,s}(x)-\Pi_{m,t}(x)|
&= |\Pi_L(x,\theta(s^{-1}d(x)))-\Pi_L(x,\theta(t^{-1}d(x)))| 
\cr&
\leq C d(x) |\theta(s^{-1}d(x))-\theta(t^{-1}d(x))|
}\leqno (17.12)
$$
by (3.20). If $d(x) \leq 3 s \leq 3t$, then
$$
d(x) |\theta(s^{-1}d(x))-\theta(t^{-1}d(y))|
\leq 2 d(x) \Big|{d(x) \over s} - {d(x) \over t}\Big|
= 2 d(x)^2 \, {|s-t| \over st}
\leq 18 |s-t|,
\leqno (17.13)
$$
and we are happy. If $d(x) \geq 2t \geq 2 s$, then
$\theta(s^{-1}d(x)) = \theta(t^{-1}d(x)) = 0$ by
(17.6), and we are happier. We are left with the case when
$3s \leq d(x) \leq 2t$; then
$$
d(x) |\theta(s^{-1}d(x))-\theta(t^{-1}d(y))|
= d(x) \theta(t^{-1}d(y)) \leq d(x) \leq 2t \leq 6 (t-s)
\leqno (17.14)
$$
by (17.6) and because $3s \leq 2t$. This takes care
of the case when $x\in L^\eta$. The other case is trivial, since
$\Pi_{m,s}(x)=\Pi_{m,t}(x) = x$ when $x\in \R^n \sm L^\eta$.
So $\Pi_{m,s}(x)$ is Lipschitz in $s$ too, and (17.4) holds.

Finally, (17.5) is a direct consequence of the fact that
$\Pi_L$ preserves the faces too, by (3.22).

We still need to prove the lemma when $r_0 < 1$; denote by
$\Pi'_{m,s}$ the mapping that we just obtained for the unit grid;
naturally, we set
$$
\Pi_{m,s}(x) = r_0 \Pi'_{m,r_0^{-1}s}(r_0^{-1} x)
\ \hbox{ for $x\in \R^n$ and }
0 \leq s \leq 10^{-1} r_0 ;
\leqno (17.15)
$$
the properties (17.2), (17.3), and (17.5) follow at once by 
conjugation, and for (17.4) a rapid inspection shows that the
two Lipschitz constants for $\Pi_{m,s}(x)$ do not even depend on $r_0$.
(We don't really need to know this, but it feels better.)
Lemma 17.1 follows.
\qed

\ms
We shall need to know that 
$$
|\Pi_{m,s}(x)-x| \leq C \Min\big(s,\dist(x,{\cal S}_m)\big)
\ \hbox{for $x\in \R^n$ and $0 \leq s \leq 10^{-1} r_0$.}
\leqno (17.16)
$$
And indeed, by (17.2) we may assume that 
$d(x) = \dist(x,{\cal S}_m) \leq 2s$, because otherwise $\Pi_{m,s}(x)=x$.
Then pick $z\in {\cal S}_m$ such that $|z-x| = d(x)$, and observe 
that
$$\eqalign{
|\Pi_{m,s}(x)-x| &\leq |\Pi_{m,s}(x)-\Pi_{m,s}(z)| + |\Pi_{m,s}(z)-x|
\cr&
= |\Pi_{m,s}(x)-\Pi_{m,s}(z)| + |z-x|
\leq C |z-x| = C d(x)
}\leqno (17.17)
$$
because $\Pi_{m,s}(z)=z$ by (17.2); (17.16) follows.

Next we want a version of Lemma Lemma 17.1 that works for all the 
dimensions $m$ at the same time; naturally we shall obtain it by 
composing mappings $\Pi_{m,s}(x)$ provided by Lemma 17.1.
We keep our usual dyadic grid of mesh $r_0$.

\ms\proclaim Lemma 17.18.
There is a mapping
$\Pi : \R^n \times [0,10^{-1}r_0] \to \R^n$, with the following
properties:
$$
|\Pi(x,s) - x| \leq C s
\ \hbox{ for $x\in \R^n$ and $0 \leq s \leq 10^{-1}r_0$},
\leqno (17.19)
$$
$$
\Pi(x,s) \in F
\ \hbox{ when $F$ is any face of the grid, $x\in \R^n$, and }
\dist(x,F) \leq C^{-1} s,
\leqno (17.20)
$$
$$
\Pi \hbox{ is $C$-Lipschitz on } \R^n \times [0,10^{-1}r_0] 
\leqno (17.21)
$$
and 
$$
\hbox{ every } \Pi(\cdot,s) 
\hbox{ preserves all the faces of our usual grid.}
\leqno (17.22)
$$

\ms
For $s \in [0,10^{-1}r_0]$, set 
$$
s_m = (6C)^{-m} s 
\ \hbox{ for } 0 \leq m \leq n-1,
\leqno (17.23)
$$
where $C$ is the constant of (17.16) (chosen so that $C \geq 1$)
and then 
$$
\Pi(x,s)
= \Pi_{0,s_0}\circ\Pi_{1,s_1}\cdots\circ\Pi_{n-1,s_{n-1}}(x)
\leqno (17.24)
$$
for $x\in \R^n$. Notice that $\Pi_{m,s_{m}}$ is well defined, because
$0 \leq s_m \leq 10^{-1}r_0$, and that (17.19) holds (with a larger 
constant $C$) by successive applications of (17.16). Also,
(17.21) follows from (17.4) and the chain rule, and (17.22)
is a consequence of (17.5). 

We are left with (17.20) to check. 
Let $F$ be a face and $x\in \R^n$ be such that 
$$
\dist(x,F) \leq  s_{n-1} = (6C)^{1-n} s;
\leqno (17.25)
$$
we want to check that $\Pi(x,s) \in F$.
Set $x_{n+1} = x_n = x$, then $x_{n-1} = \Pi_{n-1,s_{n-1}}(x)$,
and by induction 
$$
x_{k} = \Pi_{k,s_k}(x_{k+1})
\ \hbox{ for } 0 \leq k \leq n-1. 
\leqno (17.26)
$$
Thus $\Pi(x,s)= x_0$. Notice that
$$
|x_k-x| \leq C (s_{n-1} + \cdots s_k)
\leqno (17.27)
$$
by successive applications of (17.16), and where
for the few next lines $C$ will stay the same as in 
(17.16) and (17.23).

Let $m$ denote the dimension of $F$; observe that
$$
\dist(x_{m+1},F) \leq \dist(x,F) + |x_{m+1}-x| 
\leq s_{n-1} + C\sum_{k > m} s_k.
\leqno (17.28)
$$
Next denote by $l$ the smallest nonnegative integer such that
$$
\dist(x_{l+1},F') 
\leq s_{n-1} + 4C\sum_{k > l} s_k
\leqno (17.29)
$$
for some face $F' \i F$ of dimension $l$.
Thus $l \leq m$, by (17.28). Let us check that
$$
s_{n-1} + 4C\sum_{k > l} s_k \leq s_l.
\leqno (17.30)
$$
If $l = n-1$, (17.30) holds because the left-hand side is 
$s_l$. If $l < n-1$,
$$
s_{n-1} + 4C\sum_{k > l} s_k
\leq 5C\sum_{k > l} s_k
\leq 5C s_l \sum_{j \geq 1} (6C)^{-j}
\leq {5 s_l \over 6} \sum_{j \geq 0} 6^{-j}
= s_l
\leqno (17.31)
$$
because we assumed that $C \geq 1$; so (17.30) holds. Now
$$
\dist(x_{l+1},{\cal S}_l) \leq \dist(x_{l+1},F')
\leq s_{n-1} + 4C\sum_{k > l} s_k
\leq s_l 
\leqno (17.32)
$$
by (17.29) and (17.30), so (17.3) says that
$x_l = \Pi_{l,s_l}(x_{l+1})$ lies in ${\cal S}_l$. 

If $x_l \in F$, we are happy because all the later 
$\Pi_{k,s_k}$ preserve the faces, so $\Pi(x,s)= x_0$ 
lies in $F$ too.
So assume that $x_l \notin F$. Let $F''$ denote a face of 
dimension $l$ that contains $x_l$, and notice that 
$F'' \neq F'$ because $x_l \notin F'$ since $F' \i F$. 
Use (17.29) to choose $z \in F'$ such that
$$
|z-x_{l+1}| \leq s_{n-1} + 4C\sum_{k > l} s_k.
\leqno (17.33)
$$
If $l > 0$, (3.8) says that
$$
\dist(z,\d F') \leq \dist(z,F'') \leq |z-x_l|
\leq |z-x_{l+1}| + C s_l
\leq s_{n-1} + 4C\sum_{k > l} s_k
+ C s_l
\leqno (17.34)
$$
by (17.26), (17.16), and (17.33), and so
$$
\dist(x_l,\d F') \leq \dist(z,\d F') + |z-x_l|
\leq 2s_{n-1} + 8C\sum_{k > l} s_k
+ 2C s_l.
\leqno (17.35)
$$
Since $2s_{n-1} + 8C\sum_{k > l} s_k
\leq 2s_l$ by (17.30), we get that
$$
\dist(x_l,\d F') \leq 4C s_l
\leq 4C \sum_{k > l-1} s_k,
\leqno (17.36) 
$$
which contradicts the minimality of $l$, because $\d F' \i F' \i F$.
So in fact $l = 0$, and $F'$ and $F''$ are just points of the grid.
Then $F'' = \{ x_l \}$ and $F'= \{ z \}$, and these points are 
distinct. But the last part of (17.34) is still valid, and says that
$|z-x_l|$ is very small. This contradiction shows that $x_l \in F$
was the only option, and completes our proof of (17.20);
Lemma 17.18 follow.
\qed

\msi{\bf 17.b. A stable deformation for the $E_k$.}
\ms

Recall from Section 16 that we have defined a family of mappings
$g_t$, $0 \leq t \leq 1$, that satisfy the constraints (1.4)-(1.8)
with respect to our limit set $E$. We want to use the magnetic
projection given by Lemma 17.18 to modify the $g_t$ and make them
work for the $E_k$ as well.
As usual, we start in the rigid case.

Let $\varepsilon_0$ be small, to be chosen below,
and set
$$
h_t(x) = \Pi(g_t(x),s_t(x))
\ \hbox{ for $x\in U$ and } 0 \leq t \leq 1,
\leqno (17.37) 
$$
where we set
$$
s_t(x) = C \Min(\varepsilon_0, |g_t(x)-x|),
\leqno (17.38) 
$$
where $C$ is the constant of (17.20). 
We shall choose $\varepsilon_0$ much smaller than $(10C)^{-1}r_0$, 
so $h_t(x)$ is well defined. Observe that since 
$0 \leq s_t(x) \leq C \varepsilon_{0}$, (17.19) yields
$$
|h_t(x)-g_t(x)| = |\Pi(g_t(x),s_t(x))-g_t(x)| 
\leq C s_t(x) \leq C \varepsilon_0
\leqno (17.39) 
$$
for $x\in U$ and $0 \leq t \leq 1$ (and with a new constant $C$).
We are interested in the following.

\ms\proclaim Lemma 17.40. 
For $k$ large enough, the mappings $h_t$, $0 \leq t \leq 1$, 
satisfy the conditions (1.4)-(1.8), relative to $E_k$ and the ball
$B'' = \overline B(X_0, R'')$, where
$$
R'' = R_0 +  4\Lambda^2 (1+ |f|_{lip}) \delta_6 + C \Lambda \varepsilon_0.
\leqno (17.41) 
$$

\ms
We give the statement with $\Lambda$ because it will be 
valid in the Lipschitz case, but for the moment we may take 
$\Lambda = 1$.

There is still no difficulty with (1.4) and (1.8), since we merely 
composed $g_t(x)$ with continuous functions of $x$ and $t$, which
happen to be Lipschitz when $t=1$.
Notice that by (17.19),
$$
h_t(x) = g_t(x) = x
\ \hbox{ when } g_t(x) = x,
\leqno (17.42)
$$
because then $s_t(x) = 0$. Because of this,
$$
h_0(x) = x \ \hbox{ for } x\in U,
\leqno (17.43)
$$
by (16.12) and because $\varphi_0(x) = x$ for $x\in U$,
by (11.14). Let us also check that
$$
h_t(x) = x \ \hbox{ for $x\in U_{ext}$ and } 0 \leq t \leq 1,
\leqno (17.44)
$$
where $U_{ext} = \big\{ x\in \R^n \, ; \, \dist(x, \R^n \sm U) \leq 
\delta_0/2 \}$ as in (11.2). Let $x\in U_{ext}$ be given. By (11.3),
$\varphi_t(x) = x$ for  $0 \leq t \leq 1$, and in particular
$f(x) = x$. We will be finished as soon as we check that
$g(x) = f(x)$, because then $g_t(x) = x$ for all $t$, by (16.12) and 
(16.13), and we can apply (17.42).

But $\dist(x,\wh W) \geq \delta_0/2$ because 
$\delta_0 = \dist(\wh W, \R^n \sm U)$ by (11.2), so 
$\dist(x,X_1) \geq \delta_0/2$ because $X_1 \i X_0 \i \wh W$
(by (11.20)). Finally, $\delta_0/2 > \delta_6$ by (12.7), so
the first part of (16.15) says that $x\in U \sm V_1^+$, 
and then $g(x) = f(x)$ by (16.8); (17.44) follows. 

\ms
For the next verifications, we shall often need to
restrict to $E_k$ (we shall not have enough information
on the values of the $\varphi_t$ far from $E$), and
we shall find it more convenient to work on the set
$$
H = \big\{ x\in U \, ; \, \dist(x,\R^n \sm U) \geq \delta_0/2 \big\}
\supset \R^n \sm U_{ext}.
\leqno (17.45)
$$
because $H$ is a compact subset of $U$ and it will be easier to
use our assumption (10.4) (i.e., the convergence of the $E_k$ to $E$)
on that set. Indeed, set
$$
d_k = \sup_{x\in E_k \cap H} \dist(x,E);
\leqno (17.46)
$$
it is easy to deduce from (10.4)-(10.6) that 
$\lim_{k \to +\infty} d_k = 0$ (cover $H$ with a finite set of balls). 

We claim that for $k$ large,
$$
h_t(x) = x \ \hbox{ for } 0 \leq t \leq 1
\hbox{ when } x\in E_k \sm B(X_0,R_0+\varepsilon_0).
\leqno (17.47)
$$
Because of (17.44), it is enough to prove this when 
$x\in E_k \cap H \sm B(X_0,R_0+\varepsilon_0)$. Then
$$
\dist(x,E) \leq d_k < {1 \over 2} \Min(\varepsilon_0, \delta_1)
\leqno (17.48)
$$
for $k$ large. 

For each $t \in [0,1]$, the set $W_t$ of (11.13) is contained 
in $B = \overline B(X_0,R_0)$, by (1.5) and (11.1), 
so $\dist(x,E) < \varepsilon_0/2 \leq {1 \over 2} \dist(x,W_t)$
by (17.48) and because $x \notin B(X_0,R_0+\varepsilon_0)$.
This is good, because (11.12) says that then $\varphi_t(x) = x$.
[This is not a surprise; recall that we computed $\varphi_t(x)-x$ by 
Whitney-extending the values of $\varphi_t(\xi)-\xi$ on $E \cup E_{ext}$, 
which happen to vanish near $x$.] 

If we also prove that $g(x) = f(x)$,
(16.12) and (16.13) will say that $g_t(x) = x$ for all $t$, and the 
result will follow by (17.42).
For this, it is enough to prove that $x \notin V_1^+$, by (16.8).
By (17.48), we can find $z\in E$ such that 
$|z-x| < {1 \over 2} \Min(\varepsilon_0, \delta_1)$.
In particular, $z \notin B$ (because $x \notin B(X_0,R_0+\varepsilon_0)$), 
so (1.5) says that $f(z) = \varphi_1(z) = z$. 
Thus $z\in \R^n \sm W_f$ (see (11.19)), 
and so $\dist(x,\R^n \sm W_f) \leq |z-x| < \delta_1/2$, and indeed
this makes $x \in V_1^+$ impossible, by the second part of (16.15). 
This proves our claim (17.47), and (1.5) (for $E_k$ and the $h_t$,
and with a slightly larger ball) follows.

Before we prove (1.6), let us check that for $x\in U$ and
$0 \leq t \leq 1$, we can find $s \in [0,1]$ such that
$$
|g_t(x)-\varphi_s(x)| \leq 4 \Lambda^2 (1+|f|_{lip}) \delta_6.
\leqno (17.49)
$$
When $t \leq 1/2$, we just take $s = t/2$ and observe that $g_t(x) = 
\varphi_s(x)$ by (16.12). When $t \geq 1$, we take $s=1$ and
observe that $|g_t(x)-\varphi_1(x)| 
\leq |g(x)-f(x)| \leq 4 \Lambda^2 (1+|f|_{lip}) \delta_6$,
by (16.13) and (16.11). So (17.49) holds. Notice also that it implies 
that
$$
|h_t(x)-\varphi_s(x)| 
\leq 4 \Lambda^2 (1+|f|_{lip}) \delta_6 + C \varepsilon_0,
\leqno (17.50)
$$
by (17.39).

We are ready to prove (1.6). In fact, we just need to prove that for $k$ large,
$$
h_t(x) \in B'' \hbox{ when } x\in E_k \cap B(X_0,R_0+\varepsilon_0),
\leqno (17.51)
$$
where $B''$ is as in the statement of Lemma 17.40,
since (17.47) says that $h_t(x) = x \in B''$ when 
$x\in E_k \cap B'' \sm B(X_0,R_0+\varepsilon_0)$. 
Pick $z\in E$ such that $|z-x| < {4 \over 3} \dist(x,E)
\leq {4 \over 3} d_k$. For $k$ large enough,
$z\in B(X_0,R_0+2\varepsilon_0)$, and by (1.5) and (1.6) 
for $E$ and the $\varphi_t$, 
$\varphi_t(z) \in B(X_0,R_0+2\varepsilon_0)$ for $0 \leq t \leq 1$.

Also let $H'$ denote a compact neighborhood of $H$ in $U$.
The function $(y,t) \in H' \times [0,1] \to \varphi_t(x)$
is uniformly continuous; if $k$ is large enough,
$z \in H'$ because $x\in H$, and $|z-x|$ is so small that
$|\varphi_t(z)-\varphi_t(x)| \leq \varepsilon_0$. Then
$\varphi_t(x) \in B(X_0,R_0+3\varepsilon_0)$ for 
$0 \leq t \leq 1$. We then use (17.50) and get that
$h_t(x) \in B''$, if the constant $C$ in (17.41)
is large enough). This proves (1.6) for the $h_t$.

\ms
Finally we need to prove that (1.7) holds, relative to $E_k$, and for 
$k$ large. As before, we shall restrict our attention to $H$ first.
Set, for $0 \leq j \leq j_{max}$ and $k \geq 0$,
$$
d_{j,k} = \sup_{x\in L_j \cap E_k \cap H} \dist(x,E\cap L_j);
\leqno (17.52)
$$
we claim that for each $j$,
$$
\lim_{k \to +\infty} d_{j,k} = 0.
\leqno (17.53)
$$
Otherwise, we can find $j \leq j_{max}$ and a sequence
of points $x_k \in E_k \cap L_j \cap H$, for which 
$t_k = \dist(x_k,E\cap L_j)$ does not tend to $0$.
Passing to a subsequence, we may even assume that
$t_k \geq a$ for some $a > 0$, and that
$x_k$ tends to a limit $x_{\infty}$. Then
$x_\infty \in L_j \cap H$ because $L_j \cap H$ is closed,
and $x_\infty \in E$ because $E$ is closed and
$\dist(x_k, E)$ tends to $0$ by (10.4) (see near (17.46)). 
Now the fact that $x_\infty \in L_j \cap E$ contradicts
the fact that $t_k \geq a \,$; our claim (17.53) follows.
Next set
$$
\eta_{j,k} = \sup\big\{ |g_t(x)-g_t(y)| \, ; \,
x \in H, |x-y| \leq 2d_{j,k}, \hbox{ and } 
0 \leq t \leq 1 \big\};
\leqno (17.54)
$$
then $\lim_{k \to +\infty} \eta_{j,k} = 0$, by (17.53)
and because $(x,t) \to g_t(x)$ is uniformly continuous on
$H'\times [0,1]$, where $H'$ a compact neighborhood of
$H$ in $U$.

Let us check that (1.7) holds when $k$ is large enough. 
Let $j \leq j_{max}$ and $x\in E_k \cap L_j$ be given; 
we want to check that $h_t(x) \in L_j$ for $0 \leq t \leq 1$. 
We may assume that $x\in H$, because otherwise
$h_t(x) = x$ by (17.44). 
Pick $y\in E\cap L_j$ such that $|y-x| \leq 2d_{j,k}$; then
$$
\dist(g_t(x),L_j) \leq |g_t(x)-g_t(y)|
\leq \eta_{j,k} \leq \varepsilon_0
\leqno (17.55)
$$
because $g_t(y) \in L_j$ by (1.7) for the $g_t$ relative to $E$, 
and if $k$ is so large that $\eta_{j,k} < \varepsilon_0$.
Also, $\dist(g_t(x),L_j) \leq |g_t(x)-x|$ because $x\in L_j$;
altogether, 
$$
\dist(g_t(x),L_j) \leq C^{-1} s_t(x),
\leqno (17.56)
$$
by (17.38) and where $C$ is as in (17.38) and (17.20). 
Let $F$ be a face of $L_j$ such that 
$\dist(g_t(x),L_j) = \dist(g_t(x),F)$; then 
$$
h_t(x) = \Pi(g_t(x),s_t(x)) \in F
\leqno (17.57)
$$
by (17.37), (17.20), and (17.56). Thus $h_t(x) \in F \i L_j$; 
this completes our proof of (1.7), and Lemma 17.40 follows.
\qed

\ms
We shall also need to know about the analogue, for the mappings
$h_t$ and the set $E_k$, of the sets $W_t$ and $\wh W$. We claim that
for $k$ large and $x\in E_k$, 
$$
h_t(x) = x \ \hbox{ for $0 \leq t \leq 1$ when } 
\dist\big(x, \bigcup_{0 \leq t \leq 1} W_t\big) \geq 2 d_k
\leqno (17.58)
$$
(where $W_t$ is as in (11.13) and $d_k$ as in and (17.46)) and 
$$\eqalign{
\dist(h_t(x), \wh W) &\leq 4 \Lambda^2 (1+|f|_{lip}) \delta_6 
+ C \Lambda \varepsilon_0 
\ \hbox{ for $0 \leq t \leq 1$}
\cr&
\hskip2cm \hbox{ when } 
\dist\big(x, \bigcup_{0 \leq t \leq 1} W_t\big) \leq 2 d_k.
}\leqno (17.59)
$$
The proof will be almost the same as for (17.47) and (17.51).
First let $x\in E_k$ be such that 
$\dist\big(x, \bigcup_{0 \leq t \leq 1} W_t\big) \geq 2 d_k$.
If $x\in U_{ext}$, (17.44) says $h_t(x) = x$, as needed.
So we may assume that $x\in \R^n \sm U_{ext} \i H$
(by (17.45)). Then $\dist(x,E) \leq {1 \over 2} \dist(x,W_t)$ 
for all $t$, by (17.46). 
Hence, by (11.12) $\varphi_t(x) = x$ for $0 \leq t \leq 1$. 
In particular, $f(x) = x$. Next let us check that $x\notin V_1^+$. 
Let $y\in E$ be such that 
$|y-x| < {3 \over 2} \dist(x,E) \leq {3 \over 2} d_k$.
Then $\dist(y,W_1) > d_k/2$, so $f(y) = \varphi_1(y) = y$.
That is, $y \in \R^n \sm W_f$ (see the definition (11.19))
and $\dist(x,\R^n \sm W_f) \leq |y-x| \leq 2 d_k$. For $k$ large
this forces $x\notin V_1^+$, by the second part of (16.15).
Hence $g(x) = f(x) = x$, by (16.8), and
then $g_t(x) = x$ for all $t$, by (16.12) and (16.13).
Finally, $h_t(x) = x$ by (17.42), as needed for (17.58).

Now suppose that 
$\dist\big(x, \bigcup_{0 \leq t \leq 1} W_t\big) \leq 2 d_k$,
and choose $y\in W_t$ be such that 
$|y-x| < {3 \over 2} \dist(x,E) \leq 3 d_k$.
By (2.1) and (2.2), $\varphi_t(y) \in \wh W_t$ for $0 \leq t \leq 1$.
If $x\in U_{ext}$, $h_t(x)  = x$ by (17.44), 
hence $\dist(h_t(x),\wh W) \leq |x-y| 
\leq 3  d_k \leq \varepsilon_0$ (for $k$ large). Otherwise, 
(17.45) says that $x\in H$, and we use the
uniform continuity of $(y,t) \to \varphi_t(y)$ on $H' \times [0,1]$
to show that for $k$ large enough,
$\dist(\varphi_t(x),\wh W) \leq \varepsilon_0$ for
$0 \leq t \leq 1$. Then we apply (17.50) and get that
$\dist(g_t(x),\wh W) \leq 4 \Lambda^2 (1+|f|_{lip}) \delta_6 + 
C \Lambda  \varepsilon_0$ for all $t$, as needed for (17.59).

\ms $\dagger$
Let us now say how we modify all this when we work with the Lipschitz 
assumption.
We don't need to change Lemmas 17.1 and 17.18, but we need to modify 
the definition of the $h_t$. We first set
$$
h_t(x) = x
\ \hbox{ when $x\in U_{ext}$ and } 0 \leq t \leq 1;
\leqno (17.60) 
$$
this is not too shocking, because of (17.44). Let us also
define the $h_t$ on $E_k \cap H_1$, where
$$
H_1 = \big\{ x\in U \, ; \, \dist(x,\R^n \sm U) \geq \delta_0/3\big\}
\leqno (17.61) 
$$
We shall see soon that although the two sets overlap,
our two definitions coincide on their intersection
$E_k \cap H_1 \cap U_{ext}$. 
We modify our original definition by (17.37) and (17.38) and set,
for $x\in E_k \cap H_1$,
$$
\wt s_t(x) = C \Min(\lambda \varepsilon_0, 
|\wt g_t(x)- \psi(\lambda x)|), 
\leqno (17.62) 
$$
where $C$ is still as in (17.20) and $\wt g$ was defined near
(16.29), and then
$$
\wt h_t(x) = \Pi(\wt g_t(x),\wt s_t(x))
\ \hbox{ for } 0 \leq t \leq 1.
\leqno (17.63) 
$$
We naturally intend to set
$$
h_t(x) = \lambda^{-1} \psi^{-1}(\wt h_t(x))
\ \hbox{ for $x\in E_k \cap H_1$ and } 0 \leq t \leq 1,
\leqno (17.64) 
$$
but as usual we first need to check that $\wt h_t(x) \in B(0,1)$,
and this will be easier after some estimates on $\wt h_t(x)$.
Let us first describe some situations where 
$\wt h_t(x) = \psi(\lambda x)$. The analogue of (17.42)
is now that
$$
\wt h_t(x) = \wt g_t(x) = \psi(\lambda x)
\ \hbox{ when $x\in E_k \cap H_1$ and } g_t(x) = x
\leqno (17.65)
$$
(or equivalently when $\wt g_t(x) = \psi(\lambda x)$, see (16.28)
and recall that $x\in H_1 \i \R^n \sm U_{ext}$ by (17.61) and (11.2)).
The proof stays the same: just observe that $\wt s_t(x) = 0$ and 
apply (17.19). Because of this, we shall easily get that 
$\wt h_t(x) = \wt g_t(x) = \psi(\lambda x)$ in the situations
where we proved that $h_t(x) = g_t(x) = x$.

First observe that we still get that $h_t(x) = g_t(x) = x$
when $x\in E_k \cap H_1 \cap U_{ext}$, by the proof of (17.44)
(and where we replace (16.13) with (16.29)). Because of this,
our two definitions on that set coincide, which is good because
we also get that $h_t(x)$ is continuous in both variables and
Lipschitz in $x$.

We also get that 
$$\eqalign{
\wt h_t(x) = \wt g_t(x) = \psi(\lambda x)
&\hbox{ when $t=0$, when $x\in E_k \sm B(X_0,R_0+\varepsilon_0)$,} 
\cr& 
\hbox{ and when }
\dist\big(x, \bigcup_{0 \leq t \leq 1} W_t\big) \geq 2 d_k,
}\leqno (17.66)
$$
as in (17.43), in (17.47), and in (17.58). 
In all these cases, we need $k$ to be large enough
(so that we can define $h_t(x)$, but not only), and then we 
follow the proofs above (except that we replace (6.13) with (16.29) 
and (17.42) with (17.65)).
Of course in all these cases, the formula in (17.64) 
makes sense because $\wt h_t(x) \in B(0,1)$, 
and yields $h_t(x) = g_t(x) = x$.

Next we generalize the formulas (17.49) and (17.50).
We shall restrict to $x\in E_k \cap H_1$,
because for $x\in E_k \sm H_1$, (17.60) will be enough.

Let us first check that for $x\in U_{int}$ (the set defined in 
(12.72) and where $\wt f$ is nicely defined) and $0 \leq t \leq 1$, 
we can find $s \in [0,1]$ such that
$$
|\wt g_t(x)-\psi(\lambda \varphi_s(x))| 
\leq 4 \lambda \Lambda (1+|f|_{lip}) \delta_6.
\leqno (17.67)
$$
When $t \leq 1/2$, we take $s=2t$, and (17.67) holds trivially
because $\wt g_t(x) = \psi(\lambda \varphi_s(x))$ by (16.12) and (12.75). 
When $t \geq 1/2$, we take $s=1$ 
(and hence $\psi(\lambda \varphi_s(x))= \wt f(x)$)
and use (16.29) to get that 
$$
|\wt g_t(x)-\wt f(x)| \leq |\wt g(x)-\wt f(x)|
\leq 4 \lambda \Lambda (1+|f|_{lip}) \delta_6,
\leqno (17.68)
$$
by the line above (13.35), (12.8), (14.26), and (15.51)
(or faster, if you are willing to lose a factor $\Lambda^2$,
by (16.11) and (16.3)). This proves (17.67); since
$$
|\wt h_t(x) - \wt g_t(x)| \leq C \wt s_t(x) \leq C \lambda \varepsilon_0
\leqno (17.69)
$$
by (17.63), (17.19), and (17.62), we deduce from (17.67) that
(still for $x\in E_k \cap H_1$)
$$
|\wt h_t(x)-\psi(\lambda \varphi_s(x))| 
\leq 4 \lambda \Lambda (1+|f|_{lip}) \delta_6 + C \lambda 
\varepsilon_0.
\leqno (17.70)
$$

When $x\in E_k \cap H_1 \sm U_{int}$, we can even say a bit more.
For $t \leq 1/2$, we still have that
$\wt g_t(x)=\psi(\lambda \varphi_s(x))$ with $s=2t$ by (16.12).
For $t \geq 1/2$, we did not want to define $\wt g_t(x)$
directly in Section 16, and instead we set $g_t(x) = f(x)$ directly;
see the definition about nine lines below (16.30). 
This is because if $x$ is any point of $U\sm U_{int}$, 
we cannot be sure that $f(x) \in U$, and the we cannot define
$\wt f(x) = \psi(\lambda f(x))$ or $\wt g_t(x)$. But here we 
only care about points $x\in E_k \cap H_1$. Let $H'_1$ be a compact
neighborhood of $H_1$, with $H'_1 \i U$. Notice that for 
$x\in E\cap H'_1$ and $0 \leq t \leq 1$, 
$\varphi_t(x) \in \wh W \cup H'_1$, either by (2.2) or because
$\varphi_t(x) = x \in H'_1$. Then there is a compact neighborhood 
$H''$ of $E\cap H'_1$, such that $\varphi_t(H'') \i U$
for $0 \leq t \leq 1$.
Since for $k$ large, $E_k \cap H_1 \i H''$, we get that 
$\varphi_t(E_k \cap H_1) \i \i W$ for $k$ large.
That is, for $k$ large we can define $\psi(\lambda \varphi_t(x))$ 
for $x\in E_k \cap H_1$. Then we can set
$\wt g_t(x) = \psi(\lambda \varphi_1(x)) = \psi(\lambda f(x))$
for $t \geq 2$, because $g_t(x) = f(x) = \varphi_1(x)$. 
This is of course better than (17.68).
And we still have (17.69) (with the same proof) and hence (17.70)
in this case. So we get a good definition of $\wt g_t$ and $\wt h_t$ 
on $E_k \cap H_1$.

Recall that we shall not worry about $x\in E_k \sm H_1$, because 
we have the simpler formula (17.60).

Now we want to generalize (17.59) (on the set $E_k \cap H_1$). 
Let us first prove that 
$$
\dist(\wt h_t(x), \psi(\lambda \wh W)) 
\leq 4 \lambda \Lambda (1+|f|_{lip}) \delta_6 
+ C \lambda \varepsilon_0  
\ \hbox{ for $0 \leq t \leq 1$}
\leqno (17.71)
$$
when $x \in E_k \cap H_1$ is such that 
$\dist\big(x, \bigcup_{0 \leq t \leq 1} W_t\big) \leq 2 d_k$.
For such $x$ and $k$, we can find $t \in [0,1]$ and $y\in W_t$ such that
$|y-x| \leq 3d_k$; then $y \in \wh W$ by (2.2), and also
$\varphi_s(y) \in \wh W$ for $0 \leq s \leq 1$, either
because $y\in W_s$ and by (2.2), or because $y \in E \sm W_s$
and $\varphi_x(y)=y\in \wh W$. Thus
$$
\psi(\lambda \varphi_s(y)) \in \psi(\lambda \wh W)
\ \hbox{ for } 0 \leq s \leq 1.
\leqno (17.72)
$$

If $x\in E_k \cap H_1 \cap U_{ext}$, we checked that (17.60) is
still valid; it says that $h_t(x) = x$, 
so $\wt h_t(x) = \psi(\lambda x)$, and (17.71) holds because
$$
\dist(\wt h_t(x), \psi(\lambda \wh W)) 
\leq |\psi(\lambda x)-\psi(\lambda y)| 
\leq \lambda \Lambda |x-y| \leq 3\lambda \Lambda d_k
\leqno (17.73)
$$
(since $y \in \wh W$) and for $k$ large. 
By (17.45) we are left with the case when $x\in H$, and (17.70) says that 
$|\wt h_t(x)-\psi(\lambda \varphi_s(x))| 
\leq 4 \lambda \Lambda (1+|f|_{lip}) \delta_6 + C \lambda 
\varepsilon_0$ for some $s$. In addition, 
$|\psi(\lambda \varphi_s(x))-\psi(\lambda \varphi_s(y))| \leq 
\lambda \varepsilon_0$ if $k$ is large enough, because $|y-x| \leq 3d_k$
and by the usual uniform continuity argument near $H$ (see below (17.51) for 
instance). We deduce (17.71) from this and (17.72).

Notice that
$$
\dist(\psi(\lambda \wh W), \R^n \sm B(0,1))
\geq \lambda \Lambda^{-1} \dist(\wh W, \R^n \sm U)
= \lambda \Lambda^{-1} \delta_0
\geq 10 \lambda \Lambda (1+|f|_{lip}) \delta_6
\leqno (17.74)
$$
by (12.6) and (12.7); thus (17.71) implies that
$\wt h_t(x) \i B(0,1)$ when $k$ is large and $x \in E_k\cap H_1$ 
is such that $\dist\big(x, \bigcup_{0 \leq t \leq 1} W_t\big) \leq 2 d_k$.

If instead $x \in E_k\cap H_1$ but
$\dist\big(x, \bigcup_{0 \leq t \leq 1} W_t\big) > 2 d_k$,
(17.66) says that $\wt h_t(x) = \psi(\lambda x) \in B(0,1)$
(as in 17.58). So we proved that for $k$ large,
$$
\wt h_t(x) \i B(0,1) \hbox{ for } x \in E_k \cap H_1
\hbox{ and } t \in [0,1].
\leqno (17.75)
$$
This is good to know, because now we may define $h_t$ 
on $E_k \cap H_1$ by (17.64), and complete our verifications 
with a free mind. 

Notice for instance that (17.70) implies that for $x\in E_k \cap H_1$
$$
|h_t(x)-\varphi_s(x)| 
\leq 4 \Lambda^2 (1+|f|_{lip}) \delta_6 + C \Lambda \varepsilon_0
\leqno (17.76)
$$
(for the same $s$ that we found for (17.67)); this is almost as good 
as (17.50).

We already checked that (1.4) and (1.8) hold (see above (17.66), 
when we checked that our two definitions coincide on 
$E_k \cap H_1 \cap E_{ext}$), and (1.5) follows from (17.66). 
We still need to check (1.6), but we may now repeat the proof given 
near (17.51), except that when we consider points of $E_k \sm H_1$,
we just use the simpler (16.60). A similar argument applies to the 
final comment (17.59). We even added the constant $\Lambda$ in advance, 
to take care of the extra $\Lambda$ in (17.76).

We are left with the proof of (1.7), which we only need to modify 
slightly. We keep the same definition for the $d_{j,k}$ in (17.52), but 
define the modulus of continuity $\eta_{j,k}$ of (17.54) as before, 
only with respect to the $\wt g_t$; they both tend to $0$ 
for the same reasons as before.

Then suppose that $\eta_{j,k} \leq \lambda \varepsilon_0$, and
let $x\in L_j \cap E_k$ be given. If $x\in E_{ext}$,
(16.60) says that $h_t(x) = x \in L_j$ (as needed), so
we may assume that $x\in H_1$, and then we shall use 
(17.62) and (17.63). Choose $y\in L_j \cap E$ such that 
$|y-x| \leq d_k$, and observe that
$$\eqalign{
\dist(\wt g_t(x), \psi(\lambda L_j))
&\leq \Min\big(|\wt g_t(x)-\wt g_t(y)|, 
|\wt g_t(x) - \psi(\lambda x)|\big)
\cr&
\leq \Min\big(\lambda \varepsilon_0, 
|\wt g_t(x) - \psi(\lambda x)|\big)
\leq C^{-1}\wt s_t(x)
}\leqno (17.77)
$$
since $g_t(y)$ and $x$ both lie in $L_j$,
because $|\wt g_t(x)-\wt g_t(y)| 
\leq \eta_{j,k} \leq \lambda\varepsilon_0$, and by (17.62).
This allows us to use (17.20), get that
$\wt h_t(x) = \Pi(\wt g_t(x),\wt s_t(x)) \in \wt F$,
where $\wt F$ is a face of $\psi(\lambda L_j)$
that lies close to $\wt g_t(x)$, and conclude as before.
So Lemma 17.40 also holds under the Lipschitz assumption.

Notice finally that our final estimates (17.58) and (17.59)
still hold in the present context. For (17.58), 
this follows from (17.66) and (17.64), or directly (17.60);
for (17.59) we use (17.71) and (17.64) if $x\in E_k \cap H$,
and (17.60) and the fact that 
$\dist(h_t(x),\wh W) = \dist(x,\wh W) \leq 2 d_k$ otherwise.
$\dagger$

\msi{\bf 17.c. A last minute modification of our deformation.}
\ms

The family $\{ h_t \}$ that we just constructed is almost
perfect, but we would also like to make the set where
$h_1(x) \neq x$ a little smaller, so that it stays away from 
the boundary of 
$$
W_f = \big\{ x\in U \, ; \, f(x) \neq x \big\}.
\leqno (17.78)
$$
The point is that
the sets $E_k$ could have a large piece in $W_f$ very near $\d W_f$,
while the corresponding piece of $E$ lies in $\d W_f$ and is
not accounted for in some estimates; this could be bad.

So we want to replace $g_t$ and $h_t$ by mappings that
coincide with the identity very near $\d W_f$.
We want to continue the $g_t$, and then the $h_t$,
with a deformation that only moves points near $\d W_f$.
First recall from (16.13) (or (16.28) and (16.29)) and
(16.8) that
$$
g_1(x) = g(x) = f(x) 
\ \hbox{ for } x\in \R^n \sm V_1^+.
\leqno (17.79)
$$
Thus, by the second part of (16.15),
$$
g_1(x) = g(x) = f(x) 
\ \hbox{ when } \dist(x, \R^n \sm W_f) \leq \delta_1/2.
\leqno (17.80)
$$
As usual, we start with the rigid case.
Let us define $g_t$ for $1 \leq t \leq 2$.  
Let $\varepsilon_\ast$ be (much) smaller than $\delta_1/4$, 
to be chosen below. First set
$$
g_t(x) = g_1(x) 
\ \hbox{ for $1 \leq t \leq 2$ when }
\dist(x,\R^n \sm W_f) \geq 2 \varepsilon_\ast.
\leqno (17.81)
$$
When $\varepsilon_\ast \leq \dist(x,\R^n \sm W_f) \leq 2 \varepsilon_\ast$, 
first define
$$\eqalign{
g_2(x) &= {\dist(x,\R^n \sm W_f) -\varepsilon_\ast \over 
\varepsilon_\ast} \, g_1(x) 
+ {2\varepsilon_\ast -\dist(x,\R^n \sm W_f) \over \varepsilon_\ast} 
\, x
\cr&
= {\dist(x,\R^n \sm W_f) -\varepsilon_\ast \over 
\varepsilon_\ast} \, f(x) 
+ {2\varepsilon_\ast -\dist(x,\R^n \sm W_f) \over \varepsilon_\ast} 
\, x,
}\leqno (17.82)
$$
where the identity comes from (17.80); just set
$$
g_2(x) = x \ \hbox{ when } \dist(x,\R^n \sm W_f) < \varepsilon_\ast.
\leqno (17.83)
$$
Notice that $g_2$ is continuous across the two 
obvious boundaries. Now set
$$
g_t(x) = (2-t) g_1(x) + (t-1) g_2(x)
= (2-t) f(x) + (t-1) g_2(x)
\leqno (17.84)
$$
when $\dist(x,\R^n \sm W_f) \leq 2 \varepsilon_\ast$
(i.e., in the two last cases) and for $1 \leq t \leq 2$.

Now we have a complete definition of the $g_t$, $0 \leq t \leq 2$.
Again, $g_t(x)$ is a continuous function of $t$ and $x$, by
construction. We still define the $h_t$, $1 \leq t \leq 2$, 
by the formulas (17.37) and (17.38); observe that
$$
h_t(x) = h_1(x) 
\hbox{ for $1 \leq t \leq 2$ when }
\dist(x,\R^n \sm W_f) \geq 2 \varepsilon_\ast,
\leqno (17.85)
$$
by (17.81) and because $s_t(x)$, and then $h_t(x)$ depend only on the 
values $g_t(x)$. We are happier because
$$
h_2(x) = g_2(x) = x
\ \hbox{ when } \dist(x,\R^n \sm W_f) < \varepsilon_\ast,
\leqno (17.86)
$$
because $s_2(x) = 0$ by (17.83) and (17.38).
Now we want to check that
$$
\hbox{Lemma 17.40 also holds for the $h_{2t}$, $0 \leq t \leq 1$.}
\leqno (17.87)
$$

Since $h_t(x)$, $1 \leq t \leq 2$ is continuous in both variables and
Lipschitz in $x$, (1.4) and (1.8) still hold as before; 
(17.42) still holds for $t \geq 1$, for the same reasons. Then
$$
h_t(x) = g_1(x) = f(x) = x
\hbox{ for $1 \leq t \leq 2$ and }
x \in \R^n \sm W_f,  
\leqno (17.88)
$$
because $f(x) = x$ by definition of $W_f$,
$g_2(x) = x$ by (17.83), hence $g_t(x)=x$ 
by (17.84), and finally $h_t(x)=x$ by the extended (17.42).

For (1.5), it is enough to check that for $k$ large,
$h_t(x) = x$ for $x\in E_k \sm B(X_0,R_0+\varepsilon_0)$.
By (17.47), this is true for $t\leq 1$. Also, $f(x) = x$
(if $k$ is large enough): this was one of the intermediate steps
in the proof of (17.47) (we proved that $\dist(x,W_t) > 0$
when $\dist(x,E) \leq d_k$; otherwise use (17.44) and its proof). 
Then $x\in \R^n \sm W_f$, and 
we can apply (17.88). So (1.5) holds.

For (1.6) we also need to check that $h_t(x)$ does not
escape too far when $x\in E_k \cap B(X_0,R_0+\varepsilon_0)$.
We may restrict to $t>1$, since Lemma 17.40 itself takes care of 
$0 \leq t \leq 1$.

If $\dist(x,\R^n \sm W_f) \leq 2\varepsilon_\ast$,
$$\eqalign{
|g_t(x)-x| &\leq \Max(|g_1(x)-x|,|g_2(x)-x|)
= \Max(|f(x)-x|,|g_2(x)-x|) 
\cr&
\leq |f(x)-x|
\leq (1+|f|_{lip}) \dist(x,\R^n \sm W_f)
\leq 2 (1+|f|_{lip}) \varepsilon_\ast
}\leqno (17.89)
$$
by the various definitions (17.80)-(17.84),
and because $f(x)-x$ is a $(1+|f|_{lip})$-Lipschitz
mapping that vanishes on $\R^n \sm W_f$.
Then $|h_t(x)-x| \leq C (1+|f|_{lip}) \varepsilon_\ast$
too, by (17.37), (17.38), and (17.19).
We get that $h_t(x) \in B''$ if $x\in B(X_0,R_0+\varepsilon_0)$
(and if $\varepsilon_\ast$ is small enough compared to
$\varepsilon_0$ or $\delta_6$; see Remark 11.17 to check 
that we do not cheat). 

In the remaining case when 
$\dist(x,\R^n \sm W_f) \geq 2\varepsilon_\ast$,
(17.81) says that $g_t(x) = g_1(x)$ for $t \geq 1$, so 
$h_t(x) = h_1(x) \in B''$ because we use the same formula 
(17.37), and by (17.51). This proves (1.6).

The verification of (1.7) is a little easier than before.
The only case when we do not already know that $h_t(x) \in L_j$
is when $\dist(x,\R^n \sm W_f) \leq 2\varepsilon_\ast$,
and then (17.89) implies that 
$$
\dist(g_t(x),L_j) \leq |g_t(x)-x| 
\leq 2 (1+|f|_{lip}) \varepsilon_\ast
\leq \varepsilon_0
\leqno (17.90)
$$
if $\varepsilon_\ast$ is small enough compared to $\varepsilon_0$.
Then (17.55) holds and we can conclude as before (that is, $g_t(x)$
is so close to $L_j$ that the magnetic projection sends it back to 
$L_j$). This completes our proof of (17.87).

Let us record the fact that
$$
h_t(x) = x \ \hbox{ for $x\in U_{ext}$ and $0 \leq t \leq 2$;}
\leqno (17.91)
$$
for $t \leq 1$, this comes from (17.44), and the proof of
(17.44), which uses (17.42), also gives that $g(x) = f(x) = x$. Then 
(17.81)-(17.83) yield $g_t(x) = g_1(x) = x$ for $t \geq 1$, and 
hence $h_t(x) = x$ by (17.37) and (17.38) (because $s_t(x) = 0$).

Let us also say a few words about the analogue of $\wh W$ for the 
$h_t$, $0 \leq t \leq 2$. Set
$$
W_{k,t} = \big\{ x\in E_k \, ; \, h_t(x) \neq x \big\}
\hbox{ for } 0 \leq t \leq 2
\ \hbox{ and } \ 
\wh W_k = \bigcup_{0 \leq t \leq 2} W_{k,t} \cup h_t(W_{k,t}).
\leqno (17.92)
$$
We claim that 
$$
\wh W_k \i \i U \ \hbox{ for $k$ large.}
\leqno (17.93)
$$
Because of (17.91), we know that the $W_{k,t}$ do not meet
$U_{ext}$. Also, (17.58) and (17.59) give the desired control
on the $W_{k,t}$, $0 \leq t \leq 1$ (recall that
$4 \Lambda^2 (1+|f|_{lip}) \delta_6 < \delta_0/2 = \dist(\wh W, \R^n 
\sm U)/2$ by (12.6) and (12.7)).
So it is enough to consider $t > 1$. 

Let $x\in W_{k,t}$ be given. If 
$\dist(x,\R^n \sm W_f) \geq 2\varepsilon_\ast$, 
(17.81) and the fact that we still use (17.37) say that $h_t(x) = h_1(x)$, 
so $x \in W_{k,1}$ and the desired control comes from (17.58) and (17.59).
Otherwise, observe that $|g_t(x)-x| \leq 2 (1+|f|_{lip}) \varepsilon_\ast$
by (17.89), and $|h_t(x)-g_t(x)| \leq C \varepsilon_0$
because we use the formulas (17.37) and (17.38) and by the proof of
(17.39). On the other hand, $\dist(x,\R^n \sm U) \geq \delta_0/2$ because
$x \in U \sm U_{ext}$; so $\dist(h_t(x),\R^n \sm U) \geq \delta_0/4$ 
if $\varepsilon_0$ and $\varepsilon_\ast$ are small enough, and (17.93) follows.

\ms
$\dagger$
Let us do a similar construction under the Lipschitz assumption.
We still have (17.79) (by (16.8), (16.28), and (16.29); also see eight
lines below (16.30) for the definition on $U \sm U_{int}$)
and (17.80) (for the same reason). As we did in (17.60),
we first set 
$$
h_t(x) = x \ \hbox{ when $x\in U_{ext}$ and } 1 \leq t \leq 2.
\leqno (17.94)
$$
This way, it will be enough to define $h_t$ on the set $E_k \cap H_1$,
where we shall find it convenient to first define mappings $\wt g_t$, 
$1 \leq t \leq 2$. We shall also check that the two
formulas coincide on $E_k \cap H_1 \cap U_{ext}$.

Recall that when $x\in E_k \cap H_1$, 
the $h_t(x)$, $0 \leq t \leq 1$, were defined by (17.64),
in terms of functions $\wt h_t(x)$, themselves defined
by (17.62) and (17.63). 
In particular, we had first observed that
the $\wt g_t(x)$, and in particular $\wt f(x) = \wt g_1(x)$ were
well defined (in terms of $g_t(x)$ and $f(x)$). Here we proceed 
similarly. First we set
$$
\wt g_2(x) = \wt g_1(x) = \wt g(x) 
\ \hbox{ when } \dist(x,\R^n \sm W_f) \geq 2 \varepsilon_\ast
\leqno (17.95)
$$
where the fact that $\wt g_1(x) = \wt g(x)$ comes 
from (16.29), or (when $x\in \R^n \sm V_1^+$)
from the definition above (16.30).
When $\varepsilon_\ast \leq \dist(x,\R^n \sm W_f) \leq 2 \varepsilon_\ast$
(and $x\in E_k \cap H_1$), set
$$\eqalign{
\wt g_2(x) &= {\dist(x,\R^n \sm W_f) -\varepsilon_\ast \over 
\varepsilon_\ast} \, \wt g_1(x) 
+ {2\varepsilon_\ast -\dist(x,\R^n \sm W_f) \over \varepsilon_\ast} 
\, \psi(\lambda x)
\cr&
= {\dist(x,\R^n \sm W_f) -\varepsilon_\ast \over 
\varepsilon_\ast} \, \wt f(x) 
+ {2\varepsilon_\ast -\dist(x,\R^n \sm W_f) \over \varepsilon_\ast} 
\, \psi(\lambda x),
}\leqno (17.96)
$$
where the identity comes from (17.80) through the usual
change of variable. Finally set 
$$
\wt g_2(x) = \psi(\lambda x) \ \hbox{ when } \dist(x,\R^n \sm W_f) < 
\varepsilon_\ast.
\leqno (17.97)
$$
This defines the mapping $\wt g_2$ on $E_k \cap H_1$, and now we 
define the $\wt g_t(x)$, $1 \leq t \leq 2$, by
$$
\wt g_t(x) = (2-t) \wt g_1(x) + (t-1) \wt g_2(x)
= (2-t) \wt f(x) + (t-1) \wt g_2(x).
\leqno (17.98)
$$
We proceed as near (17.62), define
$$
\wt s_t(x) = C \Min(\lambda \varepsilon_0, 
|\wt g_t(x)- \psi(\lambda x)|)
\leqno (17.99) 
$$
also for $1 < t \leq 2$, and where $C$ is still as in (17.20), 
and then 
$$
\wt h_t(x) = \Pi(\wt g_t(x),\wt s_t(x))
\ \hbox{ for $x\in E_k \cap H_1$ and } 1 < t \leq 2,
\leqno (17.100) 
$$
as in (17.63). As usual, we want to set
$$
g_t(x) = \lambda^{-1}\psi^{-1}(\wt g_t(x))
\ \hbox{ and } \ 
h_t(x) = \lambda^{-1}\psi^{-1}(\wt h_t(x))
\leqno (17.101)
$$
for $x\in E_k \cap H_1$ and $1 \leq t \leq 2$,
but we shall first need to check that this makes sense.

When $\dist(x,\R^n \sm W_f) \geq 2 \varepsilon_\ast$,
(17.95) implies that $\wt g_2(x) \in B(0,1)$,
so $g_t(x) = \lambda^{-1}\psi^{-1}(\wt g_t(x))$ is well
defined, and $g_t(x)=g(x)$. Since we used the same formula
(17.63) to define $\wt h_1(x)$, we also get that
$\wt h_t(x) = \wt h_1(x)$, so $h_t(x) = \lambda^{-1}\psi^{-1}(\wt 
h_t(x))$ is well defined, and we also get that
$$
h_t(x) = h_1(x) \ \hbox{ for } 1 \leq t \leq 2
\leqno (17.102)
$$
when $x\in E_k \cap H_1$ is such that
$\dist(x,\R^n \sm W_f) \geq 2 \varepsilon_\ast$.

So suppose now that $\dist(x,\R^n \sm W_f) \leq 2\varepsilon_\ast$.
Let check that then
$$
|\wt h_t(x)-\psi(\lambda x)| 
\leq C \lambda \Lambda (1+|f|_{lip}) \varepsilon_\ast
\leqno (17.103) 
$$
for $t \geq 1$. Notice that
$$\eqalign{
|\wt g_t(x)-\psi(\lambda x)| 
&\leq \Max(|\wt g_1(x)-\psi(\lambda x)|,|\wt g_2(x)-\psi(\lambda x)|)
\cr& = \Max(|\wt f(x)-\psi(\lambda x)|,|\wt g_2(x)-\psi(\lambda x)|) 
\cr&
= |\wt f(x)-\psi(\lambda x)|
\leq \lambda \Lambda |f(x)-x|
\cr&
\leq \lambda \Lambda(1+|f|_{lip}) \dist(x,\R^n \sm W_f)
\leq 2 \lambda \Lambda (1+|f|_{lip}) \varepsilon_\ast
}\leqno (17.104)
$$
by (17.98), (17.96) or (17.97), and the fact that $f(x)-x$ 
is $(1+|f|_{lip})$-Lipschitz and vanishes on $\R^n \sm W_f$. 
Also,
$$
\wt s_t(x) \leq C |\wt g_t(x)- \psi(\lambda x)|
\leqno (17.105) 
$$
by (17.99), so
$$
|\wt h_t(x)-\wt g_t(x)| \leq C \wt s_t(x)
\leq C |\wt g_t(x)- \psi(\lambda x)|
\leq C \lambda \Lambda (1+|f|_{lip}) \varepsilon_\ast
\leqno (17.106) 
$$
by (17.100), (17.19), and (17.104). Now (17.103) follows from
(17.104) and (17.106).

We are now ready to prove that $\wt g_t(x)$ and 
$\wt h_t(x)$ lie in $B(0,1)$ when $x\in E_k \cap H_1$.
We already treated the case when 
$\dist(x,\R^n \sm W_f) \leq 2 \varepsilon_\ast$
(see (17.102)); in the other case, notice that
$\dist(\psi(\lambda x); \R^n \sm B(0,1)) \geq \Lambda^{-1} \lambda
\dist(x,\R^n \sm U) \geq \Lambda^{-1} \lambda \delta_0/3$
by (17.61), hence $\wt h_t(x)$ lie in $B(0,1)$, by (17.103),
and similarly for $\wt g_t(x)$, by (17.104).
So the definitions in (17.101) make sense.

When $x\in E_k \cap H_1 \cap U_{ext}$, (17.101)
yields the same result as (17.94), because (17.94)
says that $h_t(x) = x$, while
$h_1(x) = g_1(x) = f(x) = x$ by (17.44) and its proof by (17.42),
then $\wt g_2(x) = \psi(\lambda x)$ by (17.95)-(17.97),
then $\wt g_t(x) = \psi(\lambda x)$ by (17.98), and finally
$g_t(x) = x$ by (17.101).

We continue with the verifications that follow the definition
of the $h_t$. We still have (17.85), by (17.94) or (17.102).
Instead of (17.86), let us just check that
$$
h_2(x) = x \ \hbox{ when $x\in E_k$ and } 
\dist(x,\R^n \sm W_f) < \varepsilon_\ast.
\leqno (17.107) 
$$
When $x\in U_{ext}$, $h_2(x) = x$ by (17.94). 
When $x\in E_k \cap H_1$, (17.97) yields 
$\wt g_2(x) = \psi(\lambda x)$; then (17.98) yields
$\wt g_t(x) = \psi(\lambda x)$ for $t=2$
(yes, there is a double definition of $\wt g_2(x)$,
but the notation is acceptable because they give the same result),
and finally $\wt s_2(x) = 0$ because $\wt g_2(x) - \psi(\lambda x) = 0$
and hence $\wt h_2(x) = \wt g_2(x) = \psi(\lambda x)$ by (17.100);
(17.107) follows.

Our next verification is (17.87), the fact that Lemma 17.40
still holds for our extended family. As before, only (1.5)-(1.7)
need to be checked, because (1.4) and (1.8) can be seen from the 
definition (this is why we made sure to have an overlap when
we used two definitions, on $U_{ext}$ and on $E_k \cap H_1$).
For these verifications, we already know the desired result
for $0 \leq t \leq 1$, and we only need to worry when 
$h_t(x) \neq h_t(1)$ for some $t > 1$. Hence we may restrict to
$x\in E_k \cap H_1$ such that $\dist(x,\R^n \sm W_f) \leq 
2\varepsilon_\ast$ (see (17.95) and (17.98), and compare (17.99)-(17.101)
to (17.62)-(17.64)). 

For (1.5), we suppose in addition $x\in E_k \sm 
B(X_0,R_0+\varepsilon_0)$ and we want to know that $h_t(x)=x$
for $t>1$. But we already know that $g_1(x) = x$, so
(17.96) yields $\wt g_2(x) = \psi( \lambda x)$,
hence $\wt g_t(x) = \psi( \lambda x)$ by (17.98), 
$\wt s_t(x) = 0$ by (17.99), and finally
$h_t(x) = g_t(x) = x$ for $t> 1$, by (17.100) and (17.101).

For (1.6) we suppose that $x\in E_k \cap B(X_0,R_0+\varepsilon_0)$
and we want to prove that $h_t(x) \in B''$ for $t > 1$.
Since we may assume that $x\in E_k \cap H_1$ and 
$\dist(x,\R^n \sm W_f) \leq 2\varepsilon_\ast$
we can use (17.103), which by (17.101) implies that
$$
|h_t(x)-x| \leq C \Lambda^2 (1+|f|_{lip}) \varepsilon_\ast,
\leqno (17.108) 
$$
from which we easily deduce that $h_t(x) \in B''$, because 
$x\in B(X_0,R_0+\varepsilon_0)$ and if $\varepsilon_\ast$ is small enough.

We are left with (1.7). We are given $x\in E_k \cap L_j$, and we want
to check that $h_t(x) \in L_j$ for $t \geq 1$. Again we can assume that 
$x\in E_k \cap H_1$ and $\dist(x,\R^n \sm W_f) \leq 2\varepsilon_\ast$, 
so (17.103) holds, and hence
$$
\dist(\wt g_t(x), \wt L_j) \leq |\wt g_t(x)-\psi(\lambda x)| 
\leq C \lambda \Lambda (1+|f|_{lip}) \varepsilon_\ast
\leq \lambda \varepsilon_0
\leqno (17.109)
$$
because $x\in L_j$ and if $\varepsilon_\ast$ is small enough
compared to $\varepsilon_0$. Then $\wt s_t(x) = C |\wt 
g_t(x)-\psi(\lambda x)|$ in (17.99), and since there is a face $\wt F$
of $\wt L_j$ such that $\dist(\wt g_t(x),\wt F) \leq |\wt g_t(x)-\psi(\lambda x)|
= C^{-1}\wt s_t(x)$ (by (17.109)), we get that 
$\Pi(\wt g_t(x), \wt s_t(x)) \in \wt F$ by (17.20) and our choice of 
$C$ in (17.99); Thus $\wt h_t(x)= \Pi(\wt g_t(x), \wt s_t(x)) \in \wt F \i \wt L_j$ 
by (17.100) and $h_t(x) \in L_j$, as needed. This completes
our verification of (17.87) in the Lipschitz context.

We also want to generalize (17.91)-(17.93). In fact, (17.91) still holds, by
our definition (17.94), (17.92) is a definition, and the
proof of (17.93) also goes the same way: we just need to consider 
points $x\in E_k \cap H_1$ such that 
$\dist(x,\R^n \sm W_f) \leq 2\varepsilon_\ast$, hence for which (17.108)
holds. But $\dist(x,\R^n \sm U) \geq \delta_0/3$ because $x\in H_1$,
so $\dist(h_t(x),\R^n \sm U) \geq \delta_0/4$, by (17.108),
and (17.93) follows.

\vfill\eject 
\msi{\bf 18. The final accounting and the proof of Theorem 10.8
in the rigid case} 
\ms 

In the previous sections, we managed to construct deformations
of $E$ and the $E_k$; we are especially interested in the last
one, the family $h_{2t}$, $0 \leq t \leq 1$.
By (17.87), or the corresponding verification in the Lipschitz
case (see above (17.108), Lemma 17.40 holds for the $h_{2t}$, 
$0 \leq t \leq 1$. Also, (17.93) says that the analogue of (2.4) 
for $E_k$ holds for $k$ large enough. Thus we can apply Definition 2.3
and the quasiminimality of $E_k$; we get that for $k$ large,
$$
\H^d(E_k \cap W) \leq M \H^d(h_2(E_k \cap W)) + (R'')^d h,
\leqno (18.1)
$$
where $R''$ is as in (17.41) and
$$
W = \big\{ y \in U \, ; h_2(x) \neq x \big\}. 
\leqno (18.2)
$$
Most of this section will consist in estimating the two
sides of (18.1) (and especially the right-hand side).
We shall need to cut $U$ into small pieces, and we shall
start with the least important ones. We shall also try to
treat the rigid and Lipschitz assumptions simultaneously
when this is possible, but some estimates for the Lipschitz
cases will be done in the next section.

Return to the definition of $h_2$. 
First recall from (17.86) (or (17.107) in the Lipschitz case) 
that $h_2(x)=x$ when $x\in E_k$ is such that
$\dist(x,\R^n \sm W_f) < \varepsilon_\ast$; hence
$$
E_k \cap W \i \big\{ x\in W_f \, ; \, 
\dist(x,\R^n \sm W_f) \geq \varepsilon_\ast \big\} \i W_f. 
\leqno (18.3)
$$
{\bf The exterior skin.}
In the rigid case and on the set 
$$
A_\ast = \big\{ x\in W \, ; \,
\varepsilon_\ast \leq \dist(x,\R^n \sm W_f) \leq 2 \varepsilon_\ast 
\big\},
\leqno (18.4)
$$
we defined $g_2$ by (17.82), and then decided to use (17.37) as before
(see below (17.84)). That is
$$
h_2(x) = \Pi(g_2(x),s_2(x)),
\leqno (18.5)
$$
where $s_2(x)$ is still defined as in (17.38), and 
$g_2(x)$ is defined by (17.82). We shall need to know that
$$
h_2 \hbox{ is $C$-Lipschitz on } A_\ast.
\leqno (18.6)
$$
Since $\Pi$ is $C$-Lipschitz (by (17.21)), (15.8) and 
(17.38) say that it is enough to show that $g_2$ is 
$3(1+|f|_{lip})$-Lipschitz on $A_\ast$. And indeed, 
for $x,y\in A_\ast$,
$$\leqalignno{
|g_2(x)-g_2(y)|&= 
\Big|{\dist(x,\R^n \sm W_f) -\varepsilon_\ast \over 
\varepsilon_\ast} \, f(x) 
+{2\varepsilon_\ast -\dist(x,\R^n \sm W_f) \over \varepsilon_\ast} 
\, x
\cr& \hskip 1.8cm
- {\dist(y,\R^n \sm W_f) -\varepsilon_\ast \over 
\varepsilon_\ast} \, f(y) 
- {2\varepsilon_\ast -\dist(y,\R^n \sm W_f) \over \varepsilon_\ast} 
\, y \Big|
\cr&
\leq {\dist(x,\R^n \sm W_f) -\varepsilon_\ast \over 
\varepsilon_\ast} \, |f(x)-f(y)|
+{2\varepsilon_\ast -\dist(x,\R^n \sm W_f) \over \varepsilon_\ast}
\, |x-y|
&(18.7)
\cr& \hskip 1.8cm
+ {|\dist(x,\R^n \sm W_f) - \dist(y,\R^n \sm W_f)|\over \varepsilon_\ast} \,|f(y)-y|
\cr&
\leq |x-y| + |f(x)-f(y)| + |x-y|\, {|f(y)-y| \over \varepsilon_\ast}
\leq 3(1+|f|_{lip})\,|x-y|
}
$$
by (17.82), and because $y\in A_\ast$ and $f(y)-y$ is a $(1+|f|_{lip})$-Lipschitz
function that vanishes on $\R^n \sm W_f$. So (18.6) holds.

$\dagger$ Under the Lipschitz assumption, we start with
an analogue of (18.6). We first work on $E_k \cap H_1$
(where $H_1$ is as in (17.61)) and check that
$$
\wt g_2 \hbox{ is $C(\Lambda,f) \lambda$-Lipschitz on }
A_\ast \cap E_k \cap H_1.
\leqno (18.8)
$$
On $A_\ast \cap E_k \cap H_1$, we defined $\wt g_2$ by (17.96)
(see (18.4)). Recall that we restricted to $E_k \cap H_1$ because 
we were able to define $\wt f(x) = \psi(\lambda f(x))$ 
for $x\in E_k \cap H_1$, and compute with it. 
Then we can follow the same proof of (18.7) and get (18.8).

Then we observe that $\wt h_2$ was defined by (17.99) and (17.100),
so it is also $C(\Lambda,f)\lambda$-Lipschitz on $E_k \cap H_1 \cap A_\ast$.
Finally, 
$$
h_2 \hbox{ is $C(\Lambda,f)$-Lipschitz on }
A_\ast \cap E_k \cap H_1,
\leqno (18.9)
$$
because of (17.101). 

Recall from (17.94) that $h_2(x)=x$ on $U_{ext}$,
and that $H_1$ and $U_{ext}$ cover $\R^n$ (even, with an overlap),
by (17.61) and (11.2). So (18.9) will be good enough.
For instance, 
$$\eqalign{
\H^d(h_2(E_k \cap A_\ast))
&\leq \H^d(h_2(E_k \cap A_\ast\cap H_1)) 
+ \H^d(h_2(E_k \cap A_\ast \sm H_1))
\cr&
\leq C(\Lambda,f) \H^d(E_k \cap H_1 \cap A_\ast))
+ \H^d(E_k \cap A_\ast \sm H_1)
\cr&
\leq C(\Lambda,f) \H^d(E_k \cap A_\ast))
}\leqno (18.10)
$$
by (18.9) and the trivial estimate on $U_{ext}$. $\dagger$

In the rigid case, we also have the conclusion of (18.10), directly by (18.6).
Thus the contribution of $A_\ast$ to the 
right-hand side of (18.1) is easily controlled in both cases.
Next return to the general case, and recall that 
$$
\hbox{$W_f$ is compactly contained in $U$,}
\leqno (18.11)
$$
by (17.78) or (11.19), (11.3), and (11.2). 
Since $A_\ast \i W_f$ by (18.4), it is a compact subset of $U$. 
Then we can apply (10.14) to $A_\ast$ and get that
$$
\limsup_{k \to +\infty} \H^d(E_k \cap A_\ast))
\leq C_M \H^d(E \cap A_\ast)).
\leqno (18.12)
$$
Finally
$$
\limsup_{\varepsilon_\ast \to 0}  \H^d(E \cap A_\ast)
\leq \limsup_{\varepsilon_\ast \to 0} 
\H^d\big(\big\{ x \in E \, ; \,
0 < \dist(x,\R^n \sm W_f) \leq 2 \varepsilon_\ast \big\}\big)
= 0
\leqno (18.13)
$$
(the monotone intersection of these sets is empty), 
so we deduce from (18.10)-(18.13) that
$$
\H^d(h_2(E_k \cap A_\ast)) \leq \eta + C \H^d(E \cap A_\ast))
\leq 2 \eta
\leqno (18.14)
$$
for $k$ large, and provided that we choose $\varepsilon_\ast$ small 
enough (depending on the usual quantities).

Observe that by (18.3) and (18.4),
$$
E_k \cap W\sm A_\ast \i W_\ast := \big\{ x\in W_f \, ; \dist(x,\R^n \sm W_f) 
> 2 \varepsilon_\ast \big\},
\leqno (18.15)
$$
and so we shall now concentrate on $E_k \cap W_\ast$.
On this set, and in the rigid case, (17.85), (17.37), and (16.13)
say that
$$
h_2(x) = h_1(x) = \Pi(g_1(x),s_1(x))
= \Pi(g(x),s_1(x)),
\leqno (18.16)
$$
where by (17.38) and (16.13)
$$
s_1(x) = C \Min(\varepsilon_0, |g_1(x)-x|)
= C \Min(\varepsilon_0, |g(x)-x|).
\leqno (18.17)
$$

$\dagger$ In the Lipschitz case, either $x\in E_k \cap U_{ext}$
and (17.94) says that $h_2(x)=x$, or else $x\in E_k \cap H_1$ and
(since $x\in W_\ast$) (17.101) says that
$h_2(x) = \lambda^{-1} \psi^{-1}(\wt h_2(x))$,
where $\wt h_2(x) = \Pi(\wt g_2(x),\wt s_2(x))$
by (17.100). Since (17.95) says that
$\wt g_2(x) = \wt g_1(x) = \wt g(x)$, the definition (17.99)
yields
$$
\wt s_2(x) 
= C \Min(\lambda \varepsilon_0, |\wt g_2(x)- \psi(\lambda x)|)
= \wt s_1(x) 
= C \Min(\lambda \varepsilon_0, |\wt g(x)- \psi(\lambda x)|)
\leqno (18.18) 
$$
and hence
$$
\wt h_2(x) = \Pi(\wt g_2(x),\wt s_2(x))
= \Pi(\wt g_1(x),\wt s_1(x))
= \Pi(\wt g(x),\wt s_1(x))
\leqno (18.19)
$$
for $x\in E_k \cap H_1 \cap W_\ast$
(a good enough  analogue of (18.16) and (18.17)). $\dagger$ 

In both cases, (18.16)-(18.19), (16.10), and (17.21) yield
$$
h_2 \hbox{ is $C$-Lipschitz on } E_k \cap W_\ast
\leqno (18.20) 
$$
with a constant $C$ that depends on $M$ and $|f|_{lip}$ in particular,
and also on $\Lambda$ in the Lipschitz case.

\msi{\bf The part outside the balls.}
Our next small set is
$$
W^1_\ast = W_\ast \sm V_1^+,
\leqno (18.21) 
$$
where $V_1^+$ is as in (16.7). 
Then (18.20) yields
$$
\H^d(h_2(E_k \cap W^1_\ast)) 
\leq C \H^d(E_k \cap W^1_\ast)),
\leqno (18.22)
$$
where we don't mention the dependence on $\Lambda$ and $|f|_{lip}$.
Again $\overline W^1_\ast$ is a compact subset of $U$,
because $W_\ast \i W_f$ by (18.15), and by (18.11). 
So (10.14) yields
$$
\limsup_{k \to +\infty} \H^d(E_k \cap W^1_\ast))
\leq C_M \H^d(E \cap \overline W^1_\ast))
\leqno (18.23)
$$
and hence
$$
\H^d(h_2(E_k \cap W^1_\ast)) 
\leq \eta + C \H^d(E \cap \overline W^1_\ast))
\leqno (18.24)
$$
for $k$ large. Of course it will be interesting to control
$\H^d(E \cap \overline W^1_\ast))$, and we shall do this more
easily after the next step. 

\msi{\bf The small rings.} Next we want to control the contribution
of the small rings.

\ms\proclaim Lemma 18.25. Set
$$
R^1 = \bigcup_{j\in J_1} [\overline B_j \sm aB_j] \,,
\ \ 
R^2 = \bigcup_{j\in J_2} [\overline B_j \sm aB_j] \,,
\hbox{ and } \ 
R^3 = \hskip -0.1 cm
\bigcup_{j\in J_3 \, ; \, x \in Z(y_j)} 
[B_{j,x} \sm B_{j,x}^-] \,,
\leqno (18.26)
$$
where $B_{j,x}$ and $B_{j,x}^-$ are as in 
(15.19) and (15.17). Then
$$
\H^d(E \cap [R^1 \cup R^2 \cup R^3])
\leq C(f,\gamma) (1-a), 
\leqno (18.27)
$$
where $C(f,\gamma)$ depends on $|f|_{lip}$,
$\H^d(E \cap W_f)$, and $\gamma$.

\ms
Indeed, we know from (13.24) and the definition (11.20) that
$$
\H^d(E \cap R^1) \leq C (1-a) \H^d(E \cap W_f).
\leqno (18.28)
$$
Similarly, (14.20) says that 
$$\eqalign{
\H^d(E \cap R^2)
&\leq \sum_{j\in J_2} \H^d(E \cap \overline B_j \sm aB_j)
\leq C (1-a) \H^d(E \cap W_f).
}\leqno (18.29)
$$
We are thus left with the $B_{j,x} \sm B_{j,x}^{-}$.
First fix $j\in J_3$ and $x\in Z(y_j)$, recall the definitions
(15.16)-(15.19), and cover $(2-a)E(x) \sm aE(x)$ 
by balls $A_k$ of radius $(1-a)r_j$. 
By (15.19) and (15.20) (but we could also use the definition
(15.16) of $E(x)$), the ellipsoid $(2-a)E(x)$
is contained $B(x,{3 \over 2} \gamma^{-1} r_j)$; by the definition
(15.16), it is also contained in the $d$-plane $P_x$, and so we need less than
$C \gamma^{-d} (1-a)^{1-d}$ balls $A_k$ to cover $(2-a)E(x) \sm aE(x)$. 
The slightly larger balls $2A_k$ cover $B_{j,x} \sm B_{j,x}^{-}$, so
$$
\H^d(E \cap B_{j,x} \sm B_{j,x}^{-})
\leq \sum_{k} \H^d(E \cap 2A_k).
\leqno (18.30)
$$
Let us apply Proposition 4.1 (the local Ahlfors-regularity of $E$)
to each $2A_k$. This is allowed, because $4A_k \i W_f \i U$,
since
$$
r_j \leq {\gamma \over 4} \delta_6 \leq {\gamma \over 40} \delta_1
\leq {\gamma \over 40}\dist(x,\R^n \sm W_f)
\leqno (18.31)
$$
by (15.15), (12.7), (11.22), and because
$x \in X_9 \i X_1$ (see (15.1)). We get that
$$\eqalign{
\H^d(E \cap B_{j,x} \sm B_{j,x}^{-})
&\leq \sum_{k} \H^d(E \cap 2A_k) 
\leq C \sum_k (1-a)^d r_j^d
\leq C (1-a) \gamma^{-d} r_j^d 
\cr&
\leq C (1-a) \gamma^{-d}  |f|_{lip}^d \H^d(E \cap B(x,|f|_{lip}^{-1}r_j))
}\leqno (18.32)
$$
by (18.30) and Proposition 4.1, applied a second time but in the other
direction. The reader should not worry about $|f|_{lip}^{-1}$ being too
large: we know that $|f|_{lip} \geq 1$ because $f(z)=z$ near $\infty$.
We claim that 
$$
\hbox{ the balls $B(x,|f|_{lip}^{-1}r_j)$,
$j\in J_3$ and $x\in Z(y_j)$, are all disjoint.} 
\leqno (18.33)
$$
For different $j$, this is because 
$f(B(x,|f|_{lip}^{-1}r_j)) \i B(f(x),r_j) = D_j$,
and the $D_j$ are disjoint by (15.13).
For the same $j$ and different $x\in Z(y_j)$,
this follows from (15.6) and the fact that $y_j \in Y_9$
by (15.12), because we can safely assume that $\gamma \leq 1$. 
Now
$$\eqalign{
H^d(E \cap R^3)
&\leq \sum_{j\in J_3} \sum_{z\in Z(y_j)}
\H^d(E \cap B_{j,z} \sm B_{j,z}^{-})
\cr&
\leq C (1-a) \gamma^{-d}  |f|_{lip}^d \sum_{j\in J_3} \sum_{z\in Z(y_j)}
\H^d(E \cap B(x,|f|_{lip}^{-1}r_j))
\cr&
\leq C (1-a) \gamma^{-d} |f|_{lip}^d \H^d(E \cap W_f)
}\leqno (18.34)
$$
by (18.32), (18.33), and because $B(x,|f|_{lip}^{-1}r_j) \i W_f$
by (18.31). Lemma 18.25 follows.
\qed

\ms
By (18.20),
$$
\H^d(h_2(E_k \cap W_\ast \cap [R^1 \cup R^2 \cup R^3])) 
\leq C \H^d(E_k \cap [R^1 \cup R^2 \cup R^3])).
\leqno (18.35)
$$
But $R^1 \cup R^2 \cup R^3$ is compact because
the sets $J_1$, $J_2$, and $J_3$ are finite, and it is contained
in $U$ by construction of the $B_j$ and $B_{j,x}$, so
(10.14) and (18.27) say that
$$
\H^d(E_k \cap [R^1 \cup R^2 \cup R^3]))
\leq C \H^d(E \cap [R^1 \cup R^2 \cup R^3]))+\eta
\leq C(f,\gamma) (1-a) + C \eta
\leqno (18.36)
$$
for $k$ large. Altogether,
$$
\H^d(h_2(E_k \cap W_\ast \cap [R^1 \cup R^2 \cup R^3])) 
\leq C(f,\gamma) (1-a) + C \eta.
\leqno (18.37)
$$

\msi{\bf Return to the part outside the balls.}
We still want to estimate $\H^d(E \cap \overline W^1_\ast))$, 
but we shall even consider the set $E \cap  W^2$, where
$$
W^2 = W_f \sm {\rm int}(V_1^+).
\leqno (18.38)
$$
Let us first check that $\overline W^1_\ast \i W^2$. Let
$x\in \overline W^1_\ast$ be given; then $x$ is the 
limit of some sequence $\{ x_k \}$ in $W^1_\ast = W_\ast \sm V_1^+$
(see (18.21)).
By (18.15), $\dist(x_k,\R^n \sm W_f)  > 2 \varepsilon_\ast$ for all 
$k$, so $x\in W_f$. But also $x_k$ lies out of $V_1^+$,
so $x \notin {\rm int}(V_1^+)$. Thus $\overline W^1_\ast \i W^2$.

Now let $x\in E \cap W^2$ be given; we want
to show that it lies in one of many small sets.
Observe that $x$ lies in $X_0 = E \cap W_f$, by (11.20).
First assume that $x\in \overline B_j$ for some $j\in J_1 \cup J_2$.
Notice that ${1+a \over 2} B_j$ is contained in ${\rm int}(V_1^+)$ by (16.7), 
so it does not meet $W^2$. Hence 
$x\in \overline B_j \sm {1+a \over 2} B_j \i R^1\cup R^2$, 
and the corresponding subset of $E \cap W^2$ is small, by (18.27).

So we may assume that 
$x\in X_0 \sm \bigcup_{j\in J_1 \cup J_2} \overline B_j$.
The case when $x\notin X_9$ is controlled by (14.24), which says that 
$\H^d\big(X_0 \sm \big[X_9 \cup \big(\bigcup_{j\in J_1 \cup J_2} \overline B_j 
\big)\big]\big) \leq 7 \eta$, 
so we may even assume that $x\in X_9$.

The case when $x\in X_9 \sm X_{11}$ is covered by (15.8) and (15.11),
so we may assume that $x\in X_{11}$. In addition, (15.14) allows us to
assume that $x \in f^{-1}\big( \bigcup_{j\in J_3} \overline D_j \big)$,
and by (15.36) $x$ lies in $B_{j,z}$ for some $j\in J_3$
and $z\in Z(y_j)$. At the same time, 
$B_{j,z}^+$ is contained in $V_1^+$ by (16.7), so
its interior does not meet $W^2$.
The interior contains $B_{j,z}^-$ (see the definitions 
(15.17)-(15.19)), hence $x\in B_{j,z}\sm B_{j,z}^- \i R^3$
(compare with (18.26)).
The corresponding set is again controlled by (18.27), and altogether
$$
\H^d(E \cap W^2)) \leq C(f,\gamma) (1-a) + C \eta.
\leqno (18.39)
$$
Since $\overline W^1_\ast \i W^2$, we get that
$$
\H^d(h_2(E_k \cap W^1_\ast)) 
\leq \eta + C \H^d(E \cap \overline W^1_\ast))
\leq C(f,\gamma) (1-a) + C \eta
\leqno (18.40)
$$
for $k$ large, by (18.24) and (18.39).

\ms
We are left with the set $V_1^+$. Recall from (16.15)
that
$$
\dist(x,X_1) \leq \delta_6
\hbox{ and } \dist(x,\R^n \sm W_f) > \delta_1/2
\ \hbox{ for $x\in V_1^+$.}
\leqno (18.41)
$$

In addition, recall from (11.2) that 
$\delta_0 = \dist(\wh W,\R^n \sm U)$; 
since $X_1 \i X_0 \i \wh W$ by (11.20), we get that
$\dist(z,\R^n \sm U) \geq \delta_0$ for $z\in X_1$, and hence
$$
\dist(x,\R^n \sm U) \geq {2\delta_0 \over 3}
\ \hbox{ for $x\in V_1^+$,}
\leqno (18.42)
$$
because $\delta_6 < \delta_0/3$ by (12.7).
The definition (17.61) then immediately yields
$$
V_1^+  \i  H_1.
\leqno (18.43)
$$ 
Next we check that
$$
|g(x)-x| \geq {\delta_5 \over 2}
\ \hbox{ for $x\in V_1^+$.}
\leqno (18.44)
$$
Use (18.41) to find $z\in X_1$ such that $|z-x| \leq \delta_6$,
and notice that $|f(z)-z| \geq \delta_5$ by the definition (12.5).
Then
$$\eqalign{
|g(x)-x| &\geq |f(x)-x| - |g(x)-f(x)|
\geq |f(x)-x| - ||f-g||_\infty
\cr&
\geq|f(x)-x| - 4\Lambda^2(1+|f|_{lip}) \delta_6
\cr&
\geq |f(z)-z| - (1+|f|_{lip})|z-x| - 4\Lambda^2(1+|f|_{lip}) \delta_6
\cr&
\geq \delta_5  - (1+4\Lambda^2) (1+|f|_{lip}) \delta_6
\geq {\delta_5 \over 2}
}\leqno (18.45)
$$
by (16.11) and (12.7). So (18.44) holds. We want to use this to show
that for $k$ large, 
$$
|h_2(x)-x| \geq {\delta_5 \over 4} > 0
\ \hbox{ for } x\in E_k \cap V_1^+,
\leqno (18.46)
$$
so we want to estimate $|h_2(x)-g(x)|$.

We start in the rigid case. 
If $\varepsilon_0$ is chosen small enough
(compared to $\delta_5$), we get that 
$s_1(x) = C \Min(\varepsilon_0, |g(x)-x|) = C \varepsilon_0$
(by (18.17) and (18.45)). In this case, (18.16) simplifies to
$$
h_2(x) = \Pi(g(x),C\varepsilon_0)
\ \hbox{ for } x\in V_1^+
\leqno (18.47)
$$
and hence, by (17.19),
$$
|h_2(x)-g(x)| 
\leq C\varepsilon_0;
\leqno (18.48)
$$
in this case (18.46) follows from (18.44).

$\dagger$ Similarly, under the Lipschitz assumption, 
let $x\in E_k \cap V_1^+$ be given. First observe that
$x\in E_k \cap H_1$, by (18.43). In addition,
$x\in W_\ast$ if $\varepsilon_\ast < \delta_1/4$
(compare the definition (18.15) with (18.41)),
so (18.18) and (18.19) hold.

Recall from the remark below (16.8) that since $x\in V_1^+$,
$g(x)$ was defined in terms of a function $\wt g$, through
the relation $\wt g(x) = \psi(\lambda g(x))$; thus (18.18)
yields
$$
\wt s_1(x) = C \Min(\lambda \varepsilon_0, |\wt g(x)- \psi(\lambda x)|)
= C \Min(\lambda \varepsilon_0, |\psi(\lambda g(x))- \psi(\lambda x)|)
= C\lambda \varepsilon_0
\leqno (18.49)
$$ 
because $|\psi(\lambda g(x))- \psi(\lambda x)|
\geq \lambda \Lambda^{-1} |g(x)-x| \geq \lambda \varepsilon_0$
if $\varepsilon_0$ is small enough, by (18.44), and
(18.19) simplifies to
$$
\wt h_2(x) = \Pi(\wt g(x),\wt s_1(x)) 
= \Pi(\wt g(x),C\lambda \varepsilon_0).
\leqno (18.50)
$$
Since (17.101) says that 
$$
h_2(x) = \lambda^{-1} \psi^{-1}(\wt h_2(x)),
\leqno (18.51)
$$
we deduce from (18.50) and (17.19) that
$$
|h_2(x)-g(x)| 
= |\lambda^{-1}\psi^{-1}(\wt h_2(x))- \lambda^{-1}\psi^{-1}(\wt g(x))| 
\leq \lambda^{-1}\Lambda |\wt h_2(x)-\wt g(x)|
\leq C\varepsilon_0
\leqno (18.52)
$$
where we we don't care that $C$ depends on $\Lambda$.
That is, (18.48) still holds under the Lipschitz assumption,
and (18.46) follows from (18.44) as before.
$\dagger$

Return to the general case; (18.46) implies that 
$$
E_k \cap V_1^+ \i W,
\leqno (18.53)
$$
where $W = \big\{ y \in U \, ; h_2(x) \neq x \big\}$ 
is as in (18.2).

\msi{\bf The balls $B_j$, $j\in J_1 \cup J_2$.}
Next we estimate the contribution of $E_k \cap V_1^+$
to the right-hand side of (18.1). We start with 
$\bigcup_{j\in J_1} {1+a \over 2} B_j \i \bigcup_{j\in J_1} B_j$.
Observe that the set $R$ in (13.3) and (13.23) is contained 
in the more recent $R^1$ of (18.26), so (13.23) and (11.20) 
say that
$$
\H^d\Big(g\Big(\bigcup_{j\in J_1} B_j \sm R^1\Big)\Big)
\leq
\H^d\Big(g\Big(\bigcup_{j\in J_1} B_j \sm R\Big)\Big)
\leq C(\alpha,f) N^{-1} \H^d(E\cap W_f).
\leqno (18.54)
$$
In the rigid case, (18.47) and (17.21) imply that 
on $V_1^+$, $h_2(x)$ is a $C$-Lipschitz function of
$g(x)$; hence (18.54) yields
$$\eqalign{
\H^d\Big(h_2\Big(E_k \cap V_1^+ \cap 
\Big(\bigcup_{j\in J_1} B_j \sm R^1\Big)\Big)\Big)
&\leq
C \H^d\Big(g\Big(\bigcup_{j\in J_1} B_j \sm R^1\Big)\Big)
\cr&\leq C(\alpha,f) N^{-1} \H^d(E\cap W_f).
}\leqno (18.55)
$$
$\dagger$
Under the Lipschitz assumption, (18.50)
says that on $E_k \cap V_1^+$, 
$\wt h_2(x)$ is a $C$-Lipschitz function of $\wt g(x)$;
hence, since $\wt g(x) = \psi(\lambda g(x))$ is a 
$\lambda \Lambda$-Lipschitz function of $g(x)$
and $h_2(x) = \lambda^{-1} \psi^{-1}(\wt h_2(x))$
is a $\lambda^{-1}\Lambda$-Lipschitz function of $\wt h_2(x)$
(see (18.51)),
$h_2(x)$ is a $C\Lambda^2$-Lipschitz function of $g(x)$.
Thus (18.55) still holds,
only with a larger constant that also depends on $\Lambda$. 
$\dagger$.

We can treat 
$\bigcup_{j\in J_2} {1+a \over 2} B_j \i \bigcup_{j\in J_2} B_j$ 
almost the same way as for $J_1$. Indeed, 
$$\eqalign{
\H^d\Big(g\Big(E_k &\cap \Big(\bigcup_{j\in J_2} 
B_j \sm R^2\Big)\Big)\Big)
\leq \sum_{j\in J_2} \H^d(g(E_k \cap B_j \sm R^2))
\cr&
\leq \sum_{j\in J_2} \H^d(g(E_k \cap aB_j))
\leq C\Lambda^{2d} (1+ \Lambda |f|_{lip})^{d} \gamma \sum_{j\in J_2} r_j^d
\cr&
\leq C(\Lambda, |f|_{lip}) \gamma \, \H^d(E\cap W_f) 
= C \gamma \, \H^d(E\cap W_f)
}\leqno (18.56)
$$
for $k$ large enough, by (18.26), (14.15) or (14.37), 
the fact that $g=g_j$ on $aB_j$ (see (16.5)), and (14.19).
As for the case of $J_1$, $h_2(x)$ is a 
$C$-Lipschitz function of $g(x)$ on $E_k \cap V_1^+$ 
(by (18.47) or (18.50)), so
$$\eqalign{
\H^d\Big(h_2\Big(E_k \cap V_1^+ \cap \Big(\bigcup_{j\in J_2} 
B_j \sm R^2 \Big)\Big)\Big)
&\leq C \H^d\Big(g\Big( E_k \cap \Big(\bigcup_{j\in J_2} 
B_j \sm R^2\Big)\Big)\Big) 
\cr&
\leq C \gamma \, \H^d(E\cap W_f).
}\leqno (18.57)
$$

\msi{\bf The main contribution from the $B_{j,x}^+ \, ,$ $j\in J_3$ and 
$x\in Z(y_j)$.}
Our last piece of $V_1^+$ is the union of the $B_{j,x}^+$, and 
more precisely, of the $B_{j,x}^-$, since the rest is contained
in $R^3$ by (18.26). 
We start with the rigid case.
Let us check that for $j\in J_3$, $x\in Z(y_j)$,
and $z\in B_{j,x}^-$,
$$
h_2(z) = g(z) \in Q_j \cap D_j.
\leqno (18.58)
$$
Recall that $Q_j$ is the common value of the
$d$-planes $A_x(P_x)$, $x\in Z(y_j)$; see above (15.16).
Fix $j\in J_3$, $x\in Z(y_j)$, and $z\in B_{j,x}^-$. 
By (16.6) and (15.29), $g(z) = g_{j,z}(z) = \pi_j(f(z))$,
where $\pi_j$ denotes the orthogonal projection onto $Q_j$.
By (15.22), $f(z) \in {1 + a \over 2} D_j \i D_j$,
so $g(z) \in Q_j \cap D_j$ (recall that $Q_j$ goes through
$y_j = f(x) = A_x(x)$). Now 
$h_2(z) = \Pi(g(z), C\varepsilon_0)$ by (18.47), 
and we still need to check that $h_2(z)=g(z)$.

Denote by $F(y_j)$ the smallest face of our grid 
that contains $y_j = f(x)$, by $W(y_j)$ the affine subspace 
spanned by $F(y_j)$, and by $m$ the dimension of $F(y_j)$ and $W(y_j)$.
By Lemma 12.27, $Q_j \i W(y_j)$.
Also recall from (15.1) and (11.26)
that $x\in X_9 \i X_2 = X_{1,\delta_2}$. So (11.23)-(11.25) say
that $x\in X_{1,\delta_2}(m)$ and 
$$
\dist(y_j,{\cal S}_{m-1}) \geq \delta_2.
\leqno (18.59)
$$
But 
$$
|g(z)-y_j| \leq r_j \leq \delta_6 \leq {\delta_2 \over 10}
\leqno (18.60)
$$
because $g(z) \in D_j= B(y_j,r_j)$ and by (15.15) and (12.7).
We know that $g(z) \in Q_j \i W(y_j)$ and of course $y_j \in W(y_j)$;
hence the line segment $[y_j,g(z)]$ is contained in $W(y_j)$.
In addition, by (18.59) and (18.60), it does not meet 
$\d F(y_j) \i {\cal S}_{m-1}$, and since $y\in F(y_j)$,
we get that $g(z)\in F(y_j)$ too.
We also deduce from (18.59) and (18.60) that
$$
\dist(g(z),{\cal S}_{m-1}) \geq {\delta_2 \over 2}.
\leqno (18.61)
$$

We return to $h_2(z) = \Pi(g(z), C\varepsilon_0)$,
and use the definition of $\Pi$ in (17.23) and (17.24).
Thus
$$
h_2(z)
= \Pi_{0,s_0}\circ\Pi_{1,s_1}\cdots\circ\Pi_{n-1,s_{n-1}}(g(z)),
\leqno (18.62)
$$
where $s_j = (6C)^j (C\varepsilon_0) \leq C' \varepsilon_0$
for $0 \leq j \leq n-1$. For $j \geq m$,
$g(z)\in F(y_j) \i {\cal S}_m$ and so $\Pi_{j,s_j}(g(z)) = g(z)$
by (17.2). For $j < m$, 
$\dist(g(z),{\cal S}_{j}) \geq {\delta_2 \over 2} > 2 s_j$
by (18.61), and now $\Pi_{j,s_j}(g(z)) = g(z)$ by the second part
of (17.2). Altogether, $h_2(z) = g(z)$ and (18.58) follows.

We may now sum over $j$ and $x$:
$$\eqalign{
\H^d\Big(h_2\Big(\bigcup_{j\in J_3} \bigcup_{x\in Z(y_j)} B_{j,x}^+ \sm R^3\Big)\Big)
&\leq
\H^d\Big(h_2\Big(\bigcup_{j\in J_3} \bigcup_{x\in Z(y_j)} B_{j,x}^-)\Big)\Big)
\cr&
\leq \sum_{j\in J_3} \H^d(Q_j \cap D_j)
= \omega_d \sum_{j\in J_3} r_j^d,
}\leqno (18.63)
$$
by (18.26), (18.58), and where $\omega_d$ denotes the $\H^d$-measure 
of the unit ball in $\R^d$. We put everything together and get that
$$\eqalign{
\H^d(h_2(E_k \cap W)) 
&\leq 
C \eta + C(f,\gamma) (1-a) + C(\alpha,f) N^{-1} + C(f) \gamma
+ \omega_d \sum_{j\in J_3} r_j^d,
}\leqno (18.64)
$$
where we no longer write the dependence on $\H^d(E\cap W_f)$,
by (18.14), (18.15), (18.21) and (18.40), (18.37)
(which control everything except 
$E_k \cap V_1^+ \sm [R^1 \cup R^2 \cup R^3]$),
and then (18.55), (18.57), and (18.63) which take care of 
the rest of $E_k \cap V_1^+$.

$\dagger$ Let us say what we get easily under the Lipschitz 
assumption; additional information will be needed, but we shall
only take care of that in the next section. 
We first check that for $j\in J_3$, $x\in Z(y_j)$,
and $z\in E_k \cap B_{j,x}^-$,
$$
\wt h_2(z) = \wt g(z) \in \wt Q_j,
\leqno (18.65)
$$
where $\wt Q_j$ denotes, as above (15.48),
the common image $\wt A_x(P_x)$, $x\in Z(y_j)$.
By (16.6), $g(z) = g_j(z) = g_{j,z}(z)$; by the remark 
below (16.8), we can set $\wt g(z) = \psi(\lambda g(z))$;
by (15.47), $z\in U_{int}$; by (15.53),
$\wt g_j(z) = \psi(\lambda g_j(z))$, and hence
$\wt g(z) = \wt g_j(z)$; by the line below (15.50),
$\wt g(z)= \wt g_{j,x}(z)$; finally, by (15.48),
$\wt g(z) = \wt\pi_j(f(z))$,
where $\wt\pi_j$ denotes the orthogonal projection 
onto $\wt Q_j$. So $\wt g(z) \in \wt Q_j$.
Also, $\wt h_2(z) = \Pi(\wt g(z), C\lambda \varepsilon_0)$
by (18.50).

We still need to check that $\wt h_2(z) = \wt g(z)$.
Denote by $F(y_j)$ the smallest face of our (deformed) grid 
that contains $y_j = f(x)$; thus $\wt F = \psi(\lambda F(y_j))$
is a flat face of the usual dyadic grid, the smallest one
that contains $\wt y_j = \psi(\lambda y_j)$. Denote by
$\wt W$ the affine subspace spanned by $\wt F$, and by $m$ 
the dimension of $\wt F$ and $\wt W$.
By Lemma 12.40, $\wt Q_j \i \wt W$.
As before, $x\in X_9 \i X_2 = X_{1,\delta_2}$, so (11.23)-(11.25) say
that $x\in X_{1,\delta_2}(m)$ and 
$$
\dist(y_j,{\cal S}_{m-1}) \geq \delta_2
\leqno (18.66)
$$
as in (18.59). Now 
$$\eqalign{
|g(z)-y_j| &= |\lambda^{-1}\psi^{-1}
(\wt g(z))-\lambda^{-1}\psi^{-1}(\wt y_j)|
\leq \lambda^{-1} \Lambda |\wt g(z)-\wt y_j|
\cr&
\leq \lambda^{-1} \Lambda |\wt f(z)-\wt y_j|
= \lambda^{-1} \Lambda |\psi(\lambda f(z))-\psi(\lambda y_j)|
\cr&
\leq \Lambda^2 |f(z)-y_j| \leq \Lambda^2 r_j
\leq \Lambda^2 \delta_6 \leq {\delta_2 \over 10}
}\leqno (18.67)
$$
because $\wt g(z) = \wt\pi_j(\wt f(z))$ and
$\wt Q_j$ goes through $\wt y_j = \wt f(y_j) = \wt A_x(x)$, then
by (12.75) and (15.47), then
because $f(z) \in D_j$ by (15.22), and by (15.15) and (12.7).
Since $\wt g(z) \in \wt Q_j \i \wt W$ and $\wt y_j \in \wt W$,
the line segment $L = [\wt y_j, \wt g(z)]$ is contained in $\wt W$.
Since for $\xi \in L$, 
$$\eqalign{
|\lambda^{-1}\psi^{-1}(\xi)-y_j|
&=|\lambda^{-1}\psi^{-1}(\xi)-\lambda^{-1}\psi^{-1}(\wt y_j)|
\leq \lambda^{-1} \Lambda |\xi-\wt y_j|
\cr&
\leq \lambda^{-1} \Lambda |\wt g(z)-\wt y_j|
\leq {\delta_2 \over 10}
}\leqno (18.68)
$$
by the end of (18.67), (18.66) says that $\lambda^{-1}\psi^{-1}(L)$ does not meet 
$\d F(y_j) \i {\cal S}_{m-1}$. Since $y\in F(y_j)$, we get that 
$g(z)\in F(y_j)$ too.
We deduce from (18.66) and (18.67) that
$\dist(g(z),{\cal S}_{m-1}) \geq {\delta_2 \over 2}$,
as in (18.61), which implies that
$$
\dist(\wt g(z),\wt{\cal S}_{m-1}) 
\geq {\lambda \Lambda^{-1} \delta_2 \over 2}
\leqno (18.69)
$$
where we denote by $\wt{\cal S}_{m-1} = \psi(\lambda {\cal S}_{m-1})$
the $(m-1)$-dimensional skeleton in the standard dyadic grid.

Recall that $\wt h_2(z) = \Pi(\wt g(z), C\lambda \varepsilon_0)$
by (18.50); we may now use the definitions (17.23) and (17.24)
as above, and the same argument based on (17.2) yields that
$\wt h_2(z) = \wt g(z)$, if $\varepsilon_0$  is small enough
compared to $\delta_2$; (18.65) follows.

We'll also need to know that for $j\in J_3$ and $x\in Z(y_j)$,
and $k$ large enough
$$
h_2(E_k \cap B_{j,x}^-) \i D_j \cap \lambda^{-1}\psi^{-1}(\wt Q_j).
\leqno (18.70)
$$
Let $z\in E_k \cap B_{j,x}^-$ be given. We already know from (18.65)
that $h_2(z) = g(z) \in \lambda^{-1}\psi^{-1}(\wt Q_j)$,
so we just need to check that $g(z) \in D_j$.

By (15.22), $f(z) \in {1+a \over 2}\, D_j$. 
If $k$ is large enough, then by (10.4)
$z\in E^{\varepsilon r_j}$, and (16.6) and (15.60) say that
$|g(z) - f(z)| = |g_j(z) - f(z)| \leq \Lambda^2 (1+3|f|_{lip}) 
\varepsilon r_j$; hence  $g(z) \in D_j$ and (18.70) follows. 

We may now follow the proof of (18.63) and (18.64); we get that
for $k$  large
$$\eqalign{
\H^d\Big(h_2\Big(E_k \cap \bigcup_{j\in J_3} \bigcup_{x\in Z(y_j)} 
B_{j,x}^+ \sm R^3\Big)\Big)
&
\leq \sum_{j\in J_3} \H^d(D_j \cap \lambda^{-1}\psi^{-1}(\wt Q_j))
}\leqno (18.71)
$$
and then
$$\eqalign{
\H^d(h_2(E_k \cap W)) 
&\leq 
C \eta + C(f,\gamma) (1-a) + C(\alpha,f) N^{-1} + C(f) \gamma
\cr& \hskip2.5cm
+ \sum_{j\in J_3} \H^d(D_j \cap \lambda^{-1}\psi^{-1}(\wt Q_j)).
}\leqno (18.72)
$$
We will only see in the next section how to use our assumption (10.7)
or a weaker (but more complicated) one to control the last sum.
$\dagger$

\msi{\bf A lower bound for $H^d(f(E \cap W))$.}
We found in (18.64) or (18.72), a first upper bound
for the right-hand side of (18.1). The main term
in (18.64) is $\omega_d \sum_{j\in J_3} r_j^d$, and
we want to bound it by $H^d(f(E \cap W))$, plus small errors.
The following lemma, which is our analogue of 
Lemma 4.111 in [D2], will be useful.  

\ms\proclaim Lemma 18.73. For each $j\in J_3$,
$$
\H^d(D_j \cap f(E \cap W_f)) \geq \big(1-C(f,\gamma)\varepsilon \big) 
\, \omega_d r_j^d,
\leqno (18.74)
$$
where $C(f,\gamma)$ depends only on $|f|_{lip}$ and $\gamma$.

\ms
We shall give a different proof here, so as not to have
to construct a competitor again. Instead we shall take advantage of
the fact that we could choose extremely small balls $D_j$,
which we control by differentiability and density results.
Our proof of Lemma 18.73 will also work, with no modification,
in the Lipschitz context.

\ms
Let $j\in J_3$ be given and pick some $x\in Z(y_j)$; by 
(15.12), $y_j \in Y_{11} \i Y_9 = f(X_9)$ (see above (15.1)), so 
$Z(y_j)$ is not empty and we can choose $x\in Z(y_j)$. Then
$x\in X_9 \i X_5 = \cup_{s\in S} X_5(s)$ (see (15.1) and (11.47)), 
so there is an index $s$ such that the description 
in (11.42)-(11.46) is valid. In particular, the graph $\Gamma_s$
contains all the points $z+F_x(z)$, $z\in P_x \cap B(x,\delta_3)$,
where $F_x : P_x \to P_x^\perp$ is the $C^1$ function of (11.42).
Let $Q_j = A_x(P_x)$ be as above, denote by $\pi_j$ the 
orthogonal projection onto $Q_j$, and
define $G : P_x \cap B(x,{\delta_3 \over 2}) \to Q_j$ by
$$
G(z) = \pi_j(f(z+F_x(z))) 
\ \hbox{ for } 
z \in P_x \cap B(x,{\delta_3 \over 2}).
\leqno (18.75)
$$
Even under the Lipschitz assumption, we really want to use
$Q_j$ and the fact that $f$ itself (and not $\wt f$) is
well approximated by the affine function $A_x$, as in (11.46).
Notice that for $z \in P_x \cap B(x,{\delta_3 \over 2})$ (as in 
(18.75)),
$$
|F_x(z)| \leq ||DF_x||_\infty |z-x| \leq \varepsilon |z-x|
\leqno (18.76)
$$
by (11.42), so $z+F_x(z) \in \Gamma_s \cap B(x,\delta_3)$
by (11.43), and 
$$
|f(z+F_x(z))- A_x(z+F_x(z))| \leq \varepsilon |z+F_x(z)-x|
\leq 2 \varepsilon |z-x|
\leqno (18.77)
$$
by (11.46). Then
$$\eqalign{
|G(z)-A_x(z)| &= |\pi_j(f(z+F_x(z)))-A_x(z)|
= |\pi_j(f(z+F_x(z)))-\pi_j(A_x(z))|
\cr&
\leq |f(z+F_x(z))-A_x(z)|
\cr&
\leq |f(z+F_x(z))- A_x(z+F_x(z))| + |A_x(z+F_x(z)) - A_x(z)|
\cr&
\leq 2 \varepsilon |z-x| + |f|_{lip} |F_x(z)|
\leq (2+|f|_{lip}) \varepsilon |z-x|
}\leqno (18.78)
$$
by (18.75), because $\pi_j(A_x(z))=A_x(z)$ (since $z\in P_x$
and so $A_x(z) \in Q_j$ by definition of $Q_j$), and by (18.77),
(11.36), and (18.76).

We can apply this to $z\in 2\overline E(x)$, where 
$E(x) = P_x \cap A_x^{-1}(Q_j \cap D_j)$ is
as in (15.16), because 
$$
2E(x) \i 2B_{j,x} \i B(x,3\gamma^{-1} r_j) \i B(x,\delta_6)
\i B(x,{\delta_3 \over 10})
\leqno (18.79)
$$
by (15.19), (15.20), (15.15), and (12.7).
[Again all those things hold in the Lipschitz case; see
the remark below (15.39).] We get that
$$
|G(z)-A_x(z)| \leq (2+|f|_{lip}) \varepsilon |z-x|
\leq 3\gamma^{-1} (2+|f|_{lip}) \varepsilon r_j
\ \hbox{ for } z\in 2\overline E(x).
\leqno (18.80)
$$
Denote by $\d$ the boundary of $2E(x)$ in $P_x$;
because $x\in X_9 \i X_8$ and by the definitions (14.21) and (14.5),
$A_x$ is a bijective affine mapping from $\d$
to $Q_j \cap \d (2D_j)$. Call $\d'$  the unit sphere in
the vector $d$-space parallel to $Q_j$. For each 
$w \in Q_j \cap D_j$; the mapping 
$a_w : z \to {A_x(z)-w \over |A_x(z)-w|}$, from 
$\d$ to $\d'$, is well-defined (because $A_x(z)$ does
not take the value $w$). Its degree is the same for all
$w \in Q_j \cap D_j$, and it is equal to $\pm 1$
because this is its value at $w=y_j$.

For $w\in Q_j \cap D_j$ and $z\in \d$, 
$A_x(z) \in Q_j \cap \d (2D_j)$, so $|A_x(z)-w| \geq r_j$,
hence by (18.80) the segment 
$[A_x(z),G(z)]$ does not contain $w$.
Then 
$$
(z,t) \to a_{w,t}(z) = {(1-t)A_x(z)+tG(z)-w \over |(1-t)A_x(z)+tG(z)-w|}
\leqno (18.81)
$$
is defined and continuous on $\d \times [0,1]$. It
defines a homotopy from $a_w$ to $a_{w,1}$,
among mappings from $\d$ to $\d'$. Thus $a_{w,1}$ has the
same degree as $a_w$, namely, $\pm 1$. Then $a_{w,1}$ does not 
extend continuously as a mapping from $2 \overline E(x)$ to $\d'$, 
and this forces $G(z)-w$ to vanish at some point $z\in 2 \overline E(x)$
(otherwise, use ${G(z)-w \over |G(z)-w|}$). We just proved
that 
$$
G(2 \overline E(x)) \hbox{ contains } Q_j \cap D_j.
\leqno (18.82)
$$
Set $\lambda = 1-3\gamma^{-1} (2+|f|_{lip}) \varepsilon$.
We want to estimate the size of
$$
Y = Q_j \cap \lambda D_j \sm \pi_j[D_j \cap f(E \cap W_f)].
\leqno (18.83)
$$
Let $w\in Y$ be given, and use (18.82) to find
$z\in 2 \overline E(x)$ such that $G(z) = w$. 
Set $y = z+F_x(z)$ and observe that 
$$
|y-x| = |z+F_x(z)-x| \leq |z-x| + \varepsilon |z-x|
< 4 \gamma^{-1} r_j < {\delta_3 \over 5}
\leqno (18.84)
$$
by (18.76) and (18.79), so $y \in \Gamma_s$ by (11.43).
Notice that 
$$
w = G(z) = \pi_j(f(z+F_x(z))) = \pi_j(f(y))
\leqno (18.85)
$$
by (18.75) and other definitions.

If $y \in E \cap W_f$, we observe that
$$\eqalign{
|w-f(y)| &= |\pi_j(f(y))-f(y)| 
\leq \dist(f(y),Q_j)
\leq |f(y)-A_x(z)| 
\cr&
= |f(z+F_x(z))-A_x(z)|
\leq (2+|f|_{lip}) \varepsilon |z-x|
\leq 3\gamma^{-1} (2+|f|_{lip}) \varepsilon r_j
}\leqno (18.86)
$$
by (18.85), because $A_x(z) \in Q_j$, by the last lines of (18.78)
and the end of (18.80). So $f(y) \in D_j$ because $w\in \lambda D_j$,
and $w = \pi_j(f(y))$ lies in $\pi_j[D_j \cap f(E \cap W_f)]$, 
a contradiction with the definition of $Y$.

So $y \notin E \cap W_f$. But 
$$
y \in B(x,4 \gamma^{-1} r_j) \i B(x,\delta_6)
\i B(x,\delta_1/10) \i W_f,
\leqno (18.87)
$$
by (18.84), (15.15), (12.7), (11.22), and the fact that $x\in X_9 \i X_1$.
So $y \notin E$. In addition, $y\in \Gamma_s$
(see below (18.84)). Set $\rho = 4 \gamma^{-1} r_j$;
we just proved that $y\in B(x,\rho) \cap \Gamma_s \sm E$;
hence $w = \pi_j(f(y)) \in \pi_j(f(B(x,\rho) \cap \Gamma_s \sm E))$
and now
$$\eqalign{
\H^d(Y) &\leq \H^d(\pi_j(f(B(x,\rho) \cap \Gamma_s \sm E)))
\leq |f|_{lip}^d \H^d(B(x,\rho) \cap \Gamma_s \sm E).
}\leqno (18.88)
$$
We want to apply (11.44) to $B(x,\rho)$; this is allowed because 
$\rho = 4 \gamma^{-1} r_j \leq \delta_6 \leq \delta_3/10$
(again by (15.15) and (12.7)), and (11.44) yields
$$
\H^d(B(x,\rho) \cap [\Gamma_s \cup E] \sm X_3(s)) 
\leq \varepsilon \rho^d \leq C(\gamma) \varepsilon  r_j^d
\leqno (18.89) 
$$
Now $\Gamma_s \sm E \i \Gamma_s \cup E \sm X_3(s)$, just because
$X_3(s) \i X_0 \i E$, so (18.88) yields
$$
\H^d(Y) \leq C(\gamma) |f|_{lip}^d \varepsilon r_j^d.
\leqno (18.90) 
$$
Finally,
$$\eqalign{
\H^d(D_j \cap f(E \cap W_f)) 
&\geq \H^d(\pi_j[D_j \cap f(E \cap W_f)])
\geq \H^d(Q_j \cap \lambda D_j) - \H^d(Y)
\cr&
=\lambda^d \omega_d r_j^d - \H^d(Y)
\geq \omega_d r_j^d - C(f,\gamma) \varepsilon r_j^d
}\leqno (18.91)
$$
by (18.83), (18.90), and because  
$\lambda = 1-3\gamma^{-1} (2+|f|_{lip}) \varepsilon$.
Lemma 18.73 follows.
\qed

\msi{\bf The final estimate.}
We are now ready to conclude, at least under the rigid assumption.
We sum (18.74) over $j\in J_3$ and get that
$$\eqalign{
\sum_{j\in J_3} \omega_d r_j^d
&\leq \big(1-C(f,\gamma)\varepsilon \big)^{-1} \sum_{j\in J_3} 
\H^d(D_j \cap f(E \cap W_f))
\cr&
\leq \big(1-C(f,\gamma)\varepsilon \big)^{-1} \H^d(f(E \cap W_f)) 
\leq \H^d(f(E \cap W_f)) + C'(f,\gamma)\varepsilon
}\leqno (18.92)
$$
because the $D_j$ are disjoint (see (15.13)), 
and if $\varepsilon$ is small enough.
We compare this with (18.64) and get that
$$
\H^d(h_2(E_k \cap W)) \leq \H^d(f(E \cap W_f)) + {\cal E},
\leqno (18.93)
$$
with
$$
{\cal E} \leq C \eta + C(f,\gamma) (1-a) + C(\alpha,f) N^{-1} + C(f) \gamma
+ C'(f,\gamma)\varepsilon.
\leqno (18.94)
$$
Recall from (18.46) that $|h_2(x)-x| \geq {\delta_5 \over 4}$
for $x\in E_k \cap V_1^+$. Since $\{ E_k \}$ converges to $E$
and $h_2$ is continuous, we also get that 
$$
|h_2(x)-x| \geq {\delta_5 \over 4}
\hbox{ for } x\in E \cap {\rm int}(V_1^+), 
\leqno (18.95)
$$
and hence $E \cap {\rm int}(V_1^+) \i W$
(recall that $W = \big\{ y \in U \, ; h_2(x) \neq x \big\}$ 
by (18.2)). Hence $E \sm W \i E \sm {\rm int}(V_1^+)$, and 
$$\eqalign{
\H^d(E \cap W_f) - \H^d(E \cap W) &\leq \H^d(E \cap W_f \sm W)
\leq \H^d(E \cap W_f \sm {\rm int}(V_1^+))
\cr&
= \H^d(E \cap W^2)
\leq C(f,\gamma) (1-a) + C \eta,
}\leqno (18.96)
$$
by (18.38) and (18.39).
At the same time, $W$ is open (see the definition above), and
Theorem 10.97 (our main lower semicontinuity result)
says that for $k$ large,
$$
\H^d(E \cap W) \leq \H^d(E_k \cap W) + \eta.
\leqno (18.97)
$$
By (18.96), (18.97), and then (18.1) and (18.93),
$$\eqalign{
\H^d(E \cap W_f) 
&\leq \H^d(E \cap W) + C(f,\gamma) (1-a) + C \eta
\cr&
\leq \H^d(E_k \cap W) + C(f,\gamma) (1-a) + C \eta
\cr& 
\leq M \H^d(h_2(E_k \cap W)) + h (R'')^d + C(f,\gamma) (1-a) + C \eta
\cr&
\leq M \H^d(f(E \cap W_f)) + M {\cal E} + h (R'')^d 
+ C(f,\gamma) (1-a) + C \eta
\cr&
\leq M \H^d(f(E \cap W_f)) + {\cal E}' + h (R'')^d 
}\leqno (18.98)
$$
for $k$ large, where  ${\cal E}'$ is given by the same
sort of formula (18.94) as ${\cal E}$, and where 
$$
R'' = R_0 + 4\Lambda^2 (1+ |f|_{lip}) \delta_6 + C \Lambda \varepsilon_0
\leqno (18.99)
$$
is still as in (17.41). Now we choose our various small constants
$\gamma$, $a$, $\alpha$, $N$, $\eta$, $\varepsilon$, $\delta_6$, 
and $\varepsilon_0$ in this order (as announced in Remark 11.17), 
so as to make ${\cal E}'$ and $(R'')^d-R_0^d$ arbitrarily small. 
Since (18.98) holds for all these choices, we get that
$$
\H^d(E \cap W_f) \leq M\H^d(f(E \cap W_f)) + h R_0^d.
\leqno (18.100)
$$
This is the same as (2.5) in the circumstances that were described at 
the beginning of Section 11 : compare $E \cap W_f$ in (11.19) with 
$W_1$ in Definition 2.3, and recall that $f = \varphi_1$ (see above 
(11.18)).
So we completed the verification of (10.9), and therefore proved 
Theorem~10.8 in the special case when we have the rigid assumption.
\qed

\msi 
{\bf 19. Proof of Theorem~10.8, and variants, under the Lipschitz assumption.} 
\ms
In this section we work under the Lipschitz assumption, and try to 
prove (18.100) (almost) as in the previous section.

The only remaining difficulty is that we have (18.72)
instead of (18.64), and the difference between the two right-hand sides
is
$$
\Delta = \sum_{j\in J_3} \H^d(D_j \cap \lambda^{-1}\psi^{-1}(\wt Q_j))
- \omega_d \sum_{j\in J_3} r_j^d.
\leqno (19.1)
$$
If we follow the proof above and use (18.72) instead of (18.64),
we do not need to modify Lemma 18.73 (which is still valid
in the Lipschitz context), but we need to add $\Delta$
to $\cal E$ in (18.93) and (18.94), and then $M\Delta$ to
${\cal E}'$ in (18.98).
So, if we could prove that $\Delta$ can be made arbitrarily small
(with a choice of constants as above), then we could quietly
follow the same proof above and get the conclusion.
Before we do this, we need to modify a little our definition
of our final mapping $h_2$ in the $B_{j,x}$, $x\in Z(y_j)$,
for some $j \in J_3$.

\msi{\bf The indices $j\in J_4$ for which $y_j \in L_i$
only if it lies in its interior.}
We define the set $J_4 \i J_3$ to be the set of indices
$j\in J_3$ such that, for all $i \in [0, j_{max}]$ such that
$y_j \in L_i$, $y_j$ actually lies in the $n$-dimensional
interior of $L_i$ (that is, for the ambient topology of $\R^n$).

For these $j$, we were a little too prudent with the definition 
of $g_j$ and $h_2$, because the boundary condition (1.7) is much easier 
to fulfill in these cases. 
The truth is that we should have defined the corresponding
$g_j$ differently, but rather than modifying the construction above,
we shall fix it by continuing our deformation $\{ h_t \}$, 
$0 \leq t \leq 2$, a little further. Set
$$
h_t(z) = h_2(z)
\ \hbox{ when $2 < t \leq 3$ and }
z \in U \sm \bigcup_{j \in J_4} \bigcup_{x\in Z(y_j)} B_{j,x}^-.
\leqno (19.2)
$$
In the remaining sets $B_{j,x}^-$, $j \in J_4$ and $x\in Z(y_j)$,
we can proceed independently  (because these sets are disjoint by
(15.27) and (15.20)). Fix $j\in J_4$, denote by 
$Q_j$ the common value of the $A_x(P_x)$, $x\in Z(y_j)$, as we did 
before, and define the
$h_t$, $2 < t \leq 3$, on $B_{j,x}^-$ by
$$
h_t(z) = (1-\beta_{j,x}(z,t)) h_2(z) + \beta_{j,x}(z,t) \pi_j(h_2(z)),
\leqno (19.3)
$$
where $\pi_j$ denotes the orthogonal projection on $Q_j$
and
$$
\beta_{j,x}(z,t) = \Min\Big( 1,
{100(1+|f|_{lip}) \over (1-a) r_j}\, (t-2) \dist(z,\d B_{j,x}^-)
\Big).
\leqno (19.4)
$$
Notice that $\beta_{j,x}(z,t) = 0$, and hence $h_t(z)=h_2(z)$,
when $t=2$ and when $z\in \d B_{j,x}^-$, so we glue things
in a continuous way. The final mapping $h_3$ is even 
Lipschitz, because it is Lipschitz on each $B_{j,x}^-$
(there is a finite collection of them), on the rest of $\R^n$, 
and is continuous across the boundary.

We need to check that the $h_{3t}$, $0 \leq t \leq 3$ still
satisfy the conditions (1.4)-(1.8), relative to the set $E_k$,
as in Lemma 17.40 and (17.87) (proved below (17.107)).

The continuity condition (1.4) follows from its counterpart 
for the $h_{2t}$, and the Lipschitz property (1.8) was just discussed. 
For the other properties (1.5)-(1.7), and the constraint on the $\wh W$-set, 
we just need to worry about the only places where we change something, i.e.,
the sets $B_{j,x}^-$. 

So let $j \in J_4$, $x\in Z(y_j)$, and $z\in E_k \cap B_{j,x}^-$ be given,
and let us derive some general information on $z$ and the $h_t(z)$.
We know that $x\in X_0$, so $\varphi_1(x) \neq x$ by (11.19) and (11.20), 
which implies that $x\in B$ (by (1.5)) and that $x$ and $f(x)$ both lie in $\wh W$
(by the definition (2.2)). Since in addition 
$$
|z-x| \leq 2 \gamma^{-1} r_j \leq \delta_6
\leqno (19.5)
$$
by (15.20) and, say, (15.24), the first part implies that 
$z\in B(X_0,R_0 + \delta_6) \i B(X_0,R'')$
(see (11.1) and (17.41)). Then (1.5) holds.

Let us check that
$$
h_2(z) = g(z) \ \hbox{ for } z\in E_k \cap B_{j,x}^-.
\leqno (19.6)
$$
By (18.65), $\wt h_2(z) = \wt g(z)$.
In addition, $z\in V_1^+$ (see the definition (16.7) and (15.20)),
hence $z \in H_1$ by (18.43). By (17.101), we can compute
$h_2(z)$ and $g(z)$ in terms of $\wt h_2(z)$ and $\wt g(z)$; 
(19.6) follows.

Moreover, we get that for $z\in E_k \cap B_{j,x}^-$ and $t \geq 2$,
$$
|h_t(z)-y_j| \leq \Max(|h_2(z)-y_j|,|\pi_j(h_2(z))-y_j|)
= |h_2(z)-y_j| = |g(z)-y_j| \leq \Lambda^2 \delta_6
\leqno (19.7)
$$
by (19.3), because $y_j = f(x) = A_x(x) \in Q_j$ by definitions,
and by (18.67) (which holds because $z \in E_k \cap B_{j,x}^-$).

Return to our verifications. For (1.6), we need to check that
$h_t(z) \in B(X_0,R'')$ when $z\in E_k \cap B(X_0,R'')$,
and since we already know about $h_t$, $0 \leq t \leq 2$,
we can restrict to $z\in E_k \cap B_{j,x}^-$ as above. 
The desired estimate follows from (19.7) because $y_j = f(x) \in B$, 
by (1.5) for $f = \varphi_1$
and because $x\in X_9 \i E_0 = E \cap W_f \i B$
(see (15.1), (11.20), and (11.19)).

For the verification of (2.4), recall that $x$ and $y_j = f(x)$ both lie in 
(the old set) $\wh W$, and so $z$ and the new $h_t(z)$
lie within $\Lambda^2 \delta_6$ of $\wh W$, by (19.7) in particular.
This still puts them in a compact subset of $U$, by (12.6) and (12.7).

We are left with (1.7) to check. 
That is, we fix $0 \leq i \leq j_{max}$, and we want to check that
for $k$ large, $h_t(z) \in L_i$ for $0 \leq t \leq 3$
as soon as $z\in E_k \cap L_i$. We know this when
$0 \leq t \leq 2$, by (1.7) for the $h_{2t}$ 
(see (17.87)); so we may assume that $t > 2$.
We also know this when 
$z \in E_k \sm \bigcup_{j \in J_4} \bigcup_{x\in Z(y_j)} B_{j,x}^-$,
because (19.2) says that $h_t(z)=h_2(z) \in L_i$.
So we may fix $j\in J_4$, $x\in Z(y_j)$, and it is enough
to check that 
$$
h_t(z) \in L_i
\hbox{ for $z\in E_k \cap L_i \cap B_{j,x}^-$ and }
2 \leq t \leq 3.
\leqno (19.8)
$$
We first assume that $y_j \notin L_i$. Let us 
check that for $k$ large,
$$
E_k \cap L_i \cap B_{j,x}^- = \emptyset
\leqno (19.9)
$$
for $x\in Z(y_j)$; (19.8) will follow trivially.
Let $z\in B_{j,x}^-$ be given; we proved
below (18.68) that $g(z) \in F(y_j)$, where $F(y_j)$
is the smallest face of our grid that contains $y_j$.
In fact, since $g(z)$ is also far from the boundary
of $F(y_j)$ (by (18.69)), we see that $g(z)$ lies in
$I$, the interior of $F(y_j)$.
Now suppose that $z\in L_i$ for some $i$.
Since $g(z)=h_2(z)$, (1.7) for $h_2$ says that
$g(z) \in L_i$. Then some face $F$ of $L_i$ meets $I$,
and since $I$ is the interior of a face, $F$ contains $I$.
But $y_j \in I$, by definition of the smallest face $F(y_j)$,
so $y_j \in L_i$. Our current assumption says that this is impossible, 
so  $z\notin L_i$, and (19.9) follows. 

We are left with the case when $y_j \in L_i$. Since $j\in J_4$,
this implies that $y_j$ is an interior point of $L_i$,
and we want to deduce from this that
$$
\overline B(y_j,\Lambda^2 \delta_6) \i L_i.
\leqno (19.10)
$$ 
Since by (19.7) $|h_t(z)-y_j| \leq \Lambda^2 \delta_6$ 
for $z\in E_k \cap B_{j,x}^-$ and $t \geq 2$, (19.8)
will follow at once.

Denote by $\delta(L_i)$ the (true $n$-dimensional) boundary of $L_i$,
and set $D= \dist(y_j,\delta(L_i))$; 
we know that $D>0$, and we want to show that $D \geq \Lambda^2 \delta_6$.
Take $\xi \in \delta(L_i)$ such that $|\xi-y_j| = D$, 
and let $F$ denote the smallest face of our grid that contains $\xi$. 
Since $\delta(L_i)$ is composed of full faces too, $F$ is contained in 
$\delta(L_i)$. Hence $y_j \notin F$, and $F$ does not contain the smallest 
face $F(y_j)$ that contains $y_j$. 
First assume that $F(y_j)$ is not reduced to $\{ y_j \}$; then
(3.8) (and a conjugation by $\psi$ that makes us lose a constant 
$\Lambda^2$) yields
$$
D = |\xi-y_j| \geq \dist(y_j,F) \geq \Lambda^{-2} \dist(y_j,\d F(y_j)).
\leqno (19.11)
$$
In addition, $x\in Z(y_j)$, so (15.1) says that
$x \in X_9 \i X_2 = X_{1,\delta_2}$. By (11.23)-(11.25), this means
that $x\in X_{1,\delta_2}(m)$, where $m$ is the dimension $F(y_j)$,
and hence $\dist(y_j, {\cal S}_{m-1}) \geq \delta_2$.
Then $D \geq \Lambda^{-2} \delta_2 \geq 10 \delta_6$
by (12.7). The sad truth is that this is not enough, but this
is easy to fix: we just need to require that $\delta_6 \leq 2 
\Lambda^{-4} \delta_2$ in (addition to) (12.7), 
or replace $\Lambda^2$ with $\Lambda^4$ in (12.7),
and then we get the desired estimate which implies (19.10).
In the remaining case when $F(y_j)=\{ y_j \}$, $D$ is at least
equal to the smallest distance between a vertex of the grid
(namely $y_j$) and a face that does not contain it.
This distance is at least $\lambda^{-1} \Lambda^{-1} r_0$;
hence $D \geq \lambda^{-1} \Lambda^{-1} r_0 \geq 2 \Lambda^2 
\delta_6$, if we put an extra power of $\Lambda$ in (12.7).

This completes our verification of (1.7) for our extended family
of $h_t$, and also the verification of (1.4)-(1.8). 
We also checked (2.4) for the extended family
$\{ h_t \}$ along the way, so we have the analogue 
of (18.1) for $h_3$.

\ms
Now we want to see why $h_3$ is possibly better than $h_2$.
We do not need to modify any of our estimates, except for the
ones relative to the $B_{j,x}^-$, $x\in Z(y_j)$ and $j\in J_4$.
Fix such $j$ and $x$, and let us first check that for $k$ large,
$$
h_3 \hbox{ is $C\Lambda^2(1+|f|_{lip})$-Lipschitz on } 
E_k \cap B_{j,x}^-.
\leqno (19.12)
$$
We already know that $h_2(z)=g(z)=g_{j}(z)$ on $E_k \cap B_{j,x}^-$
(also see (16.6)). But (15.62) says that $g_j$ is 
$C\Lambda^2(1+|f|_{lip})$-Lipschitz on 
$E^{\varepsilon r_j} \cap B_{j,x}^-$, hence also on 
$E_k \cap B_{j,x}^-$ (just because 
$E_k \cap B_{j,x}^- \i E^{\varepsilon r_j}$ for $k$ large).
But we need to worry a little about the rapid fluctuations
of $\beta_{j,x}(z,3)$. As usual, pick 
$z, w \in E_k \cap B_{j,x}^-$ and write
$$\leqalignno{
h_3(z)-h_3(w) &= 
(1-\beta_{j,x}(z,3)) h_2(z) + \beta_{j,x}(z,3) \pi_j(h_2(z))
\cr&\hskip3cm
-(1-\beta_{j,x}(w,3)) h_2(w) - \beta_{j,x}(w,3) \pi_j(h_2(w))
\cr&
= - [\beta_{j,x}(z,3)-\beta_{j,x}(w,3)]\,[h_2(z)-\pi_j(h_2(z))]
&(19.13)
\cr&\hskip0.6cm
+(1- \beta_{j,x}(w,3)) [h_2(z)-h_2(w)]
+ \beta_{j,x}(w,3) [\pi_j(h_2(z))-\pi_j(h_2(w))]
}
$$
so that
$$\leqalignno{
|h_3(z)-h_3(w)| &\leq |\beta_{j,x}(z,3)-\beta_{j,x}(w,3)|\,
|h_2(z)-\pi_j(h_2(z))|
+ C\Lambda^2(1+|f|_{lip}) |z-w|
\cr&
\leq {100 \Lambda(1+|f|_{lip}) \over (1-a) r_j}
|z-w| \, |h_2(z)-\pi_j(h_2(z))|
+ C\Lambda^2(1+|f|_{lip}) |z-w|
&(19.14)
}
$$
because $h_2$ is $C\Lambda^2(1+|f|_{lip})$-Lipschitz and by (19.4).
Then we want to estimate $|h_2(z)-\pi_j(h_2(z))|$, i.e., prove
that $h_2(z)$ is close to $Q_j$. Recall from
(18.65) and the line below that 
$$
\wt h_2(z) \in \wt Q_j = \wt A_x(P_x); 
\leqno (19.15)
$$
also, we checked below (19.6) that $z\in E_k \cap H_1$,
so we may apply (17.101) and we get that
$$
\wt h_2(z) = \psi(\lambda h_2(z)).
\leqno (19.16)
$$
Then
$$\eqalign{
|\wt h_2(z)-\wt f(x)| &= |\psi(\lambda h_2(z))-\psi(\lambda f(x))|
\leq \lambda \Lambda |h_2(z)-f(x)|
\cr&
= \lambda \Lambda |h_2(z)-y_j|
\leq \lambda \Lambda r_j
}\leqno (19.17)
$$
because $h_2(z) \in D_j$ by (18.70).

Use (19.15) to write $\wt h_2(z) = \wt A_x(\xi)$, with $\xi \in P_x$;
we want to evaluate $|\xi-x|$. Recall from just above (19.3) that 
$Q_j = A_x(P_x)$. The discussion near (15.41) says that
that since $y_j \in Y_{10}$ (by (15.12)),
the restriction of $\psi$ to $Q_j$ is differentiable
at $\lambda y_j$, with a derivative $D_{\psi}$ such that
$$
\lambda D_{\psi}(DA_x(v)) = D\wt A_x(v)
\ \hbox{ for } v \in P'_x \, ,
\leqno (19.18)
$$
where $P'_x$ denotes the vector space parallel to $P_x$.
Since $\wt A_x$ is affine and $\wt A_x(x) = \wt f(x)$
(see (12.38)), we get that
$$\eqalign{
|\wt h_2(z) - \wt f(x)| &= |\wt A_x(\xi) - \wt A_x(x)|
= |D\wt A_x(\xi-x)| = |\lambda D_{\psi}(DA_x(\xi-x))|
\cr&
\geq \lambda\Lambda^{-1} |DA_x(\xi-x)|
\geq \lambda\Lambda^{-1} \gamma |\xi-x|
}\leqno (19.19)
$$
because $D_{\psi}$, just like $\psi$ itself,
is $\Lambda$-biLipschitz, and because $DA_x$ has no
contracting direction since 
$x\in Z(y_j) \i X_9 \i X_8 \i  X_6 \sm X_7$; see (15.1), (14.21),
and (14.5). We compare (19.19) to (19.17) and get that
$$
|\xi-x| \leq \Lambda^2\gamma^{-1} r_j
\leq \Lambda^2 \delta_6 \leq {\delta_3 \over 10}
\leqno (19.20)
$$
by (15.15) and (12.7). Now
$$\eqalign{
|h_2(z)-\pi_j(h_2(z))| &= \dist(h_2(z),Q_j)
\leq |h_2(z) - A_x(\xi)| 
\cr&
\leq |h_2(z) - f(\xi)| + |f(\xi) - A_x(\xi)|
}\leqno (19.21)
$$
because $A_x(\xi) \in Q_j$
(since $\xi \in P_x$ and $Q_j = A_x(P_x)$).
By (19.20), we can apply (11.46) and get that
$$
|f(\xi) - A_x(\xi)| \leq \varepsilon |\xi-x|.
\leqno (19.22)
$$
Also, (19.20) allows us to apply
(12.52), which says that $\wt f(\xi) = \psi(\lambda f(\xi))$
is well defined, and also that
$|\wt f(\xi)-\wt A_x(\xi)| \leq \lambda \varepsilon |\xi-x|$.
Then
$$\eqalign{
|h_2(z) - f(\xi)| &= |\lambda^{-1}\psi^{-1}(\wt h_2(z)) 
- \lambda^{-1}\psi^{-1}(\wt f(\xi))|
\leq \lambda^{-1} \Lambda |\wt h_2(z)-\wt f(\xi)|
\cr&
= \lambda^{-1} \Lambda |\wt A_x(\xi)-\wt f(\xi)|
\leq \Lambda \varepsilon |\xi-x|
}\leqno (19.23)
$$
by (19.16), because $\wt h_2(z) = \wt A_x(\xi)$ 
by definition of $\xi$, and by (12.52).
Altogether,
$$
|h_2(z)-\pi_j(h_2(z))| \leq (1+\Lambda) \varepsilon |\xi-x|
\leq (\Lambda^2+\Lambda^4) \gamma^{-1} \varepsilon r_j
\leqno (19.24)
$$
by (19.21), (19.22), (19.23), and (19.20).

For the record, notice that (19.3) and 
(19.24) imply that for $z\in E_k \cap B_{j,x}^{-} \,$,
$$
|h_3(z)-h_2(z)| \leq |h_2(z)-\pi_j(h_2(z))| 
\leq \varepsilon (\Lambda^4+\Lambda^2) \gamma^{-1} r_j
\leq C \varepsilon \Lambda^4 \delta_6
\leq \delta_5/10
\leqno (19.25)
$$
by (15.15), (12.7), and if $\varepsilon$ is small enough 
(depending on $\Lambda$).
We used $\varepsilon$ here just so that we don't have to put an
extra power of $\Lambda$ in the definition of $\delta_6$, but we 
could have done that too. Since we simply have that $h_3(z)=h_2(z)$
when $z$ lies in no $B_{j,x}^{-} \,$, we get that for $k$ large,
$$
|h_3(z)-h_2(z)| \leq \delta_5/10
\ \hbox{ for } z \in E_k.
\leqno (19.26)
$$

Let us return to our $z$ and $w$, plug (19.24) into (19.14), and get that
$$\leqalignno{
|h_3(z)-h_3(w)| 
&\leq {100 \Lambda(1+|f|_{lip}) \over (1-a) r_j}
|z-w| \, |h_2(z)-\pi_j(h_2(z))|
+ C\Lambda^2(1+|f|_{lip}) |z-w|
\cr&
\leq C (\Lambda^4+\Lambda^2) \gamma^{-1} \varepsilon r_j
{\Lambda(1+|f|_{lip}) \over (1-a) r_j} \, |z-w|
+ C\Lambda^2(1+|f|_{lip}) |z-w|
&(19.27)
\cr&
\leq C\Lambda^2(1+|f|_{lip}) |z-w|
}
$$
if $\varepsilon$ is small enough, depending on $\Lambda$, $\gamma$,
and $a$. This proves (19.13).

\ms
Next we take care of little rims. Set
$$
R(j,x) = \Big\{ z\in B_{j,x}^- \, ; \, 
\dist(z,\d B_{j,x}^-) \leq {(1-a) r_j \over 100(1+|f|_{lip})}
\Big\}
\leqno (19.28)
$$
for $j\in J_4$ and $x\in Z(y_j)$. By (15.17), (15.20),
and the proof of (18.32), we get that
$$
\H^d(E \cap \overline R(j,x))
\leq C(f) (1-a) \gamma^{1-d} r_j^d
\leq C(f,\gamma) (1-a) \H^d(E \cap B(x,|f|_{lip}^{-1} r_j)).
\leqno (19.29)
$$
The total contribution of these annuli to the right-hand side
of (18.1) is still small, because (as happened near (18.34))
$$\leqalignno{
\sum_{j\in J_4} \sum_{x\in Z(y_j)} \H^d(h_3(E_k \cap R(j,x)))
&\leq C \sum_{j\in J_4} \sum_{x\in Z(y_j)} \H^d(E_k \cap R(j,x))
\cr&
\leq \eta + C \sum_{j\in J_4} \sum_{x\in Z(y_j)} \H^d(E \cap \overline R(j,x))
\cr&
\leq \eta + C(f,\gamma) (1-a) \sum_{j\in J_4} \sum_{x\in Z(y_j)}
\H^d(E \cap B(x,|f|_{lip} r_j))
&(19.30)
\cr&
\leq \eta + C(f,\gamma) (1-a) H^d(E \cap W_f)
= \eta + C(f,\gamma) (1-a)
}
$$
by (19.12), for $k$ large and by (10.14), by (19.29), and because
the $B(x,|f|_{lip}^{-1} r_j)$ are disjoint by (18.33) and contained 
in $W_f$ by (18.31).

We are left with $B_{j,x}^- \sm R(j,x)$. Observe that 
$\beta_{j,x}(z,3) = 1$ for $z \in B_{j,x}^- \sm R(j,x)$,
by (19.4), so (19.3) and the first part of (18.70) yield
$$
h_3(z) = \pi_j(h_2(z)) \in Q_j \cap D_j
\hbox{ for } z \in E_k \cap B_{j,x}^- \sm R(j,x),
\leqno (19.31)
$$
at least for $k$ large. 
Again all the sets $Q_j \cap D_j = A_x(P_x) \cap D_j$, $x\in Z(y_j)$, 
coincide, and now
$$\eqalign{
\H^d\Big(h_3 \Big(E_k &\cap \bigcup_{j\in J_4}\bigcup_{x\in Z(y_j)} B_{j,x}^-\Big)\Big)
\cr&\leq \sum_{j\in J_4} \sum_{x\in Z(y_j)} \H^d(h_3(E_k \cap R(j,x)))
+ \sum_{j\in J_4} \H^d(Q_j \cap D_j)
\cr&
\leq \omega_d \sum_{j\in J_4} r_j + \eta + C(f,\gamma) (1-a). 
}\leqno (19.32)
$$
Thus the contribution of all the sets $B_{j,x}^-$ where we modified $h_2$
is just as good as in the rigid case, and we shall only need to worry
about the contribution of the indices $j\in J_3 \sm J_4$.

\msi{\bf We get rid of some small set in $Y_{11}$.}
Let us introduce a small bad set $Z_0 \i Y_{11}$.
For each $y\in U$, denote by $F(y)$ the smallest face of
our grid on $U$ that contains $y$. 
Also set $\wt y = \psi(\lambda y)$ and call
$\wt F(y) = \psi(\lambda F(y))$ the smallest
face of the usual dyadic grid that contains $\wt y$.
Finally call $\wt W(y)$ the smallest affine space that
contains $\wt F(y)$. Then set
$$\eqalign{
A_r(y) &= r^{-d} \sup\big\{ 
\H^d(B(y,r) \cap \lambda^{-1}\psi^{-1}(\wt Q)) \, ; \,
\wt Q \hbox{ is a $d$-dimensional}
\cr&\hskip 3.5cm
\hbox{affine subspace of $\wt W(y)$ that contains $\wt y$} \big\};
}\leqno (19.33)
$$
when the dimension of $\wt W(y)$ is less than $d$, just set
$A_r(y) = 0$. For $0 \leq i \leq j_{max}$, set
$L'_i = L_i \sm {\rm int}(L_i)$, where ${\rm int}(L_i)$
is really the interior of $L_i$, taken in $\R^n$ and regardless
of the dimension of the faces that compose it, and then set
$$
\wh L = \bigcup_{0 \leq i \leq j_{max}} L_i
\ \hbox{ and }\ 
\wh L' = \bigcup_{0 \leq i \leq j_{max}} L'_i.
\leqno (19.34)
$$
Still denote by $\omega_d$ the $d$-dimensional Hausdorff measure 
of the unit ball in $\R^d$. Set
$$
Z = \big\{ y \in \wh L' \, ; \, 
\limsup_{r \to 0} A_r(y) > \omega_d \big\}.
\leqno (19.35)
$$

We chose this definition because it will be easy to use,
and we chose to use the condition (10.7) because it is
not too complicated, and because it implies that
$$
H^d(Z) = 0.
\leqno (19.36)
$$
Let us check this. Let us rather use the translation of (10.7)
that is given below (10.7) itself. This condition gives an exceptional
set $Z_0$ such that $\H^d(Z_0) = 0$ and, if $y \in U \sm Z_0$
lies in $L'_i=L_i \sm {\rm int}(L_i)$ and is such that 
${\rm dimension}(F(y)) > d$, then we can find $t = t(y)>0$ such that 
the restriction of $\psi$ to $\lambda F(y) \cap B(\lambda y,t)$ is $C^1$. 

We want to show that $\H^d$-almost every $y\in Z$ lies in $Z_0$.
Let $y\in Z$ be given; then $y\in L'_i$ for some $i$,
there is a face of $L_i$ that contains $y$, and this face contains $F(y)$ 
by definition of $F(y)$ as a smallest face.

If ${\rm dimension}(F(y)) > d$ and $y\in Z \sm Z_0$,
we can find $t = t(y)>0$ such that the restriction of $\psi$ 
to $\lambda F(y) \cap B(\lambda y,t)$ is $C^1$.
Since $\psi$ is Lipschitz, this also means that
the restriction of $\psi^{-1}$ to the face 
$\wt F(y) = \psi(\lambda F(y))$ is $C^1$ in a neighborhood
of $\wt y = \psi(\lambda y)$. Recall that $y$ is an interior
point of $F(y)$, so $\wt F(y)$ coincides with $\wt W(y)$
(the affine affine space spanned by $\wt F(y)$) near $\wt y$.

With the notation above, if $\wt Q$ is a $d$-dimensional
affine subspace of $\wt W(y)$, the restriction of $\psi^{-1}$ to $\wt Q$
is also $C^1$ near $\wt y$, with uniform estimates with respect to
$\wt Q$. Then $\lambda^{-1}\psi^{-1}(\wt Q)$ is a $C^1$ surface near $y$,
and $\lim_{r \to 0} r^{-d} \H^d(B(y,r) \cap \lambda^{-1}\psi^{-1}(\wt Q))
= \omega_d$, uniformly in $\wt Q$. Thus
$\limsup_{r \to 0} A_r(y) \leq \omega_d$, which contradicts the
fact that $y\in Z$ and takes care of the case when 
${\rm dimension}(F(y)) > d$.

If ${\rm dimension}(F(y)) < d$, then by definition $A_r(y) = 0$ 
for $r > 0$ small, and $y\notin Z$ (a contradiction).

If ${\rm dimension}(F(y)) = d$, there is only one possible
choice of $\wt Q$ in the definition (19.33) of $A_r(y)$,
namely $\wt W(y)$, and 
$$\eqalign{
A_r(y) &= r^{-d} \H^d(B(y,r) \cap \lambda^{-1}\psi^{-1}(\wt W(y))
= r^{-d} \H^d(B(y,r) \cap \lambda^{-1}\psi^{-1}(\wt F(y))
\cr&
= r^{-d} \H^d(B(y,r) \cap F(y))
}\leqno (19.37)
$$
for $r$ small, because $y$ is an interior point of $F(y)$
(and hence $\wt y$ is an interior point of $\wt F(y)$).
But for each face $F$ of dimension $d$ and $\H^d$-almost-every
interior point $y\in F$, $\lim_{r \to 0} r^{-d} \H^d(B(y,r) \cap F(y)) = 
\omega_d$ (because $F$ is rectifiable), so 
$\H^d({\rm int(F)} \cap Z) = 0$. This takes care of
the case when ${\rm dimension}(F(y)) = d$.
This was our last case, and (19.36) follows.

\ms 
We now assume (19.36) (and the other assumptions of Theorem 10.8,
except perhaps (10.7)) and show that $E$ is quasiminimal as in 
(10.9). We proceed as in the last sections, with only two modifications.
The first one occurs in Step 4 (in Section 15), and we shall explain it 
now. The second one is the one that was described earlier in this 
section, and concerns the indices $j\in J_4$.

So we do not change anything up to Section 15; we also define
$Y_9$, $Y_{10}$, and $Y_{11}$ as before, but before we cover
$Y_{11}$ by disks $D_j$ (near (15.12)), we remove some small pieces.

First set $Y_{12} = Y_{11} \sm Z$. Then $\H^d(Y_{11} \sm Y_{12}) = 0$
by (19.36). Set $X_{12} = X_{11} \cap f^{-1}(Y_{12})$;
the same proof as for (15.11)
(or, more precisely, for (4.77) in [D2])  
yields that
$$
\H^d(X_{11} \sm X_{12}) = 0.
\leqno (19.38)
$$
We shall remember that by (19.35),
$$
\limsup_{r \to 0} A_r(y) \leq \omega_d
\ \hbox{ when $y \in Y_{12} \cap \wh L'$.}
\leqno (19.39)
$$
Let $\delta_{9} > 0$ be small, set
$$\eqalign{
Y_{13} = Y_{13}(\delta_9) = [Y_{12} \sm \wh L'] 
\cup \big\{ y\in Y_{12} \cap \wh L' \, ; \,
A_r(y) \leq \omega_d + \varepsilon \hbox{ for }
0 < r \leq \delta_9 \big\}
}\leqno (19.40)
$$
and then
$$
X_{13} = X_{13}(\delta_9) = X_{12} \cap f^{-1}(Y_{13}).
\leqno (19.41)
$$
Notice that $Y_{12}$ is, by (19.39), the monotone union of
the sets $Y_{13}(\delta_9)$, so $X_{12}$ is the monotone union of
the sets $X_{13}(\delta_9)$. Thus
we can choose $\delta_9 > 0$ so small that
$$
\H^d(X_{12} \sm X_{13}) \leq \eta.
\leqno (19.42)
$$
We choose $\delta_9 > 0$ like this, and then 
cover $Y_{13}$ as we did before (for $Y_{11}$)
by  balls $D_j = B(y_j,r_j)$, $j\in J_3$, so that 
$$
y_j \in Y_{13} \ \hbox{ and } \ 
0 < r_j < \Min(\delta_8,\delta_9) 
\ \hbox{ for } j\in J_3, 
\leqno (19.43)
$$
$$
\hbox{ the $\overline D_j$, $j\in J_3$, are disjoint}
\leqno (19.44)
$$
and
$$
\H^d\big(X_{13} \sm f^{-1}\big(\bigcup_{j\in J_3} \overline D_j \big)\big)
\leq \eta.
\leqno (19.45)
$$ 

Then we continue the construction as before, all the way through
Section 18, and arrive to the second modification. 
We define $J_4 \i J_3$ as we did near (19.2) and 
we continue the deformation all the way to $h_3$,
just as was explained at the beginning of this section.

We follow the proof of (18.71), but restrict to the indices
$j\in J_3\sm J_4$; we get that
$$\eqalign{
\H^d\Big(h_3\Big(E_k \cap \bigcup_{j\in J_3\sm J_4} \,\bigcup_{x\in Z(y_j)} 
B_{j,x}^+ \sm R^3\Big)\Big)
&
\leq \sum_{j\in J_3 \sm J_4} \H^d(D_j \cap \lambda^{-1}\psi^{-1}(\wt Q_j))
}\leqno (19.46)
$$
(recall that $h_3=h_2$ on these sets).
Now let $j\in J_3 \sm J_4$ be given. By definition of $J_4$, we can find
$i \in [0,j_{max}]$ such that $y_j \in L_i$, without lying on the
interior of $L_i$. That is, $y_j \in L'_i \i \wh L'$
(see (19.34) and the definition above it). Notice that $y_j \in Y_{13}$
because of our first modification.
Since $y_j \in \wh L'$, (19.40) implies that 
$$
A_r(y_j) \leq \omega_d + \varepsilon \ \hbox{ for }
0 < r \leq \delta_9 
\leqno (19.47)
$$
But $D_j = B(y_j,r_j)$ and $r_j < \delta_9$
by (19.43), so 
$$
r_j^{-d} \H^d(D_j \cap \lambda^{-1}\psi^{-1}(\wt Q_j))
\leq \omega_d + \varepsilon
\leqno (19.48)
$$
by (19.47), the definition (19.33), and because $\wt Q_j$ contains 
$y_j$ and is a $d$-dimensional subspace of the affine space spanned
by $\wt F(y_j)$ (see the discussion below (18.65)). We replace in
(19.46) and get that
$$\eqalign{
\H^d\Big(h_3\Big(E_k \cap \bigcup_{j\in J_3\sm J_4} 
\, \bigcup_{x\in Z(y_j)} B_{j,x}^+ \sm R^3\Big)\Big)
&
\leq \sum_{j\in J_3 \sm J_4} (\omega_d +\varepsilon) r_j^d.
}\leqno (19.49)
$$
Then we add this to (19.32) and get that
$$
\H^d\Big(h_3\Big(E_k \cap \bigcup_{j\in J_3} \bigcup_{x\in Z(y_j)} 
B_{j,x}^+ \sm R^3\Big)\Big)
\leq \sum_{j\in J_3} (\omega_d +\varepsilon) r_j^d
+  \eta + C(f,\gamma) (1-a)
\leqno (19.50)
$$
because $B_{j,x}^+ \sm R^3 \i B_{j,x}^-$ (see (18.26)).
The last part is an error term smaller than  
$\cal E$ in (18.93) and (18.94). We also have the small term
$$
\varepsilon \sum_{j\in J_3} r_j^d
\leq C(f) \varepsilon \sum_{j\in J_3} 
\H^d(E\cap B(x,|f|_{lip}^{-1} r_j))
\leq C(f) \, \varepsilon \, \H^d(E\cap W_f)
\leqno (19.51)
$$
by Proposition 4.1 (which we can apply because of (18.31) and
$W_f \i U$, as for the proof of (18.32)), and then the disjointness (18.33)
and the fact that $B(x,|f|_{lip}^{-1} r_j) \i W_f$ by (18.31).
This term too is dominated by $\cal E$, so (19.50) is essentially as good 
as (18.63) (the difference is controlled by $\cal E$).

We may now continue the proof as before. 
There is a last place where we need to be careful, 
when we use (18.95) to prove set inclusions in (18.97). 
Previously the sets $W$ and $W^2$ were defined in terms of $h_2$, 
and now we need the same inclusions with the sets $W_3$ and $W^2_3$
defined in terms of $h_3$. Fortunately, (19.26) says that
$|h_3(z)-h_2(z)| \leq \delta_5/10$ for $z\in E_k$;
this stays true for $z\in E$ (because $h_3-h_2$ is
continuous and $E$ is the limit of $\{ E_k \}$);
then (18.95) also holds for $h_3$, with the smaller
constant $\delta_5/10$, and we can complete the argument as before
(i.e., $E \sm W_3 \i E \sm {\rm int}(V_1^+)$, and (18.96)-(18.100)
are valid).

This finally completes our proof of Theorem 10.8
in the remaining Lipschitz case.
\qed

\msi{\bf Remark 19.52.} 
Our proof shows that in Theorem 10.8 (and under the Lipschitz
assumption), we can replace the assumption (10.7) with the 
slightly weaker (but more complicated) (19.36).

It is a little sad that the author was not able to get rid
of (10.7) or (19.36) altogether. We seem to be close to that,
but not quite close enough. It would seem natural to try the 
following modification of what we do for $d$-dimensional faces.
Notice that we just need to apply the definition of $A_r(y)$ at points
$y_j$, $j\in J_3$, and to the specific $d$-dimensional set 
$\wt Q_j = \wt A_x(P_x)$, $x\in Z(y_j)$. Modulo some additional
cutting, we could restrict to a subset of $Y_{11}$ where
$\psi(\lambda \cdot)$ coincides with a $C^1$ mapping. Then
we are supposed to go from $y_j$ to $\wt y_j = \psi(\lambda y_j)$,
get $\wt Q_j$, which is also the image by $D\psi(\lambda \cdot)$
of the tangent place to $Y_{11}$ of $f(E)$ at $y_j$, and
show that it has density $1$. For instance, we would know this
for $\wt f(E) = \psi(\lambda f(E))$, which is tangent to
$\wt Q_j$ at $\wt y_j$. But could it be that by bad luck,
$\lambda^{-1}\psi^{-1}(\wt Q_j)$ has more little wrinkles
than $f(E)$, even though they are tangent.

One option to try to  overcome this could be to try to project
points back from $\lambda^{-1}\psi^{-1}(\wt Q_j)$ to $f(E)$, or a 
flatter set, but there are difficulties because we need to do this 
in a Lipschitz way, and more importantly along the faces 
(because of (1.7)), and for instance $f(E)$ may have little holes
(although probably small because $E$ has surjective projections
at places where it is flat). Because all this seems complicated,
the author decided to leave Theorem 10.8 as it is for the moment.

\bigskip
\centerline{PART V : ALMOST MINIMAL SETS AND OTHER THEOREMS ABOUT 
LIMITS}
\ms

In this part we apply the limiting results of the
previous part to sequences of almost minimal, or even minimal sets.
The proofs will usually not be very hard, but this part should be 
useful because it is likely that the results of this paper 
will more often be applied in the almost minimal context.

In Section 20 we give three slightly different definitions of 
sliding almost minimal sets (Definition 20.2), and then show that 
the two last ones are equivalent (Proposition 20.9).
The definitions and proof are inspired of [D5]; 
the point is to unify some of the definitions, and to make it easier 
to check some assumptions. 

In Section 21 we use Theorem 10.8 (our main result about limits)
to show that limits of coral sliding almost minimal sets
(of a given type and with a given gauge function) are also
coral sliding almost minimal sets, of the same type and with the same
gauge function. See Theorem 21.3.

In Section 22 we prove an upper semicontinuity result for $\H^d$,
which says that if $\{ E_k \}$ is a convergent sequence
of coral sliding almost minimal sets in $U$ (as in Theorem 21.3), 
then for each compact set $H \i U$,
$\limsup_{k \to +\infty} \H^d(E_k \cap H) \leq \H^d(E\cap H)$.
See Theorem~22.1. We also prove Lemma 22.3, 
where we only assume that the $E_k$ lie in $GSAQ(U,M,\delta,h)$
and merely get that
$\limsup_{k \to +\infty} \H^d(E_k \cap H) \leq (1+Ch) M \H^d(E\cap H)$.

In Section 23 we consider sequences of almost minimal sets 
$E_k$ that live in domains $U_k$ and with boundary sets $L_{j,k}$
that depend slightly on $k$. We get an analogue of 
Theorems~10.8 and 21.3 that works when $U_k$ and the $L_{j,k}$
are small bilipschitz variations of the limits $U$ and the $L_j$.
See Theorem~23.8, which is proved brutally with a change of variables.

We apply this result in Section 24, to the special case
of blow-up limits. We find two sets of flatness conditions 
on the sets $L_j$ (see Definitions 24.8 and 24.29)
under which the blow-up limits at the origin of a
sliding almost minimal set are sliding minimal sets in $\R^n$,
associated to boundary sets $L_j^0$ obtained from the $L_j$ by the
same blow-up. See Theorem~24.13 and Proposition 24.35.

\msi
{\bf 20. Three notions of almost minimal sets.} 
\ms

We shall more often apply the regularity results above, 
and in particular Theorem~10.8 about limits, in the 
simpler context of almost minimal sets. 

In this section we adopt the same point of view as in
[D5],  
and introduce three types of almost minimal sets;
we shall mostly restrict to the two last ones, which
are slightly weaker, turn out to be equivalent to each other
under mild assumptions, and for which the desired limiting theorem 
will easily follow from Theorem 10.8.
The main point of this section will be the equivalence
between our second and third definitions. It is perhaps
not vital because we can hope to work with a single
definition at a time, but the author will feel much
better for not hiding a little secret under the rug.
Also, the regularity results of the previous sections translate a 
little better in terms of our second definition, while the third
one seems a little simpler.

So we shall give three different definitions of almost minimal sets,
for which we keep the same setting as in Definition 2.3.
That is, we are given an open set $U$ (equal to the unit ball
when we work under the rigid assumption, and to a bilipschitz
image of the unit ball when we work under the Lipschitz assumption), 
and boundaries $L_j$, $0 \leq j \leq j_{max}$. We give a 
special name to $\Omega = L_0$, and require that all our
sets be contained in $\Omega$ (but we can take $\Omega = U$).

Now we also give ourselves a gauge function, 
i.e., a function $h : (0,+\infty) \to [0,+\infty]$
such that
$$
h \hbox{ is continuous from the right and }
\lim_{t \to 0} h(t)  = 0;
\leqno (20.1) 
$$
let us not assume that $h$ is nondecreasing for the moment,
because we don't need this. It would also make sense, in view of
the definition below, to assume that the product $h(r) r^d$ is 
nondecreasing, but let us not do that either.

\ms\proclaim Definition 20.2.
Let $E \i \Omega \cap U$ be a relatively closed in $U$
and such that, as in (1.2),
$$
\H^d(E\cap B) < +\infty
\ \hbox{ for every compact ball $B$ such that $B \i U$.}
\leqno (20.3)
$$
We say that $E$ is an $A_+$-almost minimal set (of dimension $d$)
in $U$, with the sliding conditions given by the closed sets 
$L_j$, $0 \leq j \leq j_{max}$, and the gauge function $h$, 
if for every choice of one-parameter family $\{ \varphi_t \}$,
$0 \leq t \leq 1$, of continuous functions with the properties
(1.4)-(1.8) relative to a ball $B = \overline B(x,r)$, 
and also such that $\widehat W \i \i U$ as in (2.4), we have that
$$
\H^d(W_1) \leq (1+h(r)) \H^d(\varphi_1(W_1)),
\leqno (20.4)
$$
where as usual 
$W_1 = \big\{ y\in E \, ; \, \varphi_1(y) \neq y \big\}$.
\hfill\break
We say that $E$ is an $A$-almost minimal set
(with the sliding conditions given by the closed sets 
$L_j$, $0 \leq j \leq j_{max}$, and the gauge function $h$)
if under the same circumstances,
$$
\H^d(W_1) \leq \H^d(\varphi_1(W_1)) + h(r) r^d.
\leqno (20.5)
$$
Finally, we say that $E$ is an $A'$-almost minimal set
(with the sliding conditions given by the closed sets 
$L_j$, $0 \leq j \leq j_{max}$, and the gauge function $h$)
if under the same circumstances,
$$
\H^d(E \sm \varphi_1(E)) 
\leq \H^d(\varphi_1(E) \sm E)  + h(r) r^d.
\leqno (20.6)
$$

\ms
So the accounting in the three cases is slightly different,
but the competitors are the same (and are the same as for
the generalized quasiminimal sets in Definition 2.3).
We could also have forced the competitors to be such that
$\wh W$, instead of being merely compactly contained in $U$,
is contained in a ball of radius $r$ which itself is compactly 
contained in $U$, and this would probably not have made a big
difference in practice, but we decided to keep the same competitors
as above.

The last two definitions look slightly easier
to use. Let us also check that in Definition~20.2, 
we could replace (20.6) with
$$
\H^d(E \cap \wh W) 
\leq \H^d(\varphi_1(E) \cap \wh W)  + h(r) r^d
\leqno (20.7)
$$
and get an equivalent definition. 
Notice that $\varphi_1(E)$ coincides with
$E$ out of $\wh W$ (by the definition (2.4) of $\wh W$),
and $\H^d(E \cap \wh W) < +\infty$ (by (20.3)
and because $\wh W \i \i U$); 
then (20.7) is obtained from (20.6) by adding
$\H^d(E \cap \varphi_1(E) \cap \wh W)$ to both sides.

We shall now worry about the inclusion relations between our 
three classes of almost minimal sets.

It is fairly easy to see that if $E$ is $A_+$-almost minimal
in $U$, then $E$ is also $A$-almost minimal
in every smaller open set $U_\tau = \big\{ x\in U \, ; \,
\overline B(x,\tau) \i U \big\}$, with the same boundaries
$L_j$, but a slightly larger gauge function $\wt h$
(that depends only on $h$, $\tau$, and $n$ through local
Ahlfors-regularity constants). The proof is the 
same as in Remark 4.5 of [D5], 
and it is fairly easy once you notice that $E$ is quasiminimal,
hence locally Ahlfors-regular. This is also the reason why we
restrict to a smaller set $U_\tau$.
The converse looks like it could be wrong, but the author 
does not know for sure, even in the case without boundary.

It is also easy to see that if $E$ is $A'$-almost minimal in $U$,  
then it is $A$-almost minimal in $U$, with the same $L_j$ and the 
same gauge function $h$. To see this, let $E$ be $A'$-almost minimal, 
let the $\varphi_t$ be as in the definition, and let us deduce (20.5) from
(20.7). Set $Z = E \cap \wh W \sm W_1$ and observe that
$$\eqalign{
\H^d(W_1) &= \H^d(E \cap \wh W) - \H^d(E \cap \wh W \sm W_1)
= \H^d(E \cap \wh W) - \H^d(Z)
\cr&\leq \H^d(\varphi_1(E) \cap \wh W) - \H^d(Z) + h(r) r^d 
\cr& \leq \H^d(\varphi_1(E) \cap \wh W \sm Z) + h(r) r^d
}\leqno (20.8)
$$
because $W_1 \i E \cap \wh W$, by (2.7), and because
$Z \i \varphi_1(E) \cap \wh W$ since $\varphi_1(z)=z$ for $z\in Z$
(by definition of $W_1$). For (2.5) is is enough to check that
$\varphi_1(E) \cap \wh W \sm Z \i \varphi_1(W_1)$. So let
$y\in \varphi_1(E) \cap \wh W \sm Z$ be given, and write 
$y = \varphi_1(x)$. If $x\in W_1$, we are happy. Otherwise, 
$\varphi_1(x) = x$, hence $y = x \in E \cap \wh W \sm W_1 = Z$,
which is impossible. The $A$-minimality of $E$ follows.

Notice that if $E$ is $A$-almost minimal 
(and hence also if $E$ is $A'$-almost minimal),
then it is quasiminimal in every ball $B(x,r) \i U$, with
$M = 1$ and $h = h(r)$. So we shall be able to apply the
regularity results of the previous parts to almost minimal sets.

The fact that $A$-minimality implies $A'$-minimality
will be a little more complicated to prove, and in fact, under
the Lipschitz assumption we shall only be able to do it under
the same additional assumption (10.7) as for Theorem 10.8.
The following is a generalization of 
Proposition 4.10 in [D5]. 

\ms\proclaim Proposition 20.9.
Suppose that the rigid assumption holds, or that
the Lipschitz assumption, plus one of the two technical
conditions (10.7) or (19.36), hold.
Let $E$ be an $A$-almost minimal set in $U$, 
with the sliding conditions given by the closed sets 
$L_j$, $0 \leq j \leq j_{max}$, and the gauge function $h$.
Then $E$ is also $A'$-almost minimal in $U$, 
with the same $L_j$ and the same $h$.
The converse also holds (see above).

\ms
Without of our extra assumption (10.7), we do not
know whether, under the Lipschitz assumption alone,
$A_+$ minimality and $A$-minimality always imply $A'$-minimality.
But we do not have good reasons to think that it fails either.

Our assumption (20.1) should not bother much, but if it fails we can still
do something. The fact that $h(r)$ tends to $0$ when $r$ tends to $0$
will be used only once, at the beginning of the proof in the 
Lipschitz case, to show that $E^\ast$ is rectifiable and 
Ahlfors-regular. If we do not suppose this, we can suppose instead
that $E \in QSAQ(U,M,\delta,h)$ for a number $h> 0$ that is 
small enough (depending on $n$, $M$, and $\Lambda$) for Theorem 5.16
and Propositions~4.1 and 4.74 to apply. If $h$ is not continuous from
the right, our proof will only show that $E$ is $A'$-almost minimal
with the larger gauge function $h'(r) = \liminf_{\varepsilon \to 0^+} 
h(r+\varepsilon)$.

We shall need to revise the proof of [D5], 
because it involves a modification of a family $\{ \varphi_t \}$,
$0 \leq t \leq 1$, and we want to make sure that we do not destroy
the boundary conditions (1.7). Also, the proof in the Lipschitz case
will be a little more complicated, and will use the rectifiability
of $E$, so we shall give two different arguments,
one for the rigid case and one for the Lipschitz case. 
Of course the second argument also works in the rigid case.

\msi
{\bf Proof of Proposition 20.9 under the rigid assumption.}
Let $E$ be $A$-almost minimal; we want to prove that
$E$ is $A'$-almost minimal, so we give ourselves mappings
$\{ \varphi_t \}$, $0 \leq t \leq 1$, that satisfy
(1.4)-(1.8) relative to a ball $B = \overline B(x_0,r_0)$, 
and are also such that $\widehat W \i \i U$.
If $\varphi_1(W_1)$ were disjoint from $E\sm W_1$,
we could easily deduce (20.7) from (20.5): we would 
say that
$$\eqalign{
\H^d(E \cap \wh W) &= \H^d(E \cap \wh W \sm W_1) + \H^d(W_1)
\cr&
\leq \H^d(E \cap \wh W \sm W_1) + \H^d(\varphi_1(W_1)) + h(r_0) r_0^d
\cr&
= \H^d(\varphi_1(E \cap \wh W \sm W_1)) 
+ \H^d(\varphi_1(W_1)) + h(r_0) r_0^d
\cr&
= \H^d(\varphi_1(E \cap \wh W)) + h(r_0) r_0^d,
}\leqno (20.10)
$$
as needed. In general, we want to modify $\varphi_1$
slightly, so as to be able to almost apply the argument
above. And rather than move $\varphi_1(W_1)$, it will be more 
convenient to make $W_1$ artificially larger, by a minor 
modification that will not change (20.7) significantly,
but will make (20.5) more useful.

We first want to construct a vector-valued function $v$, defined on
$U$, which we shall see as a direction in which we are allowed to
move the points. Denote by $\cal F$ the set of faces of dimension
at least $d$ of our usual dyadic grid. For each $F \in {\cal F}$, set
$$\eqalign{
F_\tau &= \big\{ x\in F \, ; \, \dist(x,\d F) \geq \tau \} 
\ \hbox{ and }
\cr
F_\tau^+ &= \big\{ x\in \R^n \, ; \, 
\dist(x,F_{\tau}) \leq {\tau \over 10} \},
}\leqno (20.11)
$$
where the very small $\tau > 0$ will be chosen later,
and then set
$$\eqalign{
h_F(x) &= 1-10\tau^{-1} \dist(x,F_{\tau})
\ \hbox{ for } x \in F_\tau^+,
\cr
h_F(x) &= 0 \hskip2.8cm \hbox{ for } x \in \R^n \sm F_\tau^+.
}\leqno (20.12)
$$
Finally choose for each $F \in {\cal F}$
a vector $v_F$ in the vector space ${\rm Vect}(F)$ 
parallel to $F$, and such that ${1 \over 2} \leq |v_F| \leq 1$, 
and set
$$
v(x) = \sum_{F \in {\cal F}} h_F(x) \, v_F
\ \hbox{ for } x\in \R^n.
\leqno (20.13)
$$
Recall from (3.8) that if $F, G \in {\cal F}$ are
different faces, with $\dim(F) \geq \dim(G)$, then
$$
\dist(y,G) \geq \dist(y,\d F)
\ \hbox{ for } y\in F.
\leqno (20.14)
$$
In particular, $\dist(y,G) \geq \tau$ if $y\in F_\tau$,
and hence
$$
\dist(F_\tau^+,G_\tau^+) \geq {8 \tau \over 10}.
\leqno (20.15)
$$
Thus the sum in (20.13) has at most one term, and when
we compute the differential of $v$ term by term, we get
that
$$
v \hbox{ is $10\tau^{-1}$-Lipschitz.}
\leqno (20.16)
$$
Also, (20.12), (20.13), and (20.15) yield
$$
v(x) = v_F \ \hbox{ for } x\in F_\tau.
\leqno (20.17)
$$

We also need a cut-off function $\chi$. First select compact subsets
$S$ and $S'$ of $U$, such that
$$
\wh W \i S \i {\rm int}(S') \i S' \i U
\leqno (20.18)
$$
and $S' \i B(x_0,r_0+\tau)$, where $B = \overline B(x_0,r_0)$ 
is as in (1.4)-(1.8). Let us also make sure, for instance by replacing
$S'$ with a smaller compact set, that
$$
\H^d(E \cap S' \sm S) \leq \varepsilon,
\leqno (20.19)
$$ 
where the small number $\varepsilon > 0$ is chosen in advance.
Then choose a Lipschitz function $\chi$ on $U$, 
so that
$$
0 \leq \chi(x) \leq 1
\ \hbox{ everywhere, }
\chi(x) = 1 \hbox{ on $S$, and }
\chi(x) = 0 \hbox{ on } U \sm S'.
\leqno (20.20)
$$
We shall select an extremely small $t_0 > 0$ and continue 
the one parameter family $\{ \varphi_t \}$ with mappings 
$\varphi_t$, $1 \leq t \leq 1+t_0$, defined by
$$
\varphi_t(x) = \psi_t(\varphi_1(x)),
\ \hbox{ where }\ 
\psi_t(y) = y + (t-1) \chi(y) v(y).
\leqno (20.21)
$$
Our constant $t_0$ will be chosen last, depending
on $\varphi_1$, $\tau$, $\varepsilon$, $S$, $S'$, and even $\chi$
if needed, so small that $\psi_{t}$ is $2$-Lipschitz
for $1 \leq t \leq 1+t_0$, and hence
$$
\varphi_{1+t_0} \hbox{ is $2\, ||\varphi_1||_{lip}$-Lipschitz
on $E$.}
\leqno (20.22)
$$

Next we want to check that the $\varphi_{(1+t_0)t}$, 
$0 \leq t \leq 1$, satisfy the required conditions (1.4)-(1.8), 
but this time with respect to the slightly larger ball 
$B' = \overline B(x_0,r_0+\tau+t_0)$. 
We still have (1.4) and (1.8) because $v$ and $\chi$
are Lipschitz. For (1.5), we just need to worry about
$t > 1$. Observe that when $x \in E \sm B'$, 
$\chi(\varphi_1(x)) = \chi(x) = 0$ and hence
$\varphi_t(x) = \varphi_1(x)$ for $t \geq 1$;
thus (1.5) holds.
For (1.6), let $x\in E \cap B'$ be given. If
$x\in B$, (1.6) for our initial $\varphi_t$ says that
$\varphi_1(x) \in B$, and then $\varphi_t(x) \in B'$
for $t \geq 1$, by (20.21). If
$x\in B(x_0,r_0+\tau)\sm B$, then
$\varphi_1(x)=x \in B(x_0,r_0+\tau)\sm B$
by (1.5) for our initial $\varphi_t$, so
$\varphi_t(x) \in B'$ for $t \geq 1$, again by (20.21)
and because $t-1 \leq t_0$. Finally, if
$x\in B' \sm B(x_0,r_0+\tau)$, we still have
that $\varphi_1(x)=x$, and now $\chi(x) = 0$
and so $\varphi_t(x) = x$ for $t > 1$.
So (1.6) holds.

We are left with (1.7) to check. Let $j \leq j_{max}$
and $x\in E \cap L_j \cap B'$ be given; we want to check that
$\varphi_t(x) \in L_j$ for $t \geq 1$ (we already know
about $t \leq 1$, by assumption). Set $y = \varphi_1(x)$;
thus
$$
\varphi_t(x) = \psi_t(y) 
= y + (t-1) \chi(y) v(y)
\leqno (20.23)
$$
by (20.21). If $v(y) = 0$, then $\varphi_t(x) = \varphi_1(x)\in L_j$.
Otherwise, $y\in F_\tau^+$ for some $F\in {\cal F}$,
and 
$$
\varphi_t(x) = y + (t-1) \chi(y) h_F(y) v_F
\leqno (20.24)
$$
by (20.12) and (20.13) (also recall that the $F_\tau^+$
are disjoint by (20.15)). 
Let $G$ denote the smallest face of our grid that contains $y$.
We claim that $G$ contains $F$. Let $z\in F_\tau$ be such
that $|z-y| \leq \tau/10$; if $G$ does not contain $F$,
(3.8) applies and says that $\dist(z,G) \geq \dist(z, \d F)
\geq \tau$, a contradiction since $y\in G$. So $F \i G$. 

Next we claim that $\varphi_t(x) \in G$, at least if we take 
$t_0 < \tau/10$. Denote by $z'$ the orthogonal projection of 
$y$ onto the smallest affine space $W$ that contains $F$; then
$|z'-y| \leq |z-y| \leq \tau/10$ (because $z\in W$), so 
$\dist([z,z'],\d F) \geq \dist(z,\d F) - |z'-y|
\geq \dist(z, \d F) - 2\tau/10 
\geq 8\tau/10$ (because $z\in F_\tau$), and so $z'$
lies in the interior of $F$.

By (20.24), $\varphi_t(x) = y + \lambda v_F$, with 
$|\lambda| \leq t-1 \leq t_0 < \tau/10$ and $v_F \in {\rm Vect}(F)$. 
Let us compute with coordinates. The face $F$ is given by some
equations $z_i = a_i$, where the $z_i$ are coordinates of the current
point $z$, and the $a_i \in 2^{-m}\Bbb Z$ are constants, and some
inequalities $z_j \in I_j$, where each $I_j$ is a dyadic interval of
size $2^{-m}$. When we replace $y$ with 
$\varphi_t(x) = y + \lambda v_F$, we only modify some of the 
$z_j$, but since $\dist(z',\d F) \geq 8 \tau /10$, the corresponding 
coordinates stay in the interior of corresponding $I_j$. The other 
coordinates $z_i$ stay whatever they were, and altogether 
$\varphi_t(x)$ lies in exactly the same faces that contain $y$.
Since $y\in G$, we get that $\varphi_t(x) \in G$.
By definition of $G$ as the the smallest face that 
contains $y$, $G \i L_j$ because $y = \varphi_1(x)\in L_j$
(by (1.7) for $\varphi_1$). Hence $\varphi_t(x) \in L_j$,
as needed for (1.7).

We also need to check the assumption (2.4) for our extended
family. Let $x\in E$ and $t \in [0,1+t_0]$ be such that
$\varphi_t(x) \neq x$. If $t \leq 1$, we know (by assumption)
that $x\in W_t$ and $\varphi_t(x) \in \wh W$, a compact subset 
of $U$. So we may assume that $t > 1$, and also that
$\varphi_t(x) \neq \varphi_1(x)$. Set $y = \varphi_1(x)$;
since $\varphi_t(x) = y + (t-1) \chi(y) v(y)$ by (20.23),
we get that $\chi(y) \neq 0$, hence $y\in S'$ by (20.20).
Thus $\dist(\varphi_t(x),S') \leq t-1 \leq t_0$ and, if
$t_0$ is small enough, this forces $\varphi_t(x)$ to stay in a 
(fixed) compact subset of $U$. Also, either $x=y$,
and then $x\in S'$, or else $x \neq y = \varphi_1(x)$, 
hence $x\in W_1 \i S'$ too, so $W_t \i S'$ for $t \geq 1$. 
Thus (2.4) holds.

We may now use our assumption that $E$ is $A$-minimal. Set
$$
\varphi = \varphi_{1+t_0} \,\hbox{ and }\,
W = W_{1+t_0} = \big\{ x\in E \, ; \, \varphi(x) \neq x \big\};
\leqno (20.25)
$$
then
$$
\H^d(W) \leq \H^d(\varphi(W)) + h(r_0+\tau+t_0) (r_0+\tau+t_0)^d
= H^d(\varphi(W)) + h(r_1) r_1^d
\leqno (20.26)
$$
by (20.5) for the extended family, and with 
$r_1 = r_0+\tau+t_0$.

We want to say that $W$ is large. First  observe that
$$
\big\{ x\in E \, ; \, |\varphi_1(x) - x| > t_0 \big\}
\i W,
\leqno (20.27)
$$
just because $|\varphi(x)-\varphi_1(x)| 
= |\psi_{1+t_0}(\varphi_1(x))-\varphi_1(x)| \leq t_0$
by (20.21). Set
$$
A_\tau = \bigcup_{F \in {\cal F}} F_\tau,
\leqno (20.28)
$$
and notice that by (20.17), $v(x) \neq 0$ on $A_\tau$; then
$$
S \cap A_\tau \cap E \sm W_1 \i W
\leqno (20.29)
$$
because if $x\in S \cap A_\tau \cap E \sm W_1$, then
$\varphi_1(x) = x \in S \cap A_\tau$ and hence
$$
\varphi(x) 
= \varphi_{1+t_0}(x) = \psi_{1+t_0}(x)
= x + t_0 \chi(x) v(x) \neq x
\leqno (20.30)
$$
by (20.25), (20.21), because $\xi(x) = 1$ by (20.20), 
and because $v(x) \neq 0$.

Recall that we want to prove that $E$ is $A'$-almost minimal,
so we want to establish (20.7), i.e., estimate
$\H^d(E \cap \wh W)$, where
$$
\wh W  = \bigcup_{0< t \leq 1} W_t \cup \varphi_t(W_t)
\leqno (20.31)
$$
is as in (2.2). But it will be more convenient to work
with the compact set $S$ of (20.18), and estimate
$\H^d(E \cap S)$; we shall see that it makes no difference
for (2.7). We write 
$S = (S\cap W) \cup (S \sm W)$, and so
$$
\H^d(E \cap S) \leq \H^d(W) + \H^d(E \cap S \sm W)
\leq \H^d(\varphi(W)) + h(r_1) r_1^d
+ \H^d(E \cap S \sm W)
\leqno (20.32)
$$
by (20.26). Next we estimate $\H^d(E \cap S \sm W)$.
Set 
$$
Z_\tau = \R^n  \sm A_\tau
= \R^n \sm \bigcup_{F\in {\cal F}} F_\tau
\leqno (20.33)
$$
(by (20.28)). By the definition (20.11), every interior point of a 
face of dimension $\geq d$ lies in $F_\tau$ for $\tau$ small, so
$$
\bigcap_{\tau > 0} Z_\tau = {\cal S}_{d-1},
\leqno (20.34)
$$
where ${\cal S}_{d-1}$ still denotes the union of the
faces of dimension $d-1$ of our net; since the
intersection is decreasing and $\H^d(E\cap S) < +\infty$
(because $S \i \i U$), we get that
$$
\H^d(E\cap S \sm A_\tau) 
= \H^d(E\cap S \cap Z_\tau) \leq \varepsilon
\leqno (20.35)
$$
if $\tau$ is chosen small enough, and where $\varepsilon>0$
is the same small number given in advance as in (20.19).

Similarly, $W_1 = \big\{ x\in E \, ; \, \varphi_1(x) \neq x \big\}$
is the monotone union of the sets 
$\big\{ x\in E \, ; \, |\varphi_1(x) - x| > t_0 \big\}$
that show up in (20.27), so (20.27) says that if $t_0$
is small enough,
$$
\H^d(W_1 \sm W) \leq \varepsilon
\leqno (20.36)
$$
(again, this holds because $W_1 \i E \cap \wh W \i E\cap S$
and hence $\H^d (W_1) < +\infty$). Then
$$\eqalign{
\H^d(E \cap S \sm W)
&\leq \H^d(W_1 \sm W) + \H^d(E \cap S \sm (W \cup W_1))
\cr&
\leq \varepsilon + \H^d(E \cap S \sm (W \cup W_1))
\cr&
= \varepsilon + \H^d(E \cap S \sm (W \cup W_1 \cup A_\tau))
\cr&
\leq \varepsilon + \H^d(E \cap S \sm A_\tau) \leq 2 \varepsilon
}\leqno (20.37)
$$
by (20.36), because
$E \cap S \cap A_\tau \sm (W \cup W_1) = \emptyset$
by (20.29), and by (20.35). So (20.32) yields
$$
\H^d(E \cap S) 
\leq \H^d(\varphi(W)) + h(r_1) r_1^d + 2 \varepsilon
\leqno (20.38)
$$
and our next step is to estimate $\H^d(\varphi(W))$.
Recall from (20.25) and (20.21) that 
$\varphi = \psi_{1+t_0} \circ \varphi_1$, so
$$
\varphi(W) = \psi_{1+t_0}(\varphi_1(W)).
\leqno (20.39)
$$
We first consider $\psi_{1+t_0}(\varphi_1(W) \sm S)$.
Let $x\in W$ be such that $\varphi_1(x)$ lies outside of $S$; 
then $\varphi_1(x) = x$, because otherwise $x\in W_1$ and
$\varphi_1(x) \in \varphi_1(W_1) \i \wh W \i S$ by (20.31)
and (20.18). In addition, $x\in S'$ because 
otherwise $\chi(x) = 0$ by (20.20) and
$\varphi(x) = \psi_{1+t_0}(x) = x$ by (20.21); 
this is impossible because $x\in W$. So
$x\in E \cap S' \sm W_1$, and even 
$x\in E \cap S' \sm S$ because $\varphi_1(x) = x$ 
and we assumed that $\varphi_1(x)$ lies outside of $S$.
Hence $\psi_{1+t_0}(\varphi(x)) 
= \psi_{1+t_0}(x) \in \psi_{1+t_0}(E \cap S' \sm S)$. 
We just checked that $\psi_{1+t_0}(\varphi_1(W) \sm S) \i 
\psi_{1+t_0}(E \cap S' \sm S)$, and so
$$
\H^d(\psi_{1+t_0}(\varphi_1(W) \sm S))
\leq \H^d(\psi_{1+t_0}(E \cap S' \sm S))
\leq 2^d \H^d(E \cap S' \sm S)
\leq 2^d\varepsilon
\leqno (20.40)
$$
because $\psi_{1+t_0}$ is $2$-Lipschitz
(see above (20.22)), and by (20.19).

We are left with $\psi_{1+t_0}(\varphi_1(W) \cap S)$.
By (20.34), the monotone intersection
of the sets $\varphi_1(E) \cap S \cap Z_\tau$,
when $\tau$ tends to $0$, is contained in ${\cal S}_{d-1}$.
Since $\H^d(\varphi_1(E) \cap S) < +\infty$,
we get that
$$
\H^d(\varphi_1(E) \cap S \cap Z_\tau)
\leq \varepsilon
\leqno (20.41)
$$
if $\tau$ is chosen small enough (depending on $\varphi_1$).
And then
$$
\H^d(\psi_{1+t_0}(\varphi_1(W) \cap S \cap Z_\tau))
\leq 2^d \H^d(\varphi_1(W) \cap S \cap Z_\tau) 
\leq 2^d \varepsilon
\leqno (20.42)
$$
because $W \i E$ and $\psi_{1+t_0}$ is $2$-Lipschitz.
We are now left with 
$\psi_{1+t_0}(\varphi_1(W) \cap S \cap A_\tau)$.
Write
$$
\varphi_1(W) \cap S \cap A_\tau = 
\bigcup_{F \in {\cal F}} G_F,
\hbox{ with } G_F = \varphi_1(W) \cap S \cap F_\tau,
\leqno (20.43)
$$
and observe that for $y\in S \cap F_\tau$,
$\chi(y) = 1$ by (20.20), and
$$
\psi_{1+t_0}(y) = y + t_0 v(y) = y + t_0 v_F
\leqno (20.44)
$$
by (20.21) and (20.17). Hence
$$
\psi_{1+t_0}(\varphi_1(W) \cap S \cap A_\tau)
= \bigcup_{F \in {\cal F}} \psi_{1+t_0}(G_F)
= \bigcup_{F \in {\cal F}} [G_F+t_0 v_F]
\leqno (20.45)
$$
by (20.43), and 
$$\eqalign{
\H^d(\psi_{1+t_0}(\varphi_1(W) \cap S \cap A_\tau))
&\leq \sum_{F \in {\cal F}} \H^d(G_F+t_0 v_F)
\cr&
= \sum_{F \in {\cal F}} \H^d(G_F)
\leq \H^d(\varphi_1(W)))
}\leqno (20.46)
$$
because the $G_F$ are disjoint (by (20.15)) and contained in 
$\varphi_1(W)$. Altogether,
$$
\H^d(\varphi(W)) 
= \H^d(\psi_{1+t_0}(\varphi_1(W)))
\leq \H^d(\varphi_1(W)) + 2^{d+1} \varepsilon
\leqno (20.47)
$$
by (20.39), (20.40), (20.42), and (20.46). Also  recall that if
$x\in W \sm S$, then $\varphi_1(x) = x$
(because (20.18) says that $W_1 \i \wh W \i S$; also see the 
definition of $W_1$ below (20.4)) and $x\in S'$ (because
otherwise $\varphi(x) = \psi_{1+t_0} \circ \varphi_1(x)
= \psi_{1+t_0}(x) = x$ by (20.21) and because $\chi(x) = 0$
by (20.20)); hence
$$\eqalign{
\H^d(\varphi_1(W)) & \leq \H^d(\varphi_1(W \cap S)) + \H^d(\varphi_1(W \sm S))
\cr&
\leq \H^d(\varphi_1(E \cap S)) + \H^d(W \sm S)
\cr&
\leq \H^d(\varphi_1(E \cap S)) +  \H^d(E \cap S' \sm S)
\cr&
\leq \H^d(\varphi_1(E \cap S)) + \varepsilon
}\leqno (20.48)
$$
because we just saw that $\varphi_1(x) = x$ on $W \sm S$, then because
$W \i E \cap S'$ (see the definition (20.25) and recall that on $U \sm S'$, 
$\varphi_1(x) = x$ by (20.18) and (20.31) and hence 
$\varphi_{1+t_0}(x) = x$ by (20.21) and (20.20)),
and finally by (20.19). Hence
$$
\H^d(E \cap S) 
\leq \H^d(\varphi(W)) + h(r_1) r_1^d + 2 \varepsilon
\leq \H^d(\varphi_1(E \cap S)) 
+ h(r_1) r_1^d + (2^{d+1}+3) \varepsilon
\leqno (20.49)
$$
by (20.38), (20.47), and (20.48).

Recall that $\varphi_1(x) = x$ for $x\in E \sm S$
and $\varphi_1(E\cap S) \i S$ (because $\wh W \i S$);
then $E$ and $\varphi_1(E)$ coincide out of $S$, and so
$$
E \sm \varphi_1(E) = S \cap E \sm \varphi_1(E)
\ \hbox{ and } \ 
\varphi_1(E) \sm E = S \cap \varphi_1(E) \sm E.
\leqno (20.50)
$$
Since both sets have a finite measure and contain 
$E\cap \varphi_1(E) \cap S$, we get that
$$\eqalign{
\H^d(E \sm \varphi_1(E))-\H^d(\varphi_1(E) \sm E)
&= \H^d(S \cap E \sm \varphi_1(E))-\H^d(S \cap \varphi_1(E) \sm E)
\cr&
= \H^d(S \cap E) -\H^d(S \cap \varphi_1(E)) 
}\leqno (20.51)
$$
by subtracting $\H^d(E\cap \varphi_1(E) \cap S)$
from both terms. In addition, $S \cap \varphi_1(E)
= \varphi_1(E \cap S)$ because $\varphi_1(E\cap S) \i S$
and $\varphi_1(x) = x \notin S$ for $x \notin S$.
Now (20.51) and (20.49) yield
$$
\H^d(E \sm \varphi_1(E))-\H^d(\varphi_1(E) \sm E)
\leq  h(r_1) r_1^d + (2^{d+1}+3) \varepsilon.
\leqno (20.52)
$$
Recall from the line below (20.26) that $r_1 = r_0+\tau+t_0$, which
is as close to $r_0$ as we want; since $h$ is continuous from the right,
$h(r_1)$ is as close to $h(r_0)$ as we want. 
Also, $\varepsilon$  is as small as we want too, and since (20.52)
holds with all these choices, we get (20.6). This completes our
proof of Proposition 20.9 in the rigid case.
\qed

\msi{\bf Proof of Proposition 20.9 in the general case.}
The proof that we give below will use the same strategy as in the
rigid case, but will be more complicated because we have a technical
problem. When we modify $\varphi_1$, we move the points a little bit 
along the faces of our grid, and we do this because we want
to preserve the boundary conditions (1.7). If we try to do this 
with our bilipschitz faces, this small translation along the
faces may well multiply $\H^d(\varphi_1(E))$ by a factor of $2$,
even if our translation is very small, and this would of course be bad 
for our estimates. So we will have to find flatter parts of our faces
where we can translate things without increasing the measure too much,
and for this an almost-covering argument with disjoint small balls where
$\varphi_1(E)$ looks nice will be helpful.
The proof below also works in the rigid case, with some 
simplifications.

Let $E$ be $A$-almost minimal, and let the 
$\{ \varphi_t \}$, $0 \leq t \leq 1$, satisfy
(1.4)-(1.8) relative to a ball $B = \overline B(x_0,r_0)$ 
and be such that $\widehat W \i \i U$.

Observe that $E$ is rectifiable because it is quasiminimal. 
More precisely, choose $h> 0$ small, as in Theorem 5.16
with $M = 1$, and then use (20.1) to find $\delta > 0$
such that $h(r) \leq \delta$ for $0 < r \leq \delta$. 
Then $E\in QMAQ(U,M,\delta,h)$, with $M = 1$ (compare
the definitions 20.2 and 2.3). Now Theorem 5.16 says that
$E$ is rectifiable, as needed. Similarly, Propositions~4.1 and 4.74
say that $E^\ast$ (the core of $E$, defined in (3.2)), is locally
Ahlfors-regular.

Let $\varepsilon > 0$ be small, and let 
$S$ be a compact set such that 
$$
\wh W \i {\rm int}(S) \i S \i U, \ 
S \i B(x_0,r_0+\varepsilon), \, \hbox{ and } \, 
\H^d(E \cap S \sm \wh W) < \varepsilon
\leqno (20.53)
$$
(recall from (1.5), (1.6), and (2.2) that 
$\wh W \i \overline B(x_0,r_0)$).
Denote by $\mu$ the restriction of $\H^d$ to the set
$\varphi_1(E\cap S) = \varphi_1(E) \cap S$
(recall that $\varphi_1(E\cap \wh W) \i \wh W$ and
$\varphi_1(x) = x$ for $x\in E \sm \wh W$), and by $\nu$ the image 
by $\varphi_1$ of the restriction of $\H^d$ to 
$E\cap S$, defined by $\nu(A) = \H^d(E \cap S \cap \varphi_1^{-1}(A))$
for $A \i R^n$ (a Borel set). By (20.3), 
$\H^d(E\cap S) < +\infty$, hence $\mu$ and $\nu$ are finite
measures (recall that $\varphi_1$ is Lipschitz).

We want to cover a substantial part of 
$$
G_0 = \varphi_1(E^\ast\cap \wh W) \i \wh W
\leqno (20.54)
$$
(as before, the inclusion comes from the fact that 
$\varphi_1(E\cap W_1) \i \wh W$; see (2.1) and (2.2))
by a collection of disjoint balls $B_j$, and then we will continue our
mapping $\varphi_1$ by composing $\varphi_1$ by deformations defined
on the $B_j$. Our first task is to eliminate various pieces of 
$G_0$. We defined $G_0$ in terms of the core $E^\ast$, because 
it costs nothing in terms of measure (recall from (3.29)
or (8.26) on page 58 of [D4]  
that $\H^d(E \sm E^\ast) = 0$) and $E^\ast$, being locally
Ahlfors-regular, is a little easier to control. 
For instance let us check that
$$
\liminf_{r \to 0} r^{-d}\nu(B(y,r)) > 0
\ \hbox{ for } y\in G_0.
\leqno (20.55)
$$
Let $y\in G_0$ be given, pick $x\in E^\ast \cap \wh W$ such that
$\varphi_1(x) = y$, and for $r > 0$ small, set 
$\rho = (1+|\varphi_1|_{lip})^{-1} r$; then 
$\varphi_1(B(x,\rho)) \i B(y,r)$. 
By (20.53), $B(x,\rho) \i S$ for $\rho$ small, so
$\nu(B(y,r) \geq \H^d(E\cap B(x,\rho) \cap S) \geq C^{-1} \rho^{d}
\geq C^{-1} r^d$, where we don't even want to know
what $C$ depends on. This is enough for (20.55),
which we just mention because it simplifies 
the discussion below.

First we remove the points where $\nu$ is much larger than 
$\mu$. Let $M > 1$ be very large, and set
$$
Y_0 = \big\{ y\in G_0 \, ; \, 
\limsup_{r \to 0} {\nu(B(y,r)) \over \mu(B(y,r))} \geq M \big\}.
\leqno (20.56)
$$
We don't need to worry about the value of $0/0$ here, since
$\nu(B(y,r)) > 0$ by (20.53). Also set
$$
X_0 = E \cap \varphi_1^{-1}(Y_0) = E \cap \wh W \cap \varphi_1^{-1}(Y_0)
\leqno (20.57)
$$
(recall that $Y_0 \i G_0 \i \wh W$ and $\varphi_1(x) = x \notin \wh W$
when $x\in E \sm \wh W$).
By (2) in Lemma 2.13 of [Ma], 
$$
\nu(A) \geq M \mu(A)
\ \hbox{ for every Borel set } A \i Y_0,
\leqno (20.58)
$$
and in particular 
$$
\H^d(X_0) = \nu(Y_0) \geq M\mu(Y_0)
\leqno (20.59)
$$
because $X_0 \i S$ (since $\wh W \i S$ by (20.53)). 
We will not need to worry too much about $Y_0$ because 
$$
\H^d(Y_0) \leq M^{-1} \H^d(X_0)
\leq M^{-1} \H^d(E\cap S) \leq M^{-1} (1+\H^d(E \cap \wh W)),
\leqno (20.60)
$$
by (20.53). 
On the other hand, (1) in Lemma 2.13 of [Ma], 
says that
$$
\nu(A) \leq M \mu(A)
\ \hbox{ for every Borel set } A \i G_0 \sm Y_0;
\leqno (20.61)
$$
we don't intend to use the huge constant $M$, but merely the
fact that $\nu$ is absolutely continuous with respect to $\mu$
on the set $G_1 = G_0 \sm Y_0$. This will help, because we
can now remove some small sets in $G_1$ without fear of losing
a large mass in the source space.

Denote by $G_2$ the set of points $y\in G_1$ with the 
following good properties. First,
$$
\varphi_1(E^\ast) \hbox{ has an approximate tangent plane
$P(y)$ at $y$}
\leqno (20.62)
$$
and 
$$
\lim_{r \to 0} r^{-d} \H^d(\varphi_1(E^\ast) \cap B(y,r)) = \omega_d,
\leqno (20.63)
$$
where as usual $\omega_d$ is the $\H^d$-measure of the unit ball in 
$\R^d$. These properties are true for $H^d$-almost every
$y\in G_0$, because $E^\ast$ is rectifiable (with finite measure
in a neighborhood of $S$) and $\varphi_1$ is Lipschitz. 
We can also replace $E^\ast$ with $E$ in (20.62) and (20.63),
since none of these properties are sensitive to adding a set of
vanishing $H^d$-measure.

Next, if we are in the Lipschitz case, set 
$\wt \varphi_1(x) = \psi(\lambda \varphi_1(x))$ for $x\in U$
(and where $\lambda$ and $\psi$ are as in Definition 2.7),
and also set $\wt y = \psi(\lambda y)$; we require that
$$
\wt \varphi_1(E^\ast) \hbox{ has an approximate tangent plane
$\wt P(y)$ at $\wt y$.}
\leqno (20.64)
$$
In addition, denote by $F(y)$ the smallest face of our 
(twisted) net that contains $y$ and by 
${\rm dim}(F(y))$ its dimension. We demand that
$$
{\rm dim}(F(y)) \geq d,
\leqno (20.65)
$$
and also that if $\wt W(y)$ denotes the smallest affine space that
contains $\wt F(y) = \psi(\lambda F(y))$,
$$
\wt P(y)  \i \wt W(y).
\leqno (20.66)
$$
Finally we exclude the exceptional set $Z$ of (19.35). 
In other words, we demand that if $y$ lies in some boundary piece 
$L_i$, $0 \leq i \leq j_{max}$, but does not lie in its $n$-dimensional 
interior (see the definition of $L'_i$ and $\wh L'$ near (19.34)), then
$$
\limsup_{r \to 0} A_r(y) \leq \omega_d
\leqno (20.67)
$$
where $A_r(y)$ is given by (19.33). Let us check that
all these properties are true for $\H^d$-almost every $y\in G_1$,
i.e., that
$$
\H^d(G_1 \sm G_2) = 0.
\leqno (20.68)
$$
We know that (20.62) and (20.63) hold almost everywhere, and 
so does (20.64), because $\wt \varphi_1(E^\ast)$ is rectifiable
and (for the invariance of negligible sets) $\psi$ is bilipschitz.
For (20.65) we remove a set of dimension $d-1$, and (20.67) holds
almost everywhere because we assumed (19.36) or the stronger (10.7).
The fact that (10.7) implies (19.36) is proved below (19.36).
Finally, let us check that we can arrange (20.66) almost everywhere.
Let $F$ be any face, and let us say how we can get (20.66) for
almost every $y \in G_1$ such that $F(y) = y$. Set 
$A = F \cap G_1$ and $\wt A = \psi(\lambda A)$.
This last set is rectifiable (it is also a subset of 
$\wt \varphi_1(E^\ast)$), so for almost every $y\in A$,
we can find an approximate tangent $d$-plane to $\wt A$ at 
$\wt y = \psi(\lambda y)$. 
Call it $Q(y)$, and observe that by definitions it is contained in the 
smallest affine subspace that contains $\psi(\lambda F)$. 
By the almost-everywhere uniqueness of the approximate tangent plane 
to $\wt\varphi_1(E^\ast)$,
we just have to show that for almost every $y \in A$, 
$Q(y)$ is also an approximate tangent plane to
$\wt\varphi_1(E^\ast)$ (and not merely $\wt A$) at $\wt y$. But by 
Theorem 6.2 on page 89 of [Ma],  
$\lim_{r \to 0} r^{-d}\H^d(\wt\varphi_1(E^\ast) \cap B(\wt y,r)\sm \wt A) = 0$
for $\H^d$-almost every $\wt y \in \wt A$. For such $\wt y$, any 
approximate tangent plane to $\wt\varphi_1(E^\ast)$ at $\wt y$
also works for $\wt\varphi_1(E^\ast)$, as needed.
This completes the proof of (20.68).

Let us now select, for each point $y\in G_2$,
a small radius $r(y)$ with the following good properties. First, 
$$
r(y) \leq {1 \over 4\Lambda^{2}} \dist(y, U \sm S);
\leqno (20.69)
$$
this true for $r(y)$ small enough, because (20.54) and (20.53)
say that $y\in G_0 \i \wh W \i {\rm int}(S)$. We also choose 
$r(y)$ so small that
$$
r(y) \leq {1 \over 4\Lambda^4} \dist(y, \d F(y))
\leqno (20.70)
$$
(where $\dist(y, \d F(y))$, the distance to the
boundary of $F(y)$, is positive because $y$ lies in the interior of 
$F(y))$), 
$$
\omega_d - \varepsilon \leq r^{-d} \H^d(\varphi_1(E^\ast) \cap B(y,r)) 
\leq \omega_d + \varepsilon
\ \hbox{ for } 0 < r \leq r(y)
\leqno (20.71)
$$
(we use the same small $\varepsilon > 0$ as before to save notation).
We shall not need a uniform variant for the existence of a
tangent plane to $\varphi_1(E)$, because in the delicate part of the
argument, we shall work with $\wt\varphi_1(E)$. So
we use (20.64) to require that for $0 < r \leq r(y)$,
$$
\H^d\big(\big\{ z\in \wt\varphi_1(E) \cap B(\wt y,\lambda \Lambda r) \, ; \,
\dist(z, \wt P(y)) \geq \varepsilon \lambda r \big\}\big) 
\leq \varepsilon \lambda^d r^d.
\leqno (20.72)
$$
Finally we require a uniform version of (20.67), i.e., that
if $y$ lies in some boundary piece $L_i$, but not in the (true)
interior of $L_i$,
$$
A_r(y) \leq \omega_d + \varepsilon
\ \hbox{ for } 0 \leq r \leq 2r(y).
\leqno (20.73)
$$
This completes our definition of $r(y)$ when $y\in G_2$. 

We now use a consequence of Besicovitch's covering lemma.
Consider, for each $y\in G_2$, the balls $\overline B(y,r)$, 
$0 < r < \min(\varepsilon, r(y))$ (we are again using
the same $\varepsilon$ in a different role), and for which
$\mu(\d B(y,r)) = 0$ (almost every $r$ satisfies this, 
since the $\d B(y,r)$ are disjoint).
By Theorem 2.8 in [Ma] (applied to all these balls) 
we get a collection of disjoint balls $B_j= B(y_j,r_j)$, 
with the following properties:
$$
0 < r_j < \min(\varepsilon, r(y_j))
\leqno (20.74)
$$ 
and $\mu(\d B_j) = 0$ for all $j$, and 
$$
\mu(G_2 \sm \bigcup_j B_j) = 
\mu(G_2 \sm \bigcup_j \overline B_j) = 0.
\leqno (20.75)
$$

Now we can define a continuation for our family $\{ \varphi_t \}$,
with which we shall eventually apply the definition of $A$-minimality.
We want to define $\varphi_t$ for $1 \leq t \leq 2$, by 
$$
\varphi_t(x) = g_t(\varphi_1(x))
\ \hbox{ for $x\in E$ and } 1 \leq t \leq 2,
\leqno (20.76)
$$
where the functions $g_t : U \to U$ are such that
$$
g_t(y) = y  \ \hbox{ for } y\in U \setminus \bigcup_j B_j
\hbox{ and for } t = 1
\leqno (20.77)
$$
and will now be defined separately on the $B_j$.
We shall use cut-off functions $\xi_j$, defined
by $\xi_j(y) = 0$ for $y\in U \sm B_j$, and
$$
\xi_j(y) = \min\big\{ 1, (\tau r_j)^{-1} \dist(y,\d B_j)\big\}
\ \hbox{ for } y\in B_j.
\leqno (20.78)
$$
Here $\tau > 0$ is another small constant that will be chosen soon.

We start with the simpler case when $y_j$ does not lie in any
$L'_i = L_i \sm {\rm int}(L_i)$, where ${\rm int}(L_i)$ is the
$n$-dimensional interior of $L_i$. In this case we can pick
any unit vector $v_j$, and set
$$
g_t(y) = y + (t-1) \xi_j(y) \eta r_j v_j
\leqno (20.79)
$$
for $y\in B_j$ and $1 \leq t \leq 2$, where $\eta > 0$ is a
minuscule constant, to be chosen later.

When the rigid assumption holds, we define the $g_t$ by the
same formula (20.79), but we make sure to choose $v_j$ in the
vector space parallel to  the smallest face $F(y_j)$ that
contains $y_j$. This precaution will only help if $y_j$
lies in some $L'_i$.

In the remaining case when the Lipschitz assumption holds
and $y_j$ lies in some $L'_i$, we need
to be more careful and use the mapping $\psi$ of Definition~2.7.
Still denote by $\wt P(y_j)$ the approximate tangent plane to 
$\wt \varphi_1(E^\ast)$ at $\wt y_j = \psi(\lambda y_j)$,
as in (20.64), and denote by $\wt \pi_j$ the orthogonal projection 
onto $\wt P(y_j)$. Also choose a unit vector $\wt v_j$ in the vector space
parallel to $\wt P(y_j)$, and then set
$$
\wt g_t(y) = \psi(\lambda y) 
+ (t-1)\xi_j(y) [\wt\pi_j(\psi(\lambda y))-\psi(\lambda y)]
+ (t-1) \xi_j(y) \eta \lambda r_j \wt v_j
\leqno (20.80)
$$
for $y\in B_j$ and $1 \leq t \leq 2$. Notice that
$$\eqalign{
|\wt g_t(y) - \psi(\lambda y_j)|
& \leq |\psi(\lambda y)-\psi(\lambda y_j)|
+ (t-1)|\wt\pi_j(\psi(\lambda y))-\psi(\lambda y)| + (t-1)\eta \lambda r_j
\cr&
\leq 2|\psi(\lambda y)-\psi(\lambda y_j)| + \eta \lambda r_j
\leq 2\lambda \Lambda r_j + \eta \lambda r_j \leq 3 \lambda \Lambda r_j
}\leqno (20.81)
$$
because $\psi(\lambda y_j)$ lies in $\wt P(y_j)$ and if 
$\eta$ is small enough. Recall from (20.69) and (20.74) that 
$$
\dist(y_j,U \sm S) \geq 4 \Lambda^2 r(y_j) \geq 4\Lambda^2 r_j
\leqno (20.82)
$$
Since $\psi$ maps $\lambda U$ to $B(0,1)$,
this implies that $\psi(\lambda y_j)$ is at distance at least 
$4 \lambda \Lambda r_j$ from $B(0,1) \sm \psi(\lambda S)$,
hence also from $\R^n \sm \psi(\lambda S)$ (recall that $y_j \i S$,
which is compactly contained in $U$). Now (20.81) yields  
$\dist(\wt g_t(y),\R^n \sm \psi(\lambda S)) \geq \lambda \Lambda r_j$.
In particular $\wt g_t(y) \in \psi(\lambda S)$, so we can define 
$$
g_t(y) = \lambda^{-1} \psi^{-1}(\wt g_t(y)) \in S.
\leqno (20.83)
$$
Observe that by (20.81)
$$
|g_t(y) -y_j| \leq \lambda^{-1} \Lambda |\wt g_t(y) - \psi(\lambda y_j)|
\leq 3 \Lambda^2 r_j.
\leqno (20.84)
$$

This completes our definition of the $g_t$ and the $\varphi_t$.
Our next task is to show that the $\varphi_{2t}$, $0 \leq t \leq 1$,
define an acceptable competitor for $E$. There is no problem with
(1.4) and (1.8); our mappings $\varphi_t(x)$, $t \geq 1$ are 
clearly continuous in $x$ and $t$, and Lipschitz in $x$
(notice in particular that all our definitions yield 
$g_t(y) = y$ on the $\d B_j$). Also, 
$$
\varphi_t(x) = x \ \hbox{ for $0 \leq t \leq 2$ when } x\in E \sm S, 
\leqno (20.85)
$$
just because $S$ contains $\wh W$ and the $B_j$
(see (20.53), (20.69), and (20.77)). In addition, we claim that
$$
\varphi_t(x) \in S \hbox{ when } x\in E \cap S. 
\leqno (20.86)
$$
When $t \leq 1$, this comes from the fact that 
$\varphi_t(W_t) \i \wh W \i S$. When $t \geq 1$,
we know that $\varphi_1(x) \in S$, and then we just need to use
(20.76) and (20.77) or (20.83). This proves (20.86), and (1.5) and (1.6),
relative to $B(x_0,r_0+\varepsilon)$, follow (by (20.53)). 
We also get the analogue of (2.4), where we use the compact set $S \i U$
in lieu of $\wh W$.

Next we check the boundary condition (1.7). We do this under the 
Lipschitz assumption; the rigid case is just simpler. 
Let $i \leq j_{max}$ and $x\in E \cap L_i \cap \overline B(x_0,r_0+\varepsilon)$
be given; recall that we want to check that $\varphi_t(x) \in L_i$ for
all $t$. We already know this for $t \leq 1$, by (1.7) for the initial
$\varphi_t$ (and (1.5) if $x\notin B(x_0,r_0)$), so we can assume that 
$t > 1$. Set $y = \varphi_1(x)$; by (20.76) $\varphi_t(x) = g_t(y)$,
and we can assume that $y\in B_j$ for some $j$, because otherwise
$\varphi_t(x) = g_t(y) = y = \varphi_1(x) \in L_i$ by (20.76). Let
us record that (1.7) will follow as soon as we show that
$$
g_t(y) \in L_i \ \hbox{ for } 1 < t \leq 2
\leqno (20.87)
$$
when $i$, $j$, $x \in E \cap L_i$, and $y = \varphi_1(x) \in B_j$
are as above.

By (1.7) for $\varphi_1$, $y = \varphi_1(x)$ lies in $L_i$.
Let $F$ be a face of $L_i$ that contains $y$, and let us check that 
$$
F(y_j) \i F. 
\leqno (20.88)
$$
Set $\wt F = \psi(\lambda F)$ and 
$\wt F_j = \psi(\lambda F(y_j))$, and observe that
$\wt F_j$ is the smallest rigid face that contains 
$\wt y_j = \psi(\lambda y_j)$.
Suppose that (20.88) fails; then $\wt F_j$ is not contained in $\wt F$. 
If in addition $\wt F_j$ is not reduced to the point $\wt y_j$, 
(3.8) yields
$$\eqalign{
\dist(\wt y_j,\wt F) &\geq \dist(\wt y_j, \d(\wt F_j))
\geq \lambda \Lambda^{-1}\dist(y_j, \d F(y_j))
\geq 4 \lambda \Lambda^3 r(y_j)
\geq 4 \lambda \Lambda r_j
}\leqno (20.89)
$$
by (20.70) and (20.74). If  instead $\wt F_j$ is reduced to the point
$\wt y_j$, then $\wt F$ is a rigid face that does not contain the
vertex $\wt y_j$, hence $\dist(\wt y_j,\wt F)
\geq 2^{-m} \geq 4 \lambda \Lambda r_j$, by (20.74) and if 
$\varepsilon$ is small enough; so the conclusion of (20.89) holds in 
both cases.
But $\psi(\lambda y) \in \wt F$ because $y\in F$, so
$$
\dist(\wt y_j,\wt F) \leq |\psi(\lambda y_j)-\psi(\lambda y)| 
\leq \lambda \Lambda |y-y_j| \leq \lambda \Lambda r_j,
\leqno (20.90)
$$
a contradiction which proves (20.88). 

Recall that we want to check (20.87). 
We start with the most interesting case when 
$y_j \in L'_i$. Then $g_t(y)$ was defined by (20.80) and (20.83),
and (by (20.83)) it is enough to check that $\wt g_t(y) \in \wt F$.

Recall that $\wt \pi_j$ is the orthogonal projection onto the 
approximate tangent plane $\wt P(y_j)$ of (20.64), which itself is 
contained in the affine plane $\wt W_j$ spanned by $\wt F_j$,
by (20.66). Denote by $\wt W$ the affine span of $F$;
by (20.88), $\wt F_j \i \wt F$ and hence $\wt W_j \i \wt W$.
Thus the points $\wt y = \psi(\lambda y)$, $\wt \pi_j(\wt y)$, and even
$\wt w = \wt y + (t-1) \xi_j(y) [\wt \pi_j(\wt y)-\wt y]$
all lie in $\wt W$. Observe that by (20.80),
$\wt g_t(y) = \wt w + (t-1) \xi_j(y) \eta \lambda r_j \wt v_j$,
and since we chose $\wt v_j$ in the vector space parallel to 
$\wt P(y_j) \i \wt W$, we see that $\wt g_t(y) \i \wt W$.

We want to show that $\wt g_t(y)$ even lies in $\wt F$.
We start from the fact that $\wt y_j \in \wt F_j \i \wt F$ 
(by definition of $\wt F_j = \psi(\lambda F(y_j))$ and by (20.88)),
with
$$\eqalign{
\dist(\wt y_j, \d \wt F)
&\geq \dist(\wt y_j, \d \wt F_j)
\geq \lambda \Lambda^{-1}\dist(y_j, \d F(y_j))
\geq 4 \lambda \Lambda^3 r(y_j)
\geq 4 \lambda \Lambda r(y_j)
}\leqno (20.91)
$$
by (20.70) and (20.74) (that is, as in (20.89)).
But $|\wt g_t(y)-\wt y_j| \leq 3 \lambda \Lambda r_j$
by (20.81), so the line segment $[\wt y_j,\wt g_t(y)] \i \wt F$
does not meet $\d \wt F$, and $\wt g_t(y) \in \wt F$; (20.87)
follows, because $F \i L_i$ by definition, and this takes care of our 
first case.

We are left with the case when $y_j \notin L'_i$. Since 
$y_j \in F(y_j) \i F \i L_i$ by (20.88), this implies that $y_j$ 
lies in the interior of $L_i$. We want to show that in fact
$$
\overline B(y_j, 4\Lambda^2 r_j) \i L_i \, ,
\leqno (20.92)
$$
and for this we shall proceed as for (19.10).
Denote by $\delta(L_i)$ the boundary of $L_i$,
set $D= \dist(y_j,\delta(L_i)) > 0$, and pick $\xi \in \delta(L_i)$ 
such that $|\xi-y_j| = D$. Denote by $G$ the smallest face of our 
grid that contains $\xi \,$; since $\delta(L_i)$ is itself an union of 
faces, $G$ is contained in $\delta(L_i)$. Since $D > 0$,
$G$ does not contain $y_j$, and even less $F(y_j)$.

First assume that $F(y_j)$ is not reduced to $\{ y_j \}$; then
(3.8) (applied to the rigid faces $\psi(\lambda G)$ and $\wt F_j$) yields
$$\eqalign{
D &= |y_j- \xi| \geq \dist(y_j,G) 
\geq \lambda^{-1} \Lambda^{-1} \dist(\psi(\lambda y_j), \psi(\lambda G)) 
\cr&
\geq \lambda^{-1} \Lambda^{-1} \dist(\psi(\lambda y_j), \wt F_j) 
\geq \Lambda^{-2} \dist(y_j,\d F(y_j))
\geq 4 \Lambda^2 r(y_j) \geq 4 \Lambda^2 r_j
}\leqno (20.93)
$$
by (20.70) and (20.74). If instead $F(y_j)=\{ y_j \}$,
and since $D$ is the distance from the vertex $y_j$ to a face that
does not contain it, we get that 
$D \geq \lambda^{-1} \Lambda^{-1} 2^{-m} \geq 4 \Lambda^2 r_j$,
by (20.74) and if $\varepsilon$ is small enough. Thus 
$D \geq 4 \Lambda^2 r_j$ in both cases, and (20.92) follows.
In this case the fact that $\varphi_t(x)$ lies
in $L_i$ is trivial because 
$\varphi_t(x) = g_t(y) \in \overline B(y_j, 4\Lambda^2 r_j)$,
by (20.84).

This completes our proof of (1.7), and the series of
verifications for the extended family 
$\{ \varphi_t \}$, and now we can use the $A$-minimality
of $E$. This yields
$$
\H^d(W_2) \leq \H^d(\varphi_2(W_2)) + h(r_0+\varepsilon) 
(r_0+\varepsilon)^d,
\leqno (20.94)
$$
by (20.5) and where 
$W_2 = \big\{ x\in E \, ; \, \varphi_2(x) \neq x \big\}$.

Recall that $W_2 \cup \varphi_2(W_2) \i S$, by (20.85) and (20.86).
We start with an estimate of $\H^d(\varphi_2(W_2))$.
Set $A = S \cap \varphi_1(E) \setminus \bigcup_j 
[B_j \cap \varphi_1(E)]$. Then
$$\eqalign{
\H^d(\varphi_2(W_2)) &\leq \H^d(S \cap \varphi_2(E))
= \H^d(g_2(S \cap \varphi_1(E)))
\cr&
\leq \H^d(A) + \sum_j \H^d(g_2(B_j \cap \varphi_1(E)))
}\leqno (20.95)
$$
by (20.76), because $S \cap \varphi_1(E) \i 
\i A \cup \big(\bigcup_j [B_j \cap \varphi_1(E)]\big)$, and 
because $g_2(y)=y$ on $A$ (by (20.77)).
This will be compared to the fact that
$$
\H^d(S \cap \varphi_1(E)) = \H^d(A) + \sum_j \H^d(B_j \cap \varphi_1(E))
\leqno (20.96)
$$
by definition of $A$ and because the $B_j$ are disjoint.
Next we estimate the sum in (20.96). There are two types of
indices $j$; we start with the simple case when $g_t$ was 
defined by (20.79). That is,
$g_2(y) = y + \xi_j(y) \eta r_j v_j$ for $y\in B_j$.
Write $B_j = B_{j,int} \cup B_{j,ext}$, where
$B_{j,int} = \big\{y\in B_j ; \dist(y,\d B_j) \geq \tau r_j\big\}$
and $B_{j,ext} = B_j \sm B_{j,int}$. On $B_{j,int}$, (20.78)
yields $\xi_j(y)=1$ and $g_2(y) = y + \eta r_j v_j$, and hence
$$
\H^d(g_2(B_{j,int}\cap \varphi_1(E)))
= \H^d(B_{j,int}\cap \varphi_1(E))
\leq \H^d(B_j \cap \varphi_1(E)).
\leqno (20.97)
$$
On $B_{j,ext}$, although $\xi_j$ is only $(\tau r_j)^{-1}$-Lipschitz,
we can choose $\eta $ so small that $g_2$ is $2$-Lipschitz on
$B_{j,ext}$, and we get that
$$
\H^d(g_2(B_{j,ext}\cap \varphi_1(E))) 
\leq 2^d \H^d(B_{j,ext}\cap \varphi_1(E))
\leqno (20.98)
$$
But two applications of (20.71) yield
$$\eqalign{
\H^d(B_{j,int}\cap \varphi_1(E))
&= \H^d(B_{j,int}\cap \varphi_1(E^\ast))
\geq [(1 -\tau) r_j]^d (\omega_d -\varepsilon) 
\cr&
\geq (1 -\tau)^d \, {\omega_d -\varepsilon \over \omega_d +\varepsilon}\,
\H^d(B_j\cap \varphi_1(E))
}\leqno (20.99)
$$
because $\H^d(E \sm E^\ast) = 0$. Thus, if $\varepsilon$ is small 
enough, depending on $\tau$, we get that
$$
\H^d(B_{j,ext}\cap \varphi_1(E)) 
= \H^d(B_j\cap \varphi_1(E)) - \H^d(B_{j,int}\cap \varphi_1(E))
\leq C \tau \H^d(B_j\cap \varphi_1(E))
\leqno (20.100)
$$
and, by (20.97) and (20.98),
$$
\H^d(g_2(B_j\cap \varphi_1(E)))
\leq (1+C\tau) \H^d(B_j \cap \varphi_1(E)).
\leqno (20.101)
$$
Now we consider the more complicated case when we used 
(20.80)-(20.83) to define $g_2$. By (20.80),
$$
\wt g_2(y) = \psi(\lambda y) 
+ \xi_j(y) [\wt\pi_j(\psi(\lambda y))-\psi(\lambda y)]
+ \xi_j(y) \eta \lambda r_j \wt v_j
\leqno (20.102)
$$
for $y\in B_j$. We start with the good set 
$$
G(j) = \big\{ y \in B_{j,int} \cap \varphi_1(E) \, ; \, 
\dist(\psi(\lambda y),\wt P(y_j)) \leq \varepsilon \lambda r_j \big\},
\leqno (20.103)
$$
where $\wt P(y_j)$ is the tangent plane that shows up in (20.64) and 
(20.72), for instance. If $y\in G(j)$, (20.78) yields $\xi_j(y)=1$,
then (20.102) says that $\wt g_2(y) = \wt\pi_j(\psi(\lambda y)) 
+\eta \lambda r_j \wt v_j$, hence
$$
|\wt g_2(y) - \psi(\lambda y)| \leq
|\wt\pi_j(\psi(\lambda y)) - \psi(\lambda y)| + \eta \lambda r_j
\leq (\eta + \varepsilon) \lambda r_j
\leqno (20.104)
$$
(by (20.103) and because $\wt \pi_j$ is the orthogonal projection 
on $\wt P(y_j)$).
Then $|g_2(y)-y| \leq (\eta + \varepsilon) \Lambda r_j$
by (20.83), and hence
$g_2(y) \in B(y_j, (1+ \eta \Lambda + \varepsilon \Lambda) r_j)$.
Also, $\wt\pi_j(\psi(\lambda y)) \in \wt P(y_j)$, 
and since we chose $\wt v_j$ in the vector space
parallel to $\wt P(y_j)$ (see above (20.80)), we see that
$\wt g_2(y) \in \wt P(y_j)$. Thus
$$
g_2(y) \in \lambda^{-1} \psi^{-1}(\wt P(y_j)) \cap 
B(y_j, (1+ \eta \Lambda + \varepsilon \Lambda)\, r_j).
\leqno (20.105)
$$
by (20.83). By definition of $\wt P(y_j)$
(see (20.64)), $\psi(\lambda y_j) \in \wt P(y_j)$.
By (20.66), $\wt P(y_j) \i \wt W(y_j)$, the affine span of
$\wt F_j = \psi(\lambda F(y_j))$. We now deduce from (20.105) and 
the definition (19.33) of $A_r(y_j)$ that
$$\eqalign{
\H^d(g_2(G(j)))
&\leq \H^d(\lambda^{-1} \psi^{-1}(\wt P(y_j)) \cap 
B(y_j, (1+ \eta \Lambda + \varepsilon \Lambda) r_j))
\cr&
\leq (1+ \eta \Lambda + \varepsilon \Lambda)^d r_j^d
A_{(1+ \eta \Lambda + \varepsilon \Lambda) r_j}(y_j)
\cr&
\leq (\omega_d + \varepsilon)
(1+ \eta \Lambda + \varepsilon \Lambda)^d r_j^d
\cr&
\leq (1 + C\varepsilon \Lambda + C\eta \Lambda)
\,\H^d(B_j \cap \varphi_1(E))
}\leqno (20.106)
$$
by (20.73) and (20.71). Next we consider the less good set
$$
G'(j) = \big\{ y \in B_{j,ext} \cap \varphi_1(E) \, ; \, 
\dist(\psi(\lambda y),\wt P(y_j)) \leq \varepsilon \lambda r_j \big\}.
\leqno (20.107)
$$
We claim that $\wt g_2$ is $C \lambda\Lambda$-Lipschitz 
on $G'(j)$. The first term in the definition (20.102) is
$\psi(\lambda y)$, which is $\lambda \Lambda$-Lipschitz;
the third one, $\xi_j(y) \eta \lambda r_j \wt v_j$,
is $C \tau^{-1} \eta \lambda$-Lipschitz, which is much better
if $\eta$ is small enough.
Notice that $|\wt\pi_j(\psi(\lambda y))-\psi(\lambda y)|
\leq \varepsilon \lambda r_j$ on $G'(j)$, hence
the second term $\xi_j(y) [\wt\pi_j(\psi(\lambda y))-\psi(\lambda y)]$
is $C \tau^{-1} \varepsilon \lambda + C \Lambda \lambda$-Lipschitz,
our claim follows, and $g_2$ is $C \Lambda^2$-Lipschitz on $G'(j)$. 
Then
$$\eqalign{
\H^d(g_2(G'(j)))
&\leq C \Lambda^{2d} \H^d(G'(j))
\leq C \Lambda^{2d} \H^d(B_{j,ext} \cap \varphi_1(E))
\cr&
\leq C \tau \H^d(B_j \cap \varphi_1(E))
}\leqno (20.108)
$$
by (20.100), and where we no longer write the dependence on $\Lambda$
in the last line. We are left with
$$
G''(j) = \big\{ y \in B_{j} \cap \varphi_1(E) \, ; \, 
\dist(\psi(\lambda y),\wt P(y_j)) > \varepsilon \lambda r_j \big\}.
\leqno (20.109)
$$
On this set (20.102) only yields that
$\wt g_2$ is $C \lambda \Lambda \tau^{-1}$-Lipschitz,
hence by (20.83) $g_2$ is $C \Lambda^2 \tau^{-1}$-Lipschitz.
Fortunately $G''(j)$ is small. Indeed if $y\in G''(j)$,
then $\psi(\lambda y)$ lies in the bad set of (20.72),
whose measure is at most $\varepsilon \lambda^d r^d_j$
hence $\H^d(G''(j)) \leq \varepsilon \Lambda^d r^d$
and (dropping soon the dependence on $\Lambda$ and by (20.71) again),
$$\eqalign{
\H^d(g_2(G''(j)))
&\leq C \Lambda^{2d} \tau^{-d} \H^d(G''(j))
\leq C \tau^{-d} \varepsilon r_j^d
\leq C \tau^{-d} \varepsilon\H^d(B_j \cap \varphi_1(E)).
}\leqno (20.110)
$$
We add (20.106), (20.108), and (20.110) and get that
$$
\H^d(g_2(B_j\cap \varphi_1(E))) 
\leq (1 + C\eta + C\tau + C\varepsilon \tau^{-d})
\H^d(B_j \cap \varphi_1(E)).
\leqno (20.111)
$$
We had a slightly better estimate (20.101) in the first case,
so (20.111) holds in all cases, and when we compare (20.95)
to (20.96), we now get that
$$\eqalign{
\H^d(\varphi_2(W_2)) & \leq \H^d(S \cap \varphi_1(E))
+ C(\eta + \tau + \varepsilon \tau^{-d})
\sum_i \H^d(B_j \cap \varphi_1(E))
\cr&\leq \H^d(S \cap \varphi_1(E))
+ C(\eta + \tau + \varepsilon \tau^{-d})
\H^d(S \cap \varphi_1(E))
}\leqno (20.112)
$$
(recall that the $B_j$ are disjoint, and (by (20.69) and (20.74))
contained in $S$). Now we want to check that
$$
\H^d(E \cap S \sm W_2) \leq C (M^{-1} + \tau), 
\leqno (20.113)
$$
where $C$ is allowed to depend on $\H^d(E\cap \wh W)$,
and $M$ is as in (20.56) and (20.60).

Let $E \cap S \sm W_2$ be given. Let us remove a few small sets.
A first possibility is that
$x\in W = \big\{ x\in U \, ; \ \varphi_1(x) \neq x \big\}$.
Set $y=\varphi_1(x)$; then $y \neq x$ because $x\in W$. 
Since $\varphi_2(x) = x$ (because $x\notin W_2$) and 
$\varphi_2(x) = g_t(y)$ by (20.76) we get that
$g_2(y) \neq y$, and even $|g_2(y)-y| = |y-x| = |\varphi_1(x)-x|$.
Then $y$ lies in some $B_j$. If $g_2(y)$ was computed by (20.79),
this implies that $|\varphi_1(x)-x| \leq \eta r_j$. Otherwise,
(20.84) says that $|g_2(y)-y| \leq |g_2(y)-y_j|+r_j
\leq 4\Lambda^2 r_j$. In both cases,
$|\varphi_1(x)-x| \leq 4\Lambda^2 r_j \leq 4\Lambda^2 \varepsilon$
by our precaution (10.74). If $\varepsilon$ is small enough, depending
on $\tau$, we deduce from this that
$$
\H^d(E \cap W \sm W_2)
\leq \H^d\big(\{x \in E \cap W \, ; \, |\varphi_1(x)-x|
\leq 4\Lambda^2 \varepsilon \}\big) \leq \tau
\leqno (20.114)
$$
because the monotone intersection, when $\varepsilon$
tends to $0$, of the sets in (20.114) is empty, and all these
sets are contained in $E \cap W$ for which $\H^d(E \cap W) < +\infty$.
So we may restrict to $x\in E \cap S \sm [W_2 \cap W]$. 
Since $\H^d(E\cap S \sm \wh W) \leq \varepsilon$ by (20.53),
this set contributes little to (20.113), and we may assume
that $x\in \wh W$. Since $x\in E \sm W$, we get that
$\varphi_1(x) = x$, and so $x\in \varphi_1(E\cap \wh W)$.
Thus (20.54) says that $x$ almost always lies in $G_0$. 
Next we take care of $Y_0$, which by (20.60) is such that
$$
\H^d(Y_0) \leq M^{-1} (1+ \H^d(E\cap \wh W));
\leqno (20.115)
$$
this is less than the right-hand side of (20.113), 
so we may now assume that 
$x \in G_1 = G_0 \sm Y_0$ (see below (20.61)),
or even that $x$ lies in some $B_j$, because
(20.68) says that $\H^d(G_1\sm G_2) = 0$
and then (20.75) says that the $B_j$ almost cover $G_2$
(recall from the line below (20.53) that $\mu$ is the 
restriction of $\H^d$ to $\varphi_{1}(E \cap S)$).
By (20.100), 
$$
\sum_j \H^d(B_{j,ext}\cap \varphi_1(E)) 
\leq C \tau \sum_j\H^d(B_j\cap \varphi_1(E))
\leq C \tau H^d(S \cap \varphi_1(E))
\leqno (20.116)
$$
(recall again that the $B_j$ are disjoint and contained in $S$
(by (20.69) and (20.74)).
This bound is also compatible with (20.113),
so we are left with the case when $x\in B_{j,int}$.
In this case, $\xi_j(x) = 1$, and we claim that
$\varphi_2(x) = g_2(x) \neq x$. The first part follows from
the (20.76) because $\varphi_1(x) = x$. When $g_2(x)$
is given by (20.79), the second part is obvious. When
we use (20.80), projecting on $\wt P_j$ yields
$$
\wt\pi_j(\wt g_2(x)) = \wt\pi_j(\psi(\lambda x))+ \eta \lambda r_j 
\wt v_j \neq \wt\pi_j(\psi(\lambda x))
\leqno (20.117)
$$
because $\wt v_j$ was chosen to be a unit vector in the direction
of $\wt P(y_j)$. Then $\wt g_2(x) \neq \psi(\lambda x)$ and,
by (20.83), $g_2(x) \neq x$, as needed. But this is impossible,
because we assumed that $x\in E \cap S \sm W_2$. Then (20.113)
holds, and we may now put all our estimates together:
$$\leqalignno{
\H^d(E \cap \wh W) &\leq \H^d(E \cap S)
\leq \H^d(W_2) + C(M^{-1} + \tau) 
\cr&\leq \H^d(\varphi_2(W_2)) + h(r_0+\varepsilon) (r_0+\varepsilon)^d
+ C(M^{-1} + \tau) 
\cr&
\leq \H^d(S \cap \varphi_1(E))
+ C(\eta + \tau + \varepsilon \tau^{-d}) 
+ h(r_0+\varepsilon) (r_0+\varepsilon)^d
+ C(M^{-1} + \tau)
& (20.118)
}
$$
by (20.53), (20.113), (20.94), and (20.112) (where we now see
$\H^d(S \cap \varphi_1(E))$ as a constant).
Let us check that
$$
\H^d(S \cap \varphi_1(E)) \leq \H^d(\wh W \cap \varphi_1(E)) + \varepsilon.
\leqno (20.119)
$$
Suppose $y \in S \cap \varphi_1(E) \sm \wh W$, and let
$x \in E$ such that $\varphi_1(x) = y$. If $y \neq x$,
(2.2) says that $y \in \wh W$, which is impossible. So
$y=x$, and now $y \in E \cap S \sm \wh W$; (20.119)
then follows from (20.53).

When we add (20.118) and (20.119), we get that
$\H^d(E \cap \wh W) \leq \H^d(\wh W \cap \varphi_1(E)) + e$,
with 
$$
e = C(\eta + \tau + \varepsilon \tau^{-d}) 
+ h(r_0+\varepsilon) (r_0+\varepsilon)^d + C(M^{-1} + \tau).
\leqno (20.120)
$$
Of course, $C$ depends on $E$ and $\varphi_1$ in various ways,
but we can choose $\tau$, then $\varepsilon$ and $M$
(recall that we never used $M$ in the estimates, so we
can choose it as large as we want), then $\eta$
so small that $e$ is as close to $h(r_0) r_0^d$ as we want.
This proves (20.7), the $A'$-almost minimality of $E$ follows,
and so does Proposition 20.9 (in the general case).
\qed

\msi
{\bf 21. Limits of almost minimal sets and of minimizing sequences.} 
\ms  
In this section we just rewrite Theorem 10.8 in the context of almost 
minimal sets. For our first statement, we consider a gauge function 
$h : (0,+\infty) \to [0,+\infty]$ which is right-continuous, i.e., such that
$$
h(r) = \lim_{\rho \to r \, ; \, \rho > r} h(\rho)
\ \hbox{ for } r > 0, 
\leqno (21.1)
$$
and for which
$$
\lim_{r \to 0} h(r) = 0.
\leqno (21.2)
$$

\ms\proclaim Theorem 21.3.
Let an open set $U$ and boundary pieces $L_j$, $0 \leq j \leq j_{max}$,
be given, and suppose that the Lipschitz assumption holds (see 
Definition 2.7). Also suppose that the technical assumption (10.7),
or the weaker (19.36) holds (but this is not needed under the rigid 
assumption (2.6)). Let $\{ E_k \}$ be a sequence of coral (see Definition 3.1)
and relatively closed sets in $U$, that converges 
locally in $U$ to the closed set $E$ (as in (10.4)-(10.6)).
\break
{\bf 1.} If each $E_k$ is an $A_+$-almost minimal set in $U$,
with the sliding conditions given by the sets $L_j$, and the gauge 
function $h$ (see Definition 20.2), then 
$E$ is coral, and it is an $A_+$-almost minimal set in $U$,
with the sliding conditions given by the same sets $L_j$ 
and the same gauge function $h$.
\hfill\break
{\bf 2.} If each $E_k$ is an $A$-almost minimal set in $U$,
with the sliding conditions given by the sets $L_j$, and the gauge 
function $h$, then $E$ is coral, and it is an $A$-almost minimal set in $U$,
with the sliding conditions given by the same sets $L_j$ and the gauge 
function $h$.
\hfill\break
{\bf 3.} If each $E_k$ is an $A'$-almost minimal set in $U$,
with the sliding conditions given by the sets $L_j$, and the gauge 
function $h$, then $E$ is coral, and it is an $A'$-almost minimal set in $U$, 
with the sliding conditions given by the same sets $L_j$ and 
the gauge function $h$.

\msi {\bf Proof.}
We start with limits of $A_+$-almost minimal sets.
Since we want to apply Theorem~10.8, we compare Definition 20.2
with the definition 2.3 of quasiminimality.
If $E_k$ is $A_+$-almost minimal as above, then for each 
$\delta > 0$, $E_k \in GSAQ(U,M(\delta),\delta,0)$, 
with $M(\delta) = 1+h(\delta)$.

By Theorem 10.8, or its variant in Remark 19.52
where we assume (19.36) instead of (10.7), 
$E$ also satisfies this property.
Notice in particular that since here the last constant $h$ in the 
definition of $GSAQ$ is zero, the additional constraint above (10.2)
that requires $h$ to be small is automatically satisfied.
That is, for each $\delta > 0$, $E \in GSAQ(U,M(\delta),\delta,0)$.
But then $E$ is $A_+$-almost minimal with the gauge function $h'$
defined by 
$$
h'(r) = \liminf_{\delta \to r^+} \, h(\delta). 
\leqno (21.4)
$$
By (21.1), $h'=h$ and Part 1 of our result follows.

Next consider a sequence of $A$-almost minimal sets.
If the $E_k$ are as in Part 2, then for each $\delta > 0$, 
$E_k \in GSAQ(U,1,\delta,h(\delta))$ for all $k$.

We have a minor additional difficulty here, because in order to apply
Theorem 10.8, we have to assume that 
$E_k \in GSAQ(U,M,\delta,h)$ with $h$ sufficiently small,
depending on $n$, $M$ and $\Lambda$. Here this is true
for $\delta$ small, by (21.2), but maybe not for $\delta$ large.

Fortunately, as was noted below the statement of Theorem 10.8,
this assumption that $h$ be small enough is only needed to get
the right regularity and lower semicontinuity properties,
but as soon as it is satisfied for some acceptable combination of
$M,\delta,h$ (here with $M=1$ and $\delta$ so small that
$h=h(\delta)$ works), we get the limiting theorem for the other 
combinations. Thus $E \in GSAQ(U,1,\delta,h(\delta))$ for $\delta > 0$.

Then we return to Definition 20.2 and get that 
$E$ is $A$-almost minimal with the gauge function $h'$ of (21.4).
Since $h'=h$ by (21.1), Part 2 follows.

For Part 3, we just need to observe that because of Proposition 21.9, 
we do not need to distinguish between $A$-almost minimal and $A'$-almost minimal
(notice that the additional sufficient condition for the equivalence,
(10.7) or (19.36), is satisfied). Then Part 3 follows from Part 2.
\qed

\msi{\bf Remark 21.5.}
Probably we could modify our proof of Theorem 10.8 to make it work 
also for $A'$-almost minimal sets (and even with the variant of
quasiminimal sets defined with the same accounting as in (20.6)). 
We should not expect a huge simplification, and in particular we cannot
content ourselves with applying the almost minimality of
$E_k$ with any extension of our initial mapping $\varphi_1$, 
because it still could be that $\varphi_1(E_k)$ is a very bad competitor 
because it contains may parallel sheets, that could easily be merged
to produce a better competitor, while these sheets are already merged 
for $E$. 

Also, we would have to take into account the possibility
that the set $\varphi_1(E \cap W_1)$ meets $E \sm W_1$ 
(where as usual $W_1 = \big\{ x\in \R^n \, ; \, \varphi_1(x) \neq x \big\}$),
while this does not happen with $E_k$. Then $\varphi_1$ defines
a better competitor for $E$ than for $E_k$, which is also bad
for our proof. We did not pay attention to
this case in the proof of Theorem 10.8, because it did not matter
with the accounting for quasiminimal sets, but of course we could try 
to fix it, for instance by allowing a larger piece of $\wh W$
in the definition of $X_0$ in (11.20). But this becomes similar
to our proof of Proposition 21.9, so the author does not expect
to win much by trying a direct proof.

\msi{\bf Remark 21.6.}
If we did not assume (21.1), we would still have that the limit
$E$ is almost minimal, but this would be with the gauge function $h'$
defined in (21.4); this is easy to see from the proof, and 
(for Part 3) the similar comment below Proposition 20.9.

\msi{\bf Remark 21.7.}
Similarly, we do not really need to assume (21.2), but instead we can
assume that (10.2) holds, i.e., that there are constants 
$M$, $\delta$, and $\hbar$, with $\hbar$ small enough 
(depending on $n$, $\Lambda$, and $M$) such that
$E_k \in GSAQ(U,M,\delta, \hbar)$ for all $k$. Then
we can use the remark below Theorem 10.8 (as we did for Part 2)
and proceed as above, because (10.2) is enough for the regularity
results of Section 10. We shall apply this now, in the context of
local minimizing sequences.

\ms
Here is the notation for the next corollary.
We are given, as in Theorem 21.3, an open set $U$ 
and boundary pieces $L_j$, $0 \leq j \leq j_{max}$,
and we suppose that 
$$
\hbox{the Lipschitz assumption holds, as well as (10.7) or (19.36)}
\leqno (21.8)
$$
(again see Definition 2.7 and observe that (10.7) is automatic under
the rigid assumption). We are also given a sequence $\{ E_k \}$ of 
coral relatively closed sets in $U$, and we assume that
$$
\hbox{the $E_k$ converge locally in $U$ to the 
relatively closed set $E \i U$.}
\leqno (21.9)
$$
In addition, we assume that there are constants 
$M$, $\delta$, and $h$, with $h$ small enough 
(depending on $n$, $\Lambda$, and $M$) such that (10.2)
holds, i.e., 
$$
E_k \in GSAQ(U,M,\delta, h) \hbox{ for all } k.
\leqno (21.10)
$$
Finally, we assume that $\{ E_k \}$ is a locally minimizing
sequence, in the following sense. Given $\delta > 0$, we say that
one-parameter family $\{\varphi_t \}$ of functions is 
$\delta$-admissible for $E_k$ if it satisfies the conditions
(1.4)-(1.8), relative to $E_k$ and some ball $B$ of
radius $r < \delta$, and in addition the compactness condition 
(2.4) holds (relative to $E_k$). Recall that 
(2.4) says that $\wh W(E_k) \i\i U$, where we set
$$
\widehat W(E_k) = \bigcup_{0< t \leq 1} W_t(E_k) \cup 
\varphi_t(W_t(E_k)),
\leqno (21.11)
$$
with
$$
W_t(E_k) = \big\{ y \in E_k \, ; \varphi_t(y) \neq y \big\} 
\ \hbox{ for } 0 < t \leq 1.
\leqno (21.12)
$$

We shall assume that there exists $\delta > 0$ such that,
for each $\varepsilon > 0$ we can find $k_0 \geq 0$
such that
$$
\H^d(W_1(E_k)) \leq \H^d(\varphi_1(W_1(E_k))) + \varepsilon
\leqno (21.13)
$$
for every $k \geq k_0$ and every one-parameter family $\{\varphi_t \}$
which is $\delta$-admissible for $E_k$.

Or we shall assume that, with the same quantifiers,
$$
\H^d(E_k\sm\varphi_1(E_k)) \leq \H^d(\varphi_1(E_k) \sm E_k) + \varepsilon
\leqno (21.14)
$$
for every $k \geq k_0$ and every one-parameter family $\{\varphi_t \}$
which is $\delta$-admissible for $E_k$. 
This second assumption, which is more in the mode of $A'$-almost 
minimal set, is more natural in some contexts.

\ms\proclaim Corollary 21.15.
Let $U$, the boundary pieces $L_j$, and the sequence $\{ E_k \}$
of coral quasiminimal sets satisfy the conditions (21.8)-(21.14).
Then $E$ is a coral local minimizer in $U$, in the sense that
$$
\H^d(E\sm\varphi_1(E)) \leq \H^d(\varphi_1(E) \sm E)
\leqno (21.16)
$$
for every one-parameter family $\{\varphi_t \}$
which is $\delta$-admissible for $E$.

To prove the corollary, observe that for $k \geq k_0$,
(21.13) or (21.14) says that $E_k$ is $A$-almost minimal 
or $A'$-almost minimal, with
the strange gauge function $h_\varepsilon$ defined by 
$h_\varepsilon(r) = r^{-d} \varepsilon$
for $0 < r < \delta$ and $h_\varepsilon(r) = +\infty$ for $r \geq \delta$.

This function does not satisfy (21.2), but Remark 21.7 and our 
assumption (21.10) allow us to dispense with this condition. Then
by Theorem 21.3, $E$ is coral, and almost minimal with the same gauge 
function $h_\varepsilon$. Since this is true for all $\varepsilon > 0$,
we also get that $E$ is almost minimal with the gauge function $h_0$.
If we were dealing with (21.14) and $A'$-almost minimal sets, we
directly get (21.16) from this. If we were dealing with 
(21.13) and $A$-almost minimal sets, we get 
$$
\H^d(W_1(E)) \leq \H^d(\varphi_1(W_1(E)))
\leqno (21.17)
$$
instead of (21.16), but by the easy part of 
Proposition~20.9, (21.17) implies (21.16); Corollary 21.15 follows.
\qed

\msi{\bf Remark 21.18.}
In the conclusion of Corollary 21.15, we may also replace 
(21.16) with (21.17), since the two conditions are equivalent 
(by Proposition~20.9, applied with $h_0$).

\msi
{\bf 22. Upper semicontinuity of $H^d$ along sequences of almost 
minimal sets.} 
\ms

The main result of this section is the following upper semicontinuity
result.

\ms\proclaim Theorem 22.1.
Let $U$, the $L_i$, the sequence $\{ E_k \}$, and 
the set $E$ satisfy the assumptions of Theorem 21.3 (any part) 
or Corollary 21.15. Then for every compact set $H \i U$,
$$
\H^d(E \cap H) \geq \limsup_{k \to +\infty} \H^d(E_k \cap H).
\leqno (22.2)
$$

\ms
Notice that if the $E_k$ are only supposed to be quasiminimal,
the conclusion may fail, even when there is no boundary condition.
For instance, $E_k$ may coincide locally with the graph of the
function $x \to 2^{-k}\sin(2^{k} x)$, which converges to a line;
then (22.2) fails. So, for the sequences of Theorem 21.3, the
condition (21.2) is really needed this time.

The proof will only use the rectifiability of $E^\ast$, a covering 
argument, and an application of the quasiminimality (or almost 
minimality) of the $E_k$ in balls where $E$ is flat.
It is essentially a special case of the following lemma, which is a 
generalization of Lemma 3.12 on page 85 of [D5],  
and which we shall prove first.

\ms\proclaim Lemma 22.3.
Let $U$, the $L_i$, the sequence $\{ E_k \}$, and 
the set $E$ satisfy the assumptions of Theorem 10.8. 
Then for every compact set $H \i U$,
$$
(1+Ch) M \H^d(E \cap H) \geq \limsup_{k \to +\infty} \H^d(E_k \cap H),
\leqno (22.4)
$$
with a constant $C$ that depends only on $n$, $M$, and $\Lambda$.

\ms
Our proof of Lemma 22.3 will be similar to the proof of
Proposition 21.9 in the Lipschitz case. Let $\{ E_k \}$,
$E$, and $H$ be as in the statement. We first try to cover a
big piece of $E\cap H$ by small balls.

Our assumptions allow us apply the results of Section 10. 
In particular, the $E_k$ are uniformly locally Ahlfors-regular 
(by (10.10)), and $E$ is locally Ahlfors-regular (by (10.11)) 
and rectifiable (by Proposition 10.15). 

Let $\varepsilon > 0$ be given, and use the fact that $\H^d(E)$ 
is locally finite in $U$ (for instance, because $E^\ast$ is locally 
Ahlfors regular) to choose an open set $V$ such that
$$
H \i V \i\i U \ \hbox{ and } \ 
\H^d(E \cap V \sm H) \leq \varepsilon. 
\leqno (22.5)
$$
Next, the fact that $E$ is rectifiable implies that 
for $\H^d$-almost every $x\in E \cap H$,
$$
\lim_{r \to 0} r^{-d} \H^d(E \cap B(x,r)) = \omega_d,
\leqno (22.6)
$$
(see Theorem 17.6 on page 240 in [Ma]), and 
$$
E \hbox{ has a tangent plane $P(x)$ at $x$.}
\leqno (22.7)
$$
Recall that the fact that an approximate tangent plane to $E$ is a true 
tangent plane comes from the local Ahlfors-regularity of $E$; see for 
instance Exercise 41.21 on page 277 of [D4]. 

We shall assume that the Lipschitz assumption holds; the rigid case
is easier, and we could also obtain it the complicated way, by 
pretending that $U = B(0,1)$ and $\psi$ is the identity.
For $x\in E \cap H$, denote by $F(x)$ the smallest (twisted) 
face of our grid that contains $x$. We also set $\wt x = \psi(\lambda x)$ 
and $\wt F(x) = \psi(\lambda F(x))$ (a true dyadic face).
For almost every $x \in E\cap H$ such that (22.7) holds,
we also have that
$$
\wt E = \psi(\lambda E) \hbox{ has a tangent plane $\wt P(x)$ at $\wt x$,}
\leqno (22.8)
$$
because $\wt E$ is also rectifiable and locally Ahlfors-regular
(recall that $\psi$ is bilipschitz). We also want to show that
$\H^d$-almost everywhere on $E\cap H$, 
$$
\wt P(x) \hbox{ is contained in the smallest affine space that 
contains $\wt F(x)$.}
\leqno (22.9)
$$
We proceed roughly as for (20.66).
Fix a face $F$ of our twisted grid, and first 
observe that by Theorem 6.2 on page 89 of [Ma], 
$$
\lim_{r \to 0} r^{-d} \H^d(E \cap B(x,r) \sm F) = 0.
\leqno (22.10)
$$
for $\H^d$-almost every $x\in F \cap E$. Then 
$$
\lim_{\rho \to 0} \rho^{-d} 
\H^d(\wt E \cap B(\wt x,\rho) \sm \wt F) = 0,
\ \hbox{ with } \wt F = \psi(\lambda F)
\leqno (22.11)
$$
(because $\psi$ is bilipschitz). Next notice that 
$\wt E \cap \wt F$ is rectifiable; hence for $\H^d$-almost every 
$x\in F \cap E$, $\wt E \cap \wt F$ has an approximate tangent
$\wt P'(x)$ at $\wt x$, which of course can be chosen
inside the affine span of $\wt F$. When (22.11) holds,
$\wt P'(x)$ is also an approximate tangent plane to the whole
$\wt E$ (the additional part has vanishing density).
By local Ahlfors-regularity of $\wt E$, $\wt P'(x)$ is even 
a true tangent plane to $\wt E$. It is easy to see that
for local Ahlfors-regular sets, the tangent plane is unique,
so $\wt P'(x) = \wt P(x)$ almost everywhere on $F$.
Since there is only a finite number of faces to try, we get that
$\wt P(x)$ is contained in the affine span of $\wt F$ for 
$\H^d$-almost every $x\in E$ and all the faces $F$ that
contain $x$; we apply this to $F = F(x)$ and get (22.9).

We don't even need to know that the tangent plane is unique
to make the argument work, because we just need to find, for
almost every $x \in E \cap H$, a tangent plane that satisfies
(22.9); so we could use the plane $\wt P'(x)$ associated
to $F(x)$, for instance.

Observe also that the set of points $x \in E \cap H$ for which
the dimension of $F(x)$ is less than $d$ is excluded by 
(22.9); this is all right, because this set is $\H^d$-negligible.

We also exclude the exceptional set $Z$ of (19.35). That is,
let us denote by $X$ the set of points $x\in E \cap H$
that satisfy the conditions (22.6)-(22.9) above, and in addition,
if $x$ is contained in one of the sets 
$L'_i = L_i \sm {\rm int}(L_i)$, where ${\rm int}(L_i)$ denotes the
true ($n$-dimensional) interior of $L_i$, and
$$
\limsup_{r \to 0} A_r(x) \leq \omega_d,
\leqno (22.12)
$$
where $A_r(x)$ is given by (19.33). Thus, by the discussion above,
$$
\H^d(E \cap H \sm X) = 0.
\leqno (22.13)
$$

For the next stage of the proof, we select a small radius $r(x)$
for every $x\in X$.  we choose $r(x)$ so that
$$
r(x) \leq {1 \over 4\Lambda^{2}} 
\min(\dist(x, U \sm V), \dist(x,  \d F(x)))
\leqno (22.14)
$$
(which is positive because $x\in E\cap H \i V = {\rm int}(V)$ and 
$x$ lies in the (face) interior of $F(x)$),
$$
\omega_d - \varepsilon \leq r^{-d} \H^d(E \cap B(x,r)) 
\leq \omega_d + \varepsilon
\ \hbox{ for } 0 < r \leq r(x)
\leqno (22.15)
$$
(possible by (22.6)), 
$$
\dist(z, P(x)) \leq \varepsilon r
\ \hbox{ for } 
z\in E \cap B(x,2r) 
\hbox{ and } 0 < r \leq r(x),
\leqno (22.16)
$$
$$
\dist(\wt z, \wt P(x)) \leq \varepsilon \lambda r
\ \hbox{ for } 
\wt z \in \wt E \cap B(\wt x, 2\lambda \Lambda r) 
\hbox{ and } 0 < r \leq r(x),
\leqno (22.17)
$$
and, when $x$ lies in some $L'_i = L_i \sm {\rm int}(L_i)$,
$$
A_r(x) \leq \omega_d + \varepsilon
\ \hbox{ for } 0 \leq r \leq r(x).
\leqno (22.18)
$$
We add two constraints that will simplify our life
when we check the boundary condition (1.7). 
We require that for each $i\in [0,j_{max}]$,
$$
r(x) < {1 \over 4\Lambda^{2}} \dist(x, \R^n \sm L_i)
\leqno (22.19)
$$
when $x$ lies in the ($n$-dimensional) interior of $L_i$, and
on the opposite
$$
r(x) < {1 \over 2} \dist(x, L_i)
\leqno (22.20) 
$$
when $x \in U \sm L_i$.

Let us apply Theorem 2.8 in [Ma] 
to the family of balls $\overline B(x,r)$,
$x \in X$ and $0 < r < \min(r(x), \rho_0)$,
where $\rho_0$ will be chosen later, and such that
$\H^d (E\cap \d B(x,r)) = 0$. We get a collection of disjoint
$B_j = B(x_j,r_j)$, $j \in J_1$, such that
$$
0 < r_j < \min(r(y_j), \rho_0),
\leqno (22.21)
$$ 
$\H^d(E \cap \d B_j) = 0$ for all $j$, and
$$
\H^d(X \sm \bigcup_{j \in J_1} B_j) = 
\H^d(X \sm \bigcup_{j \in J_1} \overline B_j) = 0.
\leqno (22.22)
$$
Let us choose a finite subset $J$ of $J_1$, so that
$$
\H^d(X \sm \bigcup_{j \in J} B_j) \leq \varepsilon.
\leqno (22.23)
$$
Set $X_1 = E \cap H \sm \bigcup_{j \in J} B_j$; then by (22.13),
$$
\H^d(X_1) = \H^d(X \sm \bigcup_{j \in J} B_j) \leq \varepsilon
\leqno (22.24)
$$
and we can use the definition of $\H^d$ to cover $X_1$
by balls $B_i = B(x_i,r_i)$, $i \in I$, so that 
$$
\hbox{$r_i \leq \rho_0$ for $i\in I \ $ and } \ 
\sum_{i\in I} r_i^d \leq C \varepsilon.
\leqno (22.25)
$$
Because $X_1$ is compact, we can replace $I$ with a finite 
subset for which the $B_i$ still cover $X_1$ and (removing the useless
balls) each $B_i$ meets $X_1$. By definition,
$$
E \cap H \i \bigcup_{j \in I \cup J} B_j.
\leqno (22.26)
$$
Since $E \cap H$ is compact, $I \cup J$ is finite, and $\{ E_k \}$ converges to $E$,
we also get that
$$
E_k \cap H \i \bigcup_{j \in I \cup J} B_j
\ \hbox{ for $k$ large enough.}
\leqno (22.27)
$$
For each $i\in I$ pick $y_i \in E \cap B_i$. Then for $k$ large, we 
can find $y_{i,k} \in E_k \cap B_i$ for every $i\in I$, and of course
$B_i \i B(y_{i,k},2r_i)$. We shall choose 
$$
\rho_0 < {1 \over 10 \Lambda^2} 
\min(\dist(H,U \sm V), \lambda^{-1} r_0, \delta),
\leqno (22.28) 
$$
where the constants $\lambda$ and $r_0$ come from Definition 2.7
(the Lipschitz assumption), and $\delta$ comes from our 
$GSAQ(U,M,\delta,h)$ assumption. We don't care how small they are, 
the main point is that they depend only on $E$ and the sequence $\{ E_k \}$.
Then $B(y_{i,k},4r_i) \i V \i U$, and by (10.11) (the uniform local Ahlfors-regularity 
of the  $E_k$),
$$
\H^d(E_k \cap B_i) \leq \H^d(E_k \cap B(y_{i,k},2r_i)) \leq C r_i^d.
\leqno (22.29)
$$
This holds for $k$ large enough (and all $i\in I$), with a constant that 
depends only on $E$ and $\{ E_k \}$. By (22.25), this yields
$$
\H^d(E_k \cap \bigcup_{i\in I} B_i) \leq C \varepsilon
\leqno (22.30)
$$
for $k$  large, and we are left with the contributions of the balls
$B_j$, $j\in J$.

We need to use the quasiminimality of $E_k$, and for this we construct
a one parameter family of mappings $\{\varphi_{j,t} \}$, 
$0 \leq t \leq 1$, for each $j\in J$. 

Fix $j\in J$ for the moment, and define the cut-off function $\xi_j$
by $\xi_j(y) = 0$ for $y\in U \sm B_j$, and
$$
\xi_j(y) = \min\big\{ 1, (\tau r_j)^{-1} \dist(y,\d B_j)\big\}
\ \hbox{ for } y\in B_j,
\leqno (22.31)
$$
where the small constant $\tau > 0$ will be chosen later
(before $\varepsilon$ and $\rho_0$).

We start with the easier case when when $x_j$ does not lie 
in any $L'_i$. Then we pick a unit vector $v_j$ parallel to $P(x_j)$, 
and set
$$
\varphi_{j,t}(x) = x + t \xi_j(x) [\pi_j(x) - x + \eta r_j v_j]
\leqno (22.32)
$$
for $y\in B_j$ and $0 \leq t \leq 1$, where 
$\pi_j$  denotes the orthogonal projection onto $P(x_j)$
and $\eta > 0$ is a minuscule constant, to be chosen later
(depending on $\tau$  and $\varepsilon$). We do nothing on 
$U \sm B_j$, i.e., set
$$
\varphi_{j,t}(x) = x 
\ \hbox{ for $x\in U \sm B_j$ and } 0 \leq t \leq 1.
\leqno (22.33)
$$
This is also the formula that we would use under the rigid assumption,
but because of (1.7) (and under the Lipschitz assumption) we shall need to be more 
careful in our second case. 

If $x_j \in L'_i$ for some $i$, we proceed as we did near (20.80). 
Denote by $F_j = F(x_j)$ the smallest face of our twisted grid that 
contains $x_j$, set $\wt F_j = \psi(\lambda F_j)$ (a rigid face), 
call $\wt W_j$ the affine space spanned by $\wt F_j$, 
choose a unit vector $\wt v_j$ in the vector space parallel to $\wt P(x_j)$ 
(the approximate tangent plane to $\wt E = \psi(\lambda E)$ at
$\wt x_j = \psi(\lambda x_j)$). Notice that by (22.9), 
$\wt v_j$ also lies in the vector space parallel to
$\wt W_j$ and $\wt F_j$. Denote by $\wt \pi_j$
the orthogonal projection onto $\wt P_j$, and finally set
$$
\wt \varphi_{j,t}(x) = \psi(\lambda x) 
+ t \xi_j(x) \big[\wt\pi_j(\psi(\lambda x))-\psi(\lambda x)
+ \eta \lambda r_j \wt v_j \big]
\leqno (22.34)
$$
for $x\in B_j$ and $0 \leq t \leq 1$. Notice that
$$\eqalign{
|\wt \varphi_{j,t}(x) - \psi(\lambda x_j)|
& \leq |\psi(\lambda x)-\psi(\lambda x_j)|
+ t|\wt\pi_j(\psi(\lambda x))-\psi(\lambda x)| + t\eta \lambda r_j
\cr&
\leq 2|\psi(\lambda x)-\psi(\lambda x_j)| + \eta \lambda r_j
\leq 2\lambda \Lambda r_j + \eta \lambda r_j \leq 3 \lambda \Lambda r_j
}\leqno (22.35)
$$
because $\psi(\lambda x_j)$ lies in $\wt P_j$ and if
$\eta$ is small enough. Then by (22.14) and (22.21),
$$\eqalign{
\dist(\wt \varphi_{j,t}(x),\R^n \sm \psi(\lambda V)) 
&\geq \dist(\psi(\lambda x_j),\R^n \sm \psi(\lambda V)) 
- 3 \lambda \Lambda r_j
\cr&
\geq \lambda \Lambda^{-1} \dist(x_j,U \sm V) - 3 \lambda \Lambda r_j
\cr&
\geq 4 \lambda \Lambda r(x_j) - 3\lambda \Lambda r_j
\geq \lambda \Lambda r_j.
}\leqno (22.36)
$$
In particular $\wt \varphi_{j,t}(x) \in \psi(\lambda V) \i B(0,1)$
and we can define 
$$
\varphi_{j,t}(x) = \lambda^{-1} \psi^{-1}(\wt \varphi_{j,t}(x)) \in V.
\leqno (22.37)
$$
Observe that by (22.35)
$$
|\varphi_{j,t}(x) -x_j| 
\leq \lambda^{-1} \Lambda |\wt \varphi_{j,t}(x) - \psi(\lambda x_j)|
\leq 3 \Lambda^2 r_j.
\leqno (22.38)
$$

This completes our definition of the $\varphi_{j,t}$ on $B_j$, and 
naturally we keep the trivial definition (22.33) on $U \sm B_j$.
Our next task is to show that the $\varphi_{j,t}$, $0 \leq t \leq 1$,
define an acceptable competitor for $E_k$. There is no problem with
(1.4) and (1.8); our mappings $\varphi_{j,t}(x)$, $t \leq 1$, are 
clearly continuous in $x$ and $t$ and Lipschitz in $x$. 

Notice that $\varphi_{j,t}(x) = x$ when $t=0$.
By (22.33) $\varphi_{j,t}(x) = x$ when $x\in U \sm B_j$.
When $x \in B_j$ and we use (22.32), notice that
$x + t \xi_j(x) [\pi_j(x) - x] \in B_j$
(because $\pi_j(x) \in B_j$ and $B_j$ is convex), so
$\varphi_{j,t}(x) \in B(x_j,(1+\eta) r_j)$. 
When $x \in B_j$ but we use (22.34) and (22.37),
we only get that $\varphi_{j,t}(x) \in B(x_j,3\Lambda^2 r_j)$,
by (22.38). In both cases,
$$
\varphi_{j,t}(B_j) \i B(x_j,3\Lambda^2 r_j),
\leqno (22.39)
$$
and the $\varphi_{j,t}$ satisfy (1.5) and (1.6), relative to the ball
$B=\overline B(x_j,3\Lambda^2 r_j)$. They also satisfy (2.4) with
$\wh W \i B$, which is compact and contained in $V \i U$
by  (22.14) and (22.21).

Finally let us check (1.7). We do this under the 
Lipschitz assumption; the rigid case is just simpler.
Let $i \leq j_{max}$ and $x\in E_k \cap L_i$ be given.
Recall that we want to check that $\varphi_{j,t}(x) \in L_i$ for
all $t$, so we may assume that $x\in B_j$, because otherwise
$\varphi_{j,t}(x)=x \in L_i$. 

A first case is when $x_i \in {\rm int}(L_i)$. In this case,
(22.21) and (22.19) imply that $B(x_j,4\Lambda^2r_j) \i L_i$,
and then (22.39) implies that 
$\varphi_{j,t}(x) \in \varphi_{j,t}(B_j) \i L_i$,
as needed. The case when $x_j \notin L_i$ is impossible,
because (22.21) and (22.20) would imply that $B_j$ does not
meet $L_i$. So we are left with the case when $x_j \in L'_i$.
Recall that in this case $\varphi_{j,t}(x)$ was defined by 
(22.34) and (22.37).

Let $F$ be a face of $L_i$ that contains $x$, and
set $\wt F = \psi(\lambda F)$ and $\wt F(x_j) = \psi(\lambda F(x_j))$.
We first check that
$$
\wt F(x_j) \i  \wt F.
\leqno (22.40)
$$
Suppose not, first assume that $\wt F(x_j)$ is not reduced to
$\{ \psi(\lambda x_j) \}$, and apply (3.8), (22.14), and (22.21), to get that
$$\eqalign{
\dist(\psi(\lambda x_j),\wt F) &\geq \dist(\psi(\lambda x_j), \d(\wt F(x_j)))
\geq \lambda \Lambda^{-1}\dist(x_j, \d F(x_j))
\cr&\geq 4 \lambda \Lambda r(x_j)
\geq 4 \lambda \Lambda r_j.
}\leqno (22.41)
$$
If instead $\wt F(x_j) = \{ \psi(\lambda x_j) \}$, $\psi(\lambda x_j)$ 
is a vertex, so $\dist(\psi(\lambda x_j),\wt F) \geq r_0$
because $\wt F$ is a face that does not contain it, and the conclusion
of (22.41) still holds, by (22.21) and (22.28).
But $\psi(\lambda x) \in \wt F$ because $x\in F$, so
$$
\dist(\psi(\lambda x_j),\wt F) \leq |\psi(\lambda x_j)-\psi(\lambda x)| 
\leq \lambda \Lambda |x-x_j| \leq \lambda \Lambda r_j,
\leqno (22.42)
$$
a contradiction which proves (22.40). 

Return to  $\varphi_{j,t}(x)$, which was defined by (22.34) and 
(22.37); we want to show that $\varphi_{j,t}(x) \i L_i \,$, and (by 
definition of $F$) it is enough to show that $\varphi_{j,t}(x) \in F$,
or equivalently that $\wt \varphi_{j,t}(x) \in \wt F$.
Recall that $\wt \pi_j$ is the orthogonal projection onto $\wt P_j$,
which by (22.9) is contained in the affine span $\wt W_j$ of $\wt F(x_j)$,
and (by (22.40)) in the affine span $\wt W$ of $\wt F$.

Set $\wt x = \psi(\lambda x) \in \wt F$; we know that 
$\wt x \in \wt F$, and then, by (22.34),
$$
\wt \varphi_{j,t}(x) = \wt x + t \xi_j(x) \big[\wt\pi_j(\wt x)-\wt x
+ \eta \lambda r_j \wt v_j \big] \in \wt W
\leqno (22.43)
$$
because $\wt v_j$ was chosen in the vector space parallel to 
$\wt P_j \i \wt W$.
But $\psi(\lambda x_j) \in \wt F(x_j)$ and $\wt F(x_j) \i \wt F$, so
$$
\dist(\psi(\lambda x_j), \d \wt F)
\geq \dist(\psi(\lambda x_j), \d \wt F(x_j))
\geq 4 \lambda \Lambda r(x_j) \geq 4 \lambda \Lambda r_j
\leqno (22.44)
$$
by (22.14) and (22.21). Set $I = [\psi(\lambda x_j),\wt \varphi_{j,t}(x)]$.
By (22.35), its length is at most $3 \lambda \Lambda r_j$,
so (22.44) says that $\dist(I,\d \wt F) \geq \lambda \Lambda r_j$.
In particular, $I$ does not cross $\d \wt F$; since
its initial point is $\psi(\lambda x_j) \in \wt F(x_j) \i \wt F$
(by (22.40)), and $I \i \wt W$ (by (22.43)), we get that $I \i \wt F$.
Hence $\wt \varphi_{j,t}(x) \in \wt F$, as desired, and (1.7) follows.

\ms
We are now allowed to test the quasiminimality of $E_k$
on the $\varphi_{j,t}$; notice in particular that 
the radius of our ball $B=\overline B(x_j,3\Lambda^2 r_j)$
is smaller than the threshold $\delta > 0$ in our quasiminimality 
assumption (10.2), by (22.21) and (22.28). We get that
$$
\H^d(W_1(E_k)) \leq M \H^d(\varphi_{j,1} W_1(E_k)) + 3 h \Lambda^2 
r_j^d, 
\leqno (22.45)
$$
where
$$
W_1(E_k) = \big\{ x\in E_k \, ; \, \varphi_{j,1}(x) \neq x \big\}
\i E_k \cap B_j
\leqno (22.46)
$$
(by (22.33)). Next we estimate $\H^d(\varphi_{j,1} (E_k \cap B_j))$.
Let us first check that for $k$ large, 
$$
\hbox{$\varphi_{j,1}$ is $2\Lambda^2$-Lipschitz on $E_k \cap B_j$}.
\leqno (22.47)
$$
We start in the easier case when $\varphi_{j,1}$ was defined by (22.32).
Observe that for $x, y \in B_j$,
$$\eqalign{
\varphi_{j,1}(x) - \varphi_{j,1}(y)
&= x + \xi_j(x) [\pi_j(x) - x + \eta r_j v_j]
- y - \xi_j(y) [\pi_j(y) - y + \eta r_j v_j]
\cr&
= (x-y) + \xi_j(x)[\pi_j(x) - \pi_j(y) - (x-y)]
\cr&\hskip2.5cm
+ [\xi_j(x)-\xi_j(y)][ \pi_j(y) - y + \eta r_j v_j]
\cr&
= (1-\xi_j(x))(x-y) + \xi_j(x)[\pi_j(x) - \pi_j(y)]
\cr&\hskip2.5cm
+ [\xi_j(x)-\xi_j(y)][ \pi_j(y) - y + \eta r_j v_j]
}\leqno (22.48)
$$
and hence
$$\eqalign{
|\varphi_{j,1}(x) - \varphi_{j,1}(y)|
&\leq |x-y|+ |\xi_j(x)-\xi_j(y)|[|\pi_j(y) - y| + \eta r_j ]
\cr& 
\leq|x-y|+ (\tau r_j)^{-1} |x-y| [|\pi_j(y) - y| + \eta r_j ]
}\leqno (22.49)
$$
by (22.31). Since $\{ E_k \}$ converges to $E$, for 
$k$ large enough we have that for each $y\in E_k \cap B_j$, there
is a $z\in E$ such that $|z-y| \leq \varepsilon r_j$.
Then $z\in B(x_j,2r_j)$, and (since $r_j \leq r(x_j)$ by (22.21))
(22.16) says that $\dist(z, P(x)) \leq \varepsilon r_j$.
Hence $|\pi_j(y) - y| \leq |\pi_j(z) - z|+\varepsilon r_j
\leq 2\varepsilon r_j$; then (22.47) easily follows from
(22.49), if $\varepsilon$ and $\eta$ are small enough compared to 
$\tau$.

When $\varphi_{j,1}$ is given by (22.34) and (22.37),
it is enough to check that $\wt \varphi_{j,1}$ is 
$2\lambda \Lambda$-Lipschitz on $E_k \cap B_j$. 
Let $x, y \in E_k \cap B_j$ be given; by (22.34)
and the same computation as for (22.49),
$$\leqalignno{
|\wt \varphi_{j,1}(x) - \wt \varphi_{j,1}(y)|
&\leq |\psi(\lambda x)-\psi(\lambda y)|+ 
|\xi_j(x)-\xi_j(y)|\big[|\wt\pi_j(\psi(\lambda y)) - \psi(\lambda y)| 
+ \eta \lambda r_j \big]
\cr& 
\leq \lambda \Lambda |x-y|+ (\tau r_j)^{-1} |x-y| 
\big[|\wt\pi_j(\psi(\lambda y)) - \psi(\lambda y)| 
+ \eta \lambda r_j \big].
&(22.50)
}
$$
Let $z$ be, as before, a point of $E$ such that
$|z-y| \leq \varepsilon r_j$. Set $\wt z = \psi(\lambda z)$,
and notice that $|\wt z - \psi(\lambda y)| \leq \lambda\Lambda |z-y|
\leq \lambda\Lambda\varepsilon r_j$, and also
$|\wt z - \psi(\lambda x_j)| 
\leq \lambda\Lambda |z-x_j| < 2 \lambda \Lambda r_j$.
Since $\wt z \in \wt E = \psi(\lambda E)$, (22.17) yields
$\dist(\wt z, \wt P(x_j)) \leq \varepsilon \lambda r_j$. Thus
$$
\dist(\psi(\lambda y), \wt P(x_j)) 
\leq \dist(\wt z, \wt P(x_j)) + |\wt z - \psi(\lambda y)|
\leq \varepsilon \lambda r_j + \lambda\Lambda \varepsilon r_j
\leq 2\lambda\Lambda\varepsilon r_j
\leqno (22.51)
$$
for $y\in E_k \cap B_j$ (and $k$ large). Then (22.50) yields
$$
|\wt \varphi_{j,1}(x) - \wt \varphi_{j,1}(y)|
\leq \lambda \Lambda |x-y|+ (\tau r_j)^{-1} |x-y|
\big[2\lambda\Lambda\varepsilon r_j + \eta \lambda r_j \big],
\leqno (22.52)
$$
$\wt \varphi_{j,1}$ is $2\lambda \Lambda$-Lipschitz on $E_k \cap B_j$
(for $k$ large and if $\varepsilon$ and $\eta$ are small enough), 
and we get (22.47) in our last case.

Write $B_j = B_{j,int} \cup B_{j,ext}$, where
$B_{j,int} = \big\{y\in B_j ; \dist(y,\d B_j) < \tau r_j\big\}$,
$B_{j,ext} = B_j \sm B_{j,int}$, and $\tau > 0$ is still the small 
constant in the definition (22.31) of $\xi_j$. Observe that by
(22.15), 
$$\eqalign{
\H^d(E \cap B_{j,ext})
&= \H^d(E \cap B_j) - \H^d(E \cap B_{j,int})
\cr&
\leq (\omega_d + \varepsilon) r_j^d 
- (\omega_d - \varepsilon) (1-\tau)^d r_j^d
\leq C \tau r_j^d
}\leqno (22.53)
$$
if $\varepsilon$ is small enough compared to $\tau$.
Since $B_{j,ext}$ is closed, we can apply the weak lower semicontinuity
result of (10.14), and we get that for $k$ large
$$
\H^d(E_k \cap B_{j,ext}) \leq C_M\H^d(E \cap B_{j,ext}) + \tau r_j^d
\leq C \tau r_j^d;
\leqno (22.54)
$$
of course if the reader does not feel like using (10.14), he may
also use the flatness of $E\cap 2B_j$ (see (22.16)), the fact that
the $E_k$ converge to $E$, and the uniform Ahlfors regularity of the
$E_k^\ast$ near $B_j$ to get (22.54). Next, by (22.47),
$$
\H^d(\varphi_{j,1}(E_k \cap B_{j,ext})) \leq C \tau r_j^d,
\leqno (22.55)
$$
where we no longer write the dependence on $\Lambda$.

We are left with the contribution of $B_{j,int}$. 
If we used (22.32) to define $\varphi_{j,1}$, we get that 
$\varphi_{j,1}(x) = \pi_j(x) + \eta r_j v_j$ for $x\in B_{j,int}$,
because $\xi_j(x) = 1$ on $B_{j,int}$. Then
$$
\H^d(\varphi_{j,1} (E_k \cap B_{j,int})) 
=\H^d(\pi_j (E_k \cap B_{j,int})) 
\leq \H^d(P(x_j) \cap B_j)
\leq \omega_d r_j^d.
\leqno (22.56)
$$
If instead we used (22.34) and (22.37), observe that by (22.34),
$\wt \varphi_{j,1}(x) = \wt\pi_j(\psi(\lambda x)) + \eta \lambda r_j 
\wt v_j$ for $x\in E_k \cap B_{j,int}$. We chose $\wt v_j$ parallel to 
$\wt P(x_j)$, so $\wt \varphi_{j,1}(x) \in \wt P(x_j)$, and hence
$$
\varphi_{j,1}(x) \in \lambda^{-1} \psi^{-1}(\wt P(x_j)),
\leqno (22.57)
$$
by (22.37). In addition,
$$\eqalign{
|\varphi_{j,1}(x)-x| 
&\leq \lambda^{-1} \Lambda |\wt\varphi_{j,1}(x)-\psi(\lambda x)|
\leq \lambda^{-1} \Lambda \big[
|\wt\pi_j(\psi(\lambda x))-\psi(\lambda x)|+ \eta \lambda r_j \big]
\cr& \leq \lambda^{-1} \Lambda \dist(\psi(\lambda x),\wt P(x_j))
+\eta \Lambda^2 r_j 
\leq 2\Lambda^2\varepsilon r_j+ \Lambda^2 \eta r_j
}\leqno (22.58)
$$
by (22.37) again and (22.51). 

Since $x\in B_{j,int}$, we get that
$|x-x_j| \leq (1-\tau) r_j$, and then (22.58) implies that
$\varphi_{j,1}(x) \in B_j$, if $\varepsilon$ and $\eta$ are small
enough compared to $\tau$. Thus
$\varphi_{j,1}(x) \in B_j \cap \lambda^{-1} \psi^{-1}(\wt P(x_j))$.

But $\wt P(x_j)$ is a $d$-plane which is contained in the affine span 
of the face $\wt F_j = \psi(\lambda F(x_j)$
(see (22.9) or the description above (22.34)), so the definition
(19.33) yields
$$
\H^d(\varphi_{j,1} (E_k \cap B_{j,int}))
\leq \H^d(B_j \cap \lambda^{-1} \psi^{-1}(\wt P(x_j)))
\leq r_j^d A_{r_j}(x_j).
\leqno (22.59)
$$
In addition, we only used (22.34) and (22.37) when $x\in L'_i$ for 
some $i$, and then (22.18) (together with (22.21)) says that
$A_{r_j}(x_j) \leq \omega_d+\varepsilon$. Thus
$$
\H^d(\varphi_{j,1} (E_k \cap B_{j,int}))
\leq r_j^d A_{r_j}(x_j)
\leq (\omega_d+\varepsilon) r_j^d
\leq (1+ C\varepsilon)\, \omega_d r_j^d.
\leqno (22.60)
$$
Since (22.56) was a better estimate, we shall remember that (22.60) holds 
in all cases. We group this with (22.55) and get that
$$
\H^d(\varphi_{j,1}(E_k \cap B_j)) 
\leq C \tau r_j^d + (1+C\varepsilon)\, \omega_d r_j^d
\leq (1+C\tau)\omega_d r_j^d 
\leqno (22.61)
$$
(if $\varepsilon$ is small enough).

We now return to (22.45) and give a lower bound for $\H^d(W_1(E_k))$.
We claim that $E_k \cap B_{j,int} \i W_1(E_k)$. Let 
$x\in E_k \cap B_{j,int}$ be given. If we used (22.32),
$\varphi_{j,1}(x) = \pi_j(x) + \eta r_j v_j$ and
$\pi_j(\varphi_{j,1}(x)) = \pi_j(x) + \eta r_j v_j \neq \pi_j(x)$
because we chose $v_j$ parallel to $P(x_j)$. Then 
$\varphi_{j,1}(x) \neq x$ and $x\in W_1(E_k)$.
If instead we used (22.34) and (22.37), (22.34) yields
$\wt \varphi_{j,1}(x) = \wt\pi_j(\psi(\lambda x)) 
+ \eta \lambda r_j \wt v_j$, and now 
$\wt\pi_j(\wt \varphi_{j,1}(x)) = \wt\pi_j(\psi(\lambda x)) 
+ \eta \lambda r_j \wt v_j \neq \wt\pi_j(\psi(\lambda x))$
because $\wt v_j$ is parallel to $\wt P_j(x)$; in this case
$\wt \varphi_{j,1}(x) \neq \psi(\lambda x)$, hence
$\varphi_{j,1}(x) \neq x$ and $x\in W_1(E_k)$.
So $E_k \cap B_{j,int} \i W_1(E_k)$, and (22.45) yields
$$\eqalign{
\H^d(E_k \cap B_{j,int}) &\leq \H^d(W_1(E_k)) 
\leq M \H^d(\varphi_{j,1} W_1(E_k)) + 3 h \Lambda^2 r_j^d
\cr&\leq M \H^d(\varphi_{j,1}(E_k \cap B_j)) + 3 h \Lambda^2 r_j^d
}\leqno (22.62)
$$
by (22.46). By (22.54), (22.62), and (22.61),
$$\eqalign{
\H^d(E_k \cap B_j) &\leq \H^d(E_k \cap B_{j,int}) +C \tau r_j^d
\cr&\leq M \H^d(\varphi_{j,1}(E_k \cap B_j)) + 3 h \Lambda^2 r_j^d 
+ C \tau r_j^d
\cr& \leq \big[M \omega_d + 3 h \Lambda^2 + C \tau \big] r_j^d
\cr& \leq \big[M \omega_d + 3 h \Lambda^2 + C \tau \big]
(\omega_d - \varepsilon)^{-1} \H^d(E \cap B_j)
}\leqno (22.63)
$$
by (22.15). We sum over $j$ and get that
$$\eqalign{
\H^d(E_k \cap H) &\leq \H^d(E_k \cap \bigcup_{j\in I \cup J} B_j)
\leq C \varepsilon + \sum_{j\in J} \H^d(E_k \cap B_j)
\cr&
\leq C \varepsilon + \big[M \omega_d + 3 h \Lambda^2 + C \tau \big]
(\omega_d - \varepsilon)^{-1} \sum_{j\in J} \H^d(E \cap B_j)
\cr&
\leq C \varepsilon + \big[M \omega_d + 3 h \Lambda^2 + C \tau \big]
(\omega_d - \varepsilon)^{-1} \H^d(E \cap V)
\cr&
\leq C \varepsilon + \big[M \omega_d + 3 h \Lambda^2 + C \tau \big]
(\omega_d - \varepsilon)^{-1} [\H^d(E \cap H)+\varepsilon]
}\leqno (22.64)
$$
by (22.27), (22.30), (22.63), the fact that the $B_j$, $j\in J$ are
disjoint and contained in $V$ (by (22.14) and (22.21)), and finally (22.5).

This holds for all choices of small constants $\tau$ and 
$\varepsilon$, with $\varepsilon$ small enough (depending on $\tau$) 
and $k$ large enough (depending on $\tau$ and $\varepsilon$).
We let $k$ tend to $+\infty$ in (22.64), notice that we can take
$\tau$ and $\varepsilon$ as small as we want (depending on $E$,
the $E_k$, and $H$), and get (22.4).
This completes our proof of Lemma 22.3.
\qed

\msi{\bf Proof of Theorem 22.1.}
We start with the case when the $E_k$ satisfy the assumptions of 
Theorem 21.3. In both of our three cases, the sets $E_k$
lie in some fixed $GSAQ(U,M,\delta,h)$, and we shall be able to
apply Theorem 10.8 and Lemma 22.3. 

In our first case when the
$E_k$ are $A_+$-almost minimal, we can take $h=0$
and $M = 1+h(\delta)$ as close to $1$ as we want
(because (21.2) holds), only at the price of choosing
$\delta$ small. We do that, get (22.4) with $h=0$
and $M$ close to $1$, let $M$ tend to $0$, and get (22.2).

In the second case when the $E_k$ are $A$-almost minimal,
we take $M=1$ and $h=h(\delta)$, which is also as small
as we want by (21.2), apply Lemma 22.3, let $\delta$
tend to $0$, and get (22.2) as above.

In the third case, we just need to observe that the 
$E_k$ are also $A$-almost minimal, by the easy part of 
Proposition 20.9, and use the second case.

So we may assume now that the $E_k$ satisfy the assumptions of
Corollary 21.15. Because of (21.10) we can still apply the
results of Section 10, with some acceptable choice of
$M$ and $h > 0$, but we do not want to use Lemma 22.3
with these constants, because we want to get rid of $M$ and $h$.

Instead we want to use the proof of Lemma 22.3, change a little
bit the definitions and accounting at the end, and use our 
asymptotic minimality assumption (21.13) or (21.14).

The main difference will be that we want to group some
of the families $\varphi_{j,t}$. Recall that (21.13) and (21.14)
come with some $\delta > 0$. Let us pick a maximal family
$z_l$, $1 \leq l \leq L$, of points of $H$, so that the $y_l$
lie at mutual distances $\delta /10$ from each other.
Then do the construction of Lemma 22.3, with 
$\rho_0 < 10^{-2} \Lambda^{-2} \delta$ (in addition
to the constraint in (22.28)). This gives, in particular, 
a family of balls $B_j$, $j\in J$, and for $j\in J$
a family $\varphi_{j,t}$.

For each $l$, denote by $J_0(l)$ the set of indices
$j \in J$ such that $|x_j - y_l| \leq \delta /10$.
Each $j\in J$ lies in some $J_0(l)$, because otherwise
we could add $x_j \in H$ to our collection of points $y_l$
and that collection would not be maximal. We prefer to have 
disjoint sets, so we set
$$
J(l) = J_0(l) \sm \bigcup_{m < l} J(m)
\ \hbox{ for } 1 \leq l \leq L.
\leqno (22.65)
$$
Then we fix $l$, and define a new family $\{ \varphi_t \}$,
$0 \leq t \leq 1$, by 
$$
\varphi_t(x) = x \hbox{ for $0 \leq t \leq 1$ and }
x\in U \sm \bigcup_{j\in J(l)} B_j
\leqno (22.66)
$$
and
$$
\varphi_t(x) = \varphi_{j,t}(x)
\hbox{ for $0 \leq t \leq 1$ and } x\in B_j.
\leqno (22.67)
$$
The $B_j$ are disjoint, so there in no ambiguity in (22.67).
Also, the $\varphi_t$ are continuous across the natural boundaries
$\d B_j$, $j\in J(l)$, because the $\varphi_{j,t}$ are, and by 
(22.33). Because of this, the $\varphi_t$ satisfy the conditions
(1.4) and (1.8).

We first need to check that the family $\{ \varphi_t \}$
is $\delta$-admissible for $E_k$ (see the definition below (21.10)). 
Let $W_t(E_k)$ and $\wh W(E_k)$ be as in (21.11) and (21.12).
If $x\in W_t(E_k)$ for some $t\in [0,1]$, then $x \in B_j$ for
some $j\in J(l)$ (by (22.66)), and then
$\varphi_t(x) = \varphi_{j,t}(x) \in B(x_j,3\Lambda^2 r_j)$,
by (22.39). This means that
$$
\wh W(E_k) \i \bigcup_{j\in J(l)} B(x_j,3\Lambda^2 r_j)
\i B(y_l,\delta/3)
\leqno (22.68)
$$
because $r_j \leq \rho_0 \leq 10^{-2} \Lambda^{-2} \delta$
and $|x_j - y_l| \leq \delta /10$. Thus (1.4) and (1.5)
hold for the $E_k$, with respect to the ball 
$B = \overline B(y_l,\delta/3)$. Also, the 
$\overline B(x_j,3\Lambda^2 r_j)$ are still contained
in $V \i U$ (by (22.14) and (22.21)), so 
$\wh W(E_k) \i\i U$ and $\{ \varphi_t \}$ is $\delta$-admissible 
for $E_k$. We thus get that for each $\varepsilon > 0$
(we can keep the same one as before), we can find 
$k_0 \geq 0$ such that for $k \geq k_0$,
$$
\H^d(W_1(E_k)) \leq \H^d(\varphi_1(W_1(E_k))) + \varepsilon
\leqno (22.69)
$$
(as in (21.13) or
$$
\H^d(E_k\sm\varphi_1(E_k)) \leq \H^d(\varphi_1(E_k) \sm E_k) + 
\varepsilon.
\leqno (22.70)
$$
(as in (21.14).
We first assume that (22.69) holds for $k \geq k_0$, 
because it is closer to what we had for Lemma 22.3.
Then
$$\eqalign{
\H^d(\varphi_1(W_1(E_k))) 
&\leq \H^d(\varphi_1(E_k \cap \bigcup_{j\in J(l)}B_j)) 
\leq \sum_{j\in J(l)} \H^d(\varphi_1(E_k \cap B_j))
\cr&
= \sum_{j\in J(l)} \H^d(\varphi_{1,j}(E_k \cap B_j))
\leq \sum_{j\in J(l)} (1+C\tau)\omega_d r_j^d
}\leqno (22.71)
$$
because $W_1(E_k) \i \bigcup_{j\in J(l)}B_j$
by (22.66), and by (22.67) and (22.61).

On the other hand, $W_1(E_k)$ contains $E_k \cap B_{j,int}$
for $k$ large, by the proof above (22.62); hence
$$\leqalignno{
\sum_{j\in J(l)} \H^d(E_k \cap B_j)
&\leq \sum_{j\in J(l)} \big[\H^d(E_k \cap B_{j,int}) 
+ C \tau r_j^d \big]
\cr&
\leq \H^d(E_k \cap \bigcup_{j\in J(l)} B_{j,int}) 
+ C \sum_{j\in J(l)} \tau r_j^d 
\leq  \H^d(W_1(E_k))
+ C \sum_{j\in J(l)} \tau r_j^d 
\cr&
\leq \H^d(\varphi_1(W_1(E_k))) + \varepsilon
+ C \sum_{j\in J(l)} \tau r_j^d 
\leq \varepsilon + \sum_{j\in J(l)} (1+C\tau)\omega_d r_j^d
&(22.72)
}
$$
by (22.54), because the $B_j$ are disjoint, and by (22.69)
and (22.71). We sum this over $l$, use the fact that 
$J$ is the disjoint union of the $J(l)$, and get that
$$
\sum_{j\in J} \H^d(E_k \cap B_j)
\leq L \varepsilon + \sum_{j\in J} (1+C\tau)\omega_d r_j^d \, ;
\leqno (22.73)
$$
the extra $L$ will not disturb, because it is fixed as soon as
we know $\delta$, and we can choose $\varepsilon$ later.
We proceed as in (22.64):
$$\eqalign{
\H^d(E_k \cap H) &\leq \H^d(E_k \cap \bigcup_{j\in I \cup J} B_j)
\leq C \varepsilon + \sum_{j\in J} \H^d(E_k \cap B_j)
\cr&
\leq C \varepsilon + L \varepsilon + \sum_{j\in J} (1+C\tau)\omega_d r_j^d
\cr&
\leq C \varepsilon + L \varepsilon + \sum_{j\in J} 
(1+C\tau)\omega_d (\omega_d - \varepsilon)^{-1} \H^d(E \cap B_j)
\cr&
\leq C \varepsilon + L \varepsilon 
+ (1+C\tau){ \omega_d \over \omega_d - \varepsilon} \H^d(E \cap V)
\cr&
\leq C \varepsilon + L \varepsilon 
+ (1+C\tau){ \omega_d \over \omega_d - \varepsilon}
[\H^d(E \cap H)+\varepsilon]
}\leqno (22.74)
$$
by (22.27), (22.30), (22.73), (22.15), the fact that the $B_j$, 
$j\in J$ are disjoint and contained in $V$ (by (22.14) and (22.21)),
and (22.5). 

This is true for $\tau$ small enough, $\varepsilon$ small enough
(depending on $\tau$ as well), and $k$ large (depending on both).
We let $k$ tend to$+\infty$, and then let $\varepsilon$ and $\tau$
tend to $0$, and we get (22.2).

Now assume that (22.70) holds. The simplest is to show
that (22.69) holds as well. Set $W = W_1(E_k)$ and observe that
$$
W \sm \varphi_1(W) \i E_k \sm\varphi_1(E_k)
\leqno (22.75)
$$
because if $x\in W \sm \varphi_1(W)$, it lies in $E_k$, 
does not lie in $\varphi_1(W)$, and does not lie in 
$\varphi_1(E_k \sm W)$ either, because $\varphi_1(y) = y \notin W$
when $y \in E_k \sm W$. Similarly,
$$
\varphi_1(E_k) \sm E_k \i \varphi_1(W) \sm W
\leqno (22.76)
$$
because if $x\in \varphi_1(E_k) \sm E_k$ and $y\in E_k$
is such that $\varphi_1(y)=x$, then $y\in W$
(otherwise, $x = \varphi_1(y) = y \in E_k$) and of course
$x\notin W$ (because $W \i E_k$ and $x\notin E_k$). Then
$$\eqalign{
\H^d(W) &= \H^d(W \cap \varphi_1(W)) + \H^d(W \sm \varphi_1(W))
\cr& \leq \H^d(W \cap \varphi_1(W)) + \H^d(E_k \sm \varphi_1(E_k))
\cr&\leq \H^d(W \cap \varphi_1(W)) + \H^d(\varphi_1(E_k) \sm E_k) + \varepsilon
\cr&\leq \H^d(W \cap \varphi_1(W)) + \H^d(\varphi_1(W) \sm W) + \varepsilon
=  \H^d(\varphi_1(W)) + \varepsilon
}\leqno (22.77)
$$
by (22.75), (22.70), and (22.76), as needed for (22.69).
This completes our proof of (22.2) when (22.70) holds;
Theorem 22.1 follows.
\qed

\msi
{\bf 23. Limits of quasiminimal and almost minimal sets in variable domains.} 
\ms

The main result of this section is a variant of Theorem 10.8 where we 
allow the domain $U$ and the boundary pieces $L_i$ to vary slightly 
along the sequence. 
We shall not try to obtain an optimal result here
(this would probably involve following the long proof carefully),
and instead we shall state a result (Theorem 23.8)
that we can easily deduce from Theorem 10.8 by a change of variable.

Let us explain the notation for Theorem 23.8. 
We are given an open set $U$ (the limit) and boundary pieces
$L_j$, $0 \leq j \leq j_{max}$, and we assume as in (10.1) that
$$
\hbox{the Lipschitz assumption are satisfied in $U$.}
\leqno (23.1)
$$
But we also give ourselves a sequence $\{ U_k \}$ of open sets,
and for each $k \geq 0$ a collection of boundary pieces
$L_{j,k}$, $0 \leq j \leq j_{max}$. We shall make our life simpler
and assume that $U_k$ and the $L_{j,k}$ are parameterized by
$U$ and the $L_j$, using a single bilipschitz mapping
$\xi_k : U \to U_k$. That is,
$$
U_k = \xi_k(U) \ \hbox{ and } \ 
L_{j,k} = \xi_k(L_j) \ \hbox{ for } 0 \leq j \leq j_{max}.
\leqno (23.2)
$$
In addition, we assume that the $\xi_k$ become optimally bilipschitz, in 
the sense that there exist constants $\eta_k \geq 1$ such that
$$
\xi_k \hbox{ is $\eta_k$-bilipschitz on $U$ and }
\lim_{k \to +\infty} \eta_k = 1.
\leqno (23.3)
$$
We also assume that
$$
\lim_{k \to +\infty} \xi_k(x) = x
\hbox{ for } x\in U.
\leqno (23.4)
$$
These are quite strong assumptions on our sequence, but our main
example will be a blow-up sequence at a point of an initial
domain where each $L_j$ has a tangent cone, in which case the 
$\xi_k$ can be constructed by hand. See Section 24. 

Now we give ourselves a sequence $\{ E_k \}$ of quasiminimal sets.
We assume that the following properties hold for some choice of 
constants $M \geq 1$, $\delta \in (0,+\infty]$, and $h > 0$.
First, each $E_k$ is a relatively closed subset of the corresponding
set $U_k$, $E_k$ is coral (as in (10.3) and Definition 3.1), and
$$
E_k \in GSAQ(U_k,M,\delta,h),
\leqno (23.5)
$$
where of course we define $GSAQ(U_k,M,\delta,h)$ relative to the
sets $L_{j,k}$, $0 \leq j \leq j_{max}$. We assume, as in (10.4), that
for some relatively closed set $E \i U$,
$$
\lim_{k \to +\infty} E_k = E \hbox{ locally in $U$.}
\leqno (23.6)
$$
Since the $E_k$ are contained in slightly different domains $U_k$,
let us say what this means: for each compact set $H \i U$ and each
$\varepsilon >0$, we can find $k_0 \geq 0$ such that for $k \geq k_0$,
$$
\dist(x,E_k) \leq \varepsilon \hbox{ for every } x\in E \cap H
\ \hbox{ and }\,
\dist(x,E) \leq \varepsilon \hbox{ for every } x\in E_k \cap H.
\leqno (23.7)
$$

When the rigid assumption does not hold, we also assume that
(10.7) or (19.36) holds (in $U$, for the $L_j$).

\ms\proclaim Theorem 23.8. Let $U$, $\{ U_k \}$, the $L_j$,
the $L_{j,k}$, and $\{ E_k \}$ satisfy all the conditions above,
including (10.7) or (19.36) when the $L_j$ don't satisfy the
rigid assumption.
Also assume that $h$ is small enough, depending on $M$, $n$,
and the constant $\Lambda$ that comes from (23.1).
Then $E$ is coral, and 
$$
E \in GSAQ(U,M,\delta,h),
\leqno (23.9)
$$
with the same constants $M$, $\delta$, and $h$ as in (23.5).

\ms
We first prove the theorem, and then comment later. 
Since we want to reduce to a fixed domain, we consider the sets 
$$
\wt E_k = \xi_k^{-1}(E_k) \i U,
\leqno (23.10)
$$
and we want to apply Theorem 10.8 to the sequence $\{ \wt E_k \}$. 
Since $\xi_k$ is bilipschitz, $\wt E_k$ is closed in $U$ and coral.
We claim that
$$
\wt E_k \in GSAQ(U,\eta_k^{2d}M,\eta_k^{-1}\delta,\eta_k^{2d}h),
\leqno (23.11)
$$
where $GSAQ$ is defined in terms of the boundaries $L_j$.
The proof is the same as for Proposition 2.8, so we won't repeat all
the verifications, but just give the outline. 
For each one parameter family of mappings
$\wt \varphi_t : \wt E_k \to \R^n$, $0 \leq t \leq 1$, 
that satisfies (1.4)-(1.8) with a radius $\wt r < \eta_k^{-1}\delta$
and (2.4), we observe that
$\wt \varphi_t(\wt E_k) \i U$ by (2.4), and so
we can define mappings $\varphi_t : E_k \to U_k$ by
$$
\varphi_t = \xi_k\circ \wt \varphi_t \circ \xi_k^{-1}.
\leqno (23.12)
$$
It is easy to see that the properties (1.4)-(1.8) and (2.4) for
the $\varphi_t$ (relative to $U_k$ and the $L_{j,k}$, and with a
radius $r \leq \eta_k \wt r$ smaller than $\delta$) follow from their 
counterpart for the $\wt\varphi_t$, so we can apply (2.5) to $\varphi_1$; 
this yields
$$
\H^d(W_1) \leq M \H^d(\varphi_1(W_1)) + h r^d,
\leqno (23.13)
$$
with $W_1 = \big\{x\in E_k \, ; \, \varphi_1(x) \neq x \big\}$.
The analogue of $W_1$ for the $\wt \varphi_t$ is 
$$
\wt W_1 = \big\{x\in \wt E_k \, ; \, \wt \varphi_1(x) \neq x \big\}
= \xi_k^{-1}(W_1),
\leqno (23.14)
$$
so (23.13) yields
$$
\H^d(\wt W_1) \leq \eta_k^{d} \H^d(W_1)
\leq \eta_k^{d} M \H^d(\varphi_1(W_1)) + \eta_k^{d} h r^d
\leq \eta_k^{2d} M \H^d(\wt\varphi_1(\wt W_1)) 
+ \eta_k^{2d} h \wt r^d,
\leqno (23.15)
$$
as needed for (23.11). 

We shall want to check that
$$
\lim_{k \to +\infty} \wt E_k = E
\ \hbox{ locally in } U,
\leqno (23.16)
$$
but we start with simple consequences of (23.3)
and (23.4). First notice that (23.4) implies the
apparently stronger fact that
$$
\lim_{k \to +\infty} \xi_k(x) = x
\hbox{ uniformly on every compact subset of } U,
\leqno (23.17)
$$
simply because (23.3) says that the $\xi_k$ are Lipschitz
on $U$, with uniform bounds. Next we check that for every compact
set $H \i U$, we can find a compact set $K \i U$ and an integer
$k_0$ such that
$$
H \i \xi_k(K) \ \hbox{ for } k \geq k_0.
\leqno (23.18)
$$
Set $r = {1 \over 4} \dist(H,\R^n \sm U) > 0$
and cover $H$ by a finite number of balls $B(x_i,r)$. 
For each $i$, we (23.17) says that for $k$ large,
$|\xi_k(y)-y| \leq r$ for $y\in \d B(x_i,3r)$.
For such $k$, and each $x\in B(x_i,r)$, the mapping
from $\d B(x_i,3r)$ to the unit sphere which maps
$z\in \d B(x_i,3r)$ to ${\xi_k(z)-x \over |\xi_k(z)-x|}$
is well defined (because $|\xi_k(z)-x| \geq |z-x|-r \geq r$),
continuous, and of degree $1$ because we can easily find a
homotopy from this map to to the map $z \to {|z-x| \over |z-x|}$
(among continuous mappings : $\d B(x_i,3r) \to \d B(0,1)$).
Then there is no homotopy from it to a constant, which implies
that $x\in \xi_k(B(x_i,3r))$  (because otherwise we could use
the mappings $z \to {\xi_k(z_t)-x \over |\xi_k(z_t)-x|}$, with 
$z_t = t x_i + (1-t)(z-x_i)$. Set $K = \cup_i \overline B(x_i,3r)$;
we just proved that $\xi_k(K)$ contains all the $B(x_i,r)$, 
and (23.8) follows.

We are now ready to check that $\{\wt E_k\}$ converges to $E$, 
as in (23.16). Let $H \i U$ be compact, set
$d_H = \dist(H,\R^n \sm U)$, and let $\varepsilon \in (0,d_H/2)$
be given. First we want to show that if $k$ large enough,
then 
$$
\hbox{for every $x\in E \cap H$ we can find $\wt y_k\in \wt E_k$ 
such that $|x-y_k| \leq 2\varepsilon$.}
\leqno (23.19)
$$
By (23.7), we can find $y_k\in E_k$ such that 
$|x-y_k| \leq \varepsilon$. Notice that $y_k$ lies in
the compact set $H_1 = \big\{ y\in \R^n \, ; \, \dist(y,H) \leq 
d_H/2 \big\}$. By (23.18), there is a compact set $K$ such that 
$\xi_k(K)$ contains $H_1$ for $k$ large. 

Set $\wt y_k = \xi_k^{-1}(y_k)$. We already knew that
$\xi_k^{-1}(y_k)$ is defined and lies in $\wt E_k$, because
$E_k \i U_k = \xi_k(U)$ and $\wt E_k = \xi_k^{-1}(E_k)$;
now we also know that $\wt y_k \in K$ for $k$ large, and 
(23.17) says that for $k$ large, 
$|\wt y_k-y_k| = |\wt y_k-\xi_k(\wt y_k)| \leq \varepsilon$.
Thus $\wt y_k\in \wt E_k$ and $|x-y_k| \leq 2\varepsilon$;
this proves (23.19).

Conversely, if $k$ is large enough, then for each 
$x_k \in \wt E_k \cap H$ (23.17) says that 
$|\xi_k(x_k)-x_k| \leq \varepsilon$; but 
$\xi_k(x_k) \in E_k \in H_1$ (by (23.10) and the definition 
of $H_1$), and (23.7) gives (again if $k$ is large enough) 
a point $x\in E$ such that $|x-\xi_k(x_k)| \leq \varepsilon$.
Thus $|x-x_k| \leq 2\varepsilon$.
This completes our proof of (23.16).

We may now return to Theorem 23.8. If we prove that 
$$
E \in GSAQ(U,M',\delta',h')
\leqno (23.20)
$$
for every choice of $M' > M$, $\delta' < \delta$, and $h' > h$,
the desired conclusion (23.9) will follow; this comes from the way
we defined $GSAQ(U,M,\delta,h)$ in Definition 2.3. 
So it is enough to show that the $\wt E_k$ satisfy the assumptions 
of Theorem 10.8, with any such choice of $M',\delta',h'$, 
and where the limit in (10.4) is still $E$.

But (10.1) is the same as (23.1), (10.2) follows from (23.11)
(if we use the fact that $\eta_k$ tends to $1$ by (23.3), 
and restrict to $k$ large), and (10.4) follows from
(23.16). Since we also assumed (10.7) or (19.36) if the rigid 
assumption does not hold, we can apply Theorem~10.8 or Remark 19.52
and get (23.20). Theorem 23.8 follows.
\qed

\ms
We end this section with a few comments and extensions of Theorem 23.8.

\msi{\bf Remark 23.21.}
Our bilipschitz assumption (23.3) could easily be made more local. That is,
if instead of (23.3) we assume that each $\xi_k: U \to U_k$ is bilipschitz, 
and that for each compact set $K \i U$ and each $\eta > 1$ we can find 
$k_0 \geq 0$ such that 
$$
\eta^{-1} |x-y| \leq |\xi_k(x)-\xi_k(y)| \leq \eta |x-y|
\ \hbox{ for $k \geq k_0$ and } x, y \in K,
\leqno (23.22)
$$
we still get the same conclusion. This is easy, because given a
competitor (i.e., a one parameter family $\{ \varphi_t \}$ with the
usual conditions), we can restrict our attention to a relatively
compact open set $V$ such that $\wh W \i V \i\i U$, and then
apply Theorem 23.8 with the open set $V$. 

To be honest, even though we claim that the proof still works, the
situation is a little more complicated than this, 
because $V$, together with the boundary sets $L_{j} \cap V$, may not
satisfy the Lipschitz assumption as in (23.1). This gets better if we 
choose $V$, possibly a little larger, so that $\psi(\lambda V)$ 
is a ball centered at the origin (and with a radius close to $1$). 
Without saying more, we can observe that our proof of Theorem 10.8
still works in that case (we never used the fact that the radius of
$B(0,1)$ is a dyadic number). Or, on a more formal level, we can also
make sure that $\psi(\lambda V) = B(0,1-2^{-m})$ for some integer $m$,
and observe that the sets $(1-2^{-m})^{-1}\psi(L_j \cap V)$ satisfy the rigid
assumption, although perhaps with a grid with a $2^m$ times smaller 
mesh. This has a small incidence on the statement of Theorem 10.8, so 
we can get an acceptable formal proof this way.

\msi{\bf Remark 23.23.}
The lower semicontinuity of $\H^d$, as in Theorem 10.97, is still true
under the assumptions of Theorem 23.8. That is, under the assumptions
of Theorem 23.8, we also have that
$$
\H^d(E\cap V) \leq \liminf_{k \to +\infty} \, \H^d(E_k \cap V)
\ \hbox{ for every open set $V \i U$.}
\leqno (23.24)
$$
(as in (10.98)). Indeed, for this it is enough to show that
if $V \i U$ is open, and if $H$ is any compact subset of $V$,
$$
\H^d(E\cap H) \leq \liminf_{k \to +\infty} \, \H^d(E_k \cap V),
\leqno (23.25)
$$
just because we can write 
$\H^d(E\cap V) = \lim_{l \to +\infty} \H^d(E\cap H_l)$, where
$\{ H_l \}$ is an increasing sequence of compact subsets of $V$.

Let $V$ and $H$ be given, and let $W$ be an open set
that contains $H$ and is relatively compact in $V$.
We have seen in the proof of Theorem 23.8
that the sets $\wt E_k \cap V = \xi_k^{-1}(E_k) \cap V$ 
satisfy the assumptions of Theorem 10.8 and Theorem 10.97
for some acceptable choices of constants (i.e., with $h$ small enough), 
and with the same limit $E$, so Theorem 10.97 yields
$$
\H^d(E\cap H) \leq \liminf_{k \to +\infty} \, \H^d(\wt E_k \cap W).
\leqno (23.26)
$$
For $k$ large, $\xi_k(\wt E_k \cap W) = \xi_k(\wt E_k) \cap \xi_k(W) 
\i E_k \cap V$ (by (23.10) and because (23.17) says that the $\xi_k$ 
converge to the identity uniformly on $\overline W$), and
$$
\H^d(\wt E_k \cap W)
\leq \eta_k^{d} \H^d(\xi_k(\wt E_k \cap W))
\leq \eta_k^{d} \H^d(E_k \cap V)
\leqno (23.27)
$$
by (23.3); now (23.25) and then (23.24) follow by letting
$k$ tend to $+\infty$.

\msi{\bf Remark 23.28 (limits of almost minimal sets).}
Our statements about limits of almost minimal sets, namely
Theorem~21.3 and Corollary 21.15, also hold in the present setting.

We start with Theorem~21.3. If, in Theorem 23.8 above, we replace the
quasiminimality assumption (23.5) with the assumption that
$$\eqalign{
&E_k \hbox{ is an almost minimal set in $U_k$, with the sliding }
\cr& \hskip 0.2cm
\hbox{conditions given by the $L_{j,k}$ and the gauge function } h,
}\leqno (23.29)
$$
where $h$ is a gauge function such that (21.1) and (21.2) hold
(otherwise, see Remarks 21.6 and 21.7), and one type of almost 
minimal set ($A_+$, $A$, or $A'$) is chosen.
The conclusion is then, as in Theorem 21.3, that the limit set
$E$ is an almost minimal set in $U$, with the sliding 
conditions given by the $L_{j}$, the gauge function $h$,
and the same type ($A_+$, $A$, or $A'$).
The proof consists in following our short proof of 
Theorem 21.3, and applying Theorem~23.8 instead of Theorem 10.8.

Similarly, if we replace the assumption (23.5) with the assumptions
(21.10)-(21.14) (and taken, for $E_k$, in the domain $U_k$ and 
relative to the boundary sets $L_{j,k}$), we get the same conclusion 
as in Corollary 21.15, namely that $E$ is a coral local minimal set 
in $U$, that satisfies (21.16).
Once again, we just follow the proof of Corollary 21.15.

\msi{\bf Remark 23.30 (upper semicontinuity of $\H^d$).}
The upper semicontinuity results of Section 22 can also be
generalized in the context of slowly changing domains.
We start with an extension of Lemma 22.3.

\ms\proclaim Lemma 23.31.
Let $U$, the $L_i$, the $U_k$ the $L_{j,k}$, the $\{ E_k \}$, and 
the limit $E$ satisfy the assumptions of Theorem 10.8. 
Then for every compact set $H \i U$,
$$
(1+Ch) M \H^d(E \cap H) \geq \limsup_{k \to +\infty} \H^d(E_k \cap H),
\leqno (23.32)
$$
with a constant $C$ that depends only on $n$, $M$, and $\Lambda$.

\ms
We proceed and as in Theorem~23.8, to deduce Lemma 23.31 from
Lemma 22.3 and a change of variable. We already proved near (23.20)
that for each choice of $M'>M$, $\delta'<\delta$, and $h'>h$,
the end of the sequence $\{ \wt E_k \}$ satisfies the hypotheses
of Theorem~10.8 or Remark 19.52. Then (if $h$ is small enough)
they also satisfy the conclusion of Lemma~22.3: 
for each compact set $K \i U$,
$$
(1+Ch) M' \H^d(E \cap K) 
\geq \limsup_{k \to +\infty} \H^d(\wt E_k \cap K),
\leqno (23.33)
$$
with a constant $C$ which we can still bound 
in terms of $n$, $M$, and $\Lambda$ (the constant $C$ that
is associated to $2M$ also works for all $M' \leq 2M$).

Return to (23.32), let $H \i U$ be compact, and let $\varepsilon > 0$
be given. Pick a compact set $K \i U$ such that $H$ is 
contained in the interior of $K$ and
$\H^d(E \cap K) \leq \H^d(E \cap H) + \varepsilon$. 
Observe that for $k$ large, 
$$
\xi_k(E_k \cap H) = \xi_k(E_k) \cap \xi_k(H)
= \wt E_k \cap \xi_k(H) \i \wt E_k \cap K
\leqno (23.34)
$$
by (23.10) and because by (23.17) the $\xi_k$ converge 
to the identity uniformly on $H$. Thus for $k$ large,
$$\eqalign{
\H^d(E_k \cap H) 
&\leq \eta_k^d \H^d(\xi_k(E_k \cap H))
\leq \eta_k^d \H^d(\wt E_k \cap K) 
\cr&
\leq \eta_k^d [(1+Ch) M' \H^d(E \cap K) + \varepsilon]
\cr&
\leq \eta_k^d [(1+Ch) M' (\H^d(E \cap H) + \varepsilon) + \varepsilon]
}\leqno (23.35)
$$
by (23.34). We let $k$ tend to $+\infty$, use (23.3),
let $\varepsilon$ tend to $0$, and get (23.29).
We take the $limsup$ of this, recall that $\eta_k$ tends to $1$
by (23.3), then observe that $\varepsilon$ and $M'-M$ are arbitrarily 
small, and get (23.32)
\qed

\ms
Here is a generalization of Theorem 22.1.

\ms\proclaim Lemma 23.36.
Let $U$, the $L_i$, the $U_k$ the $L_{j,k}$, the $\{ E_k \}$, and 
$E$ satisfy the same assumptions as in 
Theorem~21.3 and Corollary 21.15, but modified as in 
Remark 23.28. Then for every compact set $H \i U$,
$$
\limsup_{k \to +\infty} \H^d(E_k \cap H) \leq \H^d(E \cap H).
\leqno (23.37)
$$

\ms
In the case of Theorem~21.3, we can just follow the proof of
Theorem~21.3, and replace Lemma 22.3 with Lemma 23.31 when needed.
In the case of Corollary 21.15, the simplest seems to show that
the sets $\wt E_k$ of (23.10) satisfy the assumption of
Corollary 21.15, and then compute as in Lemma 23.31.
We claim that the verifications are quite similar to what 
we did for Theorem 23.8 and Lemma 23.31, and we skip them.
\qed

\msi
{\bf 24. Blow-up limits.} 
\ms

In this section we apply the results of Section 23 to the case of 
blow-up limits of an almost minimal set $E$, at a point $x_0 \in E$
near which the boundary pieces $L_j$ behave in a roughly $C^1$
way.

We fix an open set $U \i \R^n$, boundary pieces 
$L_j$, $0 \leq j \leq j_{max}$, a point $x_0 \in U$,
and a quasiminimal (or almost minimal) set $E \i U$.
We shall systematically assume that
$$
\hbox{$U$ and the $L_j$ satisfy the Lipschitz assumption}
\leqno (24.1)
$$
(see Definition 2.7), and things will be more interesting
when $E$ is assumed to be coral and we take $x_0 \in E$.

We are also given a sequence $\{ r_k \}$, with 
$$
\lim_{k \to +\infty} r_k = 0,
\leqno (24.2)
$$
along which we want to define a blow up. We define the sets $E_k$ by
$$
E_k = r_k^{-1} (E - x_0) = \big\{ z\in \R^n \, ; \,
x_0 + r_k z \in E \big\}.
\leqno (24.3)
$$
A simple computation shows that
$$
E_k \in QSAQ(U_k,M,r_k^{-1}\delta,h)
\hbox{ if }
E \in QSAQ(U,M,\delta,h),
\leqno (24.4)
$$
with $U_k = r_k^{-1} (U - x_0)$, and where on $U_k$ we use the 
boundary pieces
$$
L_{j,k} = r_k^{-1} (L_j - x_0), \, 0 \leq j \leq j_{max}.
\leqno (24.5)
$$ 

Similarly, if $E$ is almost minimal in $U$, with the sliding 
conditions coming from the $L_j$, and the gauge function $h$,
then $E_k$ is almost minimal in $U_k$, with the sliding 
conditions coming from the $L_{j,k}$, the gauge function $h_k$
defined by $h_k(t) = h(tr_k)$, and the same type of
almost minimality ($A_+$, $A$, or $A'$) as $E$; see Definition 20.2.
Notice that, by (24.2), the sets $U_k$ converge to $\R^n$.

A \underbar{blow-up sequence} for $E$ at the point $x_0$ is any sequence
$\{ E_k \}$ defined by (24.3), under the condition (24.2). 
In some cases, the sequence converges, i.e., there is closed set 
$E_\infty \i \R^n$ such that
$$
E_\infty = \lim_{k \to +\infty} E_k
\ \hbox{ locally in $\R^n$}.
\leqno (24.6)
$$
When this is the case, i.e., when (24.6) holds
for some sequence $\{ r_k \}$ that tends to $0$, we say that
$E_\infty$ is a \underbar{blow-up limit} of $E$ at $x_0$. Of course, 
different sequences may yield different blow-up limits $E_\infty$, 
even though in some cases one can prove that $E$ has only one
blow-up limit at $x_0$.

By general compactness arguments, we can always extract, from any 
blow-up sequence for $E$ at the point $x_0$, a convergent subsequence.
Thus $E$ always has at least one blow-up limit at $x_0$.
We shall only consider the case when $x_0 \in E$ (because otherwise
$E_\infty = \emptyset$) and $E$ is coral (because otherwise the
fuzzy sets $E_k \sm E_k^\ast$ could just tend to anything).
When in addition, $E \in QSAQ(U,M,\delta,h)$ with a small enough $\delta$,
we know that $E$ is locally Ahlfors-regular, and so $E_\infty$ is
Ahlfors-regular too, so it cannot be too bad.

Our intention is to prove that under reasonable assumptions on
the $L_j$ near $x_0$, $E_\infty$ is quasiminimal and sometimes minimal, 
with sliding boundary conditions coming from limit sets $L_j^0$
that we want to define now.

Let us forget about $E$ for a moment, and restrict our attention
to the sets $L_j$. We assume that for each $0 \leq j \leq j_{max}$
such that $x_0 \in L_j$, there is a closed set $L_j^0$ such that  
$$
\hbox{ the sets $L_{j,k} = r_k^{-1} (L_j - x_0)$
converge to $L_j^0$, locally in $\R^n$.}
\leqno (24.7)
$$
This is a weak way of asking for a tangent set at $x_0$, at least
along the sequence $\{ r_k \}$. When $x_0 \notin L_j$, it is easy
to see that the $L_{j,k}$ go away to infinity, and we can always
set $L_j^0 = \emptyset$.

This condition will not be enough for us,
we want to say more uniform about the way the $L_{j,k}$ converge
to the $L_j^0$. The main point of the following definition is that
it is just what we need to apply Theorem 23.8, and its conditions
should still be easy to check in simple concrete situations. 
Otherwise it is a little too complicated, and we shall partially
address this point later.

\ms\proclaim Definition 24.8. 
Suppose the Lipschitz assumption (24.1) holds, 
and let $x_0 \in U$ and a sequence $\{ r_k \}$ that tends to $0$
(as in (24.2)) be given. 
We say that the \underbar{configuration} of $L_j$ \underbar{is flat} 
at $x_0$, \underbar{along the sequence} $\{ r_k \}$, if we can find 
closed sets $L_j^0$, $0 \leq j \leq j_{max}$,  
such that (24.7) holds when $x_0 \in L_j$, and for each radius $R > 0$, 
numbers $\eta_k \geq 1$ such that
$\lim_{k \to  +\infty} \eta_k = 1$, and bilipschitz mappings
$\xi_k : B(0,R) \to \xi_k(B(0,R)$ such that $\xi_k(0) = 0$,
$$
\rho_k^{-1} |x-y| \leq |\xi_k(x)-\xi_k(y)| \leq \rho_k |x-y|
\hbox{ for } x, y \in B(0,R),
\leqno (24.9)
$$
$$
L_{j,k} \cap B(0,\rho_k^{-1}R) \i \xi_k(L_j^0 \cap B(0,R))
\i L_{j,k} \cap B(0,\rho_kR)
\leqno (24.10)
$$
for $0 \leq j \leq j_{max}$ and $k$ large enough, and 
$$
\lim_{k \to +\infty} \xi_k(x) = x
\,\hbox{ for } x\in B(0,R).
\leqno (24.11)
$$

\ms
Probably a more reasonable notion would be that the configuration of $L_j$ 
is flat at $x_0$, meaning along any sequence $\{ r_k \}$ that tends 
to $0$. Then it would not be too hard to check that the  $L_j^0$
are cones, and don't depend on the sequence $\{ r_k \}$.
This corresponds more to the usual notion of being $C^1$ at $x_0$.

The main defect of Definition 24.8 is that it concerns the whole
configuration of the sets $L_j$, ie., both the faces that compose them
and their relative positions. For the existence of limits $L_j^0$ as
in (24.7), this is not a problem and we can check it face by face. 
That is, if for each of the faces $F$ that compose an $L_j$, we can 
find a closed set $F_0$ such that
$$
\hbox{ the $F^{k} = r_k^{-1} (F- x_0)$
converge to $F^0$, locally in $\R^n$,}
\leqno (24.12)
$$
(with $F^0 = \emptyset$ if $x_0 \notin F$), then (24.7)
holds, with a set $L_j^0$ which is just the union of the $F^0$,
where $F$  runs along the faces of $L_j$ that contain $x_0$.

To compensate this defect, we shall give in Proposition 24.35
a sufficient condition for the existence of the $\xi_k$ in
Definition 24.8, that can be checked face by face, so that we don't
have to worry about how the different faces (if they are nice)
are glued to each other.

Also notice that the flatness condition above is satisfied
trivially under the Lipschitz assumption; most of the work in this 
section will consist in dealing with the other case.

We are ready for the main statement of this section.

\ms\proclaim Theorem 24.13.
Let $E$ be a coral closed set in $U$, $x_0 \in E$, and let  
the sequence $\{ r_k \}$ tend to $0$ as in (24.2). Assume that 
(24.1) holds, that the configuration of $L_j$ is flat at $x_0$, 
along the sequence $\{ r_k \}$, and that the limit sets 
$L_j^0$ defined by (24.7) satisfy (10.7) or (19.36).
Finally assume that $E_\infty$ is a closed subset
of $\R^n$ such that (24.6) holds. \hfill\break
If $E \in QSAQ(U,M,\delta,h)$, with a constant $h \geq 0$ which is 
small enough (depending on $n$, $M$, and $\Lambda$), then
$$
E_\infty \in QSAQ(\R^n,M,+\infty,h),
\leqno (24.14)
$$
with respect to the sliding boundary conditions associated to the
$L_j^0$, $0 \leq j \leq j_{max}$, defined by (24.5). \hfill\break
If $E$ is an almost minimal set in $U$, with the sliding 
conditions coming from the $L_j$, and a gauge function $h$
such that $\lim_{r \to 0} h(r) = 0$, then 
$$\eqalign{
&E_\infty \hbox{ is locally minimal in $\R^n$, with the 
sliding boundary conditions} 
\cr& \hskip2.2cm
\hbox{ defined by the $L_j^0$, $0 \leq j \leq j_{max}$}.
}\leqno (24.15)
$$

\ms
Some comments about the definitions are in order.
Concerning (10.7) or (19.36), we will see that the $L_j^0$
have natural decompositions into faces, so the definitions make sense,
and since we expect the $L_j^0$ to be at least as nice as the $L_j$,
assuming (10.7) or (19.36) for them does not hurt more than usual.

For (24.15), we accept the three types ($A_+$, $A$, or $A'$) 
of almost minimality, and our conclusion (24.15) means that for 
each one parameter family of mappings $\varphi_t : E_\infty \to \R^n$ 
that satisfy (1.4)-(1.8) for some closed ball $B \i \R^n$,
we have the two minimality properties
$$
\H^d(E_\infty \sm \varphi_1(E_\infty)) 
\leq \H^d(\varphi_1(E_\infty) \sm E_\infty)
\leqno (24.16)
$$
as in (20.6) with $h(r) = 0$, and
$$
\H^d(\varphi_1(W_1)) \leq \H^d(W_1),
\hbox{ with } 
W_1 = \big\{ x\in E_\infty \, ; \, \varphi_1(x) \neq x \big\}
\leqno (24.17)
$$
as in (20.4) or (20.5), which coincide when $h(r) = 0$.
Notice that our conditions for the competitors (i.e, the
family $\{ \varphi_t \}$) simplify here, because (2.4) is automatic,
and (1.5)-(1.6) reduce to the fact that $\varphi_t(x) = x$
for $0 \leq t \leq 1$ when $|x|$ is large enough.

The class $QSAQ(\R^n,M,+\infty,h)$ is as in Definition 2.3, 
and it makes sense even without the Lipschitz assumption.

\ms
Let us now prove Theorem 24.13. Let $R \geq 1$ be a (large) number;
we want to apply Theorem 23.8 in a ball comparable to $B(0,R)$, and
our first task will be to check the Lipschitz assumption
for the boundary pieces $L_j^0$, with a grid that we need to construct.
Of course we shall use the grid that is provided by the Lipschitz assumption 
(24.1) for the $L_j$. In the special case of the rigid assumption, we 
don't even need to worry about this, because we can directly 
apply Theorem 10.7 or Corollary 21.15 (once we look at them in a 
small enough ball, the $L_j$ coincide with simple cones and we don't 
eve have a variable domain).

Recall that the Lipschitz assumption in Definition 2.7 comes with 
constants $\lambda > 0$, $\Lambda \geq 1$, and a $\Lambda$-bilipschitz mapping
$\psi : \lambda U \to B(0,1)$. For each large enough $k$, we
define a mapping $\psi_k : B(0,3\Lambda) \to \R^n$ by
$$
\psi_k(z) = \rho_k^{-1} \psi(\lambda (x_0 + R r_k z)) - 
\rho_k^{-1} \psi(\lambda x_0),
\leqno (24.18)
$$
where $\rho_k$ is a power of $2$ that we choose so that
$$
\lambda R r_k \leq \rho_k \leq 2 \lambda R r_k.
\leqno (24.19)
$$
Thus $\psi_k$ is defined on $B(0,3\Lambda)$ as soon as 
$B(x_0, 3\Lambda R r_k) \i U$. 
By the normalization (24.19), the $\psi_k$ are $2\Lambda$-bilipschitz, 
and since $\psi_k(0)=0$ there is a subsequence of $\{ r_k \}$ for which 
$\psi_k$ converges uniformly on $B(0,3\Lambda)$ to some 
$2\Lambda$-bilipschitz mapping 
$\psi_0 : B(0,3\Lambda) \to \psi_0(B(0,3\Lambda))$.
We replace $\{ r_k \}$ with such a subsequence; this will not alter 
the construction of a grid for the limit set.
Of course $\psi_0(B(0,2\Lambda))$ contains $B(0,1)$.
Set
$$
U_R = R \psi_0^{-1}(B(0,1))
\hbox{ and } \lambda_R = R^{-1};
\leqno (24.20)
$$
then we have a $2\Lambda$-bilipschitz mapping 
$\psi_0 : \lambda_R U_R \to B(0,1)$,
that we can use to define a grid on $U_R$ and then check the rigid 
assumption. Before we do this, record the fact that by (24.20)
and because $\psi_0$ is $2\Lambda$-Lipschitz,
$$
B(0,(2\Lambda)^{-1} R) \i U_R \i B(0,2\Lambda R).
\leqno (24.21)
$$
So we want to construct the grid. To each face $F$ of our initial 
grid $\cal G$, we associate the face $\wt F = \psi(\lambda F)$,
then for each $k$ the larger face 
$\wt F^k = \rho_k^{-1} (\wt F - \psi(\lambda x_0))$.
We took $\rho_k$ dyadic, so it is a translation of a dyadic cube,
and it contains the origin if $x_0 \in F$.

We claim that for each $F$ such that $x_0 \in F$, there is a finite 
union $\wt F^\infty$ of faces of the standard unit dyadic grid 
${\cal G}_0$, such that
$$
\wt F^k \cap B(0,2) = \wt F^\infty \cap B(0,2)
\ \hbox{ for $k$ large.}
\leqno (24.22)
$$
The simplest will be to check this with coordinates.
The face $\wt F$ is given by some equations $z_\ell = a_\ell 2^{-m}$, 
$\ell\in I_1$, and some inequalities 
$a_\ell 2^{-m} \leq z_\ell \leq (a_\ell+1) 2^{-m}$,
$\ell\in I_2$, where the $a_\ell$ are integers, $I_1$ and $I_2$ form
a partition of $\{ 1, \ldots, n \}$, and $2^{-m}$ is the scale of
our initial dyadic grid. Notice that if $\wt F$ was a 
face of a larger size, we could just adapt the argument below 
(and  work with a smaller $m$).

Denote by $b_\ell$ the $\ell$-th coordinate of $\psi(\lambda x_0)$;
then $\wt F^k = \rho_k^{-1} (\wt F - \psi(\lambda x_0))$ is given by the 
equations $x_\ell = \rho_k^{-1} [a_\ell 2^{-m}-b_\ell]$, $\ell\in I_1$,
and the inequalities 
$$
\rho_k^{-1} [a_\ell 2^{-m}-b_\ell] 
\leq x_\ell \leq \rho_k^{-1} [(a_\ell+1) 2^{-m}-b_\ell],
\leqno (24.23)
$$
$\ell\in I_2$. When $b_\ell \in 2^{-m} \Bbb Z$ (and $k$ is so large that
$\rho_k^{-1}$ is a multiple of $2^m$), the corresponding equation
or inequality has integer coefficients. We shall now check that
when $b_\ell \notin 2^{-m} \Bbb Z$, we get an inequality which is 
automatically satisfied when $|x_\ell| \leq 2$; then $\wt F^k$
coincides, in $B(0,2)$, with a union of faces of ${\cal G}_0$, 
as needed.

So let us check that we get no condition when $b_\ell \notin 2^{-m} \Bbb Z$. 
Since $x_0$ lies in $F$, we get that
$\psi(\lambda x_0) \in \wt F$. The equation $b_\ell = a_\ell 2^{-m}$
is not satisfied (because $b_\ell \notin 2^{-m} \Bbb Z$), so $\ell\in I_2$
and $a_\ell 2^{-m} \leq b_\ell \leq (a_\ell+1) 2^{-m}$. In addition, both
inequalities are strict, so there is an $\varepsilon > 0$ such that
$a_\ell 2^{-m} + \varepsilon \leq b_\ell \leq (a_\ell+1) 2^{-m} - 
\varepsilon$. Then, as soon as $k$ is so large that 
$\rho_k^{-1} \varepsilon > 2$, we get that
$\rho_k^{-1} [a_\ell 2^{-m}-b_\ell] \leq -2$ and 
$\rho_k^{-1} [(a_\ell+1) 2^{-m}-b_\ell] \geq 2$, and (24.23) 
is a tautology on $[-2,2]$, as needed for our claim (24.22).

This gives us an idea of a grid on $U_R$: simply use the cubes
$\lambda_R^{-1} \psi_0(Q \cap B(0,1)) \i U_R$,
where $Q$ runs in the (rather small) family of unit dyadic cubes that meet 
$B(0,1)$. We still need to check that our family of boundaries $L_j^0$
satisfy the Lipschitz assumption, as in Definition 2.7, and this 
means that each set $A_j = \psi_0(\lambda_R L_j^0 \cap \lambda_R U_R))$ 
is (the intersection with $B(0,1)$ of) a finite union of faces of ${\cal G}_0$.

So let $0 \leq j \leq j_{max}$ be given, and assume that 
$x_0 \in L_j$; otherwise, $L_j^0 = \emptyset$ and the result is trivial.
We claim that for $k$ large,
$$
A_j = B(0,1) \cap \Big(\bigcup_{F \in {\cal F}} \wt F^\infty \Big),
\leqno (24.24)
$$
where ${\cal F}$ is the set of faces of $L_j$ that touch $x_0$
and $\wt F^\infty$ is associated to $F$ as in (24.22).
First pick $w\in A_j$, and write $w = \psi_0(\lambda_R z)$
for some $z\in L_j^0 \cap U_R$. 
By (24.7), we can write $z = \lim_{k \to +\infty} z_k$, with 
$z_k \in L_{j,k}$. Set $x_k = x_0 + r_k z_k$; thus $x_k \in L_j$,
and since there is only a finite number of faces,
we can assume (after taking a subsequence) that all the $x_k$
lie in a same face $F$. In addition, $F$ contains $x_0$,
because otherwise the $z_k$ would go away to infinity. 
Recall from (24.20) and (24.21) that 
$\lambda_R U_R = R^{-1} U_R \i B(0,2\Lambda)$, so the
$z_k$ lie in $B(0,2\Lambda)$ for $k$ large (because $z\in 
B(0,2\Lambda)$), and hence
$$
w = \psi_0(\lambda_R z) = \lim_{k \to +\infty} \psi_0(\lambda_R z_k)
= \lim_{k \to +\infty} \psi_k(\lambda_R z_k)
\leqno (24.25)
$$
because the $\psi_k$ converge to $\psi_0$ uniformly on $B(0,3\Lambda)$.
Now
$$\eqalign{
\psi_k(\lambda_R z_k) &= \psi_k(R^{-1} z_k)
= \rho_k^{-1} \psi(\lambda (x_0 + R r_k R^{-1} z_k)) - 
\rho_k^{-1} \psi(\lambda x_0)
\cr&
= \rho_k^{-1} \psi(\lambda x_k) - 
\rho_k^{-1} \psi(\lambda x_0)
}\leqno (24.26)
$$
by (24.20) and (24.18). But $x_k \in F$, so 
$\psi(\lambda x_k) \in \wt F$ and hence (by (24.26))
$\psi_k(\lambda_R z_k) \in \rho_k^{-1}(\wt F -\psi(\lambda x_0))
= \wt F^k$. By (24.22) (and because $\lambda_R z_k \in B(0,2)$ for 
$k$ large; see above (24.25)), $\psi_k(\lambda_R z_k) \in \wt F^\infty$.
But $\wt F^\infty$ is closed (and no longer depends on $k$), 
so (24.25) implies that $w\in \wt F^\infty$, 
as needed for the first inclusion.

For the other inclusion, let $w$ lie in $B(0,1) \cap \wt F^\infty$
for some face $F$ of $L_j$ such that $x_0 \in F$.
By (24.22), $w\in \wt F^k$ for $k$ large, so we can write 
$w = \rho_k^{-1} (\psi(\lambda x_k) - \psi(\lambda x_0))$ 
for some $x_k\in F$. Then set $z_k = r_k^{-1} (x_k-x_0) \in L_{j,k}$
(by (24.5)); notice that
$$
|z_k| = r_k^{-1} |x_k-x_0| \leq r_k^{-1} \lambda^{-1} \Lambda 
|\psi(\lambda x_k) - \psi(\lambda x_0)|
= r_k^{-1} \lambda^{-1} \Lambda \rho_k |w|
\leq 2 \Lambda  R |w|
\leqno (24.27)
$$
by (24.19), so the $z_k$ lie in $B(0,2\Lambda R)$, and there is a 
subsequence for which they converge to a limit $z$.
By (24.26) (or rather its proof), $\psi_k(\lambda_R z_k) = w$.
By the uniform convergence of the $\psi_k$ on $B(0,3\Lambda)$,
$w = \lim_{k \to +\infty} \psi_0(\lambda_R z_k) = \psi_0(\lambda_R z)$.
But $z\in L_i^0$ because $z_k \in L_{j,k}$ and by (24.7), and 
$z\in U_R$ by (24.20) and because $\psi_0(\lambda_R z) = w \in B(0,1)$. 
Thus $w \in A_j$, and the converse inclusion holds.

This completes our proof of (24.24); we now know that each $A_j$
is a union of faces of ${\cal G}_0$, hence the $L_j^0$ satisfy the
Lipschitz assumption in the domain $U_R$.

\ms
Now we are ready to apply Theorem 23.8 in the domain $U_R$ and 
with the sequence $\{ E_k \}$. We apply Definition 24.8, with the
radius $3\Lambda R$; this gives, for $k$ large, a 
$\rho_k$-bilipschitz mapping $\xi_k$ defined on $B(0,3\Lambda R)$.
We are interested in the restriction of $\xi_k$ to $U_R$
(recall from (24.21) that $U_R \i B(0,3\Lambda R)$), and the domain
$U_{R,k} = \xi_k(U_R)$. We want to apply Theorem 23.8 to the
sets $E_k \cap U_{R,k}$ and the boundaries $L_{j,k} \cap U_{R,k}$,
so let us check the assumptions.

We checked (23.1) (the fact that the $L_j^0$ satisfy the Lipschitz
assumption on $U_R$), and $U_{R,k} = \xi_k(U_R)$ by definition.
For (23.2), we also need to check that 
$$
L_{j,k} \cap U_{R,k} = \xi_k(L_j^0 \cap U_R).
\leqno (24.28)
$$
First let $x\in L_{j,k} \cap U_{R,k}$ be given.
Since $U_{R,k} = \xi_k(U_R) \i B(0,2\rho_k \Lambda R)$
by (24.9) and because $\xi_k(0)=0$, (24.10) 
(with $R$ replaced with $3\Lambda R$) says that
we can write $x = \xi_k(y)$ for some 
$y \in L^0_j \cap B(0,3 \Lambda R)$. But then
$y = \xi_k^{-1}(x) \in U_R$ because $x\in U_{R,k}$
and hence $x\in \xi_k(L_j^0 \cap U_R)$. Conversely, if 
$x \in \xi_k(L_j^0 \cap U_R)$ and $y \in L_j^0 \cap U_R$
is such that $\xi_k(y)=x$, (24.10) says that $x\in L_{j,k}$, 
and obviously $x \in U_{R,k}$ because $y\in U_R$. So (24.28) holds.

The bilipschitz condition (23.3) comes from (24.9), and (23.4)
follows from (24.11). Also, we assumed that the $L_j^0$ satisfy the 
unpleasant additional condition (10.7) or (19.36), so their restriction 
to $U_R$ does too.

We start under the first assumption that $E \in QSAQ(U,M,\delta,h)$,
and by (24.4) we get that $E_k \in QSAQ(U_k,M,r_k^{-1}\delta,h)$.
For $k$ large, $U_{R,k} \i U_k = r_k^{-1}(U-x_0)$, and 
$r_k^{-1}\delta$ becomes much larger than $4 \Lambda R$ and the diameter of 
$U_{R,k}$, so when we restrict to $U_{R,k} \i U_k$, we get that
$E_k \in QSAQ(U_k,M,+\infty,h)$. That is, (23.5) holds with 
$\delta = +\infty$. Since (23.6) (the limit in $U_R$)
follows from (24.6), and if $h$ is small enough, Theorem 23.8
says that $E_\infty \cap U_{R} \in QSAQ(U_R,M,+\infty,h)$.

Now this holds for every $R > 0$, and since (by (24.21)) 
$U_R$ tends to the whole $\R^n$ (when $R \to +\infty$), we get that 
$E_\infty \in QSAQ(\R^n,M,+\infty,h)$, as promised in (24.14).

Now suppose that $E$ is $A_+$-almost minimal, with a gauge function $h$ 
that tends to $0$. Then for each $M > 1$, we can find
$\delta > 0$ such that $E \in QSAQ(U,M,\delta,0)$
(just compare with Definition 20.2); by our first case, we get that
$E_\infty \in QSAQ(\R^n,M,+\infty,0)$, and since this is true for 
each $M > 1$, we even get that $E$ is locally minimal in $\R^n$, as needed.

If $E$ is $A$-almost minimal, still with a gauge function $h$ 
that tends to $0$, Definition 20.2 says that for each small number 
$h' > 0$, we can find $\delta > 0$ such that $E \in QSAQ(U,1,\delta,h')$.
Then by our first case,  $E_\infty \in QSAQ(\R^n,1,+\infty,h')$, and
again $E$ is locally minimal in $\R^n$.

If $E$ is $A'$-almost minimal, we can still conclude as above, 
except that we now use the generalization of Theorem 21.3 that goes 
like Theorem 23.8 (but with almost minimal sets), as explained
in Remark 23.28. Observe that we can always apply this statement
with a gauge function $\wt h$ which is larger than $h$, continuous 
from the right, and still tending to $0$. Or we could use 
Remark 21.6.

This concludes our proof of Theorem 24.13.
\qed

\ms
Let us now give a slightly simpler sufficient condition for the
flatness of the configuration of the $L_i$, which depends only on
the regularity at $x_0$ of the faces of the $L_j$ (and not,
apparently, on the way they are arranged in space). As we shall see, 
the construction of the bilipschitz mapping $\xi_k$ will be simpler
than expected, because the bilipschitz property will come from the
fact that the differential stays close to the identity.

We keep the notation of the beginning of this section; that is,
$U$ is an open set, $x_0 \in U$, and the boundary pieces $L_j$
satisfy the Lipschitz assumption (as in (24.1)). Denote by
$\cal G$ the associated grid on $U$, and by $\cal F$ the
set of faces of $\cal G$ that contain $x_0$ and are contained in some
set $L_j$. Also call ${\cal F}_\ell$ the set of faces $F \in {\cal F}$
that are $\ell$-dimensional.

Let a sequence $\{ r_k \}$, that tends to $0$, be given too.
For $F\in {\cal F}$ and $k \geq 0$, set $F^k = r_k^{-1} (F-x_0)$.

\proclaim Definition 24.29.
We say that \underbar{the faces} of the $L_j$ \underbar{are flat}
at $x_0$, \underbar{along the sequence} $\{ r_k \}$, 
when for each $1 \leq \ell \leq n$ and each face $F\in {\cal F}_\ell$, 
there is a set $F^0$ that contains $0$, such that
$$\eqalign{
&F^0 \hbox{ is a closed convex $\ell$-dimensional polyhedron}
\cr&\hskip2.5cm
\hbox{ in some $\ell$-dimensional vector space $V_F$,}
}\leqno (24.30)
$$
and with the following connection with the $F^k$.
For $1 \leq R < +\infty$, there is a sequence 
$\{\varepsilon_k \}$ such that
$$
\lim_{k \to +\infty} \varepsilon_k = 0,
\leqno (24.31)
$$
and for each large enough $k$, a Lipschitz mapping 
$\psi_{F,k} : F^0 \cap B(0,R) \to F^k$, such that
$$
\psi_{F,k}(0) = 0,
\leqno (24.32)
$$
$$
|D\psi_{F,k} - I| \leq \varepsilon_k
\ \ \hbox{ $\H^\ell$-almost everywhere on } F^0 \cap B(0,R),
\leqno (24.33)
$$
and
$$
\psi_{F,k}(F^0 \cap B(0,R)) \supset F^k \cap B(0,(1-\varepsilon_k)R).
\leqno (24.34)
$$

\ms
Let us comment on the slightly strange aspects of this definition.
We shall see soon (in (24.37)) that $F^0$ is the limit of the $F^k$, 
and this is also why we require that $0 \in F^0$ (recall that $x_0 \in F$ 
for $F \in \cal F$, hence $0 \in F^k$).
We could also check that when $F$ is of dimension $\ell > 0$, 
our polyhedron $F^0$ is unbounded and has a nonempty interior 
in $V_F$. It can even fill the whole space $V_F$.
We decided not to let the $F^0$ depend on $R$ 
(this would have at least complicated the proof, maybe with
no true additional generality), even though the relation with $F^k$ 
is only stated in each ball $B(0,R)$.
Similarly, requiring that $F^0$ is a convex polyhedron will simplify
our life, and will probably not hurt in applications.
Finally, we shall not try to see whether (24.34) could, or could not, 
be deduced from the other assumptions.

\ms\proclaim Proposition 24.35.
Let $U$, $x_0$, $\{ r_k \}$, and the $L_j$ be as above.
If the the faces of the $L_j$ are flat at $x_0$
along the sequence $\{ r_k \}$, then 
the configuration of $L_j$ is flat at $x_0$
along the sequence $\{ r_k \}$.

\ms
Before we start the construction of mappings 
$\xi_k : B(0,R) \to \R^n$, as in Definition 24.8, let us
use the $\psi_{F,k}$ to control some of the geometry of the
$F^0$. Let $\ell \geq 1$ and $F\in {\cal F}_\ell$ be given.

First observe that if $\psi_{F,k}$ is as in Definition 24.29,
then for $x, y \in F^0 \cap B(0,R)$,
$$
|\psi_{F,k}(x)-\psi_{F,k}(y)-x+y| \leq \varepsilon_k |x-y|.
\leqno (24.36)
$$
Indeed, for almost all choices of $x$ and $y$ 
(for the $2\ell$-dimensional product of Lebesgue measures), 
we can compute $\psi_{F,k}(x)-\psi_{F,k}(y)$ as the integral of
$D\psi_{F,k}$ on the segment $[x,y]$; this comes from the
fact that $\psi_{F,k}$ is Lipschitz, hence absolutely integrable on
almost every line. Notice that $[x,y] \i F^0$, because $F^0$ is 
convex by (24.30).
Thus (24.36) holds for almost all choices of $x$ and $y$; the general
case follows because $\psi_{F,k}$ is continuous.

Next we check that 
$$
F^0 = \lim_{k \to +\infty} F^k \ \hbox{ in } \R^n,
\leqno (24.37)
$$
(with the same definition as in (10.4)-(10.6)).
Pick $r > 0$, and let us first check that 
$d_k = \sup\big\{ \dist(x,F^k) \, ; \, x\in F^0 \cap B(0,r) \big\}$
tends to $0$. Let the $\psi_{F,k}$ be, for $k$ large, as in Definition 24.29,
with $R=2r$, and simply observe that for $x\in F^0 \cap B(0,r)$,
$\dist(x,F^k) \leq |\psi_{F,k}(x)-x| \leq \varepsilon_k |x|$
by (24.33) and (24.32). Thus $d_k \leq \varepsilon_k r$
for $k$ large.

Then we control 
$d'_k = \sup\big\{ \dist(x,F^0) \, ; \, x\in F^k \cap B(0,r) \big\}$.
We keep the same choice of $R=2r$ and $\psi_{F,k}$, and observe that
for $x\in F^k \cap B(0,r)$, (24.34) allows us to write
$x=\psi_{F,k}(z)$ for some $z\in F^0 \cap B(0,R)$. Then
$\dist(x,F^0) \leq |x-z| = |\psi_{F,k}(z)-z| \leq \varepsilon_k |z|
\leq \varepsilon_k R$, so $d'_k \leq \varepsilon_k R$ for $k$ large,
and (24.37) follows.

Next we check that the mappings $\psi_{F,k}$ essentially preserve
the boundaries. That is, denote by $\d F^0$ the boundary of $F^0$,
seen as a subset of the vector space $V_F$. We claim that for $k$
large,
$$
\hbox{for } x\in F^0 \cap B(0,R),
\psi_{F,k}(x) \in \d F^k
\hbox{ if and only if } x\in \d F^0.
\leqno (24.38)
$$
For this we shall use a little bit of topology.
Notice that by (24.36), $\psi_{F,k}$ is a bilipschitz mapping from
$F^0 \cap B(0,R)$ to its image. We compose with the affine mapping
$\rho_k : x \to x_0 +r_k x$, and get an image
$\rho_k \circ \psi_{F,k}(F^0 \cap B(0,R)) \i F$, which is contained in 
$U$ for $k$ large. We further compose with the usual bilipschitz
mapping $z \to \psi(\lambda z)$, and get a bilipschitz mapping
$h_k : F^0 \cap B(0,R) \to h_k(F^0 \cap B(0,R)) \i \wt F$,
where now $\wt F$ is a (straight) dyadic cube of dimension $\ell$.

Let $x \in F^0 \cap B(0,R)$ be an interior point of $F^0$,
suppose that $\psi_{F,k}(x) \in \d F^l$, and derive a contradiction.
Recall that $\d F$ was in fact defined as the bilipschitz image
of $\d \wt F$, so we are assuming that $h_k(x) \in \d \wt F$.
Let $S$ be the unit sphere in $V_F$, and let us see what happens
to the mapping $f: S \to V_F$, defined by
$f(\xi) = x + t \xi$, where we choose $t > 0$ so small
that $f(S) \i F^0 \cap B(0,R)$. This map cannot be homotoped
to a constant, through mappings from $S$ to $V_F \sm \{ x \}$,
yet, for $t$ small we shall use $h_k$ to find such a homotopy.
Let $L$ denote the bilipschitz constant for $h_k$; then
$h_k \circ f(S) \i A$, where
$A = \big\{ z\in \wt F \, ; \, L^{-1}t \leq |z-h_k(x)| \leq Lt \big\}$.
Now, because $h_k(x) \in \d \wt F$ and $\wt F$ is a dyadic cube,
and if $t$ is small enough, $A$ can be contracted (to a point)
inside the slightly larger annular region
$A_C = \big\{ z\in \wt F \, ; \, (CL)^{-1}t \leq |z-h_k(x)| \leq CLt \big\}$.
That is, there is a continuous function 
$\varphi : A \times [0,1] \to A_C$ such that $\varphi(z,0) = z$
and $\varphi(z,1) = c$ for $z\in A$, and where $c \in A_C$ is a 
constant. The desired deformation is the mapping
$\wt \varphi : S \times [0,1] \to V_F$ defined by
$\wt \varphi(\xi,t) = h_k^{-1}(\varphi_t(h_k\circ f(\xi),t))$;
it is easy to see that for $t$ small, it is defined (because $A_C \i 
h_k(F^0 \cap B(0,R))$) and continuous, that it avoids the value $x$, 
and that $\wt \varphi(\xi,1) = h_k^{-1}(c)$ is constant.
This contradiction shows that $\psi_{F,k}(x)$ is an interior
point of $F^l$.

The same argument, applied with the mapping $h_k^{-1}$, shows that
if $x\in F^0 \cap B(0,R)$ and $z= \psi_{F,k}(x)$ is an interior point of 
$F^l$ (which means that $h_k(x)$ is an interior point of $\wt F$),
then $x$ lies in the interior of $F^0$. This time, we use the fact 
that since $F^0$ is a convex polyhedron, there is a constant $C$ such 
that, for $x\in \d F^0$ and $t$ small enough,
$\big\{ w\in F^0 \, ; \, L^{-1}t \leq |w-x| \leq Lt \big\}$
can be contracted inside 
$\big\{ w\in F^0 \, ; \, C^{-1}L^{-1}t \leq |w-x| \leq CLt \big\}$.
This proves (24.38).

It will also be good to know that $\d F^0$ can also be seen as 
the combinatorial boundary of $F^0$, i.e., that
$$
\d F^0 = \bigcup_{H \in {\cal F}(F)} H^0,
\leqno (24.39)
$$
where ${\cal F}(F)$ denotes the set of strict subfaces of
$F$ that meet $x_0$. Indeed, let $x\in \d F^0$ be given,
choose $R > |x|$, and notice that by (24.38),
$\psi_{F,k}(x) \in \d F^k$ for $k$ large. Then there is a strict subface
$H$ of $F$ such that $\psi_{F,k}(x) \in H^k$ for infinitely many $k$.
But $x = \lim_{k \to +\infty} \psi_{F,k}(x)$ by (24.36) with $y=0$,
so by (24.37) $x$ lies in $H^0$. Conversely, let $x\in H^0$
for some $H \in {\cal F}(F)$, and choose $\psi_{F,k}$ as above.
By (24.37), $x$ is the limit of some sequence $\{ x_k \}$,
with $x_k \in H^k$. By (24.34) and since $H^k \i F^k$, we can write 
$x_k = \psi_{F,k}(z_k)$ for some $z_k \in F^0 \cap B(0,R)$,
and by (24.36) $|z_k-x_k|$ tends to $0$. Also, 
$z_k \in \d F^0$ by (24.38), so 
$x = \lim_{k\to +\infty} x_k = \lim_{k\to +\infty} z_k$
lies in $\d F^0$ too. This completes the proof of (24.39).

Finally, we shall need the following consequence of (3.8), to control the 
geometry of the faces $F^0$. Let $F$ and $G$ be two faces of
$\cal F$, and suppose that $F$ is neither reduced to a point 
nor contained in $G$; we claim that
$$
\dist(y,\d F^0) \leq \Lambda^2 \dist(y,G^0) 
\ \hbox{ for } y\in F^0. 
\leqno (24.40)
$$
Let $y\in F^0$ be given, and use (24.37) to
choose points $y_k \in F^k$
such that $\lim_{k \to +\infty} y_k = y$. Set
$x_k = x_0+r_k y_k \in F$, $\wt x_k = \psi(\lambda x_k)$, and as usual
$\wt F = \psi(\lambda F)$ and $\wt G = \psi(\lambda G)$. Then
$$\eqalign{
\dist(y,G^0) 
&=\lim_{k \to +\infty} \dist(y,G^k)
=\lim_{k \to +\infty} \dist(y_k,G^k)
= \lim_{k \to +\infty} r_k^{-1} \dist(x_k,G)
\cr&
\geq \Lambda^{-1} \lambda^{-1}\limsup_{k \to +\infty} r_k^{-1} 
\dist(\wt x_k, \wt G)
\cr&
\geq \Lambda^{-1} \lambda^{-1}\limsup_{k \to +\infty} 
r_k^{-1}\dist(\wt x_k, \d(\wt F))
\cr&
\geq \Lambda^{-2} \limsup_{k \to +\infty} r_k^{-1} \dist(x_k, \d F)
= \Lambda^{-2} \limsup_{k \to +\infty} \dist(y_k, \d F^k)
}\leqno (24.41)
$$
by (24.37), various definitions, and (3.8).
Use (24.41) to choose a sequence of points $z_k \in \d F^k$, such that 
$$
\limsup_{k \to +\infty} \dist(y_k, z_k)
\leq \Lambda^{2} \dist(y,G^0).
\leqno (24.42)
$$
Now $\d F$ is the union of a certain number of strict
subfaces $H$, and each $\d F^k$ is the union of the corresponding
subfaces. So we can choose a subsequence for which all the $z_k$
lie in $H^k$ for the same $H$.

Of course $\dist(y,G^0) < +\infty$, and by (24.42)
$\{ z_k \}$ is bounded. So we can choose a new subsequence so
that $\{ z_k \}$ tends to a limit $z_\infty$.
Since $x_0+r_k z_k \in H$,  $\{ z_k \}$ is bounded, and $r_k$ 
tends to $0$, we get that $x_0 \in H$, hence $H\in \cal F$.
Also, $z_\infty \in H^0$, since $z_k \in H^k$ and by 
(24.37) for $H$. Hence $z_\infty \in \d F^0$, 
by the representation formula (24.39). Finally,  (24.42)
implies that $\dist(y, z_\infty) \leq \Lambda^{2} \dist(y,G^0)$,
and (24.40) follows.

\ms
We are now ready to start the construction of our mappings $\xi_k$.
Let $R > 0$ be given; 
we apply our flatness assumption to every face $F \in \cal F$, 
in the larger ball $B(0,4^{n+1} R)$, to get mappings 
$\psi_{F,k} : F^0 \cap B(0,4^{n+1} R) \to F^k$ with the
properties (24.31)-(24.34).

We shall obtain $\xi_k$ after building successive extensions, 
defined on the following collection of skeletons. Set
${\cal S}_{0} = \{ 0 \}$ and for $1 \leq \ell \leq n$,
$$
B_\ell = B(0,4^{n-\ell+1}R), \ 
{\cal S}_{\ell} = \bigcup_{F \in {\cal F}_\ell} F^0 \cap B_\ell \, ,
\ \hbox{ and } 
{\cal S}_{\ell}^+ 
= \bigcup_{m \leq \ell} {\cal S}_{m}.
\leqno (24.43)
$$
We gave a special definition to ${\cal S}_{0}$ just because, 
even when $x_0$ is not a point of the grid,
we find it more pleasant to take ${\cal S}_{0} = \{ 0 \}$.

Set $\xi_k^0(0) = 0$ to start the process.
We want to define successive extensions $\xi_k^\ell$ of
$\xi_k^0$, with the following properties. First,
$\xi_k^\ell$ is defined on ${\cal S}_{\ell}^+$, and is an
extension of $\xi_k^{\ell-1}$ if $\ell \geq 1$. 
Next, for each $F \in {\cal F}_\ell$ and $k$ large enough,
$$
F^k \cap {1 \over 2} B_\ell \i \xi_k^\ell(F^0 \cap B_\ell) 
\i F^k \cap 2 B_\ell,
\leqno (24.44)
$$
and $\xi_k^\ell$ is locally Lipschitz on $F^0 \cap B_\ell$, with 
$$
|D\xi_k^\ell - I| \leq C_\ell \varepsilon_k
\ \ \hbox{ $\H^\ell$-almost everywhere on } F^0 \cap B_\ell,
\leqno (24.45)
$$
where $I$ denotes the identity on $\R^n$ and 
$C_\ell \geq 1$ is a geometric constant that will be
computed by induction. Finally, we require that
$$
|\xi_k^\ell(x)-\xi_k^\ell(y)-x + y| 
\leq 2 (1+2\Lambda^2)^2 C_\ell \varepsilon_k |x-y|
\ \hbox{ for } x, y \in {\cal S}_{\ell}^+.
\leqno (24.46)
$$

\ms
We already have $\xi_k^0$, and let us construct $\xi_k^1$
to warm up. We choose it so that
$$
\xi_k^1 = \psi_{F,k} \ \hbox{ on } F^0 \cap B_1
\leqno (24.47)
$$
for every $F \in {\cal F}_1$.
This definition is coherent: there is no conflict with the
fact that $\xi_k^0(0)=0$, by (24.32), and similarly if 
$F$ and $G$ are different faces of ${\cal F}_1$, then
$F^0 \cap G^0 = \{ 0 \}$ (for instance by (24.40) and because 
$\d F^0$ and $\d G^0$ are contained in $\{ 0 \}$ by (24.39)), 
and $\psi_{F,k}(0) = \psi_{G,k}(0) = 0$.
Notice that it could happen that ${\cal S}_{1}$ is
empty because ${\cal F}_1$ is empty, but this does not disturb.
Next, (24.44) (i.e., the fact that
$F^k \cap {1 \over 2} B_1 \i \psi_{F,k}(F^0 \cap B_1) 
\i F^k \cap 2 B_1$ holds by definition of $\psi_{F,k}$, 
and in particular (24.34) (for the surjectivity) and 
(24.36) with $y=0$ (for the management of radii).
Also, (24.45) holds with $C_1 = 1$, by (24.33). 
We don't need to check (24.46) for $\ell = 1$,
because we shall do it in a more general case now.

\ms
Now we check that (24.46) (for some $\ell \geq 1$) 
follow from (24.45) for $\ell$ and (24.46) for $\ell-1$
(notice that (24.46) is obvious for $\ell = 0$).
First we claim that because of (24.45) for $\ell$,
we have that
$$
|\xi_k^l(x) - \xi_k^l(y) - x + y| \leq C_\ell \varepsilon_k |x-y|
\ \hbox{ for } x, y \in F^0\cap B_\ell
\leqno (24.48)
$$
for each face $F\in {\cal F}_\ell$. The verification is the same
as for (24.36): we first check this for $x, y$ in a dense subset,
using the convexity of $F^0\cap B_\ell$ and the 
absolute continuity of $\xi_k^l$ on almost all lines, and the
general case follows because $\xi_k^l$ is also continuous.

Then we check (24.46); we intend to use (24.40) to control
the position of the different faces.
Let $x, y \in {\cal S}_{\ell}^+$ be given.
By the definition (24.43), we can find 
$m_x, m_y \leq \ell$ so that $x\in {\cal S}_{m_x}$ and 
$y\in {\cal S}_{m_y}$; choose $m_x$ and $m_y$ are as small 
as possible, i.e., consider the first occurrence. Then 
use (24.43) again to choose
$F \in {\cal F}_{m_x}$ and $G \in {\cal F}_{m_y}$ such that
$x \in F^0 \cap B_{m_x}$ and $y\in G^0 \cap B_{m_y}$.

First assume that $m_x = \ell$. If $x\in G^0$, then in fact it lies
in $G^0 \cap B_{m_y}$ (because $m_y \leq m_x$), and
(24.46) follows from (24.48) for $G$. So we may assume that 
$x \notin G^0$, and then $F$ is neither reduced to a point 
(it would be $x_0$, and then $x=0 \in G^0$)
nor contained in $G$ (because if $F \i G$, then (24.37) shows
that $F^0 \i G^0$), so we may apply (24.40). We get that
$$
\dist(x,\d F^0) \leq \Lambda^2 \dist(x,G^0) \leq \Lambda^2 |x-y| 
< +\infty.
\leqno (24.49)
$$
In particular, $\d F^0$ is not empty, and we can pick $x' \in \d F^0$ 
such that $|x'-x| = \dist(x,\d F^0) \leq \Lambda^2 |x-y|$.
By the representation formula (24.39) for $\d F^0$,
and the fact that each $H^0$ is convex and contains the origin,
we see that $tx' \in \d F^0$ for $0 \leq t \leq 1$; we use this 
to see that $|x'| \leq |x|$ (otherwise, some $tx'$ is closer
to $x$). Hence $x'\in B_{\ell}$. By (24.39) again, 
$x'\in H^0$ for some $H \in \cal F$ of dimension $m \leq \ell -1$, 
so
$$
x' \in H^0 \cap B_\ell \i {\cal S}_{m}
\i {\cal S}_{\ell -1}^+.
\leqno (24.50)
$$
Since $x'\in \d F^0 \cap B_\ell \i F^0 \cap B_\ell$, we can apply (24.48)
to get that
$$
|\xi_k^\ell(x')-\xi_k^\ell(x)-x' + x| 
\leq C_\ell \varepsilon_k |x'-x| 
\leq C_\ell \Lambda^2 \varepsilon_k |x-y|.
\leqno (24.51)
$$
We will continue with the proof in a moment, but let us record some 
cases first. 

If $m_x < \ell$, we simply keep 
$x' = x \in {\cal S}_{\ell -1}^+$, and some estimates will be simpler.

If $m_y = \ell$, we can assume (as above) that $y \notin F^0$,
and then we choose $y' \in {\cal S}_{\ell -1}^+$ as we did 
for $x$; if $m_y < \ell$, we just keep $y' = y$. Now
$$\eqalign{
|\xi_k^\ell(x')-\xi_k^\ell(y')-x' + y'| 
&= |\xi_k^{\ell-1}(x')-\xi_k^{\ell-1}(y')-x' + y'|
\cr&
\leq 2 (1+2\Lambda^2) C_{\ell-1} \varepsilon_k |x'-y'|
\cr&
\leq 2 (1+2\Lambda^2)^2 C_{\ell-1} \varepsilon_k |x-y|
}\leqno (24.52)
$$
because $x', y' \in {\cal S}_{\ell -1}^+$, by
(24.46) for $\ell -1$, and because 
$|x'-y'| \leq |x-y| + |x'-x| + |y'-y| \leq (1+2\Lambda^2) |x-y|$.
If $x' \neq x$ or $y' \neq y$, we add (24.51) or its analogue for $y$,
and get that
$$\eqalign{
|\xi_k^\ell(x)-\xi_k^\ell(y)-x + y| 
&\leq 2 (1+2\Lambda^2)^2 C_{\ell-1} \varepsilon_k |x-y|
+ 2 C_\ell \Lambda^2 \varepsilon_k |x-y|
}\leqno (24.53)
$$
which implies (24.46) if we make sure to choose
$C_\ell \geq 2(1+2\Lambda^2) C_{\ell-1}$.

\ms
Next we define $\xi_k^{\ell+1}$ when $1 \leq \ell < n$, assuming
that we already have $\xi_k^{\ell}$.
We shall construct our extension $\xi_{k}^{l+1}$ independently 
on each $F^0 \cap B_{\ell+1}$, $F \in {\cal F}_{\ell+1}$, 
and of course we shall make sure to keep 
the same values as $\xi_k^{\ell}$ on 
$F^0 \cap B_{\ell+1} \cap {\cal S}_{\ell}^+$.
We claim that if we proceed this way there will be no conflict 
of definition between faces. More precisely, if   
$F, G \in {\cal F}_{\ell+1}$ are different faces, we claim that
$F^0 \cap G^0 \cap B_{\ell+1} \i {\cal S}_{\ell}^+$, 
so $\xi_{k}^{l} = \xi_{k}^{l+1}$ was in fact already defined on 
the intersection. 

To prove the claim, let
$y\in F^0 \cap G^0 \cap B_{\ell+1}$ be given.
Since $F$ is neither reduced to a point nor contained
in $G$, (24.40) says that $\dist(y, \d F^0) \leq \Lambda^2 
\dist(y,G^0) = 0$. Recall from (24.39) that $\d F^0$ 
is the finite union of the closed sets $H^0$, where $H$ is a strict
subface of $F$ that contains $x_0$. Then $y$ lies in such an $H^0$,
and by the definition (24.43), $y\in {\cal S}_{\ell}^+$;
our claim follows.

So we now fix $F \in {\cal F}_{l+1}$ and proceed to
define $\xi_{k}^{l+1}$ on $F^0 \cap B_{\ell+1}$.
By (24.39), $\d F^0 = \bigcup_{H \in {\cal F}(F)} H^0$,
and by induction assumption $\xi_k^\ell$ is defined
on $\d F^0 \cap B_\ell$, with values in $\d F^k$.
We want to extend this mapping to $F^0 \cap B_{\ell+1}$,
with values in $F^k$, and for this it will be easier to use our 
bilipschitz mapping $\psi_{F,k}$ to return to the vector space $V_F$
and work there.

By (24.48) (on the faces $H^0$, and with $y=0$), 
$$
|\xi_k^\ell(x) - x| \leq C_l \varepsilon_k |x|
\ \hbox{ for } x\in \d F^0 \cap B_{\ell} \, ;
\leqno (24.54)
$$
in particular, $\xi_k^\ell(x) \in 2 B_\ell$
and (24.34) (which we can apply in the larger ball
$B(0,4^{n+1} R)$; compare with (24.43) and recall that 
$\ell > 1$) allows us to set 
$$
h(x) = \psi_{F,k}^{-1} \circ \xi_k^\ell(x) \in F^0
\ \hbox{ for } x\in \d F^0 \cap B_{\ell}.
\leqno (24.55)
$$
In fact, (24.38) (with the same radius $4^{n+1}R$ that we used
to define $\psi_{F,k}$, and applied to $h(x)$) says that 
$h(x) \in \d F^0$. And of course, by (24.36) with $y=0$ and because
$\xi_k^\ell(x) \in 2 B_\ell \i B(0,4^{n+1} R)$, we get that 
$h(x) \in 3B_\ell$. That is,
$$
h(x) \in \d F^0 \cap 3 B_\ell
\ \hbox{ for } x\in \d F^0 \cap B_{\ell}.
\leqno (24.56)
$$

We also have good estimates on the Lipschitz
constant for $h_1 = h-I$. Indeed, for 
$x, y \in \d F^0 \cap B_\ell$,
$$
|h_1(x)-h_1(y)| 
= |\psi_{F,k}^{-1} \circ \xi_k^\ell(x) 
- \psi_{F,k}^{-1} \circ \xi_k^\ell(y) - x + y| \leq a+b,
\leqno (24.57)
$$
with 
$$
a = |\psi_{F,k}^{-1} \circ \xi_k^\ell(x) - \xi_k^\ell(x)
- \psi_{F,k}^{-1} \circ \xi_k^\ell(y) + \xi_k^\ell(y)|
\leqno (24.58)
$$
and $b = |\xi_k^\ell(x)-\xi_k^\ell(y)-x+y|$.
By (24.36) applied to $\psi_{F,k}^{-1} \circ \xi_k^\ell(x)$ and 
$\psi_{F,k}^{-1} \circ \xi_k^\ell(y)$, we get that
$$
a \leq \varepsilon_k |\psi_{F,k}^{-1} \circ \xi_k^\ell(x)
-\psi_{F,k}^{-1} \circ \xi_k^\ell(y)|
\leq 2 \varepsilon_k |\xi_k^\ell(x)-\xi_k^\ell(y)| 
\leq 4 \varepsilon_k |x-y|
\leqno (24.59)
$$
(because (24.36) and (24.46) also imply that $\psi_{F,k}$
and $\xi_k^\ell$ are $2$-bilipschitz on the domains where we work).
Also, (24.46) implies that $b \leq 2 (1+2\Lambda^2)^2 C_\ell 
\varepsilon_k |x-y|$; altogether,
$$
|h_1(x)-h_1(y)| \leq 
\big(4+2 (1+2\Lambda^2)^2 C_\ell \big) \varepsilon_k |x-y|.
\leqno (24.60)
$$
Let us take advantage that $F^0$ and the values of $h_1$
lie in the space $V_F$ to apply the Whitney extension theorem to $h_1$;
we get an extension of $h_1$ to the whole $F^0 \cap B_\ell$, 
so that
$$
h_1 : F^0 \cap B_\ell \to V_F
\hbox{ is $C'_{\ell} \varepsilon_k$-Lipschitz,} 
\leqno (24.61)
$$
with $C'_{\ell} \leq C (4+2 (1+2\Lambda^2)^2 C_\ell$.
Now set
$$
h(x) = h_1(x)+x
\hbox{ for } x \in F^0 \cap B_\ell;
\leqno (24.62)
$$
This function is an extension of the mapping in (24.55), 
it still a Lipschitz function (like $h_1$), and by (24.61)
$$
|Dh - I| \leq C'_{\ell} \varepsilon_k 
\ \hbox{ $\H^{\ell +1}$-almost everywhere on } F^0 \cap B_\ell.
\leqno (24.63)
$$
We now need some topological information, which will lead to
(24.44) for $\ell+1$. We start with a control on the restriction
of $h$ to $\d F^0$. We claim that
$$
\d F^0 \cap {1\over 3} B_\ell \i h(\d F^0 \cap B_\ell).
\leqno (24.64)
$$
Let $y\in \d F^0 \cap {1\over 3} B_\ell$, and set 
$y_k = \psi_{F,k}(y)$. Notice that 
$y_k \in F^k \cap {1\over 2} B_\ell$ by (24.36), and
$y_k \in \d F^k$ by (24.38). This means that $y_k \in H^k$
for some $\ell$-dimensional subface $H$ of $F$. Observe that then
$x_0 + y_k r_k \in H$, hence $\dist(H,x_0) \leq r_k 4^{n-\ell+1}R$.
If $k$ is large enough (depending on the finite list of faces $H$
that get close to $x_0$), this can only happen if $x_0\in H$,
i.e., if $H\in {\cal F}(F)$. By induction assumption, 
(24.44) holds for $H$, and since $y_k \in H^k \cap {1\over 2} B_\ell$, 
we get that $y_k = \xi_k^{\ell}(x)$ for some $x\in H^0 \cap B_\ell$.
Now $h(x) = \psi_{F,k}^{-1}(y_k) = y$ by (24.55), and this proves 
(24.64).

Next we use (24.64) and a connectedness argument to show that
$$
F^0 \cap  {1\over 3} B_\ell \i h(F^0 \cap B_\ell).
\leqno (24.65)
$$
Denote by $W = F^0 \sm \d F^0$ the interior of $F^0$ in the space 
$V_F$; we just need to show that $h(W\cap B_\ell)$ contains 
$W\cap {1\over 3} B_\ell$, because (24.64) takes care of 
$\d F^0 \cap {1\over 3} B_\ell$. 

Of course we may assume that $W\cap {1\over 3} B_\ell$ is not empty. 
Set $Y = W \cap {1\over 3} B_\ell \cap h(W\cap B_\ell)$; 
we want to show that $Y = W \cap {1\over 3} B_\ell$. 

First we check that if $k$ is large enough, $Y$ is not empty.
Notice that
$$
\lim_{k\to +\infty} h(z) = z \hbox{ for } z\in F^0 \cap B_\ell,
\leqno (24.66)
$$
by (24.62) and because  $\lim_{k\to +\infty} h_1(z) = 0$ by (24.61)
and because $h_1(0)=0$; thus we pick $z\in W\cap {1\over 3} B_\ell$
and observe that $h(z) \in W\cap {1\over 3} B_\ell$ for $k$ large,
hence $Y$ is not empty.

Next, $Y$ is open. This is because $h: W\cap B_\ell \to V_F$
is bilipschitz (by (24.61) and (24.62)), hence open
(for instance by degree theory, or a fixed point theorem as 
in the implicit function theorem).

Finally, $Y$ is closed in $W\cap {1\over 3} B_\ell$: 
if $\{ y_j \}$ is a sequence in $Y$, with a limit
$y \in W\cap {1\over 3} B_\ell$, then $y \in Y$ because we can
find $x_j \in W\cap B_\ell$ such that $h(x_j) = y_j$,
the $x_j$ converge to a limit $x$ because 
$h$ is bilipschitz on $F^0 \cap B_\ell$, 
$h(x) = y$ because $h$ is continuous,
$x \in B_\ell$ because all the $x_j$ lie in 
${5\over 6}B_\ell$ (again because $h$ is bilipschitz and $h(0)=0$), 
and $x\in W$ (because otherwise $x\in \d F^0$ and $h(x) = y \in \d F^0$
by (24.56), a contradiction). 

Since $W \cap {1\over 3} B_\ell$ is connected (and even convex), 
we get that $Y = W\cap {1\over 3} B_\ell$, as needed for (24.65).

Let us now deduce from (24.65) that
$$
h(F^0 \cap {1 \over 4}B_\ell) \i F^0.
\leqno (24.67)
$$
If $W \cap {1 \over 4}B_\ell = \emptyset$, i.e., if 
$F^0 \cap {1 \over 4}B_\ell \i \d F^0$,
this is a direct consequence of (24.55).
Otherwise, first observe that for $k$ large, 
$h(W \cap {1 \over 4}B_\ell)$ meets $F^0$
(pick $x\in W\cap {1 \over 4}B_\ell$, and observe that $h(x) \in W$
for $k$ large, by (24.66)). Next we claim that
$$
h(W \cap {1 \over 4}B_\ell) \hbox{ does not meet $\d F^0$.}
\leqno (24.68)
$$
Indeed, suppose that $x\in W \cap {1 \over 4}B_\ell$ is such that
$h(x) \in \d F^0$. Notice that $h(x) \in {1 \over 3}B_\ell$,
so by (24.64) we can find $z\in \d F^0 \cap B_\ell$
such that $h(z) = h(x)$. This is impossible,
because $x \in W = F^0 \sm \d F^0$ and $h$ is bilipschitz
(hence injective) on $F^0 \cap B$. So (24.68) holds.

Since $W \cap {1 \over 4}B_\ell$ is convex and 
$h(W \cap {1 \over 4}B_\ell)$ meets $F^0$, (24.68) says that
$h(W \cap {1 \over 4}B_\ell) \i W \i F^0$. But 
$h(\d F^0 \cap {1 \over 4}B_\ell) \i F^0$ by (24.55), so (24.67) holds.

\ms
We are now ready to define $\xi_k^{\ell+1}$ on $F^0$,
and check our induction assumptions. By (24.67) and the bilipschitz
property of $h$ (see (24.61) and (24.62)), 
$h(F^0 \cap {1 \over 4}B_\ell) \i F^0 \cap {1 \over 3}B_\ell$,
and we can we set 
$$
\xi_k^{\ell+1} = \psi_{F,k} \circ h
\ \hbox{ on } F^0 \cap {1 \over 4}B_\ell = F^0 \cap B_{\ell+1}.
\leqno (24.69)
$$
This gives the desired definition of $\xi_k^{\ell+1}$
(recall that we can proceed face by face).
Notice that for $x\in F^0 \cap B_{\ell+1}$,
$\xi_k^{\ell+1}(x) \in F^k$ (by definition of $\psi_{F,k}$), 
and $\xi_k^{\ell+1}(x) \in 2B_{\ell+1}$, by (24.36).
This proves the second inclusion in (24.44) (for our
$F\in {\cal F}_{\ell + 1}$).

For the first inclusion we pick $y\in F^k \cap {1 \over 2}B_{\ell+1}$
and we need to find $x\in F^0 \cap B_{\ell+1}$ such that
$\xi_k^{\ell+1}(x) = y$. First apply (24.34) to find 
$w\in F^0$ such that $\psi_{F,k}(w) = y$. By (24.36), 
$w\in {2 \over 3}B_{\ell+1} \i {1 \over 3}B_\ell$.
By (24.65), we can find $x\in F^0 \cap B_\ell$ such that
$h(x) = w$. By (24.61) and (24.62), and because 
$w\in {2 \over 3}B_{\ell+1}$, $x\in B_{\ell+1}$.
Then $\xi_k^{\ell+1}(x) = \psi_{F,k} \circ h(x) = y$
by (24.69) and our definitions, which proves 
the second inclusion in (24.44).

The estimate (24.45) on $D\xi_k^{\ell+1}$
follows from (24.69), (24.33), (24.63), and the fact
that all our mappings are Lipschitz. We just need to pick
$C_{\ell+1}$ somewhat larger that $C_\ell$.

We already checked (near (24.48)) that (24.46) follows from 
(24.45) and the induction assumption, so we completed
our induction step, and we get mappings $\xi_k^\ell$,
$1 \leq \ell \leq n$, with the properties (24.44)-(24.46).

For the verification of flatness, we need to construct a $\xi_k$,
as in Definition 24.8, which is defined on the whole $B(0,R)$,
so we need to extend $\xi_k^n$ a last time.
We proceed as before, set $g_k = \xi_k^n -I$ on
${\cal S}_n^+$, observe that by (24.46) $g_k$ is 
$2(1+2\Lambda^2)^2 C_n \varepsilon_k$-Lipschitz
on ${\cal S}_n^+$, use the Whitney extension theorem (on $\R^n$)
to extend $g_k$, and set $\xi_k = g_k + I$ on $B(0,R)$
(we shall not need the values further away).
Let us now check that $\xi_k$ satisfies the properties 
required in Definition 24.8.

First, we need to find sets $L_j^0$, $0 \leq j \leq j_{max}$,
such that (24.7) holds when $x_0 \in L_j$ (otherwise, we don't need
to find $L_j^0$, or we can take the empty set). 
Fix $j$, and let us try 
$L_j^0 = \bigcup_{F \in {\cal F} ; F \i L_j} F^0$.
From (24.37) we deduce that 
$$
L_j^0 = \lim_{k \to +\infty} \bigcup_{F \in {\cal F} ; F \i L_j} F^k.
\leqno (24.70)
$$
Recall that $F^k = r_k^{-1} (F - x_0)$ (see before Definition 24.29),
and set $L'_j = \bigcup_{F \in {\cal F} ; F \i L_j} F$ and
$L'_{j,k} = r_k^{-1} (L'_j - x_0)$; then (24.70) says that
$L_j^0 = \lim_{k \to +\infty} L'_{j,k}$, while (24.7) requires
$L_j^0 = \lim_{k \to +\infty} L_{j,k}$, i.e., for the possibly 
larger set $L_j$. That is, $L_j$ is the union of all the faces $F$
that are contained in $L_j$, and not only those that contain $x_0$,
as in the definition of $\cal F$ (above Definition 24.29).
But the difference does not hurt: there is only a finite number
of faces $F$ in $L_j$ that do not contain $x_0$, and removing them
does not change the limit of the $L_{j,k}$ because the corresponding 
sets $F^k$ go away to infinity. So our $L_j^0$ satisfy (24.7).

We know that $\xi_k$ is bilipschitz, that $\xi_k(0) = \xi_k^{0}(0)=0$,
the bilipschitz condition (24.9) holds with any 
$\rho_k \geq 1+ C (1+2\Lambda^2)^2 C_n \varepsilon_k$
(because of the small Lipschitz constant for $g_k = \xi_k - I$),
and (24.11) holds because $\lim_{k \to +\infty} g_k(x) = 0$, so we are
left with (24.10) to check.

First let $y\in L_{j,k} \cap B(0,\rho_k^{-1} R)$. If $k$ is large 
enough, and by the discussion above, $y\in L'_{j,k}$, which means that
$y\in F^k$ for some $F\in \cal F$ such that $F \i L_j$.
Let $\ell$ be such that $F\in {\cal F}_\ell$; by (24.44) and because
${1\over 2} B_\ell \supset {1\over 2} B_n \supset B(0,R)$,
we can find $x\in F^0 \cap B_\ell$ such that $\xi_k^\ell(x) = y$.
Then $x\in L_j^0$, $\xi_k(x) = \xi_k^\ell(x) = y$, and
$x \in B(0,R)$ by (24.48), and if $\rho_k \geq 1+C C_\ell \varepsilon_k$.
This yields the first part of (24.10). For the second part, we take
$x\in L_j^0 \cap B(0,R)$, choose a face $F \i {\cal F}$ such
that $F \i L_j$ and $x\in F^0$, and notice that 
$\xi_k(x) = \xi_k^\ell(x) \in F^k$, by (24.44),
and that $|\xi_k(x)| \leq (1+C_\ell \varepsilon_k)|x| 
\leq (1+C_\ell \varepsilon_k) R \leq \rho_k R$
by (24.48) and if $\rho_k \geq 1 + 2 C_\ell \varepsilon_k$. 
This proves the second inclusion in 
(24.11), the $\xi_k$ satisfy the requirements for Definition 24.8,
and this concludes our proof of Proposition 24.35.
\qed

\bigskip
\centerline{PART VI : OTHER NOTIONS OF QUASIMINIMALITY}
\ms

\msi
{\bf 25. Elliptic integrands; the main lower semicontinuity result.} 
\ms

Up to now we used the Hausdorff measure $\H^d(E)$ to measure the
size of our sets $E$, but it is natural to consider other measure like
$\int_E f(x) d\H^d(x)$ (if our space is not homogeneous), and even
with functions $h$ that depend not only on the position of $x$ in 
space, but also on the tangent plane to $E$ at $x$, to model 
nonisotropic spaces. To the author's knowledge, the question of
nonisotropic integrands like $h$, in the context of the Plateau 
problem, was raised by F. Almgren.
In [A1], 
he states his generalization of Reifenberg's theorem on the 
homological Plateau problem in terms of elliptic integrands, and even 
adds, probably to explain the use of currents and varifolds in 
[A1]: 
``It does not seem possible to extend the arguments of De Giorgi 
or of Reifenberg to general elliptic integrands. 
In particular, the orthogonal invariance of the $m$ area integrand 
$F=1$ is essential for the applicability of Reifenberg's methods".
The author of these notes does not know whether this sentence is taken 
too seriously by the specialists, but just to make sure we shall 
explain in this section why many of the results of the previous 
sections still hold when $\H^d$ is integrated against a reasonable
elliptic integrand. A good part of it is based on an adaptation of
Dal Maso, Morel, and Solimini's uniform concentration lemma, which
Y. Fang's proved to make his extension of Reifenberg's existence result 
for the homological Plateau problem work also in the 
context of elliptic integrands; see [Fa].  
As usual, we shall need to change the lemmas because of the boundary
conditions, but not the general scheme of the proof.

Let us first say what sort of elliptic integrands we shall consider,
and how we integrate them on (rectifiable) sets. Our integrands will
be Borel-measurable positive functions $f : U \times G(n,d) \to (0,+\infty)$,
where $U$ is an open set in $\R^n$ and $G(n,d)$ denotes the 
Grasssman manifold of vector $d$-planes in $\R^n$, 
and their integral on rectifiable sets $E \i U$ will be defined by
$$
J_f(E) = \int_E f(x,T_xE) \, d\H^d(x),
\leqno (25.1)
$$
where $T_x E$ denotes the non oriented vector $d$-plane
which gives the approximate tangent plane to $E$ at $x$;
thus $T_x E$ is defined $\H^d$-almost everywhere on $E$
because $E$ is rectifiable.

We shall not really need to define $J_f(E)$ when $E$ is not rectifiable, 
because we shall concentrate on quasiminimal sets, but let us mention 
here that we could do so easily with a trick: we could define an
auxiliary function $\wt f : U \to (0,+\infty)$
(possibly using the values of $f$, but not necessarily), and
then set 
$$
J_{f,\wt f}(E) = J_f(E_{rec}) + \int_{E \sm E_{rec}} \wt f(x) \, d\H^d(x)
\leqno (25.2)
$$
for Borel sets $E$ with $\H^d(E) < +\infty$, and where $E_{rec}$
denotes the rectifiable part of $E$. This may sound a little 
artificial, but the issue typically shows up when we want to state a
result connected to the Plateau problem, want to define a functional
$J$ even for sets that are not rectifiable, but know anyway that the 
minimizers (or even the very good competitors) will be rectifiable.

We will work with the following class of integrands, which is the
same as in Fang's paper [Fa]; 
Almgren [A1] and [A3] 
mentions slightly more restricted classes, but the spirit is the same.

\ms\proclaim Definition 25.3. 
Let $U \in \R^n$ be open. For $0 < a \leq b < +\infty$,
we denote by ${\cal I}(U,a,b)$ (or just ${\cal I}(a,b)$)
the set of continuous
functions $f : U \times G(n,d) \to (0,+\infty)$, such that
$$
a \leq f(x,T) \leq b
\ \hbox{ for $x\in U$ and } T \in G(n,d),
\leqno (25.4)
$$
and, for each $x\in U$, there is a radius $r(x) > 0$ and function 
$\varepsilon_x : (0,r(x)] \to [0,1]$, with 
$$
\lim_{r \to 0} \varepsilon_x(r) = 0
\leqno (25.5)
$$
(that will measure the near optimality of planes near $x$),
such that 
$$
J_f(P \cap B(x,r)) \leq J_f(S\cap B(x,r)) + \varepsilon_x(r) r^d
\leqno (25.6)
$$
when $0 < r \leq r(x)$, $P$ is a $d$-plane through $x$, 
and $S \i \overline B(x,r)$ is a compact rectifiable set which cannot 
be mapped into $P \cap \d B(x,r)$ by any Lipschitz mapping
$\psi : \overline B(x,r) \to \overline B(x,r)$ 
such that $\psi(y) = y$ for $y\in P \cap \d B(x,r)$.

\ms
This definition is probably not optimal, but something like
(25.6) is needed if we want to have existence results for the
(local) minimization of $J_f$. Let us just explain what may go 
wrong, without computing a precise example. Take $n=2$, $d=1$,
parameterize $G(2,1)$ by the angle $\theta$ of a line of $G(2,1)$
with the horizontal direction, and consider functions of the form 
$f(x,\theta) = f_1(x_1) f_2(\theta)$, where $f_1(x)$ is a function
of the horizontal variable of $\R^2$ which is minimal on the 
$x_1$-axis, and $f_2(\theta)$ is an nice function, but such that
$f_2(\pm \pi/4) < f_2(0)/2$. A good minimizing sequence 
will be composed of zig-zag curves that stay close to the $x_1$-axis, 
with for instance slopes that stay close to $\pm 1$; it will converge 
to the axis itself, which is not a minimizer because $f_2(0)$ is too large. 
Then we can easily cook up some some problems for which there is no minimizer 
because, if one existed, it would have to be the $x_1$-axis.

Our definition is a little unpleasant because it is hard to control
the list of sets $S$ that satisfy the non-retractability condition
above, but there are convexity conditions that imply that
$f\in {\cal I}(a,b)$ (assuming (25.4)).
At least we have one example : the constant $1$ 
lies in ${\cal I}(1,1)$ because the orthogonal
projection of any $S$ as above contains $P\cap B(x,r)$.
Similarly, if $f$ is a continuous function of $x$ alone
such that $a \leq f(x) \leq b$ everywhere, then
$f\in {\cal I}(a,b)$.
We are a little sorry because we do not allow 
functions $f$ that are merely lower semicontinuous; 
see Remark 25.87, Claim 25.89, and (25.96)
below for slightly more general conditions that work.

The main result of this section is the following generalization
of Theorem 10.97. 

\ms\proclaim Theorem 25.7.
Let $U$, $\{ E_k \}$, and $E$ satisfy the hypotheses (10.1), 
(10.2), (10.3), and (10.4). Also suppose that $h$ is small enough, 
depending only on $n$, $M$, and $\Lambda$.
Then 
$$
J_f(E\cap V) \leq \liminf_{k \to +\infty} \, J_f(E_k \cap V)
\leqno (25.8)
$$
for $0 < a \leq b < +\infty$, every open set $V \i U$,
and every $f \in {\cal I}(V,a,b)$.

\ms
Notice that the sets $E_k$ and $E$ are rectifiable, 
by Theorem 5.16 and Proposition 10.15,
so $J_f(E_k\cap V)$ and $J_f(E\cap V)$ are well defined by (25.1).

We shall first prove this for sets $V$ that are relatively compact in $U$.
As in [Fa], 
our argument will be inspired by 
Dal Maso, Morel, and Solimini's [DMS], 
but with some pleasant simplifications in the covering argument.
The present argument will rely directly on the rectifiability of the
limit $E$ (which gives flatness almost everywhere), rather than
the concentration lemma (which gives flatness with a more quantitative control, 
even though on balls that are not centered at the original point); 
this will allow us to apply a more standard version of the 
Vitali covering lemma, at the (small) price of a less constructive 
argument.

Let $\{ E_k \}$, $E$, and $V \i\i U$ be as in the statement. Our first
task is to find a large set $E^1$ of $E\cap V$, and lots of nice small
balls centered on $E^1$. First observe that $E$ is rectifiable, by 
Proposition 10.15, hence for $\H^d$-almost every $x\in E \cap V$,
$$
\lim_{r \to 0} r^{-d} \H^d(E\cap B(x,r)) = \omega_d
\leqno (25.9)
$$
(see Theorem 17.6 on page 240 in [Ma]), and 
$$
E \hbox{ has a tangent plane $P(x)$ at $x$}
\leqno (25.10)
$$
(here Theorem 15.19 (3) on page 212 of [Ma] 
only gives an approximate tangent plane, but since $E$ is
locally Ahlfors-regular near $x$ (by (10.11)) 
Exercise 41.21 on page 277 of [D4]  
says that this plane is a true tangent plane.

Set $E^0$ denote the set of points $x\in E \cap V$
such that (25.9) and (25.10) hold. In fact, other constraints
will appear later, which will lead us to removing other negligible
pieces of $E^0$, but this will not matter.
Let $\varepsilon > 0$ be small; we shall let 
it tend to $0$ at the end of the estimate.
For each $x\in E^0$, we select $r(x) > 0$ with various
properties, which are all true for $r$ small enough. First,
$$
\hbox{$r(x)$ is small enough, depending on $r_0$, $\lambda$, 
$\Lambda$, and $\dist(x,\R^n \sm U)$,}
\leqno (25.11)
$$
where $\lambda$, $\Lambda$, and $r_0$ (the scale of the
the dyadic cube that we use in the unit ball) are as in the
Lipschitz assumption; how small will depend on simple geometric 
constraints that will arise in our construction, and we don't
need to know this precisely.
Next,
$$
\big|\omega_d - r^{-d} \H^d(E\cap B(x,r)) \big|
\leq \varepsilon
\ \hbox{ for } 0 < r < 2 r(x),
\leqno (25.12)
$$
and
$$
\dist(y,P(x)) \leq \varepsilon |y-x|
\ \hbox{ for } y \in E \cap B(x,3r(x)). 
\leqno (25.13)
$$
We also need to control the variations of our integrand
near $x$. The function $x \to T_x E$ (where with our new notation,
$T_x E$ is the vector space parallel to $P(x)$) is Borel-measurable
on $E^0$; this is unpleasant, but not hard to check, especially 
because $E$ is locally Ahlfors-regular, but anyway we leave the proof
to the reader. Then by Lusin's theorem,
we can find a Borel set $E^1 \i E^0$, such that
$$
\H^d(E \cap V \sm E^1) \leq \varepsilon,
\leqno (25.14)
$$
and on which $T_x E$ is a continuous function of $x$. 
We use this, and the uniform continuity of $f$ on
$W \times G(n,d)$ for some neighborhood $W$ of $x$, 
to say that if $r(x) > 0$ small enough (depending on $x$),
$$
|f(y,T_yE) - f(x,T_xE)| \leq \varepsilon
\ \hbox{ for } y \in E \cap B(x,2r(x))
\leqno (25.15)
$$
and 
$$
|f(y,T_xE) - f(x,T_xE)|\leq \varepsilon
\ \hbox{ for } y \in B(x,2r(x)).
\leqno (25.16)
$$
We certainly require this from $r(x)$, but other similar constraints concerning 
the set $E^1$ and the radius $r(x)$ will show up, in relation with 
our boundary constraints; we shall find it more pleasant to
mention them later, as we use them.

For the moment, we fix $x\in E^1$ and $0 < r \leq r(x)$, 
set $B = \overline B(x,r)$, and try to evaluate the contribution of 
$B$ to the two sides of (25.8). Set $D = P(x) \cap B$;
we first compare $J_f(E \cap B)$ with $J_f(D)$. We just observe
that
$$\eqalign{
J_f(E \cap B) &= J_f(E^1 \cap B) + J_f((E \sm E^1) \cap B)
\cr&
\leq \int_{E^1 \cap B} f(y,T_yE) \, d\H^d(y) + b \H^d((E \sm E^1) \cap B)
\cr&
\leq (f(x,T_xE)+\varepsilon) \, \H^d(E^1 \cap B) 
+ b \H^d((E \sm E^1) \cap B)
}\leqno (25.17)
$$
by (25.1), because $f \in {\cal I}(a,b)$, and by (25.15).
Then 
$$
\H^d(E^1 \cap B) \leq \H^d(E \cap B) \leq (\omega_d +\varepsilon) r^d 
\leqno (25.18)
$$
by (25.12), and 
$$
f(x,T_xE) \omega_d r^d 
= \int_{D}  f(x,T_xE) \, d\H^d(y)
\leq J_f(D) + \varepsilon \omega_d r^d 
= J_f(P(x) \cap B) + \varepsilon \omega_d r^d 
\leqno (25.19)
$$
by (25.16) and because $D = P(x) \cap B$. Thus
$$\eqalign{
(f(x,T_xE)+\varepsilon) &\, \H^d(E^1 \cap B) 
\leq (\omega_d^{-1} r^{-d} J_f(P(x) \cap B) + 2\varepsilon) 
(\omega_d +\varepsilon) r^d
\cr&
\leq J_f(P(x) \cap B) + r^d 
(\varepsilon \omega_d^{-1} r^{-d} J_f(P(x) \cap B)+2\varepsilon \omega_d
+2\varepsilon^2)
\cr&
\leq J_f(P(x) \cap B) + \varepsilon r^d (b+2\omega_d+\varepsilon)
}\leqno (25.20)
$$
by (25.18), (25.19), and (25.1), and now (25.17) yields
$$\eqalign{
J_f(E \cap B)
\leq J_f(P(x) \cap B) + \varepsilon (b+2\omega_d+\varepsilon) r^d 
+ b \H^d((E \sm E^1) \cap B).
}\leqno (25.21)
$$
We also need a lower bound for $J_f(E_k \cap B)$, and for
this we shall need to introduce a set $S$
as in Definition 25.3 and use the quasiminimality of $E_k$
to show that $S$ cannot be retracted.  

First observe that for $k$ large enough (depending on $x$,
but this will not matter)
$$
E_k \cap {3\over 2} B \i H,
\ \hbox{ where }
H = \big\{ y\in \R^n \, ; \, 
\dist(y,P(x)) \leq 3 \varepsilon r \big\},
\leqno (25.22)
$$
just by (25.13) and because the $E_k$ converge to $E$
in $2B$. We want to modify $E_k$ a first time, in the set
$$
A_{20} = \big\{ y\in \R^n \, ; \, 
(1-20\varepsilon) r \leq |y-x| \leq (1+20\varepsilon) r \big\}
\leqno (25.23)
$$
because we want a set $S$ such that $S\cap \d B(x,r) \i P(x)$.
Denote by $\pi$ the orthogonal projection on $P(x)$, and set
$$
g(y) = \alpha(|y-x|) \pi(y) + (1-\alpha(|y-x|)) y
\leqno (25.24)
$$
for $y\in \R^n$, where $\alpha$ is the continuous, piecewise
affine mapping defined by $\alpha(t) = 0$ for 
$t\in [0,(1-20\varepsilon) r] \cup [(1+20\varepsilon) r,+\infty)$, 
$\alpha(t) = 1$ for $t\in [(1-10\varepsilon) r,(1+10\varepsilon) r]$,
and $\alpha$ is affine on each of the two remaining intervals
$[(1-20\varepsilon) r,(1-10\varepsilon) r]$
and $[(1+10\varepsilon) r,(1+20\varepsilon) r]$.
Notice that 
$$
|g(y)-y| \leq |\pi(y)-y| \leq 3 \varepsilon r
\ \hbox{ for } y\in H,
\leqno (25.25)
$$
and hence, if we set $A_{5} = \big\{ y\in \R^n \, ; \, 
(1-5\varepsilon) r \leq |y-x| \leq (1+5\varepsilon) r \big\}$,
$$
g(H) \cap A_5 \i P(x)
\leqno (25.26)
$$
because if $y\in H$ is such that $g(y) \in A_5$,
then $||y-x|-r| \leq 8 \varepsilon r$, and hence
$g(y) = \pi(y)$. We set $S = g(E_k) \cap B$;
the next lemma is probably the key step of the proof.  

\ms\proclaim Lemma 25.27.
There is no Lipschitz mapping $\psi : \overline B \to \overline B$
such that (as in Definition~25.3) $\psi(y) = y$ for 
$y\in P(x) \cap \d B$ and $\psi(S) \i P(x) \cap \d B$.

This will allow us to apply Definition 25.3 and get (25.6).
We want to prove the lemma by contradiction, 
suppose there exists such a $\psi$,
and use it to construct a new mapping $\varphi$ and an impossible 
competitor for $E_k$. First observe that $\pi\circ\psi$ has the same
properties as $\psi$, so, at the price of replacing $\psi$ with
$\pi\circ\psi$, we may assume that
$$
\psi(z) \in P(x) \cap B \ \hbox{ for } z\in B.
\leqno (25.28)
$$
We want to extend $\psi$ to $\R^n$,
and we do this in two steps. First we set
$$
\psi(z)=z \ \hbox{ for } 
z\in (P(x) \sm B) \cup [\R^n \sm (1+5\varepsilon)B].
\leqno (25.29)
$$
The first extension that we get this way is still Lipschitz; 
we can easily check this by hand, using the fact that $\psi(z)=z$ 
for $z\in P(x) \cap \d B$
(connect a point of $\overline B$ to any point of the
rest of a domain through a point of $P(x) \cap \d B$).
Then we extend $\psi$ to the whole $\R^n$, using for
instance the Whitney extension theorem; we can even make
sure that $\psi((1+5\varepsilon)B) \i (1+5\varepsilon)B$,
because otherwise we can compose the restriction to
$(1+5\varepsilon)B$ with the radial projection from $\R^n$
onto $(1+5\varepsilon)B$. 

If we were dealing with quasiminimal sets with no boundary 
constraints, we would use the mapping $\varphi = \psi \circ g$
as the endpoint of a one parameter family (as in Definition 1.3),
to test the quasiminimality of $E_k$ and get a contradiction.
The point is that $E_k \cap B$ is sent to $S$ by $g$,
and then to $P(x) \cap \d B$ by $\psi$, which means that
all its measure disappears. We will see that $\H^d(E_k \cap A)$
is small, and so is $\H^d(\varphi(E_k \cap A))$, so we would get that
$\H^d(\varphi(E_k \cap (A\cup B)))$ is small, while
$\H^d(E_k \cap B)$ is reasonably large, because $E_k$ is 
locally Ahlfors-regular. The ensuing contradiction would 
prove Lemma 25.27.

But we have boundary constraints coming from the $L_j$, and
so we will need to modify $\varphi$ before we use it to test
the quasiminimality of $E_k$. This will be easier if we first
add a few constraints on the set $E^1$ and the radius $r(x)$.

We first replace the sets $E^0$ and $E^1$ above by the slightly
smaller sets where we add the requirement that every point $x\in E^0$,
is a Lebesgue density point of $E \cap F$ for every face $F$ 
of our dyadic grid on $U$ that contains $x$.
This means that 
$$
\lim_{r \to 0} r^{-d} \H^d(E \cap B(x,r) \sm F) = 0
\leqno (25.30)
$$
for every such face. For each face $F$, (25.30) is true for
$\H^d$-almost every $x\in E\cap F$, so our new constraint
only removes a $\H^d$-negligible set from $E^0$ and $E^1$.
Notice that for each $x\in U$, there is a smallest face 
$F_x$ of our grid that contains $x$ (because the intersection
of two faces that contain $x$ is a face that contains $x$),
and (25.30) for $F_x$ is stronger than for the other faces.

We also put additional conditions on the radius
$r(x)$, $x\in E^1$. Namely, we require that
$$
\dist(x,F) \geq 2 r(x)
\hbox{ for every face $F$ of our grid that does not contain $x$,}
\leqno (25.31)
$$
we recall our constraint (25.11) that $r(x)$ be small enough
(constraints of that type will arise soon), and we also demand that
$$
\H^d(E \cap B(x,r) \sm F_x) \leq \varepsilon^d r^d
\ \hbox{ for } 0 < r < 3r(x). 
\leqno (25.32)
$$

Because of the boundary constraints, we shall need a retraction
on the smallest face $F = F_x$ that contains $x$. Recall from 
Lemma 3.14 and Remark 3.25 that when $F$ is a standard dyadic cube 
(i.e., under the rigid assumption), there is a natural projection 
$\pi = \pi^F$, defined on a $(r_0/3)$-neighborhood of $F$ 
(where $r_0$ is the scale of our smallest cubes). 
Under the Lipschitz assumption, we use the rigid face
$\wt F = \psi(\lambda F)$, and the projection $\pi_F$ defined by
$$
\pi_F(y) = \lambda^{-1} \psi^{-1}(\pi^{\wt F}(\psi(\lambda y)),
\leqno (25.33)
$$
which is now defined on a neighborhood of $F$
whose width near $x$ could easily be computed in terms of $r_0$, 
$\lambda$, $\Lambda$, and $\dist(x,\R^n \sm U)$.
With this projection comes a retraction, defined by
$$
\Pi_F(y,s) = \lambda^{-1} \psi^{-1}
\big(s(\pi^{\wt F}(\psi(\lambda y))+(1-s)\psi(\lambda y)\big),
\leqno (25.34)
$$
for $y$ the same neighborhood of $F$ and $0 \leq s \leq 1$.
We shall use a different time for different points; that is,
for $y\in \R^n$ and $0 \leq t \leq 1/2$, we set
$$
s(y,t) = 0 \hbox{ when } |y-x| \geq (1+30\varepsilon) r,
\leqno (25.35)
$$
and
$$
s(y,t) = 2t \, 
\min\big(1,{(1+30\varepsilon) r - |y-x| \over 10\varepsilon r}\big) 
\ \hbox{ when } |y-x| \leq (1+30\varepsilon) r
\leqno (25.36)
$$
(so that in particular $s(y,t) = 2t$ when
$y\in \overline B(x,(1+20\varepsilon) r$). Then we set
$$
s(y,t) = s(y,1/2) \hbox{ when } 1/2 \leq t \leq 1.
\leqno (25.37)
$$
We also want to interpolate between $y$ and $\varphi(y) = \psi(g(y))$,
so we set
$$
\varphi_t(y) = y \ \hbox{ for } 0 \leq t \leq 1/2
\leqno (25.38)
$$
and
$$
\varphi_t(y) = (2t-1)\varphi(y) + (2-2t)y \ \hbox{ for } 1/2 \leq t \leq 1.
\leqno (25.39)
$$
Finally, we want to use the family $\{ h_t \}$ defined by
$$
h_t(y) = \Pi_F(\varphi_t(y),s(y,t)),
\leqno (25.40)
$$
but a few verifications will be needed.

First of all, $\pi_F(z)$ and $\Pi_F(z,s)$ are well defined
when $z\in 2B$ and $0 \leq s \leq 1$. Indeed,
$x\in F$ and hence $\dist(z,F) \leq |z-x| \leq 2 r$,
and (25.11) allows us to choose $r(x)$ so small that $2B$ is
contained in the neighborhood of $F$ that was mentioned below (25.33).
Next observe that
$$
\varphi_t(y) = y \hbox{ for $0 \leq t \leq 1$ when }
|y-x| \geq (1+20\varepsilon) r, 
\leqno (25.41)
$$
since $g(y) = y$ by (25.24) and because $\alpha(|y-x|) = 0$, 
hence $\varphi(y) = \psi(y) = y$ by (25.29), and finally
$\varphi_t(y) = y$ by (25.38) and (25.39).

If moreover $(1+30\varepsilon) r \leq |y-x| \leq 2 r$, 
then $s(y,t) = 0$ by (25.35) and (25.37), so 
$h_t(y) = \varphi_t(y) = y$ by (25.40) and (25.34).

If instead $y\in \R^n \sm 2B$, and even though $\Pi_F(y,s)$ is not formally 
defined above, we can safely extend the definitions (for 
instance, set $\Pi_F(z,0) = z$) and keep $h_t(y)=y$ there. So
$$
h_t(y)=y \hbox{ for $0 \leq t \leq 1$ when } 
|y-x| \geq  (1+30\varepsilon) r,
\leqno (25.42)
$$
and nothing happens there. 

Before we continue with other regions, 
it will be useful to know that for some constant $C_\Lambda$,
that depends on the local Ahlfors regularity constant for $E$
near $x$, 
$$
\dist(y,F) \leq C_\Lambda \varepsilon r
\ \hbox{ for } y \in E \cap {29 \over 10} B.
\leqno (25.43)
$$
Indeed otherwise, $E \cap B(y,C_\Lambda \varepsilon r)$
does not meet $F$, and then $\H^d(E \cap B(x,3r) \sm F)
\geq \H^d(E \cap B(y,C_\Lambda \varepsilon r) \sm F)
\geq C^{-1} (C_\Lambda \varepsilon r)^d$ by 
local Ahlfors regularity of $E$ (see Propositions~4.1 and 4.74).
This contradicts (25.32) if $C_\Lambda$ is chosen large enough; (25.43)
follows. Since the $E_k$ converge to $E$, we also get that
for $k$ large, 
$$
\dist(y,F) \leq 2C_\Lambda \varepsilon r
\ \hbox{ for } y \in E_k \cap {28 \over 10} B.
\leqno (25.44)
$$
For the moment, we only care about $y \in E_k \cap 2B$.
Set $\wt y = \psi(\lambda y)$, and notice that for 
$0 \leq s, s' \leq 1$, 
$$\leqalignno{
|\Pi_F(y,s)- \Pi_F(y,s')| 
&= \lambda^{-1} \Big[
\psi^{-1} \big(s(\pi^{\wt F}(\wt y)+(1-s)\wt y\big)
- \psi^{-1} \big(s'(\pi^{\wt F}(\wt y)+(1-s')\wt y\big)
\Big]
\cr&
\leq \lambda^{-1} \Lambda |s-s'| |\pi^{\wt F}(\wt y)-\wt y|
\leq C \lambda^{-1} \Lambda |s-s'| \dist(\wt y, \wt F)
\cr&
\leq C \Lambda^2 |s-s'| \dist(y,F)
\leq C(\Lambda) |s-s'| \varepsilon r
& (25.45)
}
$$
by (25.34), because $\pi^{\wt F}$ is Lipschitz
and $\pi^{\wt F}(z) = z$ on $\wt F$ (see (3.6)),
and by (25.43); here and below, $C(\Lambda)$ is our notation
for a constant that depends on $\Lambda$ (but also on $M$ and 
the other usual constants).
Since $\Pi_F(y,0) = y$, we get that
$$
|\Pi_F(y,s)- y| \leq C(\Lambda) \varepsilon r
\ \hbox{ for $y \in E_k \cap {18 \over 10} B$ and $0 \leq s \leq 1$.}
\leqno (25.46)
$$

Notice also that for $y, z \in 2B$ and $0 \leq s \leq 1$, 
$$
|\Pi_F(y,s)- \Pi_F(z,s)| \leq C \Lambda^2 |y-z|,
\leqno (25.47)
$$
by (25.34) and because $\pi^{\wt F}$ is $C$-Lipschitz.

We continue our study of the $h_t$ with what happens in the region
$$
R_1 = \big\{ y\in E_k \, ; \, 
(1+10\varepsilon) r \leq |y-x| \leq (1+30\varepsilon) r \big\}.
\leqno (25.48)
$$
Let $y\in R_1$ be given; first observe that 
$|g(y)-y| \leq 3\varepsilon r$ by (25.25) and (25.22), hence
$|g(y)-x| \geq 7\varepsilon r$ (by definition of $R_1$), and 
(25.29) says that 
$$
\varphi(y) = \psi(g(y)) = g(y)
\ \hbox{ for } y\in R_1. 
\leqno (25.49)
$$
We'll need to know that 
$$
|g(y)-g(z)| \leq 2 |y-z|
\ \hbox{ for } y, z \in E_k \cap 2B
\leqno (25.50)
$$
so we return to the definition (25.24), 
write $g(y) = \alpha \pi(y) + (1-\alpha) y$,
with $\alpha = \alpha(|y-x|)$, and similarly 
$g(z) = \beta \pi(y) + (1-\beta) y$
with $\beta = \alpha(|z-x|)$, and write that
$$\eqalign{
|g(y)-g(z)| &
= |\alpha \pi(y) + (1-\alpha) y -\beta \pi(y) - (1-\beta) y|
\cr&
\leq |(\alpha-\beta)(\pi(y)-y)| +
\beta |\pi(y)-\pi(z)| + (1-\beta)|y-z|
\cr&
\leq 3 \varepsilon r |\alpha-\beta| + |y-z| \leq 2 |y-z|
}\leqno (25.51)
$$
by (25.22), and because 
$|\alpha-\beta| \leq (10\varepsilon)^{-1} |y-z|$ 
by the definition below (25.24). So (25.50) holds.

By (25.50), the definitions (25.38) and (25.39),
and the fact that $\varphi = g$ on $R_1$, 
each $\varphi_t$ is also $2$-Lipschitz on $R_1$.
In addition, notice that for $y\in R_1$,
$$
|\varphi_t(y)-y| \leq |\varphi(y)-y| 
= |g(y)-y| \leq 3\varepsilon r
\leqno (25.52)
$$
by (25.38) and (25.39), (25.49), and (25.25) and (25.22).
Thus $\dist(\varphi_t(y),F) \leq C(\Lambda) \varepsilon r$
by (25.44), and the proof of (25.45) also yields
$$
|\Pi_F(\varphi_t(y),s)- \Pi_F(\varphi_t(y),s')| 
\leq C(\Lambda) |s-s'| \varepsilon r
\leqno (25.53)
$$
for $0 \leq s,s' \leq 1$
(just replace $\dist(y,F)$ with $\dist(\varphi_t(y),F)$).
Recall that the $\varphi_t(y)$ stay in $2B$, 
where $\Pi_F$ is well defined and all the formulas that we use
make sense (because $r \leq r(x)$ and if $r(x)$ is chosen small 
enough). Similarly, we still have (as in (25.47)) that
$$
|\Pi_F(\varphi_t(y),s)- \Pi_F(\varphi_t(z),s)| 
\leq C(\Lambda) |\varphi_t(y)-\varphi_t(z)|
\leq C(\Lambda) |y-z|,
\leqno (25.54)
$$
for $y, z \in R_1$ and $0 \leq s \leq 1$ 
by (25.34), and because $\pi^{\wt F}$ is $C$-Lipschitz
and $\varphi_t$ is $2$-Lipschitz on $R_1$.
Hence, 
$$\leqalignno{
|h_t(y)-h_t(z)| &= |\Pi_F(\varphi_t(y),s(y,t))-\Pi_F(\varphi_t(z),s(y,t))|
\cr&
\leq |\Pi_F(\varphi_t(y),s(y,t))-\Pi_F(\varphi_t(y),s(z,t))| 
\cr&\hskip3.2cm
+ |\Pi_F(\varphi_t(y),s(z,t))-\Pi_F(\varphi_t(z),s(z,t))|
& (25.55)
\cr&
\leq C(\Lambda) \varepsilon r |s(y,t)-s(z,t)| + C(\Lambda) |y-z|
\leq C(\Lambda) |y-z|
}
$$
by (25.40), (25.53), (25.54), and our definition of $s(y,t)$ in (25.36).
Thus $h_t$ is $C(\Lambda)$-Lipschitz on $R_1$, and in particular
$$
\H^d(h_1(R_1)) \leq C(\Lambda) \H^d(R_1).
\leqno (25.56)
$$

Next we consider
$$
R_2 = \big\{ y\in E_k \, ; \, |y-x| < (1+10\varepsilon) r 
\hbox{ and } g(y) \notin B \big\}.
\leqno (25.57)
$$
Let us first check that 
$$
h_t(y) \in B(x,2\Lambda^2 r)
\hbox{ when $y\in E_k \cap B(x,(1+10\varepsilon) r)$ and $0 \leq t \leq 1$.}
\leqno (25.58)
$$
When $t \leq 1/2$, $\varphi_t(y) = y$ by (25.38), 
hence $h_t(y) = \Pi_F(y,s(y,t))$ by (25.40). 
Then $|h_t(y)-y| \leq C(\Lambda) \varepsilon r < r/2$ by (25.46) 
and if $\varepsilon$ is small enough, and $h_t(y) \in 2B$ (which
is better than promised). 

When $t \geq 1/2$, $s(y,t) = s(y,1/2) = 1$
by (25.37) and (the sentence below) (25.36), hence (25.40) yields
$$
h_t(y) = \Pi_F(\varphi_t(y),1) = \pi_F \circ \varphi_t(y)
\leqno (25.59)
$$
(compare (25.34) with (25.33)). Notice that $\pi_F(x) = x$
(by (25.33), because $\psi(\lambda x) \in \wt F$
(since $x\in F$), and by (3.6)); then
$$
|h_t(y)-x| = |\pi_F \circ \varphi_t(y) - \pi_F(x)|
\leq  \Lambda^2 |\varphi_t(y)-x| 
\leqno (25.60)
$$
because $\pi_F$ is $\Lambda^2$-Lipschitz by (25.33).
If $g(y) \in B$, we also get that 
$\varphi(y) = \psi \circ g(y) \in B$, by (25.28),
hence $\varphi_t(y) \in 2B$ by (25.38) and (25.39);
then (25.60) says that $|h_t(y)-x| < 2\Lambda^2 r$,
as needed for (25.58).

We are left with the case when $g(y) \notin B$, i.e., when
$y\in R_2$. We claim that
$$
\varphi(y) = \pi(y)
\ \hbox{ for } y\in R_2.
\leqno (25.61)
$$
As soon as we prove this, we will get that
$\varphi_t(y) \in 2B$ (by (25.38) and (25.39)), 
and (25.58) will follow from (25.60). 

Now we prove the claim. Let $y\in R_2$ be given.
By definition of $R_2$, $g(y) \notin B$; 
hence by (25.25) and (25.22), $|y-x| \geq (1-3 \varepsilon)r$.
Since $|y-x| \leq (1+10\varepsilon) r$ by definition of $R_2$,
we get that $\alpha(|y-x|) = 1$ (see below (25.24)), and
$g(y) = \pi(y) \in P(x)$ by (25.24). Since $g(y) \notin B$,
(25.29) yields $\varphi(y) = \psi \circ g(y) = g(y) = \pi(y)$,
as needed for (25.61).

We are a little more interested in what happens for $t=1$. Then
(25.59), (25.39), and (25.61) say that
$$
h_1(y) = \pi_F \circ \varphi_1(y) 
= \pi_F \circ \varphi(y) = \pi_F \circ \pi(y).
\leqno (25.62)
$$
This is good, because it means that $h_1$ is Lipschitz on
$R_2$, with a constant that depends on $\Lambda$, but not on
$\varepsilon$, for instance. Then
$$
\H^d(h_1(R_2)) \leq C(\Lambda) \H^d(R_2).
\leqno (25.63)
$$
We are left with the region 
$$
R_3 = \big\{ y\in E_k \, ; \, |y-x| \leq (1+10\varepsilon) r 
\hbox{ and } g(y) \in B \big\}.
\leqno (25.64)
$$
On this last region, we do not control the Lipschitz norm
of $h_1$ (because we do not control the Lipschitz norm of $\psi$),
but fortunately (25.59) and (25.39) yield 
$h_1(y) = \pi_F \circ \varphi(y) = \pi_F \circ \psi(g(y))$, so 
$$
h_1(R_3) \i \pi_F \circ \psi(g(R_3)) \i \pi_F \circ \psi(g(E_k) \cap B)
= \pi_F \circ \psi(S) \i \pi_F(P(x) \cap \d B)
\leqno (25.65)
$$
because we set $S=g(E_k) \cap B$ (above Lemma 25.27)
and by definition of $\psi$ (below that lemma).
Since $\pi_F$ is Lipschitz, we get that
$$
\H^d(h_1(R_3)) =0.
\leqno (25.66)
$$

\ms
We want to apply the definition of a quasiminimal set,
so we check that the $h_t$ satisfy the conditions (1.4)-(1.8),
relative to the ball $2\Lambda^2 B$.
The continuity and Lipschitz
conditions (1.4) and (1.8) are satisfied (all our maps are Lipschitz,
even though with possibly huge constants), and (1.5) follows 
from (25.42). 

For (1.6), we just need to check that
$h_t(y) \in B(x,2\Lambda^2 r)$ when 
$y\in E_k \cap B(x,(1+30\varepsilon)r$, because otherwise (25.42)
says that $h_t(y)=y$. When $y\in B(x,(1+10\varepsilon)r$, this follows
from (25.58), so we may assume that $y \in R_1$ (see (25.48)).
By (25.52), $|\varphi_t(y)-y| \leq 3\varepsilon r$,
and so (25.40) yields
$$\eqalign{
|h_t(y)-x| &= |\Pi_F(\varphi_t(y),s(y,t))-x|
\leq |\Pi_F(\varphi_t(y),0)-x| + C(\Lambda) \varepsilon r
\cr&
\leq |\varphi_t(y)-x| + C(\Lambda) \varepsilon r \leq 2r
}\leqno (25.67)
$$
by (25.53) with $s'=0$, because $\Pi_F(\varphi_t(y),0)=\varphi_t(y)$
and if $\varepsilon$ is small enough; (1.6) follows.

As usual, we end the verification with the boundary condition (1.7). 
Let $y\in E_k$ be given, suppose $y\in L_j$ for some $j$,
and let us check that $h_t(y) \in L_j$ for $0 \leq t \leq 1$.
There is nothing to check if $|y-x| \geq  (1+30\varepsilon) r$,
because $h_t(y)=y$ by (25.42). Otherwise, let $G$ be a face of our grid
that contains $y$.

First assume that 
$(1+20\varepsilon)r \leq |y-x| \leq (1+30\varepsilon)r$.
Then $\varphi_t(y)=y$ by (25.41), and $h_t(y) = \Pi_F(y,s(y,t)$
by (25.40). Both $\pi_F$ and $\Pi_F$ 
were designed to preserve all the faces of our grid: see Lemma 3.4 for
$\pi^{\wt F}$, observe that $s \pi^{\wt F} + (1-s) I$
preserves the faces of the usual dyadic grid too 
(by convexity of the faces), and then $\pi_F$ and $\Pi_F$
preserve the face $G$, because we conjugate with $\psi(\lambda \cdot)$ 
(see (25.33) and (25.34)). Thus $h_t(y) \in G$, as needed.

So we may assume that $|y-x| \leq (1+20\varepsilon)r$.
For $t \leq 1/2$, $\varphi_t(y)=y$ by (25.38), so $h_t(y) = 
\Pi_F(y,s(y,t))$, and we get that $h_t(y) \in G$ by the same
argument as above. So we restrict to $t \geq 1/2$. Then
$s(y,t) = s(y,1/2) = 1$, by (25.37) and (the line below)
(25.36). Then (25.40), together with (25.33) and (25.34),
yields 
$$
h_t(y) = \Pi_F(\varphi_t(y),1) = \pi_F(\varphi_t(y)).
\leqno (25.68)
$$

Our next case is when 
$(1+10\varepsilon)r \leq |y-x| \leq (1+20\varepsilon)r$;
then $y\in R_1$ (see (25.48)) and (25.52) says that
$|\varphi_t(y)-y| \leq 3 \varepsilon r$. Notice that
$\dist(\varphi_t(y),F) \leq \dist(y,F) +3 \varepsilon r
\leq C(\Lambda) \varepsilon r$ by (25.44), and hence,
if $\varepsilon$ is small enough,
$h_t(y) = \pi_F(\varphi_t(y)) \in F$ by definition of 
$\pi_F$ (see near (25.33), and then Lemma 3.4 and Remark 3.25).
But by (25.31), $G$ contains $x$; since $F$ was chosen 
(below (25.32)) as the smallest face that contains $x$, 
we get that $F \i G \i L_j$, as needed.

Next we assume that $y\in R_2$. In this case we still have 
that $h_t(y) = \pi_F(\varphi_t(y))$, by (25.59), and in addition
$\varphi(y) = \pi(y)$ by (25.61). In this case,
$|\varphi_t(y)-y| \leq |\varphi(y)-y| = \dist(y,P(x))
\leq 3 \varepsilon r$
by (25.39), the definition of $\pi$ above (25.24), and (25.22).
Thus $\dist(\varphi_t(y),F) \leq C(\Lambda) \varepsilon r$
again, and we may conclude as before.

We are left with the case when $y\in R_3$. Then
$g(y) \in B$, and by (25.29),
$\varphi(y) = \psi(g(y)) \in P(x) \cap \d B$. 
Recall that $\varphi_t(y) \in [y,\varphi(y)]$ by (25.39);
since $\dist(y,P(x)) \leq 3 \varepsilon r$ by (25.22)
and $\dist(y,B) \leq |y-g(y)| \leq 3\varepsilon r$
by (25.25) and (25.22), this yields 
$$
\dist(\varphi_t(y),P(x) \cap B) \leq 6\varepsilon r
\leqno (25.69)
$$
But we want to show that $\varphi_t(y)$ lies close to
$F$, and since $F$ is a distorted face which may not be flat, 
we shall need to show that $\varphi_t(y)$ lies close to $E_k$,
and then use (25.44).

Unfortunately, we shall need to use Lemma 9.14.
Recall that $x\in E^1 \i E$ (see the definitions 
above (25.17) and (25.9)), but since we want to apply the
lemma to the set $E_k$, we restrict to $k$ large, choose
$x_k \in E_k \cap B(x,\varepsilon r)$, and apply the
lemma with $y=x_k$, $t = r$ (so that $B(x_k,2t) \i B(x,3r)$),
and $P = P(x)$.
Recall that $E_k = E_k^\ast$ because we assumed (10.3).
The size condition (9.15) is satisfied because $r \leq r(x)$ and if 
$r(x)$ is small enough (this is allowed by (25.11)). 
If some $L_i$ meets $B(x,2r)$ and $G$ is a face of $L_i$ 
that meets $B(x,2r)$, (25.31) says that $G$ contains $x$; 
since $F$ is the smallest face that contains $x$,
we get that $F \i G \i L_i$. Thus the set $L$ of (9.16) contains $F$,
and our assumption (9.17) holds with $\eta = C(\Lambda) \varepsilon$,
by (25.44). Finally the assumption (9.18) is satisfied for $k$ large
(and with the constant $2\varepsilon$), by (25.13) and because 
the $E_k$ converge to $E$ (recall that we apply the lemma to $E_k$,
which is why we only get $2\varepsilon$).
Thus, if $\varepsilon$ is small enough, 
the lemma applies, and we get (9.19). That is,
$$
\dist(p,E_k) \leq 2\varepsilon r
\ \hbox{ for } p \in P(x) \cap B(x_k, 3r/2).
\leqno (25.70)
$$
Return to $y\in R_3$. For each $t \in [0,1]$,
(25.69) gives $p\in P(x) \cap \overline B$ such that 
$|p-\varphi_t(y)| \leq 6\varepsilon r$. Then (25.70)
gives $z \in E_k$ such that $|z-p| \leq 2\varepsilon r$.
In turn $z\in 2B$, so (25.44) says that 
$\dist(z,F) \leq C(\Lambda) \varepsilon r$.
Altogether $\dist(\varphi_t(y),F) \leq C(\Lambda) \varepsilon r$,
and (if $\varepsilon$ is small enough), (25.68) implies that
$h_t(y) \in F \i L_j$, as needed.

This completes our proof of (1.7). Notice that it was 
surprisingly easy to get, by requiring $r(x)$ to be small, the 
only difficulty was to compose with $\pi_F$
in a way that would not destroy good Lipschitz bounds on 
$R_1 \cup R_2$ (because we need (25.56) and (25.63)); 
this is where we used our good control on $E_k \cap 2B$. 
This completes also the verification of (1.4)-(1.8). We also have
(2.4), because by our proof of (1.6), the analogue of $\wh W$ 
is contained in $B(x,2\Lambda^2 r)$, which is compactly supported in 
$U$ if $r(x)$ was chosen small enough.

Anyway, the quasiminimality of $E_k$ now yields
$$
\H^d(W_1) \leq M \H^d(h_1(W_1)) + h r^d,
\leqno (25.71)
$$
as in (2.5), and where as usual 
$W_1 = \big\{ y \in E_k \cap 2\Lambda^2 B \, ; h_1(y) \neq y \big\}$. 
For each $y\in E_k \cap B(x,r)$, (25.24) says that $g(y) \in B(x,r)$;
hence $y\in R_3$ (see (25.64)). 
If in addition $y \notin h_1(R_3)$, then $h_1(y) \neq y$ 
and $y\in W_1$. Thus
$$
\H^d(W_1) \geq \H^d(E_k \cap B(x,r) \sm h_1(R_3))
= \H^d(E_k \cap B(x,r)) \geq C^{-1} r
\leqno (25.72)
$$
by (25.66), the local Ahlfors-regularity of $E_k$, and the fact that
$E_k$ meets $B(x,r/10)$ because $x\in E$. As usual, this holds for $k$
large (depending on $x$), and with a constant $C$ that may depend
on $\Lambda$, for instance, but not on $k$ or $x$.

On the other hand, $W_1 \i R_1 \cup R_2 \cup R_3$, by (25.42), hence
$$
\H^d(h_1(W_1)) \leq C \H^d(R_1 \cup R_2) 
\leqno (25.73)
$$
by (25.56), (25.63), and (25.66). Notice also that 
$|y-x| \geq (1-3\varepsilon) r$ for $y\in R_2$,
because $g(y) \notin B$ and by (25.25) and (25.22). Thus
$$
R_1 \cup R_2 \i A, \hbox{ where }
A = \big\{ x\in E_k \, ; \, (1-30\varepsilon) r \leq 
|y-x| \leq (1+30\varepsilon) r \big\}.
\leqno (25.74)
$$
By (25.22) again, $A$ is contained in the thin strip $H$
around $P(x)$, and we can cover $A$ by less than 
$C \varepsilon^{-d+1}$ balls $D_l$ of radius $\varepsilon r$, which we may
even choose centered on $A$. By the local Ahlfors regularity of
$E_k$ (and because these balls stay far from $\R^n \sm U$ 
if $r(x)$ is small enough), $\H^d(E_k \cap D_l) \leq C \varepsilon^d r^d$. 
We sum and get that
$$
\H^d(A) \leq C \varepsilon r^d
\leqno (25.75)
$$ 
and hence, by (25.73) and (25.74),
$\H^d(h_1(W_1)) \leq C\varepsilon r^d$. If $h$ is small enough 
(depending on $n$, $M$, and $\Lambda$ through 
the constants $C$ of (25.72)), and $\varepsilon$ is small enough 
(depending on our various constants, 
but not $x$ or $r$), this contradicts (25.71) or (25.72).
This contradiction proves that $\psi$ does not exist and finishes our 
proof of Lemma 25.27.
\qed

\ms
We may now return to our initial construction, with
$x\in E^1$, $r \leq r(x)$, and $S = g(E_k) \cap B$.
By Lemma 25.27 and Definition 25.3, we get that (25.6) holds, i.e.,
$$
J_f(P(x) \cap B(x,r)) \leq J_f(S\cap B(x,r)) + \varepsilon_x(r) r^d
\leqno (25.76)
$$
with $\varepsilon_x(r)$ coming from (25.5). 
But $J_f(P(x) \cap B) = J_f(P(x) \cap B(x,r))$
because $\H^d(P(x)\cap \d B) = 0$, and by (25.5)
$\varepsilon_x(r) \leq \varepsilon$ if $r(x)$ was chosen small enough.
Thus (25.76) implies that
$$
J_f(P(x) \cap B) \leq J_f(S) + \varepsilon r^d.
\leqno (25.77)
$$

Let $z\in S$ be given, and choose $y \in E_k$ such that $g(y) = z$. 
Notice that $|y-x| \leq (1+30\varepsilon)r$, because otherwise
(25.24) would yield $g(y) = y \notin B$. Then $y\in H$ by (25.22), and 
(25.25) says that $|z-y| = |g(y)-y| \leq 3 \varepsilon r$.

A first option is that $|z-x| \leq (1-23\varepsilon)r$;
then $|y-x| \leq (1-20\varepsilon)r$, and (25.24) yields
$g(y) = y$ (because $\alpha(|y-x|)=0$). Then $z \in E_k$,
and we get that
$$
J_f(S\cap B(x,(1-23\varepsilon)r))
\leq J_f(E_k\cap B(x,(1-23\varepsilon)r))
\leq J_f(E_k\cap B).
\leqno (25.78)
$$
If $|z-x| \geq (1-23\varepsilon)r$, then 
$|y-x| \geq (1-26\varepsilon)r$ and $y\in A$,
the annulus in (25.74). Thus $S\sm B(x,(1-23\varepsilon)r) \i g(A)$.
Now (25.50) says that $g$ is $2$-Lipschitz on $A$ and hence
$$
J_f(S\sm B(x,(1-23\varepsilon)r))
\leq J_f(g(A)) \leq b \H^d(g(A)) 
\leq 2^d b \H^d(A)
\leq C b \varepsilon r^d
\leqno (25.79)
$$
by (25.4) (observe also that $g(A)$ is rectifiable) and (25.75).
Hence
$$
J_f(P(x) \cap B) \leq J_f(S) + \varepsilon r^d
\leq J_f(E_k\cap B) + C b \varepsilon r^d + \varepsilon r^d
\leqno (25.80)
$$
by (25.77), (25.78), and (25.79). We compare with (25.21)
and get that
$$\eqalign{
J_f(E \cap B)
&\leq J_f(P(x) \cap B) + \varepsilon (b+2\omega_d+\varepsilon) r^d 
+ b \H^d((E \sm E^1) \cap B)
\cr&
\leq J_f(E_k\cap B) + C\varepsilon r^d  + b \H^d((E \sm E^1) \cap B),
}\leqno (25.81)
$$
where in the last line $C$ is allowed to depend on $b$ too.

This estimate is essentially what we wanted; note that for
each $x \in E^1$, it holds for $k$ large (how large depends on $x$),
and for a constant $C$ that does not depend on $x$ or $k$.
We need a covering argument to complete our estimate.

We have a set $E^1 \i E \cap V$, and for each $x\in E^1$
we have a family of closed balls $B = B(x,r)$,
$0 < r \leq r(x)$, which forms a Vitali covering of $E^1$.
Note also that $\H^d(E^1) \leq \H^d(E\cap V) < +\infty$
(because $V \i\i U$ for the moment). 
By Theorem 2.8 on page 34 in [Ma], 
we can extract from this large family of balls a disjoint
family $\{ B_i \}$, $i\in I$, so that
$\H^d(E^1 \sm \bigcup_{i\in I} B_i) = 0$.
Then we can choose a finite set $I_0 \i I$, such that
$$
\H^d(E^1 \sm \bigcup_{i\in I_0} B_i) \leq \varepsilon.
\leqno (25.82)
$$
For $k$ large enough, (25.81) holds for every $B_i$, $i\in I_0$,
and now
$$\eqalign{
J_f(E\cap V) &\leq \sum_{i\in I_0} J_f(E\cap B_i) 
+ J_f(E \cap V \sm \bigcup_{i\in I_0} B_i)
\cr& \leq \sum_{i\in I_0} J_f(E\cap B_i) 
+ b \H^d(E \cap V \sm \bigcup_{i\in I_0} B_i)
\leq \sum_{i\in I_0} J_f(E\cap B_i) + 2 b \varepsilon
\cr&
\leq 2 b \varepsilon + \sum_{i\in I_0} \big[J_f(E_k\cap B_i) 
+ C\varepsilon r_j^d  + b \H^d((E \sm E^1) \cap B_j) \big]
\cr&
\leq 2 b \varepsilon + J_f(E_k\cap V) + b \H^d(E \cap V \sm E^1)
+ C\varepsilon \sum_{i\in I_0} r_j^d
\cr&
\leq J_f(E_k\cap V) + C \varepsilon + C\varepsilon \sum_{i\in I_0} r_j^d
}\leqno (25.83)
$$
by (25.4), (25.14), (25.82), then by (25.81), 
where we set $B_j = \overline B(x_j, r_j)$,
because the $B_j$ are disjoint and contained in $V$
(if each $r(x)$ was chosen small enough, according to (25.11)),
and by (25.14) again. Since $E$ is locally Ahlfors-regular
and each $B_j$ is centered on $E$ and such that $10B_j \i V \i U$,
we get that 
$$
\sum_{i\in I_0} r_j^d \leq C \sum_{i\in I_0} \H^d(E \cap B_i)
\leq C \H^d(E \cap V)
\leqno (25.84)
$$
(because the $B_i$ are disjoint). Then (25.83) says that
$$
J_f(E\cap V) \leq J_f(E_k\cap V) + C(1+\H^d(E \cap V)) \varepsilon
\leqno (25.85)
$$
for $k$ large. Thus 
$J_f(E\cap V) \leq \liminf_{k \to +\infty} \, J_f(E_k \cap V)
+ C(1+\H^d(E \cap V)) \varepsilon$ and, since this estimate holds for 
every small $\varepsilon$, we get (25.8).

This takes care of the special case when $V$ is compactly 
contained in $U$. In the general case, we write $V$ as the increasing 
union of open sets $V_m \i \i U$, notice that 
$$
J_f(E\cap V_m) \leq \liminf_{k \to +\infty} \, J_f(E_k \cap V_m)
\leq \liminf_{k \to +\infty} \, J_f(E_k \cap V)
\leqno (25.86)
$$
for each $m$, then take the limit in $m$ and get (2.8) for $V$.
This completes our proof of Theorem 25.7.
\qed

\msi{\bf Remark 25.87.}
The author feels that it is a pity that we do not allow $f$ 
to be merely lower semicontinuous, but was not able to come up
with a clean statement, so we will just give two possible 
substitutes here.

First observe that we used the continuity of $f$ only twice, 
in (25.15), and (25.16) (but where lower semicontinuity would 
have been enough), to be used in the last line of (25.17), 
then in (25.21), to prove that for our nice balls $B$,
$J_f(E\cap B)$ is almost as small as $J_f(P(x)\cap B)$
(the measure of a nearby disk).

We want to replace our continuity assumption with 
the following one: for each $x\in U$, each $d$-plane $P$
through $x$, and each $C^1$, embedded, submanifold $\Gamma$
of dimension $d$ through $x$, which admits $P$ as a tangent plane
at $x$,
$$
\limsup_{r \to 0} {1 \over r^d} \big[
J_f(\Gamma \cap B(x,r)) - J_f(P \cap B(x,r)) \big] \leq 0.
\leqno (25.88)
$$
Notice that this goes in the direction opposite to (25.6);
this can be seen as a form of continuity is some direction,
possibly much weaker than the full continuity asked above,
but hard to think about as a lower semicontinuity property.
We could have given the same definition, where instead $\Gamma$
is the graph of some $C^1$ mapping $F : P \to P^\perp$,
with $DF(x)=0$, and the two definitions would have been equivalent.

\msi{\bf Claim 25.89.} Theorem 25.7 also holds when we replace
${\cal I}(U,a,b)$ with the class ${\cal I}_l(U,a,b)$ of functions
$f : U \times G(n,d) \to [a,b]$ that satisfy (25.5), (25.6), and 
(25.88).

\ms
That is, for the sake of Theorem 25.7, we can replace the 
continuity of $f$  by the condition (25.88) in the definition
of ${\cal I}(U,a,b)$.

Our claim will follow as soon as we show that, with a suitable 
modification of the set $E^1$ and, for $x\in E^1$, of the radius
$r(x)$, we still have (25.21) for $B = \overline B(x,r)$, 
when $x\in E^1$ and $r>0$ is small enough.
Recall that $E$ is rectifiable; thus we can write $E$ as
null set, plus a countable collection of sets $F_i$, where 
$F_i$ is contained in a $C^1$, embedded, submanifold $\Gamma_i$
of dimension $d$. We may even assume that the $F_i$ are disjoint. 
Then almost every point $x\in E$ lies in some $\Gamma_i$,
and is even a point of vanishing density for $E \sm \Gamma_i$,
i.e.,
$$
\lim_{r \to 0} r^{-d} \H^d(E \cap B(x,r)\sm \Gamma_i) = 0;
\leqno (25.90)
$$
see [Ma], Theorem 6.2 (2) on page 89. 
We leave the definition of $E^0$ and $E^1$ as it was, except that
we forget about the conditions (25.14)-(25.16) (which concerned the
continuity of $f$, and replace them by the constraints that for
$x\in E^1$, $x$ lies in some $\Gamma_i$, $i = i(x)$, and (25.90)
holds. When we choose $r(x)$, we require that for $0 < r \leq 2r(x)$,
$$
\H^d(E \cap \overline B(x,r)\sm \Gamma_i) < \varepsilon r^d
\leqno (25.91)
$$
and 
$$
J_f(\Gamma_{i(x)} \cap \overline B(x,r)) 
\leq J_f(P \cap B(x,r)) + \varepsilon r^d,
\leqno (25.92)
$$
which we obtain as a limit of (25.88) for $r'>r$, because
$\H^d(P \cap \d B(x,r)) = 0$.
Then for $x\in E^1$ and $0 < r \leq r(x)$, and if we set
$B = \overline B(x,r)$ as before, we get that
$$\eqalign{
J_f(E \cap B) &\leq J_f(E \cap \Gamma_{i(x)} \cap B)
+ b \H^d(E \cap B \sm \Gamma_{i(x)})
\cr&
\leq J_f(P(x) \cap B(x,r)) + (1+b)\varepsilon r^d,
}\leqno (25.93)
$$
which is even better than (25.21). Our claim follows.

The reader is probably worried about the $limsup$ in (25.88),
because the most logical statement would use a $liminf$.
The proof above accommodates a liminf too, if we are more careful.
Instead of having (25.92) for all the radii $r \leq 2 r(x)$,
we would only get it for a sequence of radii that tends to $0$.
Then, in the application of the Vitali covering lemma near (25.82),
we would only choose balls $B_i$ with a radius that satisfies (25.92).
We decided not to bother.

\ms
We can also try to take care of our semicontinuity issue
by extending the class ${\cal I}(U,a,b)$ after the fact.
That is, denote by ${\cal I}^+(U,a,b)$ the class of 
functions $f : U \times G(n,d)$ such that, for each compact
set $H \i U$, there is a sequence $\{ f_m \}$ in ${\cal I}(U,a,b)$, 
with $a \leq f_m \leq f$ everywhere, and 
$\lim_{m \to +\infty} f_m(x,T) = f(x,T)$
for $x\in H$ and $T\in G(n,d)$. We claim that
$$
\hbox{Theorem 25.7 also holds with ${\cal I}(U,a,b)$ replaced
by ${\cal I}^+(U,a,b)$.}
\leqno (25.94)
$$
This would be nice if we could characterize easily ${\cal I}^+(U,a,b)$
(for instance, by lower semicontinuity and the conditions 
(25.4)-(25.6)), but the truth is that the author does not know how to
manipulate (25.6) concretely.

Let us prove the claim anyway. Let $f\in {\cal I}^+(U,a,b)$, 
the sequence $\{ E_k \}$ and its limit $E$, and $V \i U$ be given.
As for Theorem 25.7 itself, it is enough to prove the conclusion (25.8)
when $V \i\i U$ (otherwise, write $V$ as an increasing union of 
open sets that are compactly contained in $U$). 
Then let $\{ f_m \}$ be as in the definition of ${\cal I}^+(U,a,b)$,
relative to $H = \overline V$; notice that for $m \geq 0$
$$
J_{f_m}(E\cap V) 
\leq \liminf_{k \to +\infty} J_{f_m}(E_k\cap V)
\leq \liminf_{k \to +\infty} J_{f}(E_k\cap V)
\leqno (25.95)
$$
because (25.8) holds for $f_m$ and $f_m \leq f$,
and that
$$
J_{f}(E\cap V) = \lim_{m \to +\infty} J_{f_m}(E\cap V) 
\leqno (25.96)
$$
by the dominated convergence theorem; (25.8) and our claim (25.94) follow.

\ms
The next lemma will help with the extension of Theorem 10.8
to $f$-quasiminimal sets. 

\ms\proclaim Lemma 25.97.
Let $f: U \times G(n,d) \to [a,b]$ satisfy the conditions 
(25.4)-(25.6). Then
$$
\liminf_{r \to 0} {1 \over r^d} \big[
J_f(\Gamma \cap B(x,r)) - J_f(P \cap B(x,r)) \big] \geq 0
\leqno (25.98)
$$
for each $x\in U$, each $d$-plane $P$ through $x$, and each $C^1$, embedded, 
submanifold $\Gamma$ of dimension $d$ through $x$, which admits $P$ 
as a tangent plane at $x$.

\ms
Notice that the conclusion is the opposite of (25.88), so we could 
replace (25.88) into the simpler (but apparently stronger)
$$
\lim_{r \to 0} {1 \over r^d} \big[
J_f(\Gamma \cap B(x,r)) - J_f(P \cap B(x,r)) \big] = 0
\leqno (25.99)
$$
in Claim 25.89, without changing the result.

The proof of Lemma 25.97 goes a little as for Lemma 25.27.
Let $x, P$, and $\Gamma$ be as in the statement, and let $r> 0$
be small. We want to  apply (25.6) to a suitable set 
$S \i \overline B(x,r)$, and since we want $S \cap \d B(x,r)$
to be contained in $P$, we use $S = B \cap g(\Gamma)$, where 
$B = \overline B(x,r)$, $g$ is defined by (25.24), 
$\pi$ denotes the orthogonal projection onto $P$, 
$\alpha$ is defined as below (25.24), 
and $\varepsilon > 0$ is a small positive number, which
will tend to  $0$ at the end of the argument.

As before, we want to show that there is no  Lipschitz mapping
$\psi : B \to B$ such that $\psi(y) = y$ for $y\in P\cap \d B$
and $\psi(S) \i P \cap \d B$. 
Let us suppose that $\psi$ exists and use it to define 
an impossible mapping $h : P\cap B \to P \cap \d B$. 
First extend $\psi$ to $P$, by setting $\psi(y)=\rho(y)$ 
for $y\in P \sm B$, where $\rho(y) = x + r {y-x \over |y-x|}$ 
is the radial projection of $y$ on $\d B$. The extension is still
Lipschitz because $\psi(y)=y$ on $P \cap \d B$.

By definition of $\Gamma$, there is a $C^1$ function $F : P \to  P^\perp$, 
with $DF(x) = 0$, such that for $r$ small, $\Gamma$ coincides with 
the graph of $F$ in $B(x,2r)$.
Also, for $r$ small enough, $|F(y)| \leq \varepsilon r$ for
$y\in P \cap 2B$. Set $\lambda = 1+2\varepsilon$ and, for
$y\in P\cap B$, set $\wt y = x + \lambda (y-x) \in \lambda B$,
and then $z = \wt y + F(\wt y)$ Thus $z \in B(x,(1+3\varepsilon) r)$
(because $|F(\wt y)| \leq \varepsilon r$).
If $g(z)\in B$, then $g(z) \in B \cap g(\Gamma) = S$, and
$\psi(g(z))$ is defined and lies in $P \cap \d B$. Otherwise, 
notice that 
$$
|g(z)-z| \leq |\pi(z)-z| = |F(\wt y)| \leq \varepsilon r
\leqno (25.100)
$$
by (25.24) and because $\wt y \in \lambda B$, and
$$
|z-y| \leq |F(\wt y)| + |\wt y - y| \leq 3 \varepsilon r.
\leqno (25.101)
$$
Hence $|z-x| \geq |g(z)-x|-\varepsilon r \geq (1-\varepsilon)r$
(because $g(z) \notin B$), and 
$|z-x| \leq |y-x| + 3 \varepsilon r \leq (1+3) \varepsilon r$,
so $\alpha(|x- z|) =1$, hence $g(z) = \pi(z) \in P$,
and again $\psi(g(z))$ is defined and lies in $P \cap \d B$
by definition of our extension $\psi$ on $P \sm B$. 
So we can define $h : P \cap B \to P \cap \d B$ by $h(y) = \psi(g(z))$,
and obviously $h$ is continuous. Also, if $y \in \d B$,
(25.101) still holds and yields
$|z-y| \leq 3 \varepsilon r$, then $\alpha(|x- z|) =1$,
and hence $g(z) = \pi(z)$ by (25.24). In addition,
$|g(z)-z| = |\pi(z)-z| = |F(\wt y)| \leq \varepsilon r$ as in 
(25.100), so $|g(z)| \geq |z|-\varepsilon r > r$ (by definition of $\lambda$),
which means that $g(z) \in P \sm B$ and $\psi(g(z)) = \rho(g(z))$.
Thus $|h(y)-y| = |\psi(g(z))-y| = |\rho(g(z))-y| \leq |g(z)-y|
\leq 4 \varepsilon r$; 
this implies that the restriction of $h$ to $P \cap \d B$
is of degree $1$, which is impossible because it has a continuous extension
from $P \cap B$ to $P \cap \d B$. 

This contradiction shows that $\psi$ does not exist, and this allows 
us to apply (25.6). That is, 
$$
J_f(P \cap B(x,r)) \leq J_f(S \cap B(x,r)) + \varepsilon_x(r) r^d
= J_f(g(\Gamma)\cap B(x,r)) + \varepsilon_x(r) r^d.
\leqno (25.102)
$$
We claim that for $r$ small,
$$
J_f(g(\Gamma)\cap B(x,r)) \leq J_f(\Gamma\cap B(x,r)) + 
C\varepsilon r^d.
\leqno (25.103)
$$
Let $z\in g(\Gamma) \cap B(x,r)$ be given, and let $y\in \Gamma$
be such that $g(y)=z$; if $|y-x| \leq (1 - 20\varepsilon) r$,
$\alpha(|y-x|) = 0$, hence $g(y) = y$. The corresponding subset
of $g(\Gamma)$ is controlled by the first term in the second hand of 
(25.103). The case when $|y-x| \geq (1 + 20\varepsilon) r$ is 
impossible, because we would have that $g(y)=y$ for the same reasons.
We are left with $g(A)$, where
$A = \big\{ y \in \Gamma \, ; \, (1 - 20\varepsilon) r \leq |y-x| 
\leq (1 - 20\varepsilon) r \big\}$. We observe that 
$\H^d(A) \leq C \varepsilon r^d$, and that $g$ is $C$-Lipschitz on $A$
(recall that $|\pi(y)-y| \leq \varepsilon r$ for $y\in A$,
and use the usual argument).
This proves (25.103), and because of (25.102) we get that the $liminf$
in (25.98) is larger than $-C\varepsilon$; since $\varepsilon > 0$
is arbitrarily small, (25.98) and Lemma 25.97 follow.
\qed

\msi
{\bf 26. Limits of $f$-quasiminimal sets associated to elliptic integrands.} 
\ms

We shall now describe a few implications of  
Theorem 25.7, in a context of quasiminimal and almost minimal
sets relative to an integrand in the class ${\cal I}(U,a,b)$.
We could also use the slightly larger classes 
${\cal I}_l(U,a,b)$ and ${\cal I}^+(U,a,b)$
defined for Claim 25.89 and (25.94), but we shall stick
to ${\cal I}(U,a,b)$ for simplicity.

We start with some simple observations on quasiminimal sets.
Let $f\in {\cal I}(U,a,b)$ be given. Since we want to define
$J_f(E)$ also for sets $E$ that are not necessarily rectifiable,
define an auxiliary function $\wt f : U \to (0,+\infty)$, also with 
$a \leq \wt f(x) \leq b$; this way we can define 
$J_{f,\wt f}(E)$ as in (25.2). Of course 
$J_{f,\wt f}(E) = J_f(E)$, as defined in (25.1), when $E$ is
rectifiable (which will be our main case).

Then, we can define the class $GSAQ_f(U,M,\delta,h)$, 
as we did in Definition 2.3, except that
we replace (2.5) with the corresponding inequality
$$
J_{f,\wt f}(W_1) \leq M J_{f,\wt f}(\varphi_1(W_1)) + h r^d.
\leqno (26.1)
$$
With our assumptions, notice that
$$
a \H^d(A) \leq J_{f,\wt f}(A) \leq b \H^d(A) 
\leqno (26.2)
$$
when $A$ is a Borel set such that $\H^d(A) < +\infty$.
Then it is easy to see that
$$
E \in GSAQ(U,a^{-1}bM,\delta,a^{-1}h)
\, \hbox{ as soon as } \,
E \in GSAQ_f(U,M,\delta,h).
\leqno (26.3)
$$
That is, quasiminimal sets relative to $f$ are also
quasiminimal relative to $1$, and if $h$ is small enough
(now depending on $a$ and $b$ as well), 
Theorem 5.16 says that $E$ is rectifiable. Then we can forget
about $\wt f$ altogether (since $E$ and also its competitors
$\varphi_1(E)$, where $\varphi_1$ satisfies (1.8), are rectifiable),
and concentrate on $f$ and the formula (25.1). In particular, our 
class $GSAQ_f$ does not depend on $\wt f$.

We don't need to worry about the regularity results for 
$E \in GSAQ_f(U,M,\delta,h)$, since we can apply the results that we 
proved for plain quasiminimal sets.

Also, Theorem 25.7 applies to quasiminimal sets 
$E\in GSAQ_f(U,M,\delta,h)$, $h$ small enough, since they are
plain quasiminimal sets.

\msi
{\bf Claim 26.4.}
Theorem 10.8 is still valid when we take $g\in {\cal I}_l(U,a,b)$,
and replace $GSAQ(U,M,\delta,h)$ with $GSAQ_g(U,M,\delta,h)$
both in the assumption (10.2) and the conclusion (10.9).

We decided to call our integrand $g$ because the letter
$f$ is  used for the Lipschitz map of Sections 11-19.
We work with the slightly larger class 
${\cal I}_l(U,a,b)$ of Claim 25.89 because the proof 
for $g\in {\cal I}(U,a,b)$ works as well with ${\cal I}_l(U,a,b)$;
we shall not attempt to see what happens when $g\in {\cal I}^+(U,a,b)$.

Also notice that under the Lipschitz assumption,
this time we restrict to the additional condition (10.7), which is easy
to use, and do not attempt to use the weaker (19.36).

Because of the length of the proof, we shall not check every detail,
so the reader is invited to use a little more caution than usual
before applying this result.

Most of the construction of stable competitors, as in Sections 10-17,
does not need to be changed (we shall just modify the definition of
the radii $r(y)$ defined in (15.4), before we define the balls 
$B_j$, $j\in J_3$). 
In particular, the estimates for all the 
small perturbation pieces will give equivalent results when we 
estimate sets with $J_g$ rather than $\H^d$, because of (25.4).
Even in Section 18, nothing much happens before the estimates
near (18.58), when we study the main contribution from the 
$B_{j,x}^+$, $j\in J_3$ and $x\in Z(y_j)$.

Under the rigid assumption, we still have (18.58), 
for the same reason as before, but instead of
(18.63) and (18.64), we use this to prove that
$$\eqalign{
J_g\Big(h_2\Big(\bigcup_{j\in J_3} \bigcup_{x\in Z(y_j)} B_{j,x}^+ \sm R^3\Big)\Big)
&\leq
J_g\Big(h_2\Big(\bigcup_{j\in J_3} \bigcup_{x\in Z(y_j)} B_{j,x}^-)\Big)\Big)
\cr&
\leq \sum_{j\in J_3} J_g(Q_j \cap D_j),
}\leqno (26.5)
$$
where, as we recall, $D_j = B(y_j,r_j)$ is a ball and $Q_j$ is a $d$-plane
through $y_j$; we then deduce from this and previous estimates that
(as in (18.64))
$$\eqalign{
J_g(h_2(E_k \cap W)) 
&\leq 
C \eta + C(f,\gamma) (1-a) + C(\alpha,f) N^{-1} + C(f) \gamma
+ \sum_{j\in J_3} J_g(Q_j \cap D_j).
}\leqno (26.6)
$$
Under the Lipschitz assumption, the 
same estimates as before lead to the following analogue of (18.72):
$$\eqalign{
J_g(h_2(E_k \cap W)) 
&\leq 
C \eta + C(f,\gamma) (1-a) + C(\alpha,f) N^{-1} + C(f) \gamma
\cr& \hskip2.5cm
+ \sum_{j\in J_3} J_g(D_j \cap \lambda^{-1}\psi^{-1}(\wt Q_j)).
}\leqno (26.7)
$$
Now we need to change the argument a little, because we want
lower bounds for $J_g(D_j \cap f(E\cap W_f))$ that fit with 
(26.6) or (26.7).

We begin with the rigid case, and explain how to modify the
definition of the $D_j$ in Section 15, so as to have additional useful
properties. We start with the sets $X_9 \i E$ and $Y_9 = f(X_9)$,
and of course $Y_9$ is rectifiable because $f$ is Lipschitz.
This means that we can find a countable collection $\{ \Gamma_m \}$,
$m \geq 0$, of $C^1$ submanifolds of dimension $d$, such that
$\H^d(Y_9 \sm \bigcup_{m \geq 0} \Gamma_m) = 0$.

Set $\Gamma'_m = \Gamma_m \sm \bigcup_{l<m} \Gamma_l$,
and then $Y_9(m) = Y_9 \cap \Gamma'_m$; thus the $Y_9(m)$
are disjoint, and almost cover $Y_9$. Then denote by
$Y'_9(m)$ the set of $y\in Y_9(m)$ such that
$$
\lim_{r \to 0} r^{-d} \H^d(B(y,r) \cap \Gamma_m \sm Y_9(m)) = 0;
\leqno (26.8)
$$
we know from [Ma], Theorem 6.2 (2) on page 89 
that (26.8) holds for $\H^d$-almost every $y\in Y_9(m)$
(because $Y_9(m) \i \Gamma_m$), so
$\H^d(Y_9 \sm \bigcup_{m \geq 0} Y'_9(m)) = 0$.
Set $X'_9(m) = X_9 \cap f^{-1}(Y'_9(m))$ for $m \geq 0$;
we claim that
$$
\H^d\big(X_9 \sm \bigcup_{m \geq 0} X'_9(m)\big) = 0.
\leqno (26.9)
$$
The justification is the same as for (15.11), relies on the
fact that $f : X_9 \to Y_9$ is at most $N$-to-$1$, and is done 
for (4.77) in [D2].  

For each $m \geq 0$ and $y\in Y'_9(m)$, there is a radius $r_1(y)$
such that
$$
\H^d(B(y,r) \cap \Gamma_m \sm Y_9(m)) \leq \eta r^d
\ \hbox{ for } 0 < r \leq r_1(y)
\leqno (26.10)
$$
(where $\eta > 0$ is the usual small number in Sections 11-18)
and, because of Lemma 25.97,
$$
J_g(Q_m(y) \cap B(y,r)) \leq J_g(\Gamma \cap B(y,r)) + \eta r^d
\ \hbox{ for } 0 < r \leq r_1(y),
\leqno (26.11)
$$
where $Q_m(y)$ denotes the tangent plane to $\Gamma_m$ at $y$.
We now modify our definition of $Y_{10}$ in Section 15. 
We replace $Y_9$ by $Y'_9 = \bigcup_{m \geq 0} Y'_9(m)$
and $X_9$ by $X'_9 = \bigcup_{m \geq 0} X'_9(m)$ (we know from (26.9)
that we don't lose any mass), and in addition to
the defining condition (15.2) on $r(y)$, we require
that $r(y) \leq r_1(y)$ for $y\in Y'_9$. Then we define $\delta_7$,
$\delta_8$, and the sets $X_{10}$ and $Y_{10}$ as before, except that
in (15.7) and (15.8) we replace $X_9$ and $Y_9$ with $X'_9$ and $Y'_9$.
This way we get the additional property that
$r_1(y) > \delta_8$ when $y\in Y_{10}$, and in particular,
once we choose the balls $D_j = B(y_j,r_j)$, $j\in J_3$, that
$r_j < r_1(y_j)$ for $j\in J_3$ (by (15.12)).

Now fix $j\in J_3$, and let $m$ be such that $y_j \in Y'_9(m)$.
We have a $d$-plane $Q_j$, which is the common value of the
$A_x(P_x)$, $x\in Z(y_j)$ (see above (15.16), and we claim that
it is also equal to the tangent plane $Q_m(y_j)$ to $\Gamma_m$
at $y_j$. Since both are $d$-dimensional, it is enough to check that 
$Q_m(y_j) \i Q_j$. Let $v$ be a unit vector in the direction of 
$Q_m(y_j)$, and let $\varepsilon > 0$ be given. For $r > 0$,
(26.8) says that $B(y_j + r v/2, \varepsilon r)$ meets
$Y_9(m)$. Since $Y_9(m) \i Y_9 = f(X_9)$, we can find
$z \in X_9$ such that $f(z) \in B(y_j + r v/2, \varepsilon r)$.
By (15.2), $z\in B(x,2\gamma^{-1}r)$ for some $x\in Z(y_j)$.
By (11.40), $|f(z)-A_x(z)| \leq \varepsilon |z-x| \leq 
2\gamma^{-1}\varepsilon r$ if $r$  is small enough (recall that 
$Z(y_j)$ is finite). At the same time, 
$\dist(z,P_x) \leq \varepsilon r$ for $r$ small, because $P_x$
is tangent to $E$ at $x$. Let $\wt z$ denote the projection of
$z$ on $P_x$; then
$$\leqalignno{
\dist(y_j + r v/2,Q_j) 
&\leq \varepsilon r +  \dist(f(z),Q_j)
= \varepsilon r +  \dist(f(z),A_x(P_x))
\cr&
\leq \varepsilon r +  \dist(f(z),A_x(\wt z))
\leq \varepsilon r + |f(z)-A_x(z)| +|A_x|_{lip} |z-\wt z|
\cr&
\leq \varepsilon r + 2\gamma^{-1}\varepsilon r + |f|_{lip} \varepsilon r
\leq C \varepsilon r
&(26.12)
}
$$
by (11.36). For each $\varepsilon > 0$, this holds for $r$ small 
enough; since $y_j \in Q_j$, it follows that $v$ lies in the vector
space parallel to $Q_j$, and $Q_m(y_j) = Q_j$, as needed.

Since $r_j < r_1(y_j)$, we can apply (26.11) and then (26.10) to get that
$$\eqalign{
J_g(Q_j \cap D_j) &= J_g(Q_m(y_j) \cap B(y_j,r_j))
\leq J_g(\Gamma_m \cap B(y_j,r_j)) + \eta r_j^d
\cr&
\leq J_g(Y_9(m)\cap B(y_j,r_j)) 
+ b \H^d(B(y_j,r_j) \cap \Gamma_m \sm Y_9(m))
+ \eta r_j^d
\cr&
\leq J_g(Y_9(m)\cap B(y_j,r_j)) + (1+b) \eta r_j^d.
}\leqno (26.13)
$$
Now $Y_9(m) \i f(E \cap W_f))$, because $Y_9(m) \i Y_9 = f(X_9)$,
and $X_9 \i X_0 = E \cap W_f$ by (11.20), so (26.13) says that
$$
J_g(D_j \cap Q_j) \leq J_g(D_j \cap f(E \cap W_f)) + (1+b) \eta r_j^d.
\leqno (26.14)
$$
We sum over $j \in J_3$, and get that
$$\eqalign{
\sum_{j\in J_3} J_g(D_j \cap Q_j) 
&\leq \sum_{j\in J_3} J_g(D_j \cap f(E \cap W_f)) 
+ (1+b) \eta \sum_{j\in J_3} r_j^d
\cr&
\leq J_g(f(E \cap W_f)) + (1+b) \eta \sum_{j\in J_3} r_j^d
}\leqno (26.15)
$$
because the $D_j$ are disjoint. Since
$J_g(D_j \cap Q_j) \geq a \H^d(D_j \cap Q_j) = a \omega_d r_j^d$
for $j\in J_3$, we deduce from this that
$$
\sum_{j\in J_3} r_j^d 
\leq a^{-1} \omega_d^{-1} \sum_{j\in J_3} J_g(D_j \cap Q_j)
\leqno (26.16)
$$
and so, by (26.15) and if $\eta$ is small enough,
$$
\sum_{j\in J_3} r_j^d \leq 2a^{-1} \omega_d^{-1} J_g(f(E \cap W_f))
\leqno (26.17)
$$
because the $D_j$ are disjoint.
We now compare (26.15) with (26.6), and get that
$$\eqalign{
J_g(h_2(E_k \cap W)) 
&\leq J_g(f(E \cap W_f)) + {\cal E},
}\leqno (26.18)
$$
where
$$
{\cal E} = C \eta + C(f,\gamma) (1-a) + C(\alpha,f) N^{-1} + C(f) \gamma
+ (1+b) \eta \sum_{j\in J_3} r_j^d
\leqno (26.19)
$$
is a small error term (observe that $\sum_{j\in J_3} r_j^d \leq C(f)$,
by (26.17)). This is a good substitute for (18.93); from there, we 
estimate the difference between $W_f$ and $W$ as in (18.96), replace 
the lower semicontinuity estimate (18.97) by (25.8), and end the 
proof as before, with $\H^d$  replaced with $J_g$. This completes the proof 
under the rigid assumption.

\ms
Now suppose that we only have the Lipschitz assumption; thus we only 
have the estimate (26.7), and as before in Section 19, we need to 
estimate the quantity
$$
\Delta = \sum_{j\in J_3} J_g(D_j \cap \lambda^{-1}\psi^{-1}(\wt Q_j))
- \sum_{j\in J_3} J_g(Q_j \cap D_j)
\leqno (26.20)
$$
(compare with (19.1)). The first part of Section 19, where for $i\in J_4$, 
we extend our one parameter family to get a final set 
which in $D_j$ is almost contained in $Q_j$, does not need to be 
modified. We get an estimate like (19.32), with $\H^d$  replaced by 
$J_g$, which gives a contribution like the one we had in (26.6), and 
the effect is that we can remove from $\Delta$ the contribution 
from the indices $j\in J_4$.

For the second part of the argument, where we get rid of some small
set $Z$, we need to change a few definitions. For $y\in U$, denote
by $F(y)$ the smallest face of our grid that contains $y$.

Our first set $Z_1$ is the set of points $y\in U$ such that 
${\rm dimension}(F(y))>d$, that lie in $L'_i = L_i \sm {\rm int}(L_i)$ 
for some $i \in [0,j_{max}]$, but for which we cannot find $t > 0$ 
such that the restriction of $\psi$ to $\lambda F(y) \cap B(\lambda y,t)$ 
is of class $C^1$. By (10.7) (or rather the translation that was given below its statement), $\H^d(Z_1) = 0$.

Next denote by $Z_2$ the union of all the faces $F$ such that
${\rm dimension}(F) < d$; again $\H^d(Z_2) = 0$.

Our third small set $Z_3$ is a subset of $Y_{11}$
(defined by (15.9)). Consider the set $Y$
of points $y\in Y_{11}$ such that ${\rm dimension}(F(y)) = d$.
This set is rectifiable, so we can find a countable collection of
$C^1$ submanifolds $\Gamma_m$, $m \geq 0$, of dimension $d$, such that
if we set $Y' = Y \cap\big(\bigcup_m \Gamma_m\big)$, then
$\H^d(Y\sm Y') = 0$. The $\Gamma'_m = \Gamma_m \sm \bigcup_{l < m} \Gamma_l$
are disjoint, and still cover $Y'$. Now for each face $F$ of dimension $d$
and each $m$, we can apply 
[Ma], Theorem 6.2 (2) on page 89 
to show that for $\H^d$-almost every $y\in Y' \cap F \cap \Gamma'_m$,
$$
\lim_{r \to 0} r^{-d}
\H^d(B(y,r) \cap \Gamma_m \sm Y') = 0
\leqno (26.21)
$$
(say that $\Gamma_m \sm Y' \i [\Gamma_m \sm \Gamma'_m] \cup 
[\Gamma'_m \sm Y']$ and observe that $\H^d(\Gamma_m)$ is locally 
finite) and 
$$
\lim_{r \to 0} r^{-d}
\H^d(B(y,r) \cap F \sm \Gamma'_m) = 0.
\leqno (26.22)
$$
We remove from $Y_{11}$ the set $Z_3$ of $y \in F_{11}$ such that
$y\in Y \sm Y'$, or $y\in Y'$ but (26.21) or (26.22) fails for 
$F = F(y)$ and $m = m(y)$, the index such that $y\in G'_m$. 
That is, we set $Y_{12} = Y_{11} \sm (Z_1 \cup Z_2 \cup Z_3)$.
Of course $\H^d(Y_{11} \sm Y_{12}) = 0$.

Now let $y\in Y_{12}$ be given. 
We want to define a radius $r_1(y)$ under which some good things happen. 
We start in the special case when
${\rm dimension}(F(y)) = d$, set $F = F(y)$, and let $m$ be such that
$y \in G'_m$; thus (26.21) and (26.22) hold, because we excluded $Z_3$.

First, we shall take $r_1(y)$ so small that 
$$
\H^d(B(y,r) \cap F \sm \Gamma_m) \leq \eta r^d
\leqno (26.23)
$$
for $0 < r \leq r_1(y)$; this is easy because 
of (26.22). For the next condition, pick any point
$x\in Z(y)$; such a point exists because $y\in Y_9 = f(X_9)$
(see (15.1) and the line before), and in addition
$y=f(x)$ and $x\in X_{11}$ (see (15.7) and (15.10)).
Set $\wt y = \psi(\lambda y) = \wt f(x)$ (see (12.36))
and $\wt F = \psi(\lambda F)$; this last is the smallest
face of the true dyadic grid that contains $\wt y$.
By Lemma 12.40, $\wt A_x(P_x) \i W(\wt f(x)) = W(\wt y)$, where
$W(\wt y)$ denotes the smallest affine subspace that contains
$\wt F$; since $F$, and hence $W(\wt y)$, are 
$d$-dimensional, we will immediately get that
$$
\wt A_x(P_x) = W(\wt y)
\leqno (26.24)
$$
as soon as we check that $\wt A_x(P_x)$ is $d$-dimensional.
We know that $A_x(P_x)$ is $d$-dimensional (compare the
definition (14.5) with (14.21) and the line below);
then our proof of (15.41) shows that shows that the restriction
of $\psi$ to $\lambda A_x(P_x)$ is differentiable in every
direction, and (15.41) gives a relation between the directional
derivatives of $A_x$ and $\wt A_x$ on $P_x$ (and at $x$), 
which proves that $D \wt A_x$ is injective (because $D\psi$
is injective since $\psi$ is bilipschitz). So 
$\wt A_x(P_x)$ is $d$-dimensional and (26.24) holds.
No $y$ lies in the interior of $F$ (by definition of $F$
as the smallest face that contains $y$, hence $\wt y$ lies
in the interior of $\wt F$, the sets 
$\wt A_x(P_x)$, $W(\wt y)$, and $\wt F$ coincide near $y$,
and (by applying the bilipschitz map $\lambda^{-1} \psi^{-1}$)
the sets $F$ and $\lambda^{-1} \psi^{-1}(\wt A_x(P_x))$
coincide near $y$. Now set $\wt Q = \wt A_x(P_x)$; we even know
from (15.40) that all the $x\in Z(y)$ give the same $\wt Q$.
If $r_1(y)$ is small enough, then for $0 < r < r_1(y)$,
$$\eqalign{
J_g(B(y,r) \cap \lambda^{-1}\psi^{-1}(\wt Q))
&= J_g(B(y,r) \cap F)
\cr&
\leq J_g(B(y,r) \cap \Gamma_m) + b \H^d(B(y,r) \cap F \sm \Gamma_m)
\cr&
\leq J_g(B(y,r) \cap \Gamma_m) + b \eta r^d
}\leqno (26.25)
$$
by (26.23). 
Denote by $P(\Gamma_m)$ the tangent $d$-plane to $\Gamma_m$
at $y$. We claim that $P(\Gamma_m) = A_x(P_x)$.
We know from (15.9) that $A_x(P_x)$ does not depend on $x$; 
since both sets are $d$-dimensional, it is enough to check that
$P(\Gamma_m) \i A_x(P_x)$. We then proceed as we did near (26.12).
Let $v$ be a unit vector in the direction of 
$P(\Gamma_m)$, and let $\varepsilon > 0$ be given. 
For $\rho > 0$, (26.21) says that 
$B(y + \rho v/2, \varepsilon \rho)$ meets
$Y' \i Y \i Y_{11}$. Since $Y_{11} = f(X_{11})$, we can find
$z \in X_{11}$ such that $f(z) \in B(y + \rho v/2, \varepsilon \rho)$.
By (15.2), $z\in B(x,2\gamma^{-1}\rho)$ for some $x\in Z(y)$.
By (11.40), $|f(z)-A_x(z)| \leq \varepsilon |z-x| \leq 
2\gamma^{-1}\varepsilon \rho$ if $\rho$  is small enough 
(don't worry, $Z(y)$ is finite). Also 
$\dist(z,P_x) \leq \varepsilon \rho$ for $\rho$ small, because $P_x$
is tangent to $E$ at $x$. Let $\wt z$ denote the projection of
$z$ on $P_x$; then
$$\leqalignno{
\dist(y + \rho v/2,A_x(P_x)) 
&\leq \varepsilon \rho +  \dist(f(z),A_x(P_x))
\leq \varepsilon \rho +  \dist(f(z),A_x(\wt z))
\cr&
\leq \varepsilon \rho + |f(z)-A_x(z)| +|A_x|_{lip} |z-\wt z|
\cr&
\leq \varepsilon \rho + 2\gamma^{-1}\varepsilon \rho + |f|_{lip} \varepsilon \rho
\leq C \gamma^{-1} \varepsilon \rho
&(26.26)
}
$$
by (11.36). For each $\varepsilon > 0$, this holds for $\rho$ small;
since $y \in A_x(P_x)$, it follows that $v$ lies in the vector
space parallel to $A_x(P_x)$, and $P(\Gamma_m) = A_x(P_x)$, as needed.

We add one more constraint to the choice of $r_1(y)$ above: we apply
the definition (25.88)-(25.89) of ${\cal I}_l(U,a,b)$, and require
that 
$$
J_g(B(y,r) \cap \Gamma_m) \leq J_g(B(y,r) \cap P(\Gamma_m)) + \eta r^n
= J_g(B(y,r) \cap A_x(P_x)) + \eta r^n
\leqno (26.27)
$$
for $0 < r \leq r_1(y)$; then by (26.25)
$$
J_g(B(y,r) \cap \lambda^{-1}\psi^{-1}(\wt Q))
\leq J_g(B(y,r) \cap A_x(P_x)) + (1+b) \eta r^d.
\leqno (26.28)
$$
We like this because if we ever
pick $y_j = y$ and $r_j \leq r_1(y)$ for some $j\in J_3$, 
we will immediately deduce from (26.28) that
$$
J_g(D_j \cap \lambda^{-1}\psi^{-1}(\wt Q_j))
\leq J_g(Q_j \cap D_j) + (1+b) \eta r_j^d,
\leqno (26.29)
$$
just because $D_j = B(y_j,r_j)$, $Q_j$ is the common value of
the $A_x(P_x)$, $x\in Z(y_j)$, and $\wt Q_j$ is the common value of
the $\wt A_x(P_x)$. 

This takes care of the definition of $r_1(y)$ in our first case
when ${\rm dimension}(F(y)) = d$. Notice that since we removed
$Z_2$, ${\rm dimension}(F(y)) < d$ is impossible. 
We are left with the case when ${\rm dimension}(F(y)) > d$.
If $y$ lies in no set $L'_i = L_i \sm {\rm int}(L_i)$,
we won't need $r_1(y)$, and we can set $r_1(y) = +\infty$. 
Finally, if $y \in L'_i$ for some $i \leq j_{max}$,
the fact that we removed the set $Z_1$ implies that
the restriction of $\psi$ to $\lambda F \cap B(\lambda y,t(y))$
is of class $C^1$.

Since $y \in Y_{11}$, and as in our first case,
Lemma 12.40 says that for $x\in Z(y)$ the
$d$-plane $\wt A_x(P_x)$ is contained in $W(\wt y)$, 
where $W(\wt y)$ is still the smallest affine subspace $W(\wt y)$
that contains $\wt F$ (and $\wt F$ is the smallest rigid dyadic face that
contains $\wt y = \wt f(x)$); but the difference is that now the dimension
of $W(\wt y)$ is larger than $d$. 
However, there is a neighborhood of $\wt y$ in $W(\wt y)$ 
where the restriction of $\lambda^{-1} \psi^{-1}$ is of class $C^1$
(in fact, the $C^1$-regularity of this inverse map is
the best definition of the $C^1$-regularity of $\psi$
on $\lambda F$). Then, if we set 
$\Gamma(y) = \lambda^{-1}\psi^{-1}(\wt A_x(P_x))$,
there is a neighborhood of $y$ in $U$ where 
$\Gamma(y)$ is is a $C^1$ submanifold of $U$.
By (25.88) again,
$$
\limsup_{r \to 0} r^{-d} \big[ J_g(B(y,r) \cap \Gamma(y))
- J_g(B(y,r) \cap P(y)) \big] \leq 0,
\leqno (26.30)
$$
where $P(y)$ denotes the tangent to $\Gamma(y)$ at $y$.
Then we need to check that 
$$
P(y) = A_x(P_x) \hbox{ for } x\in Z(y).
\leqno (26.31)
$$
The fact that all the sets $A_x(P_x)$ coincide comes from (15.9),
and for the equality with $P(y)$ we shall be able to compute. Let
$R$ denote the differential of $\lambda^{-1} \psi^{-1}$
at $\wt y$; this map is only defined on the vector space 
parallel to $\wt F$, but this will be enough. 
Also denote by $P'$ the vector space
parallel to $P_x$; we know from (11.40) that the restriction
of $DA_x$ to $P'$ is the differential of the restriction
of $f$ to $P_x$. Similarly, (12.39) says that
the restriction of $D\wt A_x$ to $P'$ is the differential 
of the restriction of $\wt f$ to $P_x$. 
We have seen that, because of Lemma 12.40,
$\wt A_x(P_x)$ is contained in the vector space parallel
to $\wt F$. Then the composition 
$R \circ D\wt A_x : P' \to \R^n$ makes sense, and we claim
that it is also the differential of the mapping $f : P_x \to \R^n$.
Indeed, for $v\in P'$, set $z = \wt f(x+tv)$ and denote by
$w$ the projection of $\wt f(x+tv)$ on $\wt A_x(P_x)$.
Then $z = \wt f(x+tv) = \wt y + t  D\wt A_x(v) + o(t)$, 
so $w = \wt y + t  D\wt A_x(v) + o(t)$ too (because 
$\wt y + t  D\wt A_x(v) \in \wt A_x(P_x)$), and finally
$$\eqalign{
f(x+tv) &= \lambda^{-1} \psi^{-1}(z)
= \lambda^{-1} \psi^{-1}(w) + O(|z-w|)
= \lambda^{-1} \psi^{-1}(w) + o(t)
\cr&
= \lambda^{-1} \psi^{-1}(\wt y) + R(w - \wt y) + o(|w - \wt y|) + o(t)
\cr&
= y + t R \circ D\wt A_x(v) + o(t)
}\leqno (26.32)
$$
because $\psi$ is Lipschitz, and as needed. Since we also have the
differential $DA_x$, we see that $DA_x = R \circ D\wt A_x$ on $P'$
and the direction of $A_x(P_x)$ is $R \circ D\wt A_x(P')$, which is 
indeed the direction of $P(y)$ (naturally obtained as the image
by $R$ of the direction of the tangent plane to $\wt A_x(P_x)$).
This proves (26.31). 

We use (26.30) and (26.31) to choose $r_1(y)$ so small 
that for $0 < r \leq r_1(y)$,
$$
J_g(B(y,r) \cap \Gamma(y))
\leq J_g(B(y,r) \cap A_x(P_x)) + \eta r^d
\leqno (26.33)
$$
for $x\in Z(y)$. This way, if we ever 
pick $y_j = y$ and $r_j \leq r_1(y)$ for some $j\in J_3$, 
we will automatically get that (for $x\in Z(y_j)$)
$$\eqalign{
J_g(D_j \cap \lambda^{-1}\psi^{-1}(\wt Q_j))
&= J_g(B(y_j,r_j) \cap \Gamma(y_j))
\cr&
\leq J_g(B(y_j,r_j) \cap A_x(P_x)) + \eta r_j^d
= J_g(Q_j \cap D_j) + \eta r_j^d,
}\leqno (26.34)
$$
which is as good as (26.29).

We may now continue the construction as suggested in Section 19;
recall that we set $Y_{12} = Y_{11} \sm Z$, and 
$\H^d(Y_{11}\sm Y_{12}) = 0$; then we set 
$X_{12} = X_{11} \cap f^{-1}(Y_{12})$, and use 
(4.77) in [D2] to get that  
$$
\H^d(X_{11}\sm X_{12}) = 0.
\leqno (26.35)
$$
Then we choose $\delta_9 > 0$, and set
$$
Y_{13} = Y_{13}(\delta_9) = \big\{ y\in Y_{11} \, ; \, 
r_1(y) < \delta_9 \big\}
\ \hbox{ and }
X_{13} = X_{13}(\delta_9) = X_{12} \cap f^{-1}(Y_{13});
\leqno (26.36)
$$
since the decreasing intersection of the $X_{13}(\delta_9)$ is
$X_{12}$, we can choose $\delta_9$ so small that
$$
\H^d(X_{11}\sm X_{13}) = \H^d(X_{12}\sm X_{13}) \leq \eta/2.
\leqno (26.37)
$$
Then we proceed as before, choose the $D_j$ as we did near (15.12)
but with the stronger constraint that (instead of (15.12))
$$
r_j < {\rm min}(\delta_8,\delta_9)
\ \hbox{ for } j\in J_3.
\leqno (26.38)
$$
We continue our construction as before, except that
we also define the modification $h_3$, which concerns the indices
$j \in J_4$, as described near (26.20).
This way, we only have to estimate the numbers 
$$
\Delta_j = J_g(D_j \cap \lambda^{-1}\psi^{-1}(\wt Q_j))
- J_g(Q_j \cap D_j)
\leqno (26.39)
$$
for $j\in J_3 \sm J_4$. For such $j$, we have defined 
a radius $r_1(y_j) > 0$, and made sure that since $y_j \in Y_{13}$,
$0 < r_j < r(y_j)$. But then $\Delta_j \leq (1+b) \eta r_j^d$,
by (26.29) or (26.34). We sum this over $j\in J_3 \sm J_4$,
and get an additional error term which is dominated by
$\sum_{j\in J_3 \sm J_4} \Delta_j 
\leq (1+b) \eta \sum_{j\in J_3 \sm J_4} \, r_j^d
\leq C(f) \eta$ by (26.17). This is small enough for us to complete 
our proof of the extension of Claim 26.4.
\qed

\ms
We may also generalize many results of Part V
to $f$-almost minimal sets, where $f$ is an elliptic integrand.
Most of the time, the proof is the same once we have
Theorem 25.7 and Claim 26.4, but we prefer to omit the details.

\msi
{\bf 27. Smooth competitors.} 
\ms

In our definition of quasiminimal sets, we used competitors 
for $E$ that come from one-parameter families $\{\varphi_t\}$,
$0 \leq t \leq 1$, for which the final mapping $\varphi_1$
is required to be Lipschitz, as in (1.8). 
We added this requirement because Almgren did it, and because 
this would not disturb in the proofs. The main advantage of this 
definition is probably that it makes it possible to show that
some types of minimal currents (typically, size minimizers) have 
supports that are almost minimal sets. The author suspects this
has been known for ages by specialists, but wrote a short proof
for this in Section 7 of [D8] anyway. 
For other classes, such as Reifenberg homological solutions
of Plateau's problem, we would not need (1.8).

To make the verification of quasiminimality easier for some other 
classes of sets, we may want to restrict the class of 
one-parameter families $\{\varphi_t\}$ (but without changing the 
main defining inequality (2.5)), typically by requiring the
final mapping $\varphi_1$ (or maybe even the whole family of
mappings $\varphi_t$) to be smoother. The issue appeared with some of 
the classes of differential chains introduced by J. Harrison, and
at some point we even produced, with J. Harrison and H. Pugh, a sketch
of proof for the some of the results in the present section
(in the special case with no boundary). The details were never written 
down, essentially because Harrison and Pugh managed to verify the
almost minimality of their supports in a different way.

The author thinks this is a reasonably interesting issue to mention,
especially because we did not find a trivial way to deal with it 
directly with density arguments, hence the present section. He wishes 
to thank J. Harrison and H. Pugh for discussions about this issue and
letting him write down this section. 

We start with some definitions.
Since all our sets (like $E$) may be thin, and we don't want 
to worry about about Whitney jets, let us agree that a function $f$ 
is of class $C^\alpha$ on the set $F$ when $f$ has a $C^\alpha$ 
extension to a neighborhood of $F$.

Let us define modified classes of quasiminimal sets. 
We keep most of the notation in Definition 2.3 as it was, and
say that the closed set $E \i \Omega$ is quasiminimal for competitors
of class $C^\alpha$, with $\alpha \in \{ 1, 2, \ldots, \infty \}$, if 
(2.5) holds for every one-parameter family 
$\{\varphi_t\}$, $0 \leq t \leq 1$,
which satisfy (1.4)-(1.8) and (2.4), and for which, in addition,
$\varphi_1$ is of class $C^\alpha$ on $E$. 
The corresponding classes are denoted
by $GSAQ(U,M,\delta,h,\varphi_1\in C^\alpha)$.

We shall also discuss the following intermediate notion of quasiminimal 
set for piecewise $C^\alpha$ competitors, denoted by
$GSAQ(U,M,\delta,h,\varphi_1\in PC^\alpha)$, where we only require
the competitor to piecewise $C^\alpha$. This last means that
the closure of $\{ x \in E \, ; \, \varphi_1(x) \neq x \}$,
which by (2.4) is required to be a compact subset of $U$, can be 
covered by a finite number of compact sets $K_l$, and 
$\varphi_1(x)$ is $C^\alpha$ on each $K_l$, with the definition
above. Thus we do not care whether the various pieces $K_l$
are smooth or not.

Let us mention yet another variant of these definitions. We say that 
$E$ is  quasiminimal for families of class $C^\alpha$, and we write
$E \in GSAQ(U,M,\delta,h,\varphi_t\in C^\alpha)$, if we only require
(2.5) for families $\{ \varphi_t \}$ that satisfy (1.4)-(1.8) and (2.4), 
and in addition define a $C^\alpha$ function on $V \times [0,1]$,
where $V$ is some neighborhood of the set $E$.
We define $GSAQ(U,M,\delta,h,\varphi_t\in PC^\alpha)$ likewise,
with piecewise $C^{\alpha}$ functions defined on $E\times [0,1]$

\ms
We shall not worry too much about the difference between the
$GSAQ(U,M,\delta,h,\varphi_t\in PC^\alpha)$ and the corresponding
$GSAQ(U,M,\delta,h,\varphi_1\in PC^\alpha)$, or directly
$GSAQ(U,M,\delta,h,\varphi_t\in C^\alpha)$ and 
$GSAQ(U,M,\delta,h,\varphi_1\in C^\alpha)$.
We shall say a few words about this in Remark 27.47 though.
We mention $GSAQ(U,M,\delta,h,\varphi_t\in PC^\alpha)$ now because 
our first result applies easily to that class with no special effort.

The main result of this section is that if our boundary pieces $L_j$ are smooth 
enough, the classes $GSAQ(U,M,\delta,h)$ and $GSAQ(U,M,\delta,h,\varphi_t\in C^1)$
are the same. See the remarks at the end of the section concerning possible further
results, in particular concerning the case of $\alpha > 1$.

\ms
We start our discussion with a basic regularity result for the
class $GSAQ(U,M,\delta,h,\varphi_t\in PC^\alpha)$. We mention
it now because it seems hard to compare our different classes before
we get the rectifiability of $E$. 
For this first result, we do not try to compare directly
$GSAQ(U,M,\delta,h,\varphi_t\in PC^\alpha)$ with $GSAQ(U,M,\delta,h)$,
but just observe that our initial proofs go through.
We start with the Lipschitz assumption.

\ms\proclaim Proposition 27.1.
For each $M \geq 1$ we can find $h>0$ and $C_M \geq 1$,
depending on the dimensions $n$ and $d$, such that if 
$E \in GSAQ(B_0,M,\delta,h,\varphi_t\in PC^\infty)$,
where $B_0 = B(0,1) \i \R^n$, and if the rigid assumption
holds, then $E$ is rectifiable and $E^\ast$ is locally 
$C_M$-Ahlfors regular in $B_0$.
\ms

The local Ahlfors regularity condition means, 
as for Proposition 4.1, that 
$$
C_M^{-1} r^d \leq \H^d(E\cap B(x,r)) \leq C_M r^d
\leqno (27.2)
$$
when $x \in E^\ast$ and $0 < r < \Min(r_0,\delta)$
are such that $B(x,2r) \i B_0$; we decided not to check that
Proposition 3.3, which says that the closed support $E^\ast$ 
of $\H^d_{\vert E}$ is also
a quasiminimal set, also holds in 
$GSAQ(B_0,M,\delta,h,\varphi_t\in PC^\infty)$,
and this is why we have to deal with $E^\ast$.

Our proof will consist in checking that modulo a few minor
modifications, all the competitors that we build to prove
Proposition 4.1 and Theorem 5.16, under the rigid assumption,
are obtained with piecewise $C^\infty$ functions. And indeed, these
mappings are compositions of Federer-Fleming projections that
can be described as follows. We fix a face $F$ of a dyadic cube,
of some dimension $m \in [d,n]$, a point $\xi$ in the interior
of $F$, and which lies outside of the image of the previous
mapping, and then compose with the mapping $p_\xi$
which sends a point $y\in F \sm\{ \xi \}$ to its radial
projection on $\d F$ (centered at $\xi$). On the other faces of the
same dimension, we set $p_\xi(y) = y$ (but then we compose with 
mappings coming from other faces). On the face $F$, we don't need
to know $\pi_\xi$ near $\xi$, which is not in the current image
of $E$, and away from $\xi$, $F$ is decomposed into a finite
collection of closed pieces $F^l$, the inverse images of the
faces of dimension $m-1$ that compose $\d F$, where after a change of
coordinates (so that $\xi = 0$ and the face of $\d F$ is contained 
in the $(m-1)$-plane with equations $x_1 = a$, and
$x_{m+1} = \ldots = x_n = 0$, the mapping $\pi_\xi$ is just given 
analytically by $\pi_{\xi,1}(y) = a$, 
$\pi_{\xi,m+1}(y) = \ldots = \pi_{\xi,n}(y) = 0$,
and $\pi_{\xi,j}(y) =  {a y_j \over y_1}$ for $2 \leq j \leq m$.
So it is easy to extend our definition of $\pi_\xi$ so that it is
defined and Lipschitz on $\R^n$, and piecewise $C^\infty$.
As before, the values of other $\pi_j$ outside of the faces don't 
matter. Recall also that for this result, we were not disturbed by
the boundary condition (1.7), because we chose to project on 
cubes parallel to our grid, so that the Federer-Fleming projections,
which preserve the faces, automatically preserve the $L_j$.
\qed

\ms
The special case of Proposition 27.1 where we use the 
stronger assumption that
$E \in GSAQ(B_0,M,\delta,h,\varphi_t\in PC^\infty)$
easily extends to the case when the Lipschitz assumption holds, 
but the bilipschitz function $\psi : \lambda U \to B(0,1)$ 
of Definition 2.7 is $C^\alpha$; we just conjugate the 
projections onto faces of dyadic cubes with $\psi$, and get 
projections on faces in $U$ that we can use as in the argument above.
It seems hard to adapt the argument to make it work when $\psi$
is merely bilipschitz, and also the author is not sure that the
notion of piecewise smooth competitors is interesting then.

\ms
Our next result will say that if the boundary pieces $L_j$ are
sufficiently smooth, the adverb ``piecewise" in the definition does not add 
anything. It will rely on the following simple extension lemma.

\ms\proclaim Lemma 27.3.
Let $F_1, F_2, \ldots, F_m$ be closed subsets of $\R^n$,
with $\dist(F_i, F_j) > 0$ for $i \neq j$, set $F = \cup_i F_i$,
and let $f: F \to \R$ be a Lipschitz mapping. Suppose 
that for $1 \leq i \leq m$, the restriction of $f$ to
$F_i$ has an extension $g_i$ to an open neighborhood
$V_i$ of $F_i$, which is both $C^{\alpha}$ and Lipschitz.
Then there $f$ has an extension to $\R^n$,
which is of class $C^{\alpha}$, and which is Lipschitz with
$$
|f|_{Lip(\R^n)} \leq C |f|_{Lip(F)} + C \sum_{i} |g_i|_{Lip(V_i)}.
\leqno (27.4)
$$

\ms
We had to mention that $g_i$ is Lipschitz on $V_i$, because
when $F_i$ is not compact this may not follow from the fact
that it is $C^\alpha$.

Let $\varepsilon > 0$ be so small that
the sets $W_i = \{ y \in \R^n \, ; \, \dist(y, F_i) \leq 3\varepsilon 
\big\}$ are disjoint and contained in the corresponding $V_i$.
Denote by $h$ the usual Lipschitz extension of $f$ to $\R^n$,
obtained from the values of $f$ on $F$ with Whitney cubes, 
as in the first pages of [St]. 
Thus $h$ is $C |f|_{Lip(F)}$-Lipschitz, and it is also 
$C^\infty$ on $\R^n \sm F$. 

For each $i$, let $\xi_i$ denote a smooth function such that
$\xi_i(y) = 1$ when $\dist(y,F_i) \leq \varepsilon$,
$\xi_i(y) = 0$ when $\dist(y,F_i) \geq 2\varepsilon$,
and $0 \leq \xi_i(y) \leq 1$ everywhere. We can choose 
$\xi_i$ so that it is $2 \varepsilon^{-1}$-Lipschitz (and we
leave the verification as an exercise).
Also set $\xi_\infty = 1 - \sum_{i} \xi_i$; notice that
$0 \leq \xi_\infty \leq 1$ by definition of $\varepsilon$.
Finally we set
$$
f = \xi_\infty h + \sum_i \xi_i g_i;
\leqno (27.5)
$$
it is clear that it is as smooth as the $g_i$, 
so we just need to check the Lipschitz bound, and since 
$f$ is smooth we just need to bound $Df$. By (27.5),
$Df = D\xi_\infty h + \sum_i D\xi_i g_i + \xi_\infty Dh
+ \sum_i \xi_i Dg_i$. The last two terms are bounded,
so we are left with
$A = D\xi_\infty h + \sum_i D\xi_i g_i$. 
Let $y\in \R^n$ be given. Notice that
$D\xi_\infty(y) + \sum_i D\xi_i(y) = 0$ 
because $\xi_\infty + \sum \xi_i$ is constant. We may
assume that $D\xi_i(y) \neq 0$ for some $i$,
because otherwise $A(y)=0$. Then $\dist(y,F_i) \leq 2 \varepsilon$,
which implies that $D\xi_j(y) = 0$ for $j \neq i$, and
$A(y) = D\xi_j(y) (g_i(y)-h(y))$. Since $g_i$ and $h$ are 
Lipschitz and coincide on $F_i$, we get that 
$|A(y)| \leq |D\xi_j(y)| (|h|_{Lip(F)} + |g_i|_{Lip(V_i)})
\dist(y,F_i) \leq C (|f|_{Lip(F)} + |g_i|_{Lip(V_i)})$,
and the lemma follows.
\qed

\ms
For the next result, we allow the Lipschitz assumption,
but require the boundary pieces $L_j$ to be sufficiently 
smooth and transverse. Of course, when we work with no boundary
pieces (or just the unique $L_0 = \R^n$), thing are much simpler
and we don't need the assumptions below, whose main goal is to
allow a re-projection on the face when we leave them.

Rather than giving simple natural conditions that works well
(the author tried to do this and did not manage), let us say what we 
will use; hopefully our assumptions will not be too brutal and will be easy to 
check in potential applications. We want local retractions on the
$L_j$, which work for all the $L_j$ at the same time. More precisely,
we shall assume that for each compact set $K \i U$, we can find
constants $\tau_0 > 0$ and $C_0 \geq 1$, so that the following holds.
For $0 < \tau \leq \tau_0$, we can find a $C^{\alpha}$ mapping 
$\pi_\tau$, defined on 
$$
K^{\tau} = \big\{ x\in U \, ; \, \dist(x,K) < \tau \big\},
\leqno (27.6)
$$
such that 
$$
|\pi_\tau(x)-x| \leq C_0 \tau 
\ \hbox{ for } x\in K^{\tau},
\leqno (27.7)
$$
$$
|\pi_\tau(x)-\pi_\tau(y)| \leq C_0 |x-y| 
\ \hbox{ for } x,y\in K^{\tau},
\leqno (27.8)
$$
$$
\pi_\tau(x)  \in L_j
\ \hbox{ for $0 \leq j \leq j_{max}$ and $x\in K^{\tau}$
such that $\dist(x,L_j) \leq C_0^{-1} \tau$,}
\leqno (27.9)
$$
but also, setting
$$
Z_j(\rho) = \big\{ x\in U \, ; \, 0 < \dist(x,L_j) < \rho \big\}
\leqno (27.10)
$$
for $0 \leq j \leq j_{max}$ and $\rho > 0$ and
$$
Z(\rho) = \bigcup_{0 \leq j \leq j_{max}} Z_j(\rho),
\leqno (27.11)
$$
such that
$$
\pi_\tau(x) = x 
\ \hbox{ for } x \in K^{\tau} \sm Z(\tau).
\leqno (27.12)
$$
Finally, we require that $\pi_\tau$ is the endpoint of a one parameter
family $\{ \pi_{\tau,t} \}$, such that $\pi_{\tau,t}(x)$ is a 
function of $x\in K^\tau$ and $t\in [0,1]$ which is both $C^\alpha$
(with no precise bound needed), but also $C_0$-Lipschitz
(we shall use this near (27.43))
and such that $\pi_{\tau,0}(x) = x$ and $\pi_{\tau,1}(x) = \pi_\tau(x)$
for $x\in K^{\tau}$, 
$$
|\pi_{\tau,t}(x)-x| \leq C_0 \tau 
\ \hbox{ for $x\in K^{\tau}$ and } 0\leq t \leq 1,
\leqno (27.13)  
$$
and
$$
\pi_{\tau,t}(x) \in L_j
\ \hbox{ for $0 \leq t \leq 1$ when } x\in L_j.
\leqno (27.14)
$$
This looks like a long list, but the reader may check that is easy
to construct such retraction, say, when the $L_j$ are two transverse 
smooth submanifolds (first project on the first one parallel the second 
one, and continue with a projection on the second one along the first 
one), or when the $L_j$ are contained in each other (retract on the
largest, then on the second largest inside the first one, etc.).

\ms\proclaim Proposition 27.15.
If the $L_j$ satisfy the assumption above, then
the two classes \break
$GSAQ(B_0,M,\delta,h,\varphi_1\in PC^\alpha)$
and $GSAQ(B_0,M,\delta,h,\varphi_1\in C^\alpha)$ are equal.

\ms
The following example shows that this result will at least be 
harder to prove if we do not assume that the $L_j$ are smooth. 
Consider, in the unit disk $D \i \R^2 \simeq \Bbb C$, 
a single boundary $L = [0,1) \cup [0,i)$ (two orthogonal intervals), 
and the set $E = L \cup J$, where $J = D \cap [0,1+i]$ is a piece 
of the first diagonal. It is easy to produce better Lipschitz, 
or piecewise $C^1$ competitors, by replacing $J$ with a shorter 
curve $\Gamma$ that ends 
somewhere else on the positive first axis, for instance
(push part of the first quadrant down and to the left). 
The obvious map $\varphi_1$ that does this (i.e., maps
$L$ to itself and $J$ to $\Gamma$) is not $C^1$, and it looks like 
there is an obstruction because we changed the angles at the origin. 
But this is not a counterexample, because we can find a smoother 
mapping, with a vanishing derivative at the origin, and which does the job
even though it destroys some angles. For instance, precompose
the function $\varphi_1$ above with the mapping $x \to |x|^2 x$.

We shall not try to extend Proposition 27.15 to such situations; this 
may be hard, and the benefit is not clear, because the class
$GSAQ(B_0,M,\delta,h,\varphi_1\in C^\alpha)$ is not too natural in
that case. Similarly, our assumptions are probably much too strong, 
but we prefer the proof to be short.

\ms
Let us prove the proposition.
We only need to show that $GSAQ(B_0,M,\delta,h,\varphi_1\in C^\alpha)
\i GSAQ(B_0,M,\delta,h,\varphi_1\in PC^\alpha)$, since the other 
inclusion is trivial. Thus we are given 
$E \in GSAQ(B_0,M,\delta,h,\varphi_1\in C^\alpha)$
and a family $\{ \varphi_t \}$ for which  $\varphi_1$ is piecewise 
$C^\alpha$, and we want to construct a modified family with a final 
map of class $C^{\alpha}$, apply the definition of 
$GSAQ(B_0,M,\delta,h,\varphi_1\in C^\alpha)$, and get (2.5) for the 
initial $\varphi_1$. 

Set $W_1 = \big\{ x \in E \, ; \, \varphi_1(x) \neq x \big\}$
and 
$$
K = \overline W_1 \cup \overline {h_1(W_1)};
\leqno (27.16)
$$ 
by (2.4), $K$ is a relatively compact subset of $U$.
We use this $K$ to apply our assumption on the $L_j$, 
with a very small constant $\tau$ that will be chosen later; 
we get mappings $\pi_\tau$ and $\pi_{\tau,s}$, $0 \leq s \leq 1$,
defined on $K^\tau$.
Since we do not want to modify the $\varphi_t$ too far from
$W_1$, we shall use a smooth cut-off function $\chi$ such that
$0 \leq \chi(x) \leq 1$ everywhere,
$$
\chi(x) = 1 \hbox{ when } \dist(x,W_1) \leq \tau/4, 
\ \chi(x) = 0 \hbox{ when } \dist(x,W_1) \geq \tau/2,
\leqno (27.17)
$$
and $|\nabla \chi| \leq C\tau^{-1}$ everywhere. 

We continue our family $\{ \varphi_t \}$ a first time. Set
$$
\varphi_t(x) = \pi_{\tau,(t-1) \chi(x)}(\varphi_1(x)) 
\ \hbox{ for $x\in E \cap K^\tau$ and } 1 \leq t \leq 2;
\leqno (27.18)
$$
notice that when $x\in W_1$, $\varphi_1(x) \in K$
and $\pi_{\tau,(t-1) \chi(x)}(\varphi_1(x))$ is well defined.
When $x\in E \cap K^\tau\sm W_1$, $\varphi_1(x) = x$ and 
$\pi_{\tau,(t-1) \chi(x)}(\varphi_1(x))$ is well defined too.
When in addition $\dist(x,W_1) \geq \tau/2$, $\chi(x) = 0$
and so $\varphi_t(x) = \varphi_1(x) = x$ (see above (27.13). 
Thus we can safely set
$$
\varphi_t(x) = x
\ \hbox{ for $x\in E \sm K^\tau$ and } 1 \leq t \leq 2,
\leqno (27.19)
$$
and $\varphi_t(x)$ is a continuous function of $x$ and $t$.
Because of what we just said, we even have that
$$
\varphi_t(x) = x
\ \hbox{ for } 1 \leq t \leq 2
\hbox{ when $x\in E$ is such that $\dist(x,W_1) \geq \tau/2$}.
\leqno (27.20)
$$
We want to continue with mappings $\varphi_t$, $2 \leq t \leq 3$,
so that the final mapping $\varphi_3$ is smooth.
Recall that $\varphi_1$ is piecewise smooth; we shall single out one piece,
$F_0 = E \sm W_1$, on which we know that $\varphi_1$
is smooth because $\varphi_1(x) = x$ there. Then, since 
$\varphi_1$ is piecewise smooth, we can cover $\overline W_1$
with a finite collection of compact sets $H_l$, $1 \leq l \leq m$, such that
$\varphi_1$ is $C^\alpha$ on some open neighborhood of $H_l$. 
Of course we may assume that $H_l \i E$.

We want to replace the $H_l$, $1 \leq l \leq m$, with slightly
smaller compact sets $F_l \i H_l$, so that
$$
\hbox{ the $F_l$, $0 \leq l \leq m$, are disjoint,}
\leqno (27.21)
$$
and, if we set
$$
F = \bigcup_{0\leq l \leq m} F_l,
\leqno (27.22)
$$
such that
$$
\H^d(E \sm F) \leq \eta
\leqno (27.23)
$$
and
$$
\dist(x,F) \leq \eta \ \hbox{ for } x\in E,
\leqno (27.24)
$$
where $\eta > 0$ is a very small constant that will be chosen later.

This is easy: we choose the $F_l$ one by one; if the $F_k$, $k<l$
have been chosen, we try $F_l = \big\{ x\in H_l \, ; \,
\dist(x,F_k) \geq a_l \hbox{ for } 0 \leq k < l \big\}$, where $a_l > 0$
will be chosen soon.
Set $H'_l = H_l \sm \big(\bigcup_{0 \leq k <l} F_k\big)$,
and observe that
$$
E \sm F \i \bigcup_{l \geq1} (H'_l \sm F_l).
\leqno (27.25)
$$
Also, for each $l$, $H'_l\sm F_l$ decreases to the empty set when $a_l$
tends to $0$, so $\H^d(H'_l \sm F_l) \leq \eta/m$
if $a_l $ is chosen small enough. Then (27.23) follows from
(27.25). In addition, if $a_l < \eta$ for $l \geq 1$ and $x\in E \sm F$,
then by (27.25) $x\in H'_l \sm F_l$ for some $l \geq 1$, and this forces
$\dist(x,F_k) \leq a_l \leq \eta$, as needed for (27.24).

Apply Lemma 27.3 to the function $\varphi_1$
and the disjoint sets $F_k$.
We get a smooth extension of the restriction of $\varphi_1$ 
to $F$, which we call $f$. Thus
$$
f(x) = \varphi_1(x) \ \hbox{ for } x\in F.
\leqno (27.26)
$$
Notice that in $f$ is Lipschitz with a norm that does not 
depend on $\eta$ or $\tau$ (because for $l \geq 1$, the
$C^\alpha$ extension of $\varphi_1$ in a neighborhood of $H_l$
that we used can be assumed to be Lipschitz too, since $H_l$
is compact). Of course $f$  may differ from $\varphi_1$ on the 
very small set $E \sm F$, but nonetheless
$$
|f(x) - \varphi_1(x)| \leq C |f|_{lip} \dist(x,F) \leq C |f|_{lip} \eta
\leq {\tau \over C_0 +2}
\leqno (27.27)
$$
for $x\in E$, by (27.24), our Lipschitz control on $f$, and if
$\eta$ is small enough compared to $\tau$.

We go from $\varphi_1$ to $f$ by the usual linear 
interpolation, i.e., set
$$
z(x,t) = (t-2)f(x)+(3-t)\varphi_1(x)
\leqno (27.28)
$$
for $x\in E$ and $2 \leq t \leq 3$, and then compose with $\pi_{\tau\chi(x)}$ 
as we did for $\varphi_2$; that is, we want to set
$$
\varphi_t(x) = \pi_{\tau,\chi(x)}(z(x,t))
\ \hbox{ for $x\in E \cap K^\tau$ and } 2 \leq t \leq 3.
\leqno (27.29)
$$
We just need to check that 
$$
z(x,t)\in K^\tau \ \hbox{ when } x\in E \cap K^\tau
\leqno (27.30)
$$
If $x \in W_1$, then $\varphi_1(x) \in K$
(by  (27.16)), and the result follows because (27.27) says that
$|f(x) - \varphi_1(x)| \leq \tau/2$. 
Otherwise, $x\in F_0 \i F$, so $f(x)=\varphi_1(x) = x \in K^\tau$, 
and the result holds too. So (27.30) holds, and  (27.29) makes sense.

On the rest of $E$, we set, as in (27.18),
$$
\varphi_t(x) = x
\ \hbox{ for $x\in E \sm K^\tau$ and } 2 \leq t \leq 3.
\leqno (27.31)
$$
When $x\in E \cap K^\tau$ but $\dist(x,W_1) \geq \tau/2$,
observe that $x\in F_0 \i F$, hence $f(x) = \varphi_1 = x$,
and since $\chi(x) = 0$ by (27.17), we get that 
$\varphi_t(x) = \pi_{\tau,0}(x) = x$. So
$$
\varphi_t(x) = x
\ \hbox{ for } 2 \leq t \leq 3
\hbox{ when $x\in E$ is such that $\dist(x,W_1) \geq \tau/2$},
\leqno (27.32)
$$
as for (27.20).

This completes our definition of the extended family 
$\{ \varphi_t \}$, $0 \leq t \leq 3$, modulo some choices of constants 
that we still need to make. We want to apply our assumption that 
$E \in GSAQ(B_0,M,\delta,h,\varphi_1\in C^\alpha)$, so let us check
the usual requirements for the $\varphi_{3t}$, $0 \leq t \leq 1$.

The continuity condition (1.4) is satisfied; in particular, for $t=2$
and $x\in  E\cap K^\tau$, (27.29) yields 
$\varphi_t(x) = \pi_{\tau,\chi(x)}(\varphi_1(x))$, 
just like (27.18). Of course $\varphi_0(x) = x$. Also,
$$
\varphi_t(x) = x
\ \hbox{ for } 0 \leq t \leq 3
\ \hbox{ when $x\in E$
is such that $\dist(x,W_1) \geq \tau/2$}, 
\leqno (27.33)
$$
by (2.1), (27.19), and (27.32). If $B = \overline B(x_0,r_0)$
was the ball for which (1.5) and (1.6) hold for the initial $\varphi_t$,
we get (1.5) for the extended family, with any ball $B'$ that contains 
$B(x_0,r_0\tau)$.

Let us now check (1.6), with the ball $B' = \overline B(x_0,r_0+(C_0+2)\tau)$.
We are given $x\in E \cap B'$, and we want to check that $\varphi_t(x) \in B'$
for all $t$. We may assume that $\dist(x,W_1) < \tau/2$, because otherwise
the result follows from (27.33), and also that $t > 1$, because we know 
(1.6) for the  $\varphi_t$, $0 \leq t \leq 1$. Let us check that
$$
|\varphi_t(x) - \varphi_1(x)| \leq (C_0+1) \tau.
\leqno (27.34)
$$
If $1 \leq t \leq 2$, $\varphi_t(x)$ is given by (27.18),
and $|\varphi_t(x) - \varphi_1(x)| \leq C_0 \tau$ by (27.13).
Otherwise, $\varphi_t(x)$ is given by (27.29); thus
$\varphi_t(x) = \pi_{\tau,\chi(x)}(z(x,t))$, where $z(x,t)$
is defined by (27.28) and lies in $[f(x),\varphi_1(x)] \cap K^\tau$ 
by (27.30). Then
$$\eqalign{
|\varphi_t(x) - \varphi_1(x)|
&= |\pi_{\tau,\chi(x)}(z(x,t)) - \varphi_1(x)|
\cr&
\leq |\pi_{\tau,\chi(x)}(z(x,t)) - z(x,t)| + |z(x,t) - \varphi_1(x)|
\cr&
\leq |\pi_{\tau,\chi(x)}(z(x,t)) - z(x,t)| + |f(x) - \varphi_1(x)|
\leq (C_0+1) \tau
}\leqno (27.35)
$$
by (27.13) and (27.27). So (27.34) holds. 
But now 
$$
\dist(\varphi_t(x),W_1) \leq \dist(x,W_1) + (C_0+1) \tau
\leq (C_0+2) \tau,
\leqno (27.36)
$$
which proves that $\varphi_t(x) \in B'$, because 
$W_1 \i B = \overline B(x_0,r_0)$.

The compactness condition (2.4) also holds, because the analogue
of $\wh W$ for the extended family lies in a $(C_0+2) \tau$-neighborhood
of $\wh W$, by (27.33) and (27.36) in particular. Of course this neighborhood
is compactly contained on $U$ if $\tau$ is small enough.

We managed to end our family with a mapping
$\varphi_3$ which is $C^\alpha$.
Indeed, (27.29) yields $\varphi_3(x) = \pi_{\tau,\chi(x)}(f(x))$
for $x\in E \cap K^\tau$, $f$ was constructed to be $C^\alpha$,
$\pi_{\tau,s}(x)$ is a $C^\alpha$ function of $s$ and $x$, and
as usual there is an overlap between the definitions by (27.29)
and (27.31), where both definitions yield $\varphi_3(x)=x$.
This takes care of the improved constraint (1.8) with $C^\alpha$.

Finally we check (1.7). We are given $x\in E \cap L_j$, and we
want to check that $\varphi_t(x) \in L_j$ for all $t$.
We can assume that $t > 1$ (otherwise, use the old (1.7)),
and that  $\dist(x,W_1) < \tau/2$ (by (27.33)).
By the old (1.7), $\varphi_1(x) \in L_j$ and now (27.18) yields
$\varphi_t(x) \in L_j$ for $1 \leq t \leq 2$, by (27.14).
So we assume that $t \geq 2$, and $\varphi_t(x)$ is given by (27.29).
If $x\in F$, then $f(x) = \varphi_1(x)$ and (27.29) yields
$\varphi_t(x) = \varphi_2(x) \in L_j$ for $t \geq 2$, as needed.
So can assume that $x\in E \sm F$. Since $x \notin F_0 = E \sm W_1$,
we get that $x\in W_1$, and then $\chi(x) = 1$
(see (27.17)), so $\varphi_t(x) = \pi_{\tau,1}(z(x,t)) 
= \pi_\tau(z(x,t))$, with $z(x,t) = (t-2)f(x)+(3-t)\varphi_1(x) \in [f(x),\varphi_1(x)]$
(see above (27.13)). By (27.27), 
$$
\dist(z,L_j) \leq |z(x,t)-\varphi_1(x)|
\leq |f(x)-\varphi_1(x)| \leq C_0^{-1} \tau,
\leqno (27.37)
$$
and now (27.9) says that $\varphi_t(x) = \pi_\tau(z(x,t)) \in L_j$.

This completes our list of verifications, and we may now apply the 
quasiminimality of $E$. Set
$W_3 = \big\{ x\in E \, ; \, \varphi_3(x) \neq x \big\}$; 
then the analogue of (2.5) for $\varphi_3$ says that
$$
\H^d(W_3) \leq M \H^d(\varphi_3(W_3)) + h \, r_1^d,
\leqno (27.38)
$$
where $r_1 = r_0+(C_0 +2)\tau$ is the radius of our ball $B'$,
and $r_0$ is the radius of our initial ball $B$.

We want to use $W_1$, so let us estimate the size of
the symmetric difference $W_1 \Delta W_3$. By definition,
$W_1 \Delta W_3 \i \Xi$, where
$$
\Xi = \big\{ x\in E \, ; \, \varphi_3(x) \neq \varphi_1(x)\big\}.
\leqno (27.39)
$$
We claim that
$$
\Xi \i (E \sm F) \cup X(\tau) \cup Y(\tau),
\leqno (27.40)
$$
where 
$$
X(\tau) = \big\{ x\in E \, ; \, \tau/4 \leq \dist(x,W_1) \leq 
\tau/2 \big\}
\leqno (27.41)
$$
and
$$
Y(\tau) = \big\{ x\in E \, ; \, \varphi_1(x) \in Z(\tau)\big\}, 
\leqno (27.42)
$$
where $Z(\tau)$ is defined by (27.10) and (27.11).
Indeed, let $x\in \Xi$ be given. If $x\in E \sm F$ we are happy, so
we may assume that $x\in F$. Then $f(x) = \varphi_1(x)$. Also,
$\dist(x,W_1) < \tau/2$, because otherwise (27.30) and (2.1) say that
$\varphi_3(x)=x =\varphi_1(x)$. Thus (27.26) applies, and
$\varphi_3(x) = \pi_{\tau, \chi(x)}(f(x)) =\pi_{\tau, \chi(x)}(\varphi_1(x))$.
Since $\varphi_3(x) \neq \varphi_1(x)$, we get that $\chi(x) \neq 0$.
If $x\in X(\tau)$, we are happy; otherwise, $\xi(x) = 1$
(see above (27.17)) and $\varphi_3(x) = \pi_{\tau,1}(\varphi_1(x))
= \pi_\tau(\varphi_1(x))$.

If $\varphi_1(x) \in Z(\tau)$, we are happy. 
Otherwise, (27.12) says that $\varphi_3(x) = \pi_\tau(\varphi_1(x))
= \varphi_1(x)$ (recall that we checked that $\varphi_1(x) \in K^\tau$ 
below (27.18)); this contradiction completes the proof of (27.40).

Our function $\varphi_3$ is $C$-Lipschitz, with a (possibly huge)
constant $C$ that does not depend on $\tau$; thus
$$\eqalign{
\H^d(\varphi_3(W_3)) 
&\leq \H^d(\varphi_3(W_3 \sm \Xi)) + C H^d(\Xi)
\cr&
\leq \H^d(\varphi_1(W_3 \sm \Xi)) + C H^d(\Xi)
\cr&
\leq \H^d(\varphi_1(W_1)) + C H^d(\Xi)
}\leqno (27.43)
$$
because $\varphi_3 = \varphi_1$ on $W_3 \sm \Xi$, and
$W_1 \Delta W_3 \i \Xi$. Also,
$\H^d(W_1) \leq \H^d(W_3) + H^d(\Xi)$ because $W_1 \Delta W_3 \i \Xi$, 
so (27.38) yields
$$\eqalign{
\H^d(W_1) &\leq \H^d(W_3) + H^d(\Xi)
\leq M \H^d(\varphi_3(W_3)) + h \, r_1^d + H^d(\Xi)
\cr&
\leq M \H^d(\varphi_1(W_1)) + C (1+M) H^d(\Xi) +  h \, r_1^d.
}\leqno (27.44)
$$
We shall now use (27.40) to estimate $H^d(\Xi)$.
Recall that we may choose $\tau$ as small as we want, and then $\eta$
even smaller. By (27.23), $\H^d(E \sm F) \leq \eta$ can be made as 
small as we want; similarly, $X(\tau) \i X'(\tau) = \big\{ x\in E \, ; \, 
0 < \dist(x,W_1) < \tau/2 \big\}$, and since all the $X'(\tau)$ have a finite 
$\H^d$-measure and their monotone intersection is empty,
$\H^d(X(\tau))$ can be made as small as we want too.
The same argument applies to $Y(\tau)$ (recall from (27.10)
and (27.11) that the monotone limit of $Z(\tau)$ is empty.

Thus (27.44) holds for arbitrarily small values of $\tau$ and hence
$H^d(\Xi)$. In addition, $r_1 = r_0+(C_0 +2)\tau$ is as close to $r_0$
as we want, we get (2.5) for the initial $\varphi_1$, and 
Proposition 27.15 follows.
\qed

\ms
When $\alpha = 1$ we can use the fact that Lipschitz functions are
not far from $C^1$ to obtain easily the main result of this section.

\ms\proclaim Corollary 27.45.
Suppose the $L_j$ satisfy the same assumption as for
Proposition 27.15, with $\alpha = 1$. Then
the two classes $GSAQ(B_0,M,\delta,h)$
and $GSAQ(B_0,M,\delta,h,\varphi_1\in C^1)$ are equal.

\ms
Recall that $GSAQ(B_0,M,\delta,h)$ is the usual class of quasiminimal 
sets that we studied in the rest of this text.

As before, one of the inclusions is trivial, and we just need to check
that if $E \in GSAQ(B_0,M,\delta,h,\varphi_1\in C^1)$,
then $E \in GSAQ(B_0,M,\delta,h)$. 
We first apply Proposition~27.15 and get that
$E \in GSAQ(B_0,M,\delta,h,\varphi_1\in PC^1)$.
This is good, because now we can apply
Proposition 27.1 to show that $E$ is rectifiable.

Let $\{ \varphi_t \}$ satisfy (1.4)-(1.8) and (2.4); we want to
copy the proof of Proposition 27.15, but we need $C^1$ mappings,
so we first pick a compact set $K$ that contains a neighborhood
of $W_1 \cup \varphi_1(W_1)$, and then use the rectifiability of $E$ 
and Theorem 3.2.29 in [Fe]  
or Theorem 15.21 in [Ma] 
to find a countable collection of $C^1$ submanifolds 
$\Gamma_j \i \R^n$, and disjoint Borel sets $F_j \i \Gamma_j$, 
so that $\H^d(E \cap K \sm \bigcup_i F_i) \leq \eta$, where
the very small $\eta >0$ will be chosen at the end of the argument.
In fact, at the price of replacing $\eta$ with $2\eta$, we can
suppose that the family is finite, and that each $F_i$ is compact.
Even more, 
Theorem 3.1.16 in [Fe] 
allows us to (make $F_j$ a tiny bit smaller and) assume that
$\varphi_1$ coincides on $F_i$ with a $C^1$ function on $\Gamma_j$.
Notice that this function can be extended into a $C^1$ function $g_i$
defined on a neighborhood of $\Gamma_i$ (and hence $F_i$).
We may now proceed as before; the only difference is that 
the small neighborhoods where we have $C^1$ extensions only 
cover $F = \cup_j F_j$ (and not $E$), but we did not use
this to apply Lemma 27.3 and define our extension $f$.
So we conclude as in Proposition 27.15.
\qed

\ms
We end this section with further results that may well be true,
but which the author was too lazy to check. 
Thus the point of the following remarks is mostly to record what
the author believes, just for the case when the potential results
may become useful. The situation for 
$GSAQ(B_0,M,\delta,h,\varphi_1\in PC^1)$ is correct,
perhaps modulo our transversality assumption for the $L_j$.
But we may feel bad about the very small difference between
Lipschitz (as in (1.8)) and $C^\alpha$.

\msi{\bf Remark 27.46.}
It is probably true that for $\alpha > 1$, 
the two classes $GSAQ(B_0,M,\delta,h)$ and 
$GSAQ(B_0,M,\delta,h,\varphi_1\in C^\alpha)$ coincide
under the same regularity condition for the $L_j$ as in Proposition 27.15.
But even if there is no boundary piece $L_j$, it seems that some nontrivial
argument is needed.

As before, and because of Proposition 27.15, it is enough to show that every
$E \in GSAQ(B_0,M,\delta,h,\varphi_1\in PC^\alpha)$
lies in $GSAQ(B_0,M,\delta,h)$. 
But this time we really need to modify our family $\{ \varphi_t \}$ on a large set, 
because $\varphi_1$ may not be smooth anywhere.

We encountered this sort of problem before, when we were dealing with
limits; we wanted to construct good competitors for sets $E_k$
that lie close to $E$, and were led to constructing stable competitors first.
Here we probably want to do something similar, and proceed roughly as 
follows. We are given our set $E \in GSAQ(B_0,M,\delta,h,\varphi_1\in PC^\alpha)$,
and a family $\{ \varphi_t \}$ that only satisfies the usual Lipschitz condition
(1.8), and we want to construct a smoother family.
We first build the stabler family of Sections 11-17, and because it is stable,
we should be able to make it smoother without making is much worse.
That is, the places where we expect the largest contributions are the
$B_{j,x}$, $j\in J_3$ and $x\in Z(y_j)$, and on these places
we composed the initial mapping $\varphi_1$ by a projection onto a $d$-plane
and compared the measure of the image with the measure of a disk.
We claim that replacing $\varphi_1$ with a smoother mapping 
before we project will not change the final estimates much.

On the other balls, or the intermediate regions (thin annuli, bad sets), 
we typically used no more than the fact that we project onto planes 
(which we still intend to do after we smooth out $\varphi_1$), 
and that our final mapping is Lipschitz with uniform
bounds (which will not be disturbed by smoothing). 
 
This description is probably enough if there is no boundary
piece $L_J$, but in the general case we would also need to compose with
retractions on the faces, as we did in Part IV and later, and for this 
the assumptions of Proposition 27.15 will probably be needed again. 
The fact that we are allowed competitors which are 
merely piecewise $C^\alpha$ may not be really needed (we can probably
glue our pieces smoothly), but is at least psychologically comforting. 
At this point the reader probably guessed why we do not want to do all this here.

\msi{\bf Remark 27.47.}
Let us say a few words about the difference between the two classes
$GSAQ(U,M,\delta,h,\varphi_t\in C^\alpha)$ and 
$GSAQ(U,M,\delta,h,\varphi_1\in C^\alpha)$.

The author believes that under the assumptions of 
Proposition 27.15, these two classes are probably equal, and that the
same thing holds for their piecewise counterparts
$GSAQ(U,M,\delta,h,\varphi_t\in PC^\alpha)$ 
and $GSAQ(U,M,\delta,h,\varphi_1\in PC^\alpha)$.

In fact, the issue may also arise with our definition of the standard
$GSAQ(U,M,\delta,h)$, where we only require the final mapping $\varphi_1$
to be Lipschitz, and we could have required instead that the whole map
$(x,t) \to \varphi_t(x)$ be Lipschitz. We expect that this yields the same class
of quasiminimal sets, but never checked.

This time, given a family $\{ \varphi_t \}$ such that $\varphi_1$
is (piecewise) smooth, we want to change the $\varphi_t$,
$0 < t < 1$, to make the family (piecewise) smooth.
If there is no boundary piece, the point is merely to find a 
(continuous and piecewise) smooth mapping on $E\times [0,1]$, with the
given boundary value for $t= 0$ and $t=1$. For this, 
the simplest is to use the definition of smoothness for
$\varphi_1$, which gives a smooth extension to an open
neighborhood of $E$, then decide brutally that
$\varphi_t = \varphi_1$ for $1-\varepsilon \leq t \leq 1$
(and similarly for $0 \leq t \leq \varepsilon$), and in the middle use 
partitions of unity and the values of $\varphi_t(x)$ on a discrete,
but rather dense set of $E \times [0,1]$ to interpolate. 
Maybe a small smooth gluing will be needed near $t=\varepsilon$
and $t=1-\varepsilon$ too (as in Lemma 27.3).
We proceeded a little like this in Section 11 when we extended
$f$, and the advantage is that the construction is rather explicit,
and in particular we can get that the new family $\{ \wt \varphi_t \}$
is such that $||\wt \varphi_t - \varphi_t||_{\infty}$ is as small as we want.

When there are boundary pieces $L_j$, the new $\wt\varphi_t$
may not respect the $L_j$, and we need to use the smooth universal
retractions $\pi_{\tau,t}$ of Proposition 27.15 (defined near
(27.6)-(27.14)) to send points back to the $L_j$. 
That is, first observe that we can choose the $\wt \varphi_t$
so that the analogues for them of the sets $W_1$, $\varphi_1(W_1)$, 
and $\wh W$ all stay within $\tau/100$ of the original $W_1$, 
$\varphi_1(W_1)$, and $\wh W$, with $\tau$ as small as we want.
Choose for $K$ the closure of $\wh W$, $\tau$ very small and in particular
such that $K^\tau \i \i U$, and use the assumptions for 
Proposition 27.15 to find the retractions $\pi_{\tau,t}$.
Then set
$$
\varphi_t^\sharp(x) = \pi_{\tau,\chi(x,t)}(\wt \varphi_t(x))
\leqno (27.48)
$$
for $x\in E$ and $0 \leq t \leq 1$, with a function $\chi$
that we still need to define. 

The point of keeping
$\varphi_t = \varphi_1$ for $1-\varepsilon \leq t \leq 1$
and similarly for $0 \leq t \leq \varepsilon$ is that we did not destroy
anything on these intervals, and so we can take
$\chi(x,t) = 0$ for $1 - \varepsilon/2 \leq t \leq 1$ and 
for $0 \leq t \leq \varepsilon/2$. We also want to take
$\chi(x,t) = 0$ when $x\in E \sm K^{\tau/2}$. We take
$\chi(x,t) = 1$ when $x\in E \cap K^{\tau/3}$ and
$\varepsilon \leq t \leq 1-\varepsilon$. We also make $\chi$
smooth, with values in $[0,1]$. We claim (but will not check)
that if $||\wt \varphi_t - \varphi_t||_{\infty} < C_0^{-1} \tau$
for all $t$, (27.48) gives a family $\{ \varphi_t^\sharp \}$
which is smooth, for which (1.7) holds, and which still satisfies (2.2) and 
(1.4)-(1.6), although perhaps with a slightly larger ball $B$ because,
to be safe, we want the new one to contain $K^\tau$.
Thus the only effect in the verification of the quasiminimality
property (2.5) is that, even though we did not change $\varphi_1$
and $W_1$, we have to replace $r^d$ in the right-hand side with a slightly larger
$r_1^d$, which does not harm much. This completes our sketch of a 
potential proof of equivalence between the classes with $\varphi_1 \in C^\alpha$ 
and $\varphi_t \in C^\alpha$.

\bigskip
\centerline{PART VII : MONOTONE DENSITY}
\ms
Monotonicity results for minimal sets or surfaces are very useful,
for instance because they usually give a good control on the
blow-up limits of these objects.

The starting point of this part is the following simple 
result (Theorem 28.4 below). 
Suppose that $E$ is a (locally) minimal set, with boundary 
pieces $L_j$ that are cones centered at $x$; then the density
$\theta(x,r) = r^{-d} \H^d(E \cap B(x,r))$ is a nondecreasing
function of $r$ (small).

We also show (in Theorem 29.1) that when in addition 
($E$ is coral and) $\theta(x,\cdot)$ is constant,
$E$ coincides with a minimal cone centered at $x$.
This, with our result of Section~24, is our way to show that
bow-up  limits of almost minimal sets are minimal cones
(Corollary 29.52 below). 

We shall establish (with essentially the same proof as for Theorem 28.4)
that $\theta(x,\cdot)$ is nearly monotone when $E$ is almost minimal
with a sufficiently small gauge function $h$ (and the $L_j$
are still cones centered at $x$); see Theorem 28.7. 
This result is extended, with a slightly different density function 
to make the computations easier, to the case when the $L_j$ 
are not exactly cones. See Remark 28.11 and Theorem 28.15.

The equality case proved in Section 29 
will allow us to show, by compactness, that (under suitable 
assumptions) if for the almost minimal set $E$,
the density $\theta(x,\cdot)$ is almost constant, then 
$E$ is close to a minimal cone, both in Hausdorff distance and in measure.
See Proposition 30.3 for a statement of approximation in an annulus, and
Proposition~30.19 for a simpler statement of approximation in a ball.
 
The results of this part clearly have some interest, but we should observe 
that it would be much better to have monotonicity results for some
quantity like $\theta(x,r)$, which would also hold when $x$ is not the
center of the $L_j$. We shall not try to prove such formulae here.

\msi
{\bf 28. Monotone density for minimal sets; almost monotone density 
in some cases} 
\ms

We start our study with the monotonicity of density for a minimal set.
We consider a coral minimal set $E$, more precisely such that
$$
E \in GSAQ(U,1,\delta,0)
\leqno (28.1)
$$
for some open set $U$, and where the boundary pieces $L_j$,
$0 \leq j \leq j_{max}$, satisfy the Lipschitz assumption. See
Definitions 2.3 and 2.7. We are also given a ball $B(x_0,r_0) \i U$,
and we assume that $r_0 \leq \delta$, and also that for
$0 \leq j \leq j_{max}$, 
$$
L_j  \hbox{ coincides, in $B(x_0,r_0)$, with a closed
cone centered at $x_0$.}
\leqno (28.2)
$$
We allow $L_j \cap B(x_0,r_0) = \emptyset$ (even though in this case
there is not much point in keeping $L_j$), and we allowed the 
Lipschitz assumption because we do not want to restrict to
plane sectors that make square or flat angles. In fact our 
proof will only use the Lipschitz assumption to make sure that
$E$ is rectifiable, and otherwise (28.2) will be enough.
Next set
$$
\theta(r) = r^{-d} \H^d(E\cap B(x_0,r))
\ \hbox{ for } 0 < r \leq r_0.
\leqno (28.3)
$$

\ms\proclaim Theorem 28.4.
Let $U$, the $L_j$, the minimal set $E$, and
$B(x_0,r_0) \i U$ satisfy the assumptions above. Then
$\theta : (0,r_0) \to \R_+$ is nondecreasing.

\ms
This should not shock the reader. The result for minimal sets
far from the boundary is classical, and relies on comparisons
of $E$ with cones, which can be obtained as limits of radial 
deformations of $E$. These deformations will preserve the boundary
pieces $L_j$, by (28.2), and we will be able to conclude from there.

We shall follow the proof of Proposition 5.16 in [D5], 
which was conveniently done in a similar context. 
We start with the integrated version of monotonicity which is stated
as Lemma 5.1 in [D5].  
In the statement of that lemma, the author required that $E$ be 
coral and that $x_0 \in E$, but this is not used in the proof
(only later in the section). The proof then used the rectifiability 
of $E$ and a radial deformation, defined near (5.3), 
and it is clear that such deformations satisfy our boundary constraint
(1.7) because of (28.2). So Lemma 5.1 in [D5] goes through. 
Once we have that lemma, the proof is a simple manipulation of 
measures and integrals, that does not use the minimality of $E$,
and it goes through as it is. Theorem 28.4 follows.
\qed

\ms
Let us now record a version of Theorem 28.4 for almost minimal sets.
We give ourselves a gauge function $h: (0,+\infty) \to [0,+\infty]$, 
which we assume to be nondecreasing and continuous on the right;
these are probably not exactly needed, but won't disturb much and the 
assumption was made in Section 4 of [D5]  
and the beginning of our Section 20. We also assume the Dini condition
$$
\int_0^{r_0} h(2t) {dt \over t} < +\infty
\leqno (28.5)
$$
(which is really used in the proof; don't mind the fact that we
have $h(2t)$, which we repeat from [D5]  
and is jut due to the fact that we wanted to estimate
$h(t)$ with an integral). We replace our minimality condition (28.1) with
the new one that
$$
E \hbox{ is an $A$-almost, or an $A'$-almost minimal set in $U$,
with gauge function $h$}
\leqno (28.6)
$$
(and with the sliding conditions given by the closed sets $L_j$). See
Definition 20.2; we don't care whether the $A$-almost or $A'$-almost
minimality is used, both are equivalently easy to use in the proof.
We now copy the analogue in the present context of
Proposition 5.24 in [D5].  

\ms\proclaim Theorem 28.7.
There exist constants $\alpha > 1$ and $\varepsilon_n > 0$, 
that depend only on the dimension $n$ and on the constant $\Lambda$
in the Lipschitz assumption, such that the following holds.
Let $U$, the $L_j$, the gauge function $h$, the almost minimal set $E$, 
and the ball $B(x_0,r_0) \i U$ satisfy the assumptions above. Suppose 
in addition that $E$ is coral, that $x_0 \in E$, and that
$h(r_0) \leq \varepsilon_n$. Then
$$
\theta(r) \exp \alpha \Big(\int_0^r h(2t) {dt \over t} \Big)
\ \hbox{ is a nondecreasing function of } r\in (0,r_0).
\leqno (28.8)
$$

\ms
Notice that the exponential tends to $1$ as $r$ tends to $0$,
so we can see it as a nice (increasing) extra term that we multiply
with $\theta(r)$ to get a nondecreasing function. The fact that
we use $h(2t)$ is just an artifact of the statement; the reader should
not worry about the case when $t$ is close to $r_0$ and $h(2t)$
may not be defined naturally: just set $h(r) = h(r_0)$ for $r \geq r_0$.

It looks strange that now we require $E$ to be coral and $x_0 \in E$;
this is because in the proof, we use a lower bound for $\theta(r)$, 
that comes from the local Ahlfors-regularity of $E$, 
to simplify a differential inequality. See Remark 28.9 though.
For the proof we proceed as we did in [D5];  
our almost minimality assumption is only used twice, once in
Lemma 5.1 as before, and once, through the local Ahlfors regularity,
in the computation of differential inequalities; so the proof goes 
through. We have to let $\varepsilon_n$ depend on $\Lambda$ because we
use our regularity theorems to prove that $E$ is rectifiable and 
locally Ahlfors-regular, and $\alpha$ depends on $\Lambda$ too, 
through the local Ahlfors-regularity bounds that we use to modify 
a differential inequality.
\qed

\msi{\bf Remark 28.9.} We may even drop our assumption that 
$E$ is coral and $x_0 \in E$ if we replace (28.6) with the stronger
$$
E \hbox{ is an $A_+$-almost minimal set in $U$, with gauge function $h$.}
\leqno (28.10)
$$
This is Proposition 5.30 in [D5], and as before 
the proof just goes through. Then we don't need the local Ahlfors-regular 
bound and $\alpha$ does not depend on $\Lambda$.

\msi{\bf Remark 28.11.} 
Theorem 28.7 can also be generalized slightly to situations
where the $L_j$ are not exactly cones.

Let us assume that, instead of (28.2), we have a
bilipschitz mapping $\xi = B(x_0,r_0) \to \xi(B(x_0,r_0)) \i \R^n$,
with the following properties. First of all, 
$$
\xi(B(x_0,r_0) \cap L_j) \hbox{ coincides, in $\xi(B(x_0,r_0))$, 
with a closed cone centered at $\xi(x_0)$;}
\leqno (28.12)
$$
this will be our replacement for (28.2). 
We want a better control (typically, of $C^1$ type)
in the smaller balls centered at $x_0$, so we assume that
for $r\in (0,r_0]$, there is a constant $\rho(r))$ such that
$$
\hbox{ the restriction of $\xi$ to $B(x_0,r)$ is 
$(1+\rho(r))$ bilipschitz}
\leqno (28.13)
$$
and
$$
\int_0^{r_0} \rho(t) {dt \over t} < +\infty.
\leqno (28.14)
$$
Then we have the following extension of Theorem 28.7.

\ms\proclaim Theorem 28.15.
There exist constants $\alpha_1 > 1$ and $\varepsilon_n > 0$, 
that depend only on the dimension $n$ and on the constant $\Lambda$
in the Lipschitz assumption, such that the following holds.
Let $U$, the $L_j$, the gauge function $h$, the coral almost minimal 
set $E$, and the ball $B(x_0,r_0) \i U$ be such that 
$h$ is nondecreasing and continuous on the right, 
(28.5) and (28.6) hold, $x_0 \in E$,
$h(r_0) \leq \varepsilon_n$, and there exists $\xi$ and $\rho$
as above, such that $\rho(r_0) \leq \varepsilon_n$ and
(28.12), (28.13), and (28.14) hold. 
Set $\wt B(r) = \xi^{-1}(B(\xi(x_0),r))$ and 
$$
\Psi(r) = \H^d(E\cap \wt B(r))
\exp\Big(\alpha_1 \Big(\int_0^r (h(9t/4)+\rho(9t/4)) {dt \over t} 
\Big)\Big)
\leqno (28.16)
$$
for $0 < r < r_0/2$; then $\Psi$ is nondecreasing
on $(0,r_0/2]$.

\ms
We were a little lazy here, because we measured the density
in terms of the slightly distorted balls $\wt B(r)$.
This way we will be able to reduce to Theorem 28.7
via a change of variable. Probably the more reasonable statement
with the the same function $\theta$ as above also holds, but for
this it seems that we would have to follow the proof above,
and in due time modify the proof of Lemma 5.1 in [D5].  
That is, we would obtain some estimate on the measure of thin annuli
by constructing directly a competitor that expands the annulus
and contracts the inside disk. Since our initial radial competitor
probably does not satisfy (1.7) (because the $L_j$ are no longer
cones), we could conjugate by $\xi$ and apply a radial transformation
in the new variables. We decided to use the function $\Psi$ above
and avoid the computations.

So we try to deduce Theorem 28.15 from Theorem 28.7 and
a change of variable. Without loss of generality, we may assume
that $\xi(x_0) = x_0 = 0$. Set $B_0 = B(x_0,r_0)$ and 
$\wt E = \xi(E\cap B_0)$; we would like to say that $\wt E$ 
is almost minimal in $\xi(B_0)$, but this is probably wrong, because
we did not assume $\xi$ to be asymptotically conformal near
each point of $B_0$ (this would be a very strong assumption
to make!), but only at the point $x_0$. So we need to be a little
careful with our assertions.

Let us first prove that $\wt E$ is quasiminimal in $B_1 = \xi(B_0)$, 
with $M=1$, $\delta=2r_0$, and $h = C\varepsilon_n$. 
That is, with the notation of Definition 2.3, that
$$
\wt E \in GSAQ(B_1,1,2r_0,C\varepsilon_n),
\leqno (28.17)
$$
where on $B_1$, we use the boundary pieces
$\wt L_j = B_1 \cap \xi(L_j)$. We put $\delta=2r_0$ as a way to
imply that we put no constraint on the size of the analogue of $\wh W$
for competitors of $\wt E$.  

The proof will be easy.
Let the $\wt \varphi_t$, $0 \leq t \leq 1$, define a competitor
for $\wt E$ in $B_1$; that is, assume that they satisfy
the analogue of (1.4)-(1.8), in a ball $\wt B$ of radius $\wt r$,
and (2.4). Then define $\varphi_t$, $0 \leq t \leq 1$, 
by $\varphi_t(x) = \xi^{-1}(\wt \varphi_t(\xi(x)))$
for $x\in \xi^{-1}(B_1)$ and $\varphi_t(x)=x$
otherwise. It is easy to see that the $\varphi_t$
define a competitor for $E$ (i.e, satisfy the usual
constraints (1.4)-(1.8) and (2.4)), in a ball $B$
that contains $\xi^{-1}(\wt B)$; we can choose $B$ 
or radius $r = \min(r_0,(1+\rho(r_0)) \wt r)$.
In addition,
$\wt W = \big\{ x\in \wt E \cap B_1 \, ; \,
\wt\varphi_1(x) \neq x \big\}$ is equal to $\xi(W)$,
where $W = \big\{ x\in E \cap B_0 \, ; \, \varphi_1(x) \neq x \big\}$. 
We may assume that $E$ is $A$-almost minimal 
(because $A'$-almost minimality implies $A$-almost minimality 
with the same gauge function $h$; see near (20.8)), 
and the defining property (20.5) yields
$$\eqalign{
\H^d(\wt W) &= \H^d(\xi(W)) \leq (1+\rho(r_0))^d \H^d(W))
\cr&
\leq (1+\rho(r_0))^d \big[\H^d(\varphi_1(W))+ h(r) r^d \big]
\cr&
\leq (1+\rho(r_0))^d 
\big[(1+\rho(r_0))^d \H^d(\wt\varphi_1(\wt W))+ h(r) r^d \big].
}\leqno (28.18)
$$
If $\H^d(\wt W) \leq \H^d(\wt\varphi_1(\wt W))$, 
we are happy. Otherwise,
$$
\H^d(\wt\varphi_1(\wt W)) \leq \H^d(\wt W) 
\leq \H^d(\wt E \cap \wt B)
\leq (1+\rho(r_0))^d \H^d(E \cap B) \leq C r^d \leq C \wt r^d
\leqno (28.19)
$$
because $E$ is locally Ahlfors-regular in $B_0$
(if $h(r_0) < \varepsilon_n$ is small enough). Then (28.18) yields
$$\eqalign{
\H^d(\wt W) 
&\leq (1+\rho(r_0))^{2d} \H^d(\wt\varphi_1(\wt W))
+ (1+\rho(r_0))^d h(r) r^d 
\cr&
\leq \H^d(\wt\varphi_1(\wt W)) + [C\rho(r_0) + h(r)] \, \wt r^d
}\leqno (28.20)
$$
because $r = (1+\rho(r_0)) \wt r$; (28.17) follows because 
$h(r) \leq h(r_0) < \varepsilon_n$ and $\rho(r_0) \leq \varepsilon_n$.

When the ball $B$ is contained in $B(0,t)$ for some
$t \leq 8r_0/9$, we can use the better bilipschitz control provided
by (28.13), and the proof of (28.20) yields
$$
\H^d(\wt W) 
\leq \H^d(\wt\varphi_1(\wt W)) + [C\rho(9t/8) + h(9t/8)] \, \wt r^d.
\leqno (28.21)
$$
Consequently,
$$
\wt E \cap B(0,t) \in GSAQ(B(0,t),1,2r_0,C\rho(9t/8) + h(9t/8)).
\leqno (28.22)
$$
We may now follow our proof of Theorem 28.7, which we want to apply
to $\wt E$.
For this proof we need to know that $\wt E$ is rectifiable and
locally Ahlfors regular in $B_1$. Since $\xi$ is bilipschitz, we just need
to check that $E$ is rectifiable and locally Ahlfors regular in $B_0$,
and this follows from our assumption (28.10), together with the fact
that $h(r_0) \leq \varepsilon_n$, and as usual 
Propositions 4.1 and 4.74 and Theorem 5.16.
We even get bounds on $r^{-d}\H^d(\wt E \cap B(x,r))$,
for $x\in E$ and $B(x,2r) \i B_1$, that depend only on
$n$ and $\Lambda$, as long as $\varepsilon_n$ is chosen small enough.
Incidentally, it is a little easier to proceed this way here, rather than
trying to use (28.17) directly, because this way we don't need
to worry about the fact that on $B_1$, the boundary pieces
$\wt L_j = B_1 \cap \xi(L_j)$ may not satisfy the Lipschitz assumption
exactly as it was stated.

Then we turn to the main ingredient of the proof, which is 
the comparison argument in [D5], Lemma 5.1,  
that we already used for Theorem 28.4.
This lemma uses a radial deformation of the set (here, this means $\wt E$). 
The fact that the boundary sets $\xi(L_j)$ are conical allows us to 
use the same deformation, and then apply (28.20) or (28.21).
That is, we do not need to use the full almost minimality
of $\wt E$, because we just need to compare with a single competitor
that lives in a small ball $B$ centered at the origin.
Now this is the only place where we use the almost minimality of $E$
in [D5] (see (5.4) there); 
after this, the same argument as in Theorem 28.7 applies, and 
again does not use the full almost minimality (or the Lipschitz assumption),
but just the rectifiability and local Ahlfors regularity of $\wt E$, plus
measure theory. We get that $H(r)$ is a nondecreasing function of 
$r \in (0,8r_0/3)$, where 
$$
H(r) = r^{-d} \H^d(\wt E \cap B(0,r))
 \exp \alpha \Big(\int_0^r (C \rho(9t/4)+ h(9t/4)) {dt \over t} \Big)
\leqno (28.23)
$$
is the analogue for $\wt E$ of the function in (28.8).
Of course there is a small difference between $\Psi$ and $H$,
because $\H^d(E \cap \xi^{-1}(B(0,r))$ is not the same
as $\H^d(\wt E \cap (B(0,r))$ (one set is the image of the
other one by $\xi$), but this will not disturb much for the monotonicity. 
Once again, we want to avoid a computation, so let us recall
how (28.23) is obtained in [D5].  
We write $H(r) = l(r) g(r)$, with $l(r) = \H^d(\wt E \cap B(0,r))$
and $g(r) = r^{-d} e^{\alpha A(r)}$, with 
$A(r) = \int_0^r (C \rho(9t/4)+ h(9t/4)) {dt \over t}$.
Then we find out that those functions have a derivative almost 
everywhere, and that it is enough to show that $H'(r) \geq 0$ almost 
everywhere. A computation shows that
$$
H'(r)= l'(r) g(r) + l(r) g'(r) = g(r) \big\{ l'(r) 
- { l(r) \over r} \,\big[d - \alpha (C \rho(9r/4) + h(9r/4)) \big]\big\}
\leqno (28.24)
$$
(see (5.27) in [D5]),  
and it tuns out that because of Lemma 5.1, the right-hand side
of (28.24) is nonnegative almost everywhere. That is, setting
$a = C \rho(9r/4) + h(9r/4)$,
$$
l'(r) \geq { l(r) \over r} (d-\alpha a)
\leqno (28.25)
$$

Now we want to replace $l(r)$ with $l_1(r) = \H^d(E\cap \wt B(r))
= \H^d(\xi^{-1}(\wt E \cap B(0,r))$, and prove the same inequality,
but with $\alpha$ replaced by a larger $\alpha_1$.
We observe that, by a change of variable,
$l_1(r) \leq (1+\rho(9r/4))^d l(r)$ and 
$l'_1(r) \geq  (1+\rho(9r/4))^{-d} l'(r)$. Therefore
$$\eqalign{
l_1'(r) &\geq  (1+\rho(9r/4))^{-d} l'(r)
\geq (1+\rho(9r/4))^{-d} { l(r) \over r} (d-\alpha a)
\cr&
\geq (1+\rho(9r/4))^{-2d} \, { l_1(r) \over r} \, (d-\alpha a).
}\leqno (28.26)
$$
Now $\rho(9r/4) < a$ and $a$ is as small as we want,
so $(1+\rho(9r/4))^{-2d} (d-\alpha a) \geq d-\alpha_1 a$
for some new constant $\alpha_1 > \alpha$,
and (28.26) gives an analogue of (28.25) for $l_1$ and $\alpha_1$,
which implies that $\Psi'(r) \geq 0$ almost everywhere, and then
that $\Psi$ is monotone, by the same argument as in [D5].  
\qed

\msi{\bf Remark 28.27.}
When $E$ is an $A_+$-almost minimal set (as in (28.10)),
we do not need to assume that $E$ is coral and that $x_0 \in E$.
The reason is the same as for Remark 28.9.

\msi
{\bf 29. Minimal sets with constant density are cones} 
\ms

Our goal for this section is to prove that under the assumptions
of Theorem 28.4, if in addition the density function $\theta$ is 
constant on some interval, the set $E$ (almost) coincides with
a minimal cone on the corresponding annulus. As a corollary we will
get that under mild assumptions, blow-up limits of coral almost minimal sets 
are minimal cones. See Corollary~29.53. 

\ms\proclaim Theorem 29.1.
Let the open set $U$ and the boundary pieces $L_j$, $0 \leq j \leq 
j_{max}$ satisfy the Lipschitz assumption, let $E$ be a coral local 
minimal set in $U$, with $E \in GSAQ(U,1,\delta,0)$ for some $\delta > 0$,
and let $B(x_0,r_0) \i U$, with $0 < r_0 \leq \delta$ be given.
Assume (as in (28.2)) that each $L_j$ coincides in $B(x_0,r_0)$ 
with a cone centered at $x_0$ and that we can find 
constants $a$, $b$, and $\theta$ such that $0 \leq a < b \leq r_0$ and
$$
\H^d(E\cap B(x_0,r)) = \theta r^d
\ \hbox{ for } a < r < b.
\leqno (29.2)
$$
Then there is a closed coral minimal cone $\cal C$ centered at $x_0$,
such that $E \cap B(x_0,b) \sm B(x_0,a) = {\cal C}\cap B(x_0,b) \sm 
B(x_0,a)$.

\ms
Let $\wh L_j$ denote the cone that coincides with $L_j$ in 
$B(x_0,r_0)$; the fact that $\cal C$ is a minimal cone in $\R^n$, associated
to the boundary pieces $\wh L_j$, is a fairly easy consequence of 
local minimality of $E$. The argument is given in more detail in
[D5], pages 125-126, but we sketch it here for the  
convenience of the reader. In the ball $B(x_0,r_0)$, 
the cone $\cal C$ is a competitor for $E$, i.e., can be obtained
as $\varphi_1(E)$ for some family $\{ \varphi_t \}$ of mappings
that contract part of $E$ along the rays through $x_0$. In 
particular, the constraint (1.7) is satisfied by (28.2) and 
because we move points along rays.
The cone also has the same measure as $E$ in $B(x_0,r_0)$, by our assumption 
of constant density. Then every competitor for the cone gives rise
to a competitor for $E$ (by scale invariance, we may assume that the
modifications only occur in $B(x_0,r_0/2)$, and then we just compose 
our two deformations); the minimality  of $E$ then implies the 
minimality of $\cal C$. Similarly, $\cal C$ is coral because $E$ is
Ahlfors-regular.

We required $E$ to be coral to have a cleaner statement; if we 
don't, we just get that conclusion that the difference between 
$E \cap B(x_0,b) \sm B(x_0,a)$ and ${\cal C} \cap B(x_0,b) \sm B(x_0,a)$
is $\H^d$-negligible. Recall that the core of a minimal set is a 
coral minimal set; we would use this to replace an initial minimal
cone $\cal C$ with a coral one.

\ms
Let us now prove the theorem. Just for convenience and by translation invariance, 
let us assume that $x_0 = 0$; notice that we work with the Lipschitz assumption, 
so the origin does not have a special position in our grid, except for 
the fact that it is the center of our cones. 

Set $A = B(0,b) \sm \overline B(0,a)$. We first follow carefully our proof of 
monotonicity for $\theta$ (Theorem 28.4), use the fact that all the
inequalities are identities almost everywhere, and get that for 
$\H^d$-almost every $x\in E\cap A$, $E$ has a tangent plane $P(x)$ at $x$, 
which goes through the origin. The existence of a tangent plane 
follows from the rectifiability and local Ahlfors regularity of $E$;
only the fact that $0 \in P(x)$ is new, and we get it essentially because 
otherwise the measure of $E$ in a thin annulus that contains $x$ 
would be too large, compared to $\H^{d-1}(\d B(0,|x|)$. 
See (6.5) in [D5] (and its translation one page later);  
the same proof (of measure theory only) applies here.

Our next stage is to show that for $\H^d$-almost every $y\in E\cap A$,
$$
E \hbox{ contains the line interval } A \cap L(y),
\leqno (29.3)
$$
where $L(y) = \big\{ \lambda y \, ; \, \lambda > 0 \big\}$
is the open half line through $y$. Once we prove this, the conclusion 
will follow, because (29.3) then also holds for all $y\in E \cap A$;
see [D5] (below (6.12)) for the easy verification.  

We can thus restrict our attention to the points $y\in E\cap A$
for which the tangent plane $P(y)$ exists and contains the origin,
but we also add the following density constraint, which is valid $\H^d$-almost 
everywhere (see [Ma], Theorem 6.2 (2) on page 89). 
For $y\in U$, denote by ${\cal F}(y)$ the collection of all the 
faces $F$ of our grid that contain $y$. We require that for 
every face $F \in {\cal F}(y)$, $y$ be a density point of 
$E\cap F$ in $F$, i.e., that
$$
\lim_{r \to 0} r^{-d} \H^d(B(y,r) \cap E \sm F) = 0.
\leqno (29.4)
$$

Let $y\in E$ be such a point, and assume that we can find
$x\in A \cap L(y)\sm E$; we want to apply the proof of 
Proposition 6.11 in [D5],  
with a few modifications, to get a contradiction.

The construction will use two radii $r_y$ and $r_x$, 
and will work as soon as $r_y$ is small enough (depending on 
$y$) and then $r_x$ is small enough (depending on $r_y$ and 
the position of $x$, $y$, and in particular on the ratio
$|y|/|x|$, which may be large). Various smallness conditions
will arise along the proof, but let us mention the first ones.
First set
$$
B_y = B(y,r_y), \ B_x = B(x, r_x), \ \hbox{ and } P = P(y)
\leqno (29.5)
$$
to simplify the notation.
As in (6.13) in [D5], we require that  
for some small $\varepsilon_0$ (that will be chosen near the end),
$$
\dist(z,P) \leq \varepsilon_0 r_y
\ \hbox{ for } z \in E \cap B(y,3r_y);
\leqno (29.6)
$$
this is true for $r_y$ small because $P = P(y)$ is a true tangent plane,
by Ahlfors-regularity 
(see Exercise 41.21 on page 277 of [D4]). 
We also demand that for each face $F$ that
contains $y$, and in particular the smallest one,
$$
\H^d(B(y,3r_y) \cap E \sm F) \leq \varepsilon_1^d r_y^d,
\leqno (29.7)
$$
where $\varepsilon_1$ is another, even smaller, positive constant.
Again small $r_y$ satisfy this, by (29.4).
Let us deduce from this that
$$
\dist(z,L_j) \leq C \varepsilon_0 r_y
\ \hbox{ when $L_j$ contains $y$ and } z \in P \cap B_y.
\leqno (29.8)
$$
We want to apply Lemma 9.14 to $E$ and the ball $B_y$. If $r_y$ is 
small enough, the first condition (9.15) on the size of $B_y$ is 
satisfied. Also, every $L_i$ that meets $3B_y$ contains $y$, 
so the set $L$ of (9.16) is the intersection of all the $L_i$ that contain $y$. 
The smallest face $F$ that contains $y$ is contained in all these $L_i$, 
hence $F \i L$ and (29.7) says that  
$\H^d(B(y,3r_y) \cap E \sm L) \leq \varepsilon_1^d r_y^d$.
Since $E$ is locally Ahlfors-regular, this implies that
$$
\dist(z,L) \leq C \varepsilon_1 r_y
\ \hbox{ for } z\in E \cap 2B_y.
\leqno (29.9)
$$
That is, the constraint (9.17) is satisfied if $\varepsilon_1$
is small enough.

Next, the flatness condition (9.18) holds (for the same $P$
and with $\varepsilon = \varepsilon_0$); if $\varepsilon_0$ is small enough, 
we can apply Lemma 9.14 and we also get that
$$
\dist(p,E) \leq \varepsilon_0 r_y
\ \hbox{ for } p\in P \cap B(y,3r_y/2),
\leqno (29.10)
$$
as in (9.19). Now (29.8) follows from this and (29.9), if $\varepsilon_1$ is 
smaller than $\varepsilon_0$.

Let us try to describe the construction of [D5],  
without entering into too much detail. For simplicity, all the
references of the type (6.x) will refer to the corresponding
number in [D5]. 

The construction of [D5] 
starts with the choice of a set $T\i B_x$, which is defined near 
(6.19); the fact that $T \i B_x$ follows from 
Lemma 6.15 (of [D5]). 
Then we consider the cone $\wh T$ over $T$, and more precisely the 
part that lives near $[x,y]$. That is, if $x_1$ and $x_2$ lie in
the half line $L(y)$ through $y$, we denote by $V(x_1,x_2)$ the
set of points of $\R^n$ whose orthogonal projection on the line
that contains $L(y)$ lies between $x_1$ and $x_2$.
We shall use a lot the piece of tube
$$
T_0 = \wh T \cap V(x_1,y)
\leqno (29.11)
$$
where $x_1$ is a point of $L(y)$ that lies quite close to $x$ 
(on the other side of $x$ as $y$, so that $V(x,y) \i V(x_1,y)$);
see (6.22) for the definition of $x_1$, which was called $x_0$
in [D5]  
(but we want to avoid a conflict with the center of our main ball). 
Notice that
$$
\dist(z,[x,y]) \leq (1+|x|^{-1}|y|) \, r_x
\ \hbox{ for } z\in T_0,
\leqno (29.12)
$$
essentially because $T \i B_x$. We shall use $T_0$
to connect $B_x$ to $B_y$; the point of the specific
choice in Lemma 6.15 is to find a tube, inside $T_0$,
that does not meet $T_0$, but we shall not need to know this here.
Now set
$$
Z = T_0 \cup \overline B_y,
\leqno (29.13)
$$
as in (6.43); this is the place where most of the transformations
will take place. That is, a few successive mappings $f_1, \ldots, f_5$
are constructed in [D5], and the  
composition $f = f_5 \circ f_4 \circ f_3 \circ f_2 \circ f_1$
is used to define a competitor for $E$, and eventually get a 
contradiction. The point of these mappings is to use the hole 
that we have in $B_x$ (that is, there is no $E$ there) 
to push $E \cap B_y$ onto a set of much smaller
measure, essentially contained in $\d B_y$.
We will not need to know here how these mappings are constructed, 
but for the information of the reader, let us rapidly
say how it works. 
Our first mapping $f_1$ uses the hole in $\wh T$ to send points of $E$ 
(vertically, i.e., in hyperplanes perpendicular to $L(y)$)
to something that looks like a thin double tunnel (say, when $d=2$ and 
$n=3$), with a common flat floor that is contained in $P$,
and which we shall use to communicate between $B_x$
(which does not meet $E$) and $P \cap B_y$
(which we want to kill). Then $f_2$ acts in $B_y$, 
where it pushes points of $E_1 = f_1(E)$ to $B_y \cap P$
(a large disk, with an entrance that was the floor of the tunnel),
plus some small part of $\d B_y$ near $P$.
Then we compose with a mapping $f_3$ that pushes the points of the 
floor, starting from the empty part $P \cap B_x$, and eventually sends 
every point of the floor to a point of $P \cap \d B_y$. This is 
good, because this is how we get rid of most of the measure of $E \cap B_y$.
In codimension larger than $1$, since $\d B_y$ has an infinite 
measure, we need to compose with two additional Federer-Fleming projections 
on $d$-dimensional skeletons, one near the boundary of our tunnel 
that is not the floor, and one near $\d B_y$, to get some 
good control on the image of the part of 
$E_3 = f_3 \circ f_2 \circ f_1(E)$ that lives there.

We return to the final mapping $f$, and record some of its
properties. There is a slight enlargement $Z^\ast$ of $Z$, 
that will be described in a minute, such that
$$
f(z) = z \hbox{ for $z\in E \sm Z^\ast$, and } f(Z^\ast) \i Z^\ast.
\leqno (29.14)
$$
This set is defined just above (6.90), and we just need to know 
that
$$
Z^\ast = Z \cup Z_1^\ast \cap Z_2^\ast
\leqno (29.15)
$$
(see just above (6.90)) with sets $Z_1^\ast$ and $Z_2^\ast$
with the following properties. First
$$
\dist(z,T_0) \leq C r_x
\ \hbox{ for } z \in Z_1^\ast;
\leqno (29.16)
$$
see above (6.79). Next $Z_2^\ast$ is defined as a 
neighborhood, roughly of width $Cr_x$ (see (6.58))
of a set $Z_2$, which itself is defined by (6.51) and contained 
in $T_0$ (because $V(x,y) \i V(x_1,y)$). And similarly
$$
\dist(z,\d B_y \cap P)
\leq 3\varepsilon_0 r_y + C r_x 
\ \hbox{ for } z \in Z_2^\ast;
\leqno (29.17)
$$
this time see the line above (6.88) for the definition
of $Z_2^\ast$ in terms of $Z_2$, (6.51) for the definition
of $Z_2$, and (6.39) for $H$.

The main nice thing about $f$ is that, as in (6.92),
$$
\H^d(E_5 \cap Z^\ast) \leq C \varepsilon_0 r_y^d + C_{x,y} r_x^{d-1},
\leqno (29.18)
$$
where $E_5 = f(E)$, and the constant $C_{x,y}$ does not 
depend on $r_x$, so that we can choose $r_x$ very small 
to make $\H^d(E_5 \cap Z^\ast)$ as small as we want 
compared to $r_y^d$.

It is easy to find a continuous one parameter family of functions
$\varphi_t : E \to \R^n$, such that $\varphi_0(z) = z$ and 
$\varphi_1(z) = f(z)$ for $z\in E$,
$\varphi_t(z) = z$ for $z\in E \sm Z^\ast$, and 
$\varphi_t(Z^\ast) \i Z^\ast$.
This is not exactly what was done in [D5],  
where a brutal linear interpolation was enough, but did not yield 
$\varphi_t(Z^\ast) \i Z^\ast$ because $Z^\ast$ is not convex.
Here we want to proceed with just a little more care, and the
most natural thing to do is construct a family that goes from 
the identity to $f_1$, then from $f_1$ to $f_2 \circ f_1$, 
then to $f_3 \circ f_2 \circ f_1$, and so on. 
For the first two times, we proceed by brutal linear interpolation.
When we go from $f_2 \circ f_1$ to $f_3 \circ f_2 \circ f_1$,
we proceed slightly differently.

The main part of the definition of $f_3$ is (6.47), which is its 
definition on the floor $F = P \cap Z \cap V(x,b)$, where
$b$ is a point of $L(y)$ a little further from $x$ than $B_y$;
see the bottom of page 114 in [D5],  
and the statement of Theorem 6.2 for the definition of $b$
(which essentially lies on $\d B(0,r_0)$ here). We don't care
about the definition of $f_3$ on the rest of the interior of $Z$,
because the intersection of $E_2 = f_2 \circ f_1(E)$ with 
the interior of $Z$ is reduced to $F$; see (6.46), and 
compare with (6.43) or see the comment three lines above (6.52) 
for a hint. And on $\overline{\R^n \sm Z}$, we set $f_3(z) = z$
(see (6.48)). So for the interior of $Z$, the only interesting
piece is the definition of $f_3$ on the floor $F$, and it is
obtained by conjugating with a bilipschitz mapping $\psi$
that goes from $F$ to a cylinder a radial mapping on the cylinder
that maps it to its boundary. To define the intermediate mappings
on $F$, we just interpolate linearly the radial projection, and
conjugate. Of course we keep the identity on $\overline{\R^n \sm Z}$,
and the way we extend to the rest of $Z$ does not matter anyway.

Even though this is not important (and a brutal interpolation would
do), for the last two segments where we compose with $f_4$ and $f_5$,
it is more natural to decompose the deformation into successive 
Federer-Fleming (radial) projections, and for each one interpolate 
linearly. This way we are sure not to leave the sets $Z_1^\ast$
and $Z_2^\ast$ introduced in (29.15).

Anyway, this allows us to define a one parameter family 
$\{ \varphi_t \}$ that goes from the identity to $f$, but
because of the boundary condition (1.7) we cannot use this
family directly, and we shall compose it with retractions.

We shall use the mapping $\Pi$ that was constructed for Lemma 17.18,
except that here we can work under the Lipschitz assumption, 
in which case we conjugate it with our usual mapping 
$\psi(\lambda \cdot)$ to make it work on a different grid. 
This only makes the constant $C$ in (17.19)-(17.21)
larger, and also forces us to define $\Pi(z,s)$ only when
$0 \leq s \leq C^{-1}$ (instead of $0 \leq s \leq 10^{-1}$),
but this will not matter.
We define a first part of our path by
$$
g_t(z) = \Pi(z,t \chi(z))
\ \hbox{ for $z\in E$ and } 0 \leq t \leq 1,
\leqno (29.19)
$$
where the cut-off function $\chi = \chi_1 + \chi_2$ is defined as 
follows. We start with $\chi_1$, for which we use the small scale
$\tau_1 = C_1 r_x$, where the geometric $C_1$ will
be chosen later and may also depend on $|y|/|x|$, and set
$$
\chi_1(z) = \Big[2\tau_1 - \dist(z,T_0)\Big]_+ 
\ \hbox{ for } z\in \R^n.
\leqno (29.20)
$$
Here $T_0$ is the truncated cone of (29.11)
and $a_+$ denotes the positive part of $a\in \R$.
Similarly, we choose $\tau_2 = C_2 \varepsilon_0 r_y$,
where the large $C_2$ will also be chosen later, and we set
$$
\chi_2(z) = \Big[2\tau_2 - \dist(z,B_y)\Big]_+ 
\ \hbox{ for } z\in\R^n.
\leqno (29.21)
$$
Thus, by (17.19), $g_t(z) = z$ for $0 \leq t \leq 1$
unless $z$ lies in a $2\tau_1$-neighborhood of $T_0$,
or a $2\tau_2$-neighborhood of $B_y$. At the end of this first stage,
we are left with $g_1(z) = \Pi(z,\chi(z))$.
For our next stage, we use the $\varphi_t$ above and set
$$
g_t(z) = \Pi(\varphi_{t-1}(z),\chi(\varphi_{t-1}(z)))
\ \hbox{ for $z\in E$ and } 1 \leq t \leq 2.
\leqno (29.22)
$$
We want to check now that the $g_{2t}$, $0 \leq t \leq 1$, define an
acceptable competitor for $E$. We start with (2.4). Set 
$$
Z_+ = \big\{ z\in E \, ; \, \chi(z) > 0 \big\}
= \big\{ z\in E \, ; \, \dist(z,T_0) < 2 \tau_1
\hbox{ or } \dist(z,B_y) < 2 \tau_2 \big\}.
\leqno (29.23)
$$
If $C_1$ and $C_2$ are large enough, $Z_+$ contains $Z^\ast$, 
(by (19.13), (29.16), and (29.17)).
If $z\in E \sm Z_+$, (29.14) and the definition of the intermediate
mappings say that $\varphi_s(z) = z$ for all $s$, 
then $\chi(\varphi_{s}(z)) = \chi(z) = 0$,
and  by (17.19) $\Pi(\varphi_{s}(z),\chi(\varphi_{s}(z))
= \varphi_{s}(z) = z$ for all $s$; thus (29.19) and (29.22) yield 
$$
g_t(z) = z \ \hbox{ when $z\in E \sm Z_+$ and $0 \leq t \leq 2$.}
\leqno (29.24)
$$

Next suppose that $z \in Z_+ \sm Z^\ast$. Since $\varphi_s(z) = z$
for $0 \leq s \leq 1$ (again by (29.14) and the definition of the intermediate
mappings), we get that $g_t(z) = \Pi(z,s)$ for some $s\in [0,\chi(z)]$. 
Since $\chi(z) \leq 2 \tau_1 + 2 \tau_2$, (17.19) yields
$$
\dist(g_t(z), Z_+) \leq |g_t(z)-z| \leq C(\tau_1+\tau_2) 
\leq C(r_x+ \varepsilon_0 r_y).
\leqno (29.25)
$$

Finally, assume that $z \in Z^\ast$; then all the $\varphi_s(z)$ 
lie in $Z^\ast$ (by (29.14)), and in this case 
$$
\dist(g_t(z),Z^\ast) 
\leq |g_t(z)-\varphi_s(z)|\leq C(\tau_1+\tau_2) \leq C(r_x+\varepsilon_0 r_y)
\leqno (29.26)
$$
for some $s$, by the proof of (29.25).

Thus, if $r_y$ and then $r_x$ are chosen small enough, $z$ and
the $g_t(z)$ lie in a compact subset of $B(x_0,r_0)$ when
$z\in Z_+$; (2.4) follows, and also (1.5) and (1.6).
Here we can take for $B$ a compact ball that is almost as wide as
$B(x_0,r_0)$, and we do not care if its radius is quite large 
(provided that it stays smaller than $r_0$), because 
there is no price to pay in (2.5) when $B$ is large, 
since $E$ is minimal.
The constraints (1.4) and (1.8) (continuity and Lipschitzness)
hold by construction, so we are left with (1.7) to check.

Let $z \in E \cap L_j$ be given. We may assume that $z\in Z_+$, 
because otherwise $g_t(z) = z \in L_j$ for all $t$.
If $z\in Z_+ \sm Z^\ast$, (29.14) says that $\varphi_s(z) = z$
for all $s$, and then (29.19) or (22.22) says that every
$g_t(z)$ is of the form $\Pi(z,s)$, which by (17.20)
lies in any face of $L_j$ that contains $z$. We are left
with the case when $z\in Z^\ast$.

Even in this case, (29.19) and (17.20) say that 
$g_t(z) = \Pi(z,t\chi(z)) \in L_j$ for $0 \leq t \leq 1$,
so it is enough to show that
$$
g_t(z) \in L_j \hbox{ when $z \in L_j \cap Z^\ast$ and } t > 1.
\leqno (29.27)
$$
Let $t > 1$ be given, and set $w = \varphi_{t-1}(z)$; thus 
$$
g_t(z) = \Pi(w, \chi(w))
\leqno (29.28)
$$ 
by (29.22). As we shall see, the mapping $\Pi$ tends to send to $L_j$
points that lie sufficiently close to $L_j$, so we want to show that
$w$ lies close to $L_j$. 
Let us first check that
$$
\dist(\xi, L_j) \leq C (1+ |y|/||x|) r_x
\ \hbox{ for } \xi \in T_0.
\leqno (29.29)
$$ 
Set $z' = z |y|/|z|$; first observe that $z' \in L_j$
because $z\in L_j$ and by the cone property (28.2).
Also, $|z'-y| \leq C r_y + C (1+ |y|/||x|) r_x$, 
because $z\in Z^\ast \i Z_+$, and by (29.15), (29.13),
(29.12) and, for $Z_1^\ast$ and $Z_2^\ast$,
(29.16) and (29.17). If $r_x$ and $r_y$ are chosen small
enough, this forces $y\in L_j$. That is, we choose $r_x$ and $r_y$
so small that the ball centered at $y$ and with radius 
$C r_y + C (1+ |y|/||x|) r_x$ does not meet any $L_j$ that does
not already contain $y$. Then $[x,y] \i L_j$ too,
by the cone property (28.2), and now (29.29) follows from
(29.12). Also recall from (29.8) that
$$
\dist(\xi, L_j) \leq C \varepsilon_0 r_y
\ \hbox{ for } \xi \in P \cap B_y.
\leqno (29.30)
$$
Let us return to $w = \varphi_{t-1}(z)$ and use this
to evaluate its distance to $L_j$. We will unfortunately
need to distinguish between cases.
First assume that 
$$
\dist(w,T_1) \leq C(1+ |y|/||x|) r_x,
\leqno (29.31)
$$
where $T_1 = \wh T \cap V(x_1,y_0)$ is defined like  $T_0$ in (20.11), 
but with a point $y_0$ that lies on $L(y)$, but at distance $r_y/5$ 
from $y$, in the direction opposite to $x$. Thus $T_1$ is a little larger
than $T_0$ (in the direction of $y$), but not much.
By (29.29), $\dist(w,L_j) \leq C(1+ |y|/||x|) r_x$ too.
If $\tau_1$ is large enough (compared to $(1+ |y|/||x|) r_x$), this implies 
that $\chi(w) \geq \chi_1(w) \geq \tau_1$ (by (29.20)), and that 
$\tau_1$ is much larger than $\dist(w,L_j)$. Then (17.20)
(applied to any face of $L_j$ that lies near $w$) implies 
that $g_t(z) = \pi(w,\chi(w)) \in L_j$, by (29.28), and as needed.
Similarly, if 
$$
\dist(w,P\cap B_y) \leq C \varepsilon_0 r_y,
\leqno (29.32)
$$
(29.30) says that a similar estimate holds for
$\dist(w,L_j)$; then, if $\tau_2$ is large enough compared
to $\varepsilon_0 r_y$, and by (29.21), 
$\chi(z) \geq \chi_2(z) \geq \tau_2$ and (17.20) implies 
that $g_t(z) = \pi(w,\chi(w)) \in L_j$.

Finally, if $w = z$, we know that $g_t(z) = \pi(z,\chi(z)) \in L_j$,
because $z\in L_j$ and by (17.20) again.

We now want to check that we always fall in one 
of these three cases. We start with the case when $t$
comes from the first part of the construction of $f$
in [D5], when we  
go from the identity to the first map $f_1$ that is
defined on pages 107-110 of [D5]. 
Let us say that we define these intermediate functions $\varphi_t$
by linear interpolation, i.e., set
$\varphi_t(z) = 2t f_1(z) +(1-2t) z$, for $0 \leq t \leq 1/2$.
Also recall that $f_1$ was obtained by moving points in vertical
hyperplanes (i.e., in directions perpendicular to the line $L(y)$
through $x$ and $y$), and inside the tube $T_1$. That is,
$f_1(z) = z$ for $z \in E \sm T_0$ (by (6.25) in [D5]), 
and $f_1(T_0) \i T_0$ by (6.26) 
(also compare our definition of $y_0$ with (6.22)). 
The same thing holds for $\varphi_t(z)$,
$0 \leq t \leq 1/2$, and we are happy because $w=z$ or (29.31) holds.

At the end of this first stage, all the points $z\in E$ are sent
to $E_1 = f_1(E)$, and we now apply to them our second mapping $f_2$,
defined on pages 111-114 of [D5]. 
This map only moves points of $B_y$
(see (6.37)), and moves them like a radial projection, centered
on a $(n-d-1)$-dimensional sphere $S_y$, onto a part of $P \cap \d B_y$ 
(recall that $Q$ on pages 106 and 111 of [D5] 
is $(n-d)$-dimensional). We do not care about the details here, 
we just need to know that if $\xi \in E_1 \i \overline B_y$ 
(the only place where we may move something) then by (6.30),
$$
\dist(\xi,P) \leq 2 \varepsilon_0 |\xi - y| \leq 2 \varepsilon_0 r_y.
\leqno (29.33)
$$
The function $f_2$ maps $\overline B_y$ to itself (see above (6.38))
and maps $E_1 \cap \overline B_y$ to a 
$2 \varepsilon_0 r_y$-neighborhood of $P$, by (6.41) and the definition
(6.39). Again we interpolate linearly, i.e. set
$\varphi_t(z) = (4t-2) f_2 \circ f_1(z) +(3-4t) f_1(z)$ 
for $1/2 \leq t \leq 3/4$, and are happy because (29.32) holds as soon as 
$w \neq f_1(z)$. (In the other case, we already new that $w = f_1(z)$ 
satisfies (29.31)).

We now consider the case when $\varphi_t$ comes between 
$f_2 \circ f_1$ and $f_3 \circ f_2 \circ f_1$, where $f_3$
is described on pages 114-116 of [D5]. 
This is the place where we said below (29.18) that we do not
interpolate linearly, but only after a conjugation with a biLipschitz
mapping. There is only one part of $E_2 = f_2 \circ f_1(E)$
where $f_3$ moves points, which is the floor $F \i P$.
See our discussion below (29.32), or directly (6.47), (6.48), and 
(6.46) in [D5]. 
But we choose the $\varphi_t$ so that they move points of $F$ along $F$,
which means that (29.31) or (29.32) holds. Thus we are happy, up to 
the stage of $f_3 \circ f_2 \circ f_1$. 
As was explained below (6.52) in [D5],  
in codimension $1$ we could stop here, but in higher codimensions
we have to complete our construction with Federer-Fleming projections
that act near $Z$. For these we can interpolate linearly between
each elementary projection and the (composition with the) next one.
By construction we never leave the sets $Z_1^\ast \cup Z_2^\ast$ of 
(29.15) (see (6.58), the proof of (6.78) (the intermediate projection 
also stay close), (6.81), and the proof of (6.96)); then the desired 
conclusion follows from (29.16) and (29.17), which show that (29.31)
or (29.32) holds whenever we move a point.

This completes the proof of the boundary constraint (1.7)
for our family $\{ g_t \}$. Thus we are allowed to use (2.5),
which says that
$$
\H^d(W) \leq \H^d(g_2(W))
\leqno (29.34)
$$
(recall that we work with minimizers here), with
$W = \big\{ z\in E \, ; \, g_2(z) \neq z \big\}$.
Recall that (29.22) yields
$g_2(z) = \Pi(f(z),\chi(f(z))$, because $\varphi_1 = f$.
Also, $\Pi(w,\chi(w))$ is a $C$-Lipschitz function of $w$,
by (29.20), (29.21), and (17.21), and with a constant $C$
that may now depend on the bilipschitz constant $\Lambda$
when we work under the Lipschitz assumption, but not on
$\varepsilon_0$ or $r_x$, for instance.
Because of this, we can easily take care of $g_2(E \cap Z^\ast)$, 
since (6.95) and (6.97) in [D5] yield  
$$
\H^d(g_2(E \cap Z^\ast)) \leq C \H^d(f(E \cap Z^\ast))
\leq C \H^d(E_5 \cap Z^\ast) 
\leq C \varepsilon_0 r_y^d + C_{x,y} r^{d-1}_x;
\leqno (29.35)
$$
note that the $r_y^2$ in (6.97) is a misprint, but that
we really get $r_x^{d-1}$ because the points of $T_0$ stay 
$C_{x,y}$-close to a line.

Because of (29.24), $W$ is contained in $Z_+$. So we still need
to worry about the set $W_0 = W \sm Z^\ast \i Z_+ \sm Z^\ast$. 
On this set $f(z)=z$ by (29.14), so $g_2(z) = \Pi(z,\chi(z)) \neq z$,
where the last part holds by definition of $W$. 
Now we follow the construction of $\Pi(z,\chi(z))$, and notice that
by (17.24)
$$
\Pi(z,\chi(z))
= \Pi_{0,s_0}\circ\Pi_{1,s_1}\cdots\circ\Pi_{n-1,s_{n-1}}(z)
\leqno (29.36)
$$
where the $\Pi_{k,s_k}$ come from Lemma 17.1 and we set
$$
s_m = (6C)^{-m} \chi(z) 
\ \hbox{ for } 0 \leq m \leq n-1,
\leqno (29.37)
$$
as in (17.23). Also denote by $x_k$, $n+1 \geq k \geq 0$, 
the successive images of $z= x_{n+1}$, defined  as (above (17.26) and)
in (17.26) by $x_k = \Pi_{k,s_k}(x_{k+1})$. By the property
(17.2) of $\Pi_{k,s_k}$, we see that
$x_k = x_{k+1}$ unless $0 < \dist(x_{k+1},{\cal S}_k) < 2s_k$.
That is, unless $x_k$ lies very close to some face of the
grid, without actually lying on that face. By an easy induction, 
we see that $\Pi(z,\chi(z)) = z$ unless
$$
0 < \dist(z, F) \leq 2 \chi(z)
\ \hbox{ for some face $F$ of our grid.}
\leqno (29.38)
$$
First assume that $z\in 2B_y$. By choosing $r_y$ small enough,
we can ensure that $2B_y$ only meets the faces $F$
that already contain $y$. But (29.38) means that 
$z \in E \cap 2B_y \sm F$. Then of course
$z \in E \cap 2B_y \sm F_0$, where $F_0$ is the smallest face that
contains $y$, and by (29.7)
$$
\H^d(W_0 \cap 2B_y) \leq \H^d(E \cap 2B_y \sm F_0) \leq 
\varepsilon_1^d r_y^d.
\leqno (29.39)
$$
By (29.23), we are left with 
$$
W_0 \sm 2B_y \i E \cap Z_+ \sm 2B_y
\i \big\{ z\in E \, ; \, \dist(z,T_0) < 2 \tau_1 \big\}.
\leqno (29.40)
$$
Because $T_0$ stays so close to $[x,y]$ (see (29.12)), we can cover 
$W_0 \sm 2B_y$ by less than $C_{x,y} r_x^{-1}$ balls centered on $E$
and with radius $Cr_x$. Then by local Ahlfors regularity
(Propositions~4.1 and 4.74), 
$\H^d(W_0 \cap 2B_y) \leq C_{x,y} r_x^{d-1}$.
Again $g_2(z) = \Pi(z,\chi(z))$ is a $C$-Lipschitz function of 
$z\in W_0$, so
$$
\H^d(g_2(W_0)) \leq C \H^d(W_0) 
\leq C \varepsilon_1^d r_y^d + C_{x,y} r_x^{d-1}.
\leqno (29.41)
$$
We add this to (29.35) and get that
$$
\H^d(g_2(W)) \leq C \varepsilon_0 r_y^d
+ C \varepsilon_1^d r_y^d + C_{x,y} r_x^{d-1}.
\leqno (29.42)
$$
On the other hand, we claim that $W$ is large because
it contains most of $E \cap B_y$. More precisely,
if $z\in E \cap B_y \sm W$, then $g_2(z) = z$
and hence $z \in g_2(E\cap B_y) \i g_2(E \cap Z^\ast)$ 
(by (29.13) and because $Z \i Z^\ast$). Thus
$$
\H^d(E \cap B_y \sm W) \leq \H^d(g_2(E \cap Z^\ast))
\leq C \varepsilon_0 r_y^d + C_{x,y} r^{d-1}_x
\leqno (29.43)
$$
by (29.35). But $\H^d(E \cap B_y) \geq C^{-1} r_y^d$
because $E$ is locally Ahlfors regular (by Propositions~4.1 and 4.74),
hence
$$
\H^d(W) \geq \H^d(E \cap B_y) - \H^d(E \cap B_y \sm W)
\geq C^{-1} r_y^d - C \varepsilon_0 r_y^d - C_{x,y} r^{d-1}_x.
\leqno (29.44)
$$
If $\varepsilon_0$ and $\varepsilon_1$ are small enough and 
$r_x$ is small enough 
(depending also on the position of $x$ and $y$ through $|y|/|x|$),
(29.42) and (29.44) contradict (29.34), and this concludes our proof 
of Theorem 29.1.
\qed

\ms
We complete this section with a simple consequence of Theorem 29.1
and the results of Section 10 on limits, revised in Section 24 so that they
apply to blow-up limits.

Let us list the assumptions for the next theorem. Most of them are
the same as for our Theorem 24.13 on blow-up limits, which we intend 
to combine with Theorem 29.1.
We are given a coral almost minimal set $E$ in the open set $U$,
an origin $x_0 \in E$, and a sequence $\{ r_k \}$, with
$$
\lim_{k \to +\infty} r_k = 0.
\leqno (29.45)
$$
We assume that
$$
\hbox{$U$ and the $L_j$ satisfy the Lipschitz assumption,}
\leqno (29.46)
$$
as in Definition 2.27, and that
$$
\hbox{the configuration of $L_j$ is flat at $x_0$, 
along the sequence $\{ r_k \}$.}
\leqno (29.47)
$$
See Definition 24.8, but also recall that we have a simpler condition,
the flatness of the faces of the $L_j$ along the sequence, which is 
introduced in Definition 24.29 and implies it; see Proposition 24.35.
Recall that (29.47) comes with a collection of limit sets $L_j^0$,
the natural blow-up limits of the $L_j$, that are defined by (24.7).
We assume that
$$
\hbox{the $L_j$ satisfy  (10.7) or (19.36),}
\leqno (29.48)
$$
the additional assumptions that we used for our theorems on limits,
and that 
$$\eqalign{
&\hbox{$E$ is a coral almost minimal set in $U$, with }
\cr& \hskip.2cm
\hbox{ sliding conditions coming from the $L_j$,}
}\leqno (29.49)
$$
and with a gauge function $h$ such that 
$$
\lim_{r \to 0} h(r) = 0.
\leqno (29.50)
$$
For this, we accept the three types ($A_+$, $A$, or $A'$) 
of almost minimality; see Definition~20.2. Also see Definition 3.1
for corality.

We also give ourselves a closed set $E_\infty \i U$,
and we assume that
$$
E_\infty = \lim_{k \to +\infty} r_k^{-1} (E- x_0) 
\ \hbox{ locally in } \R^n
\leqno (29.51)
$$
(see near (10.5) for the definition).
Finally, we suppose that the following limit exists:
$$
\theta(x_0) =\lim_{\rho \to 0} \rho^{-d} \H^d(E\cap B(x_0,\rho)).
\leqno (29.52)
$$
Notice that the existence of $\theta(x_0)$ follows from Theorem 28.7
when $h$ satisfies the Dini condition (28.5) and the $L_j$ are cones.
One may also use Remark 28.11 to prove this when the $L_j$ are almost 
cones (but with conditions stronger than (29.47)).

\ms\proclaim Corollary 29.53. 
Let the coral almost minimal set $E$ in $U$, the point $x_0 \in E$, 
the sequence $\{ r_k \}$, and the set $E_\infty$ satisfy the 
conditions above. Then $E_\infty$ is a coral minimal cone, with 
the sliding boundary conditions defined by the sets $L^0_j$, 
$0 \leq j \leq j_{max}$, defined by (24.7), and with the constant density
$$
r^{-d} \H^d(E_\infty \cap B(x_0,r)) = \theta(x_0). 
\leqno (29.54)
$$

\ms
The first assumptions allow us to apply Theorem 24.13,
which says that $E_\infty$ is a coral minimal set, with 
the sliding boundary conditions defined by the sets $L^0_j$.
Recall that this is defined by (24.16) or (24.17), as the reader 
prefers.

We still need to check that $E_\infty$ is a cone and that (29.54)
holds, and naturally we start with (29.54). 
For this shall need to say more about how Theorem 24.13 is proved.
We consider the same sets $E_k$ (compare with (24.3)), and the main 
point of the proof is to show that Theorem 23.8 can be applied.
A long first part consists in showing that for any fixed
large radius $R \geq 1$, the $L_j^0$ satisfy the 
Lipschitz assumption on some appropriate domain $U_R$
(defined by (24.18)-(24.20)). Once this is done, we apply 
Theorem 23.8 to the domains $U_{R,k} = \xi_k(U_R)$ and 
the sets $E_k\cap U_{R,k}$. 
In turn Theorem 23.8 consists in applying Theorem 10.8
to a single domain (with single boundary sets), but the
different sets $\wt E_k = \xi_k^{-1}(E_k)$, where the 
bilipschitz mappings $\xi_k$ come from the condition (29.47),
and satisfy the asymptotic conditions (23.3) and (23.4).
Eventually, one proves that these sets satisfy the assumptions
of Theorem 10.8 (see (23.16)-(23.20)). 
Their limit is still
$E_\infty$ (by (23.3) and mostly (23.4)), and Theorem 10.97,
which has the same assumptions as Theorem 10.8, shows that for 
$0 < \rho < \rho_1 < R/2$,
$$\eqalign{
\H^d(E_\infty\cap B(x_0,\rho)) 
&\leq \liminf_{k \to +\infty} \H^d(\wt E_k\cap B(x_0,\rho))
\leq \liminf_{k \to +\infty} \H^d(E_k\cap B(x_0,\rho_1))
\cr&
= \liminf_{k \to +\infty} r_k^{-d}\H^d(E\cap B(x_0, r_k\rho_1))
= \rho_1^d \theta(x_0)
}\leqno (29.55)
$$
where we used the asymptotic bilipschitz property (23.3) for the
change of variable to control the measures, and then the scale 
invariance and (29.56). Since we may take $\rho_1$ as close to $\rho$ 
as we want, this gives the upper bound in (29.54).

Similarly, we can apply Lemma 22.2 (whose assumptions are the 
same as for Theorem~10.8) and with any choice of $M > 1$ and 
$h> 0$; this yields, for $0 < \rho_1 < \rho_2 < \rho < R/2$,
$$\eqalign{
(1+Ch) M \H^d(E_\infty\cap B(x_0,\rho)) 
& \geq (1+Ch) M \H^d(E_\infty\cap \overline B(x_0,\rho_2)) 
\cr&
\geq \limsup_{k \to +\infty} \H^d(\wt E_k\cap \overline B(x_0,\rho_2))
\cr& 
\geq \limsup_{k \to +\infty} \H^d(E_k\cap B(x_0,\rho_1))
\cr& 
= \limsup_{k \to +\infty} r_k^{-d}\H^d(E\cap B(x_0, r_k\rho_1))
= \rho_1^d \theta(x_0)
}\leqno (29.56)
$$
by (22.4) and the same sort of computation as above.
This yields the rest of (29.54).

Once we know (29.54), we can apply Theorem 29.1, we get that
$E_\infty$ is a cone, and this completes the proof of Corollary 29.53.
\qed

\msi
{\bf 30. Nearly constant density and approximation by minimal cones} 
\ms
In this  section we use Theorem 29.1 and the results of 
Section 10 on limits to give sufficient conditions, in terms of density,
for an almost minimal set to be very close to a minimal cone.

For the main statement, we give ourselves a fixed ball $B_0 = B(x_0,r_0)$,
boundary pieces $L_j^0$, $0 \leq j \leq j_{max}$,
and we suppose that the Lipschitz assumption of Definition 2.7 holds
for $B_0$ and the $L_j^0$ and (in the non rigid case), 
that the $L_j^0$ satisfy the technical assumption (10.7), 
or the weaker (19.36). 
We also suppose that $0 \in B_0$ and that, for $0 \leq j \leq j_{max}$,
$$
L_j^0 \hbox{ coincides with a cone in } B_0.
\leqno (30.1) 
$$
We see $B_0$ and the collection of $L_j^0$ as a model for
domains $U$, endowed with boundary pieces $L_j$, 
and such that there is a  bilipschitz mapping $\xi$ such that 
$$
\xi(B_0) = U \ \hbox{ and } \ 
L_j  = \xi(L_j^0)
\hbox{ for } 0 \leq j \leq j_{max}.
\leqno (30.2)
$$
We want to say that when the bilipschitz constant of $\xi$ is
small, $E$ is a coral quasiminimal set in $U$ with constants
$M$ close enough to $1$ and $h$ small enough,
and the density ratios of $E$ in two 
different balls centered at $\xi(0)$ are close enough, then 
$E$ looks a lot like a minimal cone in the corresponding annulus. 
The statement will be a little complicated, but later on, 
in Proposition 30.19, we shall consider the simpler case when the annulus 
is just a ball centered at the origin.

\ms\proclaim Proposition 30.3. Let $B_0$ and the 
$L_j^0$, $0 \leq j \leq j_{max}$, be as above.
In particular, assume that we have (30.1), the Lipschitz assumption, 
and (10.7) or (19.36). 
For each $\tau > 0$, we can find $\varepsilon > 0$ such that the 
following holds.
Let $U$, the $L_j$, and $\xi$ satisfy (30.2), and assume that
$$
\xi \hbox{ is $(1+\varepsilon)$-bilipschitz.}
\leqno (30.4)
$$
Let $E$ be a coral quasiminimal set in $U$, with
$$
E \in GSAQ(U,1+\varepsilon,2r_0,\varepsilon),
\leqno (30.5)
$$
set $x_0 = \xi(0) \in U$, and assume that
for some choice of radii $0 < r_1 < r_2 \leq r_0$,
$$
r_2^{-d} \H^d(E \cap B(x_0,r_2)) 
\leq r_1^{-d} \H^d(E \cap B(x_0,r_1)) + \varepsilon.
\leqno (30.6)
$$
Then there is a minimal cone $T$ centered at $x_0$ such that
$$
\dist(y,T) \leq \tau r_0 \hbox{ for } 
y\in E \cap B(x_0,r_2-\tau) \sm B(x_0,r_1+\tau)
\leqno (30.7)
$$
$$
\dist(z,E) \leq \tau r_0 \hbox{ for } 
z\in T \cap B(x_0,r_2-\tau) \sm B(x_0,r_1+\tau),
\leqno (30.8)
$$
$$\eqalign{
\big| \H^d(T \cap B(y,t)) - \H^d(E \cap B(y,t)) \big| &\leq \tau r_0^d
\cr
\hbox{ for $y$ and $t$ such that }
B(y,t) \i B(x_0,r_2-\tau) &\sm B(x_0,r_1+\tau),
}\leqno (30.9)
$$
and 
$$
\big| H^d(E \cap B(x_0,r)) - H^d(T \cap B(x_0,r)) \big|
\leq \tau r_0^d
\ \hbox{ for } r_1+\tau \leq r \leq r_2-\tau.
\leqno (30.10)
$$

\ms
Let us comment this statement before we prove it.

Proposition 30.3 is a generalization of Proposition 7.1 
in [D5]. 

In (30.6), it could happen that $B(x_0,r_0))$ goes slightly out of $U$, 
so we could have written $\H^d(E \cap U \cap B(x_0,r_2))$ 
instead of $\H^d(E \cap B(x_0,r_2))$ to be more explicit
(but the result is the same since $E \i U$).

Of course (30.6) is only meaningful if $r_1$ is not too close to 
$r_2$, but otherwise the conclusion is empty anyway.

We can always apply the result to domains $U' \supset U$
and quasiminimal sets $E'$ in $U'$, since it is easy to check
that $E'\cap U \in GSAQ(U,1+\varepsilon,2r_0,\varepsilon)$
as soon as $E' \in GSAQ(U',1+\varepsilon,2r_0,\varepsilon)$.

We are lucky because we don't need the Dini condition in (28.5), 
or even the fact that $E$ is almost minimal. Thus we may not be
in the situation where we know for sure that the density 
$r^{-d} \H^d(E\cap B(x,r))$ is almost nondecreasing 
(as in Theorems 28.7 and 28.15).
On the other hand, we only state a result at a fixed scale, 
not an asymptotic result, and we will not compute $\varepsilon$
in terms of $\tau$, but just apply a compactness argument.

The author is not too happy about the statement
of Proposition 30.3, because we let $\varepsilon$
depend on the $L_j$. This is because the proof below,
just like the limiting result of Section 23, for instance, 
uses the very fine bilipschitz convergence on domains,
for which it seems too hard to extract converging subsequences. 
Both for Section 23 and here, there are probably ways to improve the 
statement, but the author is not really sure of what would 
be needed, and hopes that in practice Proposition 30.3 will often
be enough. Anyway, we still put the scale invariant factors $r_0$ and
$r_0^d$ in (30.7)-(30.10), even though our statement allows 
$\varepsilon$ to depend on $r_0$ through $B_0$ and the $L_j$.

\ms
We shall prove the proposition by compactness. It will be easier to 
take limits of sets like $\wt E = \xi^{-1}(E)$, because they live in
the fixed domain $B_0$.
A simple computation, using (30.4) and (30.5), 
shows that if $E$ is as in the statement,
$$
\wt E \in GSAQ(B_0,1+C\varepsilon,r_0,C\varepsilon),
\leqno (30.11)
$$
where the associated boundary pieces on $B_0$ are still the $L_j^0$.

Now let us fix $B_0$, the $L_j^0$, and $\tau$, and
assume that we cannot find $\varepsilon > 0$ as in the statement.
Let $\xi_k$, $U_k=\xi_k(B_0)$, the sets $L_j^k$, $E_k$, and the radii
$r_{1,k}$ and $r_{2,k}$ provide a counterexample, 
associated with $\varepsilon_k = 2^{-k}$. 
By translation invariance, we may assume that $\xi_k(0) = 0$.
Also set $\wt E_k = \xi_k^{-1}(E_k)$.

Recall that $\xi_k$ is defined on $B_0$ and 
$(1+2^{-k})$-bilipschitz (by (30.4)). Since $\xi_k(0) = 0$ for all $k$, 
we see that modulo extracting a subsequence, we may assume that the $\xi_k$ 
converge, uniformly on $B_0$, to a mapping $\eta$ which in addition is 
$1$-bilipschitz. That is, $\eta$ is the restriction of a linear isometry of 
$\R^n$ that fixes $0$. Let us replace $\xi_k$ with $\eta^{-1}\circ \xi_k$,
$E_k$ with $\eta^{-1}(E_k)$, and so on; we get a new counterexample for 
which $\eta$ is the identity. So we may assume that the $\xi_k$ converge,
uniformly on $B_0$, to the identity. Because of our bilipschitz
property (30.4), the $\xi_k^{-1}$ also converge, uniformly on compact
subsets of $B_0$, to the identity.

Modulo extracting a new subsequence, we may assume that the sets
$\wt E_k$ converge in $B_0$ to a limit $F$; then
the sets $E_k = \xi_k(\wt E_k)$ converge,
locally in $B_0$, to the same limit $F$.
By (30.11), $\wt E_k \in GSAQ(B_0,1+C 2^{-k},r_0,C2^{-k})$.
Then for each small $\delta > 0$, we can apply 
Theorem 10.8 to the (end of the) sequence $\{ \wt E_k \}$, in the
domain $B_0$. We get that $F$ is a coral quasiminimal set (associated
to the boundary pieces $L_j^0$), with $F\in GSAQ(B_0,1+\delta,r_0,\delta)$.
Since this holds for every $\delta > 0$, $F$ is a minimal set in 
$B_0$.

Now we want to take care of the measures. 
Let $B$ be an open ball, with $B \i B(0,r_0-\tau/3)$ 
(we may assume that $\tau << r_0$, so we don't lose much). 
By the lower semicontinuity property (10.98),
$$
\H^d(F \cap B) \leq \liminf_{k \to +\infty} \H^d(\wt E_k \cap B).
\leqno (30.12)
$$
By Lemma 22.3 (applied with $h$ and $M-1$ as small as we want), we 
also get that
$$
\H^d(F \cap \overline B) \geq \limsup_{k \to +\infty} \H^d(\wt E_k 
\cap \overline B).
\leqno (30.13)
$$

We want to compare this with our density assumption (30.6).
Let us replace our sequence with a subsequence for which $r_{1,k}$ 
tends to a limit $r_1$ and $r_{2,k}$ tends to a limit $r_2$. 
Notice that $r_{2,k} \geq r_{1,k} + 2\tau$ for all $k$ 
(otherwise the conclusion (30.7)-(30.10) would be trivially true,
which is impossible for a counterexample), so $r_2 \geq r_1 + 2\tau$.
Then take $\delta > 0$ very small, and notice that for $k$ large,
$$
\wt E_k \cap B(0,r_2-\delta) = \xi_k^{-1}(E_k) \cap B(0,r_2-\delta))
\i \xi_k^{-1}(E_k \cap B(0,r_2))
\leqno (30.14)
$$
(because (30.4) says that $\xi_k$ is $(1+2^{-k})$-bilipschitz).
Then apply (30.12) to $B = B(0,r_2-\delta)$, and get that
$$\eqalign{
\H^d(F \cap B(0,r_2-\delta))
&\leq \liminf_{k \to +\infty} \H^d(\wt E_k \cap B(0,r_2-\delta))
\cr&
\leq \liminf_{k \to +\infty} \H^d(\xi_k^{-1}(E_k \cap B(0,r_2)))
\cr&
\leq \liminf_{k \to +\infty} (1+2^{-k})^d \H^d(E_k \cap B(0,r_2))
\cr&
= \liminf_{k \to +\infty} \H^d(E_k \cap B(0,r_2))
\cr&
\leq (r_1/r_2)^{-d} \liminf_{k \to +\infty} \H^d(E_k \cap B(0,r_1))
}\leqno (30.15)
$$
by (30.4) and because (30.6) holds for $E_k$, with $\varepsilon = 2^{-k}$.
Similarly,
$$
\wt E_k \cap \overline B(0,r_1+\delta) 
= \xi_k^{-1}(E_k) \cap \overline B(0,r_1+\delta))
\supset \xi_k^{-1}(E_k \cap B(0,r_1))
\leqno (30.16)
$$
for $k$ large, so when we apply (30.13) to $B= B(0,r_1+\delta)$,
we get that
$$\eqalign{
\H^d(F \cap \overline B(0,r_1+\delta))
&\geq \limsup_{k \to +\infty} \H^d(\wt E_k \cap \overline B(0,r_1+\delta))
\cr&
\geq \limsup_{k \to +\infty} \H^d(\xi_k^{-1}(E_k \cap B(0,r_1)))
\cr&
\geq \limsup_{k \to +\infty} (1+2^{-k})^{-d} \,\H^d(E_k \cap B(0,r_1))
\cr&
= \limsup_{k \to +\infty} \H^d(E_k \cap B(0,r_1)).
}\leqno (30.17)
$$
We compare (30.15) and (30.17), let $\delta$ tend to $0$, and get that
$$
\H^d(F \cap B(0,r_2)) \leq (r_1/r_2)^{-d} \H^d(F \cap \overline B(0,r_1)).
\leqno (30.18)
$$
By Theorem 28.4 (and because $F$ is minimal in the full $B_0$)
$\theta(r) = r^{-d} \H^d(F\cap B(0,r))$ is nondecreasing on 
$(0,r_0)$. But (30.18) says that $\theta(r_2) \leq \lim_{r \to r_1^+} \theta(r)$,
so $\theta$ is constant on $(r_1,r_0)$.
We apply Theorem 29.1 and get that $F$ coincides, in the
annulus $B(0,r_2) \sm \overline B(0,r_1)$, with a coral minimal cone $T$.

We shall now prove that the approximation properties (30.7)-(30.10)
are satisfied for $k$ large (and the cone $T$ that we just found), 
and this will give the desired contradiction with the definition of $E_k$.

First notice that (30.7) and (30.8) hold, because we observed earlier that
$F$ is also the limit of the $E_k$ in compact subsets of $B_0$. For (30.9) and 
(30.10), we deduce them from (30.12) and (30.13).
The details of the verification were done in
[D5], pages 128-129, 
so we refer to that and merely mention
the two minor difficulties that may worry the reader.
For (30.10), it is easy to deduce it, for a single radius, from 
(30.12), (30.13), and the fact that for the cone $T$, 
$\H^d(T \cap \d B(0,r)) = 0$. But it is enough to check 
(30.10) for a finite collection of radii, because
$\H^d(T \cap B(0,r))$ is a continuous function of $r$, while each 
$\H^d(E_k \cap B(0,r)) = 0$ is nondecreasing. For (30.9), we can 
proceed similarly, but we also need the less obvious fact that
$\H^d(T\cap \d B(y,t)) = 0$ for every ball $B(y,t)$. This is
(7.14) in [D5], 
and the verification, done as Lemma 7.34 [D5], 
only uses the fact that $\H^d(T \cap B(0,1)) < +\infty$ and
a little bit of geometric measure theory, but not the fact that
$T$ is minimal (which is good, because here minimal is merely meant
with additional boundary constraints, so our cone $T$ is probably not minimal
as in [D5]).  
This concludes our proof of Proposition 30.3 by contradiction 
and compactness.
\qed

\ms
Let us now state the analogue of Proposition 30.3
for the density in a ball (i.e., with $r_1=0$).  

\ms\proclaim Proposition 30.19. Let $0 < r_0$ be given, and
let $B_0$ and the $L_j^0$, $0 \leq j \leq j_{max}$, be as in the statement
of Proposition 30.3. In particular, assume that we have (30.1), 
the Lipschitz assumption, and (10.7) or (19.36). 
For each $\tau > 0$, we can find $\varepsilon > 0$ such that the 
following holds.
Let $\xi$, $U$, and the $L_j$ satisfy (30.2), and assume that
$$
\xi \hbox{ is $(1+\varepsilon)$-bilipschitz.}
\leqno (30.20)
$$
Let $E \i U$ be a coral quasiminimal set, with
$$
E \in GSAQ(U,1+\varepsilon, 2r_0,\varepsilon),
\leqno (30.21)
$$
set $x_0 = \xi(0) \in U$, and assume that for some $r_2 \in (0,r_0]$,
$$
r_2^{-d} \H^d(E \cap B(x_0,r_2)) 
\leq \varepsilon + \inf_{0 < r < 10^{-3} r_0} r^{-d} \H^d(E \cap B(x_0,r)).
\leqno (30.22)
$$
Then there is a minimal cone $T$ centered at $x_0$ such that
$$
\dist(y,T) \leq \tau r_0 \hbox{ for } 
y\in E \cap B(x_0,r_2-\tau)
\leqno (30.23)
$$
$$
\dist(z,E) \leq \tau r_0 \hbox{ for } 
z\in T \cap B(x_0,r_2-\tau),
\leqno (30.24)
$$
$$\eqalign{
&\big| \H^d(T \cap B(y,t)) - \H^d(E \cap B(y,t)) \big| \leq \tau r_0^d
\cr&\hskip 1.5cm
\hbox{ for $y$ and $t$ such that } B(y,t) \i B(x_0,r_2-\tau),
}\leqno (30.25)
$$
and in particular
$$
\big| H^d(E \cap B(x_0,r)) - H^d(T \cap B(x_0,r)) \big|
\leq \tau r_0^d
\ \hbox{ for } 0 < r \leq r_2-\tau.
\leqno (30.26)
$$

\ms
This is now the generalization of Proposition 7.24 in [D5]. 
We write (30.22) in this strange way because, since we do not assume
$E$ to be almost minimal with a small gauge function, we do not know 
that $\lim_{r \to 0} r^{-d} \H^d(E \cap B(x_0,r))$ exists
and gives a good lower bound on $r^{-d} \H^d(E \cap B(x_0,r))$
for $r$ small. Of course, with suitable additional assumptions on the
almost minimality of $E$, we could use Theorem 28.15 and replace
the infimum in (30.22) with the density 
$\lim_{r \to 0} r^{-d} \H^d(E \cap B(x_0,r))$.

We repeat the proof of Proposition 30.3 because we don't want to 
worry about the way $\varepsilon$ depends on $r_1$, and also because
in (30.25) we allow balls $B(y,t)$ that contain $x_0$. So we suppose that
for $k \geq 0$, we have a counterexample $E_k$ to the statement 
with $\varepsilon = 2^{-k}$, extract suitable subsequences, and find
a minimal set $F$, which is the joint limit in $B_0$ of the sequences
$\{ E_k \}$ and $\{ \wt E_k \}$. 

For each small $r_1 \in (0,10^{-3} r_0)$, we can repeat the argument near
(3.12)-(3.18), and get (3.18). This is why we required an infimum in (30.22).
Then $\theta(r) = r^{-d} \H^d(F\cap B(0,r))$ is constant on 
$(r_1,r_2)$, and $F$ coincides with a minimal cone $T$ on 
$B(0,r_2) \sm \overline B(0,r_1)$. We let $r_1$ tend to $0$ and get that 
$F \cap B_0 = T \cap B_0$. 

The desired contradiction, i.e., the fact that (30.23)-(30.26)
actually hold for $k$ large, is then proved as before. In particular,
since now $F$ coincides with $T$ near the origin, we are allowed 
balls $B(y,t)$ that contain it.
Proposition 30.19 follows.
\qed

\msi
{\bf 31. Where should we go now?} 
\ms

There were two main reasons for the present paper.
The first one was to obtain some boundary regularity results for
quasiminimal sets and almost minimal sets. 

As far as quasiminimal sets are concerned, the author believes that
the results in Parts I-III are not so far from being optimal,
because of the bilipschitz invariance. That is, Lipschitz graphs are
quasiminimal, and uniformly rectifiable sets are not so far the
Lipschitz regularity. Of course it would be nice to know
that the strange dimensional condition (6.2) can be removed,
especially because this would mean that we found another
proof of uniform rectifiability than the very complicated stopping time
argument coming from [D1].  
Also, the uniform rectifiability result of Theorem 6.1 barely
contains more information than the fact that $E$ is locally
uniformly rectifiable away from the boundary pieces $L_j$,
plus the uniform rectifiability of the pieces themselves.

The situation is quite different for the almost minimal sets
(typically, minimizers of a functional like $\int_E f(x) d\H^d(x)$,
maybe plus some lower order terms, and where $f : \R^n \to [1,M]$
is continuous). For these sets, we expect much more regularity
than what we obtained so far. Sufficiently flat sets are a little 
easier to control, because of Allard's theorem [All], 
but we should not expect precise general results, because we know
that a general description is already hard away from the
boundary. Recall that J. Taylor [Ta]  
gave a very good local description of the $2$-dimensional 
almost minimal sets in $\R^3$. A similar, but already
much less precise description is available for $2$-dimensional 
almost minimal sets in $\R^n$ (see [D5,6]),  
and there are even some first descriptions of 
$3$-dimensional almost minimal sets in $\R^4$ near special
(but non flat) points ([Lu1,2]), 
but we expect a lot of trouble except in very small dimensions.
Maybe see [D7] for a rapid description. 

Because this is always a good way to start, a first step consists
in studying the blow-up limits of our minimal sets at a point
of the boundary, and Corollary 29.53 says that in the
reasonable situations, these are sliding minimal cones 
associated to conical boundary pieces. So it seems interesting
to study (find a list of) the minimal cones in some simple situations.
Even for $2$-dimensional minimal cones in $\R^4$, with no boundary
piece, the list of minimal cones is not known. 
It was recently shown [Li2] 
that the almost orthogonal union of two planes in $\R^4$ is minimal
(partially answering a conjecture of F. Morgan [Mo2]), 
and that the orthogonal product of two one dimensional sets $Y$
in $\R^4$ is minimal too [Li3],  
but there may be lots of other $2$-dimensional minimal cones in $\R^4$
that we did not guess.
To the author's best knowledge, the list of 
$2$-dimensional minimal cones in $\R^3$, with a unique
boundary $L_1$ which is a line, is not known either, and this would
be a very good start for some versions of the most classical Plateau
problem in $3$-space.

Once we have an almost minimal set $E$, with a known blow-up limit
at some point $x$ of the boundary, we can try to give a good 
description of $E$ near $x$, for instance a nice parameterization
by the blow-up limit, like for the J. Taylor theorem [Ta]. 
This should not be too hard, in a very limited number of
situations.

Of course it would be good to have a substitute for the monotonicity
of the density $\theta(x,r) = r^{-d} \H^d(E\cap B(x,r))$
when $x\in E$ lies close to a boundary piece, but not on it.
The monotonicity fails stupidly when $E$ is a half plane
(bounded by the line $L_1$) and $x\in E \sm L_1$, but
one may dream to use more clever functions of $r$, that would 
probably need to depend on the approximate shape of $E$
and the distance to $L_j$.

A second motivation is that a good local knowledge of the
sliding almost minimal sets near the boundary should help
attacking some existence problems. A typical one is the
version of Plateau's problem that was described in the introduction:
take a simple curve $\Gamma$ in $\R^3$ and an initial set
$E_0$ bounded by $\Gamma$. You do not have to know what this means, 
but if you pick a wrong $E_0$, the problem will probably have a trivial 
solution (like a point). Then try to minimize $\H^2(E)$
(or a variant) among all the competitors for $E$
(as defined in Definition~1.3). In [D3], the author 
proposed to use sequences of quasiminimal sets, together
with the concentration lemma of [DMS]  
and a construction of adapted polyhedral networks by V. Feuvrier
[Fv1], [Fv2] 
to find existence results for problem of this type.
Some existence results were indeed found
(see [Fv3], [Li1], [Fa]), 
but often avoiding complicated problems at the
boundary. For the problem above, for instance, it would be 
good to know that for the limit $E$ of the minimizing sequence 
that we construct, and which is a sliding minimal set
by Theorem 10.8, there is a Lipschitz retraction of a neighborhood
of $E$ onto $E$. Also see [D8] for variants of this problem, 
probably not all easy to solve. As was mentioned in the introduction,
since other categories (such as size minimizing currents) also yield
sliding minimal sets, boundary regularity results could be useful 
there too.

\bigskip
REFERENCES

\smallskip 
\item {[All]} W. K.  Allard, On the first variation of a varifold. 
Ann. of Math. (2) 95 (1972), 417-491.
\smallskip 
\item {[A1]} F. J. Almgren, Existence and regularity 
almost everywhere of solutions to elliptic variational problems 
among surfaces of varying topological type and singularity structure. 
Ann. of Math. (2) 87 1968 321-391.
\smallskip
\item {[A2]} F. J. Almgren, Existence and regularity almost everywhere 
of solutions to elliptic variational problems with constraints, 
Memoirs of the Amer. Math. Soc. 165, volume 4 (1976), i-199.
\smallskip 
\item {[A3]} F. J. Almgren, The structure of limit varifolds 
associated with minimizing sequences of mappings. 
Symposia Mathematica, Vol. XIV (Convegno di Teoria Geometrica 
dell'Integrazione e Variet\`a Minimali, INDAM, Rome, 1973), pp. 413-428. 
Academic Press, London, 1974.
\smallskip 
\item {[A4]} F. J. Almgren, Plateau's problem. 
An invitation to varifold geometry. Corrected reprint of the 1966 original. 
With forewords by Jean E. Taylor and Robert Gunning, and Hugo Rossi. 
Student Mathematical Library, 13. 
American Mathematical Society, Providence, RI, 2001. xvi+78 pp. 
\smallskip
\item {[DMS]} G. Dal Maso, J.-M. Morel, and S. Solimini, 
A variational method in image segmentation: Existence and
approximation results, Acta Math. 168 (1992), no. 1-2, 89--151.
\smallskip 
\item {[D1]} G. David, Morceaux de graphes lipschitziens et int\'{e}grales 
singuli\`{e}res sur une surface, 
Rev. Mat. Iberoamericana 4 (1988), no. 1, 73--114.
\smallskip 
\item {[D2]} G. David, Limits of Almgren-quasiminimal sets, 
Proceedings of the conference on Harmonic Analysis, 
Mount Holyoke, A.M.S. Contemporary Mathematics series, Vol. 320 
(2003), 119-145.
\item {[D3]} G. David, 
Quasiminimal sets for Hausdorff measures, 
in Recent Developments in Nonlinear PDEs,
Proceeding of the second symposium on analysis and PDEs
(June 7-10, 2004), Purdue University, 
D. Danielli editor, 81--99, Contemp. Math. 439,
Amer. Math. Soc., Providence, RI, 2007.
\smallskip
\item {[D4]} G. David, \underbar{Singular sets of minimizers for 
the Mumford-Shah functional},
Progress in Mathematics 233 (581p.), Birkh\"auser 2005.
\smallskip
\item {[D5]} G. David, 
H\"older regularity of two-dimensional Almost-minimal sets in $\Bbb R^n$,
Annales de la Facult\'{e} des Sciences de Toulouse, 
Vol 18, 1 (2009), 65--246.
\smallskip
\item {[D6]} G. David, 
$C^{1+a}$-regularity for two-dimensional almost-minimal sets 
in $\R^n$, Journal of Geometric Analysis 20 (2010),  no. 4, 837-954. 
\smallskip  
\item {[D7]} G. David, Regularity of minimal and almost minimal sets 
and cones: J. Taylor's theorem for beginners,
Centre de Recherches Math\'{e}matiques, 
Lecture notes of the 50th S\'{e}minaire de math\'{e}matiques 
sup\'{e}rieures, Montreal 2011, G. Dafni, R. Mc Cann, A. Stancu editors,
CRM Proceedings and lecture notes, Volume 56, 2013, 67-117. 
\smallskip  
\item {[D8]} G. David, Should we solve Plateau's problem again? 
To be published, Advances in Analysis: The Legacy of Elias M. Stein. 
Edited by Charles Fefferman, Alexandru D. Ionescu, D.H. Phong, 
and Stephen Wainger.
\smallskip
\item {[DS1]} G. David and S. Semmes, 
\underbar{Singular integrals and rectifiable sets in $\Bbb R\sp n$: 
au-del\`{a} des} \underbar{ graphes lipschitziens,}
Ast\'{e}risque No. 193 (1991).
\smallskip
\item {[DS2]} G. David and S. Semmes, 
Quantitative rectifiability and Lipschitz mappings, 
Trans. Amer. Math. Soc. 337 (1993), no. 2, 855--889.
\smallskip
\item {[DS3]} G. David and S. Semmes,
\underbar{Analysis of and on uniformly rectifiable sets},
Mathematical Surveys and Monographs 38, 
American Mathematical Society, Providence, RI, 1993.
\smallskip
\item {[DS4]} G. David and S. Semmes, Uniform rectifiability and 
quasiminimizing sets of arbitrary codimension, 
Memoirs of the A.M.S. Number 687, volume 144,  2000.
\smallskip
\item {[DeCL]} E. De Giorgi, M. Carriero, and A. Leaci, 
Existence theorem for a minimum problem
with free discontinuity set, Arch. Rational Mech. Anal. 108 (1989), 195-218.
\smallskip 
\item {[De]} T. De Pauw, Size minimizing surfaces. 
Ann. Sci. \'{E}c. Norm. Sup\'er. (4) 42 (2009), no. 1, 37-101.
\smallskip
\item {[Fa]} Yangqin Fang,  
Existence of Minimizers for the Riefenberg Plateau problem,
ArXiv, October 2013.
\smallskip
\item {[Fe]} H. Federer, \underbar{Geometric measure theory}, 
Grundlehren der Mathematishen Wissenschaf-ten 
153, Springer Verlag 1969.
\smallskip 
\item {[Fv1]} V. Feuvrier, Un r\'{e}sultat d'existence pour les 
ensembles minimaux par optimisation sur des grilles polyh\'{e}drales,
Th\`{e}se Paris-Sud (Orsay), September 2008.
\smallskip 
\item {[Fv2]} V. Feuvrier, Remplissage de l'espace euclidien 
par des complexes polyh\'edriques d'orientation impos\'ee et de rotondit\'e uniforme. 
(French. English, French summary) [Filling Euclidean space using polyhedral 
complexes of given orientation and uniform rotundity] 
Bull. Soc. Math. France 140 (2012), no. 2, 163-235.
\smallskip 
\item {[Fv3]} V. Feuvrier, 
Condensation of polyhedral structures onto soap films.
Preprint.
arXiv:0906.3505.
\smallskip 
\item {[HP]} J. Harrison and H. Pugh,
Existence and Soap Film Regularity of Solutions to Plateau's Problem.
Preprint.
ArXiv:1310.0508.
\smallskip 
\item {[LM]} G. Lawlor and F. Morgan, 
Curvy slicing proves that triple junctions locally minimize area.
J. Diff. Geom. 44 (1996), 514-528.
\smallskip 
\item {[Li1]} X. Liang, Topological minimal sets and existence results. 
Calc. Var. Partial Differential Equations 47 (2013), no. 3-4, 523-546.
\smallskip 
\item {[Li2]} X. Liang, 
Almgren-minimality of unions of two almost orthogonal planes in $R^4$. 
Proc. Lond. Math. Soc. (3) 106 (2013), no. 5, 1005-1059.
\item {[Li3]} X. Liang, 
Almgren and topological minimality for the set $\Bbb Y \times \Bbb Y$,
Preprint. ArXiv:1203.0564.
\smallskip 
\item {[Lu1]} T. Luu, 
H\"older regularity of three-dimensional minimal cones in $\R^4$. 
Preprint.
\smallskip 
\item {[Lu2]} T. Luu, 
On some properties of three-dimensional minimal sets in $\R^4$.
Preprint.
\smallskip
\item {[Ma]}  P. Mattila, \underbar{Geometry of sets and 
measures in Euclidean space}, Cambridge Studies in
Advanced Mathematics 44, Cambridge University Press l995.
\smallskip
\item {[Mo1]} F. Morgan, Size-minimizing rectifiable currents, 
Invent. Math. 96 (1989), no. 2, 333-348.
\smallskip
\item {[Mo2]} F. Morgan, 
Soap films and mathematics. 
Differential geometry: partial differential equations on manifolds (Los Angeles, CA, 1990), 375-380, 
Proc. Sympos. Pure Math., 54, Part 1, Amer. Math. Soc., Providence, RI, 1993.
\smallskip
\item {[Mo3]} F. Morgan, Geometric measure theory. A beginner's guide. 
Fourth edition. Elsevier/Academic Press, Amsterdam, 2009. viii+249 pp. 
\smallskip
\item {[R1]} E. R. Reifenberg, Solution of the Plateau Problem for 
$m$-dimensional surfaces of varying topological type,
Acta Math. 104, 1960, 1--92.
\smallskip
\item {[R2]} E. R. Reifenberg, An epiperimetric inequality related to 
the analyticity of minimal surfaces, Ann. of Math. (2) 80, 1964, 1--14.
\smallskip
\item {[St]}	E. M. Stein, \underbar{Singular integrals and 
differentiability properties of functions},
Princeton university press 1970.
\smallskip
\item {[Ta]} J. Taylor, The structure of singularities in 
soap-bubble-like and soap-film-like minimal surfaces, 
Ann. of Math. (2) 103 (1976), no. 3, 489--539.
\smallskip

\bigskip
\vfill \vfill \vfill\vfill
\noindent Guy David,  
\smallskip\noindent 
Math\'{e}matiques, B\^atiment 425,
\smallskip\noindent 
Universit\'{e} de Paris-Sud, 
\smallskip\noindent 
91405 Orsay Cedex, France
\smallskip\noindent 
guy.david@math.u-psud.fr

\bye